%% file: main.tex
\documentclass[preprint,11pt,authoryear,nonatbib]{elsarticle}

\makeatletter
\let\c@author\relax
\makeatother

\usepackage[utf8]{inputenc}
\usepackage[backend=biber,style=apa]{biblatex}
\usepackage{setspace}
\usepackage{titlesec}
\normalsize  

\titlespacing{\subsection}{0pt}{5pt}{5pt}
\titlespacing{\section}{0pt}{5pt}{5pt}
\titlespacing{\subsubsection}{0pt}{5pt}{5pt}

\usepackage{graphicx} 
\usepackage{amsmath}
\usepackage[left=0.94in, right=0.94in, top=0.94in, bottom=0.94in]{geometry}
\usepackage[toc,page]{appendix}
\usepackage{booktabs}
\usepackage{multirow}
\usepackage[final]{changes} 
\newcommand{\rev}[1]{\textcolor{black}{#1}}

\usepackage{comment}
\usepackage{float}
\usepackage{subcaption}
 \usepackage{adjustbox}
\usepackage{lscape}
\usepackage{todonotes}
\usepackage{makecell}
\usepackage{standalone}

\usepackage{amsthm}
\newtheorem{theorem}{Theorem}
\usepackage{enumitem}
\setlist[itemize]{parsep=0cm, itemsep=2pt, leftmargin=12pt} 
\usepackage{longtable}
\usepackage[hidelinks]{hyperref} 
\usepackage{cleveref}

\usepackage{xcolor}
\usepackage{color} 
\usepackage{colortbl}
\usepackage[linesnumbered,ruled,vlined]{algorithm2e}
\usepackage{booktabs}
\usepackage[normalem]{ulem}
\usepackage{ulem}
\usepackage{subfiles}
\usepackage{tikz}
\usepackage{pgfplots} 
\pgfplotsset{compat=1.18}

\usetikzlibrary{shapes.geometric, arrows, patterns,pgfplots.groupplots,pgfplots.statistics}
\tikzstyle{startstop} = [rectangle, rounded corners, 
minimum width=3cm, 
minimum height=1cm,
text centered, 
draw=black
]

\tikzstyle{io} = [trapezium, 
trapezium stretches=true, 
trapezium left angle=70, 
trapezium right angle=110, 
minimum width=3cm, 
minimum height=1cm, text centered, 
draw=black
]

\tikzstyle{process} = [rectangle, 
minimum width=3cm, 
minimum height=1cm, 
text centered, 
text width=3cm, 
draw=black
]

\tikzstyle{decision} = [diamond, 
minimum width=3cm, 
minimum height=1cm, 
text centered, 
draw=black
]

\tikzstyle{arrow} = [thick,->,>=stealth]

\allowdisplaybreaks

\SetCommentSty{mycommfont}

\usepackage{accents}

\usepackage{lscape}
\colorlet{darkgreen}{green!70!black}

\input{customCommands}

\begin{document}

\begin{frontmatter}
\author[a]{Siv Marie Cartland Hansen$^{*}$}
\author[b]{Rowan Hoogervorst}
\author[a]{Otto Anker Nielsen}
\author[a]{Richard Martin Lusby}
\author[c]{Evelien van der Hurk$^{*}$}

\cortext[cor1]{Corresponding author.}
\affiliation[a]{organization={Department of Technology, Management and Economics, Technical University of Denmark}, addressline={Akademivej 358}, city={Kongens Lyngby}, postcode={2800}, country={Denmark}}
\affiliation[b]{organization={EDHEC Business School}, addressline={24 avenue Gustave Delory}, postcode={59057}, city={Roubaix}, country={France}}
\affiliation[c]{organization={CGI Netherlands}, addressline={George Hintzenweg 89}, postcode={3068 AX}, city={Rotterdam}, country={Netherlands}}

\title{
An Exact Algorithm for Public Transport Line Planning Considering Passenger and Operational Costs and Lost Demand
}

\date{2025}

\begin{abstract}
\added{Line planning in public transport is the strategic problem of selecting lines and their operating frequencies. This problem is important as it defines the passenger service, based on available connections and expected travel times, and drives operational cost in terms of the number of vehicles required.}
This paper presents a line planning model that minimizes the weighted sum of passenger travel time, including in-vehicle time and frequency-dependent waiting and transfer times, and operating costs for the public transport agency. Unlike traditional approaches that assume demand to be fixed, our approach requires a minimum service level for demand to be captured, ensuring that services are provided only when they are attractive to users and cost-efficient to operate. \added{The introduced capacity constraints ensure sufficient capacity on the lines and help guide the trade-off between expected demand on selected lines and their frequencies.} The resulting mixed-integer program presents a challenging combinatorial problem as the number of passenger paths grows rapidly in relation to the number of lines and frequencies considered. To address this, we propose an iterative exact algorithm that utilizes a reduced problem representation and dynamically expands it with additional frequencies and paths. 
Evaluated on four networks with varying complexity and cost trade-offs, our method achieves significant speed-ups and tighter bounds compared to solving the complete model directly by CPLEX, particularly when operator and passenger costs are more evenly balanced in the objective. 
Furthermore, we demonstrate how accounting for lost demand leads to more efficient resource use from an overall perspective.
\end{abstract} 

\end{frontmatter}

\input{01IntroductionEJOR}
\input{02LitEJOR}
\input{03Problem}
\input{04MethodsEJOR}

\input{05CaseSetup}
\input{06Results}

\input{07Discussions}

\section*{Acknowledgment}
This research was funded by DFF (Independent Research Fund Denmark) as part of the project \emph{INTRAURBAN -- Integrating behavior modelling with operations research for user-driven public transport network design (Grant number 2033-00206B)}

\printbibliography

\appendix
\renewcommand{\thesection}{Appendix~\Alph{section}}
\input{08Appendix}

\newpage

\end{document}

%% file: customCommands.tex
\newcommand{\problem}{LPP\xspace}
\newcommand{\frameworkNameLP}{Dynamic Frequency Refinement Algorithm\xspace}

\newcommand{\elasticDemand}{service-dependent demand\xspace}

\newcommand{\grid}{\textsc{grid}\xspace} 
\newcommand{\sioux}{\textsc{sioux-falls}\xspace} 
\newcommand{\lowersaxony}{\textsc{lower-saxony}\xspace} 
\newcommand{\athens}{\textsc{athens}\xspace}

\newcommand{\linePool}{\mathcal{L}}
\newcommand{\frequencySet}{\mathcal{F}}
\newcommand{\headwaySet}{\mathcal{H}}

\newcommand{\vehicleSet}{\mathcal{V}}

\newcommand{\stkout}[1]{\ifmmode\text{\sout{\ensuremath{#1}}}\else\sout{#1}\fi}

\newcommand{\frameworkNameShortLP}{DFRA\xspace}

\newcommand{\stopSet}{\ensuremath{N}\xspace}
\newcommand{\stopSetCG}{\ensuremath{\mathcal{N}}\xspace}

\newcommand{\stopSetCGA}{\ensuremath{\mathcal{N}^{A}}\xspace}
\newcommand{\stopSetCGE}{\ensuremath{\mathcal{N}^{E}}\xspace}
\newcommand{\stopSetCGT}{\ensuremath{\mathcal{N}^{T}}\xspace}
\newcommand{\stopSetCGL}{\ensuremath{\mathcal{N}^{L}}\xspace}
\newcommand{\stopn}{\ensuremath{n}\xspace}
\newcommand{\linePT}{\ensuremath{l}\xspace}

\newcommand{\ODset}{\ensuremath{\mathcal{D}}\xspace}
\newcommand{\ODor}{\ensuremath{s_d}\xspace}
\newcommand{\ODde}{\ensuremath{t_d}\xspace}
\newcommand{\ODw}{\ensuremath{w_d}\xspace}
\newcommand{\od}{\ensuremath{d}\xspace}

\newcommand{\Pset}{\ensuremath{\mathcal{P}}\xspace}
\newcommand{\PsetOD}{\ensuremath{\mathcal{P}_d}\xspace}

\newcommand{\PsetODedge}{\ensuremath{\mathcal{P}_d(a)}\xspace}
\newcommand{\PSetLineFreq}{\ensuremath{\mathcal{P}_d(l,f)}\xspace}

\newcommand{\arc}{\ensuremath{a}\xspace}

\newcommand{\arcSet}{\ensuremath{\mathcal{A}}\xspace}
\newcommand{\arcSetDrive}{\ensuremath{\mathcal{A}^{IVT}}\xspace}
\newcommand{\arcSetTransfer}{\ensuremath{\mathcal{A}^{T}}\xspace}
\newcommand{\arcSetAcces}{\ensuremath{\mathcal{A}^{A}}\xspace}
\newcommand{\arcSetEgress}{\ensuremath{\mathcal{A}^{E}}\xspace}
\newcommand{\arcSetOptOut}{\ensuremath{\mathcal{A}^{O}}\xspace}
\newcommand{\linePoolSol}{\ensuremath{\mathcal{L}^{*}}\xspace}

\newcommand{\altModePath}{\ensuremath{p'}\xspace}
\newcommand{\minIVTd}{\ensuremath{t_\od^{\min}}\xspace}
\newcommand{\vlpVI}{\ensuremath{v_{lp}^{*}}\xspace}
\newcommand{\wttThresOD}{\ensuremath{\bar{c}_{d}}\xspace}
\newcommand{\kSim}{\ensuremath{\kappa}\xspace}

\newcommand{\frequencySetAlg}[1]{\ensuremath{\mathcal{F}^{#1}}\xspace}
\newcommand{\vectorSol}[1]{\ensuremath{\left( \mathbf{y}^*, \mathbf{x}^*, \mathbf{z}^*\right)^{(#1)}}\xspace}
\newcommand{\headwaySetAlg}[1]{\ensuremath{\mathcal{H}^{#1}}\xspace}
\newcommand{\PsetAlg}[1]{\ensuremath{\mathcal{P}^{#1}}\xspace}

%% file: 01IntroductionEJOR.tex
\section{Introduction}

This paper studies the Line Planning Problem (\problem) in public transport, which involves \added{the strategic level problem of} selecting which routes (lines) to operate and determining their operating frequency.
A line defines the sequence of stops served by a vehicle, thus determining the available connections and the required transfers for passengers.
The frequency of a line, expressed as the number of services per hour, influences expected waiting and transfer times. 
It is assumed that a predefined set of candidate lines (a line pool) and their possible frequencies are available based on characteristics of the physical infrastructure.
Together, the selected lines and frequencies define a \emph{line concept} as introduced in \textcite{schobel_line_2012}, which determines the expected travel times between stops and provides an estimate of operating costs and capacity via the number of vehicles required.
Although exact travel times and operating costs depend on later planning stages, particularly timetabling and vehicle scheduling, incorporating frequency setting in the line planning phase provides a better indication of passenger service, especially since operating frequency is a key factor in passenger satisfaction and ridership \autocite{ingvardson_how_2018}.

This article presents a passenger-oriented line planning model that minimizes a weighted sum of operator costs, including ticket fare revenue, and weighted travel time, accounting for in-vehicle time, as well as frequency-dependent waiting and transfer times. 
\added{A distinctive feature of this work is that it proposes an exact solution framework for a line planning problem where:}
    \begin{itemize}
        \item \added{the frequency of a line is chosen from a large discrete set of possible frequencies, thus also enforcing a minimum frequency per line in case the line is opened.}
        \item \added{ passenger travel time depends on the selected frequency of all lines in their selected travel path(s) in terms of initial waiting time, in-vehicle time, and frequency-dependent transfer time, where multiple potential travel paths are generated in advance.}
        \item \added{passenger assignment is guided by the objective to minimize overall weighted travel time. However, passenger preferences are accounted for by limiting the maximum weighted travel time they are willing to accept and pre-generating the set of candidate geographical paths. When no suitable path is available, the demand is lost. }
        \item \added{capacity constraints based on vehicle type, frequency, and line are included, as well as a budget constraint limiting the operational resources. 
        These are intended to guide the trade-off between network coverage and provided frequency, and do not intend to represent the results of overcrowding. In fact, overcrowding is a feature that one may wish to prevent as much as possible at the strategic level. }
    \end{itemize}
    \added{Input to this line planning problem thus consists of a line pool and candidate frequencies per line, an operational budget (possibly per vehicle type), a demand estimate (e.g., a peak-hour demand), and a set of candidate geographical paths (stop and line-based) for each group of passengers that adhere to a maximum weighted travel time for this group of passengers.}

The \added{proposed} model is an adaptation of that in \textcite{hurk_shuttle_2016}, which is well-suited to include weighted travel time, 
but was not originally scalable to general line planning. 
Due to the combinatorial nature of the problem, where the number of passenger paths grows rapidly with the number of lines and frequencies considered, the resulting mixed-integer linear program (MILP) can become computationally intractable if solved directly by a commercial solver. 
To address this, we propose the \frameworkNameLP (\frameworkNameShortLP), which utilizes a reduced network representation and dynamically expands this representation. 
In particular, this algorithm solves a MILP based on a reduced network representation containing a subset of the possible frequencies and paths in each iteration, which is constructed such that it results in a lower bound to the original problem.
The obtained solution then guides the algorithm in refining the network representation, corresponding to adding additional frequencies and paths, so that the lower bound is improved over the iterations.
This process is repeated until the found solution can be proven to be optimal, which normally happens well before all frequency options and paths in the original problem are considered.
Hence, the method is particularly well-suited for considering a large set of frequencies per line, enabling a more detailed representation of both passenger and operator costs. 

We evaluate the model and solution method on instances of varying sizes that are constructed from four \added{commonly used} data sets in the LinTim software package \autocite{SchieweSchoebelJaegeretal.2024}.
The results demonstrate that the proposed algorithm can achieve significant improvements compared to the commercial solver CPLEX in terms of computational efficiency and solution quality.
Here, we find that the algorithm especially outperforms the commercial solver when the operator costs play a significant role in the trade-off in our objective.
Moreover, we show that explicitly considering lost demand leads to a better use of available resources, generally allowing for better service to the passengers who are willing to travel at reduced operating costs compared to a fixed-demand approach.

This paper presents three main contributions. 
First, we present a line planning model that accounts for lost demand, where demand is only captured if the available paths provide an acceptable service. 
In contrast to most existing approaches, this reflects that passenger demand depends on the service and allows to estimate the amount of travel demand captured. 
Second, we introduce the \frameworkNameShortLP for line planning. 
This exact algorithm is tailored to our model's definition of capacity assignments to lines, but can be adapted to different passenger-oriented formulations for the \problem.
In the context of this algorithm, we also propose customized valid inequalities,
accelerating convergence and computation time. 
Third, we conduct extensive experiments to evaluate the effectiveness of our approach.
Our analysis shows that the \frameworkNameShortLP leads to a significant speedup compared to solving the model directly by a commercial solver, especially when the operating costs take an important role in the objective trade-off.
Moreover, our findings show that we can obtain a more efficient use of resources when recognizing that demand depends on the service in optimization. 

The remainder of this paper is structured as follows. \autoref{sec:Lit} reviews related work on line planning. 
A formal problem definition and model description are given in \autoref{sec:pasgAndDemElast}. 
\autoref{sec:Methods} introduces the \frameworkNameShortLP. 
\autoref{sec:ExpDesign} describes the setup of our numerical experiments, and computational results are presented in \autoref{sec:Results}. Finally,  \autoref{sec:Discussion} summarizes our findings \added{and suggests different directions for future research}.

%% file: 02LitEJOR.tex
\section{Related Literature} \label{sec:Lit}
Line planning has been widely studied, with various models and solution methods proposed. 
For an overview of modeling approaches, we refer the reader to the surveys by \textcite{schobel_line_2012} and \textcite{schmidt_planning_2024}. This review focuses on how studies model passenger service, the scalability of exact methods for line planning, and their approach to demand (fixed or flexible).

The \problem is part of transit network design, which includes the sub-problems of stop and edge selection, line generation, line selection, and frequency setting. The \problem addresses line and frequency selection given a line pool. It is closely related to the Transit Route Network Design Problem (TRNDP), which involves line generation and selection, and the Transit Route Network Design and Frequency Setting Problem (TRNDFS), which additionally incorporates frequency setting.  The latter is also referred to as line planning with Full Pool, as the solution space is not limited by a predefined set of candidate lines \autocite[see][]{heinrich_algorithms_2022}. For an overview of studies on the TRNDP and TRNDFS, we refer to the survey by \textcite{duran-micco_survey_2022}. 

\subsection{Modeling of Passenger Service}
The \problem objective is either related to the passenger service, operator costs, or a mixture of these \autocite{schobel_line_2012}. Cost-oriented models focus on the operator costs subject to constraints on the service. 
\textcite{claessens_cost_1998} minimize total cost for the Dutch railway system while passengers are distributed to specific connections based on shortest paths before line and frequency selection. This {\em system-split approach} is also used by \textcite{bussieck_optimal_1997}, and~\textcite{goossens_branch-and-cut_2004,goossens_solving_2006}. 

Passenger-oriented models focus on metrics related to passenger service, such as travel time or transfers. 
\textcite{bussieck_optimal_1997} prioritize reducing transfers by maximizing direct connections, assuming fixed traffic loads with passengers traveling on the shortest path. 
\textcite{goerigk_line_2017} argue that fixed paths are not suitable when minimizing travel time as the time depends not only on distance but also on transfer time. 
\added{In this regard, \textcite{borndorfer_column-generation_2007} present a model that minimizes in-vehicle time, where passenger paths are iteratively generated through column generation.
Here, passengers may be assigned to arbitrarily long paths.
\textcite{borndorfer_direct_2012} extend this model by distinguishing between direct passengers and those requiring at least one transfer. 
More recently, \textcite{gattSolvingLinePlanning2025} further extended this work by incorporating service-level considerations through bounds on the passengers' travel time and requirements on the share of passengers with direct connections.}
\textcite{guan_simultaneous_2006} \added{present a model that minimizes} operator costs and travel time, including a fixed transfer penalty, which captures the general inconvenience of transferring. 
\textcite{schobel_line_2006} introduce the change \& go network (CGN), a graph representation that models travel time and transfer penalties to minimize passenger travel time under budget constraints. They suggest that frequency-dependent passenger path costs can be modeled with a CGN, an approach adopted in this work. 
\textcite{goerigk_line_2017} propose a bi-level model with operator objectives at the upper level and passenger path choice at the lower, capturing free path choice rather than assigning passengers to paths according to a system optimum. They note that adding capacity constraints significantly increases complexity, making the problem difficult to solve even when linearized.  
\textcite{Gatt2022} also apply a bi-level model for user-optimal routing for the frequency setting problem.  
The interaction between passengers and line concept is further examined by \textcite{schiewe_line_2019}, who propose a game-theoretical approach. 

In this work, we consider a model that uses a path-based passenger assignment with frequency-dependent operator and passenger path costs. 
\added{The candidate path set is pre-generated per OD pair, i.e., group of passengers with similar origin and destination, and limited to paths within a maximum weighted travel time. While we acknowledge that traditionally demand elasticity towards travel time is assumed to be more gradual, and our model results in a hard upper bound in maximum weighted travel time per OD pair, we still believe this is an interesting first step towards including interaction between demand volume and supply in the context of an exact solution method for line planning.}

\subsection{Demand Considerations}
Most studies assume fixed demand and that all passengers must be served \autocite{schmidt_planning_2024}.
However, since alternative modes usually exist, passenger demand depends on passengers' willingness to use a service. Only a few studies consider how much demand is attracted in the line planning stage. \textcite{klier_urban_2015} maximize ridership while considering that the probability of a passenger traveling in the public transport network depends on the frequency. They use a CGN to generate all frequency combinations of the lines, which allows them to precompute frequency-dependent transfer and travel time, similar to our approach. The passenger choice is determined using a discrete choice model. Passengers choose a path freely, regardless of the capacity of vehicles. 
\textcite{hartleb_modeling_2023} maximize revenue and integrate path choice while at the same time ensuring that sufficient capacity is available. 
 For each OD pair, a set of reasonable paths is precomputed, considering driving time and a fixed transfer penalty. The model assumes that passengers distribute over the paths in the choice set made available by the line plan without accounting for transfer and waiting time. 
The same objective is used by \textcite{bertsimas_data-driven_2021}, who assume that the relative number of passengers traveling depends on the frequency and whether a transfer is necessary. They do not assign passengers to explicit paths and, therefore, do not measure frequency-dependent travel time but assume that operating a higher frequency will increase modal share by a predefined percentage. 
In this work, not all stops and demand need to be served. Unlike other studies, our model minimizes weighted passenger and operating costs, including fare revenue,
and accounts for demand lost when weighted travel time exceeds OD-specific thresholds. 

\subsection{Scalability of Solution Methods}
Integrating passenger behavior with network design results in challenging mathematical models.
\textcite{cancela_mathematical_2015} note that, at the time, heuristics and metaheuristic approaches were the only practical methods to solve large-scale line planning problems. 
\textcite{bertsimas_data-driven_2021} also recognize that advanced techniques like mathematical optimization models are not widely used in transit network design due to scalability issues. 
\textcite{guan_simultaneous_2006} consider line selection with integrated routing through a candidate path set for problems of up to 49 stops (preprocessed to 9). \textcite{vermeir_exact_2021} consider line selection and line generation and solve Mandl's Swiss Network with 15 stations using a branch-and-bound algorithm. 
\textcite{cancela_mathematical_2015} incorporate frequencies, solving instances of up to 84 nodes, 143 edges, 378 ODs and 48 candidate lines within 12$\%$ of optimality. 
\textcite{goerigk_line_2017} integrate passenger choice in the line planning procedure, but the formulation scales to only ten stops, relying on a genetic algorithm for larger instances.  
Column generation has been effectively applied by \textcite{borndorfer_column-generation_2007} to a large-scale TRNDFS problem considering 410 stops and 891 edges.  
Their method finds LP-optimal solutions within 
minutes, while integer solutions are obtained with a greedy algorithm. 
\textcite{borndorfer_direct_2012} extend this framework and solve cases with 23 stops, 106 edges, 420 ODs, and up to 2,679 lines to optimality; larger cases with over 7,000 ODs are solved to near-optimality. 
\textcite{bertsimas_data-driven_2021} use a similar approach to the passenger modeling with a one-transfer limit but include line generation. A preprocessing step accelerates their subproblem by a factor ten with minimal impact on the solution, allowing them to solve problems of 233 stations, 2042 edges, and 597 ODs. \added{\textcite{gattSolvingLinePlanning2025} propose a tailored column generation-based heuristic and apply it to a medium-sized French case study with 80 stations and 2,797 OD pairs.}

Studies considering discrete frequencies often limit the analysis to only a few options, e.g., \textcite{klier_urban_2015}, \textcite{bertsimas_data-driven_2021}, \textcite{cancela_mathematical_2015}, and \textcite{borndorfer_direct_2012} consider one to four options. 
In contrast, the proposed method handles 19 options for networks of up to 51 stops and 52 edges, providing a more precise estimate of costs and service. 

Our approach shares structural similarities with the Dynamic Discretization Discovery (DDD) algorithm of~\textcite{boland_continuous-time_2017}, which addresses the Continuous Time Service Network Design Problem. The DDD algorithm iteratively solves a reduced integer program over a partially time-expanded network, dynamically adding time points so that the integer program ultimately yields an optimal solution to the continuous-time problem. The integer program is constructed such that each iteration provides a lower bound on the continuous-time objective, under the assumption that the reduced time dimension, e.g., the holding cost, does not impact the optimal solution. 
The \frameworkNameShortLP similarly solves a sequence of smaller integer programs but differs in (i) how it compresses them, (ii) how it extends from one iteration to the next while guaranteeing a non-decreasing lower bound, and (iii) not assuming that (part of) the reduced dimension does not impact the objective. 

\added{Integrating more planning stages, such as line planning and timetabling, and ideally vehicle and crew scheduling, poses a significant challenge. Each of these problems is computationally difficult on its own, and this is the main reason they are traditionally solved in sequence~\autocite{CATS2025100160}. While integrated approaches can enable better solutions, the mathematical models involved allow only very small instances to be solved to optimality and are often intractable, see e.g.,~\textcite{schiewe2020a, SCHIEWE2022100073}, who consider the integration of line planning, timetabling, and vehicle scheduling. 
Practical methods rely on iterative coordination of the respective problems, see e.g.,~\textcite{SCHOBEL2017348},~\textcite{fuchsEnhancingInteractionRailway2022}, and~\textcite{ZHANG2022240}, and heuristics, see e.g.,~\textcite{kaspi2013a}. The latter two works focus on the integration of line planning and timetabling in a rail setting, where sequentially solving the two problems may even lead to infeasibility. Therefore, a more scalable exact approach could be an interesting building block in these iterative methods.}

%% file: 03Problem.tex
\section{The Line Planning Problem with Service-Dependent Demand}\label{sec:pasgAndDemElast}
Let the travel demand \added{for a given time period} be described by an OD matrix \ODset, where every OD pair $\od \in \ODset$ is described by an origin stop \ODor, a destination stop \ODde, and a number of travelers \ODw. 
\added{As is common in the line planning literature, we assume that this OD matrix represents the demand during a representative time period, typically the morning peak hour.}
\added{The demand given by \ODset{}} will only be captured by the public transport system if the provided service in terms of weighted travel time offered by the selection of lines, frequency, and vehicles provides a sufficiently good service as well as sufficient capacity for the OD pair. Otherwise, demand is lost, and a penalty for this is incurred. We term this interaction between line concept and demand as \emph{service-dependent demand}, reflecting that passenger demand depends on the weighted travel time between any two stops.

\subsection{Passenger flows and the change \& go network}
Passenger flows are defined as paths in a CGN, which is derived from a public transport network (PTN), a candidate line pool, and a set of potential frequencies per line\added{ that are defined relative to the time period represented by the OD matrix}. 
The PTN $G = (\stopSet,E)$ is defined by a set of nodes \stopSet representing the potential stops and a set of edges $E$ indicating which pairs of stops can be connected directly, as well as the travel time required for this journey. For ease of notation, we assume that all edges are bi-directional and travel time is symmetric. 
The line pool $\linePool$ contains a set of lines, where each line $\linePT \in \linePool$ is a simple path $(e_{l 1}, \dots, e_{l m})$ in $G$, i.e., we assume that each line represents a bi-directional service between two terminal stops.
We denote by $E_l = \{e_{l 1}, \dots, e_{l m}\} \subseteq E$ the edges passed by line $\linePT \in \linePool$.
Each line $\linePT \in \linePool$ is associated with a set of frequencies $\frequencySet_l$. The frequency defines the minimum number of vehicles required to operate the line. 
We denote by $\frequencySet{} = \bigcup_{\linePT \in L} \frequencySet_l$ the set of all potential frequencies in the public transport system.

\begin{figure}[htbp]
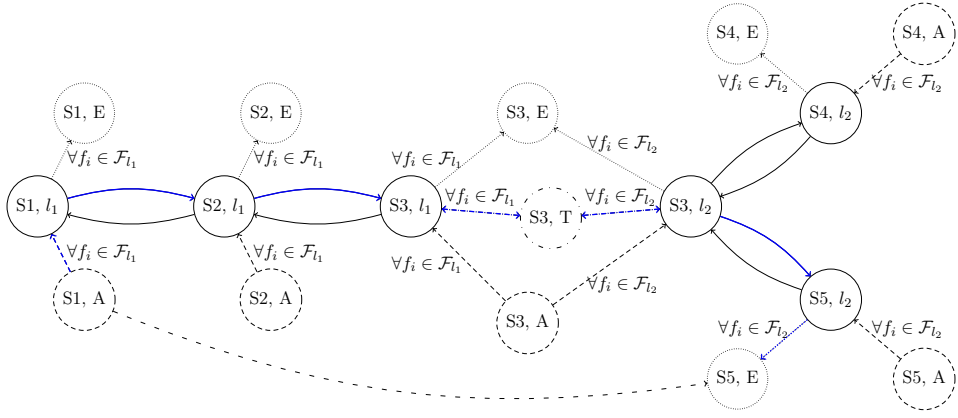

    \centering

    \includestandalone[width=0.75\textwidth]{figures/CNG/CorrectChangeAndGoGraphDense}    
    
    \caption{An example CGN with two lines ($l_1,l_2$), each with three stops. A path for OD pair (1,5) is shown in blue, as well as a direct path from the access node of stop 1 to the egress node of stop 5, representing an alternative mode.  }
    \label{fig:changeAndGo}
\end{figure}

The CGN = (\stopSetCG, \arcSet), as first introduced by \textcite{schobel_line_2006} and illustrated in \autoref{fig:changeAndGo}, is a directed multigraph that allows modeling of both travel and transfer decisions.
The node set is defined as $\stopSetCG:= \stopSetCGA  \cup \stopSetCGE  \cup  \stopSetCGT \cup \stopSetCGL $, where:
\begin{itemize}
    \item \stopSetCGA and \stopSetCGE contain a node for every stop $\stopn\in \stopSet$, that is an origin or destination stop for any OD pair $\od\in \ODset$, thus $\ODor \in \stopSetCGA$ and $\ODde \in \stopSetCGE$. These nodes represent the start and end of the journey of every OD pair, respectively. 
    \item \stopSetCGT contains a node for every stop $\stopn\in \stopSet$ in the PTN $G$ when at least two different lines in the line pool contain this stop. These nodes represent transferring passengers and can only be accessed by passengers already in the system. 
    \item \stopSetCGL contains a node for every line and every stop it contains in its path. 
\end{itemize} 

The arcset contains five subsets, i.e., $\arcSet: = \arcSetDrive \cup \arcSetAcces \cup \arcSetEgress \cup \arcSetTransfer\cup \arcSetOptOut$. 
\begin{itemize}[]
    \item \arcSetDrive contains arcs representing in-vehicle time. Two such arcs are defined for each edge $e \in E_l$ of each line $\linePT \in \linePool$, one for each direction. Arcs specific to line $l \in \linePool$ are denoted by  $\arcSetDrive_l \subseteq \arcSetDrive$.   
    \item \arcSetAcces contains arcs that represent passengers accessing a specific line at a specific frequency at a stop, thus connecting nodes in \stopSetCGA with nodes in \stopSetCGL when they are associated with the same stop $\stopn\in \stopSet$, and having one arc for each potential frequency.
    \item \arcSetEgress contains arcs that represent passengers arriving at a stop after at least one trip on a line. It contains an arc from a line-and-stop node in \stopSetCGL to an egress stop node in \stopSetCGE for every stop $\stopn\in \stopSet$, and every line $\linePT \in \linePool$ that stops at this stop, and every possible frequency $f \in \frequencySet_{l}$.
    \item \arcSetTransfer contains arcs that represent passengers transferring at a stop. Therefore, we define arcs from every node in \stopSetCGL to the transfer node of that stop in \stopSetCGT, a separate arc for each frequency. 
    \item \arcSetOptOut contains arcs that represent travelers choosing an alternative mode of transport. These arcs connect nodes in \stopSetCGA to nodes in \stopSetCGE with one direct arc for each unique OD pair $(\ODor, \ODde)\in \ODset$. 
\end{itemize}  

\autoref{fig:changeAndGo} illustrates the CGN = (\stopSetCG, \arcSet) for two lines: Line $l_1$ serves stops 1 to 3, and line $l_2$ serves stops 4, 3 and 5. This makes stop 3 the only transfer stop. 
Different node and arc types are distinguished by style: access nodes \stopSetCGA are dashed circles; egress nodes \stopSetCGE are dotted; line nodes \stopSetCGL are solid; and transfer nodes \stopSetCGT are dash-dotted. Arcs follow a similar visual scheme.
To increase readability, all frequency-dependent access, egress, and transfer arcs have been displayed once with a label ``$\forall f_i \in \mathcal{F}_{lx}$''. Note that feasible paths should have matching frequency arcs for access, egress, or transfer from a single line. This can either be enforced in a postprocessing step, as a constraint, or alternatively enforced in the graph by redefining the line nodes to be frequency-specific. These are implementation decisions that do not affect the principles of the model or the proposed algorithm. Furthermore, only a single arc in \arcSetOptOut has been included in \autoref{fig:changeAndGo}, namely representing traveling by an alternative mode for OD pair (1,5). 

\subsection{Passenger Demand}
Passengers travel along a simple path $p$ in the CGN, which is defined by a sequence of arcs $\arc \in \arcSet$. Passenger service is expressed by a function \added{$g$} dependent on the arcs in the path: 
\vspace{-0.5\baselineskip}
\begin{align}
    \added{g}(p):= \sum_{\arc \in p \cap \arcSetAcces} c^{A}_{\arc} + \sum_{\arc \in p \cap \arcSetDrive} c^{IVT}_{\arc} +  \sum_{\arc \in  p \cap \arcSetTransfer} c^{T}_{\arc} + \sum_{\arc \in p \cap \arcSetEgress} c^{E}_{\arc} 
    \label{eq: cost allocation}
\end{align}
Here, $c^{A}_{\arc}, c^{IVT}_{\arc}, c^{T}_{\arc}$, and $c^{E}_{\arc}$ represent the cost of an arc $\arc \in \arcSet$ that may depend on the type of arc, a possible distance or travel time associated to the arc, and the frequency associated to the arc. This function must be non-increasing in the frequency, that is, a higher frequency cannot incur a higher cost in the path than a lower frequency. Such cost functions allow for assigning waiting time a higher cost than in-vehicle time, making the in-vehicle time cost dependent on the line (mode), and assigning both a fixed penalty as well as a frequency-dependent penalty on transfers. 
Thus, the cost of a path depends on the frequency of the line through access arcs \arcSetAcces and egress arcs \arcSetEgress, on the length of the path through in-vehicle arcs \arcSetDrive, and both on the number of transfers, as well as some function of the frequency of the feeder and receiver line, through \arcSetTransfer.

We model the interaction between line concept and demand by setting an OD pair-dependent maximum on the cost of the path expressed by $\added{g}(p)$. All paths for each OD pair $\od\in\ODset$ must have a cost below a threshold \wttThresOD. When passengers have access to a path in the resulting public transport system with a cost lower than \wttThresOD, then we assume all demand \ODw will be captured by the network if there is sufficient capacity. Otherwise, the demand is lost, meaning the operator loses the fare revenue, and passengers must take an alternative, which may reduce social welfare. 
Lost demand will be assigned to an arc in $\arcSetOptOut$ representing an alternative-mode path $\altModePath_\od$ with a cost at least as high as the cost of the maximum acceptable weighted travel time path for OD pair $\od$.

Unique to our approach is that the threshold for demand captured by the public transport system depends on a general function $\added{g}$ that explicitly depends on the frequency of the selected lines. 
This relies on a detailed choice of frequencies, one of our model's distinguishing features. 
If a more gradual relation between $\added{g}(p)$ and captured demand is desired, one could split up the OD pair \od into multiple OD pairs with individual thresholds. This introduces a smoother relationship between service quality (as measured by $\added{g}(p)$) and the proportion of demand served, which is a step towards approximating a demand curve and integrating how demand might vary with service levels.
In the following, we will assume that this path set per OD pair \od can be pre-generated, which has notational benefits for the presentation of the model and the solution method. 

\input{Tables/Notation/Notation}

\subsection{Model} \label{sec:model}
The \problem selects a subset of lines $\linePoolSol \subseteq \linePool$ together with a frequency $\added{f^*}\in \frequencySet_l$ for each selected line $l \in \linePoolSol$. Its objective is to minimize the weighted sum of passenger cost in terms of weighted travel time and operating cost dependent on the number of lines, the \added{number of vehicles assigned to each line}, and the ticket fare revenue. 
Input consists of a line pool $\linePool$, a set of candidate frequencies for each line $\frequencySet_{l}$, and a travel demand \ODset, where each OD pair is associated with a candidate path set \PsetOD. Each path is associated with a cost $c_p$, as previously defined. Additional input consists of a set of vehicle types $\vehicleSet{}$, with $\vehicleSet{}_l$ the set of vehicle types available to operate line $\linePT \in \linePool$. Different vehicle types could represent longer versus shorter vehicles, and also allow to model different modes. The driving speed per line $\linePT \in \linePool$ is fixed, thus a line would typically be associated to one mode, but the solution could select lines that represent a multi-modal network. Let $\phi_{lf}$ denote the minimum number of vehicles needed to operate line $l \in \linePool$ at frequency $f \in \frequencySet_l$, and let $\delta_{lv}$ denote the maximum number of passengers that can be transported by one vehicle of type $v \in \vehicleSet{}_l$ assigned to line $l$ on any edge of that line. $\delta_{lv}$ can be derived from the size of the vehicle as well as the number of round trips \added{over the considered period} 
per vehicle on line $\linePT \in \linePool$. 
$c_v$ denotes the cost of using a single vehicle of type $v \in \vehicleSet{}$, and $h_{lf}$ denotes the fixed cost of opening line $\linePT \in \linePool$ at frequency $f\in\frequencySet_l$. Finally, $R$ represents the operator’s revenue per passenger, which can be extended to account for distance-based or zonal pricing depending on the OD pair. 

The decision on whether to open a line $l$ at a given frequency $f$ is represented by the binary variable $y_{lf}$. 
The cost, capacity, and frequency of a line depend on the number and type of vehicles assigned to a line. This is represented by decision variables $z_{lv}$, which can be either continuous or integer. 
The passenger flow is represented by continuous variables $x_{pd}$, which denote the fraction of passengers of OD pair $\od \in \ODset$ traveling on path $p$. 
The following MILP model is an adaptation of the shuttle planning model in \textcite{hurk_shuttle_2016}, where we allow for lost demand, model passenger flows as fractions instead of actual demand, and account for fare revenue in the objective. \vspace{-1.2\baselineskip}
\input{Model/MILP} 
\vspace{-1.5\baselineskip}

\added{The objective of the model combines passenger travel disutility and operator costs into a single welfare function. Passenger costs reflect generalized travel costs on the assigned paths, while operator costs represent fleet usage, line operation, and revenues generated from served demand, which capture the operator’s incentive to provide service. The resulting objective represents a system-optimal trade-off between service quality and operational efficiency and is not intended to model individual passenger choice behavior.} 
Constraints \eqref{con1} ensure that all passengers are assigned a path\added{, which could either be a path in the public transport network or a path representing the use of an alternative mode. 
Moreover}, constraints \eqref{con2} ensure that the vehicle capacity is not exceeded 
for any line\added{ by requiring that, for each segment of the line, the passenger demand 
does not exceed the available seating capacity within the given time period}. 
Constraints \eqref{con4} ensure that paths can only be used if the associated frequency of the line is operated. 
\added{Passengers are assigned to specific lines, which means that if several lines in a solution \linePoolSol cover the same edges, passengers would realistically take the earliest service. As a result, the model may overestimate the actual waiting times. }
Constraints \eqref{con5} enforce that each line operates at most one frequency. 
Not assigning any frequency means that the line is not operated. 
Constraints \eqref{con6} ensure that if a line $l$ is opened at frequency $f$, then at least $\phi_{lf}$ vehicles must be assigned. 
Constraint \eqref{conBudget} limits the operational resources to be within a given budget $B$. 
Finally, constraints \eqref{domain_x}, \eqref{domain_z}, and \eqref{domain_y} specify the domains of the variables. 

\added{As can be seen from the formulation, our model uses a system-optimum assignment of passengers to paths, where passengers belonging to the same OD can be split over different paths. 
This could mean that some passengers travel in the public transport network while others use an alternate mode, or that passengers are split over paths covering different lines or different segments of the same line.
While an all-or-nothing assignment is possible, splitting the passengers can be advantageous when the increased passenger cost is smaller than the additional cost that is needed to increase the capacity on some lines accordingly.
We note, though, that the number of vehicles for a line can be chosen beyond the lower bound given by the number of vehicles needed to achieve the frequency selected for this line, such that the above does generally not happen often.}

\added{In terms of capacity, it can be seen that constraints \eqref{con2} consider the time-aggregated seating capacity provided on each segment of a line over the full planning period, which matches the assumption that \ODset represents the demand for a representative period, usually the peak hour. 
As commonly assumed in line planning literature, the found line plan can then afterwards be adjusted to better match demand in other demand periods, e.g., by reducing the frequency of some lines during off-peak hours.
Moreover, uneven passenger demand during the peak hour could be further taken into account in the timetabling step.
Note that the found line plan admits a timetable requiring the same number of vehicles, but that further adjustments of the timetable to time-varying demand or other timetable constraints, e.g., headway constraints in a rail setting, might increase the required number of vehicles.}

\added{Finally, passenger comfort is modeled through hard capacity constraints rather than explicit crowding effects. This modeling choice aligns with the literature on line planning where capacity constraints are used to guarantee feasibility at an aggregate level, while a more detailed representation of passenger comfort is typically addressed in later stages when a line concept and timetable are fixed. }

%% file: figures/CNG/CorrectChangeAndGoGraphDense.tex
\begin{tikzpicture}
	\tikzstyle{cgNode} = [circle, draw, minimum size=0.75cm]
	\tikzstyle{accessNode} = [cgNode, densely dashed]
	\tikzstyle{egressNode} = [cgNode, densely dotted]
	\tikzstyle{transferNode} = [cgNode, loosely dashdotted]

    \node[cgNode] (S1L1) at (0,2) {S1, $l_1$};
    \node[accessNode] (S1A) at (1,0) {S1, A};
    \node[egressNode] (S1E) at (1,4) {S1, E};

    \node[cgNode] (S2L1) at (4,2) {S2, $l_1$};
    \node[accessNode] (S2A) at (5,0) {S2, A};
    \node[egressNode] (S2E) at (5,4) {S2, E};
    
    \node[cgNode] (S3L1) at (8,2) {S3, $l_1$};
    \node[transferNode] (S3T) at (11,1.75) {S3, T};
    \node[accessNode] (S3A) at (10.5,-0.5) {S3, A};
    \node[egressNode] (S3E) at (10.5,4) {S3, E};

    \node[cgNode] (S3L3) at (14,2) {S3, $l_2$};

    \node[cgNode] (S4L1) at (17,4) {S4, $l_2$};
    \node[accessNode] (S4A) at (19,5.7) {S4, A};
    \node[egressNode] (S4E) at (15,5.7) {S4, E};

    \node[cgNode] (S5L1) at (17,0) {S5, $l_2$};
    \node[accessNode] (S5A) at (19,-1.7) {S5, A};
    \node[egressNode] (S5E) at (15,-1.7) {S5, E};

	\tikzstyle{accessArc} = [->, densely dashed]
	\tikzstyle{egressArc} = [->, densely dotted]
	\tikzstyle{transferArc} = [<->, densely dashdotted]

    \path[->] (S1L1) edge [bend left=15] (S2L1);
    \path[->] (S2L1) edge [bend left=15] (S1L1);    
    \draw[->] (S2L1) edge [bend left=15] (S3L1);
    \draw[->] (S3L1) edge [bend left=15] (S2L1);

    \path[->] (S4L1) edge [bend left=15] (S3L3);
    \path[->] (S3L3) edge [bend left=15] (S4L1);    
    \draw[->] (S3L3) edge [bend left=15] (S5L1);
    \draw[->] (S5L1) edge [bend left=15] (S3L3);

     \foreach \station in {S1, S2} {
        \draw[accessArc] (\station A) -- (\station L1) node[midway, right]{$\forall f_i \in \mathcal{F}_{l_1}$};

        \draw[egressArc] (\station L1) -- (\station E) node[midway, right]{$\forall f_i \in \mathcal{F}_{l_1}$};
    }
    \foreach \station in {S4, S5} {
        \draw[accessArc] (\station A) -- (\station L1) node[near end, right]{$\forall f_i \in \mathcal{F}_{l_2}$};

        \draw[egressArc] (\station L1) -- (\station E) node[near start, left]{$\forall f_i \in \mathcal{F}_{l_2}$};
    }
    \foreach \line in {L1} {
        \draw[accessArc] (S3A) -- (S3\line) node[midway, left]{$\forall f_i \in \mathcal{F}_{l_1}$};

        \draw[egressArc] (S3\line) -- (S3E) node[midway, left]{$\forall f_i \in \mathcal{F}_{l_1}$};
    }
    \foreach \line in {L3} {
        \draw[accessArc] (S3A) -- (S3\line) node[near start, right]{$\forall f_i \in \mathcal{F}_{l_2}$};

        \draw[egressArc] (S3\line) -- (S3E) node[near end, right]{$\forall f_i \in \mathcal{F}_{l_2}$};
    }

    \draw[transferArc] (S3L1) -- (S3T) node[midway, above]{$\forall f_i \in \mathcal{F}_{l_1}$};
    \draw[transferArc] (S3L3) -- (S3T) node[midway, above]{$\forall f_i \in \mathcal{F}_{l_2}$};

    path = {S1, S2, S3, S5}
    \tikzstyle{pathAccessArc} = [accessArc, blue, thick] 
    \tikzstyle{pathEgressArc} = [egressArc, blue, thick] 
    \tikzstyle{pathTransferArc} = [transferArc, blue, thick]
    \draw[pathAccessArc] (S1A) -- (S1L1) node[midway, right]{};
    \path[->, blue, thick] (S1L1) edge [bend left=15] (S2L1){};
    \draw[->, blue, thick] (S2L1) edge [bend left=15] (S3L1){};
    \draw[pathTransferArc] (S3L1) -- (S3T) node[midway, above]{};
    \draw[pathTransferArc] (S3L3) -- (S3T) node[midway, above]{};
    \draw[->, blue, thick] (S3L3) edge [bend left=15] (S5L1){};
    \draw[pathEgressArc] (S5L1) -- (S5E) node[midway, right]{};
    \draw[->,loosely dashed] (S1A) edge [bend right=15] (S5E);

\end{tikzpicture}

%% file: Tables/Notation/Notation.tex
\begin{table}[htbp!] 
		\centering
        \caption{Notation and terminology}
            \renewcommand{\arraystretch}{0.97}
		\begin{tabular}{ll}
			\toprule
         	\textbf{Sets}& Description\\
            \midrule
           $\stopSet$ & set of candidate stops, with stop $\stopn\in \stopSet$\\
            $E$ & set of bidirectional edges, edge $e \in E$\\
            $\linePool$ & line pool of potential lines     \\ 
            $E_l \subseteq E$ & set of edges passed by line $l \in \linePool$\\
			$\frequencySet$ & set of frequencies   \\ 
			$\frequencySet_l \subseteq \frequencySet$ & set of frequencies that can be operated for line $l \in \linePool$ \\  
            $\arcSet$ &  set of directed arcs in the CGN \\ 
            $\arcSetDrive \subseteq \arcSet$ & set of arcs representing in-vehicle time\\
            $\arcSetDrive_l\subseteq \arcSetDrive$ & set of arcs associated with line $l \in \linePool$. Every $e \in E_l$ is associated with two directed arcs\\
            $\vehicleSet{}$ & set of vehicle types \\
            $\vehicleSet{}_l \subseteq \vehicleSet$ & set of vehicle types accepted for line $l \in \linePool$\\
			\ODset & set of OD pairs    \\ 
            \Pset & set of paths in the CGN-graph \\
            $\PsetOD \subseteq \Pset$ & set of paths for OD pair $\od \in \ODset$\\
            $\altModePath_\od \in \PsetOD$ &  alternative mode path for OD pair $\od \in \ODset$ \\ 
            $\PsetODedge \subseteq \PsetOD $ & paths traversing arc $a \in \arcSetDrive$ for OD pair \od \\
            $\PSetLineFreq \subseteq \PsetOD$& set of paths for OD pair $\od \in \ODset$ that use line $l \in \linePool$ at frequency $f \in \frequencySet_l$\\
           \midrule
           \textbf{Parameters} & Description\\
           \midrule
            $c_v$ & flexible cost per vehicle of type $v \in \vehicleSet$ \\
            $\phi_{lf}$ & minimum number of vehicles required to operate line $\linePT \in \linePool$ at frequency $f \in \frequencySet_l$ \\
             $\delta_{lv}$ & 
             maximum number of passengers per vehicle $v \in \vehicleSet$ on line $l\in \linePool$\\
            \ODw &  number of passengers in OD pair $\od  \in \ODset$\\
            $\ODor \in \stopSet$ & origin stop of OD pair \od \\
            $\ODde \in \stopSet$ & destination stop of OD pair \od \\
			$c_p$ & cost of path $p \in \mathcal{P}$\\
            $\wttThresOD$ & maximum acceptable path cost for passengers in OD pair $\od$ \\
            $h_{lf}$ & fixed cost of opening line $l \in \linePool$ at frequency $f\in  \frequencySet_l$ \\
            $B$ & operational budget \\
            $R$ & the fare revenue from a single journey \\
           $\lambda$  & a weighing parameter on the passenger costs  \\
           \midrule
           \textbf{Variables}& Description\\
           \midrule
            $y_{lf}$ & decision to open line $l$ at frequency $f$\\
            $x_{p \od}$ & the fraction of passengers of OD pair \od on path $p$\\
            $z_{lv}$ & number of vehicles of type $v \in \vehicleSet$ assigned to line $l \in \linePool$ \\
            \bottomrule
	\end{tabular}
    \label{tab:notationAndTerm}
\end{table}

%% file: Model/MILP.tex
\begin{align}
    \min \;
    \lambda \left(
    \sum_{\od\in \ODset} 
    \sum_{p \in \PsetOD}
    c_p \ODw x_{p\od}\right)  +
    \left( 
    \sum_{l \in \linePool} 
    \sum_{v \in \vehicleSet_l}
    c_v z_{lv} + 
    \sum_{l \in \linePool} 
    \sum_{f \in \frequencySet_l} 
    h_{lf} y_{lf} - 
     \sum_{\od\in \ODset} 
    \sum_{p \in \PsetOD \setminus \altModePath_\od}
    R \ODw x_{p\od} 
    \right) \label{obj} 
     \end{align}
     \vspace{-0.5cm}
subject to \vspace{-0.5cm}
     \begin{align}
     \sum_{p \in \PsetOD} x_{p\od} &= 1 && \forall \od \in \ODset \label{con1}  \\
        \sum_{\od\in \ODset} \sum_{p \in \PsetODedge} \ODw x_{p\od} &\leq \sum_{v \in \mathcal{V}_l} \delta_{lv} z_{lv}   && \forall l \in \linePool, \forall a \in \arcSetDrive_l  \label{con2}  \\
         \sum_{p \in \PSetLineFreq} x_{p\od} &\leq y_{lf}  && \forall \od \in \ODset, \forall l \in \linePool, \forall f \in \mathcal{F}_l \label{con4} \\
       \sum_{f \in \mathcal{F}_l} y_{lf} &\leq 1 && \forall l \in \linePool{} \label{con5}\\
       \sum_{v \in \vehicleSet_l} z_{lv} &\geq y_{lf} \phi_{lf} && \forall l \in \linePool, \forall f \in \mathcal{F}_l \label{con6}\\
      \sum_{l\in \linePool}\left(\sum_{v\in \mathcal{V}_l} c_v z_{lv} + \sum_{f\in\mathcal{F}_l} h_l y_{lf}\right) &\leq B  && \label{conBudget}    \\
       x_{p\od}&\geq 0 && \forall \od \in \ODset, \forall p \in \PsetOD \label{domain_x}\\   
       z_{lv} &\geq 0 && \forall l \in \linePool, \forall v \in \vehicleSet_l \label{domain_z} \\  
       y_{lf} &\in \{0,1\} && \forall l \in \linePool, \forall f \in \mathcal{F}_l \label{domain_y}  
\end{align}

%% file: 04MethodsEJOR.tex
\section{The \frameworkNameLP} \label{sec:Methods}
This section presents the \frameworkNameLP (\frameworkNameShortLP) for solving the \problem. The \frameworkNameShortLP iteratively generates non-decreasing lower bounds for the \problem by dynamically refining the relationship between frequency, operating cost, and travel time for each line. 
We refer to the problem as \problem{$(\frequencySet$)}, since the \frameworkNameShortLP works by refining the set $\frequencySet$, while $\Pset$ is derived from fixed input and the current frequency set $\frequencySet$. The remaining inputs ($\ODset$, $\linePool$, and $\vehicleSet$) remain fixed. For clarity, we describe the algorithm assuming a single vehicle type and denote vehicle assignments by $z_l$ per line $l\in\linePool$. However, the approach extends to multiple types of vehicles.

The \frameworkNameShortLP first generates a compressed frequency set \frequencySetAlg{0}, which contains a {\em single} frequency per line and is defined such that it provides a lower bound on both weighted travel time and operating cost. The initial path set \PsetAlg{0} depends on this adjusted frequency set. The algorithm then solves \problem{$(\frequencySetAlg{0})$} in which the path set is limited to \PsetAlg{0}. The resulting solution provides a lower bound on the solution quality of \problem{$(\frequencySet)$}.  If this solution is feasible for \problem{$(\frequencySet)$}, i.e., considering the full frequency set and the full path set \Pset, the algorithm terminates with an optimal solution. Otherwise, the set of frequencies is refined to a new subset \frequencySetAlg{i} that systematically eliminates infeasibilities, and the process continues. This reduction enables solving initially a very compact model, which gradually grows in complexity as the refinement progresses. 
We define the vector of solutions for $y_{lf}$, $x_{p\od}$ and $z_{l}$ in a given iteration $i$ 
as \vectorSol{i}. The \frameworkNameShortLP in an iteration $i$ first solves \problem{$(\frequencySetAlg{i})$}, next refines that set, and repeats until \vectorSol{i} is feasible for \problem{$(\mathcal{F})$}. \added{In iteration $i$, this procedure provides a selection of lines $\linePoolSol = \{ l \in \linePool \mid \exists f \in \frequencySetAlg{i}_l: y_{lf}^{*} = 1 \}$ and their corresponding frequencies. While the set of candidate lines stays the same, both the line selection and the frequency assignment may change over the course of the \frameworkNameShortLP as the frequency representations of the lines are refined.}

An upper bound can relatively easily be constructed from any solution $\vectorSol{i}$ as we will describe in the following. 
Although these upper bounds are not necessary for the \frameworkNameShortLP, they can be valuable when the
procedure is terminated based on a time limit, as they provide a feasible solution to \problem{$(\mathcal{F})$}. Moreover, these upper bounds can be used to track the progress of the algorithm throughout the solve by comparing against the lower bounds found.
In the following, we discuss the steps of the algorithm in detail. We then present a set of valid inequalities relevant in the context of lost demand. 

\subsection{Initialization}
We note that there is a one-to-one correspondence between the candidate frequency set $\frequencySet_l$ and the candidate headway set $\headwaySet_{l} = \{\frac{1}{f}\mid f \in \frequencySet_l\}$ for each line $l\in \linePool$. For notational convenience, we explore the initialization and refinement of $\frequencySet$ through the equivalent representation $\headwaySet$.  Let function \( f_v \) take a given headway as input and return the required number of vehicles, and let function \( f_h \) do the reverse and return the achievable headway for a given number of vehicles.  We define $\phi_{lh}$ to be the required number of vehicles to operate headway $h$ of line $l$. The set $\headwaySet_{l}$ is defined such that each headway $h\in \headwaySet_{l}$ is associated with a unique number of required vehicles, selecting the minimum headway that can be achieved with any integer number of vehicles.

The headway of a line can be characterized by a tuple $(h,\phi_{lh})$ describing the headway that can be operated and the number of vehicles required to operate the line at headway $h$.  The path set $\Pset$ is a function of the headways represented in $\frequencySet$, as defined by Constraints~\eqref{con4}, while $\phi_{lh}$ defines the minimum number of vehicles that must be operated, as defined by Constraints~\eqref{con6}.  In the full problem, the discrete headway set is defined by the complete set $\headwaySet = \{(h,\phi_{lh})\mid h \in \headwaySet_l, l \in \linePool \}$.

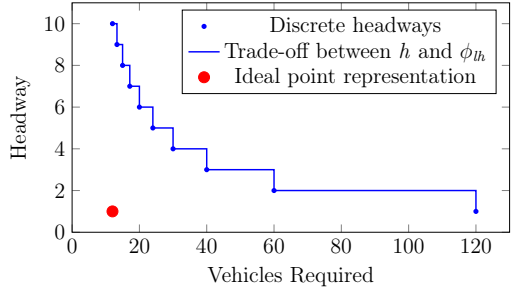
\begin{figure}
\centering
    \resizebox{0.4\textwidth}{!}{%
\begin{tikzpicture}
\begin{axis}[
    width=10cm,
    height=6cm,
    xlabel={Vehicles Required},
    ylabel={Headway},
    legend style={font=\large},
    label style={font=\large},
    tick label style={font=\large},
    xtick distance=20,
    ytick distance=2,
    ymin=0, ymax=11,
    xmin=0, xmax=130, 
    legend pos=north east,
]

\addplot[
    blue,
    thick,
    mark=*,
    mark options={fill=blue},
    only marks, 
    mark size = 1pt,
]
coordinates {
    (120.0, 1)
    (60.0, 2)
    (40.0, 3)
    (30.0, 4)
    (24.0, 5)
    (20.0, 6)
    (17.14285714, 7)
    (15.0, 8)
    (13.33333333, 9)
    (12.0, 10)
};
\addlegendentry{Discrete headways}

\addplot[
    blue,
    thick,
    const plot,
    const plot mark right,
    mark=none,
]
coordinates {
    (120.0, 1)
    (60.0, 2)
    (40.0, 3)
    (30.0, 4)
    (24.0, 5)
    (20.0, 6)
    (17.14285714, 7)
    (15.0, 8)
    (13.33333333, 9)
    (12.0, 10)
};
\addlegendentry{Trade-off between \rev{$\mathit{h}$} and \rev{$\mathit{\phi_{lh}}$}}

\addplot[
    red,
    only marks,
    mark=*,
    mark size=3pt,
]
coordinates {
    (12.0, 1)
};
\addlegendentry{Ideal point representation}

\end{axis}
\end{tikzpicture}}
\caption{The relationship between headway \rev{$h$} and required vehicles \rev{$\phi_{lh}$} for a line $l\in \linePool$ with 10 candidate frequencies.}
\label{fig:TradeOff}
\end{figure}

In the \frameworkNameShortLP, this discretization is replaced by a set of {\em ideal points} per line, which respectively provide a lower bound on the line's passenger and operator costs given the line's possible headways $\headwaySet_l$.  To examine how this lower bound is constructed, we examine the relationship between headways and costs: from a passenger perspective, shorter headways improve service, while from an operator perspective, reducing the number of vehicles used is preferred. 
In our model, we allow the assignment of more vehicles to a line than is strictly required by $\phi_{lh}$. The relationship between service and costs is, therefore, most accurately captured by the step-function in \autoref{fig:TradeOff}.
The \frameworkNameShortLP aims to represent this relationship more efficiently by recognizing that it might not be necessary to explicitly consider all discrete headways in the model. Instead of including the full set of headways, we initially define a single ideal point per line, as illustrated in \autoref{fig:TradeOff}. 
Each line $l$ is modeled under the assumption that the lowest possible headway, $h^l_{\min} = \min(\headwaySet_{l})$, can be achieved using the minimum number of vehicles, $v^l_{\min} = f_v(h^l_{\max})$, where $h^l_{\max}=\max(\headwaySet_{l})$. Consequently, the initial headway set is defined as $\headwaySetAlg{0} = \{(h^l_{\min}, v_{\min}^l)\} \mid l\in \linePool\}$.
The path set $\Pset^0$ can be calculated using the CGN derived from $\linePool, \headwaySetAlg{0}$, and the PTN, or it can be obtained from the full path set $\Pset$ as: $\Pset^0=\{p \in \PsetOD(l,h) \mid \od \in \ODset, l \in \linePool, h\in \headwaySetAlg{0}_l\}$, where $\PsetOD(l,h)$ defines the paths for OD pair $\od\in\ODset$ using line $l\in\linePool$ at headway $h\in\headwaySet_l$. 

\subsection{Solving and Refining the MILP}
Solving the initial ideal point representation of the problem, i.e., \problem{$(\frequencySetAlg{0})$}, yields an optimal solution $\vectorSol{0}$, where $z_l^*\geq v^l_{\min}$ for $l\in\linePool$. When solving \problem{$(\frequencySetAlg{0})$}, $\phi_{lh}$ in
  Constraint~\eqref{con6} is replaced by $v^l_{min}$. Let $h^*_l$ be the headway selected for line $l\in \linePoolSol$.  In reality,
$\phi_{lh_{l}^{*}}$ vehicles must be assigned to operate the promised headway. However, since Constraint~\eqref{con2} ensures that the
number of vehicles must be sufficient to transport the assigned passengers, the model may be required to assign more vehicles than
$v^l_{\min}$. In the case that $z_l^*\geq \phi_{lh_{l}^{*}}$, enough vehicles have been assigned to operate the promised headway, and the solution is therefore feasible to the full problem.

If $z_l^*<\phi_{lh_{l}^{*}}$, then the model is refined to eliminate the current solution and strengthen the lower bound. This is achieved by adjusting the set $\headwaySetAlg{0}$ so that in the next iteration $\vectorSol{0}$ is no longer feasible. To obtain this, we define the set $\headwaySetAlg{1}$, increasing the vehicles required to operate the minimum headway $h^l_{\min}$ and including an additional headway to preserve the lower bound representation on $\headwaySet$. Headway $h^l_{\min}$ is updated to require the number of vehicles necessary to operate a one-step lower headway than the headway that could be operated given $z_l^*$ vehicles, i.e., $\headwaySetAlg{1} \gets \{(h^l_{\min},f_v(f_h(z_l^*)-1))\}$. Here, we assume that headways are discretized with 1-minute intervals, but this can easily be extended to handle any given set of sorted candidate headways. The headway that can actually be operated with the number of vehicles allocated, $f_h(z^*_l)$, is also included in $\headwaySetAlg{1}$ and is set to require the minimum possible vehicles $v^l_{\min}$, i.e. $\headwaySetAlg{1} \gets \headwaySetAlg{1}\cup\{(f_h(z^*_l),v^l_{\min})\}$.

In general, for iteration $i>1$, we continue updating the representation of any line $l\in\linePoolSol$, with selected headway  $h_l^*$, if too few vehicles have been assigned in the optimal solution \vectorSol{i} to actually operate this headway. In particular, we modify the headway to $(h_l^*,f_v(f_h(z^*_l)-1)$, which still provides a lower bound but forces a higher vehicle investment for frequency $h_l^*$ than in the previous iteration. We also introduce the new headway $f_h(z^*_l)$ and set it to require a vehicle investment that is one step higher than the lowest headway in $\headwaySet^{i}$ that operates at a higher headway than $h_l^*$. If no such headway exists in $\headwaySet^{i}$, $v_{\min}^l$ is used. We therefore introduce headway $(f_h(z^*_l), \max\{v_{\min}^l, \min_{h_l \in \headwaySetAlg{i}| h_l> h_l^* } f_v(h_l-1) \})$. Note that in each iteration, new variables in terms of headways and paths are added to the model, as well as the constraints associated with these variables. The algorithmic steps of the \frameworkNameShortLP are formalized in~\autoref{MainPseudocode}, and a small worked example is provided in \autoref{app_Example}. For a given line $l$ and headway $h$, we denote the lower bound on the number of vehicles needed in iteration $i$ as $\underline{\phi}_{lh}^i$. Each iteration of the algorithm provides a lower bound for \problem{$(\frequencySet)$}, as formally stated in Theorem~\ref{thm:lb} and proved in \ref{sec:proof1}.
\begin{algorithm}[htbp]
\caption{\frameworkNameLP for the \problem}
\label{MainPseudocode}
\SetAlgoLined
\SetAlgoNoEnd
\KwIn{PTN $G$, OD matrix $\ODset$, path set $\Pset$, line pool $\linePool$, frequencies $\frequencySet$, vehicle type $v$}
\KwResult{Line concept, vehicle assignment, and passenger flow through \vectorSol{i}} 
\tcc{initialize} 
$\headwaySetAlg{0} \gets \emptyset$\;

\ForEach{line $l\in \linePool$}
{
    $\headwaySetAlg{0} \gets \headwaySetAlg{0} \cup \{(h_{\min}^l, v_{\min}^l)\}$\;
}
$\PsetAlg{0} \gets \text{calcPaths}(\headwaySetAlg{0}$)\;
$i \gets 0$\;
\tcc{build MILP and solve}
\vectorSol{i} $\gets$ solve \problem{(\headwaySetAlg{0})}\; 
\tcc{update line representation and resolve until feasible solution found}
\While{\vectorSol{i} is not feasible}
{
$i \gets i + 1$\;
$\headwaySet^{i}\gets\headwaySet^{i-1}$\;
  \ForEach{line $l\in \linePoolSol$}
    {
      \uIf{$z_{l}^*<f_v(h_{l}^{*})$}{
        \tcc{update vehicle requirements for selected headway}
        $\headwaySet^{i} \gets \headwaySet^{i}\setminus \{(h^*_l,\underline{\phi}_{lh_l^*}^{i-1})\}$\;
        $\headwaySet^{i} \gets \headwaySet^{i}\cup\{(h_l^*, f_v(f_h(z^*_l)-1)\}$\;
        \tcc{include new headway} 
        $\headwaySet^{i} \gets \headwaySet^{i}\cup \left(f_h(z_l^*), \max\left\{v_{\min}^l, \min_{h_l \in \headwaySetAlg{i-1}\mid h_l> h_l^* } f_v(h_l-1)\right\}\right)$\; 
       }
    }
    $\PsetAlg{i} \gets \text{calcPaths}(\headwaySetAlg{i}$)\;
    $\vectorSol{i} \gets $ solve \problem{(\headwaySetAlg{i})}\;    
}
\end{algorithm}

\begin{theorem}
\label{thm:lb}
The optimal solution at any iteration of the algorithm provides a lower bound on the optimal objective value of \problem{$(\mathcal{F})$}.
\end{theorem}
\begin{proof}
  See~\ref{sec:proof1}
\end{proof}

\subsection{Repairing the Solution}
A simple heuristic can be used to repair the solution \vectorSol{i} at iteration $i$ to obtain a valid upper bound if the solution is infeasible for \problem{$(\mathcal{F})$}.
If the solution is not feasible, then it implies that we have assigned fewer vehicles than are required to operate the promised headway. To correct this, we recalculate the true headway $f_h(z_l^*)$ for every $l \in \linePoolSol$ and update the passenger costs based on the true headway. This accounts for the possibility that some OD pairs $d\in\ODset$ might not have any acceptable paths available and must therefore be routed on $\altModePath_\od$. While this repair step is not essential for the
algorithm’s convergence, it bounds the solution quality and
ensures that a feasible solution is available if the procedure is
terminated before convergence.

\subsection{Termination}
If an optimal solution \vectorSol{i} at iteration $i$ is feasible for
\problem{$(\mathcal{F})$}, then it is also optimal for
\problem{$(\mathcal{F})$}. Theorem~\ref{thm:lb} guarantees that
\vectorSol{i} provides a lower bound on the objective value of
\problem{$(\mathcal{F})$}. If \vectorSol{i} is feasible for
\problem{$(\mathcal{F})$}, then it also provides an upper bound. As
the lower bound and the upper bound are the same, the solution must be
optimal. The algorithm therefore terminates as soon as an iteration
yields a feasible solution to \problem{$(\mathcal{F})$}.

\subsection{Convergence}
The algorithm converges to an optimal solution to \problem{$(\mathcal{F})$} in a finite number of steps. In the worst
case, we enumerate the complete set of headways and obtain exactly the
relationship visualized by the step function in \autoref{fig:TradeOff}. However,
as we demonstrate in \autoref{sec:Results}, a feasible solution is
often found early in the search process, and the algorithm
terminates before having enumerated the full set of headways. The
efficiency of the proposed algorithm depends on the number of
iterations that are needed and the computational effort required to
solve the MILP model at each iteration.  Since frequencies are
discrete, the worst-case number is bounded by the number of lines and
frequencies, as formally described by Theorem~\ref{thm:finite} and proved in \ref{sec:proof2}.
\begin{theorem}
\label{thm:finite}
  Assuming an optimal solution to \problem{$(\mathcal{F})$} exists, then
the algorithm converges to an optimal solution in a finite number of
steps bounded by $\sum_{l\in\linePool}|\frequencySet_l|$.
\begin{proof}
  See~\ref{sec:proof2}
\end{proof}
\end{theorem}

\subsection{Valid Inequalities}\label{subsec:ineq}
The main idea of the \frameworkNameShortLP is to represent the line
frequencies using an idealized lower bound that avoids the explicit enumeration of all frequencies and paths. In
a given iteration $i$, a path $p$ associated with a
frequency $f\in \frequencySet^i$ is assigned a cost, which by construction
is less than or equal to the true cost $c_p$ of the path in \problem{($\frequencySet$)}. 
This can lead to passengers being assigned paths where the true cost would exceed $\wttThresOD$, meaning OD pair $\od\in\ODset$ 
would not travel in the public transport system, resulting in weak lower bounds, as we not only underestimate the path
cost but also misrepresent the demand captured. To address this, we introduce a family of valid inequalities that ensure that selected paths are operated at
a sufficiently high frequency. 

Consider a path $p\in \PsetOD \setminus \altModePath_\od$ for OD pair $\od \in \ODset$ 
using lines $\linePool(p)\subseteq \linePool$. If a single line is used, we compute the minimum service level, i.e., the
maximum allowable headway, that ensures that $c_p$ is within the threshold $\wttThresOD$. Since the path cost moves in steps
depending on the headway $h$ of line $l$, we define $h_{lp}^*$ as the highest headway where $c_p \leq \wttThresOD$. 
The corresponding vehicle requirements $\vlpVI = f_v(h_{lp}^*)$ specify the investment required to operate an acceptable service for $\od$.  
For paths with multiple lines, $h_{lp}^*$ and  $\vlpVI$ are specified per line, assuming that the others operate at minimum headway.  The constraints are defined as follows:
\vspace{-0.5\baselineskip}
\begin{align}
    \vlpVI \delta_p  &\leq z_l &&\forall  d\in \ODset, p \in \PsetOD \setminus \altModePath_\od, l \in \linePool(p) \label{eq:validIneq} \\
    x_{pd} &\leq \delta_p &&\forall d\in \ODset, p \in \PsetOD\setminus \altModePath_\od \label{eq:delta}
\end{align}
Here, $\delta_p$ is activated if path $p$ is
selected. Since $x_{pd}\in[0,1]$ are continuous assignment variables, the auxiliary binary variables $\delta$ ensure that \eqref{eq:validIneq} remain valid when $x_{pd}<1$.

When $\vlpVI = \phi_{l,\min(\mathcal{F}_l)}$, the constraints are redundant, as even the minimum frequency allows an
acceptable service. However, if $\vlpVI > \phi_{l,\min(\mathcal{F}_l)}$, constraints~\eqref{eq:validIneq} and~\eqref{eq:delta} ensure
that enough vehicles are allocated such that the path is a valid part of the choice set. 

Although not strictly necessary for convergence, these valid inequalities improve the algorithm’s performance by tightening bounds and reducing the number of iterations, as demonstrated in \autoref{sec:ComputationalConsiderations}.  However, as they significantly increase model size, a practical strategy is to include only a carefully selected set.  Specifically, we introduce a headway threshold $h_{T}$ and include only those constraints where $h_{lp}^*\leq h_T$, that is, where the frequency requirement is particularly high.

%% file: 05CaseSetup.tex
\section{Case Studies and Experimental Design}\label{sec:ExpDesign} 

We now describe the instances, parameters, and the experimental design of our numerical study.

\subsection{Data Sets} \label{sec:instances}
\begin{figure}[htbp]
\centering
\subfloat[\sioux]{\includegraphics[width=0.23\textwidth]{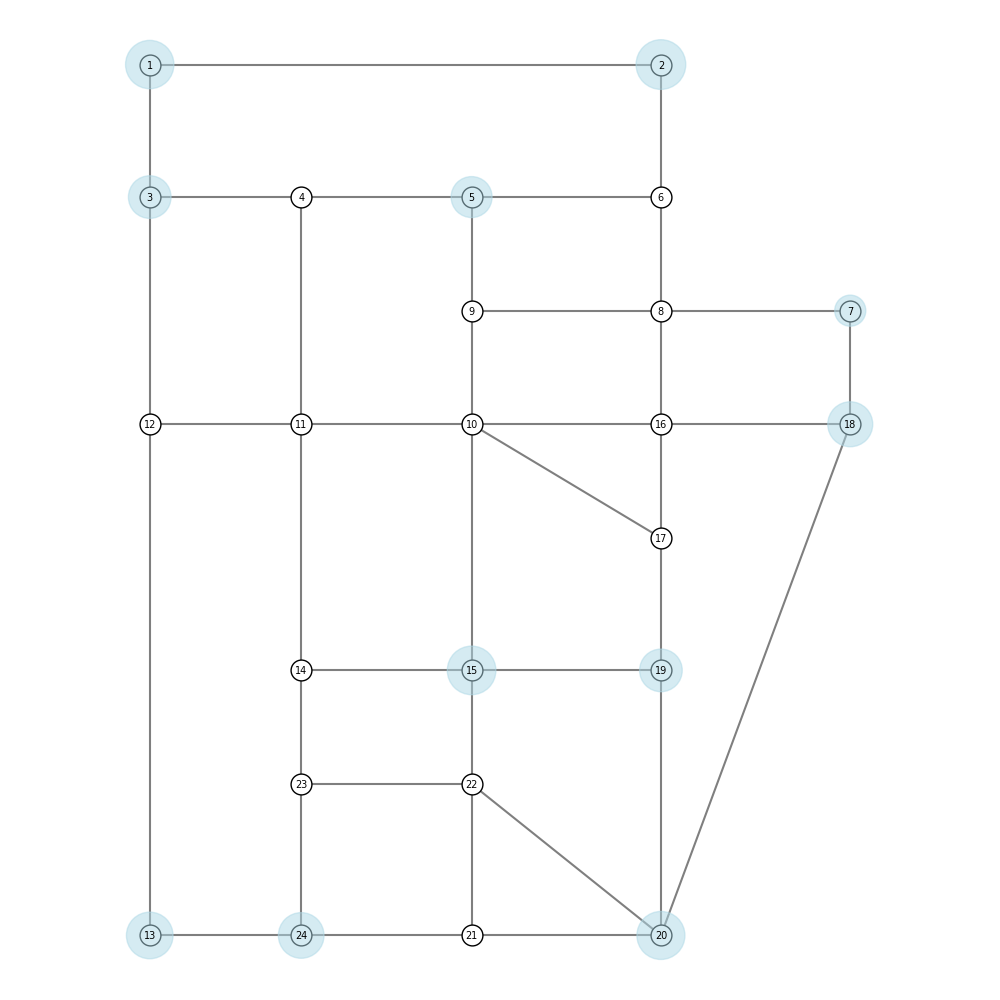}}
\subfloat[\grid]{\includegraphics[width=0.20\textwidth]{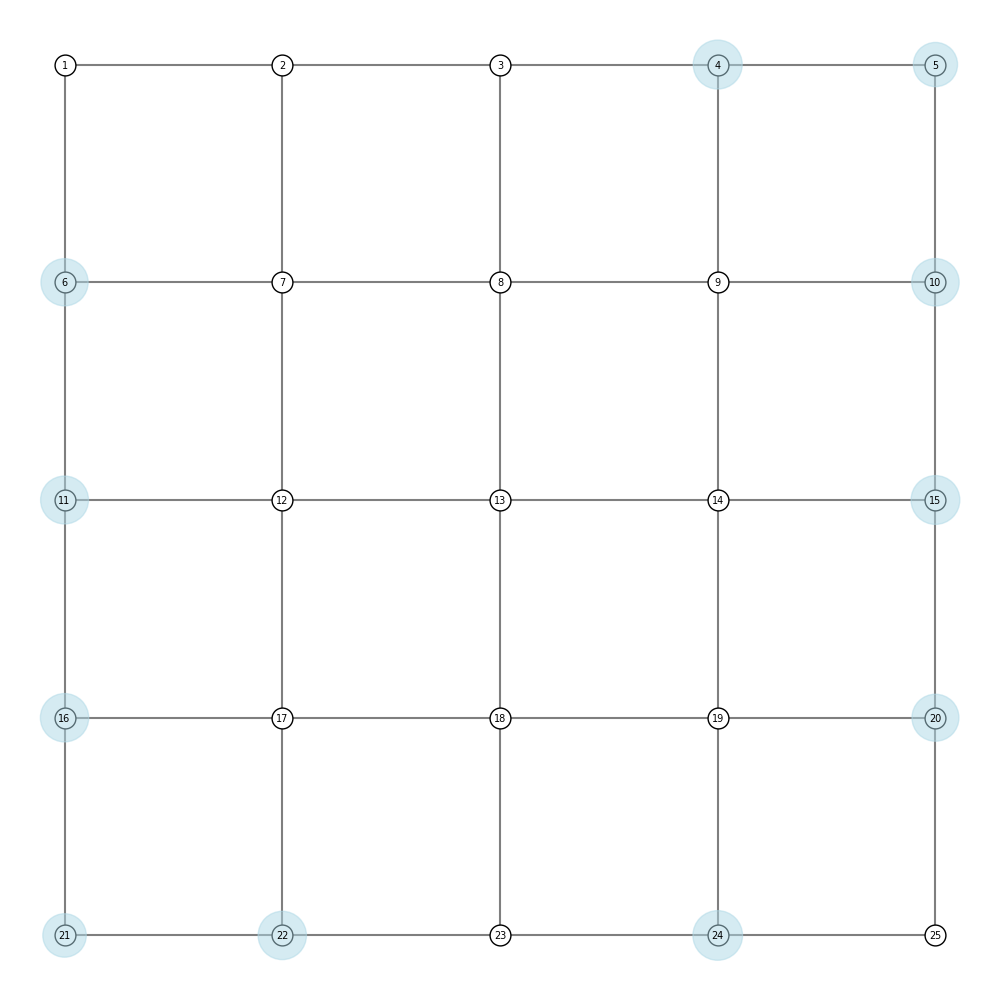}}
\subfloat[\lowersaxony]{\includegraphics[width=0.25\textwidth]{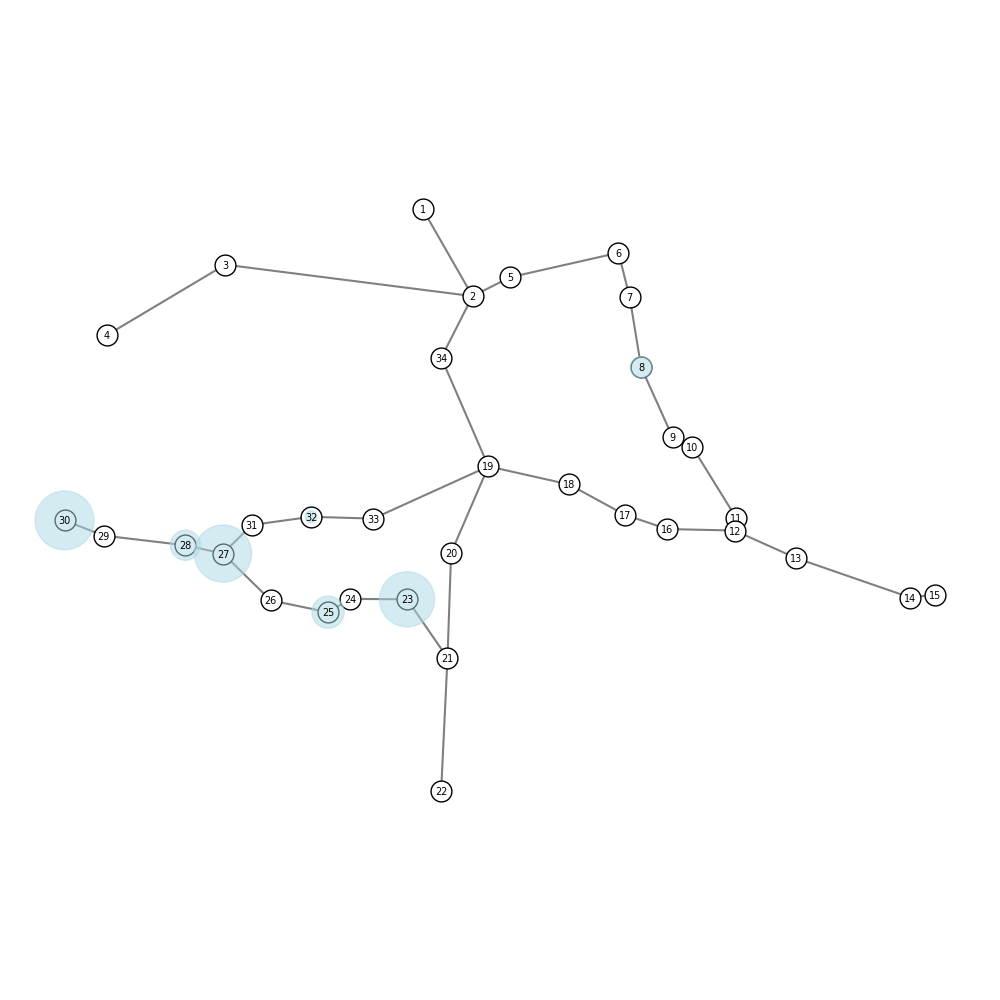}}
\subfloat[\athens]{\includegraphics[width=0.25\textwidth]{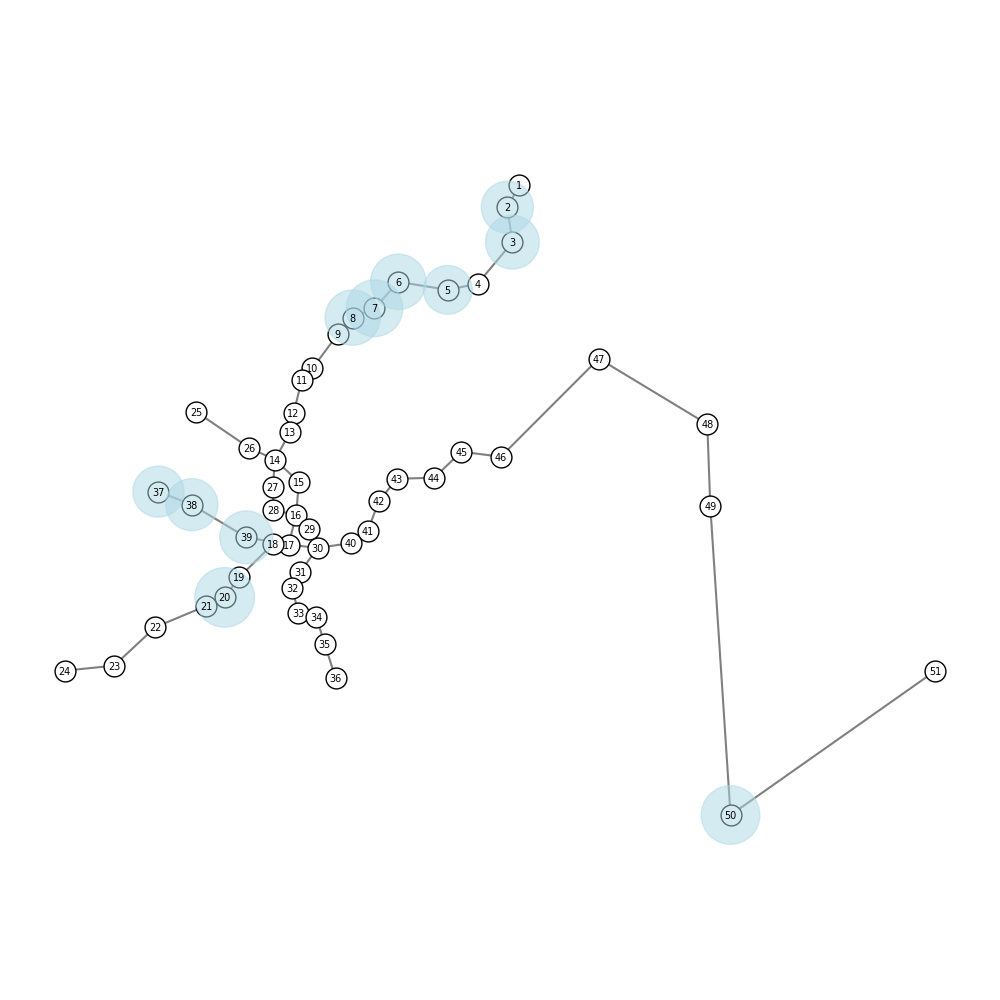}}
\caption{PTNs from LinTim. The ten stops with the highest demand for transportation are highlighted.}\label{networksWithDemand}
\end{figure}

We consider instances based on four data sets available through the LinTim toolbox \autocite{SchieweSchoebelJaegeretal.2024}.
These correspond to one \added{well-established benchmark} network (\grid) and three networks \added{based on} real-life public transport systems \added{and demand data}: the bus network in the city Sioux Falls (\sioux), the regional railway network of the German region of Lower Saxony (\lowersaxony), and the metro network in Athens (\athens). 
\autoref{networksWithDemand} illustrates the underlying structure of the PTNs while \autoref{tab:networks} summarizes their size.
It can be seen that the \grid and \sioux networks have highly symmetrical grid-like structures, while the other two are more star-shaped. 
\added{The studied networks have been widely used in the literature to benchmark cutting-edge methods for transit network design and related problems.
The \grid network is considered by, among others, \textcite{friedrich_integrating_2017}, \textcite{schiewe_line_2019}, and \textcite{schiewe_periodic_2020}, while the \sioux network is considered by \textcite{gattSolvingLinePlanning2025} and \textcite{borndorfer_direct_2012}.
	\textcite{schiewe_periodic_2020} also consider the \lowersaxony and \athens network, which also come back in various recent studies (e.g., \textcite{grafe2021solving}).}

In addition to the PTN and average travel times of the edges, each data set also provides travel demand in the form of OD matrix \ODset.  
The OD data is given as an estimated hourly demand for transportation for the \sioux  (4,115 travelers over 552 ODs) and \grid (2,546 travelers over 567 ODs) instances. 
However, for the \lowersaxony and \athens instances, demand is given on a more aggregate level.
Therefore, we scale the demand by a factor of $1/20$ (from 63,323 to 2,957 travelers over 2,385 ODs) and $1/200$ (from 359,818 to 1,799 travelers over 395 ODs) for the \lowersaxony and \athens data sets, respectively.

\subsection{Instances and Parameters}
The \problem{} relies on the fact that a line pool $\linePool$, a set of possible frequencies $\frequencySet$, and a candidate set of paths for each  OD pair $\od \in \ODset$ is available. 
Furthermore, operational cost and capacity parameters need to be specified.
To obtain a variety of instances for the experiments, we consider different subsets of ODs and line pools for each data set. 
By considering the four data sets, three OD sets, and three line pools we end up with a total of 36 instances. 

For each network, we consider OD sets of three different sizes. 
The available OD pairs are sorted according to their demand $w_{\od}$, and we form sets containing the 100, 200, and 300 OD pairs with the highest demand. 
\autoref{tab:networks} outlines the total hourly travel demand covered by each OD set and the fraction of the original demand this represents. 
While the selection of a subset of ODs, in our case, is primarily for creating a varied set of instances, we note that many  PTNs exhibit a strong inequality in the demand distribution where a few ODs account for a majority of the trips.
Hence, considering a selected set of ODs based on travel demand or using another OD matrix decomposition technique \autocite[see, e.g.,][]{schiewe_periodic_2020,hurk_shuttle_2016}  can provide good networks that serve most travelers effectively. 

We generate three line pools for each network. 
Here, we use the line pool generation heuristic proposed by \textcite{gattermann_line_2017}, that is available in LinTim, to generate a first set of lines. This heuristic has been shown to produce line pools that are not too large and lead to good line plans. 
To generate line pools of different sizes, we also apply a k-shortest path heuristic that generates lines such that all OD pairs should be able to travel on their shortest path. 
We selected lines from the line pools generated by both methods in different ways to obtain a small, medium, and large line pool for each instance. Their sizes are reported in \autoref{tab:networks}.
Each line can be opened at different frequencies, which can differ per mode or specific line. 
In our experiments, we consider all integer headways between two and 20 minutes, i.e., $\frequencySet_{l} = \{3, \frac{60}{19},\dots,20,30\}$ for each $l \in \linePool$.

\begin{table}[htbp]
\centering
\caption{Problem characteristics in terms of stops $|\stopSet|$, (undirected) edges $|E|$, volume of travelers covered by the different OD subsets, and number of candidate lines $|\linePool|$.}
\begin{tabular}{@{}lllllllll@{}} \toprule
\textbf{Data set} & \multicolumn{2}{l}{\textbf{PTN}} & \multicolumn{3}{l}{\textbf{Hourly travel demand (\% of total)}} & \multicolumn{3}{l}{\textbf{Line pool size}} \\ 
\cmidrule(l){2-3} \cmidrule(l){4-6} \cmidrule(l){7-9} 
& $|\stopSet|$ & $|E|$ & $|\ODset| = 100$ & $|\ODset|=200$ & $|\ODset| = 300$ & S & M & L \\ \midrule
\sioux & 24 & 38 & 2272 (55.2\%)  & 3238 (78.7\%)  & 3691 (89.6\%)  & 23 & 46 & 92 \\
\grid & 25 & 40 & 1524 (59.9\%) & 1955 (76.8\%)  & 2213 (86.9\%)  & 30 & 60 & 123 \\
\athens & 51 & 52 & 915 (31.0\%) & 1264 (42.8\%) & 1502 (50.8\%) & 32 & 64 & 129 \\
\lowersaxony & 34 & 35 & 1525 (84.7\%) & 1713 (95.2\%) & 1780 (98.9\%) & 41 & 82 & 120 \\ \bottomrule
\end{tabular}
\label{tab:networks}
\end{table}
\vspace{-0.5\baselineskip}

\subsubsection{Passenger Costs} \label{pasCosts}
The first element of the objective function \eqref{obj} concerns the passenger costs.
For a single passenger in OD pair $\od\in \ODset$ assigned to a path $p\in \PsetOD$, this cost is defined through the path cost $c_{p}$.
The path cost is decomposed into the cost of the individual arcs as shown in \eqref{eq: cost allocation}, allowing to weigh the elements of initial waiting, in-vehicle, and transfer time within a journey independently.
In our experiments, we let the cost $c_p$ represent the monetary value of travel time for a path $p$, where we choose the parameters $c^{A}_{\arc}, c^{IVT}_{\arc}, c^{T}_{\arc},$ and $c^{E}_{\arc}$ in line with the value proposed by the TERESA model for the value of travel time in Denmark \autocite{incentive__dtu_transport_teresa_2024}.
In particular, we let each minute of in-vehicle traveling time on a line arc $a \in \arcSetDrive$ corresponds to a cost of $\frac{119}{60}$ Danish kroner (DKK). We assume the expected waiting time of an access arc $a \in \arcSetAcces$ to be half of the headway of the outgoing line, thus assuming an equal distribution of arriving passengers over time. Studies \autocite{ingvardson_passenger_2018} suggest that passengers tend to time their arrival when headways are above a certain level. Based on this, we treat the first five minutes as perceived waiting time, assigning a weight of DKK $\frac{238}{60}$ per minute. Any additional waiting time is considered hidden waiting time, weighted at DKK $\frac{95}{60}$ per minute.
The cost of any transfer arc $a \in \arcSetTransfer$ includes a fixed penalty and the expected time waiting for the next vehicle, which we again choose to be half the headway of the outgoing line due to the absence of timetable information at this planning stage.
We use a fixed transfer penalty of 12 DKK and a value of $\frac{179}{60}$ DKK for each minute of transfer time. 
In our experiments, we do not put a weight on the egress arcs \arcSetEgress.

\subsubsection{Passenger Path Generation} \label{pasRoutes}
Determining the path set $\Pset$ is a crucial step since it defines which paths are seen as acceptable by passengers in the model and, hence, if passengers are willing to travel in the network given a set of operated lines and associated frequencies.
To determine the acceptable paths, we set a threshold relative to the shortest possible journey time $\minIVTd$ for each OD pair $\od \in \ODset$, defined by the shortest path between origin $s_d$ and destination stop $t_d$ in the PTN.
Note that as $\minIVTd$ is derived from the PTN, it does not include transfer and waiting time.

The main path set $\Pset$ in our numerical experiments reflects \elasticDemand{}.
Here, we see a path as acceptable for an OD pair $\od\in\ODset$ if its weighted travel time is within threshold $\wttThresOD$, which we define as
$\wttThresOD = \min(3\minIVTd, 1.25\minIVTd + \tau)$, 
where  $\tau$ represents the combined cost of a single transfer and the initial waiting time. The value of $\tau$ is computed under a worst-case scenario for the passengers in which the lines operate at the minimum frequency.
This setting was determined through experimentation but can be based on available passenger alternatives in case studies.

We also consider a second benchmark path set $\PsetOD^{\textrm{rigid}}$ in our experiments, which resembles a setting with ``rigid demand'' in which passengers are almost always willing to travel in the network.
Here, we assume that if an ordered set of edges  $(e_1, e_2,..., e_m) \in E$ is deemed an acceptable path for OD $\od$ at any combination of frequencies of the lines traveled on, then all frequency combinations $(e_{l_1f_1}, e_{l_2f_2},...,  e_{l_nf_m}) \in \arcSet \text{ with } f_i\in \frequencySet_{l}, f_i=f_j \text{ for }l_i=l_j $ are also part of the path set.
A path for OD $\od \in \ODset$ is acceptable if it itself, or any other frequency combination, is within a tolerance coefficient $\alpha$ of the shortest travel time $\minIVTd$. We define $\alpha$ to allow only paths within 15 minutes of \minIVTd and permitting a single transfer and initial waiting time. 
In this calculation, the cost of transferring and waiting is calculated assuming the highest frequency of a line.

We apply a $k$-shortest path algorithm for each OD pair with a sufficiently high $k$ to ensure all simple paths in the CGN with a cost within the threshold are found. 
For the rigid path set $\PsetOD^{\textrm{rigid}}$, where paths are included irrespective of frequency, we first generate the acceptable paths assuming the highest frequency for each line and then create all other frequency combination versions of that path.
Moreover, we further process both path sets to remove impractical or redundant paths.
For example, if a direct line exists between $\ODor$ and $\ODde$ for an OD $\od$, transferring along the same line is unnecessary, provided the direct path does not involve a significantly longer detour. 
To further reduce the path set, we introduce the following dominance rule. Considering two paths $r_i$ and $r_i$, we say that $r_i$ dominates $r_j$ if $r_i$ has equal or lower cost than $r_j$, $r_i$ uses a subset of the lines used by $r_j$, and $r_i$ uses the same frequencies on the subset of common lines. 
In addition, we restrict the generation of paths to only include paths using a single transfer. 
This is for both practical reasons, as it has often been observed that having more than one transfer can deter people from using public transport, and computational reasons, as the number of paths increases quickly when considering more than one transfer.
We note that this assumption follows common practice (see, e.g., \cite{borndorfer_direct_2012, bertsimas_data-driven_2021, hartleb_modeling_2023}).

\begin{table}[htbp]
\centering
\caption{Average number of paths per OD for a given line pool $\linePool$}
\begin{tabular}{@{}lcccc@{}}
\toprule
\textbf{} & \multicolumn{1}{c}{\begin{tabular}[c]{@{}c@{}}\sioux \\ $|\linePool|=23$\end{tabular}} & \multicolumn{1}{c}{\begin{tabular}[c]{@{}c@{}}\grid \\ $|\linePool|=30$\end{tabular}} & \multicolumn{1}{c}{\begin{tabular}[c]{@{}c@{}}\athens \\ $|\linePool|=32$\end{tabular}} & \multicolumn{1}{c}{\begin{tabular}[c]{@{}c@{}}\lowersaxony \\ $|\linePool|=41$\end{tabular}}  \\ \midrule
Unique frequency-independent paths & 9 & 9 & 20 & 36 \\
Number of paths in $\PsetOD$ & 1,085 & 1,077 & 778 &  3,742\\
Number of paths in $\PsetOD^{\textrm{rigid}}$ & 1,235 & 1,130 & 2,385 & 5,020 \\
 \bottomrule
\end{tabular}
\label{tab:PathSetSizeSmallLinePool}
\end{table}

Lastly, for each OD $\od\in\ODset$ we include a path $\altModePath_\od$ in both sets $\PsetOD^{\textrm{rigid}}$ and $\PsetOD$ that represents the passenger not being able to or choosing not to travel in the PTN. This path is not associated with any line or frequency and has a cost equal to the cost threshold in either set. 
The resulting size of the path sets is given in \autoref{tab:PathSetSizeSmallLinePool}, which denotes the number of paths in both path sets as well as the number of unique frequency-independent paths when not considering the frequency expansion in the CGN.
Note how the number of paths increases sharply through the consideration of frequencies and how the size of the path set $\PsetOD$ is always smaller than the path set $\PsetOD^{\textrm{rigid}}$.

\subsubsection{Operating Costs and Capacity}  
The second part of the objective considers the operating cost, which is the difference between the cost of operating the selected lines at the selected frequencies and the ticket revenue.
First, this includes the cost of the operated vehicles, where we consider a single vehicle type with a capacity of $\delta_{v} = 50$ passengers and an hourly cost of 880 DKK, which is based on the average hourly cost of operating a bus for the Danish Traffic Companies \autocite{nielsen_region_2024}. 
Second, we consider the cost of opening a line, representing the management and overhead cost associated to a line, which we set equal to the cost of one vehicle (880 DKK).
Moreover, we assume that each transported passenger generates a ticket fare of 22 DKK, which aligns with ticket fares within Danish cities.

Since planners typically operate within predefined resource constraints, the model includes a budget on overall operating costs.  However, since no predefined budget is available for the networks we study, we base the operating budget on the operational cost that we get for the problem with the rigid path set ($\PsetOD^{\textrm{rigid}}$) without enforcing a budget constraint. 
Since the path set, and hence the solution space is larger, the solution is a lower bound on the problem with frequency-dependent paths.
We then use the operating costs in the best-found solution related to operating the lines, i.e., the sum of vehicle and line cost, as input to the model when solving it with the path set $\PsetOD$.

\begin{table}[htbp]
\centering
\begin{minipage}[t]{0.45\linewidth} 
\caption{Operating costs and capacities}
\centering
\begin{tabular}{@{}ll@{}}
\toprule
Parameter & Value\\ \midrule
Vehicle cost $c_v$ (DKK/hour) & 880 \\
Line cost $h_{lf}$ (DKK) & $1c_v$ \\
Vehicle capacity $\delta_v$ (Pas./vehicle) & 50 \\
Ticket fare $R$ (DKK/pas.) & 22 \\
\bottomrule
\end{tabular}\label{tab:parameterSettingOC} 
\end{minipage}
\hspace{0.05\linewidth} 
\begin{minipage}[t]{0.45\linewidth} 
\caption{Passenger costs}
\centering
\begin{tabular}{@{}ll@{}} 
\toprule
Parameter & Value \\ \midrule
Fixed transfer penalty (DKK) & 12 \\
In-vehicle time (DKK/hour) & 119 \\
Waiting time (DKK/hour) & 238 \\
Hidden waiting time (DKK/hour) & 95 \\
Transferring time (DKK/hour) & 179 \\
\bottomrule
\end{tabular} \label{tab:parameterSettingPI}
\end{minipage}
\end{table}

\autoref{tab:parameterSettingOC} and \autoref{tab:parameterSettingPI} summarize the key parameters. While using the same parameter settings across different networks is a simplification, our focus is to demonstrate the algorithm’s functionality and the effect of considering lost demand rather than optimizing each network individually.

%% file: 06Results.tex
\section{Computational Results}\label{sec:Results}
We evaluate the proposed model and solution method on the 36 instances described in \autoref{sec:ExpDesign} 
for $\lambda\in\{0.1, 0.25, 0.5, 0.75, 1.0\}$. 
\autoref{sec:Computational} analyzes the computational performance and scalability of the \frameworkNameShortLP.  
\autoref{sec:DemandElasticity} examines the impact of \elasticDemand with respect to passenger service, the share of demand captured, and the total cost. 
The MILPs in the \frameworkNameShortLP are implemented with a restricted subset of the valid inequalities (\ref{eq:validIneq})--(\ref{eq:delta}), using a threshold $h_T$ of 10 minutes -- this setting is further explored in \autoref{sec:ComputationalConsiderations} on algorithm acceleration strategies. 
Finally, \autoref{sec:RouteAssignment} compares the model’s passenger assignment assumptions with alternative approaches.

All experiments are conducted on an 8-core Intel Xeon Gold 6226R processor with 200GB of internal memory allowing up to ten hours to solve each instance. The \frameworkNameShortLP is implemented in Java version 17.0.6 (OpenJDK 17) and all optimization models were solved using the commercial solver IBM ILOG CPLEX Optimizer 22.1.1.

\input{06_ComputationScalability}
\input{06_DemandElasticity}
\input{06_AccelerationTechniques}
\input{06_RouteAssignment}

%% file: 06_ComputationScalability.tex
\subsection{Computational Performance} \label{sec:Computational}
\begin{figure}[htbp]
	\centering
		\subfloat[t][Average computational time (lines, left axis) and number of instances solved to optimality (bars, right axis)]{
        \resizebox{!}{0.32\textwidth}{\input{figures/Comp/TimeVSLambda_withCounts.pgf}}\label{fig:TimeVsL}
        }        \hspace{\fill}
		\subfloat[t][Average passenger and operator costs (excluding revenue)]{\resizebox{!}{0.32\textwidth}
{\input{figures/Comp/PIandOCVSLambda.pgf}}
			\label{fig:costVsL}} 
		\caption{Results for different values of $\lambda$}
	\label{fig:VS_lambda}
\end{figure}
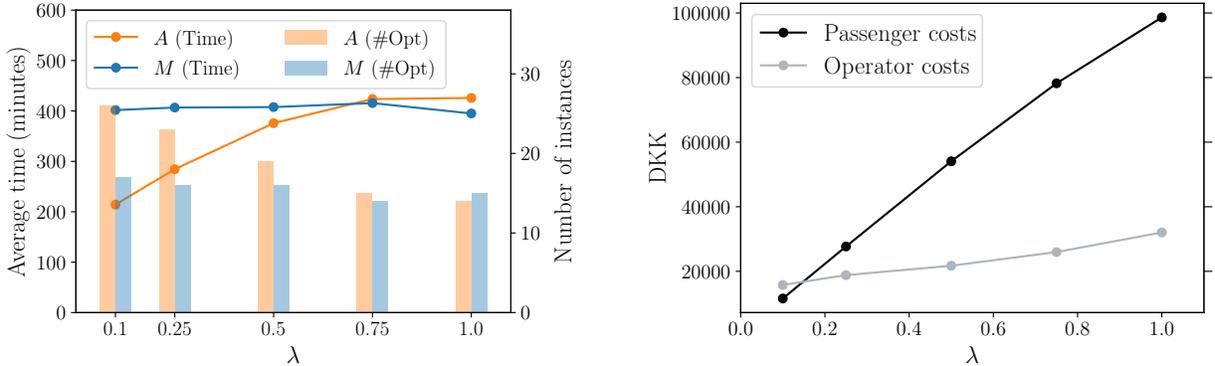

We first evaluate the computational efficiency of the proposed \frameworkNameShortLP by comparing its performance with that of directly solving the full model with CPLEX. 
For brevity, we refer to the former as {\em A} and the latter as {\em M}. 
\autoref{fig:TimeVsL} provides an overview of the results and illustrates the average computation time and the number of instances solved to optimality for each of the considered $\lambda$ values. This includes instances that are terminated at the time limit, affecting the reported average.  Detailed results on individual instances can be found in \autoref{sec:TablesSup}. 

Recall that the refinement step in
\frameworkNameShortLP is designed under the assumption that capacity -- and, consequently, operator costs -- are important for the solution. 
\autoref{fig:costVsL} shows that passenger costs tend to dominate
operator costs; here, $\lambda$ can act as a tuning parameter that allows public transport agencies to adjust the trade-off between the two according to their priorities. 
In \autoref{sec:appSocialOptimum}, we provide a more detailed breakdown of the objective components across different $\lambda$ values.

\autoref{fig:TimeVsL} confirms that the \frameworkNameShortLP\ becomes faster on average as $\lambda$ decreases, while {\em M} is largely unaffected by $\lambda$. 
The gains in performance are mainly driven by the increased effectiveness of the frequency-refinement step of the \frameworkNameShortLP\ at lower $\lambda$ levels, which enables more instances to be solved to optimality. While $M$ solves to optimality respectively 17, 16, 16, 14, and 15 instances for $\lambda$ values 0.1, 0.25, 0.5, 0.75, and 1, these counts for \frameworkNameShortLP are 26, 23, 19, 15, and 14. This refinement step assumes that the assigned frequency in an iteration is a good estimate for the optimal frequency, which is less often the case for high $\lambda$. Future research could investigate if alternative approaches to refining the frequency representation could increase the speed of the algorithm for these instances.

\begin{figure}[htbp]
	\centering
		\subfloat[{\rev{Number of instances solved to optimality over time}}]{\includegraphics[width=0.45\textwidth]{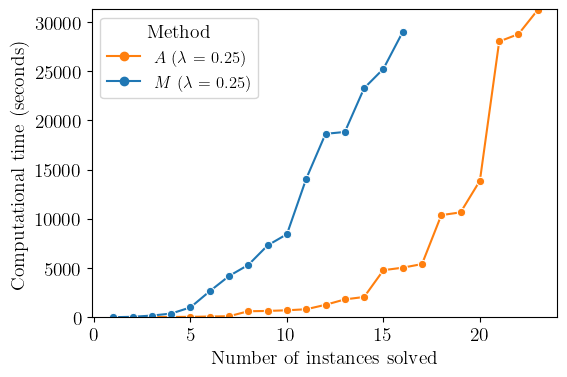}\label{fig:NrInst}}
		\subfloat[{\rev{Number of instances with a certain  gap or better at termination}}]{\includegraphics[width=0.45\textwidth]{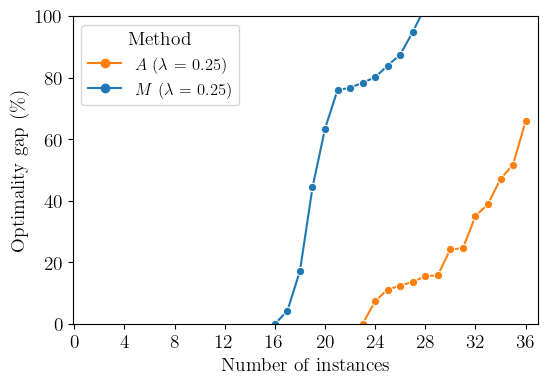}
			\label{fig:Gaps}} \\
		\caption{\rev{Results for $\lambda = 0.25$ }}
	\label{fig:InstOverTime}
\end{figure}

To analyze the performance of the \frameworkNameShortLP under a more balanced weighting of passenger and operator costs, we focus on a setting with $\lambda=0.25$. 
\autoref{fig:NrInst} shows the number of instances solved to optimality in relation to the computation time needed. Not only does the \frameworkNameShortLP solve more instances to optimality at this setting, but it does so much faster.
The number of solved instances is significantly higher for \lowersaxony, where the algorithm is able to solve all but two of the instances, while the full model can only solve three of the instances. 
\autoref{fig:Gaps} highlights
that $A$ also consistently provides tighter bounds
than $M$; it shows the cumulative number of instances that achieve a certain optimality gap or better within the time limit. 
For $A$, 27 instances have gaps within 15$\%$, while only two have gaps higher than 50\% (both for \sioux). 
In contrast, when $M$ fails to find an optimal solution within the time limit, the resulting bounds are generally weaker. 
Only 17 instances achieve a gap below 15$\%$, and in four cases, $M$ fails to compute any reasonable lower bound. 
In these instances, the model assigns no passengers to the public transport system, resulting in technically feasible but practically meaningless solutions. 
Similar plots for the other $\lambda$ values are provided in \autoref{sec:PlotSup}. These show that $A$ achieves faster solution times primarily at lower $\lambda$ values and consistently produces stronger bounds across the entire range of $\lambda$.

To obtain further insights into the performance of the methods over the different instances, they can be classified into three groups: 1) instances that can be solved to optimality by either method, 2) instances where $A$ solves the problem to optimality while $M$ reaches the time limit, and 3) instances where both methods reach the time limit.  
\autoref{tab:CompAverage} summarizes the results for each of these groups \added{in terms of computation time, optimality gap, and size of the models.
The relative difference} between the respective upper bounds obtained with each method \added{(RUB)} is computed as $\frac{z_{UB}^M -
  z_{UB}^A}{|z_{UB}^M|}$, where $z_{UB}^{i}$ denotes the upper bound obtained with method $i\in\{A,M\}$. 
 When the relative upper bound is greater than zero, it indicates that $A$ provides a tighter bound. 
\added{The ``Speed-up'' column reports the average per-instance speed-up between the two methods.}
For the 16 instances that can be solved using either method, the algorithm requires, on average, 61.9 minutes, compared to 165.3 minutes for $M$, where the largest speed-up occurs for the \lowersaxony data set (1.8 vs 314 minutes).
This group primarily consists of instances with smaller line pools and OD sets. 
For the instances solved to optimality only by $A$, $M$ produces relatively poor bounds, with an average optimality gap of $91.5\%$. 
We note that the difference in bounds is mostly related to an improvement in the lower bound, as the relative difference in upper bound is only 3.3\%.  
The last group consists of 13 instances that predominantly use either the largest line pool or the largest set of ODs. For these instances, $A$ provides solutions within 28\% of optimality on average compared to $>$100\% for $M$. Notably, both lower and upper bounds show substantial improvement with $A$. 

The performance difference can be largely attributed to $A$ generating only a fraction of the variables and constraints in $M$. 
\autoref{tab:CompAverage} reports the average percentage of variables and constraints generated by $A$ relative to $M$ in the final iteration\added{, as shown in the Var.\@ (\%) and Const.\@ (\%) columns}.
The most significant reduction in the number of variables is observed for the \lowersaxony instances, where only 1.9--3.5\% of variables are generated. This reduction likely contributes to the relatively large performance gap between the methods in this group of instances. 
In contrast, up to 19.7\% of variables are generated by $A$ for the \sioux instances. Notably, this data set contains more instances (especially with larger line pools) that \emph{A} is unable to solve to optimality. 
Regarding the number of constraints, the reduction is somewhat smaller for the \athens instances, where the \frameworkNameShortLP in a single case even generates more constraints than $M$ due to the inclusion of the valid inequalities.  
\begin{table}[htbp]
\centering
\caption{Average computational performance for instances solved by both $A$ and $M$, by either, or by both methods.}
\label{tab:CompAverage}\resizebox{0.95\linewidth}{!}{
\begin{tabular}{@{}lrrrrrrrrr@{}}
\toprule
 & \multicolumn{1}{c}{} & \multicolumn{3}{c}{Time (minutes)} & \multicolumn{2}{c}{Gap (\%)}& \multicolumn{1}{c}{RUB ($\%$)} & \multicolumn{1}{c}{Var.\@ ($\%$)}& \multicolumn{1}{c}{Const.\@ ($\%$)}\\   \cmidrule(l){3-5} \cmidrule(l){6-7} \cmidrule(l){8-8}\cmidrule(l){9-9}\cmidrule(l){10-10} 
 & \multicolumn{1}{c}{\#Inst} &  \multicolumn{1}{c}{$A$} & \multicolumn{1}{c}{$M$} & \multicolumn{1}{c}{Speed-up} & \multicolumn{1}{c}{$A$} & \multicolumn{1}{c}{$M$} & \multicolumn{1}{c}{$A$}& 
 \multicolumn{1}{c}{$A$} & \multicolumn{1}{c}{$A$}  \\  \midrule
Both optimal & 16 & 61.9 & 165.3 & 33.1 & 0 & 0 & 0  & 7.0 & 52.2\\
$A$ optimal  & 7 & 207.1 & $>600$& 11.8 & 0 & 91.5 & 3.3 & 6.4 & 38.9\\
Neither optimal &13 & $>600$ &$>600$ & 1.0 & 28.0 & $>100$ & 20.1 & 4.2 & 35.6 \\\bottomrule
\end{tabular}}
\end{table}

\begin{figure}[htbp]
    \centering
    \includegraphics[width=1.0\linewidth]{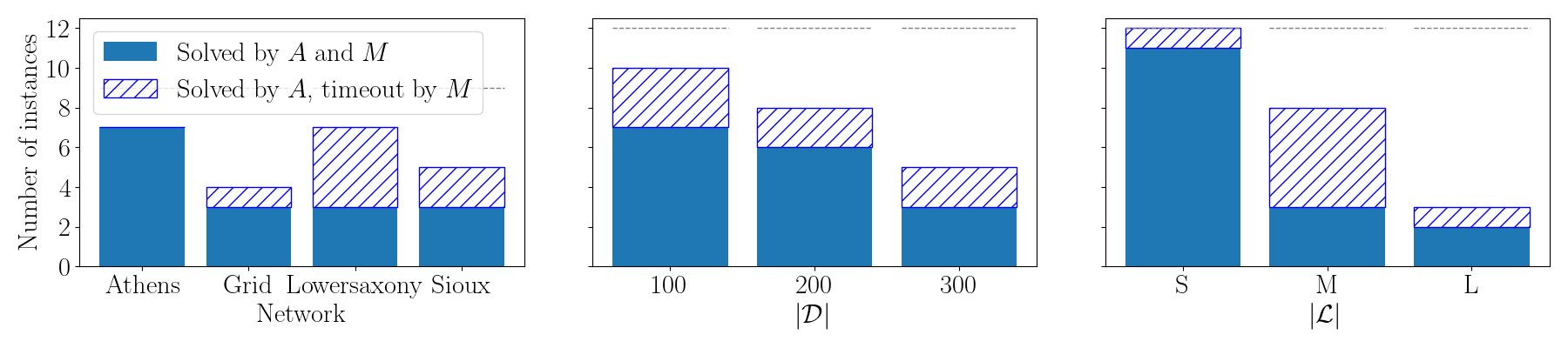}
    \caption{Number of instances solved to optimality for $\lambda = 0.25$ grouped by network, $|\ODset|$, and $|\linePool|$. The total number of instances in each group is marked with a dotted line.}
    \label{fig:DimensionsL25}
\end{figure}

We also quantify the performance in relation to the overall instance characteristics.
\autoref{fig:DimensionsL25} shows the number of instances solved to optimality over the different networks, line pools, and OD sets. 
These figures nicely show that $A$ can solve a range of problems -- both large and small -- which $M$ cannot. 
While we have focused on an in-depth analysis for $\lambda=0.25$, more notable differences in favor of $A$ are observed for $\lambda=0.1$. 
For higher values of $\lambda\in\{0.5,0.75,1\}$, the differences in computation time are less pronounced, but a consistent advantage of $A$ across all settings is its ability to deliver significantly tighter bounds. 
The respective results can be found in \autoref{sec:problemDims}. 

%% file: figures/Comp/TimeVSLambda_withCounts.pgf
\begingroup%
\makeatletter%
\begin{pgfpicture}%
\pgfpathrectangle{\pgfpointorigin}{\pgfqpoint{6.000000in}{4.000000in}}%
\pgfusepath{use as bounding box, clip}%
\begin{pgfscope}%
\pgfsetbuttcap%
\pgfsetmiterjoin%
\definecolor{currentfill}{rgb}{1.000000,1.000000,1.000000}%
\pgfsetfillcolor{currentfill}%
\pgfsetlinewidth{0.000000pt}%
\definecolor{currentstroke}{rgb}{1.000000,1.000000,1.000000}%
\pgfsetstrokecolor{currentstroke}%
\pgfsetdash{}{0pt}%
\pgfpathmoveto{\pgfqpoint{0.000000in}{0.000000in}}%
\pgfpathlineto{\pgfqpoint{6.000000in}{0.000000in}}%
\pgfpathlineto{\pgfqpoint{6.000000in}{4.000000in}}%
\pgfpathlineto{\pgfqpoint{0.000000in}{4.000000in}}%
\pgfpathlineto{\pgfqpoint{0.000000in}{0.000000in}}%
\pgfpathclose%
\pgfusepath{fill}%
\end{pgfscope}%
\begin{pgfscope}%
\pgfsetbuttcap%
\pgfsetmiterjoin%
\definecolor{currentfill}{rgb}{1.000000,1.000000,1.000000}%
\pgfsetfillcolor{currentfill}%
\pgfsetlinewidth{0.000000pt}%
\definecolor{currentstroke}{rgb}{0.000000,0.000000,0.000000}%
\pgfsetstrokecolor{currentstroke}%
\pgfsetstrokeopacity{0.000000}%
\pgfsetdash{}{0pt}%
\pgfpathmoveto{\pgfqpoint{0.868690in}{0.731763in}}%
\pgfpathlineto{\pgfqpoint{5.167416in}{0.731763in}}%
\pgfpathlineto{\pgfqpoint{5.167416in}{3.722737in}}%
\pgfpathlineto{\pgfqpoint{0.868690in}{3.722737in}}%
\pgfpathlineto{\pgfqpoint{0.868690in}{0.731763in}}%
\pgfpathclose%
\pgfusepath{fill}%
\end{pgfscope}%
\begin{pgfscope}%
\pgfsetbuttcap%
\pgfsetroundjoin%
\definecolor{currentfill}{rgb}{0.000000,0.000000,0.000000}%
\pgfsetfillcolor{currentfill}%
\pgfsetlinewidth{0.803000pt}%
\definecolor{currentstroke}{rgb}{0.000000,0.000000,0.000000}%
\pgfsetstrokecolor{currentstroke}%
\pgfsetdash{}{0pt}%
\pgfsys@defobject{currentmarker}{\pgfqpoint{0.000000in}{-0.048611in}}{\pgfqpoint{0.000000in}{0.000000in}}{%
\pgfpathmoveto{\pgfqpoint{0.000000in}{0.000000in}}%
\pgfpathlineto{\pgfqpoint{0.000000in}{-0.048611in}}%
\pgfusepath{stroke,fill}%
}%
\begin{pgfscope}%
\pgfsys@transformshift{1.259483in}{0.731763in}%
\pgfsys@useobject{currentmarker}{}%
\end{pgfscope}%
\end{pgfscope}%
\begin{pgfscope}%
\definecolor{textcolor}{rgb}{0.000000,0.000000,0.000000}%
\pgfsetstrokecolor{textcolor}%
\pgfsetfillcolor{textcolor}%
\pgftext[x=1.259483in,y=0.634541in,,top]{\color{textcolor}\rmfamily\fontsize{14.000000}{16.800000}\selectfont 0.1}%
\end{pgfscope}%
\begin{pgfscope}%
\pgfsetbuttcap%
\pgfsetroundjoin%
\definecolor{currentfill}{rgb}{0.000000,0.000000,0.000000}%
\pgfsetfillcolor{currentfill}%
\pgfsetlinewidth{0.803000pt}%
\definecolor{currentstroke}{rgb}{0.000000,0.000000,0.000000}%
\pgfsetstrokecolor{currentstroke}%
\pgfsetdash{}{0pt}%
\pgfsys@defobject{currentmarker}{\pgfqpoint{0.000000in}{-0.048611in}}{\pgfqpoint{0.000000in}{0.000000in}}{%
\pgfpathmoveto{\pgfqpoint{0.000000in}{0.000000in}}%
\pgfpathlineto{\pgfqpoint{0.000000in}{-0.048611in}}%
\pgfusepath{stroke,fill}%
}%
\begin{pgfscope}%
\pgfsys@transformshift{1.845673in}{0.731763in}%
\pgfsys@useobject{currentmarker}{}%
\end{pgfscope}%
\end{pgfscope}%
\begin{pgfscope}%
\definecolor{textcolor}{rgb}{0.000000,0.000000,0.000000}%
\pgfsetstrokecolor{textcolor}%
\pgfsetfillcolor{textcolor}%
\pgftext[x=1.845673in,y=0.634541in,,top]{\color{textcolor}\rmfamily\fontsize{14.000000}{16.800000}\selectfont 0.25}%
\end{pgfscope}%
\begin{pgfscope}%
\pgfsetbuttcap%
\pgfsetroundjoin%
\definecolor{currentfill}{rgb}{0.000000,0.000000,0.000000}%
\pgfsetfillcolor{currentfill}%
\pgfsetlinewidth{0.803000pt}%
\definecolor{currentstroke}{rgb}{0.000000,0.000000,0.000000}%
\pgfsetstrokecolor{currentstroke}%
\pgfsetdash{}{0pt}%
\pgfsys@defobject{currentmarker}{\pgfqpoint{0.000000in}{-0.048611in}}{\pgfqpoint{0.000000in}{0.000000in}}{%
\pgfpathmoveto{\pgfqpoint{0.000000in}{0.000000in}}%
\pgfpathlineto{\pgfqpoint{0.000000in}{-0.048611in}}%
\pgfusepath{stroke,fill}%
}%
\begin{pgfscope}%
\pgfsys@transformshift{2.822656in}{0.731763in}%
\pgfsys@useobject{currentmarker}{}%
\end{pgfscope}%
\end{pgfscope}%
\begin{pgfscope}%
\definecolor{textcolor}{rgb}{0.000000,0.000000,0.000000}%
\pgfsetstrokecolor{textcolor}%
\pgfsetfillcolor{textcolor}%
\pgftext[x=2.822656in,y=0.634541in,,top]{\color{textcolor}\rmfamily\fontsize{14.000000}{16.800000}\selectfont 0.5}%
\end{pgfscope}%
\begin{pgfscope}%
\pgfsetbuttcap%
\pgfsetroundjoin%
\definecolor{currentfill}{rgb}{0.000000,0.000000,0.000000}%
\pgfsetfillcolor{currentfill}%
\pgfsetlinewidth{0.803000pt}%
\definecolor{currentstroke}{rgb}{0.000000,0.000000,0.000000}%
\pgfsetstrokecolor{currentstroke}%
\pgfsetdash{}{0pt}%
\pgfsys@defobject{currentmarker}{\pgfqpoint{0.000000in}{-0.048611in}}{\pgfqpoint{0.000000in}{0.000000in}}{%
\pgfpathmoveto{\pgfqpoint{0.000000in}{0.000000in}}%
\pgfpathlineto{\pgfqpoint{0.000000in}{-0.048611in}}%
\pgfusepath{stroke,fill}%
}%
\begin{pgfscope}%
\pgfsys@transformshift{3.799640in}{0.731763in}%
\pgfsys@useobject{currentmarker}{}%
\end{pgfscope}%
\end{pgfscope}%
\begin{pgfscope}%
\definecolor{textcolor}{rgb}{0.000000,0.000000,0.000000}%
\pgfsetstrokecolor{textcolor}%
\pgfsetfillcolor{textcolor}%
\pgftext[x=3.799640in,y=0.634541in,,top]{\color{textcolor}\rmfamily\fontsize{14.000000}{16.800000}\selectfont 0.75}%
\end{pgfscope}%
\begin{pgfscope}%
\pgfsetbuttcap%
\pgfsetroundjoin%
\definecolor{currentfill}{rgb}{0.000000,0.000000,0.000000}%
\pgfsetfillcolor{currentfill}%
\pgfsetlinewidth{0.803000pt}%
\definecolor{currentstroke}{rgb}{0.000000,0.000000,0.000000}%
\pgfsetstrokecolor{currentstroke}%
\pgfsetdash{}{0pt}%
\pgfsys@defobject{currentmarker}{\pgfqpoint{0.000000in}{-0.048611in}}{\pgfqpoint{0.000000in}{0.000000in}}{%
\pgfpathmoveto{\pgfqpoint{0.000000in}{0.000000in}}%
\pgfpathlineto{\pgfqpoint{0.000000in}{-0.048611in}}%
\pgfusepath{stroke,fill}%
}%
\begin{pgfscope}%
\pgfsys@transformshift{4.776623in}{0.731763in}%
\pgfsys@useobject{currentmarker}{}%
\end{pgfscope}%
\end{pgfscope}%
\begin{pgfscope}%
\definecolor{textcolor}{rgb}{0.000000,0.000000,0.000000}%
\pgfsetstrokecolor{textcolor}%
\pgfsetfillcolor{textcolor}%
\pgftext[x=4.776623in,y=0.634541in,,top]{\color{textcolor}\rmfamily\fontsize{14.000000}{16.800000}\selectfont 1.0}%
\end{pgfscope}%
\begin{pgfscope}%
\definecolor{textcolor}{rgb}{0.000000,0.000000,0.000000}%
\pgfsetstrokecolor{textcolor}%
\pgfsetfillcolor{textcolor}%
\pgftext[x=3.018053in,y=0.401208in,,top]{\color{textcolor}\rmfamily\fontsize{16.000000}{19.200000}\selectfont \(\displaystyle \lambda\)}%
\end{pgfscope}%
\begin{pgfscope}%
\pgfsetbuttcap%
\pgfsetroundjoin%
\definecolor{currentfill}{rgb}{0.000000,0.000000,0.000000}%
\pgfsetfillcolor{currentfill}%
\pgfsetlinewidth{0.803000pt}%
\definecolor{currentstroke}{rgb}{0.000000,0.000000,0.000000}%
\pgfsetstrokecolor{currentstroke}%
\pgfsetdash{}{0pt}%
\pgfsys@defobject{currentmarker}{\pgfqpoint{-0.048611in}{0.000000in}}{\pgfqpoint{-0.000000in}{0.000000in}}{%
\pgfpathmoveto{\pgfqpoint{-0.000000in}{0.000000in}}%
\pgfpathlineto{\pgfqpoint{-0.048611in}{0.000000in}}%
\pgfusepath{stroke,fill}%
}%
\begin{pgfscope}%
\pgfsys@transformshift{0.868690in}{0.731763in}%
\pgfsys@useobject{currentmarker}{}%
\end{pgfscope}%
\end{pgfscope}%
\begin{pgfscope}%
\definecolor{textcolor}{rgb}{0.000000,0.000000,0.000000}%
\pgfsetstrokecolor{textcolor}%
\pgfsetfillcolor{textcolor}%
\pgftext[x=0.673553in, y=0.662319in, left, base]{\color{textcolor}\rmfamily\fontsize{14.000000}{16.800000}\selectfont \(\displaystyle {0}\)}%
\end{pgfscope}%
\begin{pgfscope}%
\pgfsetbuttcap%
\pgfsetroundjoin%
\definecolor{currentfill}{rgb}{0.000000,0.000000,0.000000}%
\pgfsetfillcolor{currentfill}%
\pgfsetlinewidth{0.803000pt}%
\definecolor{currentstroke}{rgb}{0.000000,0.000000,0.000000}%
\pgfsetstrokecolor{currentstroke}%
\pgfsetdash{}{0pt}%
\pgfsys@defobject{currentmarker}{\pgfqpoint{-0.048611in}{0.000000in}}{\pgfqpoint{-0.000000in}{0.000000in}}{%
\pgfpathmoveto{\pgfqpoint{-0.000000in}{0.000000in}}%
\pgfpathlineto{\pgfqpoint{-0.048611in}{0.000000in}}%
\pgfusepath{stroke,fill}%
}%
\begin{pgfscope}%
\pgfsys@transformshift{0.868690in}{1.230259in}%
\pgfsys@useobject{currentmarker}{}%
\end{pgfscope}%
\end{pgfscope}%
\begin{pgfscope}%
\definecolor{textcolor}{rgb}{0.000000,0.000000,0.000000}%
\pgfsetstrokecolor{textcolor}%
\pgfsetfillcolor{textcolor}%
\pgftext[x=0.477722in, y=1.160815in, left, base]{\color{textcolor}\rmfamily\fontsize{14.000000}{16.800000}\selectfont \(\displaystyle {100}\)}%
\end{pgfscope}%
\begin{pgfscope}%
\pgfsetbuttcap%
\pgfsetroundjoin%
\definecolor{currentfill}{rgb}{0.000000,0.000000,0.000000}%
\pgfsetfillcolor{currentfill}%
\pgfsetlinewidth{0.803000pt}%
\definecolor{currentstroke}{rgb}{0.000000,0.000000,0.000000}%
\pgfsetstrokecolor{currentstroke}%
\pgfsetdash{}{0pt}%
\pgfsys@defobject{currentmarker}{\pgfqpoint{-0.048611in}{0.000000in}}{\pgfqpoint{-0.000000in}{0.000000in}}{%
\pgfpathmoveto{\pgfqpoint{-0.000000in}{0.000000in}}%
\pgfpathlineto{\pgfqpoint{-0.048611in}{0.000000in}}%
\pgfusepath{stroke,fill}%
}%
\begin{pgfscope}%
\pgfsys@transformshift{0.868690in}{1.728755in}%
\pgfsys@useobject{currentmarker}{}%
\end{pgfscope}%
\end{pgfscope}%
\begin{pgfscope}%
\definecolor{textcolor}{rgb}{0.000000,0.000000,0.000000}%
\pgfsetstrokecolor{textcolor}%
\pgfsetfillcolor{textcolor}%
\pgftext[x=0.477722in, y=1.659310in, left, base]{\color{textcolor}\rmfamily\fontsize{14.000000}{16.800000}\selectfont \(\displaystyle {200}\)}%
\end{pgfscope}%
\begin{pgfscope}%
\pgfsetbuttcap%
\pgfsetroundjoin%
\definecolor{currentfill}{rgb}{0.000000,0.000000,0.000000}%
\pgfsetfillcolor{currentfill}%
\pgfsetlinewidth{0.803000pt}%
\definecolor{currentstroke}{rgb}{0.000000,0.000000,0.000000}%
\pgfsetstrokecolor{currentstroke}%
\pgfsetdash{}{0pt}%
\pgfsys@defobject{currentmarker}{\pgfqpoint{-0.048611in}{0.000000in}}{\pgfqpoint{-0.000000in}{0.000000in}}{%
\pgfpathmoveto{\pgfqpoint{-0.000000in}{0.000000in}}%
\pgfpathlineto{\pgfqpoint{-0.048611in}{0.000000in}}%
\pgfusepath{stroke,fill}%
}%
\begin{pgfscope}%
\pgfsys@transformshift{0.868690in}{2.227250in}%
\pgfsys@useobject{currentmarker}{}%
\end{pgfscope}%
\end{pgfscope}%
\begin{pgfscope}%
\definecolor{textcolor}{rgb}{0.000000,0.000000,0.000000}%
\pgfsetstrokecolor{textcolor}%
\pgfsetfillcolor{textcolor}%
\pgftext[x=0.477722in, y=2.157806in, left, base]{\color{textcolor}\rmfamily\fontsize{14.000000}{16.800000}\selectfont \(\displaystyle {300}\)}%
\end{pgfscope}%
\begin{pgfscope}%
\pgfsetbuttcap%
\pgfsetroundjoin%
\definecolor{currentfill}{rgb}{0.000000,0.000000,0.000000}%
\pgfsetfillcolor{currentfill}%
\pgfsetlinewidth{0.803000pt}%
\definecolor{currentstroke}{rgb}{0.000000,0.000000,0.000000}%
\pgfsetstrokecolor{currentstroke}%
\pgfsetdash{}{0pt}%
\pgfsys@defobject{currentmarker}{\pgfqpoint{-0.048611in}{0.000000in}}{\pgfqpoint{-0.000000in}{0.000000in}}{%
\pgfpathmoveto{\pgfqpoint{-0.000000in}{0.000000in}}%
\pgfpathlineto{\pgfqpoint{-0.048611in}{0.000000in}}%
\pgfusepath{stroke,fill}%
}%
\begin{pgfscope}%
\pgfsys@transformshift{0.868690in}{2.725746in}%
\pgfsys@useobject{currentmarker}{}%
\end{pgfscope}%
\end{pgfscope}%
\begin{pgfscope}%
\definecolor{textcolor}{rgb}{0.000000,0.000000,0.000000}%
\pgfsetstrokecolor{textcolor}%
\pgfsetfillcolor{textcolor}%
\pgftext[x=0.477722in, y=2.656301in, left, base]{\color{textcolor}\rmfamily\fontsize{14.000000}{16.800000}\selectfont \(\displaystyle {400}\)}%
\end{pgfscope}%
\begin{pgfscope}%
\pgfsetbuttcap%
\pgfsetroundjoin%
\definecolor{currentfill}{rgb}{0.000000,0.000000,0.000000}%
\pgfsetfillcolor{currentfill}%
\pgfsetlinewidth{0.803000pt}%
\definecolor{currentstroke}{rgb}{0.000000,0.000000,0.000000}%
\pgfsetstrokecolor{currentstroke}%
\pgfsetdash{}{0pt}%
\pgfsys@defobject{currentmarker}{\pgfqpoint{-0.048611in}{0.000000in}}{\pgfqpoint{-0.000000in}{0.000000in}}{%
\pgfpathmoveto{\pgfqpoint{-0.000000in}{0.000000in}}%
\pgfpathlineto{\pgfqpoint{-0.048611in}{0.000000in}}%
\pgfusepath{stroke,fill}%
}%
\begin{pgfscope}%
\pgfsys@transformshift{0.868690in}{3.224241in}%
\pgfsys@useobject{currentmarker}{}%
\end{pgfscope}%
\end{pgfscope}%
\begin{pgfscope}%
\definecolor{textcolor}{rgb}{0.000000,0.000000,0.000000}%
\pgfsetstrokecolor{textcolor}%
\pgfsetfillcolor{textcolor}%
\pgftext[x=0.477722in, y=3.154797in, left, base]{\color{textcolor}\rmfamily\fontsize{14.000000}{16.800000}\selectfont \(\displaystyle {500}\)}%
\end{pgfscope}%
\begin{pgfscope}%
\pgfsetbuttcap%
\pgfsetroundjoin%
\definecolor{currentfill}{rgb}{0.000000,0.000000,0.000000}%
\pgfsetfillcolor{currentfill}%
\pgfsetlinewidth{0.803000pt}%
\definecolor{currentstroke}{rgb}{0.000000,0.000000,0.000000}%
\pgfsetstrokecolor{currentstroke}%
\pgfsetdash{}{0pt}%
\pgfsys@defobject{currentmarker}{\pgfqpoint{-0.048611in}{0.000000in}}{\pgfqpoint{-0.000000in}{0.000000in}}{%
\pgfpathmoveto{\pgfqpoint{-0.000000in}{0.000000in}}%
\pgfpathlineto{\pgfqpoint{-0.048611in}{0.000000in}}%
\pgfusepath{stroke,fill}%
}%
\begin{pgfscope}%
\pgfsys@transformshift{0.868690in}{3.722737in}%
\pgfsys@useobject{currentmarker}{}%
\end{pgfscope}%
\end{pgfscope}%
\begin{pgfscope}%
\definecolor{textcolor}{rgb}{0.000000,0.000000,0.000000}%
\pgfsetstrokecolor{textcolor}%
\pgfsetfillcolor{textcolor}%
\pgftext[x=0.477722in, y=3.653293in, left, base]{\color{textcolor}\rmfamily\fontsize{14.000000}{16.800000}\selectfont \(\displaystyle {600}\)}%
\end{pgfscope}%
\begin{pgfscope}%
\definecolor{textcolor}{rgb}{0.000000,0.000000,0.000000}%
\pgfsetstrokecolor{textcolor}%
\pgfsetfillcolor{textcolor}%
\pgftext[x=0.422166in,y=2.227250in,,bottom,rotate=90.000000]{\color{textcolor}\rmfamily\fontsize{16.000000}{19.200000}\selectfont Average time (minutes)}%
\end{pgfscope}%
\begin{pgfscope}%
\pgfpathrectangle{\pgfqpoint{0.868690in}{0.731763in}}{\pgfqpoint{4.298726in}{2.990974in}}%
\pgfusepath{clip}%
\pgfsetrectcap%
\pgfsetroundjoin%
\pgfsetlinewidth{1.505625pt}%
\definecolor{currentstroke}{rgb}{1.000000,0.498039,0.054902}%
\pgfsetstrokecolor{currentstroke}%
\pgfsetdash{}{0pt}%
\pgfpathmoveto{\pgfqpoint{1.259483in}{1.799257in}}%
\pgfpathlineto{\pgfqpoint{1.845673in}{2.149790in}}%
\pgfpathlineto{\pgfqpoint{2.822656in}{2.604983in}}%
\pgfpathlineto{\pgfqpoint{3.799640in}{2.843758in}}%
\pgfpathlineto{\pgfqpoint{4.776623in}{2.854309in}}%
\pgfusepath{stroke}%
\end{pgfscope}%
\begin{pgfscope}%
\pgfpathrectangle{\pgfqpoint{0.868690in}{0.731763in}}{\pgfqpoint{4.298726in}{2.990974in}}%
\pgfusepath{clip}%
\pgfsetbuttcap%
\pgfsetroundjoin%
\definecolor{currentfill}{rgb}{1.000000,0.498039,0.054902}%
\pgfsetfillcolor{currentfill}%
\pgfsetlinewidth{1.003750pt}%
\definecolor{currentstroke}{rgb}{1.000000,0.498039,0.054902}%
\pgfsetstrokecolor{currentstroke}%
\pgfsetdash{}{0pt}%
\pgfsys@defobject{currentmarker}{\pgfqpoint{-0.041667in}{-0.041667in}}{\pgfqpoint{0.041667in}{0.041667in}}{%
\pgfpathmoveto{\pgfqpoint{0.000000in}{-0.041667in}}%
\pgfpathcurveto{\pgfqpoint{0.011050in}{-0.041667in}}{\pgfqpoint{0.021649in}{-0.037276in}}{\pgfqpoint{0.029463in}{-0.029463in}}%
\pgfpathcurveto{\pgfqpoint{0.037276in}{-0.021649in}}{\pgfqpoint{0.041667in}{-0.011050in}}{\pgfqpoint{0.041667in}{0.000000in}}%
\pgfpathcurveto{\pgfqpoint{0.041667in}{0.011050in}}{\pgfqpoint{0.037276in}{0.021649in}}{\pgfqpoint{0.029463in}{0.029463in}}%
\pgfpathcurveto{\pgfqpoint{0.021649in}{0.037276in}}{\pgfqpoint{0.011050in}{0.041667in}}{\pgfqpoint{0.000000in}{0.041667in}}%
\pgfpathcurveto{\pgfqpoint{-0.011050in}{0.041667in}}{\pgfqpoint{-0.021649in}{0.037276in}}{\pgfqpoint{-0.029463in}{0.029463in}}%
\pgfpathcurveto{\pgfqpoint{-0.037276in}{0.021649in}}{\pgfqpoint{-0.041667in}{0.011050in}}{\pgfqpoint{-0.041667in}{0.000000in}}%
\pgfpathcurveto{\pgfqpoint{-0.041667in}{-0.011050in}}{\pgfqpoint{-0.037276in}{-0.021649in}}{\pgfqpoint{-0.029463in}{-0.029463in}}%
\pgfpathcurveto{\pgfqpoint{-0.021649in}{-0.037276in}}{\pgfqpoint{-0.011050in}{-0.041667in}}{\pgfqpoint{0.000000in}{-0.041667in}}%
\pgfpathlineto{\pgfqpoint{0.000000in}{-0.041667in}}%
\pgfpathclose%
\pgfusepath{stroke,fill}%
}%
\begin{pgfscope}%
\pgfsys@transformshift{1.259483in}{1.799257in}%
\pgfsys@useobject{currentmarker}{}%
\end{pgfscope}%
\begin{pgfscope}%
\pgfsys@transformshift{1.845673in}{2.149790in}%
\pgfsys@useobject{currentmarker}{}%
\end{pgfscope}%
\begin{pgfscope}%
\pgfsys@transformshift{2.822656in}{2.604983in}%
\pgfsys@useobject{currentmarker}{}%
\end{pgfscope}%
\begin{pgfscope}%
\pgfsys@transformshift{3.799640in}{2.843758in}%
\pgfsys@useobject{currentmarker}{}%
\end{pgfscope}%
\begin{pgfscope}%
\pgfsys@transformshift{4.776623in}{2.854309in}%
\pgfsys@useobject{currentmarker}{}%
\end{pgfscope}%
\end{pgfscope}%
\begin{pgfscope}%
\pgfpathrectangle{\pgfqpoint{0.868690in}{0.731763in}}{\pgfqpoint{4.298726in}{2.990974in}}%
\pgfusepath{clip}%
\pgfsetrectcap%
\pgfsetroundjoin%
\pgfsetlinewidth{1.505625pt}%
\definecolor{currentstroke}{rgb}{0.121569,0.466667,0.705882}%
\pgfsetstrokecolor{currentstroke}%
\pgfsetdash{}{0pt}%
\pgfpathmoveto{\pgfqpoint{1.259483in}{2.734449in}}%
\pgfpathlineto{\pgfqpoint{1.845673in}{2.759523in}}%
\pgfpathlineto{\pgfqpoint{2.822656in}{2.764075in}}%
\pgfpathlineto{\pgfqpoint{3.799640in}{2.803807in}}%
\pgfpathlineto{\pgfqpoint{4.776623in}{2.700992in}}%
\pgfusepath{stroke}%
\end{pgfscope}%
\begin{pgfscope}%
\pgfpathrectangle{\pgfqpoint{0.868690in}{0.731763in}}{\pgfqpoint{4.298726in}{2.990974in}}%
\pgfusepath{clip}%
\pgfsetbuttcap%
\pgfsetroundjoin%
\definecolor{currentfill}{rgb}{0.121569,0.466667,0.705882}%
\pgfsetfillcolor{currentfill}%
\pgfsetlinewidth{1.003750pt}%
\definecolor{currentstroke}{rgb}{0.121569,0.466667,0.705882}%
\pgfsetstrokecolor{currentstroke}%
\pgfsetdash{}{0pt}%
\pgfsys@defobject{currentmarker}{\pgfqpoint{-0.041667in}{-0.041667in}}{\pgfqpoint{0.041667in}{0.041667in}}{%
\pgfpathmoveto{\pgfqpoint{0.000000in}{-0.041667in}}%
\pgfpathcurveto{\pgfqpoint{0.011050in}{-0.041667in}}{\pgfqpoint{0.021649in}{-0.037276in}}{\pgfqpoint{0.029463in}{-0.029463in}}%
\pgfpathcurveto{\pgfqpoint{0.037276in}{-0.021649in}}{\pgfqpoint{0.041667in}{-0.011050in}}{\pgfqpoint{0.041667in}{0.000000in}}%
\pgfpathcurveto{\pgfqpoint{0.041667in}{0.011050in}}{\pgfqpoint{0.037276in}{0.021649in}}{\pgfqpoint{0.029463in}{0.029463in}}%
\pgfpathcurveto{\pgfqpoint{0.021649in}{0.037276in}}{\pgfqpoint{0.011050in}{0.041667in}}{\pgfqpoint{0.000000in}{0.041667in}}%
\pgfpathcurveto{\pgfqpoint{-0.011050in}{0.041667in}}{\pgfqpoint{-0.021649in}{0.037276in}}{\pgfqpoint{-0.029463in}{0.029463in}}%
\pgfpathcurveto{\pgfqpoint{-0.037276in}{0.021649in}}{\pgfqpoint{-0.041667in}{0.011050in}}{\pgfqpoint{-0.041667in}{0.000000in}}%
\pgfpathcurveto{\pgfqpoint{-0.041667in}{-0.011050in}}{\pgfqpoint{-0.037276in}{-0.021649in}}{\pgfqpoint{-0.029463in}{-0.029463in}}%
\pgfpathcurveto{\pgfqpoint{-0.021649in}{-0.037276in}}{\pgfqpoint{-0.011050in}{-0.041667in}}{\pgfqpoint{0.000000in}{-0.041667in}}%
\pgfpathlineto{\pgfqpoint{0.000000in}{-0.041667in}}%
\pgfpathclose%
\pgfusepath{stroke,fill}%
}%
\begin{pgfscope}%
\pgfsys@transformshift{1.259483in}{2.734449in}%
\pgfsys@useobject{currentmarker}{}%
\end{pgfscope}%
\begin{pgfscope}%
\pgfsys@transformshift{1.845673in}{2.759523in}%
\pgfsys@useobject{currentmarker}{}%
\end{pgfscope}%
\begin{pgfscope}%
\pgfsys@transformshift{2.822656in}{2.764075in}%
\pgfsys@useobject{currentmarker}{}%
\end{pgfscope}%
\begin{pgfscope}%
\pgfsys@transformshift{3.799640in}{2.803807in}%
\pgfsys@useobject{currentmarker}{}%
\end{pgfscope}%
\begin{pgfscope}%
\pgfsys@transformshift{4.776623in}{2.700992in}%
\pgfsys@useobject{currentmarker}{}%
\end{pgfscope}%
\end{pgfscope}%
\begin{pgfscope}%
\pgfsetrectcap%
\pgfsetmiterjoin%
\pgfsetlinewidth{0.803000pt}%
\definecolor{currentstroke}{rgb}{0.000000,0.000000,0.000000}%
\pgfsetstrokecolor{currentstroke}%
\pgfsetdash{}{0pt}%
\pgfpathmoveto{\pgfqpoint{0.868690in}{0.731763in}}%
\pgfpathlineto{\pgfqpoint{0.868690in}{3.722737in}}%
\pgfusepath{stroke}%
\end{pgfscope}%
\begin{pgfscope}%
\pgfsetrectcap%
\pgfsetmiterjoin%
\pgfsetlinewidth{0.803000pt}%
\definecolor{currentstroke}{rgb}{0.000000,0.000000,0.000000}%
\pgfsetstrokecolor{currentstroke}%
\pgfsetdash{}{0pt}%
\pgfpathmoveto{\pgfqpoint{5.167416in}{0.731763in}}%
\pgfpathlineto{\pgfqpoint{5.167416in}{3.722737in}}%
\pgfusepath{stroke}%
\end{pgfscope}%
\begin{pgfscope}%
\pgfsetrectcap%
\pgfsetmiterjoin%
\pgfsetlinewidth{0.803000pt}%
\definecolor{currentstroke}{rgb}{0.000000,0.000000,0.000000}%
\pgfsetstrokecolor{currentstroke}%
\pgfsetdash{}{0pt}%
\pgfpathmoveto{\pgfqpoint{0.868690in}{0.731763in}}%
\pgfpathlineto{\pgfqpoint{5.167416in}{0.731763in}}%
\pgfusepath{stroke}%
\end{pgfscope}%
\begin{pgfscope}%
\pgfsetrectcap%
\pgfsetmiterjoin%
\pgfsetlinewidth{0.803000pt}%
\definecolor{currentstroke}{rgb}{0.000000,0.000000,0.000000}%
\pgfsetstrokecolor{currentstroke}%
\pgfsetdash{}{0pt}%
\pgfpathmoveto{\pgfqpoint{0.868690in}{3.722737in}}%
\pgfpathlineto{\pgfqpoint{5.167416in}{3.722737in}}%
\pgfusepath{stroke}%
\end{pgfscope}%
\begin{pgfscope}%
\pgfsetbuttcap%
\pgfsetmiterjoin%
\definecolor{currentfill}{rgb}{1.000000,1.000000,1.000000}%
\pgfsetfillcolor{currentfill}%
\pgfsetfillopacity{0.800000}%
\pgfsetlinewidth{1.003750pt}%
\definecolor{currentstroke}{rgb}{0.800000,0.800000,0.800000}%
\pgfsetstrokecolor{currentstroke}%
\pgfsetstrokeopacity{0.800000}%
\pgfsetdash{}{0pt}%
\pgfpathmoveto{\pgfqpoint{1.014523in}{2.947738in}}%
\pgfpathlineto{\pgfqpoint{4.483577in}{2.947738in}}%
\pgfpathquadraticcurveto{\pgfqpoint{4.525244in}{2.947738in}}{\pgfqpoint{4.525244in}{2.989404in}}%
\pgfpathlineto{\pgfqpoint{4.525244in}{3.576904in}}%
\pgfpathquadraticcurveto{\pgfqpoint{4.525244in}{3.618570in}}{\pgfqpoint{4.483577in}{3.618570in}}%
\pgfpathlineto{\pgfqpoint{1.014523in}{3.618570in}}%
\pgfpathquadraticcurveto{\pgfqpoint{0.972857in}{3.618570in}}{\pgfqpoint{0.972857in}{3.576904in}}%
\pgfpathlineto{\pgfqpoint{0.972857in}{2.989404in}}%
\pgfpathquadraticcurveto{\pgfqpoint{0.972857in}{2.947738in}}{\pgfqpoint{1.014523in}{2.947738in}}%
\pgfpathlineto{\pgfqpoint{1.014523in}{2.947738in}}%
\pgfpathclose%
\pgfusepath{stroke,fill}%
\end{pgfscope}%
\begin{pgfscope}%
\pgfsetrectcap%
\pgfsetroundjoin%
\pgfsetlinewidth{1.505625pt}%
\definecolor{currentstroke}{rgb}{1.000000,0.498039,0.054902}%
\pgfsetstrokecolor{currentstroke}%
\pgfsetdash{}{0pt}%
\pgfpathmoveto{\pgfqpoint{1.056190in}{3.458154in}}%
\pgfpathlineto{\pgfqpoint{1.264523in}{3.458154in}}%
\pgfpathlineto{\pgfqpoint{1.472857in}{3.458154in}}%
\pgfusepath{stroke}%
\end{pgfscope}%
\begin{pgfscope}%
\pgfsetbuttcap%
\pgfsetroundjoin%
\definecolor{currentfill}{rgb}{1.000000,0.498039,0.054902}%
\pgfsetfillcolor{currentfill}%
\pgfsetlinewidth{1.003750pt}%
\definecolor{currentstroke}{rgb}{1.000000,0.498039,0.054902}%
\pgfsetstrokecolor{currentstroke}%
\pgfsetdash{}{0pt}%
\pgfsys@defobject{currentmarker}{\pgfqpoint{-0.041667in}{-0.041667in}}{\pgfqpoint{0.041667in}{0.041667in}}{%
\pgfpathmoveto{\pgfqpoint{0.000000in}{-0.041667in}}%
\pgfpathcurveto{\pgfqpoint{0.011050in}{-0.041667in}}{\pgfqpoint{0.021649in}{-0.037276in}}{\pgfqpoint{0.029463in}{-0.029463in}}%
\pgfpathcurveto{\pgfqpoint{0.037276in}{-0.021649in}}{\pgfqpoint{0.041667in}{-0.011050in}}{\pgfqpoint{0.041667in}{0.000000in}}%
\pgfpathcurveto{\pgfqpoint{0.041667in}{0.011050in}}{\pgfqpoint{0.037276in}{0.021649in}}{\pgfqpoint{0.029463in}{0.029463in}}%
\pgfpathcurveto{\pgfqpoint{0.021649in}{0.037276in}}{\pgfqpoint{0.011050in}{0.041667in}}{\pgfqpoint{0.000000in}{0.041667in}}%
\pgfpathcurveto{\pgfqpoint{-0.011050in}{0.041667in}}{\pgfqpoint{-0.021649in}{0.037276in}}{\pgfqpoint{-0.029463in}{0.029463in}}%
\pgfpathcurveto{\pgfqpoint{-0.037276in}{0.021649in}}{\pgfqpoint{-0.041667in}{0.011050in}}{\pgfqpoint{-0.041667in}{0.000000in}}%
\pgfpathcurveto{\pgfqpoint{-0.041667in}{-0.011050in}}{\pgfqpoint{-0.037276in}{-0.021649in}}{\pgfqpoint{-0.029463in}{-0.029463in}}%
\pgfpathcurveto{\pgfqpoint{-0.021649in}{-0.037276in}}{\pgfqpoint{-0.011050in}{-0.041667in}}{\pgfqpoint{0.000000in}{-0.041667in}}%
\pgfpathlineto{\pgfqpoint{0.000000in}{-0.041667in}}%
\pgfpathclose%
\pgfusepath{stroke,fill}%
}%
\begin{pgfscope}%
\pgfsys@transformshift{1.264523in}{3.458154in}%
\pgfsys@useobject{currentmarker}{}%
\end{pgfscope}%
\end{pgfscope}%
\begin{pgfscope}%
\definecolor{textcolor}{rgb}{0.000000,0.000000,0.000000}%
\pgfsetstrokecolor{textcolor}%
\pgfsetfillcolor{textcolor}%
\pgftext[x=1.639523in,y=3.385237in,left,base]{\color{textcolor}\rmfamily\fontsize{15.000000}{18.000000}\selectfont \(\displaystyle A\) (Time)}%
\end{pgfscope}%
\begin{pgfscope}%
\pgfsetrectcap%
\pgfsetroundjoin%
\pgfsetlinewidth{1.505625pt}%
\definecolor{currentstroke}{rgb}{0.121569,0.466667,0.705882}%
\pgfsetstrokecolor{currentstroke}%
\pgfsetdash{}{0pt}%
\pgfpathmoveto{\pgfqpoint{1.056190in}{3.153987in}}%
\pgfpathlineto{\pgfqpoint{1.264523in}{3.153987in}}%
\pgfpathlineto{\pgfqpoint{1.472857in}{3.153987in}}%
\pgfusepath{stroke}%
\end{pgfscope}%
\begin{pgfscope}%
\pgfsetbuttcap%
\pgfsetroundjoin%
\definecolor{currentfill}{rgb}{0.121569,0.466667,0.705882}%
\pgfsetfillcolor{currentfill}%
\pgfsetlinewidth{1.003750pt}%
\definecolor{currentstroke}{rgb}{0.121569,0.466667,0.705882}%
\pgfsetstrokecolor{currentstroke}%
\pgfsetdash{}{0pt}%
\pgfsys@defobject{currentmarker}{\pgfqpoint{-0.041667in}{-0.041667in}}{\pgfqpoint{0.041667in}{0.041667in}}{%
\pgfpathmoveto{\pgfqpoint{0.000000in}{-0.041667in}}%
\pgfpathcurveto{\pgfqpoint{0.011050in}{-0.041667in}}{\pgfqpoint{0.021649in}{-0.037276in}}{\pgfqpoint{0.029463in}{-0.029463in}}%
\pgfpathcurveto{\pgfqpoint{0.037276in}{-0.021649in}}{\pgfqpoint{0.041667in}{-0.011050in}}{\pgfqpoint{0.041667in}{0.000000in}}%
\pgfpathcurveto{\pgfqpoint{0.041667in}{0.011050in}}{\pgfqpoint{0.037276in}{0.021649in}}{\pgfqpoint{0.029463in}{0.029463in}}%
\pgfpathcurveto{\pgfqpoint{0.021649in}{0.037276in}}{\pgfqpoint{0.011050in}{0.041667in}}{\pgfqpoint{0.000000in}{0.041667in}}%
\pgfpathcurveto{\pgfqpoint{-0.011050in}{0.041667in}}{\pgfqpoint{-0.021649in}{0.037276in}}{\pgfqpoint{-0.029463in}{0.029463in}}%
\pgfpathcurveto{\pgfqpoint{-0.037276in}{0.021649in}}{\pgfqpoint{-0.041667in}{0.011050in}}{\pgfqpoint{-0.041667in}{0.000000in}}%
\pgfpathcurveto{\pgfqpoint{-0.041667in}{-0.011050in}}{\pgfqpoint{-0.037276in}{-0.021649in}}{\pgfqpoint{-0.029463in}{-0.029463in}}%
\pgfpathcurveto{\pgfqpoint{-0.021649in}{-0.037276in}}{\pgfqpoint{-0.011050in}{-0.041667in}}{\pgfqpoint{0.000000in}{-0.041667in}}%
\pgfpathlineto{\pgfqpoint{0.000000in}{-0.041667in}}%
\pgfpathclose%
\pgfusepath{stroke,fill}%
}%
\begin{pgfscope}%
\pgfsys@transformshift{1.264523in}{3.153987in}%
\pgfsys@useobject{currentmarker}{}%
\end{pgfscope}%
\end{pgfscope}%
\begin{pgfscope}%
\definecolor{textcolor}{rgb}{0.000000,0.000000,0.000000}%
\pgfsetstrokecolor{textcolor}%
\pgfsetfillcolor{textcolor}%
\pgftext[x=1.639523in,y=3.081071in,left,base]{\color{textcolor}\rmfamily\fontsize{15.000000}{18.000000}\selectfont \(\displaystyle M\) (Time)}%
\end{pgfscope}%
\begin{pgfscope}%
\pgfsetbuttcap%
\pgfsetmiterjoin%
\definecolor{currentfill}{rgb}{1.000000,0.498039,0.054902}%
\pgfsetfillcolor{currentfill}%
\pgfsetfillopacity{0.400000}%
\pgfsetlinewidth{0.000000pt}%
\definecolor{currentstroke}{rgb}{0.000000,0.000000,0.000000}%
\pgfsetstrokecolor{currentstroke}%
\pgfsetstrokeopacity{0.400000}%
\pgfsetdash{}{0pt}%
\pgfpathmoveto{\pgfqpoint{2.930185in}{3.385237in}}%
\pgfpathlineto{\pgfqpoint{3.346851in}{3.385237in}}%
\pgfpathlineto{\pgfqpoint{3.346851in}{3.531070in}}%
\pgfpathlineto{\pgfqpoint{2.930185in}{3.531070in}}%
\pgfpathlineto{\pgfqpoint{2.930185in}{3.385237in}}%
\pgfpathclose%
\pgfusepath{fill}%
\end{pgfscope}%
\begin{pgfscope}%
\definecolor{textcolor}{rgb}{0.000000,0.000000,0.000000}%
\pgfsetstrokecolor{textcolor}%
\pgfsetfillcolor{textcolor}%
\pgftext[x=3.513518in,y=3.385237in,left,base]{\color{textcolor}\rmfamily\fontsize{15.000000}{18.000000}\selectfont \(\displaystyle A\) (\(\displaystyle \#\)Opt)}%
\end{pgfscope}%
\begin{pgfscope}%
\pgfsetbuttcap%
\pgfsetmiterjoin%
\definecolor{currentfill}{rgb}{0.121569,0.466667,0.705882}%
\pgfsetfillcolor{currentfill}%
\pgfsetfillopacity{0.400000}%
\pgfsetlinewidth{0.000000pt}%
\definecolor{currentstroke}{rgb}{0.000000,0.000000,0.000000}%
\pgfsetstrokecolor{currentstroke}%
\pgfsetstrokeopacity{0.400000}%
\pgfsetdash{}{0pt}%
\pgfpathmoveto{\pgfqpoint{2.930185in}{3.081071in}}%
\pgfpathlineto{\pgfqpoint{3.346851in}{3.081071in}}%
\pgfpathlineto{\pgfqpoint{3.346851in}{3.226904in}}%
\pgfpathlineto{\pgfqpoint{2.930185in}{3.226904in}}%
\pgfpathlineto{\pgfqpoint{2.930185in}{3.081071in}}%
\pgfpathclose%
\pgfusepath{fill}%
\end{pgfscope}%
\begin{pgfscope}%
\definecolor{textcolor}{rgb}{0.000000,0.000000,0.000000}%
\pgfsetstrokecolor{textcolor}%
\pgfsetfillcolor{textcolor}%
\pgftext[x=3.513518in,y=3.081071in,left,base]{\color{textcolor}\rmfamily\fontsize{15.000000}{18.000000}\selectfont \(\displaystyle M\) (\(\displaystyle \#\)Opt)}%
\end{pgfscope}%
\begin{pgfscope}%
\pgfpathrectangle{\pgfqpoint{0.868690in}{0.731763in}}{\pgfqpoint{4.298726in}{2.990974in}}%
\pgfusepath{clip}%
\pgfsetbuttcap%
\pgfsetmiterjoin%
\definecolor{currentfill}{rgb}{1.000000,0.498039,0.054902}%
\pgfsetfillcolor{currentfill}%
\pgfsetfillopacity{0.400000}%
\pgfsetlinewidth{0.000000pt}%
\definecolor{currentstroke}{rgb}{0.000000,0.000000,0.000000}%
\pgfsetstrokecolor{currentstroke}%
\pgfsetstrokeopacity{0.400000}%
\pgfsetdash{}{0pt}%
\pgfpathmoveto{\pgfqpoint{1.103166in}{0.731763in}}%
\pgfpathlineto{\pgfqpoint{1.259483in}{0.731763in}}%
\pgfpathlineto{\pgfqpoint{1.259483in}{2.778219in}}%
\pgfpathlineto{\pgfqpoint{1.103166in}{2.778219in}}%
\pgfpathlineto{\pgfqpoint{1.103166in}{0.731763in}}%
\pgfpathclose%
\pgfusepath{fill}%
\end{pgfscope}%
\begin{pgfscope}%
\pgfpathrectangle{\pgfqpoint{0.868690in}{0.731763in}}{\pgfqpoint{4.298726in}{2.990974in}}%
\pgfusepath{clip}%
\pgfsetbuttcap%
\pgfsetmiterjoin%
\definecolor{currentfill}{rgb}{1.000000,0.498039,0.054902}%
\pgfsetfillcolor{currentfill}%
\pgfsetfillopacity{0.400000}%
\pgfsetlinewidth{0.000000pt}%
\definecolor{currentstroke}{rgb}{0.000000,0.000000,0.000000}%
\pgfsetstrokecolor{currentstroke}%
\pgfsetstrokeopacity{0.400000}%
\pgfsetdash{}{0pt}%
\pgfpathmoveto{\pgfqpoint{1.689356in}{0.731763in}}%
\pgfpathlineto{\pgfqpoint{1.845673in}{0.731763in}}%
\pgfpathlineto{\pgfqpoint{1.845673in}{2.542090in}}%
\pgfpathlineto{\pgfqpoint{1.689356in}{2.542090in}}%
\pgfpathlineto{\pgfqpoint{1.689356in}{0.731763in}}%
\pgfpathclose%
\pgfusepath{fill}%
\end{pgfscope}%
\begin{pgfscope}%
\pgfpathrectangle{\pgfqpoint{0.868690in}{0.731763in}}{\pgfqpoint{4.298726in}{2.990974in}}%
\pgfusepath{clip}%
\pgfsetbuttcap%
\pgfsetmiterjoin%
\definecolor{currentfill}{rgb}{1.000000,0.498039,0.054902}%
\pgfsetfillcolor{currentfill}%
\pgfsetfillopacity{0.400000}%
\pgfsetlinewidth{0.000000pt}%
\definecolor{currentstroke}{rgb}{0.000000,0.000000,0.000000}%
\pgfsetstrokecolor{currentstroke}%
\pgfsetstrokeopacity{0.400000}%
\pgfsetdash{}{0pt}%
\pgfpathmoveto{\pgfqpoint{2.666339in}{0.731763in}}%
\pgfpathlineto{\pgfqpoint{2.822656in}{0.731763in}}%
\pgfpathlineto{\pgfqpoint{2.822656in}{2.227250in}}%
\pgfpathlineto{\pgfqpoint{2.666339in}{2.227250in}}%
\pgfpathlineto{\pgfqpoint{2.666339in}{0.731763in}}%
\pgfpathclose%
\pgfusepath{fill}%
\end{pgfscope}%
\begin{pgfscope}%
\pgfpathrectangle{\pgfqpoint{0.868690in}{0.731763in}}{\pgfqpoint{4.298726in}{2.990974in}}%
\pgfusepath{clip}%
\pgfsetbuttcap%
\pgfsetmiterjoin%
\definecolor{currentfill}{rgb}{1.000000,0.498039,0.054902}%
\pgfsetfillcolor{currentfill}%
\pgfsetfillopacity{0.400000}%
\pgfsetlinewidth{0.000000pt}%
\definecolor{currentstroke}{rgb}{0.000000,0.000000,0.000000}%
\pgfsetstrokecolor{currentstroke}%
\pgfsetstrokeopacity{0.400000}%
\pgfsetdash{}{0pt}%
\pgfpathmoveto{\pgfqpoint{3.643322in}{0.731763in}}%
\pgfpathlineto{\pgfqpoint{3.799640in}{0.731763in}}%
\pgfpathlineto{\pgfqpoint{3.799640in}{1.912411in}}%
\pgfpathlineto{\pgfqpoint{3.643322in}{1.912411in}}%
\pgfpathlineto{\pgfqpoint{3.643322in}{0.731763in}}%
\pgfpathclose%
\pgfusepath{fill}%
\end{pgfscope}%
\begin{pgfscope}%
\pgfpathrectangle{\pgfqpoint{0.868690in}{0.731763in}}{\pgfqpoint{4.298726in}{2.990974in}}%
\pgfusepath{clip}%
\pgfsetbuttcap%
\pgfsetmiterjoin%
\definecolor{currentfill}{rgb}{1.000000,0.498039,0.054902}%
\pgfsetfillcolor{currentfill}%
\pgfsetfillopacity{0.400000}%
\pgfsetlinewidth{0.000000pt}%
\definecolor{currentstroke}{rgb}{0.000000,0.000000,0.000000}%
\pgfsetstrokecolor{currentstroke}%
\pgfsetstrokeopacity{0.400000}%
\pgfsetdash{}{0pt}%
\pgfpathmoveto{\pgfqpoint{4.620306in}{0.731763in}}%
\pgfpathlineto{\pgfqpoint{4.776623in}{0.731763in}}%
\pgfpathlineto{\pgfqpoint{4.776623in}{1.833701in}}%
\pgfpathlineto{\pgfqpoint{4.620306in}{1.833701in}}%
\pgfpathlineto{\pgfqpoint{4.620306in}{0.731763in}}%
\pgfpathclose%
\pgfusepath{fill}%
\end{pgfscope}%
\begin{pgfscope}%
\pgfpathrectangle{\pgfqpoint{0.868690in}{0.731763in}}{\pgfqpoint{4.298726in}{2.990974in}}%
\pgfusepath{clip}%
\pgfsetbuttcap%
\pgfsetmiterjoin%
\definecolor{currentfill}{rgb}{0.121569,0.466667,0.705882}%
\pgfsetfillcolor{currentfill}%
\pgfsetfillopacity{0.400000}%
\pgfsetlinewidth{0.000000pt}%
\definecolor{currentstroke}{rgb}{0.000000,0.000000,0.000000}%
\pgfsetstrokecolor{currentstroke}%
\pgfsetstrokeopacity{0.400000}%
\pgfsetdash{}{0pt}%
\pgfpathmoveto{\pgfqpoint{1.259483in}{0.731763in}}%
\pgfpathlineto{\pgfqpoint{1.415801in}{0.731763in}}%
\pgfpathlineto{\pgfqpoint{1.415801in}{2.069831in}}%
\pgfpathlineto{\pgfqpoint{1.259483in}{2.069831in}}%
\pgfpathlineto{\pgfqpoint{1.259483in}{0.731763in}}%
\pgfpathclose%
\pgfusepath{fill}%
\end{pgfscope}%
\begin{pgfscope}%
\pgfpathrectangle{\pgfqpoint{0.868690in}{0.731763in}}{\pgfqpoint{4.298726in}{2.990974in}}%
\pgfusepath{clip}%
\pgfsetbuttcap%
\pgfsetmiterjoin%
\definecolor{currentfill}{rgb}{0.121569,0.466667,0.705882}%
\pgfsetfillcolor{currentfill}%
\pgfsetfillopacity{0.400000}%
\pgfsetlinewidth{0.000000pt}%
\definecolor{currentstroke}{rgb}{0.000000,0.000000,0.000000}%
\pgfsetstrokecolor{currentstroke}%
\pgfsetstrokeopacity{0.400000}%
\pgfsetdash{}{0pt}%
\pgfpathmoveto{\pgfqpoint{1.845673in}{0.731763in}}%
\pgfpathlineto{\pgfqpoint{2.001991in}{0.731763in}}%
\pgfpathlineto{\pgfqpoint{2.001991in}{1.991121in}}%
\pgfpathlineto{\pgfqpoint{1.845673in}{1.991121in}}%
\pgfpathlineto{\pgfqpoint{1.845673in}{0.731763in}}%
\pgfpathclose%
\pgfusepath{fill}%
\end{pgfscope}%
\begin{pgfscope}%
\pgfpathrectangle{\pgfqpoint{0.868690in}{0.731763in}}{\pgfqpoint{4.298726in}{2.990974in}}%
\pgfusepath{clip}%
\pgfsetbuttcap%
\pgfsetmiterjoin%
\definecolor{currentfill}{rgb}{0.121569,0.466667,0.705882}%
\pgfsetfillcolor{currentfill}%
\pgfsetfillopacity{0.400000}%
\pgfsetlinewidth{0.000000pt}%
\definecolor{currentstroke}{rgb}{0.000000,0.000000,0.000000}%
\pgfsetstrokecolor{currentstroke}%
\pgfsetstrokeopacity{0.400000}%
\pgfsetdash{}{0pt}%
\pgfpathmoveto{\pgfqpoint{2.822656in}{0.731763in}}%
\pgfpathlineto{\pgfqpoint{2.978974in}{0.731763in}}%
\pgfpathlineto{\pgfqpoint{2.978974in}{1.991121in}}%
\pgfpathlineto{\pgfqpoint{2.822656in}{1.991121in}}%
\pgfpathlineto{\pgfqpoint{2.822656in}{0.731763in}}%
\pgfpathclose%
\pgfusepath{fill}%
\end{pgfscope}%
\begin{pgfscope}%
\pgfpathrectangle{\pgfqpoint{0.868690in}{0.731763in}}{\pgfqpoint{4.298726in}{2.990974in}}%
\pgfusepath{clip}%
\pgfsetbuttcap%
\pgfsetmiterjoin%
\definecolor{currentfill}{rgb}{0.121569,0.466667,0.705882}%
\pgfsetfillcolor{currentfill}%
\pgfsetfillopacity{0.400000}%
\pgfsetlinewidth{0.000000pt}%
\definecolor{currentstroke}{rgb}{0.000000,0.000000,0.000000}%
\pgfsetstrokecolor{currentstroke}%
\pgfsetstrokeopacity{0.400000}%
\pgfsetdash{}{0pt}%
\pgfpathmoveto{\pgfqpoint{3.799640in}{0.731763in}}%
\pgfpathlineto{\pgfqpoint{3.955957in}{0.731763in}}%
\pgfpathlineto{\pgfqpoint{3.955957in}{1.833701in}}%
\pgfpathlineto{\pgfqpoint{3.799640in}{1.833701in}}%
\pgfpathlineto{\pgfqpoint{3.799640in}{0.731763in}}%
\pgfpathclose%
\pgfusepath{fill}%
\end{pgfscope}%
\begin{pgfscope}%
\pgfpathrectangle{\pgfqpoint{0.868690in}{0.731763in}}{\pgfqpoint{4.298726in}{2.990974in}}%
\pgfusepath{clip}%
\pgfsetbuttcap%
\pgfsetmiterjoin%
\definecolor{currentfill}{rgb}{0.121569,0.466667,0.705882}%
\pgfsetfillcolor{currentfill}%
\pgfsetfillopacity{0.400000}%
\pgfsetlinewidth{0.000000pt}%
\definecolor{currentstroke}{rgb}{0.000000,0.000000,0.000000}%
\pgfsetstrokecolor{currentstroke}%
\pgfsetstrokeopacity{0.400000}%
\pgfsetdash{}{0pt}%
\pgfpathmoveto{\pgfqpoint{4.776623in}{0.731763in}}%
\pgfpathlineto{\pgfqpoint{4.932940in}{0.731763in}}%
\pgfpathlineto{\pgfqpoint{4.932940in}{1.912411in}}%
\pgfpathlineto{\pgfqpoint{4.776623in}{1.912411in}}%
\pgfpathlineto{\pgfqpoint{4.776623in}{0.731763in}}%
\pgfpathclose%
\pgfusepath{fill}%
\end{pgfscope}%
\begin{pgfscope}%
\pgfsetbuttcap%
\pgfsetroundjoin%
\definecolor{currentfill}{rgb}{0.000000,0.000000,0.000000}%
\pgfsetfillcolor{currentfill}%
\pgfsetlinewidth{0.803000pt}%
\definecolor{currentstroke}{rgb}{0.000000,0.000000,0.000000}%
\pgfsetstrokecolor{currentstroke}%
\pgfsetdash{}{0pt}%
\pgfsys@defobject{currentmarker}{\pgfqpoint{0.000000in}{0.000000in}}{\pgfqpoint{0.048611in}{0.000000in}}{%
\pgfpathmoveto{\pgfqpoint{0.000000in}{0.000000in}}%
\pgfpathlineto{\pgfqpoint{0.048611in}{0.000000in}}%
\pgfusepath{stroke,fill}%
}%
\begin{pgfscope}%
\pgfsys@transformshift{5.167416in}{0.731763in}%
\pgfsys@useobject{currentmarker}{}%
\end{pgfscope}%
\end{pgfscope}%
\begin{pgfscope}%
\definecolor{textcolor}{rgb}{0.000000,0.000000,0.000000}%
\pgfsetstrokecolor{textcolor}%
\pgfsetfillcolor{textcolor}%
\pgftext[x=5.264638in, y=0.662319in, left, base]{\color{textcolor}\rmfamily\fontsize{14.000000}{16.800000}\selectfont \(\displaystyle {0}\)}%
\end{pgfscope}%
\begin{pgfscope}%
\pgfsetbuttcap%
\pgfsetroundjoin%
\definecolor{currentfill}{rgb}{0.000000,0.000000,0.000000}%
\pgfsetfillcolor{currentfill}%
\pgfsetlinewidth{0.803000pt}%
\definecolor{currentstroke}{rgb}{0.000000,0.000000,0.000000}%
\pgfsetstrokecolor{currentstroke}%
\pgfsetdash{}{0pt}%
\pgfsys@defobject{currentmarker}{\pgfqpoint{0.000000in}{0.000000in}}{\pgfqpoint{0.048611in}{0.000000in}}{%
\pgfpathmoveto{\pgfqpoint{0.000000in}{0.000000in}}%
\pgfpathlineto{\pgfqpoint{0.048611in}{0.000000in}}%
\pgfusepath{stroke,fill}%
}%
\begin{pgfscope}%
\pgfsys@transformshift{5.167416in}{1.518862in}%
\pgfsys@useobject{currentmarker}{}%
\end{pgfscope}%
\end{pgfscope}%
\begin{pgfscope}%
\definecolor{textcolor}{rgb}{0.000000,0.000000,0.000000}%
\pgfsetstrokecolor{textcolor}%
\pgfsetfillcolor{textcolor}%
\pgftext[x=5.264638in, y=1.449417in, left, base]{\color{textcolor}\rmfamily\fontsize{14.000000}{16.800000}\selectfont \(\displaystyle {10}\)}%
\end{pgfscope}%
\begin{pgfscope}%
\pgfsetbuttcap%
\pgfsetroundjoin%
\definecolor{currentfill}{rgb}{0.000000,0.000000,0.000000}%
\pgfsetfillcolor{currentfill}%
\pgfsetlinewidth{0.803000pt}%
\definecolor{currentstroke}{rgb}{0.000000,0.000000,0.000000}%
\pgfsetstrokecolor{currentstroke}%
\pgfsetdash{}{0pt}%
\pgfsys@defobject{currentmarker}{\pgfqpoint{0.000000in}{0.000000in}}{\pgfqpoint{0.048611in}{0.000000in}}{%
\pgfpathmoveto{\pgfqpoint{0.000000in}{0.000000in}}%
\pgfpathlineto{\pgfqpoint{0.048611in}{0.000000in}}%
\pgfusepath{stroke,fill}%
}%
\begin{pgfscope}%
\pgfsys@transformshift{5.167416in}{2.305960in}%
\pgfsys@useobject{currentmarker}{}%
\end{pgfscope}%
\end{pgfscope}%
\begin{pgfscope}%
\definecolor{textcolor}{rgb}{0.000000,0.000000,0.000000}%
\pgfsetstrokecolor{textcolor}%
\pgfsetfillcolor{textcolor}%
\pgftext[x=5.264638in, y=2.236516in, left, base]{\color{textcolor}\rmfamily\fontsize{14.000000}{16.800000}\selectfont \(\displaystyle {20}\)}%
\end{pgfscope}%
\begin{pgfscope}%
\pgfsetbuttcap%
\pgfsetroundjoin%
\definecolor{currentfill}{rgb}{0.000000,0.000000,0.000000}%
\pgfsetfillcolor{currentfill}%
\pgfsetlinewidth{0.803000pt}%
\definecolor{currentstroke}{rgb}{0.000000,0.000000,0.000000}%
\pgfsetstrokecolor{currentstroke}%
\pgfsetdash{}{0pt}%
\pgfsys@defobject{currentmarker}{\pgfqpoint{0.000000in}{0.000000in}}{\pgfqpoint{0.048611in}{0.000000in}}{%
\pgfpathmoveto{\pgfqpoint{0.000000in}{0.000000in}}%
\pgfpathlineto{\pgfqpoint{0.048611in}{0.000000in}}%
\pgfusepath{stroke,fill}%
}%
\begin{pgfscope}%
\pgfsys@transformshift{5.167416in}{3.093058in}%
\pgfsys@useobject{currentmarker}{}%
\end{pgfscope}%
\end{pgfscope}%
\begin{pgfscope}%
\definecolor{textcolor}{rgb}{0.000000,0.000000,0.000000}%
\pgfsetstrokecolor{textcolor}%
\pgfsetfillcolor{textcolor}%
\pgftext[x=5.264638in, y=3.023614in, left, base]{\color{textcolor}\rmfamily\fontsize{14.000000}{16.800000}\selectfont \(\displaystyle {30}\)}%
\end{pgfscope}%
\begin{pgfscope}%
\definecolor{textcolor}{rgb}{0.000000,0.000000,0.000000}%
\pgfsetstrokecolor{textcolor}%
\pgfsetfillcolor{textcolor}%
\pgftext[x=5.599358in,y=2.227250in,,top,rotate=90.000000]{\color{textcolor}\rmfamily\fontsize{16.000000}{19.200000}\selectfont Number of instances}%
\end{pgfscope}%
\begin{pgfscope}%
\pgfsetrectcap%
\pgfsetmiterjoin%
\pgfsetlinewidth{0.803000pt}%
\definecolor{currentstroke}{rgb}{0.000000,0.000000,0.000000}%
\pgfsetstrokecolor{currentstroke}%
\pgfsetdash{}{0pt}%
\pgfpathmoveto{\pgfqpoint{0.868690in}{0.731763in}}%
\pgfpathlineto{\pgfqpoint{0.868690in}{3.722737in}}%
\pgfusepath{stroke}%
\end{pgfscope}%
\begin{pgfscope}%
\pgfsetrectcap%
\pgfsetmiterjoin%
\pgfsetlinewidth{0.803000pt}%
\definecolor{currentstroke}{rgb}{0.000000,0.000000,0.000000}%
\pgfsetstrokecolor{currentstroke}%
\pgfsetdash{}{0pt}%
\pgfpathmoveto{\pgfqpoint{5.167416in}{0.731763in}}%
\pgfpathlineto{\pgfqpoint{5.167416in}{3.722737in}}%
\pgfusepath{stroke}%
\end{pgfscope}%
\begin{pgfscope}%
\pgfsetrectcap%
\pgfsetmiterjoin%
\pgfsetlinewidth{0.803000pt}%
\definecolor{currentstroke}{rgb}{0.000000,0.000000,0.000000}%
\pgfsetstrokecolor{currentstroke}%
\pgfsetdash{}{0pt}%
\pgfpathmoveto{\pgfqpoint{0.868690in}{0.731763in}}%
\pgfpathlineto{\pgfqpoint{5.167416in}{0.731763in}}%
\pgfusepath{stroke}%
\end{pgfscope}%
\begin{pgfscope}%
\pgfsetrectcap%
\pgfsetmiterjoin%
\pgfsetlinewidth{0.803000pt}%
\definecolor{currentstroke}{rgb}{0.000000,0.000000,0.000000}%
\pgfsetstrokecolor{currentstroke}%
\pgfsetdash{}{0pt}%
\pgfpathmoveto{\pgfqpoint{0.868690in}{3.722737in}}%
\pgfpathlineto{\pgfqpoint{5.167416in}{3.722737in}}%
\pgfusepath{stroke}%
\end{pgfscope}%
\end{pgfpicture}%
\makeatother%
\endgroup%

%% file: figures/Comp/PIandOCVSLambda.pgf
\begingroup%
\makeatletter%
\begin{pgfpicture}%
\pgfpathrectangle{\pgfpointorigin}{\pgfqpoint{6.000000in}{4.000000in}}%
\pgfusepath{use as bounding box, clip}%
\begin{pgfscope}%
\pgfsetbuttcap%
\pgfsetmiterjoin%
\definecolor{currentfill}{rgb}{1.000000,1.000000,1.000000}%
\pgfsetfillcolor{currentfill}%
\pgfsetlinewidth{0.000000pt}%
\definecolor{currentstroke}{rgb}{1.000000,1.000000,1.000000}%
\pgfsetstrokecolor{currentstroke}%
\pgfsetdash{}{0pt}%
\pgfpathmoveto{\pgfqpoint{0.000000in}{0.000000in}}%
\pgfpathlineto{\pgfqpoint{6.000000in}{0.000000in}}%
\pgfpathlineto{\pgfqpoint{6.000000in}{4.000000in}}%
\pgfpathlineto{\pgfqpoint{0.000000in}{4.000000in}}%
\pgfpathlineto{\pgfqpoint{0.000000in}{0.000000in}}%
\pgfpathclose%
\pgfusepath{fill}%
\end{pgfscope}%
\begin{pgfscope}%
\pgfsetbuttcap%
\pgfsetmiterjoin%
\definecolor{currentfill}{rgb}{1.000000,1.000000,1.000000}%
\pgfsetfillcolor{currentfill}%
\pgfsetlinewidth{0.000000pt}%
\definecolor{currentstroke}{rgb}{0.000000,0.000000,0.000000}%
\pgfsetstrokecolor{currentstroke}%
\pgfsetstrokeopacity{0.000000}%
\pgfsetdash{}{0pt}%
\pgfpathmoveto{\pgfqpoint{1.128611in}{0.731763in}}%
\pgfpathlineto{\pgfqpoint{5.706667in}{0.731763in}}%
\pgfpathlineto{\pgfqpoint{5.706667in}{3.790000in}}%
\pgfpathlineto{\pgfqpoint{1.128611in}{3.790000in}}%
\pgfpathlineto{\pgfqpoint{1.128611in}{0.731763in}}%
\pgfpathclose%
\pgfusepath{fill}%
\end{pgfscope}%
\begin{pgfscope}%
\pgfsetbuttcap%
\pgfsetroundjoin%
\definecolor{currentfill}{rgb}{0.000000,0.000000,0.000000}%
\pgfsetfillcolor{currentfill}%
\pgfsetlinewidth{0.803000pt}%
\definecolor{currentstroke}{rgb}{0.000000,0.000000,0.000000}%
\pgfsetstrokecolor{currentstroke}%
\pgfsetdash{}{0pt}%
\pgfsys@defobject{currentmarker}{\pgfqpoint{0.000000in}{-0.048611in}}{\pgfqpoint{0.000000in}{0.000000in}}{%
\pgfpathmoveto{\pgfqpoint{0.000000in}{0.000000in}}%
\pgfpathlineto{\pgfqpoint{0.000000in}{-0.048611in}}%
\pgfusepath{stroke,fill}%
}%
\begin{pgfscope}%
\pgfsys@transformshift{1.128611in}{0.731763in}%
\pgfsys@useobject{currentmarker}{}%
\end{pgfscope}%
\end{pgfscope}%
\begin{pgfscope}%
\definecolor{textcolor}{rgb}{0.000000,0.000000,0.000000}%
\pgfsetstrokecolor{textcolor}%
\pgfsetfillcolor{textcolor}%
\pgftext[x=1.128611in,y=0.634541in,,top]{\color{textcolor}\rmfamily\fontsize{14.000000}{16.800000}\selectfont \(\displaystyle {0.0}\)}%
\end{pgfscope}%
\begin{pgfscope}%
\pgfsetbuttcap%
\pgfsetroundjoin%
\definecolor{currentfill}{rgb}{0.000000,0.000000,0.000000}%
\pgfsetfillcolor{currentfill}%
\pgfsetlinewidth{0.803000pt}%
\definecolor{currentstroke}{rgb}{0.000000,0.000000,0.000000}%
\pgfsetstrokecolor{currentstroke}%
\pgfsetdash{}{0pt}%
\pgfsys@defobject{currentmarker}{\pgfqpoint{0.000000in}{-0.048611in}}{\pgfqpoint{0.000000in}{0.000000in}}{%
\pgfpathmoveto{\pgfqpoint{0.000000in}{0.000000in}}%
\pgfpathlineto{\pgfqpoint{0.000000in}{-0.048611in}}%
\pgfusepath{stroke,fill}%
}%
\begin{pgfscope}%
\pgfsys@transformshift{1.960985in}{0.731763in}%
\pgfsys@useobject{currentmarker}{}%
\end{pgfscope}%
\end{pgfscope}%
\begin{pgfscope}%
\definecolor{textcolor}{rgb}{0.000000,0.000000,0.000000}%
\pgfsetstrokecolor{textcolor}%
\pgfsetfillcolor{textcolor}%
\pgftext[x=1.960985in,y=0.634541in,,top]{\color{textcolor}\rmfamily\fontsize{14.000000}{16.800000}\selectfont \(\displaystyle {0.2}\)}%
\end{pgfscope}%
\begin{pgfscope}%
\pgfsetbuttcap%
\pgfsetroundjoin%
\definecolor{currentfill}{rgb}{0.000000,0.000000,0.000000}%
\pgfsetfillcolor{currentfill}%
\pgfsetlinewidth{0.803000pt}%
\definecolor{currentstroke}{rgb}{0.000000,0.000000,0.000000}%
\pgfsetstrokecolor{currentstroke}%
\pgfsetdash{}{0pt}%
\pgfsys@defobject{currentmarker}{\pgfqpoint{0.000000in}{-0.048611in}}{\pgfqpoint{0.000000in}{0.000000in}}{%
\pgfpathmoveto{\pgfqpoint{0.000000in}{0.000000in}}%
\pgfpathlineto{\pgfqpoint{0.000000in}{-0.048611in}}%
\pgfusepath{stroke,fill}%
}%
\begin{pgfscope}%
\pgfsys@transformshift{2.793359in}{0.731763in}%
\pgfsys@useobject{currentmarker}{}%
\end{pgfscope}%
\end{pgfscope}%
\begin{pgfscope}%
\definecolor{textcolor}{rgb}{0.000000,0.000000,0.000000}%
\pgfsetstrokecolor{textcolor}%
\pgfsetfillcolor{textcolor}%
\pgftext[x=2.793359in,y=0.634541in,,top]{\color{textcolor}\rmfamily\fontsize{14.000000}{16.800000}\selectfont \(\displaystyle {0.4}\)}%
\end{pgfscope}%
\begin{pgfscope}%
\pgfsetbuttcap%
\pgfsetroundjoin%
\definecolor{currentfill}{rgb}{0.000000,0.000000,0.000000}%
\pgfsetfillcolor{currentfill}%
\pgfsetlinewidth{0.803000pt}%
\definecolor{currentstroke}{rgb}{0.000000,0.000000,0.000000}%
\pgfsetstrokecolor{currentstroke}%
\pgfsetdash{}{0pt}%
\pgfsys@defobject{currentmarker}{\pgfqpoint{0.000000in}{-0.048611in}}{\pgfqpoint{0.000000in}{0.000000in}}{%
\pgfpathmoveto{\pgfqpoint{0.000000in}{0.000000in}}%
\pgfpathlineto{\pgfqpoint{0.000000in}{-0.048611in}}%
\pgfusepath{stroke,fill}%
}%
\begin{pgfscope}%
\pgfsys@transformshift{3.625732in}{0.731763in}%
\pgfsys@useobject{currentmarker}{}%
\end{pgfscope}%
\end{pgfscope}%
\begin{pgfscope}%
\definecolor{textcolor}{rgb}{0.000000,0.000000,0.000000}%
\pgfsetstrokecolor{textcolor}%
\pgfsetfillcolor{textcolor}%
\pgftext[x=3.625732in,y=0.634541in,,top]{\color{textcolor}\rmfamily\fontsize{14.000000}{16.800000}\selectfont \(\displaystyle {0.6}\)}%
\end{pgfscope}%
\begin{pgfscope}%
\pgfsetbuttcap%
\pgfsetroundjoin%
\definecolor{currentfill}{rgb}{0.000000,0.000000,0.000000}%
\pgfsetfillcolor{currentfill}%
\pgfsetlinewidth{0.803000pt}%
\definecolor{currentstroke}{rgb}{0.000000,0.000000,0.000000}%
\pgfsetstrokecolor{currentstroke}%
\pgfsetdash{}{0pt}%
\pgfsys@defobject{currentmarker}{\pgfqpoint{0.000000in}{-0.048611in}}{\pgfqpoint{0.000000in}{0.000000in}}{%
\pgfpathmoveto{\pgfqpoint{0.000000in}{0.000000in}}%
\pgfpathlineto{\pgfqpoint{0.000000in}{-0.048611in}}%
\pgfusepath{stroke,fill}%
}%
\begin{pgfscope}%
\pgfsys@transformshift{4.458106in}{0.731763in}%
\pgfsys@useobject{currentmarker}{}%
\end{pgfscope}%
\end{pgfscope}%
\begin{pgfscope}%
\definecolor{textcolor}{rgb}{0.000000,0.000000,0.000000}%
\pgfsetstrokecolor{textcolor}%
\pgfsetfillcolor{textcolor}%
\pgftext[x=4.458106in,y=0.634541in,,top]{\color{textcolor}\rmfamily\fontsize{14.000000}{16.800000}\selectfont \(\displaystyle {0.8}\)}%
\end{pgfscope}%
\begin{pgfscope}%
\pgfsetbuttcap%
\pgfsetroundjoin%
\definecolor{currentfill}{rgb}{0.000000,0.000000,0.000000}%
\pgfsetfillcolor{currentfill}%
\pgfsetlinewidth{0.803000pt}%
\definecolor{currentstroke}{rgb}{0.000000,0.000000,0.000000}%
\pgfsetstrokecolor{currentstroke}%
\pgfsetdash{}{0pt}%
\pgfsys@defobject{currentmarker}{\pgfqpoint{0.000000in}{-0.048611in}}{\pgfqpoint{0.000000in}{0.000000in}}{%
\pgfpathmoveto{\pgfqpoint{0.000000in}{0.000000in}}%
\pgfpathlineto{\pgfqpoint{0.000000in}{-0.048611in}}%
\pgfusepath{stroke,fill}%
}%
\begin{pgfscope}%
\pgfsys@transformshift{5.290480in}{0.731763in}%
\pgfsys@useobject{currentmarker}{}%
\end{pgfscope}%
\end{pgfscope}%
\begin{pgfscope}%
\definecolor{textcolor}{rgb}{0.000000,0.000000,0.000000}%
\pgfsetstrokecolor{textcolor}%
\pgfsetfillcolor{textcolor}%
\pgftext[x=5.290480in,y=0.634541in,,top]{\color{textcolor}\rmfamily\fontsize{14.000000}{16.800000}\selectfont \(\displaystyle {1.0}\)}%
\end{pgfscope}%
\begin{pgfscope}%
\definecolor{textcolor}{rgb}{0.000000,0.000000,0.000000}%
\pgfsetstrokecolor{textcolor}%
\pgfsetfillcolor{textcolor}%
\pgftext[x=3.417639in,y=0.401208in,,top]{\color{textcolor}\rmfamily\fontsize{16.000000}{19.200000}\selectfont \(\displaystyle \lambda\)}%
\end{pgfscope}%
\begin{pgfscope}%
\pgfsetbuttcap%
\pgfsetroundjoin%
\definecolor{currentfill}{rgb}{0.000000,0.000000,0.000000}%
\pgfsetfillcolor{currentfill}%
\pgfsetlinewidth{0.803000pt}%
\definecolor{currentstroke}{rgb}{0.000000,0.000000,0.000000}%
\pgfsetstrokecolor{currentstroke}%
\pgfsetdash{}{0pt}%
\pgfsys@defobject{currentmarker}{\pgfqpoint{-0.048611in}{0.000000in}}{\pgfqpoint{-0.000000in}{0.000000in}}{%
\pgfpathmoveto{\pgfqpoint{-0.000000in}{0.000000in}}%
\pgfpathlineto{\pgfqpoint{-0.048611in}{0.000000in}}%
\pgfusepath{stroke,fill}%
}%
\begin{pgfscope}%
\pgfsys@transformshift{1.128611in}{1.139788in}%
\pgfsys@useobject{currentmarker}{}%
\end{pgfscope}%
\end{pgfscope}%
\begin{pgfscope}%
\definecolor{textcolor}{rgb}{0.000000,0.000000,0.000000}%
\pgfsetstrokecolor{textcolor}%
\pgfsetfillcolor{textcolor}%
\pgftext[x=0.541812in, y=1.070344in, left, base]{\color{textcolor}\rmfamily\fontsize{14.000000}{16.800000}\selectfont \(\displaystyle {20000}\)}%
\end{pgfscope}%
\begin{pgfscope}%
\pgfsetbuttcap%
\pgfsetroundjoin%
\definecolor{currentfill}{rgb}{0.000000,0.000000,0.000000}%
\pgfsetfillcolor{currentfill}%
\pgfsetlinewidth{0.803000pt}%
\definecolor{currentstroke}{rgb}{0.000000,0.000000,0.000000}%
\pgfsetstrokecolor{currentstroke}%
\pgfsetdash{}{0pt}%
\pgfsys@defobject{currentmarker}{\pgfqpoint{-0.048611in}{0.000000in}}{\pgfqpoint{-0.000000in}{0.000000in}}{%
\pgfpathmoveto{\pgfqpoint{-0.000000in}{0.000000in}}%
\pgfpathlineto{\pgfqpoint{-0.048611in}{0.000000in}}%
\pgfusepath{stroke,fill}%
}%
\begin{pgfscope}%
\pgfsys@transformshift{1.128611in}{1.778521in}%
\pgfsys@useobject{currentmarker}{}%
\end{pgfscope}%
\end{pgfscope}%
\begin{pgfscope}%
\definecolor{textcolor}{rgb}{0.000000,0.000000,0.000000}%
\pgfsetstrokecolor{textcolor}%
\pgfsetfillcolor{textcolor}%
\pgftext[x=0.541812in, y=1.709077in, left, base]{\color{textcolor}\rmfamily\fontsize{14.000000}{16.800000}\selectfont \(\displaystyle {40000}\)}%
\end{pgfscope}%
\begin{pgfscope}%
\pgfsetbuttcap%
\pgfsetroundjoin%
\definecolor{currentfill}{rgb}{0.000000,0.000000,0.000000}%
\pgfsetfillcolor{currentfill}%
\pgfsetlinewidth{0.803000pt}%
\definecolor{currentstroke}{rgb}{0.000000,0.000000,0.000000}%
\pgfsetstrokecolor{currentstroke}%
\pgfsetdash{}{0pt}%
\pgfsys@defobject{currentmarker}{\pgfqpoint{-0.048611in}{0.000000in}}{\pgfqpoint{-0.000000in}{0.000000in}}{%
\pgfpathmoveto{\pgfqpoint{-0.000000in}{0.000000in}}%
\pgfpathlineto{\pgfqpoint{-0.048611in}{0.000000in}}%
\pgfusepath{stroke,fill}%
}%
\begin{pgfscope}%
\pgfsys@transformshift{1.128611in}{2.417253in}%
\pgfsys@useobject{currentmarker}{}%
\end{pgfscope}%
\end{pgfscope}%
\begin{pgfscope}%
\definecolor{textcolor}{rgb}{0.000000,0.000000,0.000000}%
\pgfsetstrokecolor{textcolor}%
\pgfsetfillcolor{textcolor}%
\pgftext[x=0.541812in, y=2.347809in, left, base]{\color{textcolor}\rmfamily\fontsize{14.000000}{16.800000}\selectfont \(\displaystyle {60000}\)}%
\end{pgfscope}%
\begin{pgfscope}%
\pgfsetbuttcap%
\pgfsetroundjoin%
\definecolor{currentfill}{rgb}{0.000000,0.000000,0.000000}%
\pgfsetfillcolor{currentfill}%
\pgfsetlinewidth{0.803000pt}%
\definecolor{currentstroke}{rgb}{0.000000,0.000000,0.000000}%
\pgfsetstrokecolor{currentstroke}%
\pgfsetdash{}{0pt}%
\pgfsys@defobject{currentmarker}{\pgfqpoint{-0.048611in}{0.000000in}}{\pgfqpoint{-0.000000in}{0.000000in}}{%
\pgfpathmoveto{\pgfqpoint{-0.000000in}{0.000000in}}%
\pgfpathlineto{\pgfqpoint{-0.048611in}{0.000000in}}%
\pgfusepath{stroke,fill}%
}%
\begin{pgfscope}%
\pgfsys@transformshift{1.128611in}{3.055986in}%
\pgfsys@useobject{currentmarker}{}%
\end{pgfscope}%
\end{pgfscope}%
\begin{pgfscope}%
\definecolor{textcolor}{rgb}{0.000000,0.000000,0.000000}%
\pgfsetstrokecolor{textcolor}%
\pgfsetfillcolor{textcolor}%
\pgftext[x=0.541812in, y=2.986542in, left, base]{\color{textcolor}\rmfamily\fontsize{14.000000}{16.800000}\selectfont \(\displaystyle {80000}\)}%
\end{pgfscope}%
\begin{pgfscope}%
\pgfsetbuttcap%
\pgfsetroundjoin%
\definecolor{currentfill}{rgb}{0.000000,0.000000,0.000000}%
\pgfsetfillcolor{currentfill}%
\pgfsetlinewidth{0.803000pt}%
\definecolor{currentstroke}{rgb}{0.000000,0.000000,0.000000}%
\pgfsetstrokecolor{currentstroke}%
\pgfsetdash{}{0pt}%
\pgfsys@defobject{currentmarker}{\pgfqpoint{-0.048611in}{0.000000in}}{\pgfqpoint{-0.000000in}{0.000000in}}{%
\pgfpathmoveto{\pgfqpoint{-0.000000in}{0.000000in}}%
\pgfpathlineto{\pgfqpoint{-0.048611in}{0.000000in}}%
\pgfusepath{stroke,fill}%
}%
\begin{pgfscope}%
\pgfsys@transformshift{1.128611in}{3.694718in}%
\pgfsys@useobject{currentmarker}{}%
\end{pgfscope}%
\end{pgfscope}%
\begin{pgfscope}%
\definecolor{textcolor}{rgb}{0.000000,0.000000,0.000000}%
\pgfsetstrokecolor{textcolor}%
\pgfsetfillcolor{textcolor}%
\pgftext[x=0.443896in, y=3.625274in, left, base]{\color{textcolor}\rmfamily\fontsize{14.000000}{16.800000}\selectfont \(\displaystyle {100000}\)}%
\end{pgfscope}%
\begin{pgfscope}%
\definecolor{textcolor}{rgb}{0.000000,0.000000,0.000000}%
\pgfsetstrokecolor{textcolor}%
\pgfsetfillcolor{textcolor}%
\pgftext[x=0.388341in,y=2.260882in,,bottom,rotate=90.000000]{\color{textcolor}\rmfamily\fontsize{16.000000}{19.200000}\selectfont DKK}%
\end{pgfscope}%
\begin{pgfscope}%
\pgfpathrectangle{\pgfqpoint{1.128611in}{0.731763in}}{\pgfqpoint{4.578056in}{3.058237in}}%
\pgfusepath{clip}%
\pgfsetrectcap%
\pgfsetroundjoin%
\pgfsetlinewidth{1.505625pt}%
\definecolor{currentstroke}{rgb}{0.000000,0.000000,0.000000}%
\pgfsetstrokecolor{currentstroke}%
\pgfsetdash{}{0pt}%
\pgfpathmoveto{\pgfqpoint{1.544798in}{0.870774in}}%
\pgfpathlineto{\pgfqpoint{2.169078in}{1.385068in}}%
\pgfpathlineto{\pgfqpoint{3.209545in}{2.228034in}}%
\pgfpathlineto{\pgfqpoint{4.250013in}{2.999232in}}%
\pgfpathlineto{\pgfqpoint{5.290480in}{3.650989in}}%
\pgfusepath{stroke}%
\end{pgfscope}%
\begin{pgfscope}%
\pgfpathrectangle{\pgfqpoint{1.128611in}{0.731763in}}{\pgfqpoint{4.578056in}{3.058237in}}%
\pgfusepath{clip}%
\pgfsetbuttcap%
\pgfsetroundjoin%
\definecolor{currentfill}{rgb}{0.000000,0.000000,0.000000}%
\pgfsetfillcolor{currentfill}%
\pgfsetlinewidth{1.003750pt}%
\definecolor{currentstroke}{rgb}{0.000000,0.000000,0.000000}%
\pgfsetstrokecolor{currentstroke}%
\pgfsetdash{}{0pt}%
\pgfsys@defobject{currentmarker}{\pgfqpoint{-0.041667in}{-0.041667in}}{\pgfqpoint{0.041667in}{0.041667in}}{%
\pgfpathmoveto{\pgfqpoint{0.000000in}{-0.041667in}}%
\pgfpathcurveto{\pgfqpoint{0.011050in}{-0.041667in}}{\pgfqpoint{0.021649in}{-0.037276in}}{\pgfqpoint{0.029463in}{-0.029463in}}%
\pgfpathcurveto{\pgfqpoint{0.037276in}{-0.021649in}}{\pgfqpoint{0.041667in}{-0.011050in}}{\pgfqpoint{0.041667in}{0.000000in}}%
\pgfpathcurveto{\pgfqpoint{0.041667in}{0.011050in}}{\pgfqpoint{0.037276in}{0.021649in}}{\pgfqpoint{0.029463in}{0.029463in}}%
\pgfpathcurveto{\pgfqpoint{0.021649in}{0.037276in}}{\pgfqpoint{0.011050in}{0.041667in}}{\pgfqpoint{0.000000in}{0.041667in}}%
\pgfpathcurveto{\pgfqpoint{-0.011050in}{0.041667in}}{\pgfqpoint{-0.021649in}{0.037276in}}{\pgfqpoint{-0.029463in}{0.029463in}}%
\pgfpathcurveto{\pgfqpoint{-0.037276in}{0.021649in}}{\pgfqpoint{-0.041667in}{0.011050in}}{\pgfqpoint{-0.041667in}{0.000000in}}%
\pgfpathcurveto{\pgfqpoint{-0.041667in}{-0.011050in}}{\pgfqpoint{-0.037276in}{-0.021649in}}{\pgfqpoint{-0.029463in}{-0.029463in}}%
\pgfpathcurveto{\pgfqpoint{-0.021649in}{-0.037276in}}{\pgfqpoint{-0.011050in}{-0.041667in}}{\pgfqpoint{0.000000in}{-0.041667in}}%
\pgfpathlineto{\pgfqpoint{0.000000in}{-0.041667in}}%
\pgfpathclose%
\pgfusepath{stroke,fill}%
}%
\begin{pgfscope}%
\pgfsys@transformshift{1.544798in}{0.870774in}%
\pgfsys@useobject{currentmarker}{}%
\end{pgfscope}%
\begin{pgfscope}%
\pgfsys@transformshift{2.169078in}{1.385068in}%
\pgfsys@useobject{currentmarker}{}%
\end{pgfscope}%
\begin{pgfscope}%
\pgfsys@transformshift{3.209545in}{2.228034in}%
\pgfsys@useobject{currentmarker}{}%
\end{pgfscope}%
\begin{pgfscope}%
\pgfsys@transformshift{4.250013in}{2.999232in}%
\pgfsys@useobject{currentmarker}{}%
\end{pgfscope}%
\begin{pgfscope}%
\pgfsys@transformshift{5.290480in}{3.650989in}%
\pgfsys@useobject{currentmarker}{}%
\end{pgfscope}%
\end{pgfscope}%
\begin{pgfscope}%
\pgfpathrectangle{\pgfqpoint{1.128611in}{0.731763in}}{\pgfqpoint{4.578056in}{3.058237in}}%
\pgfusepath{clip}%
\pgfsetrectcap%
\pgfsetroundjoin%
\pgfsetlinewidth{1.505625pt}%
\definecolor{currentstroke}{rgb}{0.678431,0.709804,0.741176}%
\pgfsetstrokecolor{currentstroke}%
\pgfsetdash{}{0pt}%
\pgfpathmoveto{\pgfqpoint{1.544798in}{1.004590in}}%
\pgfpathlineto{\pgfqpoint{2.169078in}{1.102074in}}%
\pgfpathlineto{\pgfqpoint{3.209545in}{1.194949in}}%
\pgfpathlineto{\pgfqpoint{4.250013in}{1.330131in}}%
\pgfpathlineto{\pgfqpoint{5.290480in}{1.524518in}}%
\pgfusepath{stroke}%
\end{pgfscope}%
\begin{pgfscope}%
\pgfpathrectangle{\pgfqpoint{1.128611in}{0.731763in}}{\pgfqpoint{4.578056in}{3.058237in}}%
\pgfusepath{clip}%
\pgfsetbuttcap%
\pgfsetroundjoin%
\definecolor{currentfill}{rgb}{0.678431,0.709804,0.741176}%
\pgfsetfillcolor{currentfill}%
\pgfsetlinewidth{1.003750pt}%
\definecolor{currentstroke}{rgb}{0.678431,0.709804,0.741176}%
\pgfsetstrokecolor{currentstroke}%
\pgfsetdash{}{0pt}%
\pgfsys@defobject{currentmarker}{\pgfqpoint{-0.041667in}{-0.041667in}}{\pgfqpoint{0.041667in}{0.041667in}}{%
\pgfpathmoveto{\pgfqpoint{0.000000in}{-0.041667in}}%
\pgfpathcurveto{\pgfqpoint{0.011050in}{-0.041667in}}{\pgfqpoint{0.021649in}{-0.037276in}}{\pgfqpoint{0.029463in}{-0.029463in}}%
\pgfpathcurveto{\pgfqpoint{0.037276in}{-0.021649in}}{\pgfqpoint{0.041667in}{-0.011050in}}{\pgfqpoint{0.041667in}{0.000000in}}%
\pgfpathcurveto{\pgfqpoint{0.041667in}{0.011050in}}{\pgfqpoint{0.037276in}{0.021649in}}{\pgfqpoint{0.029463in}{0.029463in}}%
\pgfpathcurveto{\pgfqpoint{0.021649in}{0.037276in}}{\pgfqpoint{0.011050in}{0.041667in}}{\pgfqpoint{0.000000in}{0.041667in}}%
\pgfpathcurveto{\pgfqpoint{-0.011050in}{0.041667in}}{\pgfqpoint{-0.021649in}{0.037276in}}{\pgfqpoint{-0.029463in}{0.029463in}}%
\pgfpathcurveto{\pgfqpoint{-0.037276in}{0.021649in}}{\pgfqpoint{-0.041667in}{0.011050in}}{\pgfqpoint{-0.041667in}{0.000000in}}%
\pgfpathcurveto{\pgfqpoint{-0.041667in}{-0.011050in}}{\pgfqpoint{-0.037276in}{-0.021649in}}{\pgfqpoint{-0.029463in}{-0.029463in}}%
\pgfpathcurveto{\pgfqpoint{-0.021649in}{-0.037276in}}{\pgfqpoint{-0.011050in}{-0.041667in}}{\pgfqpoint{0.000000in}{-0.041667in}}%
\pgfpathlineto{\pgfqpoint{0.000000in}{-0.041667in}}%
\pgfpathclose%
\pgfusepath{stroke,fill}%
}%
\begin{pgfscope}%
\pgfsys@transformshift{1.544798in}{1.004590in}%
\pgfsys@useobject{currentmarker}{}%
\end{pgfscope}%
\begin{pgfscope}%
\pgfsys@transformshift{2.169078in}{1.102074in}%
\pgfsys@useobject{currentmarker}{}%
\end{pgfscope}%
\begin{pgfscope}%
\pgfsys@transformshift{3.209545in}{1.194949in}%
\pgfsys@useobject{currentmarker}{}%
\end{pgfscope}%
\begin{pgfscope}%
\pgfsys@transformshift{4.250013in}{1.330131in}%
\pgfsys@useobject{currentmarker}{}%
\end{pgfscope}%
\begin{pgfscope}%
\pgfsys@transformshift{5.290480in}{1.524518in}%
\pgfsys@useobject{currentmarker}{}%
\end{pgfscope}%
\end{pgfscope}%
\begin{pgfscope}%
\pgfsetrectcap%
\pgfsetmiterjoin%
\pgfsetlinewidth{0.803000pt}%
\definecolor{currentstroke}{rgb}{0.000000,0.000000,0.000000}%
\pgfsetstrokecolor{currentstroke}%
\pgfsetdash{}{0pt}%
\pgfpathmoveto{\pgfqpoint{1.128611in}{0.731763in}}%
\pgfpathlineto{\pgfqpoint{1.128611in}{3.790000in}}%
\pgfusepath{stroke}%
\end{pgfscope}%
\begin{pgfscope}%
\pgfsetrectcap%
\pgfsetmiterjoin%
\pgfsetlinewidth{0.803000pt}%
\definecolor{currentstroke}{rgb}{0.000000,0.000000,0.000000}%
\pgfsetstrokecolor{currentstroke}%
\pgfsetdash{}{0pt}%
\pgfpathmoveto{\pgfqpoint{5.706667in}{0.731763in}}%
\pgfpathlineto{\pgfqpoint{5.706667in}{3.790000in}}%
\pgfusepath{stroke}%
\end{pgfscope}%
\begin{pgfscope}%
\pgfsetrectcap%
\pgfsetmiterjoin%
\pgfsetlinewidth{0.803000pt}%
\definecolor{currentstroke}{rgb}{0.000000,0.000000,0.000000}%
\pgfsetstrokecolor{currentstroke}%
\pgfsetdash{}{0pt}%
\pgfpathmoveto{\pgfqpoint{1.128611in}{0.731763in}}%
\pgfpathlineto{\pgfqpoint{5.706667in}{0.731763in}}%
\pgfusepath{stroke}%
\end{pgfscope}%
\begin{pgfscope}%
\pgfsetrectcap%
\pgfsetmiterjoin%
\pgfsetlinewidth{0.803000pt}%
\definecolor{currentstroke}{rgb}{0.000000,0.000000,0.000000}%
\pgfsetstrokecolor{currentstroke}%
\pgfsetdash{}{0pt}%
\pgfpathmoveto{\pgfqpoint{1.128611in}{3.790000in}}%
\pgfpathlineto{\pgfqpoint{5.706667in}{3.790000in}}%
\pgfusepath{stroke}%
\end{pgfscope}%
\begin{pgfscope}%
\pgfsetbuttcap%
\pgfsetmiterjoin%
\definecolor{currentfill}{rgb}{1.000000,1.000000,1.000000}%
\pgfsetfillcolor{currentfill}%
\pgfsetfillopacity{0.800000}%
\pgfsetlinewidth{1.003750pt}%
\definecolor{currentstroke}{rgb}{0.800000,0.800000,0.800000}%
\pgfsetstrokecolor{currentstroke}%
\pgfsetstrokeopacity{0.800000}%
\pgfsetdash{}{0pt}%
\pgfpathmoveto{\pgfqpoint{1.284167in}{2.963303in}}%
\pgfpathlineto{\pgfqpoint{3.472596in}{2.963303in}}%
\pgfpathquadraticcurveto{\pgfqpoint{3.517041in}{2.963303in}}{\pgfqpoint{3.517041in}{3.007747in}}%
\pgfpathlineto{\pgfqpoint{3.517041in}{3.634444in}}%
\pgfpathquadraticcurveto{\pgfqpoint{3.517041in}{3.678889in}}{\pgfqpoint{3.472596in}{3.678889in}}%
\pgfpathlineto{\pgfqpoint{1.284167in}{3.678889in}}%
\pgfpathquadraticcurveto{\pgfqpoint{1.239722in}{3.678889in}}{\pgfqpoint{1.239722in}{3.634444in}}%
\pgfpathlineto{\pgfqpoint{1.239722in}{3.007747in}}%
\pgfpathquadraticcurveto{\pgfqpoint{1.239722in}{2.963303in}}{\pgfqpoint{1.284167in}{2.963303in}}%
\pgfpathlineto{\pgfqpoint{1.284167in}{2.963303in}}%
\pgfpathclose%
\pgfusepath{stroke,fill}%
\end{pgfscope}%
\begin{pgfscope}%
\pgfsetrectcap%
\pgfsetroundjoin%
\pgfsetlinewidth{1.505625pt}%
\definecolor{currentstroke}{rgb}{0.000000,0.000000,0.000000}%
\pgfsetstrokecolor{currentstroke}%
\pgfsetdash{}{0pt}%
\pgfpathmoveto{\pgfqpoint{1.328611in}{3.501111in}}%
\pgfpathlineto{\pgfqpoint{1.550833in}{3.501111in}}%
\pgfpathlineto{\pgfqpoint{1.773056in}{3.501111in}}%
\pgfusepath{stroke}%
\end{pgfscope}%
\begin{pgfscope}%
\pgfsetbuttcap%
\pgfsetroundjoin%
\definecolor{currentfill}{rgb}{0.000000,0.000000,0.000000}%
\pgfsetfillcolor{currentfill}%
\pgfsetlinewidth{1.003750pt}%
\definecolor{currentstroke}{rgb}{0.000000,0.000000,0.000000}%
\pgfsetstrokecolor{currentstroke}%
\pgfsetdash{}{0pt}%
\pgfsys@defobject{currentmarker}{\pgfqpoint{-0.041667in}{-0.041667in}}{\pgfqpoint{0.041667in}{0.041667in}}{%
\pgfpathmoveto{\pgfqpoint{0.000000in}{-0.041667in}}%
\pgfpathcurveto{\pgfqpoint{0.011050in}{-0.041667in}}{\pgfqpoint{0.021649in}{-0.037276in}}{\pgfqpoint{0.029463in}{-0.029463in}}%
\pgfpathcurveto{\pgfqpoint{0.037276in}{-0.021649in}}{\pgfqpoint{0.041667in}{-0.011050in}}{\pgfqpoint{0.041667in}{0.000000in}}%
\pgfpathcurveto{\pgfqpoint{0.041667in}{0.011050in}}{\pgfqpoint{0.037276in}{0.021649in}}{\pgfqpoint{0.029463in}{0.029463in}}%
\pgfpathcurveto{\pgfqpoint{0.021649in}{0.037276in}}{\pgfqpoint{0.011050in}{0.041667in}}{\pgfqpoint{0.000000in}{0.041667in}}%
\pgfpathcurveto{\pgfqpoint{-0.011050in}{0.041667in}}{\pgfqpoint{-0.021649in}{0.037276in}}{\pgfqpoint{-0.029463in}{0.029463in}}%
\pgfpathcurveto{\pgfqpoint{-0.037276in}{0.021649in}}{\pgfqpoint{-0.041667in}{0.011050in}}{\pgfqpoint{-0.041667in}{0.000000in}}%
\pgfpathcurveto{\pgfqpoint{-0.041667in}{-0.011050in}}{\pgfqpoint{-0.037276in}{-0.021649in}}{\pgfqpoint{-0.029463in}{-0.029463in}}%
\pgfpathcurveto{\pgfqpoint{-0.021649in}{-0.037276in}}{\pgfqpoint{-0.011050in}{-0.041667in}}{\pgfqpoint{0.000000in}{-0.041667in}}%
\pgfpathlineto{\pgfqpoint{0.000000in}{-0.041667in}}%
\pgfpathclose%
\pgfusepath{stroke,fill}%
}%
\begin{pgfscope}%
\pgfsys@transformshift{1.550833in}{3.501111in}%
\pgfsys@useobject{currentmarker}{}%
\end{pgfscope}%
\end{pgfscope}%
\begin{pgfscope}%
\definecolor{textcolor}{rgb}{0.000000,0.000000,0.000000}%
\pgfsetstrokecolor{textcolor}%
\pgfsetfillcolor{textcolor}%
\pgftext[x=1.950833in,y=3.423333in,left,base]{\color{textcolor}\rmfamily\fontsize{16.000000}{19.200000}\selectfont Passenger costs}%
\end{pgfscope}%
\begin{pgfscope}%
\pgfsetrectcap%
\pgfsetroundjoin%
\pgfsetlinewidth{1.505625pt}%
\definecolor{currentstroke}{rgb}{0.678431,0.709804,0.741176}%
\pgfsetstrokecolor{currentstroke}%
\pgfsetdash{}{0pt}%
\pgfpathmoveto{\pgfqpoint{1.328611in}{3.176651in}}%
\pgfpathlineto{\pgfqpoint{1.550833in}{3.176651in}}%
\pgfpathlineto{\pgfqpoint{1.773056in}{3.176651in}}%
\pgfusepath{stroke}%
\end{pgfscope}%
\begin{pgfscope}%
\pgfsetbuttcap%
\pgfsetroundjoin%
\definecolor{currentfill}{rgb}{0.678431,0.709804,0.741176}%
\pgfsetfillcolor{currentfill}%
\pgfsetlinewidth{1.003750pt}%
\definecolor{currentstroke}{rgb}{0.678431,0.709804,0.741176}%
\pgfsetstrokecolor{currentstroke}%
\pgfsetdash{}{0pt}%
\pgfsys@defobject{currentmarker}{\pgfqpoint{-0.041667in}{-0.041667in}}{\pgfqpoint{0.041667in}{0.041667in}}{%
\pgfpathmoveto{\pgfqpoint{0.000000in}{-0.041667in}}%
\pgfpathcurveto{\pgfqpoint{0.011050in}{-0.041667in}}{\pgfqpoint{0.021649in}{-0.037276in}}{\pgfqpoint{0.029463in}{-0.029463in}}%
\pgfpathcurveto{\pgfqpoint{0.037276in}{-0.021649in}}{\pgfqpoint{0.041667in}{-0.011050in}}{\pgfqpoint{0.041667in}{0.000000in}}%
\pgfpathcurveto{\pgfqpoint{0.041667in}{0.011050in}}{\pgfqpoint{0.037276in}{0.021649in}}{\pgfqpoint{0.029463in}{0.029463in}}%
\pgfpathcurveto{\pgfqpoint{0.021649in}{0.037276in}}{\pgfqpoint{0.011050in}{0.041667in}}{\pgfqpoint{0.000000in}{0.041667in}}%
\pgfpathcurveto{\pgfqpoint{-0.011050in}{0.041667in}}{\pgfqpoint{-0.021649in}{0.037276in}}{\pgfqpoint{-0.029463in}{0.029463in}}%
\pgfpathcurveto{\pgfqpoint{-0.037276in}{0.021649in}}{\pgfqpoint{-0.041667in}{0.011050in}}{\pgfqpoint{-0.041667in}{0.000000in}}%
\pgfpathcurveto{\pgfqpoint{-0.041667in}{-0.011050in}}{\pgfqpoint{-0.037276in}{-0.021649in}}{\pgfqpoint{-0.029463in}{-0.029463in}}%
\pgfpathcurveto{\pgfqpoint{-0.021649in}{-0.037276in}}{\pgfqpoint{-0.011050in}{-0.041667in}}{\pgfqpoint{0.000000in}{-0.041667in}}%
\pgfpathlineto{\pgfqpoint{0.000000in}{-0.041667in}}%
\pgfpathclose%
\pgfusepath{stroke,fill}%
}%
\begin{pgfscope}%
\pgfsys@transformshift{1.550833in}{3.176651in}%
\pgfsys@useobject{currentmarker}{}%
\end{pgfscope}%
\end{pgfscope}%
\begin{pgfscope}%
\definecolor{textcolor}{rgb}{0.000000,0.000000,0.000000}%
\pgfsetstrokecolor{textcolor}%
\pgfsetfillcolor{textcolor}%
\pgftext[x=1.950833in,y=3.098874in,left,base]{\color{textcolor}\rmfamily\fontsize{16.000000}{19.200000}\selectfont Operator costs}%
\end{pgfscope}%
\begin{pgfscope}%
\pgfsetbuttcap%
\pgfsetroundjoin%
\definecolor{currentfill}{rgb}{0.000000,0.000000,0.000000}%
\pgfsetfillcolor{currentfill}%
\pgfsetlinewidth{0.803000pt}%
\definecolor{currentstroke}{rgb}{0.000000,0.000000,0.000000}%
\pgfsetstrokecolor{currentstroke}%
\pgfsetdash{}{0pt}%
\pgfsys@defobject{currentmarker}{\pgfqpoint{0.000000in}{0.000000in}}{\pgfqpoint{0.083333in}{0.000000in}}{%
\pgfpathmoveto{\pgfqpoint{0.000000in}{0.000000in}}%
\pgfpathlineto{\pgfqpoint{0.083333in}{0.000000in}}%
\pgfusepath{stroke,fill}%
}%
\begin{pgfscope}%
\pgfsys@transformshift{5.706667in}{1.139788in}%
\pgfsys@useobject{currentmarker}{}%
\end{pgfscope}%
\end{pgfscope}%
\begin{pgfscope}%
\pgfsetbuttcap%
\pgfsetroundjoin%
\definecolor{currentfill}{rgb}{0.000000,0.000000,0.000000}%
\pgfsetfillcolor{currentfill}%
\pgfsetlinewidth{0.803000pt}%
\definecolor{currentstroke}{rgb}{0.000000,0.000000,0.000000}%
\pgfsetstrokecolor{currentstroke}%
\pgfsetdash{}{0pt}%
\pgfsys@defobject{currentmarker}{\pgfqpoint{0.000000in}{0.000000in}}{\pgfqpoint{0.083333in}{0.000000in}}{%
\pgfpathmoveto{\pgfqpoint{0.000000in}{0.000000in}}%
\pgfpathlineto{\pgfqpoint{0.083333in}{0.000000in}}%
\pgfusepath{stroke,fill}%
}%
\begin{pgfscope}%
\pgfsys@transformshift{5.706667in}{1.778521in}%
\pgfsys@useobject{currentmarker}{}%
\end{pgfscope}%
\end{pgfscope}%
\begin{pgfscope}%
\pgfsetbuttcap%
\pgfsetroundjoin%
\definecolor{currentfill}{rgb}{0.000000,0.000000,0.000000}%
\pgfsetfillcolor{currentfill}%
\pgfsetlinewidth{0.803000pt}%
\definecolor{currentstroke}{rgb}{0.000000,0.000000,0.000000}%
\pgfsetstrokecolor{currentstroke}%
\pgfsetdash{}{0pt}%
\pgfsys@defobject{currentmarker}{\pgfqpoint{0.000000in}{0.000000in}}{\pgfqpoint{0.083333in}{0.000000in}}{%
\pgfpathmoveto{\pgfqpoint{0.000000in}{0.000000in}}%
\pgfpathlineto{\pgfqpoint{0.083333in}{0.000000in}}%
\pgfusepath{stroke,fill}%
}%
\begin{pgfscope}%
\pgfsys@transformshift{5.706667in}{2.417253in}%
\pgfsys@useobject{currentmarker}{}%
\end{pgfscope}%
\end{pgfscope}%
\begin{pgfscope}%
\pgfsetbuttcap%
\pgfsetroundjoin%
\definecolor{currentfill}{rgb}{0.000000,0.000000,0.000000}%
\pgfsetfillcolor{currentfill}%
\pgfsetlinewidth{0.803000pt}%
\definecolor{currentstroke}{rgb}{0.000000,0.000000,0.000000}%
\pgfsetstrokecolor{currentstroke}%
\pgfsetdash{}{0pt}%
\pgfsys@defobject{currentmarker}{\pgfqpoint{0.000000in}{0.000000in}}{\pgfqpoint{0.083333in}{0.000000in}}{%
\pgfpathmoveto{\pgfqpoint{0.000000in}{0.000000in}}%
\pgfpathlineto{\pgfqpoint{0.083333in}{0.000000in}}%
\pgfusepath{stroke,fill}%
}%
\begin{pgfscope}%
\pgfsys@transformshift{5.706667in}{3.055986in}%
\pgfsys@useobject{currentmarker}{}%
\end{pgfscope}%
\end{pgfscope}%
\begin{pgfscope}%
\pgfsetbuttcap%
\pgfsetroundjoin%
\definecolor{currentfill}{rgb}{0.000000,0.000000,0.000000}%
\pgfsetfillcolor{currentfill}%
\pgfsetlinewidth{0.803000pt}%
\definecolor{currentstroke}{rgb}{0.000000,0.000000,0.000000}%
\pgfsetstrokecolor{currentstroke}%
\pgfsetdash{}{0pt}%
\pgfsys@defobject{currentmarker}{\pgfqpoint{0.000000in}{0.000000in}}{\pgfqpoint{0.083333in}{0.000000in}}{%
\pgfpathmoveto{\pgfqpoint{0.000000in}{0.000000in}}%
\pgfpathlineto{\pgfqpoint{0.083333in}{0.000000in}}%
\pgfusepath{stroke,fill}%
}%
\begin{pgfscope}%
\pgfsys@transformshift{5.706667in}{3.694718in}%
\pgfsys@useobject{currentmarker}{}%
\end{pgfscope}%
\end{pgfscope}%
\begin{pgfscope}%
\pgfsetrectcap%
\pgfsetmiterjoin%
\pgfsetlinewidth{0.803000pt}%
\definecolor{currentstroke}{rgb}{0.000000,0.000000,0.000000}%
\pgfsetstrokecolor{currentstroke}%
\pgfsetdash{}{0pt}%
\pgfpathmoveto{\pgfqpoint{1.128611in}{0.731763in}}%
\pgfpathlineto{\pgfqpoint{1.128611in}{3.790000in}}%
\pgfusepath{stroke}%
\end{pgfscope}%
\begin{pgfscope}%
\pgfsetrectcap%
\pgfsetmiterjoin%
\pgfsetlinewidth{0.803000pt}%
\definecolor{currentstroke}{rgb}{0.000000,0.000000,0.000000}%
\pgfsetstrokecolor{currentstroke}%
\pgfsetdash{}{0pt}%
\pgfpathmoveto{\pgfqpoint{5.706667in}{0.731763in}}%
\pgfpathlineto{\pgfqpoint{5.706667in}{3.790000in}}%
\pgfusepath{stroke}%
\end{pgfscope}%
\begin{pgfscope}%
\pgfsetrectcap%
\pgfsetmiterjoin%
\pgfsetlinewidth{0.803000pt}%
\definecolor{currentstroke}{rgb}{0.000000,0.000000,0.000000}%
\pgfsetstrokecolor{currentstroke}%
\pgfsetdash{}{0pt}%
\pgfpathmoveto{\pgfqpoint{1.128611in}{0.731763in}}%
\pgfpathlineto{\pgfqpoint{5.706667in}{0.731763in}}%
\pgfusepath{stroke}%
\end{pgfscope}%
\begin{pgfscope}%
\pgfsetrectcap%
\pgfsetmiterjoin%
\pgfsetlinewidth{0.803000pt}%
\definecolor{currentstroke}{rgb}{0.000000,0.000000,0.000000}%
\pgfsetstrokecolor{currentstroke}%
\pgfsetdash{}{0pt}%
\pgfpathmoveto{\pgfqpoint{1.128611in}{3.790000in}}%
\pgfpathlineto{\pgfqpoint{5.706667in}{3.790000in}}%
\pgfusepath{stroke}%
\end{pgfscope}%
\end{pgfpicture}%
\makeatother%
\endgroup%

%% file: 06_DemandElasticity.tex
\newcommand{\DE}{{\problem}$^{(T)}$\xspace}
\newcommand{\DEpi}{{\problem}$^{(T,P)}$\xspace}
\newcommand{\Rigid}{{\problem}$^{(R)}$\xspace}
\newcommand{\RigidResolved}{{\problem}$^{(R,T)}$\xspace}

\subsection{Impact of considering service-dependent demand} \label{sec:DemandElasticity}
This section evaluates how incorporating {\elasticDemand} affects the solution quality of the {\problem}. 
We compare two model variants: one that accounts for passenger sensitivity to service frequency and one with a fixed-demand assumption. The frequency-dependent model {\problem{($\Pset$)}} restricts passenger assignment to paths where WTT does not exceed a predefined threshold, which is a function of both route and frequency. 
In contrast, the fixed-demand variant, {\problem{($\Pset^{\text{rigid}}$)}}, assumes that passenger assignment is independent of frequency, i.e., paths are available at any frequency and passenger demand is fixed by excluding alternative mode paths $\altModePath_{\od}$ for each $\od\in\ODset$ unless no feasible PTN path exists. 
We denote these models as \DE and \Rigid, respectively. 

To give a fair comparison between the two models, we apply a post-processing step to \Rigid. In this step, we estimate the number of passengers that would decide to travel in the network using the same threshold-based assumptions on the WTT as \DE. The resulting assignment is then used to compute the passenger costs, operator costs, and revenue. We refer to this solution as \RigidResolved. Furthermore, to ensure we compare using the same resource constraints, we first solve \Rigid without budget constraints to determine the required operating costs. This budget is then imposed on \DE, ensuring that both models are compared under the same operational conditions and consistent passenger behavior. Since \DE optimizes the trade-off between operator and passenger costs, we solve it using both the objective defined in \eqref{obj} and an objective using the passenger component of \eqref{obj}, which we refer to as \DEpi.

\begin{figure}[htbp]
  \centering
  \begin{tabular}{cc}
    \begin{subfigure}[b]{0.5\textwidth}
      \centering
      \resizebox{\linewidth}{!}{\input{figures/DemandPlots/TightLayout/DE_DemandCapturedPct.pgf}}
      \caption{Demand captured}
            \label{fig:DE_demand}
    \end{subfigure}  & 
        \begin{subfigure}[b]{0.5\textwidth} 
      \centering
      \resizebox{\linewidth}{!}{\input{figures/DemandPlots/TightLayout/DE_avgWTT_rescaled.pgf}} 
      \caption{Weighted travel time ($\lambda$-normalized)}
      \label{fig:DE_WTT}
    \end{subfigure} \\ [-1.0ex]
          \begin{subfigure}[b]{0.5\textwidth}
      \centering
      \resizebox{\linewidth}{!}{\input{figures/DemandPlots/TightLayout/DE_UB_rescaled.pgf}}
      \caption{Objective value ($\lambda$-normalized)}
            \label{fig:DE_UB}
    \end{subfigure}& 
    \begin{subfigure}[b]{0.5\textwidth}
      \centering
      \resizebox{\linewidth}{!}{\input{figures/DemandPlots/TightLayout/DE_PI_rescaled.pgf}}
      \caption{Passenger cost ($\lambda$-normalized)} 
                  \label{fig:DE_PI}
    \end{subfigure} 
  \end{tabular}
  \caption{Results over different values of $\lambda$ for instances solved to optimality}
  \label{fig:DE_lambda}
\end{figure}
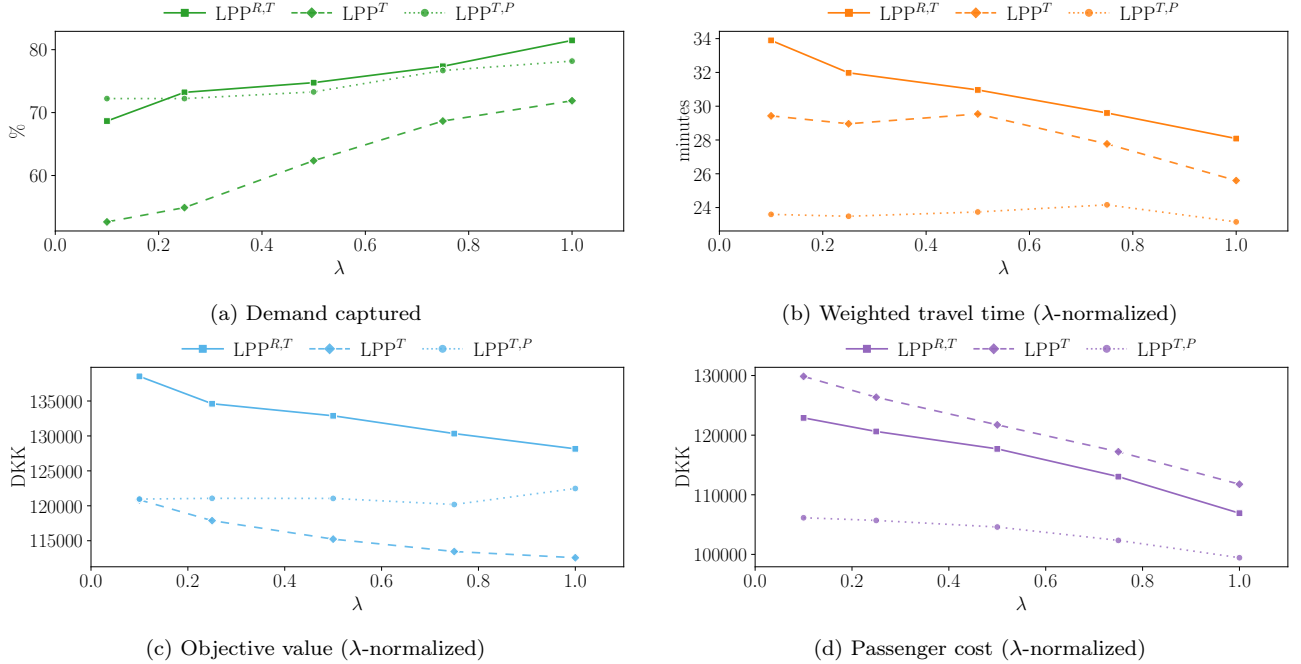

Given an OD matrix and a set of candidate paths, we evaluate solutions based on their objective components: passenger cost (including both routed and lost demand), operator cost, the WTT for served passengers, and total ridership captured by the PTN. 
\autoref{fig:DE_lambda} presents these key metrics for each model across $\lambda\in \{0.1, 0.25, 0.5, 0.75, 1.0\}$.
Passenger costs have been normalized to a fixed $\lambda$ value and the objective has been recalculated for \DEpi to include operator costs for a consistent comparison.
The results are averaged over a subset of seven instances -- including \athens, \grid, and \sioux instances -- that can be solved to optimality by all models. These instances feature the smallest line pool and between 100 and 300 ODs.  

\autoref{fig:DE_demand} clearly shows that \RigidResolved captures far less demand than the networks are designed for. 
Despite the network's capacity to serve all passengers, only an average of 69--81\% of demand is captured.
This gap suggests that under fixed-demand assumptions, the supply and actual demand are not aligned. 
The model \DEpi captures slightly less demand, except at the lowest $\lambda$ value. 
However, while the demand differs by no more than 4\% percentage points between \RigidResolved and \DEpi, the passenger service in terms of weighted travel time is consistently better for \DEpi, as seen in \autoref{fig:DE_WTT}. 
On average, WTT is reduced by 5 to 10 minutes  (from $\lambda=1.0$ to $\lambda=0.1$), corresponding to an overall improvement in the objective of 4--13\%, as illustrated in \autoref{fig:DE_UB}. 
Both models use exactly the same resources in terms of operating costs. 
The results suggest that, even with identical resources, passenger-centric planning can deliver substantially improved service levels to those realistically willing to travel by allocating resources more efficiently.

If we consider the \DE model, we observe that it serves the least demand on average yet achieves the lowest total objective value. 
This reflects the system-optimal focus of the objective: in scenarios where routing additional passengers results in disproportionately high costs, the model opts to serve fewer passengers. 
As shown in \autoref{fig:DE_UB}, this leads to total cost reductions of 12--13\% relative to \RigidResolved despite lower ridership. 
This outcome is partly driven by the cost parameters associated with lost demand. Specifically, the penalty for not serving passengers is relatively modest compared to the cost of routing them, which reduces the incentive to maximize ridership. If the penalty for lost demand were increased, the ridership achieved by \DE would more closely resemble that of \RigidResolved. Overall, these results highlight the importance of accounting for \elasticDemand. This improvement is most noticeable when compared against a benchmark model that places less emphasis on passenger service, i.e., lower values of $\lambda$. 

This analysis is limited to instances that can be solved to optimality. 
To broaden the scope, \autoref{sec:DE_appendix} presents results that include instances solved within a 10\% optimality gap, which allows for including slightly larger instances. The results remain largely consistent with those shown in \autoref{fig:DE_lambda}. Notably, on average, \DEpi has the largest ridership for all values of $\lambda$ except at $\lambda = 1$, while at the same time allowing lower travel time. This improvement is partly due to the inclusion of larger instances with more ODs, which makes it more favorable to offer more services. 

A closer examination of the individual solutions reveals that the line concepts produced using \Rigid tend to open a larger number of lines. However, once passenger flows are re-evaluated with more realistic assumptions -- specifically, imposing an OD-specific threshold on the acceptable travel time -- many lines fail to attract sufficient ridership to justify their operational cost. 
In contrast, the line plans generated by \DE and \DEpi emphasize operating fewer lines at higher frequencies.  
This approach prioritizes providing a better service where demand is concentrated, capturing more passengers on fewer, well-utilized lines rather than spreading service thinly over the network. 

Finally, we note that our approach assumes that demand responds to service levels and that travel times can be reasonably estimated using a detailed set of frequencies. However, not all frequencies are equally relevant. A sensitivity analysis in \autoref{sec:FreqAnalysis} shows that expanding the number of available frequencies from four options to all integer values between two and 20 minutes improves the objective by 5\% and increases captured demand by 19\%. The benefits of introducing more finely spaced headways within the same range diminish beyond 10 frequencies, yielding only a 1.3\% improvement in the objective and a 2\% increase in demand. This suggests that a well-chosen, limited set of practical headways may suffice for modeling purposes.

%% file: figures/DemandPlots/TightLayout/DE_DemandCapturedPct.pgf
\begingroup%
\makeatletter%
\begin{pgfpicture}%
\pgfpathrectangle{\pgfpointorigin}{\pgfqpoint{7.775300in}{3.651310in}}%
\pgfusepath{use as bounding box, clip}%
\begin{pgfscope}%
\pgfsetbuttcap%
\pgfsetmiterjoin%
\definecolor{currentfill}{rgb}{1.000000,1.000000,1.000000}%
\pgfsetfillcolor{currentfill}%
\pgfsetlinewidth{0.000000pt}%
\definecolor{currentstroke}{rgb}{1.000000,1.000000,1.000000}%
\pgfsetstrokecolor{currentstroke}%
\pgfsetdash{}{0pt}%
\pgfpathmoveto{\pgfqpoint{0.000000in}{0.000000in}}%
\pgfpathlineto{\pgfqpoint{7.775300in}{0.000000in}}%
\pgfpathlineto{\pgfqpoint{7.775300in}{3.651310in}}%
\pgfpathlineto{\pgfqpoint{0.000000in}{3.651310in}}%
\pgfpathlineto{\pgfqpoint{0.000000in}{0.000000in}}%
\pgfpathclose%
\pgfusepath{fill}%
\end{pgfscope}%
\begin{pgfscope}%
\pgfsetbuttcap%
\pgfsetmiterjoin%
\definecolor{currentfill}{rgb}{1.000000,1.000000,1.000000}%
\pgfsetfillcolor{currentfill}%
\pgfsetlinewidth{0.000000pt}%
\definecolor{currentstroke}{rgb}{0.000000,0.000000,0.000000}%
\pgfsetstrokecolor{currentstroke}%
\pgfsetstrokeopacity{0.000000}%
\pgfsetdash{}{0pt}%
\pgfpathmoveto{\pgfqpoint{0.686263in}{0.679475in}}%
\pgfpathlineto{\pgfqpoint{7.675300in}{0.679475in}}%
\pgfpathlineto{\pgfqpoint{7.675300in}{3.135897in}}%
\pgfpathlineto{\pgfqpoint{0.686263in}{3.135897in}}%
\pgfpathlineto{\pgfqpoint{0.686263in}{0.679475in}}%
\pgfpathclose%
\pgfusepath{fill}%
\end{pgfscope}%
\begin{pgfscope}%
\pgfsetbuttcap%
\pgfsetroundjoin%
\definecolor{currentfill}{rgb}{0.000000,0.000000,0.000000}%
\pgfsetfillcolor{currentfill}%
\pgfsetlinewidth{0.803000pt}%
\definecolor{currentstroke}{rgb}{0.000000,0.000000,0.000000}%
\pgfsetstrokecolor{currentstroke}%
\pgfsetdash{}{0pt}%
\pgfsys@defobject{currentmarker}{\pgfqpoint{0.000000in}{-0.048611in}}{\pgfqpoint{0.000000in}{0.000000in}}{%
\pgfpathmoveto{\pgfqpoint{0.000000in}{0.000000in}}%
\pgfpathlineto{\pgfqpoint{0.000000in}{-0.048611in}}%
\pgfusepath{stroke,fill}%
}%
\begin{pgfscope}%
\pgfsys@transformshift{0.686263in}{0.679475in}%
\pgfsys@useobject{currentmarker}{}%
\end{pgfscope}%
\end{pgfscope}%
\begin{pgfscope}%
\definecolor{textcolor}{rgb}{0.000000,0.000000,0.000000}%
\pgfsetstrokecolor{textcolor}%
\pgfsetfillcolor{textcolor}%
\pgftext[x=0.686263in,y=0.582253in,,top]{\color{textcolor}\rmfamily\fontsize{16.000000}{19.200000}\selectfont \(\displaystyle {0.0}\)}%
\end{pgfscope}%
\begin{pgfscope}%
\pgfsetbuttcap%
\pgfsetroundjoin%
\definecolor{currentfill}{rgb}{0.000000,0.000000,0.000000}%
\pgfsetfillcolor{currentfill}%
\pgfsetlinewidth{0.803000pt}%
\definecolor{currentstroke}{rgb}{0.000000,0.000000,0.000000}%
\pgfsetstrokecolor{currentstroke}%
\pgfsetdash{}{0pt}%
\pgfsys@defobject{currentmarker}{\pgfqpoint{0.000000in}{-0.048611in}}{\pgfqpoint{0.000000in}{0.000000in}}{%
\pgfpathmoveto{\pgfqpoint{0.000000in}{0.000000in}}%
\pgfpathlineto{\pgfqpoint{0.000000in}{-0.048611in}}%
\pgfusepath{stroke,fill}%
}%
\begin{pgfscope}%
\pgfsys@transformshift{1.956997in}{0.679475in}%
\pgfsys@useobject{currentmarker}{}%
\end{pgfscope}%
\end{pgfscope}%
\begin{pgfscope}%
\definecolor{textcolor}{rgb}{0.000000,0.000000,0.000000}%
\pgfsetstrokecolor{textcolor}%
\pgfsetfillcolor{textcolor}%
\pgftext[x=1.956997in,y=0.582253in,,top]{\color{textcolor}\rmfamily\fontsize{16.000000}{19.200000}\selectfont \(\displaystyle {0.2}\)}%
\end{pgfscope}%
\begin{pgfscope}%
\pgfsetbuttcap%
\pgfsetroundjoin%
\definecolor{currentfill}{rgb}{0.000000,0.000000,0.000000}%
\pgfsetfillcolor{currentfill}%
\pgfsetlinewidth{0.803000pt}%
\definecolor{currentstroke}{rgb}{0.000000,0.000000,0.000000}%
\pgfsetstrokecolor{currentstroke}%
\pgfsetdash{}{0pt}%
\pgfsys@defobject{currentmarker}{\pgfqpoint{0.000000in}{-0.048611in}}{\pgfqpoint{0.000000in}{0.000000in}}{%
\pgfpathmoveto{\pgfqpoint{0.000000in}{0.000000in}}%
\pgfpathlineto{\pgfqpoint{0.000000in}{-0.048611in}}%
\pgfusepath{stroke,fill}%
}%
\begin{pgfscope}%
\pgfsys@transformshift{3.227731in}{0.679475in}%
\pgfsys@useobject{currentmarker}{}%
\end{pgfscope}%
\end{pgfscope}%
\begin{pgfscope}%
\definecolor{textcolor}{rgb}{0.000000,0.000000,0.000000}%
\pgfsetstrokecolor{textcolor}%
\pgfsetfillcolor{textcolor}%
\pgftext[x=3.227731in,y=0.582253in,,top]{\color{textcolor}\rmfamily\fontsize{16.000000}{19.200000}\selectfont \(\displaystyle {0.4}\)}%
\end{pgfscope}%
\begin{pgfscope}%
\pgfsetbuttcap%
\pgfsetroundjoin%
\definecolor{currentfill}{rgb}{0.000000,0.000000,0.000000}%
\pgfsetfillcolor{currentfill}%
\pgfsetlinewidth{0.803000pt}%
\definecolor{currentstroke}{rgb}{0.000000,0.000000,0.000000}%
\pgfsetstrokecolor{currentstroke}%
\pgfsetdash{}{0pt}%
\pgfsys@defobject{currentmarker}{\pgfqpoint{0.000000in}{-0.048611in}}{\pgfqpoint{0.000000in}{0.000000in}}{%
\pgfpathmoveto{\pgfqpoint{0.000000in}{0.000000in}}%
\pgfpathlineto{\pgfqpoint{0.000000in}{-0.048611in}}%
\pgfusepath{stroke,fill}%
}%
\begin{pgfscope}%
\pgfsys@transformshift{4.498465in}{0.679475in}%
\pgfsys@useobject{currentmarker}{}%
\end{pgfscope}%
\end{pgfscope}%
\begin{pgfscope}%
\definecolor{textcolor}{rgb}{0.000000,0.000000,0.000000}%
\pgfsetstrokecolor{textcolor}%
\pgfsetfillcolor{textcolor}%
\pgftext[x=4.498465in,y=0.582253in,,top]{\color{textcolor}\rmfamily\fontsize{16.000000}{19.200000}\selectfont \(\displaystyle {0.6}\)}%
\end{pgfscope}%
\begin{pgfscope}%
\pgfsetbuttcap%
\pgfsetroundjoin%
\definecolor{currentfill}{rgb}{0.000000,0.000000,0.000000}%
\pgfsetfillcolor{currentfill}%
\pgfsetlinewidth{0.803000pt}%
\definecolor{currentstroke}{rgb}{0.000000,0.000000,0.000000}%
\pgfsetstrokecolor{currentstroke}%
\pgfsetdash{}{0pt}%
\pgfsys@defobject{currentmarker}{\pgfqpoint{0.000000in}{-0.048611in}}{\pgfqpoint{0.000000in}{0.000000in}}{%
\pgfpathmoveto{\pgfqpoint{0.000000in}{0.000000in}}%
\pgfpathlineto{\pgfqpoint{0.000000in}{-0.048611in}}%
\pgfusepath{stroke,fill}%
}%
\begin{pgfscope}%
\pgfsys@transformshift{5.769199in}{0.679475in}%
\pgfsys@useobject{currentmarker}{}%
\end{pgfscope}%
\end{pgfscope}%
\begin{pgfscope}%
\definecolor{textcolor}{rgb}{0.000000,0.000000,0.000000}%
\pgfsetstrokecolor{textcolor}%
\pgfsetfillcolor{textcolor}%
\pgftext[x=5.769199in,y=0.582253in,,top]{\color{textcolor}\rmfamily\fontsize{16.000000}{19.200000}\selectfont \(\displaystyle {0.8}\)}%
\end{pgfscope}%
\begin{pgfscope}%
\pgfsetbuttcap%
\pgfsetroundjoin%
\definecolor{currentfill}{rgb}{0.000000,0.000000,0.000000}%
\pgfsetfillcolor{currentfill}%
\pgfsetlinewidth{0.803000pt}%
\definecolor{currentstroke}{rgb}{0.000000,0.000000,0.000000}%
\pgfsetstrokecolor{currentstroke}%
\pgfsetdash{}{0pt}%
\pgfsys@defobject{currentmarker}{\pgfqpoint{0.000000in}{-0.048611in}}{\pgfqpoint{0.000000in}{0.000000in}}{%
\pgfpathmoveto{\pgfqpoint{0.000000in}{0.000000in}}%
\pgfpathlineto{\pgfqpoint{0.000000in}{-0.048611in}}%
\pgfusepath{stroke,fill}%
}%
\begin{pgfscope}%
\pgfsys@transformshift{7.039933in}{0.679475in}%
\pgfsys@useobject{currentmarker}{}%
\end{pgfscope}%
\end{pgfscope}%
\begin{pgfscope}%
\definecolor{textcolor}{rgb}{0.000000,0.000000,0.000000}%
\pgfsetstrokecolor{textcolor}%
\pgfsetfillcolor{textcolor}%
\pgftext[x=7.039933in,y=0.582253in,,top]{\color{textcolor}\rmfamily\fontsize{16.000000}{19.200000}\selectfont \(\displaystyle {1.0}\)}%
\end{pgfscope}%
\begin{pgfscope}%
\definecolor{textcolor}{rgb}{0.000000,0.000000,0.000000}%
\pgfsetstrokecolor{textcolor}%
\pgfsetfillcolor{textcolor}%
\pgftext[x=4.180781in,y=0.313349in,,top]{\color{textcolor}\rmfamily\fontsize{18.000000}{21.600000}\selectfont \(\displaystyle \lambda\)}%
\end{pgfscope}%
\begin{pgfscope}%
\pgfsetbuttcap%
\pgfsetroundjoin%
\definecolor{currentfill}{rgb}{0.000000,0.000000,0.000000}%
\pgfsetfillcolor{currentfill}%
\pgfsetlinewidth{0.803000pt}%
\definecolor{currentstroke}{rgb}{0.000000,0.000000,0.000000}%
\pgfsetstrokecolor{currentstroke}%
\pgfsetdash{}{0pt}%
\pgfsys@defobject{currentmarker}{\pgfqpoint{-0.048611in}{0.000000in}}{\pgfqpoint{-0.000000in}{0.000000in}}{%
\pgfpathmoveto{\pgfqpoint{-0.000000in}{0.000000in}}%
\pgfpathlineto{\pgfqpoint{-0.048611in}{0.000000in}}%
\pgfusepath{stroke,fill}%
}%
\begin{pgfscope}%
\pgfsys@transformshift{0.686263in}{1.361851in}%
\pgfsys@useobject{currentmarker}{}%
\end{pgfscope}%
\end{pgfscope}%
\begin{pgfscope}%
\definecolor{textcolor}{rgb}{0.000000,0.000000,0.000000}%
\pgfsetstrokecolor{textcolor}%
\pgfsetfillcolor{textcolor}%
\pgftext[x=0.368904in, y=1.278518in, left, base]{\color{textcolor}\rmfamily\fontsize{16.000000}{19.200000}\selectfont \(\displaystyle {60}\)}%
\end{pgfscope}%
\begin{pgfscope}%
\pgfsetbuttcap%
\pgfsetroundjoin%
\definecolor{currentfill}{rgb}{0.000000,0.000000,0.000000}%
\pgfsetfillcolor{currentfill}%
\pgfsetlinewidth{0.803000pt}%
\definecolor{currentstroke}{rgb}{0.000000,0.000000,0.000000}%
\pgfsetstrokecolor{currentstroke}%
\pgfsetdash{}{0pt}%
\pgfsys@defobject{currentmarker}{\pgfqpoint{-0.048611in}{0.000000in}}{\pgfqpoint{-0.000000in}{0.000000in}}{%
\pgfpathmoveto{\pgfqpoint{-0.000000in}{0.000000in}}%
\pgfpathlineto{\pgfqpoint{-0.048611in}{0.000000in}}%
\pgfusepath{stroke,fill}%
}%
\begin{pgfscope}%
\pgfsys@transformshift{0.686263in}{2.136085in}%
\pgfsys@useobject{currentmarker}{}%
\end{pgfscope}%
\end{pgfscope}%
\begin{pgfscope}%
\definecolor{textcolor}{rgb}{0.000000,0.000000,0.000000}%
\pgfsetstrokecolor{textcolor}%
\pgfsetfillcolor{textcolor}%
\pgftext[x=0.368904in, y=2.052751in, left, base]{\color{textcolor}\rmfamily\fontsize{16.000000}{19.200000}\selectfont \(\displaystyle {70}\)}%
\end{pgfscope}%
\begin{pgfscope}%
\pgfsetbuttcap%
\pgfsetroundjoin%
\definecolor{currentfill}{rgb}{0.000000,0.000000,0.000000}%
\pgfsetfillcolor{currentfill}%
\pgfsetlinewidth{0.803000pt}%
\definecolor{currentstroke}{rgb}{0.000000,0.000000,0.000000}%
\pgfsetstrokecolor{currentstroke}%
\pgfsetdash{}{0pt}%
\pgfsys@defobject{currentmarker}{\pgfqpoint{-0.048611in}{0.000000in}}{\pgfqpoint{-0.000000in}{0.000000in}}{%
\pgfpathmoveto{\pgfqpoint{-0.000000in}{0.000000in}}%
\pgfpathlineto{\pgfqpoint{-0.048611in}{0.000000in}}%
\pgfusepath{stroke,fill}%
}%
\begin{pgfscope}%
\pgfsys@transformshift{0.686263in}{2.910318in}%
\pgfsys@useobject{currentmarker}{}%
\end{pgfscope}%
\end{pgfscope}%
\begin{pgfscope}%
\definecolor{textcolor}{rgb}{0.000000,0.000000,0.000000}%
\pgfsetstrokecolor{textcolor}%
\pgfsetfillcolor{textcolor}%
\pgftext[x=0.368904in, y=2.826985in, left, base]{\color{textcolor}\rmfamily\fontsize{16.000000}{19.200000}\selectfont \(\displaystyle {80}\)}%
\end{pgfscope}%
\begin{pgfscope}%
\definecolor{textcolor}{rgb}{0.000000,0.000000,0.000000}%
\pgfsetstrokecolor{textcolor}%
\pgfsetfillcolor{textcolor}%
\pgftext[x=0.313349in,y=1.907686in,,bottom,rotate=90.000000]{\color{textcolor}\rmfamily\fontsize{18.000000}{21.600000}\selectfont \%}%
\end{pgfscope}%
\begin{pgfscope}%
\pgfpathrectangle{\pgfqpoint{0.686263in}{0.679475in}}{\pgfqpoint{6.989037in}{2.456422in}}%
\pgfusepath{clip}%
\pgfsetrectcap%
\pgfsetroundjoin%
\pgfsetlinewidth{1.505625pt}%
\definecolor{currentstroke}{rgb}{0.172549,0.627451,0.172549}%
\pgfsetstrokecolor{currentstroke}%
\pgfsetdash{}{0pt}%
\pgfpathmoveto{\pgfqpoint{1.321630in}{2.032116in}}%
\pgfpathlineto{\pgfqpoint{2.274680in}{2.384946in}}%
\pgfpathlineto{\pgfqpoint{3.863098in}{2.504399in}}%
\pgfpathlineto{\pgfqpoint{5.451515in}{2.705699in}}%
\pgfpathlineto{\pgfqpoint{7.039933in}{3.024241in}}%
\pgfusepath{stroke}%
\end{pgfscope}%
\begin{pgfscope}%
\pgfpathrectangle{\pgfqpoint{0.686263in}{0.679475in}}{\pgfqpoint{6.989037in}{2.456422in}}%
\pgfusepath{clip}%
\pgfsetbuttcap%
\pgfsetmiterjoin%
\definecolor{currentfill}{rgb}{0.172549,0.627451,0.172549}%
\pgfsetfillcolor{currentfill}%
\pgfsetlinewidth{0.752812pt}%
\definecolor{currentstroke}{rgb}{1.000000,1.000000,1.000000}%
\pgfsetstrokecolor{currentstroke}%
\pgfsetdash{}{0pt}%
\pgfsys@defobject{currentmarker}{\pgfqpoint{-0.041667in}{-0.041667in}}{\pgfqpoint{0.041667in}{0.041667in}}{%
\pgfpathmoveto{\pgfqpoint{-0.041667in}{-0.041667in}}%
\pgfpathlineto{\pgfqpoint{0.041667in}{-0.041667in}}%
\pgfpathlineto{\pgfqpoint{0.041667in}{0.041667in}}%
\pgfpathlineto{\pgfqpoint{-0.041667in}{0.041667in}}%
\pgfpathlineto{\pgfqpoint{-0.041667in}{-0.041667in}}%
\pgfpathclose%
\pgfusepath{stroke,fill}%
}%
\begin{pgfscope}%
\pgfsys@transformshift{1.321630in}{2.032116in}%
\pgfsys@useobject{currentmarker}{}%
\end{pgfscope}%
\begin{pgfscope}%
\pgfsys@transformshift{2.274680in}{2.384946in}%
\pgfsys@useobject{currentmarker}{}%
\end{pgfscope}%
\begin{pgfscope}%
\pgfsys@transformshift{3.863098in}{2.504399in}%
\pgfsys@useobject{currentmarker}{}%
\end{pgfscope}%
\begin{pgfscope}%
\pgfsys@transformshift{5.451515in}{2.705699in}%
\pgfsys@useobject{currentmarker}{}%
\end{pgfscope}%
\begin{pgfscope}%
\pgfsys@transformshift{7.039933in}{3.024241in}%
\pgfsys@useobject{currentmarker}{}%
\end{pgfscope}%
\end{pgfscope}%
\begin{pgfscope}%
\pgfpathrectangle{\pgfqpoint{0.686263in}{0.679475in}}{\pgfqpoint{6.989037in}{2.456422in}}%
\pgfusepath{clip}%
\pgfsetbuttcap%
\pgfsetroundjoin%
\pgfsetlinewidth{1.505625pt}%
\definecolor{currentstroke}{rgb}{0.172549,0.627451,0.172549}%
\pgfsetstrokecolor{currentstroke}%
\pgfsetstrokeopacity{0.900000}%
\pgfsetdash{{7.500000pt}{7.500000pt}}{0.000000pt}%
\pgfpathmoveto{\pgfqpoint{1.321630in}{0.791131in}}%
\pgfpathlineto{\pgfqpoint{2.274680in}{0.966992in}}%
\pgfpathlineto{\pgfqpoint{3.863098in}{1.544349in}}%
\pgfpathlineto{\pgfqpoint{5.451515in}{2.032116in}}%
\pgfpathlineto{\pgfqpoint{7.039933in}{2.282083in}}%
\pgfusepath{stroke}%
\end{pgfscope}%
\begin{pgfscope}%
\pgfpathrectangle{\pgfqpoint{0.686263in}{0.679475in}}{\pgfqpoint{6.989037in}{2.456422in}}%
\pgfusepath{clip}%
\pgfsetbuttcap%
\pgfsetmiterjoin%
\definecolor{currentfill}{rgb}{0.172549,0.627451,0.172549}%
\pgfsetfillcolor{currentfill}%
\pgfsetfillopacity{0.900000}%
\pgfsetlinewidth{0.752812pt}%
\definecolor{currentstroke}{rgb}{1.000000,1.000000,1.000000}%
\pgfsetstrokecolor{currentstroke}%
\pgfsetdash{}{0pt}%
\pgfsys@defobject{currentmarker}{\pgfqpoint{-0.058926in}{-0.058926in}}{\pgfqpoint{0.058926in}{0.058926in}}{%
\pgfpathmoveto{\pgfqpoint{0.000000in}{-0.058926in}}%
\pgfpathlineto{\pgfqpoint{0.058926in}{0.000000in}}%
\pgfpathlineto{\pgfqpoint{0.000000in}{0.058926in}}%
\pgfpathlineto{\pgfqpoint{-0.058926in}{0.000000in}}%
\pgfpathlineto{\pgfqpoint{0.000000in}{-0.058926in}}%
\pgfpathclose%
\pgfusepath{stroke,fill}%
}%
\begin{pgfscope}%
\pgfsys@transformshift{1.321630in}{0.791131in}%
\pgfsys@useobject{currentmarker}{}%
\end{pgfscope}%
\begin{pgfscope}%
\pgfsys@transformshift{2.274680in}{0.966992in}%
\pgfsys@useobject{currentmarker}{}%
\end{pgfscope}%
\begin{pgfscope}%
\pgfsys@transformshift{3.863098in}{1.544349in}%
\pgfsys@useobject{currentmarker}{}%
\end{pgfscope}%
\begin{pgfscope}%
\pgfsys@transformshift{5.451515in}{2.032116in}%
\pgfsys@useobject{currentmarker}{}%
\end{pgfscope}%
\begin{pgfscope}%
\pgfsys@transformshift{7.039933in}{2.282083in}%
\pgfsys@useobject{currentmarker}{}%
\end{pgfscope}%
\end{pgfscope}%
\begin{pgfscope}%
\pgfpathrectangle{\pgfqpoint{0.686263in}{0.679475in}}{\pgfqpoint{6.989037in}{2.456422in}}%
\pgfusepath{clip}%
\pgfsetbuttcap%
\pgfsetroundjoin%
\pgfsetlinewidth{1.505625pt}%
\definecolor{currentstroke}{rgb}{0.172549,0.627451,0.172549}%
\pgfsetstrokecolor{currentstroke}%
\pgfsetstrokeopacity{0.800000}%
\pgfsetdash{{1.500000pt}{4.500000pt}}{0.000000pt}%
\pgfpathmoveto{\pgfqpoint{1.321630in}{2.308628in}}%
\pgfpathlineto{\pgfqpoint{2.274680in}{2.308628in}}%
\pgfpathlineto{\pgfqpoint{3.863098in}{2.390476in}}%
\pgfpathlineto{\pgfqpoint{5.451515in}{2.653715in}}%
\pgfpathlineto{\pgfqpoint{7.039933in}{2.769850in}}%
\pgfusepath{stroke}%
\end{pgfscope}%
\begin{pgfscope}%
\pgfpathrectangle{\pgfqpoint{0.686263in}{0.679475in}}{\pgfqpoint{6.989037in}{2.456422in}}%
\pgfusepath{clip}%
\pgfsetbuttcap%
\pgfsetroundjoin%
\definecolor{currentfill}{rgb}{0.172549,0.627451,0.172549}%
\pgfsetfillcolor{currentfill}%
\pgfsetfillopacity{0.800000}%
\pgfsetlinewidth{0.752812pt}%
\definecolor{currentstroke}{rgb}{1.000000,1.000000,1.000000}%
\pgfsetstrokecolor{currentstroke}%
\pgfsetdash{}{0pt}%
\pgfsys@defobject{currentmarker}{\pgfqpoint{-0.041667in}{-0.041667in}}{\pgfqpoint{0.041667in}{0.041667in}}{%
\pgfpathmoveto{\pgfqpoint{0.000000in}{-0.041667in}}%
\pgfpathcurveto{\pgfqpoint{0.011050in}{-0.041667in}}{\pgfqpoint{0.021649in}{-0.037276in}}{\pgfqpoint{0.029463in}{-0.029463in}}%
\pgfpathcurveto{\pgfqpoint{0.037276in}{-0.021649in}}{\pgfqpoint{0.041667in}{-0.011050in}}{\pgfqpoint{0.041667in}{0.000000in}}%
\pgfpathcurveto{\pgfqpoint{0.041667in}{0.011050in}}{\pgfqpoint{0.037276in}{0.021649in}}{\pgfqpoint{0.029463in}{0.029463in}}%
\pgfpathcurveto{\pgfqpoint{0.021649in}{0.037276in}}{\pgfqpoint{0.011050in}{0.041667in}}{\pgfqpoint{0.000000in}{0.041667in}}%
\pgfpathcurveto{\pgfqpoint{-0.011050in}{0.041667in}}{\pgfqpoint{-0.021649in}{0.037276in}}{\pgfqpoint{-0.029463in}{0.029463in}}%
\pgfpathcurveto{\pgfqpoint{-0.037276in}{0.021649in}}{\pgfqpoint{-0.041667in}{0.011050in}}{\pgfqpoint{-0.041667in}{0.000000in}}%
\pgfpathcurveto{\pgfqpoint{-0.041667in}{-0.011050in}}{\pgfqpoint{-0.037276in}{-0.021649in}}{\pgfqpoint{-0.029463in}{-0.029463in}}%
\pgfpathcurveto{\pgfqpoint{-0.021649in}{-0.037276in}}{\pgfqpoint{-0.011050in}{-0.041667in}}{\pgfqpoint{0.000000in}{-0.041667in}}%
\pgfpathlineto{\pgfqpoint{0.000000in}{-0.041667in}}%
\pgfpathclose%
\pgfusepath{stroke,fill}%
}%
\begin{pgfscope}%
\pgfsys@transformshift{1.321630in}{2.308628in}%
\pgfsys@useobject{currentmarker}{}%
\end{pgfscope}%
\begin{pgfscope}%
\pgfsys@transformshift{2.274680in}{2.308628in}%
\pgfsys@useobject{currentmarker}{}%
\end{pgfscope}%
\begin{pgfscope}%
\pgfsys@transformshift{3.863098in}{2.390476in}%
\pgfsys@useobject{currentmarker}{}%
\end{pgfscope}%
\begin{pgfscope}%
\pgfsys@transformshift{5.451515in}{2.653715in}%
\pgfsys@useobject{currentmarker}{}%
\end{pgfscope}%
\begin{pgfscope}%
\pgfsys@transformshift{7.039933in}{2.769850in}%
\pgfsys@useobject{currentmarker}{}%
\end{pgfscope}%
\end{pgfscope}%
\begin{pgfscope}%
\pgfsetrectcap%
\pgfsetmiterjoin%
\pgfsetlinewidth{0.803000pt}%
\definecolor{currentstroke}{rgb}{0.000000,0.000000,0.000000}%
\pgfsetstrokecolor{currentstroke}%
\pgfsetdash{}{0pt}%
\pgfpathmoveto{\pgfqpoint{0.686263in}{0.679475in}}%
\pgfpathlineto{\pgfqpoint{0.686263in}{3.135897in}}%
\pgfusepath{stroke}%
\end{pgfscope}%
\begin{pgfscope}%
\pgfsetrectcap%
\pgfsetmiterjoin%
\pgfsetlinewidth{0.803000pt}%
\definecolor{currentstroke}{rgb}{0.000000,0.000000,0.000000}%
\pgfsetstrokecolor{currentstroke}%
\pgfsetdash{}{0pt}%
\pgfpathmoveto{\pgfqpoint{7.675300in}{0.679475in}}%
\pgfpathlineto{\pgfqpoint{7.675300in}{3.135897in}}%
\pgfusepath{stroke}%
\end{pgfscope}%
\begin{pgfscope}%
\pgfsetrectcap%
\pgfsetmiterjoin%
\pgfsetlinewidth{0.803000pt}%
\definecolor{currentstroke}{rgb}{0.000000,0.000000,0.000000}%
\pgfsetstrokecolor{currentstroke}%
\pgfsetdash{}{0pt}%
\pgfpathmoveto{\pgfqpoint{0.686263in}{0.679475in}}%
\pgfpathlineto{\pgfqpoint{7.675300in}{0.679475in}}%
\pgfusepath{stroke}%
\end{pgfscope}%
\begin{pgfscope}%
\pgfsetrectcap%
\pgfsetmiterjoin%
\pgfsetlinewidth{0.803000pt}%
\definecolor{currentstroke}{rgb}{0.000000,0.000000,0.000000}%
\pgfsetstrokecolor{currentstroke}%
\pgfsetdash{}{0pt}%
\pgfpathmoveto{\pgfqpoint{0.686263in}{3.135897in}}%
\pgfpathlineto{\pgfqpoint{7.675300in}{3.135897in}}%
\pgfusepath{stroke}%
\end{pgfscope}%
\begin{pgfscope}%
\pgfsetrectcap%
\pgfsetroundjoin%
\pgfsetlinewidth{1.505625pt}%
\definecolor{currentstroke}{rgb}{0.172549,0.627451,0.172549}%
\pgfsetstrokecolor{currentstroke}%
\pgfsetdash{}{0pt}%
\pgfpathmoveto{\pgfqpoint{2.028215in}{3.363810in}}%
\pgfpathlineto{\pgfqpoint{2.278215in}{3.363810in}}%
\pgfpathlineto{\pgfqpoint{2.528215in}{3.363810in}}%
\pgfusepath{stroke}%
\end{pgfscope}%
\begin{pgfscope}%
\pgfsetbuttcap%
\pgfsetmiterjoin%
\definecolor{currentfill}{rgb}{0.172549,0.627451,0.172549}%
\pgfsetfillcolor{currentfill}%
\pgfsetlinewidth{1.003750pt}%
\definecolor{currentstroke}{rgb}{0.172549,0.627451,0.172549}%
\pgfsetstrokecolor{currentstroke}%
\pgfsetdash{}{0pt}%
\pgfsys@defobject{currentmarker}{\pgfqpoint{-0.041667in}{-0.041667in}}{\pgfqpoint{0.041667in}{0.041667in}}{%
\pgfpathmoveto{\pgfqpoint{-0.041667in}{-0.041667in}}%
\pgfpathlineto{\pgfqpoint{0.041667in}{-0.041667in}}%
\pgfpathlineto{\pgfqpoint{0.041667in}{0.041667in}}%
\pgfpathlineto{\pgfqpoint{-0.041667in}{0.041667in}}%
\pgfpathlineto{\pgfqpoint{-0.041667in}{-0.041667in}}%
\pgfpathclose%
\pgfusepath{stroke,fill}%
}%
\begin{pgfscope}%
\pgfsys@transformshift{2.278215in}{3.363810in}%
\pgfsys@useobject{currentmarker}{}%
\end{pgfscope}%
\end{pgfscope}%
\begin{pgfscope}%
\definecolor{textcolor}{rgb}{0.000000,0.000000,0.000000}%
\pgfsetstrokecolor{textcolor}%
\pgfsetfillcolor{textcolor}%
\pgftext[x=2.653215in,y=3.276310in,left,base]{\color{textcolor}\rmfamily\fontsize{18.000000}{21.600000}\selectfont LPP\(\displaystyle ^{R,T}\)}%
\end{pgfscope}%
\begin{pgfscope}%
\pgfsetbuttcap%
\pgfsetroundjoin%
\pgfsetlinewidth{1.505625pt}%
\definecolor{currentstroke}{rgb}{0.172549,0.627451,0.172549}%
\pgfsetstrokecolor{currentstroke}%
\pgfsetstrokeopacity{0.900000}%
\pgfsetdash{{5.550000pt}{2.400000pt}}{0.000000pt}%
\pgfpathmoveto{\pgfqpoint{3.589165in}{3.363810in}}%
\pgfpathlineto{\pgfqpoint{3.839165in}{3.363810in}}%
\pgfpathlineto{\pgfqpoint{4.089165in}{3.363810in}}%
\pgfusepath{stroke}%
\end{pgfscope}%
\begin{pgfscope}%
\pgfsetbuttcap%
\pgfsetmiterjoin%
\definecolor{currentfill}{rgb}{0.172549,0.627451,0.172549}%
\pgfsetfillcolor{currentfill}%
\pgfsetfillopacity{0.900000}%
\pgfsetlinewidth{1.003750pt}%
\definecolor{currentstroke}{rgb}{0.172549,0.627451,0.172549}%
\pgfsetstrokecolor{currentstroke}%
\pgfsetstrokeopacity{0.900000}%
\pgfsetdash{}{0pt}%
\pgfsys@defobject{currentmarker}{\pgfqpoint{-0.058926in}{-0.058926in}}{\pgfqpoint{0.058926in}{0.058926in}}{%
\pgfpathmoveto{\pgfqpoint{0.000000in}{-0.058926in}}%
\pgfpathlineto{\pgfqpoint{0.058926in}{0.000000in}}%
\pgfpathlineto{\pgfqpoint{0.000000in}{0.058926in}}%
\pgfpathlineto{\pgfqpoint{-0.058926in}{0.000000in}}%
\pgfpathlineto{\pgfqpoint{0.000000in}{-0.058926in}}%
\pgfpathclose%
\pgfusepath{stroke,fill}%
}%
\begin{pgfscope}%
\pgfsys@transformshift{3.839165in}{3.363810in}%
\pgfsys@useobject{currentmarker}{}%
\end{pgfscope}%
\end{pgfscope}%
\begin{pgfscope}%
\definecolor{textcolor}{rgb}{0.000000,0.000000,0.000000}%
\pgfsetstrokecolor{textcolor}%
\pgfsetfillcolor{textcolor}%
\pgftext[x=4.214165in,y=3.276310in,left,base]{\color{textcolor}\rmfamily\fontsize{18.000000}{21.600000}\selectfont LPP\(\displaystyle ^{T}\)}%
\end{pgfscope}%
\begin{pgfscope}%
\pgfsetbuttcap%
\pgfsetroundjoin%
\pgfsetlinewidth{1.505625pt}%
\definecolor{currentstroke}{rgb}{0.172549,0.627451,0.172549}%
\pgfsetstrokecolor{currentstroke}%
\pgfsetstrokeopacity{0.800000}%
\pgfsetdash{{1.500000pt}{2.475000pt}}{0.000000pt}%
\pgfpathmoveto{\pgfqpoint{4.979197in}{3.363810in}}%
\pgfpathlineto{\pgfqpoint{5.229197in}{3.363810in}}%
\pgfpathlineto{\pgfqpoint{5.479197in}{3.363810in}}%
\pgfusepath{stroke}%
\end{pgfscope}%
\begin{pgfscope}%
\pgfsetbuttcap%
\pgfsetroundjoin%
\definecolor{currentfill}{rgb}{0.172549,0.627451,0.172549}%
\pgfsetfillcolor{currentfill}%
\pgfsetfillopacity{0.800000}%
\pgfsetlinewidth{1.003750pt}%
\definecolor{currentstroke}{rgb}{0.172549,0.627451,0.172549}%
\pgfsetstrokecolor{currentstroke}%
\pgfsetstrokeopacity{0.800000}%
\pgfsetdash{}{0pt}%
\pgfsys@defobject{currentmarker}{\pgfqpoint{-0.041667in}{-0.041667in}}{\pgfqpoint{0.041667in}{0.041667in}}{%
\pgfpathmoveto{\pgfqpoint{0.000000in}{-0.041667in}}%
\pgfpathcurveto{\pgfqpoint{0.011050in}{-0.041667in}}{\pgfqpoint{0.021649in}{-0.037276in}}{\pgfqpoint{0.029463in}{-0.029463in}}%
\pgfpathcurveto{\pgfqpoint{0.037276in}{-0.021649in}}{\pgfqpoint{0.041667in}{-0.011050in}}{\pgfqpoint{0.041667in}{0.000000in}}%
\pgfpathcurveto{\pgfqpoint{0.041667in}{0.011050in}}{\pgfqpoint{0.037276in}{0.021649in}}{\pgfqpoint{0.029463in}{0.029463in}}%
\pgfpathcurveto{\pgfqpoint{0.021649in}{0.037276in}}{\pgfqpoint{0.011050in}{0.041667in}}{\pgfqpoint{0.000000in}{0.041667in}}%
\pgfpathcurveto{\pgfqpoint{-0.011050in}{0.041667in}}{\pgfqpoint{-0.021649in}{0.037276in}}{\pgfqpoint{-0.029463in}{0.029463in}}%
\pgfpathcurveto{\pgfqpoint{-0.037276in}{0.021649in}}{\pgfqpoint{-0.041667in}{0.011050in}}{\pgfqpoint{-0.041667in}{0.000000in}}%
\pgfpathcurveto{\pgfqpoint{-0.041667in}{-0.011050in}}{\pgfqpoint{-0.037276in}{-0.021649in}}{\pgfqpoint{-0.029463in}{-0.029463in}}%
\pgfpathcurveto{\pgfqpoint{-0.021649in}{-0.037276in}}{\pgfqpoint{-0.011050in}{-0.041667in}}{\pgfqpoint{0.000000in}{-0.041667in}}%
\pgfpathlineto{\pgfqpoint{0.000000in}{-0.041667in}}%
\pgfpathclose%
\pgfusepath{stroke,fill}%
}%
\begin{pgfscope}%
\pgfsys@transformshift{5.229197in}{3.363810in}%
\pgfsys@useobject{currentmarker}{}%
\end{pgfscope}%
\end{pgfscope}%
\begin{pgfscope}%
\definecolor{textcolor}{rgb}{0.000000,0.000000,0.000000}%
\pgfsetstrokecolor{textcolor}%
\pgfsetfillcolor{textcolor}%
\pgftext[x=5.604197in,y=3.276310in,left,base]{\color{textcolor}\rmfamily\fontsize{18.000000}{21.600000}\selectfont LPP\(\displaystyle ^{T, P}\)}%
\end{pgfscope}%
\end{pgfpicture}%
\makeatother%
\endgroup%

%% file: figures/DemandPlots/TightLayout/DE_avgWTT_rescaled.pgf
\begingroup%
\makeatletter%
\begin{pgfpicture}%
\pgfpathrectangle{\pgfpointorigin}{\pgfqpoint{7.775300in}{3.651310in}}%
\pgfusepath{use as bounding box, clip}%
\begin{pgfscope}%
\pgfsetbuttcap%
\pgfsetmiterjoin%
\definecolor{currentfill}{rgb}{1.000000,1.000000,1.000000}%
\pgfsetfillcolor{currentfill}%
\pgfsetlinewidth{0.000000pt}%
\definecolor{currentstroke}{rgb}{1.000000,1.000000,1.000000}%
\pgfsetstrokecolor{currentstroke}%
\pgfsetdash{}{0pt}%
\pgfpathmoveto{\pgfqpoint{0.000000in}{0.000000in}}%
\pgfpathlineto{\pgfqpoint{7.775300in}{0.000000in}}%
\pgfpathlineto{\pgfqpoint{7.775300in}{3.651310in}}%
\pgfpathlineto{\pgfqpoint{0.000000in}{3.651310in}}%
\pgfpathlineto{\pgfqpoint{0.000000in}{0.000000in}}%
\pgfpathclose%
\pgfusepath{fill}%
\end{pgfscope}%
\begin{pgfscope}%
\pgfsetbuttcap%
\pgfsetmiterjoin%
\definecolor{currentfill}{rgb}{1.000000,1.000000,1.000000}%
\pgfsetfillcolor{currentfill}%
\pgfsetlinewidth{0.000000pt}%
\definecolor{currentstroke}{rgb}{0.000000,0.000000,0.000000}%
\pgfsetstrokecolor{currentstroke}%
\pgfsetstrokeopacity{0.000000}%
\pgfsetdash{}{0pt}%
\pgfpathmoveto{\pgfqpoint{0.686263in}{0.679475in}}%
\pgfpathlineto{\pgfqpoint{7.675300in}{0.679475in}}%
\pgfpathlineto{\pgfqpoint{7.675300in}{3.135897in}}%
\pgfpathlineto{\pgfqpoint{0.686263in}{3.135897in}}%
\pgfpathlineto{\pgfqpoint{0.686263in}{0.679475in}}%
\pgfpathclose%
\pgfusepath{fill}%
\end{pgfscope}%
\begin{pgfscope}%
\pgfsetbuttcap%
\pgfsetroundjoin%
\definecolor{currentfill}{rgb}{0.000000,0.000000,0.000000}%
\pgfsetfillcolor{currentfill}%
\pgfsetlinewidth{0.803000pt}%
\definecolor{currentstroke}{rgb}{0.000000,0.000000,0.000000}%
\pgfsetstrokecolor{currentstroke}%
\pgfsetdash{}{0pt}%
\pgfsys@defobject{currentmarker}{\pgfqpoint{0.000000in}{-0.048611in}}{\pgfqpoint{0.000000in}{0.000000in}}{%
\pgfpathmoveto{\pgfqpoint{0.000000in}{0.000000in}}%
\pgfpathlineto{\pgfqpoint{0.000000in}{-0.048611in}}%
\pgfusepath{stroke,fill}%
}%
\begin{pgfscope}%
\pgfsys@transformshift{0.686263in}{0.679475in}%
\pgfsys@useobject{currentmarker}{}%
\end{pgfscope}%
\end{pgfscope}%
\begin{pgfscope}%
\definecolor{textcolor}{rgb}{0.000000,0.000000,0.000000}%
\pgfsetstrokecolor{textcolor}%
\pgfsetfillcolor{textcolor}%
\pgftext[x=0.686263in,y=0.582253in,,top]{\color{textcolor}\rmfamily\fontsize{16.000000}{19.200000}\selectfont \(\displaystyle {0.0}\)}%
\end{pgfscope}%
\begin{pgfscope}%
\pgfsetbuttcap%
\pgfsetroundjoin%
\definecolor{currentfill}{rgb}{0.000000,0.000000,0.000000}%
\pgfsetfillcolor{currentfill}%
\pgfsetlinewidth{0.803000pt}%
\definecolor{currentstroke}{rgb}{0.000000,0.000000,0.000000}%
\pgfsetstrokecolor{currentstroke}%
\pgfsetdash{}{0pt}%
\pgfsys@defobject{currentmarker}{\pgfqpoint{0.000000in}{-0.048611in}}{\pgfqpoint{0.000000in}{0.000000in}}{%
\pgfpathmoveto{\pgfqpoint{0.000000in}{0.000000in}}%
\pgfpathlineto{\pgfqpoint{0.000000in}{-0.048611in}}%
\pgfusepath{stroke,fill}%
}%
\begin{pgfscope}%
\pgfsys@transformshift{1.956997in}{0.679475in}%
\pgfsys@useobject{currentmarker}{}%
\end{pgfscope}%
\end{pgfscope}%
\begin{pgfscope}%
\definecolor{textcolor}{rgb}{0.000000,0.000000,0.000000}%
\pgfsetstrokecolor{textcolor}%
\pgfsetfillcolor{textcolor}%
\pgftext[x=1.956997in,y=0.582253in,,top]{\color{textcolor}\rmfamily\fontsize{16.000000}{19.200000}\selectfont \(\displaystyle {0.2}\)}%
\end{pgfscope}%
\begin{pgfscope}%
\pgfsetbuttcap%
\pgfsetroundjoin%
\definecolor{currentfill}{rgb}{0.000000,0.000000,0.000000}%
\pgfsetfillcolor{currentfill}%
\pgfsetlinewidth{0.803000pt}%
\definecolor{currentstroke}{rgb}{0.000000,0.000000,0.000000}%
\pgfsetstrokecolor{currentstroke}%
\pgfsetdash{}{0pt}%
\pgfsys@defobject{currentmarker}{\pgfqpoint{0.000000in}{-0.048611in}}{\pgfqpoint{0.000000in}{0.000000in}}{%
\pgfpathmoveto{\pgfqpoint{0.000000in}{0.000000in}}%
\pgfpathlineto{\pgfqpoint{0.000000in}{-0.048611in}}%
\pgfusepath{stroke,fill}%
}%
\begin{pgfscope}%
\pgfsys@transformshift{3.227731in}{0.679475in}%
\pgfsys@useobject{currentmarker}{}%
\end{pgfscope}%
\end{pgfscope}%
\begin{pgfscope}%
\definecolor{textcolor}{rgb}{0.000000,0.000000,0.000000}%
\pgfsetstrokecolor{textcolor}%
\pgfsetfillcolor{textcolor}%
\pgftext[x=3.227731in,y=0.582253in,,top]{\color{textcolor}\rmfamily\fontsize{16.000000}{19.200000}\selectfont \(\displaystyle {0.4}\)}%
\end{pgfscope}%
\begin{pgfscope}%
\pgfsetbuttcap%
\pgfsetroundjoin%
\definecolor{currentfill}{rgb}{0.000000,0.000000,0.000000}%
\pgfsetfillcolor{currentfill}%
\pgfsetlinewidth{0.803000pt}%
\definecolor{currentstroke}{rgb}{0.000000,0.000000,0.000000}%
\pgfsetstrokecolor{currentstroke}%
\pgfsetdash{}{0pt}%
\pgfsys@defobject{currentmarker}{\pgfqpoint{0.000000in}{-0.048611in}}{\pgfqpoint{0.000000in}{0.000000in}}{%
\pgfpathmoveto{\pgfqpoint{0.000000in}{0.000000in}}%
\pgfpathlineto{\pgfqpoint{0.000000in}{-0.048611in}}%
\pgfusepath{stroke,fill}%
}%
\begin{pgfscope}%
\pgfsys@transformshift{4.498465in}{0.679475in}%
\pgfsys@useobject{currentmarker}{}%
\end{pgfscope}%
\end{pgfscope}%
\begin{pgfscope}%
\definecolor{textcolor}{rgb}{0.000000,0.000000,0.000000}%
\pgfsetstrokecolor{textcolor}%
\pgfsetfillcolor{textcolor}%
\pgftext[x=4.498465in,y=0.582253in,,top]{\color{textcolor}\rmfamily\fontsize{16.000000}{19.200000}\selectfont \(\displaystyle {0.6}\)}%
\end{pgfscope}%
\begin{pgfscope}%
\pgfsetbuttcap%
\pgfsetroundjoin%
\definecolor{currentfill}{rgb}{0.000000,0.000000,0.000000}%
\pgfsetfillcolor{currentfill}%
\pgfsetlinewidth{0.803000pt}%
\definecolor{currentstroke}{rgb}{0.000000,0.000000,0.000000}%
\pgfsetstrokecolor{currentstroke}%
\pgfsetdash{}{0pt}%
\pgfsys@defobject{currentmarker}{\pgfqpoint{0.000000in}{-0.048611in}}{\pgfqpoint{0.000000in}{0.000000in}}{%
\pgfpathmoveto{\pgfqpoint{0.000000in}{0.000000in}}%
\pgfpathlineto{\pgfqpoint{0.000000in}{-0.048611in}}%
\pgfusepath{stroke,fill}%
}%
\begin{pgfscope}%
\pgfsys@transformshift{5.769199in}{0.679475in}%
\pgfsys@useobject{currentmarker}{}%
\end{pgfscope}%
\end{pgfscope}%
\begin{pgfscope}%
\definecolor{textcolor}{rgb}{0.000000,0.000000,0.000000}%
\pgfsetstrokecolor{textcolor}%
\pgfsetfillcolor{textcolor}%
\pgftext[x=5.769199in,y=0.582253in,,top]{\color{textcolor}\rmfamily\fontsize{16.000000}{19.200000}\selectfont \(\displaystyle {0.8}\)}%
\end{pgfscope}%
\begin{pgfscope}%
\pgfsetbuttcap%
\pgfsetroundjoin%
\definecolor{currentfill}{rgb}{0.000000,0.000000,0.000000}%
\pgfsetfillcolor{currentfill}%
\pgfsetlinewidth{0.803000pt}%
\definecolor{currentstroke}{rgb}{0.000000,0.000000,0.000000}%
\pgfsetstrokecolor{currentstroke}%
\pgfsetdash{}{0pt}%
\pgfsys@defobject{currentmarker}{\pgfqpoint{0.000000in}{-0.048611in}}{\pgfqpoint{0.000000in}{0.000000in}}{%
\pgfpathmoveto{\pgfqpoint{0.000000in}{0.000000in}}%
\pgfpathlineto{\pgfqpoint{0.000000in}{-0.048611in}}%
\pgfusepath{stroke,fill}%
}%
\begin{pgfscope}%
\pgfsys@transformshift{7.039933in}{0.679475in}%
\pgfsys@useobject{currentmarker}{}%
\end{pgfscope}%
\end{pgfscope}%
\begin{pgfscope}%
\definecolor{textcolor}{rgb}{0.000000,0.000000,0.000000}%
\pgfsetstrokecolor{textcolor}%
\pgfsetfillcolor{textcolor}%
\pgftext[x=7.039933in,y=0.582253in,,top]{\color{textcolor}\rmfamily\fontsize{16.000000}{19.200000}\selectfont \(\displaystyle {1.0}\)}%
\end{pgfscope}%
\begin{pgfscope}%
\definecolor{textcolor}{rgb}{0.000000,0.000000,0.000000}%
\pgfsetstrokecolor{textcolor}%
\pgfsetfillcolor{textcolor}%
\pgftext[x=4.180781in,y=0.313349in,,top]{\color{textcolor}\rmfamily\fontsize{18.000000}{21.600000}\selectfont \(\displaystyle \lambda\)}%
\end{pgfscope}%
\begin{pgfscope}%
\pgfsetbuttcap%
\pgfsetroundjoin%
\definecolor{currentfill}{rgb}{0.000000,0.000000,0.000000}%
\pgfsetfillcolor{currentfill}%
\pgfsetlinewidth{0.803000pt}%
\definecolor{currentstroke}{rgb}{0.000000,0.000000,0.000000}%
\pgfsetstrokecolor{currentstroke}%
\pgfsetdash{}{0pt}%
\pgfsys@defobject{currentmarker}{\pgfqpoint{-0.048611in}{0.000000in}}{\pgfqpoint{-0.000000in}{0.000000in}}{%
\pgfpathmoveto{\pgfqpoint{-0.000000in}{0.000000in}}%
\pgfpathlineto{\pgfqpoint{-0.048611in}{0.000000in}}%
\pgfusepath{stroke,fill}%
}%
\begin{pgfscope}%
\pgfsys@transformshift{0.686263in}{0.967613in}%
\pgfsys@useobject{currentmarker}{}%
\end{pgfscope}%
\end{pgfscope}%
\begin{pgfscope}%
\definecolor{textcolor}{rgb}{0.000000,0.000000,0.000000}%
\pgfsetstrokecolor{textcolor}%
\pgfsetfillcolor{textcolor}%
\pgftext[x=0.368904in, y=0.884280in, left, base]{\color{textcolor}\rmfamily\fontsize{16.000000}{19.200000}\selectfont \(\displaystyle {24}\)}%
\end{pgfscope}%
\begin{pgfscope}%
\pgfsetbuttcap%
\pgfsetroundjoin%
\definecolor{currentfill}{rgb}{0.000000,0.000000,0.000000}%
\pgfsetfillcolor{currentfill}%
\pgfsetlinewidth{0.803000pt}%
\definecolor{currentstroke}{rgb}{0.000000,0.000000,0.000000}%
\pgfsetstrokecolor{currentstroke}%
\pgfsetdash{}{0pt}%
\pgfsys@defobject{currentmarker}{\pgfqpoint{-0.048611in}{0.000000in}}{\pgfqpoint{-0.000000in}{0.000000in}}{%
\pgfpathmoveto{\pgfqpoint{-0.000000in}{0.000000in}}%
\pgfpathlineto{\pgfqpoint{-0.048611in}{0.000000in}}%
\pgfusepath{stroke,fill}%
}%
\begin{pgfscope}%
\pgfsys@transformshift{0.686263in}{1.383328in}%
\pgfsys@useobject{currentmarker}{}%
\end{pgfscope}%
\end{pgfscope}%
\begin{pgfscope}%
\definecolor{textcolor}{rgb}{0.000000,0.000000,0.000000}%
\pgfsetstrokecolor{textcolor}%
\pgfsetfillcolor{textcolor}%
\pgftext[x=0.368904in, y=1.299994in, left, base]{\color{textcolor}\rmfamily\fontsize{16.000000}{19.200000}\selectfont \(\displaystyle {26}\)}%
\end{pgfscope}%
\begin{pgfscope}%
\pgfsetbuttcap%
\pgfsetroundjoin%
\definecolor{currentfill}{rgb}{0.000000,0.000000,0.000000}%
\pgfsetfillcolor{currentfill}%
\pgfsetlinewidth{0.803000pt}%
\definecolor{currentstroke}{rgb}{0.000000,0.000000,0.000000}%
\pgfsetstrokecolor{currentstroke}%
\pgfsetdash{}{0pt}%
\pgfsys@defobject{currentmarker}{\pgfqpoint{-0.048611in}{0.000000in}}{\pgfqpoint{-0.000000in}{0.000000in}}{%
\pgfpathmoveto{\pgfqpoint{-0.000000in}{0.000000in}}%
\pgfpathlineto{\pgfqpoint{-0.048611in}{0.000000in}}%
\pgfusepath{stroke,fill}%
}%
\begin{pgfscope}%
\pgfsys@transformshift{0.686263in}{1.799042in}%
\pgfsys@useobject{currentmarker}{}%
\end{pgfscope}%
\end{pgfscope}%
\begin{pgfscope}%
\definecolor{textcolor}{rgb}{0.000000,0.000000,0.000000}%
\pgfsetstrokecolor{textcolor}%
\pgfsetfillcolor{textcolor}%
\pgftext[x=0.368904in, y=1.715709in, left, base]{\color{textcolor}\rmfamily\fontsize{16.000000}{19.200000}\selectfont \(\displaystyle {28}\)}%
\end{pgfscope}%
\begin{pgfscope}%
\pgfsetbuttcap%
\pgfsetroundjoin%
\definecolor{currentfill}{rgb}{0.000000,0.000000,0.000000}%
\pgfsetfillcolor{currentfill}%
\pgfsetlinewidth{0.803000pt}%
\definecolor{currentstroke}{rgb}{0.000000,0.000000,0.000000}%
\pgfsetstrokecolor{currentstroke}%
\pgfsetdash{}{0pt}%
\pgfsys@defobject{currentmarker}{\pgfqpoint{-0.048611in}{0.000000in}}{\pgfqpoint{-0.000000in}{0.000000in}}{%
\pgfpathmoveto{\pgfqpoint{-0.000000in}{0.000000in}}%
\pgfpathlineto{\pgfqpoint{-0.048611in}{0.000000in}}%
\pgfusepath{stroke,fill}%
}%
\begin{pgfscope}%
\pgfsys@transformshift{0.686263in}{2.214756in}%
\pgfsys@useobject{currentmarker}{}%
\end{pgfscope}%
\end{pgfscope}%
\begin{pgfscope}%
\definecolor{textcolor}{rgb}{0.000000,0.000000,0.000000}%
\pgfsetstrokecolor{textcolor}%
\pgfsetfillcolor{textcolor}%
\pgftext[x=0.368904in, y=2.131423in, left, base]{\color{textcolor}\rmfamily\fontsize{16.000000}{19.200000}\selectfont \(\displaystyle {30}\)}%
\end{pgfscope}%
\begin{pgfscope}%
\pgfsetbuttcap%
\pgfsetroundjoin%
\definecolor{currentfill}{rgb}{0.000000,0.000000,0.000000}%
\pgfsetfillcolor{currentfill}%
\pgfsetlinewidth{0.803000pt}%
\definecolor{currentstroke}{rgb}{0.000000,0.000000,0.000000}%
\pgfsetstrokecolor{currentstroke}%
\pgfsetdash{}{0pt}%
\pgfsys@defobject{currentmarker}{\pgfqpoint{-0.048611in}{0.000000in}}{\pgfqpoint{-0.000000in}{0.000000in}}{%
\pgfpathmoveto{\pgfqpoint{-0.000000in}{0.000000in}}%
\pgfpathlineto{\pgfqpoint{-0.048611in}{0.000000in}}%
\pgfusepath{stroke,fill}%
}%
\begin{pgfscope}%
\pgfsys@transformshift{0.686263in}{2.630471in}%
\pgfsys@useobject{currentmarker}{}%
\end{pgfscope}%
\end{pgfscope}%
\begin{pgfscope}%
\definecolor{textcolor}{rgb}{0.000000,0.000000,0.000000}%
\pgfsetstrokecolor{textcolor}%
\pgfsetfillcolor{textcolor}%
\pgftext[x=0.368904in, y=2.547137in, left, base]{\color{textcolor}\rmfamily\fontsize{16.000000}{19.200000}\selectfont \(\displaystyle {32}\)}%
\end{pgfscope}%
\begin{pgfscope}%
\pgfsetbuttcap%
\pgfsetroundjoin%
\definecolor{currentfill}{rgb}{0.000000,0.000000,0.000000}%
\pgfsetfillcolor{currentfill}%
\pgfsetlinewidth{0.803000pt}%
\definecolor{currentstroke}{rgb}{0.000000,0.000000,0.000000}%
\pgfsetstrokecolor{currentstroke}%
\pgfsetdash{}{0pt}%
\pgfsys@defobject{currentmarker}{\pgfqpoint{-0.048611in}{0.000000in}}{\pgfqpoint{-0.000000in}{0.000000in}}{%
\pgfpathmoveto{\pgfqpoint{-0.000000in}{0.000000in}}%
\pgfpathlineto{\pgfqpoint{-0.048611in}{0.000000in}}%
\pgfusepath{stroke,fill}%
}%
\begin{pgfscope}%
\pgfsys@transformshift{0.686263in}{3.046185in}%
\pgfsys@useobject{currentmarker}{}%
\end{pgfscope}%
\end{pgfscope}%
\begin{pgfscope}%
\definecolor{textcolor}{rgb}{0.000000,0.000000,0.000000}%
\pgfsetstrokecolor{textcolor}%
\pgfsetfillcolor{textcolor}%
\pgftext[x=0.368904in, y=2.962852in, left, base]{\color{textcolor}\rmfamily\fontsize{16.000000}{19.200000}\selectfont \(\displaystyle {34}\)}%
\end{pgfscope}%
\begin{pgfscope}%
\definecolor{textcolor}{rgb}{0.000000,0.000000,0.000000}%
\pgfsetstrokecolor{textcolor}%
\pgfsetfillcolor{textcolor}%
\pgftext[x=0.313349in,y=1.907686in,,bottom,rotate=90.000000]{\color{textcolor}\rmfamily\fontsize{18.000000}{21.600000}\selectfont minutes}%
\end{pgfscope}%
\begin{pgfscope}%
\pgfpathrectangle{\pgfqpoint{0.686263in}{0.679475in}}{\pgfqpoint{6.989037in}{2.456422in}}%
\pgfusepath{clip}%
\pgfsetrectcap%
\pgfsetroundjoin%
\pgfsetlinewidth{1.505625pt}%
\definecolor{currentstroke}{rgb}{1.000000,0.498039,0.054902}%
\pgfsetstrokecolor{currentstroke}%
\pgfsetdash{}{0pt}%
\pgfpathmoveto{\pgfqpoint{1.321630in}{3.024241in}}%
\pgfpathlineto{\pgfqpoint{2.274680in}{2.625542in}}%
\pgfpathlineto{\pgfqpoint{3.863098in}{2.415006in}}%
\pgfpathlineto{\pgfqpoint{5.451515in}{2.131621in}}%
\pgfpathlineto{\pgfqpoint{7.039933in}{1.816683in}}%
\pgfusepath{stroke}%
\end{pgfscope}%
\begin{pgfscope}%
\pgfpathrectangle{\pgfqpoint{0.686263in}{0.679475in}}{\pgfqpoint{6.989037in}{2.456422in}}%
\pgfusepath{clip}%
\pgfsetbuttcap%
\pgfsetmiterjoin%
\definecolor{currentfill}{rgb}{1.000000,0.498039,0.054902}%
\pgfsetfillcolor{currentfill}%
\pgfsetlinewidth{0.752812pt}%
\definecolor{currentstroke}{rgb}{1.000000,1.000000,1.000000}%
\pgfsetstrokecolor{currentstroke}%
\pgfsetdash{}{0pt}%
\pgfsys@defobject{currentmarker}{\pgfqpoint{-0.041667in}{-0.041667in}}{\pgfqpoint{0.041667in}{0.041667in}}{%
\pgfpathmoveto{\pgfqpoint{-0.041667in}{-0.041667in}}%
\pgfpathlineto{\pgfqpoint{0.041667in}{-0.041667in}}%
\pgfpathlineto{\pgfqpoint{0.041667in}{0.041667in}}%
\pgfpathlineto{\pgfqpoint{-0.041667in}{0.041667in}}%
\pgfpathlineto{\pgfqpoint{-0.041667in}{-0.041667in}}%
\pgfpathclose%
\pgfusepath{stroke,fill}%
}%
\begin{pgfscope}%
\pgfsys@transformshift{1.321630in}{3.024241in}%
\pgfsys@useobject{currentmarker}{}%
\end{pgfscope}%
\begin{pgfscope}%
\pgfsys@transformshift{2.274680in}{2.625542in}%
\pgfsys@useobject{currentmarker}{}%
\end{pgfscope}%
\begin{pgfscope}%
\pgfsys@transformshift{3.863098in}{2.415006in}%
\pgfsys@useobject{currentmarker}{}%
\end{pgfscope}%
\begin{pgfscope}%
\pgfsys@transformshift{5.451515in}{2.131621in}%
\pgfsys@useobject{currentmarker}{}%
\end{pgfscope}%
\begin{pgfscope}%
\pgfsys@transformshift{7.039933in}{1.816683in}%
\pgfsys@useobject{currentmarker}{}%
\end{pgfscope}%
\end{pgfscope}%
\begin{pgfscope}%
\pgfpathrectangle{\pgfqpoint{0.686263in}{0.679475in}}{\pgfqpoint{6.989037in}{2.456422in}}%
\pgfusepath{clip}%
\pgfsetbuttcap%
\pgfsetroundjoin%
\pgfsetlinewidth{1.505625pt}%
\definecolor{currentstroke}{rgb}{1.000000,0.498039,0.054902}%
\pgfsetstrokecolor{currentstroke}%
\pgfsetstrokeopacity{0.900000}%
\pgfsetdash{{7.500000pt}{7.500000pt}}{0.000000pt}%
\pgfpathmoveto{\pgfqpoint{1.321630in}{2.095892in}}%
\pgfpathlineto{\pgfqpoint{2.274680in}{1.998632in}}%
\pgfpathlineto{\pgfqpoint{3.863098in}{2.119083in}}%
\pgfpathlineto{\pgfqpoint{5.451515in}{1.750878in}}%
\pgfpathlineto{\pgfqpoint{7.039933in}{1.299926in}}%
\pgfusepath{stroke}%
\end{pgfscope}%
\begin{pgfscope}%
\pgfpathrectangle{\pgfqpoint{0.686263in}{0.679475in}}{\pgfqpoint{6.989037in}{2.456422in}}%
\pgfusepath{clip}%
\pgfsetbuttcap%
\pgfsetmiterjoin%
\definecolor{currentfill}{rgb}{1.000000,0.498039,0.054902}%
\pgfsetfillcolor{currentfill}%
\pgfsetfillopacity{0.900000}%
\pgfsetlinewidth{0.752812pt}%
\definecolor{currentstroke}{rgb}{1.000000,1.000000,1.000000}%
\pgfsetstrokecolor{currentstroke}%
\pgfsetdash{}{0pt}%
\pgfsys@defobject{currentmarker}{\pgfqpoint{-0.058926in}{-0.058926in}}{\pgfqpoint{0.058926in}{0.058926in}}{%
\pgfpathmoveto{\pgfqpoint{0.000000in}{-0.058926in}}%
\pgfpathlineto{\pgfqpoint{0.058926in}{0.000000in}}%
\pgfpathlineto{\pgfqpoint{0.000000in}{0.058926in}}%
\pgfpathlineto{\pgfqpoint{-0.058926in}{0.000000in}}%
\pgfpathlineto{\pgfqpoint{0.000000in}{-0.058926in}}%
\pgfpathclose%
\pgfusepath{stroke,fill}%
}%
\begin{pgfscope}%
\pgfsys@transformshift{1.321630in}{2.095892in}%
\pgfsys@useobject{currentmarker}{}%
\end{pgfscope}%
\begin{pgfscope}%
\pgfsys@transformshift{2.274680in}{1.998632in}%
\pgfsys@useobject{currentmarker}{}%
\end{pgfscope}%
\begin{pgfscope}%
\pgfsys@transformshift{3.863098in}{2.119083in}%
\pgfsys@useobject{currentmarker}{}%
\end{pgfscope}%
\begin{pgfscope}%
\pgfsys@transformshift{5.451515in}{1.750878in}%
\pgfsys@useobject{currentmarker}{}%
\end{pgfscope}%
\begin{pgfscope}%
\pgfsys@transformshift{7.039933in}{1.299926in}%
\pgfsys@useobject{currentmarker}{}%
\end{pgfscope}%
\end{pgfscope}%
\begin{pgfscope}%
\pgfpathrectangle{\pgfqpoint{0.686263in}{0.679475in}}{\pgfqpoint{6.989037in}{2.456422in}}%
\pgfusepath{clip}%
\pgfsetbuttcap%
\pgfsetroundjoin%
\pgfsetlinewidth{1.505625pt}%
\definecolor{currentstroke}{rgb}{1.000000,0.498039,0.054902}%
\pgfsetstrokecolor{currentstroke}%
\pgfsetstrokeopacity{0.800000}%
\pgfsetdash{{1.500000pt}{4.500000pt}}{0.000000pt}%
\pgfpathmoveto{\pgfqpoint{1.321630in}{0.883817in}}%
\pgfpathlineto{\pgfqpoint{2.274680in}{0.860501in}}%
\pgfpathlineto{\pgfqpoint{3.863098in}{0.913380in}}%
\pgfpathlineto{\pgfqpoint{5.451515in}{1.001100in}}%
\pgfpathlineto{\pgfqpoint{7.039933in}{0.791131in}}%
\pgfusepath{stroke}%
\end{pgfscope}%
\begin{pgfscope}%
\pgfpathrectangle{\pgfqpoint{0.686263in}{0.679475in}}{\pgfqpoint{6.989037in}{2.456422in}}%
\pgfusepath{clip}%
\pgfsetbuttcap%
\pgfsetroundjoin%
\definecolor{currentfill}{rgb}{1.000000,0.498039,0.054902}%
\pgfsetfillcolor{currentfill}%
\pgfsetfillopacity{0.800000}%
\pgfsetlinewidth{0.752812pt}%
\definecolor{currentstroke}{rgb}{1.000000,1.000000,1.000000}%
\pgfsetstrokecolor{currentstroke}%
\pgfsetdash{}{0pt}%
\pgfsys@defobject{currentmarker}{\pgfqpoint{-0.041667in}{-0.041667in}}{\pgfqpoint{0.041667in}{0.041667in}}{%
\pgfpathmoveto{\pgfqpoint{0.000000in}{-0.041667in}}%
\pgfpathcurveto{\pgfqpoint{0.011050in}{-0.041667in}}{\pgfqpoint{0.021649in}{-0.037276in}}{\pgfqpoint{0.029463in}{-0.029463in}}%
\pgfpathcurveto{\pgfqpoint{0.037276in}{-0.021649in}}{\pgfqpoint{0.041667in}{-0.011050in}}{\pgfqpoint{0.041667in}{0.000000in}}%
\pgfpathcurveto{\pgfqpoint{0.041667in}{0.011050in}}{\pgfqpoint{0.037276in}{0.021649in}}{\pgfqpoint{0.029463in}{0.029463in}}%
\pgfpathcurveto{\pgfqpoint{0.021649in}{0.037276in}}{\pgfqpoint{0.011050in}{0.041667in}}{\pgfqpoint{0.000000in}{0.041667in}}%
\pgfpathcurveto{\pgfqpoint{-0.011050in}{0.041667in}}{\pgfqpoint{-0.021649in}{0.037276in}}{\pgfqpoint{-0.029463in}{0.029463in}}%
\pgfpathcurveto{\pgfqpoint{-0.037276in}{0.021649in}}{\pgfqpoint{-0.041667in}{0.011050in}}{\pgfqpoint{-0.041667in}{0.000000in}}%
\pgfpathcurveto{\pgfqpoint{-0.041667in}{-0.011050in}}{\pgfqpoint{-0.037276in}{-0.021649in}}{\pgfqpoint{-0.029463in}{-0.029463in}}%
\pgfpathcurveto{\pgfqpoint{-0.021649in}{-0.037276in}}{\pgfqpoint{-0.011050in}{-0.041667in}}{\pgfqpoint{0.000000in}{-0.041667in}}%
\pgfpathlineto{\pgfqpoint{0.000000in}{-0.041667in}}%
\pgfpathclose%
\pgfusepath{stroke,fill}%
}%
\begin{pgfscope}%
\pgfsys@transformshift{1.321630in}{0.883817in}%
\pgfsys@useobject{currentmarker}{}%
\end{pgfscope}%
\begin{pgfscope}%
\pgfsys@transformshift{2.274680in}{0.860501in}%
\pgfsys@useobject{currentmarker}{}%
\end{pgfscope}%
\begin{pgfscope}%
\pgfsys@transformshift{3.863098in}{0.913380in}%
\pgfsys@useobject{currentmarker}{}%
\end{pgfscope}%
\begin{pgfscope}%
\pgfsys@transformshift{5.451515in}{1.001100in}%
\pgfsys@useobject{currentmarker}{}%
\end{pgfscope}%
\begin{pgfscope}%
\pgfsys@transformshift{7.039933in}{0.791131in}%
\pgfsys@useobject{currentmarker}{}%
\end{pgfscope}%
\end{pgfscope}%
\begin{pgfscope}%
\pgfsetrectcap%
\pgfsetmiterjoin%
\pgfsetlinewidth{0.803000pt}%
\definecolor{currentstroke}{rgb}{0.000000,0.000000,0.000000}%
\pgfsetstrokecolor{currentstroke}%
\pgfsetdash{}{0pt}%
\pgfpathmoveto{\pgfqpoint{0.686263in}{0.679475in}}%
\pgfpathlineto{\pgfqpoint{0.686263in}{3.135897in}}%
\pgfusepath{stroke}%
\end{pgfscope}%
\begin{pgfscope}%
\pgfsetrectcap%
\pgfsetmiterjoin%
\pgfsetlinewidth{0.803000pt}%
\definecolor{currentstroke}{rgb}{0.000000,0.000000,0.000000}%
\pgfsetstrokecolor{currentstroke}%
\pgfsetdash{}{0pt}%
\pgfpathmoveto{\pgfqpoint{7.675300in}{0.679475in}}%
\pgfpathlineto{\pgfqpoint{7.675300in}{3.135897in}}%
\pgfusepath{stroke}%
\end{pgfscope}%
\begin{pgfscope}%
\pgfsetrectcap%
\pgfsetmiterjoin%
\pgfsetlinewidth{0.803000pt}%
\definecolor{currentstroke}{rgb}{0.000000,0.000000,0.000000}%
\pgfsetstrokecolor{currentstroke}%
\pgfsetdash{}{0pt}%
\pgfpathmoveto{\pgfqpoint{0.686263in}{0.679475in}}%
\pgfpathlineto{\pgfqpoint{7.675300in}{0.679475in}}%
\pgfusepath{stroke}%
\end{pgfscope}%
\begin{pgfscope}%
\pgfsetrectcap%
\pgfsetmiterjoin%
\pgfsetlinewidth{0.803000pt}%
\definecolor{currentstroke}{rgb}{0.000000,0.000000,0.000000}%
\pgfsetstrokecolor{currentstroke}%
\pgfsetdash{}{0pt}%
\pgfpathmoveto{\pgfqpoint{0.686263in}{3.135897in}}%
\pgfpathlineto{\pgfqpoint{7.675300in}{3.135897in}}%
\pgfusepath{stroke}%
\end{pgfscope}%
\begin{pgfscope}%
\pgfsetrectcap%
\pgfsetroundjoin%
\pgfsetlinewidth{1.505625pt}%
\definecolor{currentstroke}{rgb}{1.000000,0.498039,0.054902}%
\pgfsetstrokecolor{currentstroke}%
\pgfsetdash{}{0pt}%
\pgfpathmoveto{\pgfqpoint{2.028215in}{3.363810in}}%
\pgfpathlineto{\pgfqpoint{2.278215in}{3.363810in}}%
\pgfpathlineto{\pgfqpoint{2.528215in}{3.363810in}}%
\pgfusepath{stroke}%
\end{pgfscope}%
\begin{pgfscope}%
\pgfsetbuttcap%
\pgfsetmiterjoin%
\definecolor{currentfill}{rgb}{1.000000,0.498039,0.054902}%
\pgfsetfillcolor{currentfill}%
\pgfsetlinewidth{1.003750pt}%
\definecolor{currentstroke}{rgb}{1.000000,0.498039,0.054902}%
\pgfsetstrokecolor{currentstroke}%
\pgfsetdash{}{0pt}%
\pgfsys@defobject{currentmarker}{\pgfqpoint{-0.041667in}{-0.041667in}}{\pgfqpoint{0.041667in}{0.041667in}}{%
\pgfpathmoveto{\pgfqpoint{-0.041667in}{-0.041667in}}%
\pgfpathlineto{\pgfqpoint{0.041667in}{-0.041667in}}%
\pgfpathlineto{\pgfqpoint{0.041667in}{0.041667in}}%
\pgfpathlineto{\pgfqpoint{-0.041667in}{0.041667in}}%
\pgfpathlineto{\pgfqpoint{-0.041667in}{-0.041667in}}%
\pgfpathclose%
\pgfusepath{stroke,fill}%
}%
\begin{pgfscope}%
\pgfsys@transformshift{2.278215in}{3.363810in}%
\pgfsys@useobject{currentmarker}{}%
\end{pgfscope}%
\end{pgfscope}%
\begin{pgfscope}%
\definecolor{textcolor}{rgb}{0.000000,0.000000,0.000000}%
\pgfsetstrokecolor{textcolor}%
\pgfsetfillcolor{textcolor}%
\pgftext[x=2.653215in,y=3.276310in,left,base]{\color{textcolor}\rmfamily\fontsize{18.000000}{21.600000}\selectfont LPP\(\displaystyle ^{R,T}\)}%
\end{pgfscope}%
\begin{pgfscope}%
\pgfsetbuttcap%
\pgfsetroundjoin%
\pgfsetlinewidth{1.505625pt}%
\definecolor{currentstroke}{rgb}{1.000000,0.498039,0.054902}%
\pgfsetstrokecolor{currentstroke}%
\pgfsetstrokeopacity{0.900000}%
\pgfsetdash{{5.550000pt}{2.400000pt}}{0.000000pt}%
\pgfpathmoveto{\pgfqpoint{3.589165in}{3.363810in}}%
\pgfpathlineto{\pgfqpoint{3.839165in}{3.363810in}}%
\pgfpathlineto{\pgfqpoint{4.089165in}{3.363810in}}%
\pgfusepath{stroke}%
\end{pgfscope}%
\begin{pgfscope}%
\pgfsetbuttcap%
\pgfsetmiterjoin%
\definecolor{currentfill}{rgb}{1.000000,0.498039,0.054902}%
\pgfsetfillcolor{currentfill}%
\pgfsetfillopacity{0.900000}%
\pgfsetlinewidth{1.003750pt}%
\definecolor{currentstroke}{rgb}{1.000000,0.498039,0.054902}%
\pgfsetstrokecolor{currentstroke}%
\pgfsetstrokeopacity{0.900000}%
\pgfsetdash{}{0pt}%
\pgfsys@defobject{currentmarker}{\pgfqpoint{-0.058926in}{-0.058926in}}{\pgfqpoint{0.058926in}{0.058926in}}{%
\pgfpathmoveto{\pgfqpoint{0.000000in}{-0.058926in}}%
\pgfpathlineto{\pgfqpoint{0.058926in}{0.000000in}}%
\pgfpathlineto{\pgfqpoint{0.000000in}{0.058926in}}%
\pgfpathlineto{\pgfqpoint{-0.058926in}{0.000000in}}%
\pgfpathlineto{\pgfqpoint{0.000000in}{-0.058926in}}%
\pgfpathclose%
\pgfusepath{stroke,fill}%
}%
\begin{pgfscope}%
\pgfsys@transformshift{3.839165in}{3.363810in}%
\pgfsys@useobject{currentmarker}{}%
\end{pgfscope}%
\end{pgfscope}%
\begin{pgfscope}%
\definecolor{textcolor}{rgb}{0.000000,0.000000,0.000000}%
\pgfsetstrokecolor{textcolor}%
\pgfsetfillcolor{textcolor}%
\pgftext[x=4.214165in,y=3.276310in,left,base]{\color{textcolor}\rmfamily\fontsize{18.000000}{21.600000}\selectfont LPP\(\displaystyle ^{T}\)}%
\end{pgfscope}%
\begin{pgfscope}%
\pgfsetbuttcap%
\pgfsetroundjoin%
\pgfsetlinewidth{1.505625pt}%
\definecolor{currentstroke}{rgb}{1.000000,0.498039,0.054902}%
\pgfsetstrokecolor{currentstroke}%
\pgfsetstrokeopacity{0.800000}%
\pgfsetdash{{1.500000pt}{2.475000pt}}{0.000000pt}%
\pgfpathmoveto{\pgfqpoint{4.979197in}{3.363810in}}%
\pgfpathlineto{\pgfqpoint{5.229197in}{3.363810in}}%
\pgfpathlineto{\pgfqpoint{5.479197in}{3.363810in}}%
\pgfusepath{stroke}%
\end{pgfscope}%
\begin{pgfscope}%
\pgfsetbuttcap%
\pgfsetroundjoin%
\definecolor{currentfill}{rgb}{1.000000,0.498039,0.054902}%
\pgfsetfillcolor{currentfill}%
\pgfsetfillopacity{0.800000}%
\pgfsetlinewidth{1.003750pt}%
\definecolor{currentstroke}{rgb}{1.000000,0.498039,0.054902}%
\pgfsetstrokecolor{currentstroke}%
\pgfsetstrokeopacity{0.800000}%
\pgfsetdash{}{0pt}%
\pgfsys@defobject{currentmarker}{\pgfqpoint{-0.041667in}{-0.041667in}}{\pgfqpoint{0.041667in}{0.041667in}}{%
\pgfpathmoveto{\pgfqpoint{0.000000in}{-0.041667in}}%
\pgfpathcurveto{\pgfqpoint{0.011050in}{-0.041667in}}{\pgfqpoint{0.021649in}{-0.037276in}}{\pgfqpoint{0.029463in}{-0.029463in}}%
\pgfpathcurveto{\pgfqpoint{0.037276in}{-0.021649in}}{\pgfqpoint{0.041667in}{-0.011050in}}{\pgfqpoint{0.041667in}{0.000000in}}%
\pgfpathcurveto{\pgfqpoint{0.041667in}{0.011050in}}{\pgfqpoint{0.037276in}{0.021649in}}{\pgfqpoint{0.029463in}{0.029463in}}%
\pgfpathcurveto{\pgfqpoint{0.021649in}{0.037276in}}{\pgfqpoint{0.011050in}{0.041667in}}{\pgfqpoint{0.000000in}{0.041667in}}%
\pgfpathcurveto{\pgfqpoint{-0.011050in}{0.041667in}}{\pgfqpoint{-0.021649in}{0.037276in}}{\pgfqpoint{-0.029463in}{0.029463in}}%
\pgfpathcurveto{\pgfqpoint{-0.037276in}{0.021649in}}{\pgfqpoint{-0.041667in}{0.011050in}}{\pgfqpoint{-0.041667in}{0.000000in}}%
\pgfpathcurveto{\pgfqpoint{-0.041667in}{-0.011050in}}{\pgfqpoint{-0.037276in}{-0.021649in}}{\pgfqpoint{-0.029463in}{-0.029463in}}%
\pgfpathcurveto{\pgfqpoint{-0.021649in}{-0.037276in}}{\pgfqpoint{-0.011050in}{-0.041667in}}{\pgfqpoint{0.000000in}{-0.041667in}}%
\pgfpathlineto{\pgfqpoint{0.000000in}{-0.041667in}}%
\pgfpathclose%
\pgfusepath{stroke,fill}%
}%
\begin{pgfscope}%
\pgfsys@transformshift{5.229197in}{3.363810in}%
\pgfsys@useobject{currentmarker}{}%
\end{pgfscope}%
\end{pgfscope}%
\begin{pgfscope}%
\definecolor{textcolor}{rgb}{0.000000,0.000000,0.000000}%
\pgfsetstrokecolor{textcolor}%
\pgfsetfillcolor{textcolor}%
\pgftext[x=5.604197in,y=3.276310in,left,base]{\color{textcolor}\rmfamily\fontsize{18.000000}{21.600000}\selectfont LPP\(\displaystyle ^{T, P}\)}%
\end{pgfscope}%
\end{pgfpicture}%
\makeatother%
\endgroup%

%% file: figures/DemandPlots/TightLayout/DE_UB_rescaled.pgf
\begingroup%
\makeatletter%
\begin{pgfpicture}%
\pgfpathrectangle{\pgfpointorigin}{\pgfqpoint{7.782018in}{3.651310in}}%
\pgfusepath{use as bounding box, clip}%
\begin{pgfscope}%
\pgfsetbuttcap%
\pgfsetmiterjoin%
\definecolor{currentfill}{rgb}{1.000000,1.000000,1.000000}%
\pgfsetfillcolor{currentfill}%
\pgfsetlinewidth{0.000000pt}%
\definecolor{currentstroke}{rgb}{1.000000,1.000000,1.000000}%
\pgfsetstrokecolor{currentstroke}%
\pgfsetdash{}{0pt}%
\pgfpathmoveto{\pgfqpoint{0.000000in}{0.000000in}}%
\pgfpathlineto{\pgfqpoint{7.782018in}{0.000000in}}%
\pgfpathlineto{\pgfqpoint{7.782018in}{3.651310in}}%
\pgfpathlineto{\pgfqpoint{0.000000in}{3.651310in}}%
\pgfpathlineto{\pgfqpoint{0.000000in}{0.000000in}}%
\pgfpathclose%
\pgfusepath{fill}%
\end{pgfscope}%
\begin{pgfscope}%
\pgfsetbuttcap%
\pgfsetmiterjoin%
\definecolor{currentfill}{rgb}{1.000000,1.000000,1.000000}%
\pgfsetfillcolor{currentfill}%
\pgfsetlinewidth{0.000000pt}%
\definecolor{currentstroke}{rgb}{0.000000,0.000000,0.000000}%
\pgfsetstrokecolor{currentstroke}%
\pgfsetstrokeopacity{0.000000}%
\pgfsetdash{}{0pt}%
\pgfpathmoveto{\pgfqpoint{1.126536in}{0.679475in}}%
\pgfpathlineto{\pgfqpoint{7.682018in}{0.679475in}}%
\pgfpathlineto{\pgfqpoint{7.682018in}{3.135897in}}%
\pgfpathlineto{\pgfqpoint{1.126536in}{3.135897in}}%
\pgfpathlineto{\pgfqpoint{1.126536in}{0.679475in}}%
\pgfpathclose%
\pgfusepath{fill}%
\end{pgfscope}%
\begin{pgfscope}%
\pgfsetbuttcap%
\pgfsetroundjoin%
\definecolor{currentfill}{rgb}{0.000000,0.000000,0.000000}%
\pgfsetfillcolor{currentfill}%
\pgfsetlinewidth{0.803000pt}%
\definecolor{currentstroke}{rgb}{0.000000,0.000000,0.000000}%
\pgfsetstrokecolor{currentstroke}%
\pgfsetdash{}{0pt}%
\pgfsys@defobject{currentmarker}{\pgfqpoint{0.000000in}{-0.048611in}}{\pgfqpoint{0.000000in}{0.000000in}}{%
\pgfpathmoveto{\pgfqpoint{0.000000in}{0.000000in}}%
\pgfpathlineto{\pgfqpoint{0.000000in}{-0.048611in}}%
\pgfusepath{stroke,fill}%
}%
\begin{pgfscope}%
\pgfsys@transformshift{1.126536in}{0.679475in}%
\pgfsys@useobject{currentmarker}{}%
\end{pgfscope}%
\end{pgfscope}%
\begin{pgfscope}%
\definecolor{textcolor}{rgb}{0.000000,0.000000,0.000000}%
\pgfsetstrokecolor{textcolor}%
\pgfsetfillcolor{textcolor}%
\pgftext[x=1.126536in,y=0.582253in,,top]{\color{textcolor}\rmfamily\fontsize{16.000000}{19.200000}\selectfont \(\displaystyle {0.0}\)}%
\end{pgfscope}%
\begin{pgfscope}%
\pgfsetbuttcap%
\pgfsetroundjoin%
\definecolor{currentfill}{rgb}{0.000000,0.000000,0.000000}%
\pgfsetfillcolor{currentfill}%
\pgfsetlinewidth{0.803000pt}%
\definecolor{currentstroke}{rgb}{0.000000,0.000000,0.000000}%
\pgfsetstrokecolor{currentstroke}%
\pgfsetdash{}{0pt}%
\pgfsys@defobject{currentmarker}{\pgfqpoint{0.000000in}{-0.048611in}}{\pgfqpoint{0.000000in}{0.000000in}}{%
\pgfpathmoveto{\pgfqpoint{0.000000in}{0.000000in}}%
\pgfpathlineto{\pgfqpoint{0.000000in}{-0.048611in}}%
\pgfusepath{stroke,fill}%
}%
\begin{pgfscope}%
\pgfsys@transformshift{2.318442in}{0.679475in}%
\pgfsys@useobject{currentmarker}{}%
\end{pgfscope}%
\end{pgfscope}%
\begin{pgfscope}%
\definecolor{textcolor}{rgb}{0.000000,0.000000,0.000000}%
\pgfsetstrokecolor{textcolor}%
\pgfsetfillcolor{textcolor}%
\pgftext[x=2.318442in,y=0.582253in,,top]{\color{textcolor}\rmfamily\fontsize{16.000000}{19.200000}\selectfont \(\displaystyle {0.2}\)}%
\end{pgfscope}%
\begin{pgfscope}%
\pgfsetbuttcap%
\pgfsetroundjoin%
\definecolor{currentfill}{rgb}{0.000000,0.000000,0.000000}%
\pgfsetfillcolor{currentfill}%
\pgfsetlinewidth{0.803000pt}%
\definecolor{currentstroke}{rgb}{0.000000,0.000000,0.000000}%
\pgfsetstrokecolor{currentstroke}%
\pgfsetdash{}{0pt}%
\pgfsys@defobject{currentmarker}{\pgfqpoint{0.000000in}{-0.048611in}}{\pgfqpoint{0.000000in}{0.000000in}}{%
\pgfpathmoveto{\pgfqpoint{0.000000in}{0.000000in}}%
\pgfpathlineto{\pgfqpoint{0.000000in}{-0.048611in}}%
\pgfusepath{stroke,fill}%
}%
\begin{pgfscope}%
\pgfsys@transformshift{3.510347in}{0.679475in}%
\pgfsys@useobject{currentmarker}{}%
\end{pgfscope}%
\end{pgfscope}%
\begin{pgfscope}%
\definecolor{textcolor}{rgb}{0.000000,0.000000,0.000000}%
\pgfsetstrokecolor{textcolor}%
\pgfsetfillcolor{textcolor}%
\pgftext[x=3.510347in,y=0.582253in,,top]{\color{textcolor}\rmfamily\fontsize{16.000000}{19.200000}\selectfont \(\displaystyle {0.4}\)}%
\end{pgfscope}%
\begin{pgfscope}%
\pgfsetbuttcap%
\pgfsetroundjoin%
\definecolor{currentfill}{rgb}{0.000000,0.000000,0.000000}%
\pgfsetfillcolor{currentfill}%
\pgfsetlinewidth{0.803000pt}%
\definecolor{currentstroke}{rgb}{0.000000,0.000000,0.000000}%
\pgfsetstrokecolor{currentstroke}%
\pgfsetdash{}{0pt}%
\pgfsys@defobject{currentmarker}{\pgfqpoint{0.000000in}{-0.048611in}}{\pgfqpoint{0.000000in}{0.000000in}}{%
\pgfpathmoveto{\pgfqpoint{0.000000in}{0.000000in}}%
\pgfpathlineto{\pgfqpoint{0.000000in}{-0.048611in}}%
\pgfusepath{stroke,fill}%
}%
\begin{pgfscope}%
\pgfsys@transformshift{4.702253in}{0.679475in}%
\pgfsys@useobject{currentmarker}{}%
\end{pgfscope}%
\end{pgfscope}%
\begin{pgfscope}%
\definecolor{textcolor}{rgb}{0.000000,0.000000,0.000000}%
\pgfsetstrokecolor{textcolor}%
\pgfsetfillcolor{textcolor}%
\pgftext[x=4.702253in,y=0.582253in,,top]{\color{textcolor}\rmfamily\fontsize{16.000000}{19.200000}\selectfont \(\displaystyle {0.6}\)}%
\end{pgfscope}%
\begin{pgfscope}%
\pgfsetbuttcap%
\pgfsetroundjoin%
\definecolor{currentfill}{rgb}{0.000000,0.000000,0.000000}%
\pgfsetfillcolor{currentfill}%
\pgfsetlinewidth{0.803000pt}%
\definecolor{currentstroke}{rgb}{0.000000,0.000000,0.000000}%
\pgfsetstrokecolor{currentstroke}%
\pgfsetdash{}{0pt}%
\pgfsys@defobject{currentmarker}{\pgfqpoint{0.000000in}{-0.048611in}}{\pgfqpoint{0.000000in}{0.000000in}}{%
\pgfpathmoveto{\pgfqpoint{0.000000in}{0.000000in}}%
\pgfpathlineto{\pgfqpoint{0.000000in}{-0.048611in}}%
\pgfusepath{stroke,fill}%
}%
\begin{pgfscope}%
\pgfsys@transformshift{5.894159in}{0.679475in}%
\pgfsys@useobject{currentmarker}{}%
\end{pgfscope}%
\end{pgfscope}%
\begin{pgfscope}%
\definecolor{textcolor}{rgb}{0.000000,0.000000,0.000000}%
\pgfsetstrokecolor{textcolor}%
\pgfsetfillcolor{textcolor}%
\pgftext[x=5.894159in,y=0.582253in,,top]{\color{textcolor}\rmfamily\fontsize{16.000000}{19.200000}\selectfont \(\displaystyle {0.8}\)}%
\end{pgfscope}%
\begin{pgfscope}%
\pgfsetbuttcap%
\pgfsetroundjoin%
\definecolor{currentfill}{rgb}{0.000000,0.000000,0.000000}%
\pgfsetfillcolor{currentfill}%
\pgfsetlinewidth{0.803000pt}%
\definecolor{currentstroke}{rgb}{0.000000,0.000000,0.000000}%
\pgfsetstrokecolor{currentstroke}%
\pgfsetdash{}{0pt}%
\pgfsys@defobject{currentmarker}{\pgfqpoint{0.000000in}{-0.048611in}}{\pgfqpoint{0.000000in}{0.000000in}}{%
\pgfpathmoveto{\pgfqpoint{0.000000in}{0.000000in}}%
\pgfpathlineto{\pgfqpoint{0.000000in}{-0.048611in}}%
\pgfusepath{stroke,fill}%
}%
\begin{pgfscope}%
\pgfsys@transformshift{7.086065in}{0.679475in}%
\pgfsys@useobject{currentmarker}{}%
\end{pgfscope}%
\end{pgfscope}%
\begin{pgfscope}%
\definecolor{textcolor}{rgb}{0.000000,0.000000,0.000000}%
\pgfsetstrokecolor{textcolor}%
\pgfsetfillcolor{textcolor}%
\pgftext[x=7.086065in,y=0.582253in,,top]{\color{textcolor}\rmfamily\fontsize{16.000000}{19.200000}\selectfont \(\displaystyle {1.0}\)}%
\end{pgfscope}%
\begin{pgfscope}%
\definecolor{textcolor}{rgb}{0.000000,0.000000,0.000000}%
\pgfsetstrokecolor{textcolor}%
\pgfsetfillcolor{textcolor}%
\pgftext[x=4.404277in,y=0.313349in,,top]{\color{textcolor}\rmfamily\fontsize{18.000000}{21.600000}\selectfont \(\displaystyle \lambda\)}%
\end{pgfscope}%
\begin{pgfscope}%
\pgfsetbuttcap%
\pgfsetroundjoin%
\definecolor{currentfill}{rgb}{0.000000,0.000000,0.000000}%
\pgfsetfillcolor{currentfill}%
\pgfsetlinewidth{0.803000pt}%
\definecolor{currentstroke}{rgb}{0.000000,0.000000,0.000000}%
\pgfsetstrokecolor{currentstroke}%
\pgfsetdash{}{0pt}%
\pgfsys@defobject{currentmarker}{\pgfqpoint{-0.048611in}{0.000000in}}{\pgfqpoint{-0.000000in}{0.000000in}}{%
\pgfpathmoveto{\pgfqpoint{-0.000000in}{0.000000in}}%
\pgfpathlineto{\pgfqpoint{-0.048611in}{0.000000in}}%
\pgfusepath{stroke,fill}%
}%
\begin{pgfscope}%
\pgfsys@transformshift{1.126536in}{1.000327in}%
\pgfsys@useobject{currentmarker}{}%
\end{pgfscope}%
\end{pgfscope}%
\begin{pgfscope}%
\definecolor{textcolor}{rgb}{0.000000,0.000000,0.000000}%
\pgfsetstrokecolor{textcolor}%
\pgfsetfillcolor{textcolor}%
\pgftext[x=0.368904in, y=0.916994in, left, base]{\color{textcolor}\rmfamily\fontsize{16.000000}{19.200000}\selectfont \(\displaystyle {115000}\)}%
\end{pgfscope}%
\begin{pgfscope}%
\pgfsetbuttcap%
\pgfsetroundjoin%
\definecolor{currentfill}{rgb}{0.000000,0.000000,0.000000}%
\pgfsetfillcolor{currentfill}%
\pgfsetlinewidth{0.803000pt}%
\definecolor{currentstroke}{rgb}{0.000000,0.000000,0.000000}%
\pgfsetstrokecolor{currentstroke}%
\pgfsetdash{}{0pt}%
\pgfsys@defobject{currentmarker}{\pgfqpoint{-0.048611in}{0.000000in}}{\pgfqpoint{-0.000000in}{0.000000in}}{%
\pgfpathmoveto{\pgfqpoint{-0.000000in}{0.000000in}}%
\pgfpathlineto{\pgfqpoint{-0.048611in}{0.000000in}}%
\pgfusepath{stroke,fill}%
}%
\begin{pgfscope}%
\pgfsys@transformshift{1.126536in}{1.430495in}%
\pgfsys@useobject{currentmarker}{}%
\end{pgfscope}%
\end{pgfscope}%
\begin{pgfscope}%
\definecolor{textcolor}{rgb}{0.000000,0.000000,0.000000}%
\pgfsetstrokecolor{textcolor}%
\pgfsetfillcolor{textcolor}%
\pgftext[x=0.368904in, y=1.347161in, left, base]{\color{textcolor}\rmfamily\fontsize{16.000000}{19.200000}\selectfont \(\displaystyle {120000}\)}%
\end{pgfscope}%
\begin{pgfscope}%
\pgfsetbuttcap%
\pgfsetroundjoin%
\definecolor{currentfill}{rgb}{0.000000,0.000000,0.000000}%
\pgfsetfillcolor{currentfill}%
\pgfsetlinewidth{0.803000pt}%
\definecolor{currentstroke}{rgb}{0.000000,0.000000,0.000000}%
\pgfsetstrokecolor{currentstroke}%
\pgfsetdash{}{0pt}%
\pgfsys@defobject{currentmarker}{\pgfqpoint{-0.048611in}{0.000000in}}{\pgfqpoint{-0.000000in}{0.000000in}}{%
\pgfpathmoveto{\pgfqpoint{-0.000000in}{0.000000in}}%
\pgfpathlineto{\pgfqpoint{-0.048611in}{0.000000in}}%
\pgfusepath{stroke,fill}%
}%
\begin{pgfscope}%
\pgfsys@transformshift{1.126536in}{1.860662in}%
\pgfsys@useobject{currentmarker}{}%
\end{pgfscope}%
\end{pgfscope}%
\begin{pgfscope}%
\definecolor{textcolor}{rgb}{0.000000,0.000000,0.000000}%
\pgfsetstrokecolor{textcolor}%
\pgfsetfillcolor{textcolor}%
\pgftext[x=0.368904in, y=1.777329in, left, base]{\color{textcolor}\rmfamily\fontsize{16.000000}{19.200000}\selectfont \(\displaystyle {125000}\)}%
\end{pgfscope}%
\begin{pgfscope}%
\pgfsetbuttcap%
\pgfsetroundjoin%
\definecolor{currentfill}{rgb}{0.000000,0.000000,0.000000}%
\pgfsetfillcolor{currentfill}%
\pgfsetlinewidth{0.803000pt}%
\definecolor{currentstroke}{rgb}{0.000000,0.000000,0.000000}%
\pgfsetstrokecolor{currentstroke}%
\pgfsetdash{}{0pt}%
\pgfsys@defobject{currentmarker}{\pgfqpoint{-0.048611in}{0.000000in}}{\pgfqpoint{-0.000000in}{0.000000in}}{%
\pgfpathmoveto{\pgfqpoint{-0.000000in}{0.000000in}}%
\pgfpathlineto{\pgfqpoint{-0.048611in}{0.000000in}}%
\pgfusepath{stroke,fill}%
}%
\begin{pgfscope}%
\pgfsys@transformshift{1.126536in}{2.290830in}%
\pgfsys@useobject{currentmarker}{}%
\end{pgfscope}%
\end{pgfscope}%
\begin{pgfscope}%
\definecolor{textcolor}{rgb}{0.000000,0.000000,0.000000}%
\pgfsetstrokecolor{textcolor}%
\pgfsetfillcolor{textcolor}%
\pgftext[x=0.368904in, y=2.207497in, left, base]{\color{textcolor}\rmfamily\fontsize{16.000000}{19.200000}\selectfont \(\displaystyle {130000}\)}%
\end{pgfscope}%
\begin{pgfscope}%
\pgfsetbuttcap%
\pgfsetroundjoin%
\definecolor{currentfill}{rgb}{0.000000,0.000000,0.000000}%
\pgfsetfillcolor{currentfill}%
\pgfsetlinewidth{0.803000pt}%
\definecolor{currentstroke}{rgb}{0.000000,0.000000,0.000000}%
\pgfsetstrokecolor{currentstroke}%
\pgfsetdash{}{0pt}%
\pgfsys@defobject{currentmarker}{\pgfqpoint{-0.048611in}{0.000000in}}{\pgfqpoint{-0.000000in}{0.000000in}}{%
\pgfpathmoveto{\pgfqpoint{-0.000000in}{0.000000in}}%
\pgfpathlineto{\pgfqpoint{-0.048611in}{0.000000in}}%
\pgfusepath{stroke,fill}%
}%
\begin{pgfscope}%
\pgfsys@transformshift{1.126536in}{2.720998in}%
\pgfsys@useobject{currentmarker}{}%
\end{pgfscope}%
\end{pgfscope}%
\begin{pgfscope}%
\definecolor{textcolor}{rgb}{0.000000,0.000000,0.000000}%
\pgfsetstrokecolor{textcolor}%
\pgfsetfillcolor{textcolor}%
\pgftext[x=0.368904in, y=2.637664in, left, base]{\color{textcolor}\rmfamily\fontsize{16.000000}{19.200000}\selectfont \(\displaystyle {135000}\)}%
\end{pgfscope}%
\begin{pgfscope}%
\definecolor{textcolor}{rgb}{0.000000,0.000000,0.000000}%
\pgfsetstrokecolor{textcolor}%
\pgfsetfillcolor{textcolor}%
\pgftext[x=0.313349in,y=1.907686in,,bottom,rotate=90.000000]{\color{textcolor}\rmfamily\fontsize{18.000000}{21.600000}\selectfont DKK}%
\end{pgfscope}%
\begin{pgfscope}%
\pgfpathrectangle{\pgfqpoint{1.126536in}{0.679475in}}{\pgfqpoint{6.555482in}{2.456422in}}%
\pgfusepath{clip}%
\pgfsetrectcap%
\pgfsetroundjoin%
\pgfsetlinewidth{1.505625pt}%
\definecolor{currentstroke}{rgb}{0.337255,0.705882,0.913725}%
\pgfsetstrokecolor{currentstroke}%
\pgfsetdash{}{0pt}%
\pgfpathmoveto{\pgfqpoint{1.722489in}{3.024241in}}%
\pgfpathlineto{\pgfqpoint{2.616418in}{2.687408in}}%
\pgfpathlineto{\pgfqpoint{4.106300in}{2.539098in}}%
\pgfpathlineto{\pgfqpoint{5.596183in}{2.319832in}}%
\pgfpathlineto{\pgfqpoint{7.086065in}{2.132221in}}%
\pgfusepath{stroke}%
\end{pgfscope}%
\begin{pgfscope}%
\pgfpathrectangle{\pgfqpoint{1.126536in}{0.679475in}}{\pgfqpoint{6.555482in}{2.456422in}}%
\pgfusepath{clip}%
\pgfsetbuttcap%
\pgfsetmiterjoin%
\definecolor{currentfill}{rgb}{0.337255,0.705882,0.913725}%
\pgfsetfillcolor{currentfill}%
\pgfsetlinewidth{0.752812pt}%
\definecolor{currentstroke}{rgb}{1.000000,1.000000,1.000000}%
\pgfsetstrokecolor{currentstroke}%
\pgfsetdash{}{0pt}%
\pgfsys@defobject{currentmarker}{\pgfqpoint{-0.041667in}{-0.041667in}}{\pgfqpoint{0.041667in}{0.041667in}}{%
\pgfpathmoveto{\pgfqpoint{-0.041667in}{-0.041667in}}%
\pgfpathlineto{\pgfqpoint{0.041667in}{-0.041667in}}%
\pgfpathlineto{\pgfqpoint{0.041667in}{0.041667in}}%
\pgfpathlineto{\pgfqpoint{-0.041667in}{0.041667in}}%
\pgfpathlineto{\pgfqpoint{-0.041667in}{-0.041667in}}%
\pgfpathclose%
\pgfusepath{stroke,fill}%
}%
\begin{pgfscope}%
\pgfsys@transformshift{1.722489in}{3.024241in}%
\pgfsys@useobject{currentmarker}{}%
\end{pgfscope}%
\begin{pgfscope}%
\pgfsys@transformshift{2.616418in}{2.687408in}%
\pgfsys@useobject{currentmarker}{}%
\end{pgfscope}%
\begin{pgfscope}%
\pgfsys@transformshift{4.106300in}{2.539098in}%
\pgfsys@useobject{currentmarker}{}%
\end{pgfscope}%
\begin{pgfscope}%
\pgfsys@transformshift{5.596183in}{2.319832in}%
\pgfsys@useobject{currentmarker}{}%
\end{pgfscope}%
\begin{pgfscope}%
\pgfsys@transformshift{7.086065in}{2.132221in}%
\pgfsys@useobject{currentmarker}{}%
\end{pgfscope}%
\end{pgfscope}%
\begin{pgfscope}%
\pgfpathrectangle{\pgfqpoint{1.126536in}{0.679475in}}{\pgfqpoint{6.555482in}{2.456422in}}%
\pgfusepath{clip}%
\pgfsetbuttcap%
\pgfsetroundjoin%
\pgfsetlinewidth{1.505625pt}%
\definecolor{currentstroke}{rgb}{0.337255,0.705882,0.913725}%
\pgfsetstrokecolor{currentstroke}%
\pgfsetstrokeopacity{0.900000}%
\pgfsetdash{{7.500000pt}{7.500000pt}}{0.000000pt}%
\pgfpathmoveto{\pgfqpoint{1.722489in}{1.503660in}}%
\pgfpathlineto{\pgfqpoint{2.616418in}{1.247907in}}%
\pgfpathlineto{\pgfqpoint{4.106300in}{1.020004in}}%
\pgfpathlineto{\pgfqpoint{5.596183in}{0.867450in}}%
\pgfpathlineto{\pgfqpoint{7.086065in}{0.791131in}}%
\pgfusepath{stroke}%
\end{pgfscope}%
\begin{pgfscope}%
\pgfpathrectangle{\pgfqpoint{1.126536in}{0.679475in}}{\pgfqpoint{6.555482in}{2.456422in}}%
\pgfusepath{clip}%
\pgfsetbuttcap%
\pgfsetmiterjoin%
\definecolor{currentfill}{rgb}{0.337255,0.705882,0.913725}%
\pgfsetfillcolor{currentfill}%
\pgfsetfillopacity{0.900000}%
\pgfsetlinewidth{0.752812pt}%
\definecolor{currentstroke}{rgb}{1.000000,1.000000,1.000000}%
\pgfsetstrokecolor{currentstroke}%
\pgfsetdash{}{0pt}%
\pgfsys@defobject{currentmarker}{\pgfqpoint{-0.058926in}{-0.058926in}}{\pgfqpoint{0.058926in}{0.058926in}}{%
\pgfpathmoveto{\pgfqpoint{0.000000in}{-0.058926in}}%
\pgfpathlineto{\pgfqpoint{0.058926in}{0.000000in}}%
\pgfpathlineto{\pgfqpoint{0.000000in}{0.058926in}}%
\pgfpathlineto{\pgfqpoint{-0.058926in}{0.000000in}}%
\pgfpathlineto{\pgfqpoint{0.000000in}{-0.058926in}}%
\pgfpathclose%
\pgfusepath{stroke,fill}%
}%
\begin{pgfscope}%
\pgfsys@transformshift{1.722489in}{1.503660in}%
\pgfsys@useobject{currentmarker}{}%
\end{pgfscope}%
\begin{pgfscope}%
\pgfsys@transformshift{2.616418in}{1.247907in}%
\pgfsys@useobject{currentmarker}{}%
\end{pgfscope}%
\begin{pgfscope}%
\pgfsys@transformshift{4.106300in}{1.020004in}%
\pgfsys@useobject{currentmarker}{}%
\end{pgfscope}%
\begin{pgfscope}%
\pgfsys@transformshift{5.596183in}{0.867450in}%
\pgfsys@useobject{currentmarker}{}%
\end{pgfscope}%
\begin{pgfscope}%
\pgfsys@transformshift{7.086065in}{0.791131in}%
\pgfsys@useobject{currentmarker}{}%
\end{pgfscope}%
\end{pgfscope}%
\begin{pgfscope}%
\pgfpathrectangle{\pgfqpoint{1.126536in}{0.679475in}}{\pgfqpoint{6.555482in}{2.456422in}}%
\pgfusepath{clip}%
\pgfsetbuttcap%
\pgfsetroundjoin%
\pgfsetlinewidth{1.505625pt}%
\definecolor{currentstroke}{rgb}{0.337255,0.705882,0.913725}%
\pgfsetstrokecolor{currentstroke}%
\pgfsetstrokeopacity{0.800000}%
\pgfsetdash{{1.500000pt}{4.500000pt}}{0.000000pt}%
\pgfpathmoveto{\pgfqpoint{1.722489in}{1.513247in}}%
\pgfpathlineto{\pgfqpoint{2.616418in}{1.522575in}}%
\pgfpathlineto{\pgfqpoint{4.106300in}{1.521395in}}%
\pgfpathlineto{\pgfqpoint{5.596183in}{1.446436in}}%
\pgfpathlineto{\pgfqpoint{7.086065in}{1.643674in}}%
\pgfusepath{stroke}%
\end{pgfscope}%
\begin{pgfscope}%
\pgfpathrectangle{\pgfqpoint{1.126536in}{0.679475in}}{\pgfqpoint{6.555482in}{2.456422in}}%
\pgfusepath{clip}%
\pgfsetbuttcap%
\pgfsetroundjoin%
\definecolor{currentfill}{rgb}{0.337255,0.705882,0.913725}%
\pgfsetfillcolor{currentfill}%
\pgfsetfillopacity{0.800000}%
\pgfsetlinewidth{0.752812pt}%
\definecolor{currentstroke}{rgb}{1.000000,1.000000,1.000000}%
\pgfsetstrokecolor{currentstroke}%
\pgfsetdash{}{0pt}%
\pgfsys@defobject{currentmarker}{\pgfqpoint{-0.041667in}{-0.041667in}}{\pgfqpoint{0.041667in}{0.041667in}}{%
\pgfpathmoveto{\pgfqpoint{0.000000in}{-0.041667in}}%
\pgfpathcurveto{\pgfqpoint{0.011050in}{-0.041667in}}{\pgfqpoint{0.021649in}{-0.037276in}}{\pgfqpoint{0.029463in}{-0.029463in}}%
\pgfpathcurveto{\pgfqpoint{0.037276in}{-0.021649in}}{\pgfqpoint{0.041667in}{-0.011050in}}{\pgfqpoint{0.041667in}{0.000000in}}%
\pgfpathcurveto{\pgfqpoint{0.041667in}{0.011050in}}{\pgfqpoint{0.037276in}{0.021649in}}{\pgfqpoint{0.029463in}{0.029463in}}%
\pgfpathcurveto{\pgfqpoint{0.021649in}{0.037276in}}{\pgfqpoint{0.011050in}{0.041667in}}{\pgfqpoint{0.000000in}{0.041667in}}%
\pgfpathcurveto{\pgfqpoint{-0.011050in}{0.041667in}}{\pgfqpoint{-0.021649in}{0.037276in}}{\pgfqpoint{-0.029463in}{0.029463in}}%
\pgfpathcurveto{\pgfqpoint{-0.037276in}{0.021649in}}{\pgfqpoint{-0.041667in}{0.011050in}}{\pgfqpoint{-0.041667in}{0.000000in}}%
\pgfpathcurveto{\pgfqpoint{-0.041667in}{-0.011050in}}{\pgfqpoint{-0.037276in}{-0.021649in}}{\pgfqpoint{-0.029463in}{-0.029463in}}%
\pgfpathcurveto{\pgfqpoint{-0.021649in}{-0.037276in}}{\pgfqpoint{-0.011050in}{-0.041667in}}{\pgfqpoint{0.000000in}{-0.041667in}}%
\pgfpathlineto{\pgfqpoint{0.000000in}{-0.041667in}}%
\pgfpathclose%
\pgfusepath{stroke,fill}%
}%
\begin{pgfscope}%
\pgfsys@transformshift{1.722489in}{1.513247in}%
\pgfsys@useobject{currentmarker}{}%
\end{pgfscope}%
\begin{pgfscope}%
\pgfsys@transformshift{2.616418in}{1.522575in}%
\pgfsys@useobject{currentmarker}{}%
\end{pgfscope}%
\begin{pgfscope}%
\pgfsys@transformshift{4.106300in}{1.521395in}%
\pgfsys@useobject{currentmarker}{}%
\end{pgfscope}%
\begin{pgfscope}%
\pgfsys@transformshift{5.596183in}{1.446436in}%
\pgfsys@useobject{currentmarker}{}%
\end{pgfscope}%
\begin{pgfscope}%
\pgfsys@transformshift{7.086065in}{1.643674in}%
\pgfsys@useobject{currentmarker}{}%
\end{pgfscope}%
\end{pgfscope}%
\begin{pgfscope}%
\pgfsetrectcap%
\pgfsetmiterjoin%
\pgfsetlinewidth{0.803000pt}%
\definecolor{currentstroke}{rgb}{0.000000,0.000000,0.000000}%
\pgfsetstrokecolor{currentstroke}%
\pgfsetdash{}{0pt}%
\pgfpathmoveto{\pgfqpoint{1.126536in}{0.679475in}}%
\pgfpathlineto{\pgfqpoint{1.126536in}{3.135897in}}%
\pgfusepath{stroke}%
\end{pgfscope}%
\begin{pgfscope}%
\pgfsetrectcap%
\pgfsetmiterjoin%
\pgfsetlinewidth{0.803000pt}%
\definecolor{currentstroke}{rgb}{0.000000,0.000000,0.000000}%
\pgfsetstrokecolor{currentstroke}%
\pgfsetdash{}{0pt}%
\pgfpathmoveto{\pgfqpoint{7.682018in}{0.679475in}}%
\pgfpathlineto{\pgfqpoint{7.682018in}{3.135897in}}%
\pgfusepath{stroke}%
\end{pgfscope}%
\begin{pgfscope}%
\pgfsetrectcap%
\pgfsetmiterjoin%
\pgfsetlinewidth{0.803000pt}%
\definecolor{currentstroke}{rgb}{0.000000,0.000000,0.000000}%
\pgfsetstrokecolor{currentstroke}%
\pgfsetdash{}{0pt}%
\pgfpathmoveto{\pgfqpoint{1.126536in}{0.679475in}}%
\pgfpathlineto{\pgfqpoint{7.682018in}{0.679475in}}%
\pgfusepath{stroke}%
\end{pgfscope}%
\begin{pgfscope}%
\pgfsetrectcap%
\pgfsetmiterjoin%
\pgfsetlinewidth{0.803000pt}%
\definecolor{currentstroke}{rgb}{0.000000,0.000000,0.000000}%
\pgfsetstrokecolor{currentstroke}%
\pgfsetdash{}{0pt}%
\pgfpathmoveto{\pgfqpoint{1.126536in}{3.135897in}}%
\pgfpathlineto{\pgfqpoint{7.682018in}{3.135897in}}%
\pgfusepath{stroke}%
\end{pgfscope}%
\begin{pgfscope}%
\pgfsetrectcap%
\pgfsetroundjoin%
\pgfsetlinewidth{1.505625pt}%
\definecolor{currentstroke}{rgb}{0.337255,0.705882,0.913725}%
\pgfsetstrokecolor{currentstroke}%
\pgfsetdash{}{0pt}%
\pgfpathmoveto{\pgfqpoint{2.251710in}{3.363810in}}%
\pgfpathlineto{\pgfqpoint{2.501710in}{3.363810in}}%
\pgfpathlineto{\pgfqpoint{2.751710in}{3.363810in}}%
\pgfusepath{stroke}%
\end{pgfscope}%
\begin{pgfscope}%
\pgfsetbuttcap%
\pgfsetmiterjoin%
\definecolor{currentfill}{rgb}{0.337255,0.705882,0.913725}%
\pgfsetfillcolor{currentfill}%
\pgfsetlinewidth{1.003750pt}%
\definecolor{currentstroke}{rgb}{0.337255,0.705882,0.913725}%
\pgfsetstrokecolor{currentstroke}%
\pgfsetdash{}{0pt}%
\pgfsys@defobject{currentmarker}{\pgfqpoint{-0.041667in}{-0.041667in}}{\pgfqpoint{0.041667in}{0.041667in}}{%
\pgfpathmoveto{\pgfqpoint{-0.041667in}{-0.041667in}}%
\pgfpathlineto{\pgfqpoint{0.041667in}{-0.041667in}}%
\pgfpathlineto{\pgfqpoint{0.041667in}{0.041667in}}%
\pgfpathlineto{\pgfqpoint{-0.041667in}{0.041667in}}%
\pgfpathlineto{\pgfqpoint{-0.041667in}{-0.041667in}}%
\pgfpathclose%
\pgfusepath{stroke,fill}%
}%
\begin{pgfscope}%
\pgfsys@transformshift{2.501710in}{3.363810in}%
\pgfsys@useobject{currentmarker}{}%
\end{pgfscope}%
\end{pgfscope}%
\begin{pgfscope}%
\definecolor{textcolor}{rgb}{0.000000,0.000000,0.000000}%
\pgfsetstrokecolor{textcolor}%
\pgfsetfillcolor{textcolor}%
\pgftext[x=2.876710in,y=3.276310in,left,base]{\color{textcolor}\rmfamily\fontsize{18.000000}{21.600000}\selectfont LPP\(\displaystyle ^{R,T}\)}%
\end{pgfscope}%
\begin{pgfscope}%
\pgfsetbuttcap%
\pgfsetroundjoin%
\pgfsetlinewidth{1.505625pt}%
\definecolor{currentstroke}{rgb}{0.337255,0.705882,0.913725}%
\pgfsetstrokecolor{currentstroke}%
\pgfsetstrokeopacity{0.900000}%
\pgfsetdash{{5.550000pt}{2.400000pt}}{0.000000pt}%
\pgfpathmoveto{\pgfqpoint{3.812661in}{3.363810in}}%
\pgfpathlineto{\pgfqpoint{4.062661in}{3.363810in}}%
\pgfpathlineto{\pgfqpoint{4.312661in}{3.363810in}}%
\pgfusepath{stroke}%
\end{pgfscope}%
\begin{pgfscope}%
\pgfsetbuttcap%
\pgfsetmiterjoin%
\definecolor{currentfill}{rgb}{0.337255,0.705882,0.913725}%
\pgfsetfillcolor{currentfill}%
\pgfsetfillopacity{0.900000}%
\pgfsetlinewidth{1.003750pt}%
\definecolor{currentstroke}{rgb}{0.337255,0.705882,0.913725}%
\pgfsetstrokecolor{currentstroke}%
\pgfsetstrokeopacity{0.900000}%
\pgfsetdash{}{0pt}%
\pgfsys@defobject{currentmarker}{\pgfqpoint{-0.058926in}{-0.058926in}}{\pgfqpoint{0.058926in}{0.058926in}}{%
\pgfpathmoveto{\pgfqpoint{0.000000in}{-0.058926in}}%
\pgfpathlineto{\pgfqpoint{0.058926in}{0.000000in}}%
\pgfpathlineto{\pgfqpoint{0.000000in}{0.058926in}}%
\pgfpathlineto{\pgfqpoint{-0.058926in}{0.000000in}}%
\pgfpathlineto{\pgfqpoint{0.000000in}{-0.058926in}}%
\pgfpathclose%
\pgfusepath{stroke,fill}%
}%
\begin{pgfscope}%
\pgfsys@transformshift{4.062661in}{3.363810in}%
\pgfsys@useobject{currentmarker}{}%
\end{pgfscope}%
\end{pgfscope}%
\begin{pgfscope}%
\definecolor{textcolor}{rgb}{0.000000,0.000000,0.000000}%
\pgfsetstrokecolor{textcolor}%
\pgfsetfillcolor{textcolor}%
\pgftext[x=4.437661in,y=3.276310in,left,base]{\color{textcolor}\rmfamily\fontsize{18.000000}{21.600000}\selectfont LPP\(\displaystyle ^{T}\)}%
\end{pgfscope}%
\begin{pgfscope}%
\pgfsetbuttcap%
\pgfsetroundjoin%
\pgfsetlinewidth{1.505625pt}%
\definecolor{currentstroke}{rgb}{0.337255,0.705882,0.913725}%
\pgfsetstrokecolor{currentstroke}%
\pgfsetstrokeopacity{0.800000}%
\pgfsetdash{{1.500000pt}{2.475000pt}}{0.000000pt}%
\pgfpathmoveto{\pgfqpoint{5.202692in}{3.363810in}}%
\pgfpathlineto{\pgfqpoint{5.452692in}{3.363810in}}%
\pgfpathlineto{\pgfqpoint{5.702692in}{3.363810in}}%
\pgfusepath{stroke}%
\end{pgfscope}%
\begin{pgfscope}%
\pgfsetbuttcap%
\pgfsetroundjoin%
\definecolor{currentfill}{rgb}{0.337255,0.705882,0.913725}%
\pgfsetfillcolor{currentfill}%
\pgfsetfillopacity{0.800000}%
\pgfsetlinewidth{1.003750pt}%
\definecolor{currentstroke}{rgb}{0.337255,0.705882,0.913725}%
\pgfsetstrokecolor{currentstroke}%
\pgfsetstrokeopacity{0.800000}%
\pgfsetdash{}{0pt}%
\pgfsys@defobject{currentmarker}{\pgfqpoint{-0.041667in}{-0.041667in}}{\pgfqpoint{0.041667in}{0.041667in}}{%
\pgfpathmoveto{\pgfqpoint{0.000000in}{-0.041667in}}%
\pgfpathcurveto{\pgfqpoint{0.011050in}{-0.041667in}}{\pgfqpoint{0.021649in}{-0.037276in}}{\pgfqpoint{0.029463in}{-0.029463in}}%
\pgfpathcurveto{\pgfqpoint{0.037276in}{-0.021649in}}{\pgfqpoint{0.041667in}{-0.011050in}}{\pgfqpoint{0.041667in}{0.000000in}}%
\pgfpathcurveto{\pgfqpoint{0.041667in}{0.011050in}}{\pgfqpoint{0.037276in}{0.021649in}}{\pgfqpoint{0.029463in}{0.029463in}}%
\pgfpathcurveto{\pgfqpoint{0.021649in}{0.037276in}}{\pgfqpoint{0.011050in}{0.041667in}}{\pgfqpoint{0.000000in}{0.041667in}}%
\pgfpathcurveto{\pgfqpoint{-0.011050in}{0.041667in}}{\pgfqpoint{-0.021649in}{0.037276in}}{\pgfqpoint{-0.029463in}{0.029463in}}%
\pgfpathcurveto{\pgfqpoint{-0.037276in}{0.021649in}}{\pgfqpoint{-0.041667in}{0.011050in}}{\pgfqpoint{-0.041667in}{0.000000in}}%
\pgfpathcurveto{\pgfqpoint{-0.041667in}{-0.011050in}}{\pgfqpoint{-0.037276in}{-0.021649in}}{\pgfqpoint{-0.029463in}{-0.029463in}}%
\pgfpathcurveto{\pgfqpoint{-0.021649in}{-0.037276in}}{\pgfqpoint{-0.011050in}{-0.041667in}}{\pgfqpoint{0.000000in}{-0.041667in}}%
\pgfpathlineto{\pgfqpoint{0.000000in}{-0.041667in}}%
\pgfpathclose%
\pgfusepath{stroke,fill}%
}%
\begin{pgfscope}%
\pgfsys@transformshift{5.452692in}{3.363810in}%
\pgfsys@useobject{currentmarker}{}%
\end{pgfscope}%
\end{pgfscope}%
\begin{pgfscope}%
\definecolor{textcolor}{rgb}{0.000000,0.000000,0.000000}%
\pgfsetstrokecolor{textcolor}%
\pgfsetfillcolor{textcolor}%
\pgftext[x=5.827692in,y=3.276310in,left,base]{\color{textcolor}\rmfamily\fontsize{18.000000}{21.600000}\selectfont LPP\(\displaystyle ^{T, P}\)}%
\end{pgfscope}%
\end{pgfpicture}%
\makeatother%
\endgroup%

%% file: figures/DemandPlots/TightLayout/DE_PI_rescaled.pgf
\begingroup%
\makeatletter%
\begin{pgfpicture}%
\pgfpathrectangle{\pgfpointorigin}{\pgfqpoint{7.782018in}{3.651310in}}%
\pgfusepath{use as bounding box, clip}%
\begin{pgfscope}%
\pgfsetbuttcap%
\pgfsetmiterjoin%
\definecolor{currentfill}{rgb}{1.000000,1.000000,1.000000}%
\pgfsetfillcolor{currentfill}%
\pgfsetlinewidth{0.000000pt}%
\definecolor{currentstroke}{rgb}{1.000000,1.000000,1.000000}%
\pgfsetstrokecolor{currentstroke}%
\pgfsetdash{}{0pt}%
\pgfpathmoveto{\pgfqpoint{0.000000in}{0.000000in}}%
\pgfpathlineto{\pgfqpoint{7.782018in}{0.000000in}}%
\pgfpathlineto{\pgfqpoint{7.782018in}{3.651310in}}%
\pgfpathlineto{\pgfqpoint{0.000000in}{3.651310in}}%
\pgfpathlineto{\pgfqpoint{0.000000in}{0.000000in}}%
\pgfpathclose%
\pgfusepath{fill}%
\end{pgfscope}%
\begin{pgfscope}%
\pgfsetbuttcap%
\pgfsetmiterjoin%
\definecolor{currentfill}{rgb}{1.000000,1.000000,1.000000}%
\pgfsetfillcolor{currentfill}%
\pgfsetlinewidth{0.000000pt}%
\definecolor{currentstroke}{rgb}{0.000000,0.000000,0.000000}%
\pgfsetstrokecolor{currentstroke}%
\pgfsetstrokeopacity{0.000000}%
\pgfsetdash{}{0pt}%
\pgfpathmoveto{\pgfqpoint{1.126536in}{0.679475in}}%
\pgfpathlineto{\pgfqpoint{7.682018in}{0.679475in}}%
\pgfpathlineto{\pgfqpoint{7.682018in}{3.135897in}}%
\pgfpathlineto{\pgfqpoint{1.126536in}{3.135897in}}%
\pgfpathlineto{\pgfqpoint{1.126536in}{0.679475in}}%
\pgfpathclose%
\pgfusepath{fill}%
\end{pgfscope}%
\begin{pgfscope}%
\pgfsetbuttcap%
\pgfsetroundjoin%
\definecolor{currentfill}{rgb}{0.000000,0.000000,0.000000}%
\pgfsetfillcolor{currentfill}%
\pgfsetlinewidth{0.803000pt}%
\definecolor{currentstroke}{rgb}{0.000000,0.000000,0.000000}%
\pgfsetstrokecolor{currentstroke}%
\pgfsetdash{}{0pt}%
\pgfsys@defobject{currentmarker}{\pgfqpoint{0.000000in}{-0.048611in}}{\pgfqpoint{0.000000in}{0.000000in}}{%
\pgfpathmoveto{\pgfqpoint{0.000000in}{0.000000in}}%
\pgfpathlineto{\pgfqpoint{0.000000in}{-0.048611in}}%
\pgfusepath{stroke,fill}%
}%
\begin{pgfscope}%
\pgfsys@transformshift{1.126536in}{0.679475in}%
\pgfsys@useobject{currentmarker}{}%
\end{pgfscope}%
\end{pgfscope}%
\begin{pgfscope}%
\definecolor{textcolor}{rgb}{0.000000,0.000000,0.000000}%
\pgfsetstrokecolor{textcolor}%
\pgfsetfillcolor{textcolor}%
\pgftext[x=1.126536in,y=0.582253in,,top]{\color{textcolor}\rmfamily\fontsize{16.000000}{19.200000}\selectfont \(\displaystyle {0.0}\)}%
\end{pgfscope}%
\begin{pgfscope}%
\pgfsetbuttcap%
\pgfsetroundjoin%
\definecolor{currentfill}{rgb}{0.000000,0.000000,0.000000}%
\pgfsetfillcolor{currentfill}%
\pgfsetlinewidth{0.803000pt}%
\definecolor{currentstroke}{rgb}{0.000000,0.000000,0.000000}%
\pgfsetstrokecolor{currentstroke}%
\pgfsetdash{}{0pt}%
\pgfsys@defobject{currentmarker}{\pgfqpoint{0.000000in}{-0.048611in}}{\pgfqpoint{0.000000in}{0.000000in}}{%
\pgfpathmoveto{\pgfqpoint{0.000000in}{0.000000in}}%
\pgfpathlineto{\pgfqpoint{0.000000in}{-0.048611in}}%
\pgfusepath{stroke,fill}%
}%
\begin{pgfscope}%
\pgfsys@transformshift{2.318442in}{0.679475in}%
\pgfsys@useobject{currentmarker}{}%
\end{pgfscope}%
\end{pgfscope}%
\begin{pgfscope}%
\definecolor{textcolor}{rgb}{0.000000,0.000000,0.000000}%
\pgfsetstrokecolor{textcolor}%
\pgfsetfillcolor{textcolor}%
\pgftext[x=2.318442in,y=0.582253in,,top]{\color{textcolor}\rmfamily\fontsize{16.000000}{19.200000}\selectfont \(\displaystyle {0.2}\)}%
\end{pgfscope}%
\begin{pgfscope}%
\pgfsetbuttcap%
\pgfsetroundjoin%
\definecolor{currentfill}{rgb}{0.000000,0.000000,0.000000}%
\pgfsetfillcolor{currentfill}%
\pgfsetlinewidth{0.803000pt}%
\definecolor{currentstroke}{rgb}{0.000000,0.000000,0.000000}%
\pgfsetstrokecolor{currentstroke}%
\pgfsetdash{}{0pt}%
\pgfsys@defobject{currentmarker}{\pgfqpoint{0.000000in}{-0.048611in}}{\pgfqpoint{0.000000in}{0.000000in}}{%
\pgfpathmoveto{\pgfqpoint{0.000000in}{0.000000in}}%
\pgfpathlineto{\pgfqpoint{0.000000in}{-0.048611in}}%
\pgfusepath{stroke,fill}%
}%
\begin{pgfscope}%
\pgfsys@transformshift{3.510347in}{0.679475in}%
\pgfsys@useobject{currentmarker}{}%
\end{pgfscope}%
\end{pgfscope}%
\begin{pgfscope}%
\definecolor{textcolor}{rgb}{0.000000,0.000000,0.000000}%
\pgfsetstrokecolor{textcolor}%
\pgfsetfillcolor{textcolor}%
\pgftext[x=3.510347in,y=0.582253in,,top]{\color{textcolor}\rmfamily\fontsize{16.000000}{19.200000}\selectfont \(\displaystyle {0.4}\)}%
\end{pgfscope}%
\begin{pgfscope}%
\pgfsetbuttcap%
\pgfsetroundjoin%
\definecolor{currentfill}{rgb}{0.000000,0.000000,0.000000}%
\pgfsetfillcolor{currentfill}%
\pgfsetlinewidth{0.803000pt}%
\definecolor{currentstroke}{rgb}{0.000000,0.000000,0.000000}%
\pgfsetstrokecolor{currentstroke}%
\pgfsetdash{}{0pt}%
\pgfsys@defobject{currentmarker}{\pgfqpoint{0.000000in}{-0.048611in}}{\pgfqpoint{0.000000in}{0.000000in}}{%
\pgfpathmoveto{\pgfqpoint{0.000000in}{0.000000in}}%
\pgfpathlineto{\pgfqpoint{0.000000in}{-0.048611in}}%
\pgfusepath{stroke,fill}%
}%
\begin{pgfscope}%
\pgfsys@transformshift{4.702253in}{0.679475in}%
\pgfsys@useobject{currentmarker}{}%
\end{pgfscope}%
\end{pgfscope}%
\begin{pgfscope}%
\definecolor{textcolor}{rgb}{0.000000,0.000000,0.000000}%
\pgfsetstrokecolor{textcolor}%
\pgfsetfillcolor{textcolor}%
\pgftext[x=4.702253in,y=0.582253in,,top]{\color{textcolor}\rmfamily\fontsize{16.000000}{19.200000}\selectfont \(\displaystyle {0.6}\)}%
\end{pgfscope}%
\begin{pgfscope}%
\pgfsetbuttcap%
\pgfsetroundjoin%
\definecolor{currentfill}{rgb}{0.000000,0.000000,0.000000}%
\pgfsetfillcolor{currentfill}%
\pgfsetlinewidth{0.803000pt}%
\definecolor{currentstroke}{rgb}{0.000000,0.000000,0.000000}%
\pgfsetstrokecolor{currentstroke}%
\pgfsetdash{}{0pt}%
\pgfsys@defobject{currentmarker}{\pgfqpoint{0.000000in}{-0.048611in}}{\pgfqpoint{0.000000in}{0.000000in}}{%
\pgfpathmoveto{\pgfqpoint{0.000000in}{0.000000in}}%
\pgfpathlineto{\pgfqpoint{0.000000in}{-0.048611in}}%
\pgfusepath{stroke,fill}%
}%
\begin{pgfscope}%
\pgfsys@transformshift{5.894159in}{0.679475in}%
\pgfsys@useobject{currentmarker}{}%
\end{pgfscope}%
\end{pgfscope}%
\begin{pgfscope}%
\definecolor{textcolor}{rgb}{0.000000,0.000000,0.000000}%
\pgfsetstrokecolor{textcolor}%
\pgfsetfillcolor{textcolor}%
\pgftext[x=5.894159in,y=0.582253in,,top]{\color{textcolor}\rmfamily\fontsize{16.000000}{19.200000}\selectfont \(\displaystyle {0.8}\)}%
\end{pgfscope}%
\begin{pgfscope}%
\pgfsetbuttcap%
\pgfsetroundjoin%
\definecolor{currentfill}{rgb}{0.000000,0.000000,0.000000}%
\pgfsetfillcolor{currentfill}%
\pgfsetlinewidth{0.803000pt}%
\definecolor{currentstroke}{rgb}{0.000000,0.000000,0.000000}%
\pgfsetstrokecolor{currentstroke}%
\pgfsetdash{}{0pt}%
\pgfsys@defobject{currentmarker}{\pgfqpoint{0.000000in}{-0.048611in}}{\pgfqpoint{0.000000in}{0.000000in}}{%
\pgfpathmoveto{\pgfqpoint{0.000000in}{0.000000in}}%
\pgfpathlineto{\pgfqpoint{0.000000in}{-0.048611in}}%
\pgfusepath{stroke,fill}%
}%
\begin{pgfscope}%
\pgfsys@transformshift{7.086065in}{0.679475in}%
\pgfsys@useobject{currentmarker}{}%
\end{pgfscope}%
\end{pgfscope}%
\begin{pgfscope}%
\definecolor{textcolor}{rgb}{0.000000,0.000000,0.000000}%
\pgfsetstrokecolor{textcolor}%
\pgfsetfillcolor{textcolor}%
\pgftext[x=7.086065in,y=0.582253in,,top]{\color{textcolor}\rmfamily\fontsize{16.000000}{19.200000}\selectfont \(\displaystyle {1.0}\)}%
\end{pgfscope}%
\begin{pgfscope}%
\definecolor{textcolor}{rgb}{0.000000,0.000000,0.000000}%
\pgfsetstrokecolor{textcolor}%
\pgfsetfillcolor{textcolor}%
\pgftext[x=4.404277in,y=0.313349in,,top]{\color{textcolor}\rmfamily\fontsize{18.000000}{21.600000}\selectfont \(\displaystyle \lambda\)}%
\end{pgfscope}%
\begin{pgfscope}%
\pgfsetbuttcap%
\pgfsetroundjoin%
\definecolor{currentfill}{rgb}{0.000000,0.000000,0.000000}%
\pgfsetfillcolor{currentfill}%
\pgfsetlinewidth{0.803000pt}%
\definecolor{currentstroke}{rgb}{0.000000,0.000000,0.000000}%
\pgfsetstrokecolor{currentstroke}%
\pgfsetdash{}{0pt}%
\pgfsys@defobject{currentmarker}{\pgfqpoint{-0.048611in}{0.000000in}}{\pgfqpoint{-0.000000in}{0.000000in}}{%
\pgfpathmoveto{\pgfqpoint{-0.000000in}{0.000000in}}%
\pgfpathlineto{\pgfqpoint{-0.048611in}{0.000000in}}%
\pgfusepath{stroke,fill}%
}%
\begin{pgfscope}%
\pgfsys@transformshift{1.126536in}{0.832087in}%
\pgfsys@useobject{currentmarker}{}%
\end{pgfscope}%
\end{pgfscope}%
\begin{pgfscope}%
\definecolor{textcolor}{rgb}{0.000000,0.000000,0.000000}%
\pgfsetstrokecolor{textcolor}%
\pgfsetfillcolor{textcolor}%
\pgftext[x=0.368904in, y=0.748754in, left, base]{\color{textcolor}\rmfamily\fontsize{16.000000}{19.200000}\selectfont \(\displaystyle {100000}\)}%
\end{pgfscope}%
\begin{pgfscope}%
\pgfsetbuttcap%
\pgfsetroundjoin%
\definecolor{currentfill}{rgb}{0.000000,0.000000,0.000000}%
\pgfsetfillcolor{currentfill}%
\pgfsetlinewidth{0.803000pt}%
\definecolor{currentstroke}{rgb}{0.000000,0.000000,0.000000}%
\pgfsetstrokecolor{currentstroke}%
\pgfsetdash{}{0pt}%
\pgfsys@defobject{currentmarker}{\pgfqpoint{-0.048611in}{0.000000in}}{\pgfqpoint{-0.000000in}{0.000000in}}{%
\pgfpathmoveto{\pgfqpoint{-0.000000in}{0.000000in}}%
\pgfpathlineto{\pgfqpoint{-0.048611in}{0.000000in}}%
\pgfusepath{stroke,fill}%
}%
\begin{pgfscope}%
\pgfsys@transformshift{1.126536in}{1.566266in}%
\pgfsys@useobject{currentmarker}{}%
\end{pgfscope}%
\end{pgfscope}%
\begin{pgfscope}%
\definecolor{textcolor}{rgb}{0.000000,0.000000,0.000000}%
\pgfsetstrokecolor{textcolor}%
\pgfsetfillcolor{textcolor}%
\pgftext[x=0.368904in, y=1.482933in, left, base]{\color{textcolor}\rmfamily\fontsize{16.000000}{19.200000}\selectfont \(\displaystyle {110000}\)}%
\end{pgfscope}%
\begin{pgfscope}%
\pgfsetbuttcap%
\pgfsetroundjoin%
\definecolor{currentfill}{rgb}{0.000000,0.000000,0.000000}%
\pgfsetfillcolor{currentfill}%
\pgfsetlinewidth{0.803000pt}%
\definecolor{currentstroke}{rgb}{0.000000,0.000000,0.000000}%
\pgfsetstrokecolor{currentstroke}%
\pgfsetdash{}{0pt}%
\pgfsys@defobject{currentmarker}{\pgfqpoint{-0.048611in}{0.000000in}}{\pgfqpoint{-0.000000in}{0.000000in}}{%
\pgfpathmoveto{\pgfqpoint{-0.000000in}{0.000000in}}%
\pgfpathlineto{\pgfqpoint{-0.048611in}{0.000000in}}%
\pgfusepath{stroke,fill}%
}%
\begin{pgfscope}%
\pgfsys@transformshift{1.126536in}{2.300446in}%
\pgfsys@useobject{currentmarker}{}%
\end{pgfscope}%
\end{pgfscope}%
\begin{pgfscope}%
\definecolor{textcolor}{rgb}{0.000000,0.000000,0.000000}%
\pgfsetstrokecolor{textcolor}%
\pgfsetfillcolor{textcolor}%
\pgftext[x=0.368904in, y=2.217112in, left, base]{\color{textcolor}\rmfamily\fontsize{16.000000}{19.200000}\selectfont \(\displaystyle {120000}\)}%
\end{pgfscope}%
\begin{pgfscope}%
\pgfsetbuttcap%
\pgfsetroundjoin%
\definecolor{currentfill}{rgb}{0.000000,0.000000,0.000000}%
\pgfsetfillcolor{currentfill}%
\pgfsetlinewidth{0.803000pt}%
\definecolor{currentstroke}{rgb}{0.000000,0.000000,0.000000}%
\pgfsetstrokecolor{currentstroke}%
\pgfsetdash{}{0pt}%
\pgfsys@defobject{currentmarker}{\pgfqpoint{-0.048611in}{0.000000in}}{\pgfqpoint{-0.000000in}{0.000000in}}{%
\pgfpathmoveto{\pgfqpoint{-0.000000in}{0.000000in}}%
\pgfpathlineto{\pgfqpoint{-0.048611in}{0.000000in}}%
\pgfusepath{stroke,fill}%
}%
\begin{pgfscope}%
\pgfsys@transformshift{1.126536in}{3.034625in}%
\pgfsys@useobject{currentmarker}{}%
\end{pgfscope}%
\end{pgfscope}%
\begin{pgfscope}%
\definecolor{textcolor}{rgb}{0.000000,0.000000,0.000000}%
\pgfsetstrokecolor{textcolor}%
\pgfsetfillcolor{textcolor}%
\pgftext[x=0.368904in, y=2.951291in, left, base]{\color{textcolor}\rmfamily\fontsize{16.000000}{19.200000}\selectfont \(\displaystyle {130000}\)}%
\end{pgfscope}%
\begin{pgfscope}%
\definecolor{textcolor}{rgb}{0.000000,0.000000,0.000000}%
\pgfsetstrokecolor{textcolor}%
\pgfsetfillcolor{textcolor}%
\pgftext[x=0.313349in,y=1.907686in,,bottom,rotate=90.000000]{\color{textcolor}\rmfamily\fontsize{18.000000}{21.600000}\selectfont DKK}%
\end{pgfscope}%
\begin{pgfscope}%
\pgfpathrectangle{\pgfqpoint{1.126536in}{0.679475in}}{\pgfqpoint{6.555482in}{2.456422in}}%
\pgfusepath{clip}%
\pgfsetrectcap%
\pgfsetroundjoin%
\pgfsetlinewidth{1.505625pt}%
\definecolor{currentstroke}{rgb}{0.580392,0.403922,0.741176}%
\pgfsetstrokecolor{currentstroke}%
\pgfsetdash{}{0pt}%
\pgfpathmoveto{\pgfqpoint{1.722489in}{2.513148in}}%
\pgfpathlineto{\pgfqpoint{2.616418in}{2.345713in}}%
\pgfpathlineto{\pgfqpoint{4.106300in}{2.131626in}}%
\pgfpathlineto{\pgfqpoint{5.596183in}{1.789695in}}%
\pgfpathlineto{\pgfqpoint{7.086065in}{1.339426in}}%
\pgfusepath{stroke}%
\end{pgfscope}%
\begin{pgfscope}%
\pgfpathrectangle{\pgfqpoint{1.126536in}{0.679475in}}{\pgfqpoint{6.555482in}{2.456422in}}%
\pgfusepath{clip}%
\pgfsetbuttcap%
\pgfsetmiterjoin%
\definecolor{currentfill}{rgb}{0.580392,0.403922,0.741176}%
\pgfsetfillcolor{currentfill}%
\pgfsetlinewidth{0.752812pt}%
\definecolor{currentstroke}{rgb}{1.000000,1.000000,1.000000}%
\pgfsetstrokecolor{currentstroke}%
\pgfsetdash{}{0pt}%
\pgfsys@defobject{currentmarker}{\pgfqpoint{-0.041667in}{-0.041667in}}{\pgfqpoint{0.041667in}{0.041667in}}{%
\pgfpathmoveto{\pgfqpoint{-0.041667in}{-0.041667in}}%
\pgfpathlineto{\pgfqpoint{0.041667in}{-0.041667in}}%
\pgfpathlineto{\pgfqpoint{0.041667in}{0.041667in}}%
\pgfpathlineto{\pgfqpoint{-0.041667in}{0.041667in}}%
\pgfpathlineto{\pgfqpoint{-0.041667in}{-0.041667in}}%
\pgfpathclose%
\pgfusepath{stroke,fill}%
}%
\begin{pgfscope}%
\pgfsys@transformshift{1.722489in}{2.513148in}%
\pgfsys@useobject{currentmarker}{}%
\end{pgfscope}%
\begin{pgfscope}%
\pgfsys@transformshift{2.616418in}{2.345713in}%
\pgfsys@useobject{currentmarker}{}%
\end{pgfscope}%
\begin{pgfscope}%
\pgfsys@transformshift{4.106300in}{2.131626in}%
\pgfsys@useobject{currentmarker}{}%
\end{pgfscope}%
\begin{pgfscope}%
\pgfsys@transformshift{5.596183in}{1.789695in}%
\pgfsys@useobject{currentmarker}{}%
\end{pgfscope}%
\begin{pgfscope}%
\pgfsys@transformshift{7.086065in}{1.339426in}%
\pgfsys@useobject{currentmarker}{}%
\end{pgfscope}%
\end{pgfscope}%
\begin{pgfscope}%
\pgfpathrectangle{\pgfqpoint{1.126536in}{0.679475in}}{\pgfqpoint{6.555482in}{2.456422in}}%
\pgfusepath{clip}%
\pgfsetbuttcap%
\pgfsetroundjoin%
\pgfsetlinewidth{1.505625pt}%
\definecolor{currentstroke}{rgb}{0.580392,0.403922,0.741176}%
\pgfsetstrokecolor{currentstroke}%
\pgfsetstrokeopacity{0.900000}%
\pgfsetdash{{7.500000pt}{7.500000pt}}{0.000000pt}%
\pgfpathmoveto{\pgfqpoint{1.722489in}{3.024241in}}%
\pgfpathlineto{\pgfqpoint{2.616418in}{2.767342in}}%
\pgfpathlineto{\pgfqpoint{4.106300in}{2.427563in}}%
\pgfpathlineto{\pgfqpoint{5.596183in}{2.096456in}}%
\pgfpathlineto{\pgfqpoint{7.086065in}{1.695650in}}%
\pgfusepath{stroke}%
\end{pgfscope}%
\begin{pgfscope}%
\pgfpathrectangle{\pgfqpoint{1.126536in}{0.679475in}}{\pgfqpoint{6.555482in}{2.456422in}}%
\pgfusepath{clip}%
\pgfsetbuttcap%
\pgfsetmiterjoin%
\definecolor{currentfill}{rgb}{0.580392,0.403922,0.741176}%
\pgfsetfillcolor{currentfill}%
\pgfsetfillopacity{0.900000}%
\pgfsetlinewidth{0.752812pt}%
\definecolor{currentstroke}{rgb}{1.000000,1.000000,1.000000}%
\pgfsetstrokecolor{currentstroke}%
\pgfsetdash{}{0pt}%
\pgfsys@defobject{currentmarker}{\pgfqpoint{-0.058926in}{-0.058926in}}{\pgfqpoint{0.058926in}{0.058926in}}{%
\pgfpathmoveto{\pgfqpoint{0.000000in}{-0.058926in}}%
\pgfpathlineto{\pgfqpoint{0.058926in}{0.000000in}}%
\pgfpathlineto{\pgfqpoint{0.000000in}{0.058926in}}%
\pgfpathlineto{\pgfqpoint{-0.058926in}{0.000000in}}%
\pgfpathlineto{\pgfqpoint{0.000000in}{-0.058926in}}%
\pgfpathclose%
\pgfusepath{stroke,fill}%
}%
\begin{pgfscope}%
\pgfsys@transformshift{1.722489in}{3.024241in}%
\pgfsys@useobject{currentmarker}{}%
\end{pgfscope}%
\begin{pgfscope}%
\pgfsys@transformshift{2.616418in}{2.767342in}%
\pgfsys@useobject{currentmarker}{}%
\end{pgfscope}%
\begin{pgfscope}%
\pgfsys@transformshift{4.106300in}{2.427563in}%
\pgfsys@useobject{currentmarker}{}%
\end{pgfscope}%
\begin{pgfscope}%
\pgfsys@transformshift{5.596183in}{2.096456in}%
\pgfsys@useobject{currentmarker}{}%
\end{pgfscope}%
\begin{pgfscope}%
\pgfsys@transformshift{7.086065in}{1.695650in}%
\pgfsys@useobject{currentmarker}{}%
\end{pgfscope}%
\end{pgfscope}%
\begin{pgfscope}%
\pgfpathrectangle{\pgfqpoint{1.126536in}{0.679475in}}{\pgfqpoint{6.555482in}{2.456422in}}%
\pgfusepath{clip}%
\pgfsetbuttcap%
\pgfsetroundjoin%
\pgfsetlinewidth{1.505625pt}%
\definecolor{currentstroke}{rgb}{0.580392,0.403922,0.741176}%
\pgfsetstrokecolor{currentstroke}%
\pgfsetstrokeopacity{0.800000}%
\pgfsetdash{{1.500000pt}{4.500000pt}}{0.000000pt}%
\pgfpathmoveto{\pgfqpoint{1.722489in}{1.283188in}}%
\pgfpathlineto{\pgfqpoint{2.616418in}{1.250276in}}%
\pgfpathlineto{\pgfqpoint{4.106300in}{1.169264in}}%
\pgfpathlineto{\pgfqpoint{5.596183in}{1.004473in}}%
\pgfpathlineto{\pgfqpoint{7.086065in}{0.791131in}}%
\pgfusepath{stroke}%
\end{pgfscope}%
\begin{pgfscope}%
\pgfpathrectangle{\pgfqpoint{1.126536in}{0.679475in}}{\pgfqpoint{6.555482in}{2.456422in}}%
\pgfusepath{clip}%
\pgfsetbuttcap%
\pgfsetroundjoin%
\definecolor{currentfill}{rgb}{0.580392,0.403922,0.741176}%
\pgfsetfillcolor{currentfill}%
\pgfsetfillopacity{0.800000}%
\pgfsetlinewidth{0.752812pt}%
\definecolor{currentstroke}{rgb}{1.000000,1.000000,1.000000}%
\pgfsetstrokecolor{currentstroke}%
\pgfsetdash{}{0pt}%
\pgfsys@defobject{currentmarker}{\pgfqpoint{-0.041667in}{-0.041667in}}{\pgfqpoint{0.041667in}{0.041667in}}{%
\pgfpathmoveto{\pgfqpoint{0.000000in}{-0.041667in}}%
\pgfpathcurveto{\pgfqpoint{0.011050in}{-0.041667in}}{\pgfqpoint{0.021649in}{-0.037276in}}{\pgfqpoint{0.029463in}{-0.029463in}}%
\pgfpathcurveto{\pgfqpoint{0.037276in}{-0.021649in}}{\pgfqpoint{0.041667in}{-0.011050in}}{\pgfqpoint{0.041667in}{0.000000in}}%
\pgfpathcurveto{\pgfqpoint{0.041667in}{0.011050in}}{\pgfqpoint{0.037276in}{0.021649in}}{\pgfqpoint{0.029463in}{0.029463in}}%
\pgfpathcurveto{\pgfqpoint{0.021649in}{0.037276in}}{\pgfqpoint{0.011050in}{0.041667in}}{\pgfqpoint{0.000000in}{0.041667in}}%
\pgfpathcurveto{\pgfqpoint{-0.011050in}{0.041667in}}{\pgfqpoint{-0.021649in}{0.037276in}}{\pgfqpoint{-0.029463in}{0.029463in}}%
\pgfpathcurveto{\pgfqpoint{-0.037276in}{0.021649in}}{\pgfqpoint{-0.041667in}{0.011050in}}{\pgfqpoint{-0.041667in}{0.000000in}}%
\pgfpathcurveto{\pgfqpoint{-0.041667in}{-0.011050in}}{\pgfqpoint{-0.037276in}{-0.021649in}}{\pgfqpoint{-0.029463in}{-0.029463in}}%
\pgfpathcurveto{\pgfqpoint{-0.021649in}{-0.037276in}}{\pgfqpoint{-0.011050in}{-0.041667in}}{\pgfqpoint{0.000000in}{-0.041667in}}%
\pgfpathlineto{\pgfqpoint{0.000000in}{-0.041667in}}%
\pgfpathclose%
\pgfusepath{stroke,fill}%
}%
\begin{pgfscope}%
\pgfsys@transformshift{1.722489in}{1.283188in}%
\pgfsys@useobject{currentmarker}{}%
\end{pgfscope}%
\begin{pgfscope}%
\pgfsys@transformshift{2.616418in}{1.250276in}%
\pgfsys@useobject{currentmarker}{}%
\end{pgfscope}%
\begin{pgfscope}%
\pgfsys@transformshift{4.106300in}{1.169264in}%
\pgfsys@useobject{currentmarker}{}%
\end{pgfscope}%
\begin{pgfscope}%
\pgfsys@transformshift{5.596183in}{1.004473in}%
\pgfsys@useobject{currentmarker}{}%
\end{pgfscope}%
\begin{pgfscope}%
\pgfsys@transformshift{7.086065in}{0.791131in}%
\pgfsys@useobject{currentmarker}{}%
\end{pgfscope}%
\end{pgfscope}%
\begin{pgfscope}%
\pgfsetrectcap%
\pgfsetmiterjoin%
\pgfsetlinewidth{0.803000pt}%
\definecolor{currentstroke}{rgb}{0.000000,0.000000,0.000000}%
\pgfsetstrokecolor{currentstroke}%
\pgfsetdash{}{0pt}%
\pgfpathmoveto{\pgfqpoint{1.126536in}{0.679475in}}%
\pgfpathlineto{\pgfqpoint{1.126536in}{3.135897in}}%
\pgfusepath{stroke}%
\end{pgfscope}%
\begin{pgfscope}%
\pgfsetrectcap%
\pgfsetmiterjoin%
\pgfsetlinewidth{0.803000pt}%
\definecolor{currentstroke}{rgb}{0.000000,0.000000,0.000000}%
\pgfsetstrokecolor{currentstroke}%
\pgfsetdash{}{0pt}%
\pgfpathmoveto{\pgfqpoint{7.682018in}{0.679475in}}%
\pgfpathlineto{\pgfqpoint{7.682018in}{3.135897in}}%
\pgfusepath{stroke}%
\end{pgfscope}%
\begin{pgfscope}%
\pgfsetrectcap%
\pgfsetmiterjoin%
\pgfsetlinewidth{0.803000pt}%
\definecolor{currentstroke}{rgb}{0.000000,0.000000,0.000000}%
\pgfsetstrokecolor{currentstroke}%
\pgfsetdash{}{0pt}%
\pgfpathmoveto{\pgfqpoint{1.126536in}{0.679475in}}%
\pgfpathlineto{\pgfqpoint{7.682018in}{0.679475in}}%
\pgfusepath{stroke}%
\end{pgfscope}%
\begin{pgfscope}%
\pgfsetrectcap%
\pgfsetmiterjoin%
\pgfsetlinewidth{0.803000pt}%
\definecolor{currentstroke}{rgb}{0.000000,0.000000,0.000000}%
\pgfsetstrokecolor{currentstroke}%
\pgfsetdash{}{0pt}%
\pgfpathmoveto{\pgfqpoint{1.126536in}{3.135897in}}%
\pgfpathlineto{\pgfqpoint{7.682018in}{3.135897in}}%
\pgfusepath{stroke}%
\end{pgfscope}%
\begin{pgfscope}%
\pgfsetrectcap%
\pgfsetroundjoin%
\pgfsetlinewidth{1.505625pt}%
\definecolor{currentstroke}{rgb}{0.580392,0.403922,0.741176}%
\pgfsetstrokecolor{currentstroke}%
\pgfsetdash{}{0pt}%
\pgfpathmoveto{\pgfqpoint{2.251710in}{3.363810in}}%
\pgfpathlineto{\pgfqpoint{2.501710in}{3.363810in}}%
\pgfpathlineto{\pgfqpoint{2.751710in}{3.363810in}}%
\pgfusepath{stroke}%
\end{pgfscope}%
\begin{pgfscope}%
\pgfsetbuttcap%
\pgfsetmiterjoin%
\definecolor{currentfill}{rgb}{0.580392,0.403922,0.741176}%
\pgfsetfillcolor{currentfill}%
\pgfsetlinewidth{1.003750pt}%
\definecolor{currentstroke}{rgb}{0.580392,0.403922,0.741176}%
\pgfsetstrokecolor{currentstroke}%
\pgfsetdash{}{0pt}%
\pgfsys@defobject{currentmarker}{\pgfqpoint{-0.041667in}{-0.041667in}}{\pgfqpoint{0.041667in}{0.041667in}}{%
\pgfpathmoveto{\pgfqpoint{-0.041667in}{-0.041667in}}%
\pgfpathlineto{\pgfqpoint{0.041667in}{-0.041667in}}%
\pgfpathlineto{\pgfqpoint{0.041667in}{0.041667in}}%
\pgfpathlineto{\pgfqpoint{-0.041667in}{0.041667in}}%
\pgfpathlineto{\pgfqpoint{-0.041667in}{-0.041667in}}%
\pgfpathclose%
\pgfusepath{stroke,fill}%
}%
\begin{pgfscope}%
\pgfsys@transformshift{2.501710in}{3.363810in}%
\pgfsys@useobject{currentmarker}{}%
\end{pgfscope}%
\end{pgfscope}%
\begin{pgfscope}%
\definecolor{textcolor}{rgb}{0.000000,0.000000,0.000000}%
\pgfsetstrokecolor{textcolor}%
\pgfsetfillcolor{textcolor}%
\pgftext[x=2.876710in,y=3.276310in,left,base]{\color{textcolor}\rmfamily\fontsize{18.000000}{21.600000}\selectfont LPP\(\displaystyle ^{R,T}\)}%
\end{pgfscope}%
\begin{pgfscope}%
\pgfsetbuttcap%
\pgfsetroundjoin%
\pgfsetlinewidth{1.505625pt}%
\definecolor{currentstroke}{rgb}{0.580392,0.403922,0.741176}%
\pgfsetstrokecolor{currentstroke}%
\pgfsetstrokeopacity{0.900000}%
\pgfsetdash{{5.550000pt}{2.400000pt}}{0.000000pt}%
\pgfpathmoveto{\pgfqpoint{3.812661in}{3.363810in}}%
\pgfpathlineto{\pgfqpoint{4.062661in}{3.363810in}}%
\pgfpathlineto{\pgfqpoint{4.312661in}{3.363810in}}%
\pgfusepath{stroke}%
\end{pgfscope}%
\begin{pgfscope}%
\pgfsetbuttcap%
\pgfsetmiterjoin%
\definecolor{currentfill}{rgb}{0.580392,0.403922,0.741176}%
\pgfsetfillcolor{currentfill}%
\pgfsetfillopacity{0.900000}%
\pgfsetlinewidth{1.003750pt}%
\definecolor{currentstroke}{rgb}{0.580392,0.403922,0.741176}%
\pgfsetstrokecolor{currentstroke}%
\pgfsetstrokeopacity{0.900000}%
\pgfsetdash{}{0pt}%
\pgfsys@defobject{currentmarker}{\pgfqpoint{-0.058926in}{-0.058926in}}{\pgfqpoint{0.058926in}{0.058926in}}{%
\pgfpathmoveto{\pgfqpoint{0.000000in}{-0.058926in}}%
\pgfpathlineto{\pgfqpoint{0.058926in}{0.000000in}}%
\pgfpathlineto{\pgfqpoint{0.000000in}{0.058926in}}%
\pgfpathlineto{\pgfqpoint{-0.058926in}{0.000000in}}%
\pgfpathlineto{\pgfqpoint{0.000000in}{-0.058926in}}%
\pgfpathclose%
\pgfusepath{stroke,fill}%
}%
\begin{pgfscope}%
\pgfsys@transformshift{4.062661in}{3.363810in}%
\pgfsys@useobject{currentmarker}{}%
\end{pgfscope}%
\end{pgfscope}%
\begin{pgfscope}%
\definecolor{textcolor}{rgb}{0.000000,0.000000,0.000000}%
\pgfsetstrokecolor{textcolor}%
\pgfsetfillcolor{textcolor}%
\pgftext[x=4.437661in,y=3.276310in,left,base]{\color{textcolor}\rmfamily\fontsize{18.000000}{21.600000}\selectfont LPP\(\displaystyle ^{T}\)}%
\end{pgfscope}%
\begin{pgfscope}%
\pgfsetbuttcap%
\pgfsetroundjoin%
\pgfsetlinewidth{1.505625pt}%
\definecolor{currentstroke}{rgb}{0.580392,0.403922,0.741176}%
\pgfsetstrokecolor{currentstroke}%
\pgfsetstrokeopacity{0.800000}%
\pgfsetdash{{1.500000pt}{2.475000pt}}{0.000000pt}%
\pgfpathmoveto{\pgfqpoint{5.202692in}{3.363810in}}%
\pgfpathlineto{\pgfqpoint{5.452692in}{3.363810in}}%
\pgfpathlineto{\pgfqpoint{5.702692in}{3.363810in}}%
\pgfusepath{stroke}%
\end{pgfscope}%
\begin{pgfscope}%
\pgfsetbuttcap%
\pgfsetroundjoin%
\definecolor{currentfill}{rgb}{0.580392,0.403922,0.741176}%
\pgfsetfillcolor{currentfill}%
\pgfsetfillopacity{0.800000}%
\pgfsetlinewidth{1.003750pt}%
\definecolor{currentstroke}{rgb}{0.580392,0.403922,0.741176}%
\pgfsetstrokecolor{currentstroke}%
\pgfsetstrokeopacity{0.800000}%
\pgfsetdash{}{0pt}%
\pgfsys@defobject{currentmarker}{\pgfqpoint{-0.041667in}{-0.041667in}}{\pgfqpoint{0.041667in}{0.041667in}}{%
\pgfpathmoveto{\pgfqpoint{0.000000in}{-0.041667in}}%
\pgfpathcurveto{\pgfqpoint{0.011050in}{-0.041667in}}{\pgfqpoint{0.021649in}{-0.037276in}}{\pgfqpoint{0.029463in}{-0.029463in}}%
\pgfpathcurveto{\pgfqpoint{0.037276in}{-0.021649in}}{\pgfqpoint{0.041667in}{-0.011050in}}{\pgfqpoint{0.041667in}{0.000000in}}%
\pgfpathcurveto{\pgfqpoint{0.041667in}{0.011050in}}{\pgfqpoint{0.037276in}{0.021649in}}{\pgfqpoint{0.029463in}{0.029463in}}%
\pgfpathcurveto{\pgfqpoint{0.021649in}{0.037276in}}{\pgfqpoint{0.011050in}{0.041667in}}{\pgfqpoint{0.000000in}{0.041667in}}%
\pgfpathcurveto{\pgfqpoint{-0.011050in}{0.041667in}}{\pgfqpoint{-0.021649in}{0.037276in}}{\pgfqpoint{-0.029463in}{0.029463in}}%
\pgfpathcurveto{\pgfqpoint{-0.037276in}{0.021649in}}{\pgfqpoint{-0.041667in}{0.011050in}}{\pgfqpoint{-0.041667in}{0.000000in}}%
\pgfpathcurveto{\pgfqpoint{-0.041667in}{-0.011050in}}{\pgfqpoint{-0.037276in}{-0.021649in}}{\pgfqpoint{-0.029463in}{-0.029463in}}%
\pgfpathcurveto{\pgfqpoint{-0.021649in}{-0.037276in}}{\pgfqpoint{-0.011050in}{-0.041667in}}{\pgfqpoint{0.000000in}{-0.041667in}}%
\pgfpathlineto{\pgfqpoint{0.000000in}{-0.041667in}}%
\pgfpathclose%
\pgfusepath{stroke,fill}%
}%
\begin{pgfscope}%
\pgfsys@transformshift{5.452692in}{3.363810in}%
\pgfsys@useobject{currentmarker}{}%
\end{pgfscope}%
\end{pgfscope}%
\begin{pgfscope}%
\definecolor{textcolor}{rgb}{0.000000,0.000000,0.000000}%
\pgfsetstrokecolor{textcolor}%
\pgfsetfillcolor{textcolor}%
\pgftext[x=5.827692in,y=3.276310in,left,base]{\color{textcolor}\rmfamily\fontsize{18.000000}{21.600000}\selectfont LPP\(\displaystyle ^{T, P}\)}%
\end{pgfscope}%
\end{pgfpicture}%
\makeatother%
\endgroup%

%% file: 06_AccelerationTechniques.tex
\subsection{Acceleration Methods}\label{sec:ComputationalConsiderations} 
This section analyzes two configuration choices of the algorithm: different strategies for including the valid inequalities introduced in \autoref{subsec:ineq} and an extension to the refinement procedure aimed at speeding up convergence. All experiments in this section use a setting of $\lambda=0.25$.

\subsubsection{Valid inequalities}
Valid inequalities~\eqref{eq:validIneq} and~\eqref{eq:delta} are important for the convergence of the algorithm when allowing for lost demand. We evaluate three settings: no inequalities ($h_T=0$), a restricted set for paths requiring headways of 50$\%$ or less of the maximum ($h_T = 10$), and the full set of inequalities $(h_T = 20)$. \autoref{tab:validIneqNewSplitTime} summarizes the average performance across these settings, along with a reference comparison to $M$, using the same instance categorization as in \autoref{tab:CompAverage}. 
Using the restricted set ($h_T = 10$) achieves the best balance, reducing the average computation time to 61.9 minutes, which is faster than the settings with no inequalities (194.8 minutes) and the full set of inequalities (71.9 minutes). 
Although $h_T = 20$ achieves tighter bounds for the instances that cannot be solved to optimality, it is slower and solves two fewer instances to optimality compared to $h_T = 10$.  
\begin{table}[htbp]
\centering
\caption{Average computation time and gap ($\lambda=0.25$)}
\label{tab:validIneqNewSplitTime}
  \resizebox{\textwidth}{!}{
\begin{tabular}{lrrrrrrrr}
\toprule
 & \multicolumn{4}{c}{ Time (min.)} & \multicolumn{4}{c}{Gap ($\%$)} 
 \\ \cmidrule(l){2-5}  \cmidrule(l){6-9}
 & \multicolumn{1}{c}{$A (h_t=0)$} & \multicolumn{1}{c}{$A (h_t=10)$} & \multicolumn{1}{c}{$A (h_t=20)$} &
 \multicolumn{1}{c}{$M$} & \multicolumn{1}{c}{$A (h_t=0)$} & \multicolumn{1}{c}{$A (h_t=10)$} & \multicolumn{1}{c}{$A (h_t=20)$} & 
 \multicolumn{1}{c}{$M$} \\  \midrule
$A$ and $M$ optimal & 194.8 & 61.9 & 71.9 & 165.3 & 4.0 & 0 & 0 & 0 \\
$A$ optimal & 464.8 & 207.1 & 227.2 & 600 & 41.6 & 0 & 1.4 & 91.5 \\
Neither optimal & 599.5 & 600 & 600 & 600 & 46.9 & 28.0 & 25.6 & $>100$  \\
\bottomrule
\end{tabular}}
\end{table}
The valid inequalities present a trade-off: they increase the size of the subproblems but reduce the number of iterations required for the algorithm to converge, from an average of 57.2 iterations with $h_T=0$ to 15.6 iterations with $h_T=20$. A midpoint setting ($h_T = 10$) strikes the best balance with almost the same reduction in iterations while being faster on average and solving the most instances to optimality.
This is further illustrated in \autoref{sec:AppendixValidIneq}, which presents the number of instances solved over time and the cumulative number of instances with a certain gap.  
We also tested adding the inequalities as lazy constraints in a CPLEX callback but this proved inefficient due to the repeated process of solving the MILPs to integer optimality and adding constraints. 
Given its consistent advantages, we adopt the setting of $h_T = 10$ throughout the experiments. 

\subsubsection{Adding multiple frequencies in the refinement step} 
The \frameworkNameShortLP iteratively adds frequencies to eliminate infeasible solutions. By default, it adds one per selected line.
However, when many parallel lines exist in $\linePool$, it may cycle through similar solutions, requiring more iterations to converge. 
To address this, we extend the refinement step to also update the $\kSim$ most similar lines to each selected line $l\in \linePoolSol$. 
We measure similarity by the size of the intersection of trips $e\in E_l$ between lines and define $\linePool_{l\kSim}$ as the set of $\kSim$ most similar lines to line $l$. 
For each line $l\in \linePoolSol$ with insufficient vehicles assigned to operate the selected frequency $f\in \frequencySet_l$, we update the representation of both $l$ and $l_1, l_2, ..., l_\kSim\in \linePool_{l\kSim}$ to exclude solutions where frequency $f$ could be selected with $\underline{\phi}_{lf}$ vehicles. 
While this increases the model size, it reduces the number of iterations required by addressing infeasibilities before they occur. 
There are multiple ways to update the lines $\linePool_{l\kSim}$, as different lines have distinct characteristics, making the best choice of frequency to update non-trivial. Our approach prioritizes updating the lines and frequencies that are most similar to those already selected.

The proposed strategy can accelerate convergence and find better solutions in cases where the procedure is terminated. However, as $\kSim$ needs to be tuned specifically for each instance, it is not systematically applied in the other sections of this work.  
However, in specific cases, the computation time is significantly reduced. \autoref{tab:similarLines} presents two such examples, reporting the computation time, iterations required, and the number of variables and constraints relative to $M$.
For the \sioux instance, adopting this strategy with $\kSim=1$ speeds up computation time by 32\% compared to $\kSim=0$ and reduces the number of iterations from 30 to 19, with only a modest increase in model size. 
For the \lowersaxony instance, where convergence generally requires fewer iterations, higher values of $\kSim$ perform better, with $\kSim=2$ reducing computation time by $34\%$.

\begin{table}[htbp]
\centering
\caption{Impact of adding $\kSim$ similar lines in refinement step.}
\resizebox{0.8\textwidth}{!}{
\begin{tabular}{@{}lrrrrrrrr@{}}
\toprule
& 
\multicolumn{4}{c}{\begin{tabular}[c]{@{}c@{}}\sioux\\ $|\linePool|$ = 23, $|\ODset|$=300\end{tabular}} & 
\multicolumn{4}{c}{\begin{tabular}[c]{@{}c@{}}\lowersaxony\\ $|\linePool|$ = 41, $|\ODset|$=300\end{tabular}}
 \\ \cmidrule(l){2-5}  \cmidrule(l){6-9}
 &  $\kSim = 0$ & $\kSim = 1$ & $\kSim = 2 $&$ \kSim = 3$ &  $\kSim = 0$ & $\kSim = 1$ & $\kSim = 2 $&$ \kSim = 3$ \\ \midrule
Time (minutes) &  467.15   &   317.7 &353.5 & 375.4 & 84.0 & 75.2 & 55.4 & 62.4 \\
 Iterations &  30 & 19 & 16 & 15 & 26 & 18 & 13 & 12 \\
Variables (\%) & 19.7 &
23.5 &
26.8 &
30.0 & 3.5 &
4.1 &
5.0 &
5.8  \\
 Constraints (\%) & 47.5 &
51.4 &
54.6 &
59.0 & 29.5 &
32.1 &
35.4 &
37.8 \\ \bottomrule
\end{tabular}}
\label{tab:similarLines}
\end{table}
While this strategy does not enable solving additional instances within the time limit compared to $\kSim=0$, it can yield tighter bounds.  For instance, for \sioux with 46 lines and 300 ODs, we can tighten the optimality gap from $11.2\%$ to $4.8\%$ by also updating the $\kSim=2$ most similar lines. The lower and upper bounds are improved by $1.7\%$ and $5.2\%$, respectively. However, in some cases, setting $\kSim>0$ increases the problem size such that an optimal solution can no longer be found within the time limit.
This approach enables exploring the trade-off between model size and convergence speed, and although not universally applicable, it can significantly improve the performance of the \frameworkNameShortLP on certain instances.

%% file: 06_RouteAssignment.tex
\subsection{Path Assignment} \label{sec:RouteAssignment}
The model aims to find a balance between operator costs and passenger service, and while we enforce a limit on the travel time, the assigned path is not necessarily the preferred path. Additionally, capacity constraints may lead to passengers being distributed across different paths -- an adjustment they may not be able to anticipate or willing to accommodate in practice. We therefore compare the path assignment from the \problem to two different assignment methods: 1) a probabilistic assignment based on a logit model, and 2) a shortest path assignment. Both assignments are computed using the generated path set. Unlike the path assignment of the model, these assignments ignore capacity constraints. As a result, the shortest path assignment serves as a lower bound for the passenger costs for a given selection of lines and frequencies.
For the logit assignment, we calculate the probability of a passenger choosing a path $p$ in the set of available paths given the selected lines and frequencies.   
The probability of choosing path $p_i$ is defined as $P(p_i)=\frac{e^\theta c_{p_i}}{\sum_{j =1}^n e^\theta c_{p_j}} \forall i = 1, \dots, n$, where $\theta = -0.2$. 
\begin{table}[htbp]
\centering
\caption{Path Assignment}
\begin{tabular}{@{}llll@{}}
\toprule
Path Assign. & Pas. costs & Demand captured & WTT \\ \midrule
Model & 13,125.6 & 58.74\% & 21.45 \\ 
Logit & 13,183.4 & 60.84\%& 21.57 \\ 
Shortest path & 13,043.6 & 60.84\% & 21.44 \\ \bottomrule
\end{tabular}
\label{tab:routeAssignmenDEstar}
\end{table}

\autoref{tab:routeAssignmenDEstar} shows that the passenger service and demand captured are not very sensitive to the path assignment. The results, averaged over the 23 instances solved to optimality using $\lambda = 0.25$, indicate only small variations across assignments.  
The shortest path assignment yields the lowest passenger costs, as this represents the optimal routing for passengers. Compared to the model assignment, the costs are only $0.6\%$ lower, which can be translated to 0.01 weighted minutes on average per passenger. This indicates that most passengers are assigned a near-optimal path. In contrast, the probabilistic assignment results in passenger costs $0.4\%$ higher than the model assignment, as passengers are distributed over all possible paths. 
The total transported demand is slightly higher under logit and shortest path assignments, which ignore capacity constraints, suggesting additional passengers would use the system if sufficient capacity was available. 

\added{An examination of the load profiles supports this interpretation. The load profiles of the operated lines show that only a few lines reach their maximum capacity, and when this occurs, it is limited to specific arcs. \autoref{apx:Capacity} presents the load profiles and vehicle utilization rates across different values of $\lambda$. On average, less than $1.4\%$ of the arcs selected in the solution operate at maximum capacity. Moreover, vehicle utilization weighted by passenger-minutes remains modest, ranging from 55\% at $\lambda = 0.1$ to 38\% at $\lambda = 1$. Overall, these findings indicate that the designed networks have limited crowding.}

\begin{figure}[htbp]
    \centering
    \begin{subfigure}{0.45\textwidth}
        \centering
        \begin{tikzpicture}
        \begin{axis}[
            xlabel={$\lambda$},
            ylabel={Fraction of OD pairs (\%)},
            legend pos=north west,
            grid=both,
            width=\textwidth,
            height=6cm,
            ylabel style={text width=5cm, align=center}
        ]
        \addplot[only marks, mark=*, mark size=1.5pt, opacity=0.5, color=blue, forget plot] 
        table [x=lambda, y=pct_split_any, col sep=comma] {Appendix/routing_data/pct_ods_split_any_vs_lambda.csv};
        \addplot[thick, color=red]
        table [x=lambda, y=pct_split_any, col sep=comma] {Appendix/routing_data/pct_ods_split_any_vs_lambda_avg.csv};
        \addlegendentry{Average}
        \end{axis}
        \end{tikzpicture}
        \caption{Fraction of OD pairs routed on more than a single path (lines show means, scatter points show individual instances).} \label{pctODsSplit}
    \end{subfigure}%
    \hfill
    \begin{subfigure}{0.45\textwidth}
        \centering
        \begin{tikzpicture}
        \begin{axis}[
            ybar,
            xlabel={Difference from shortest path (minutes)},
            ylabel={Fraction of passengers (\%)},
            width=\textwidth,
            height=5.8cm,
            grid=both,
            bar width=4pt,
            ylabel style={text width=5cm, align=center}
        ]
        \addplot+[
            fill=blue!50
        ] table [x=bin_centers, y=pct_rel_to_routed, col sep=comma] {Appendix/routing_data/extra_cost_distribution_norm_PT_histogram_L0.25.csv};
        \end{axis}
        \end{tikzpicture}
        \caption{Distribution of passenger cost increase (weighted travel time). Averaged over all instances solved to optimality with $\lambda = 0.25$}
        \label{extraPTtimeDist}
    \end{subfigure}
            \caption{\added{Insights on route assignment}}
        \label{fig:routing_insights}
\end{figure}
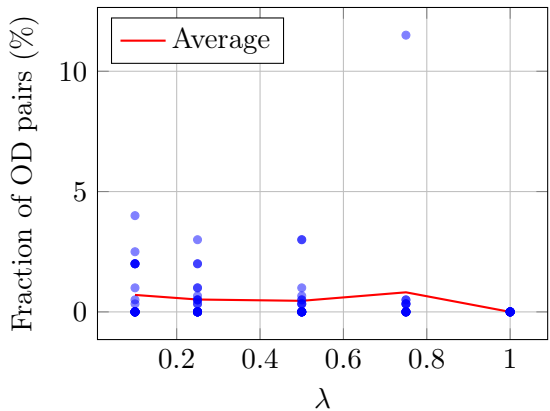
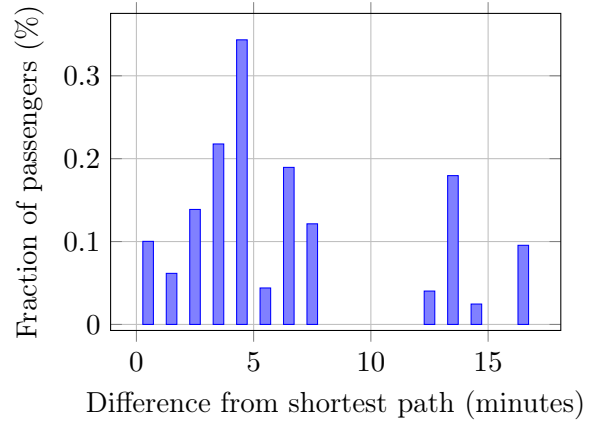

\added{
\autoref{pctODsSplit} shows the fraction of OD pairs routed on more than one path, including both passengers using multiple PT paths and those combining PT with an alternative mode. As expected, the fraction of split OD pairs decreases with higher values of $\lambda$, reflecting greater sensitivity to travel time. Nevertheless, even for lower values of $\lambda$, more than 99\% of OD pairs are, on average, assigned to a single path.
\autoref{extraPTtimeDist} shows the distribution of additional travel time per passenger compared to the shortest available path, weighted by the number of passengers routed. Passengers routed on the shortest path are excluded. While $98.4\%$ of passengers are routed on the shortest available path, the remaining are routed on paths with detours of up to 16 additional minutes, demonstrating that the model occasionally assigns paths longer than the shortest path to satisfy the capacity restrictions.  
}

%% file: 07Discussions.tex
\section{Conclusion and Discussion}\label{sec:Discussion}
This paper presents an exact solution method for line planning for public transport that integrates routing, capacity, and demand considerations. 
The proposed algorithm iteratively refines the frequency representation and associated path sets to generate non-decreasing lower bounds until a feasible solution is found, exploiting problem structure to outperform standard MILP formulations. 

We evaluate the \frameworkNameShortLP on 36 instances of various characteristics. Computational results show that our implementation can speed up computation time significantly and provide stronger bounds on the solution quality compared to solving the problem as a standard MILP. 
The method is tested over different values of the trade-off parameter $\lambda$, which controls the relative importance of passenger versus operator costs. 
We show that it performs particularly well for lower values of $\lambda$ where the passenger and operator costs contribute more evenly to the objective. In these cases, the algorithm solves more instances to optimality and achieves substantial speed-ups. Across the full range of $\lambda$ values, it consistently provides tighter bounds than solving the full formulation.

While this study represents a promising application of the \frameworkNameShortLP on instances with up to 129 lines and 300 ODs, scalability to larger, more realistic networks is limited. 
Scalability could be improved through network and OD matrix preprocessing or dynamic path generation methods such as column generation. Additionally, the \frameworkNameShortLP could be integrated into a matheuristic framework.
\added{Moreover, our results show that for passenger-focused objectives, the algorithm requires more iterations to converge as the refinement step is guided by operating costs, which becomes less relevant for the solution. 
This suggests that refinement strategies focused on the passenger service dimension could accelerate convergence and represent a promising direction for future research.}

The formulation relies on the assumption that a set of reasonable paths for each OD pair can be determined. This enables two main ideas: first, it allows for a detailed modeling of passenger preferences and utility; second, it accounts for the fact that a passenger’s willingness to accept a path depends on the service. This is represented by a threshold on the acceptable weighted travel time for each OD pair. Our results demonstrate that considering this interaction between demand and line concept during the line planning stage can lead to solutions with lower total costs and enhanced passenger service, assuming comparable passenger behavior and operator resources.
The current approach assumes that a path is either available to all passengers in an OD or not at all. 
Future work could explore more realistic approaches to demand, where the likelihood of passengers choosing a path is a function of weighted travel time rather than a binary option.
Additionally, the expected waiting and transfer times could account for the impact of parallel lines\added{, and the weighted travel time could be extended to reflect passenger discomfort due to crowding}. 
\added{Moreover, future work could look at more explicit ways to handle time-varying demand, e.g., by determining related peak and non-peak hour line plans within a single model.}

Finally, the described methodology can be applied to other network design problems containing discretized non-linear trade-offs between variables such as between passenger and operator costs in line planning. 
Future research could explore its application in other relevant contexts.

%% file: 08Appendix.tex
\input{Appendix/Example}
\input{Appendix/Proofs}
\input{Appendix/Tables}
\input{Appendix/SocialMaximum}

\input{Appendix/Plots}
\input{Appendix/TableSummary}
\input{Appendix/Demand}
\input{Appendix/ValidInequalities}
\input{Appendix/FrequencySet}
\input{Appendix/Capacity}

%% file: Appendix/Example.tex
\newpage
\setcounter{page}{1}
\section{A single line example}\label{app_Example}

To give further
insight into how the algorithm works, we consider a
simple PTN consisting of five sequential stops and a single
OD pair traveling from one end to the other. A line runs back and
forth between these stops, and each vehicle can complete one round trip per hour. We consider five different headways
$h\in \headwaySet$ with associated vehicle requirements $\phi_{lh}$
and corresponding weighted travel time (WTT) as specified in
\autoref{tab:smallExample}, which also describes the parameter setting
in terms of capacity, operating costs, and demand. We assume travelers
will use any frequency of the line, setting $\wttThresOD$ arbitrarily high. For this example, we assume that the weighted travel time translates directly to the equivalent monetary value.

\begin{table}[htbp]
\centering
\caption{Specifications of operator and passenger costs for a single line with five candidate frequencies.}
\label{tab:smallExample}
\begin{minipage}{0.3\textwidth}
    \centering
    \begin{tabular}{ccc} \toprule
        Headway & Veh.\ req.\ $\phi$ & WTT \\ \midrule
        $h_1 = 5$  & 12 & 35 \\
        $h_2 = 10$ & 6  & 40 \\
        $h_3 = 15$ & 4  & 45 \\
        $h_4 = 20$ & 3  & 50 \\
        $h_5 = 30$ & 2  & 55 \\ \bottomrule
    \end{tabular}
\end{minipage}
\hspace{1cm}  
\begin{minipage}{0.3\textwidth}
    \centering
    \begin{tabular}{cc} \toprule
        Parameters & Value \\ \midrule
        $\lambda$ & 1 \\
        $\delta_{lv} $ & 50 \\
        $c_v$ & 2000 \\
        $h_l$ & 0 \\
        \ODw & 150 \\
        $B$ & $\infty$ \\ 
        $R$ & $0$ \\ \bottomrule
    \end{tabular}
\end{minipage}
\end{table}

\begin{figure}[H]
 \centering
 \subfloat[Initialization and first solve] {\includegraphics[width=0.25\textwidth]{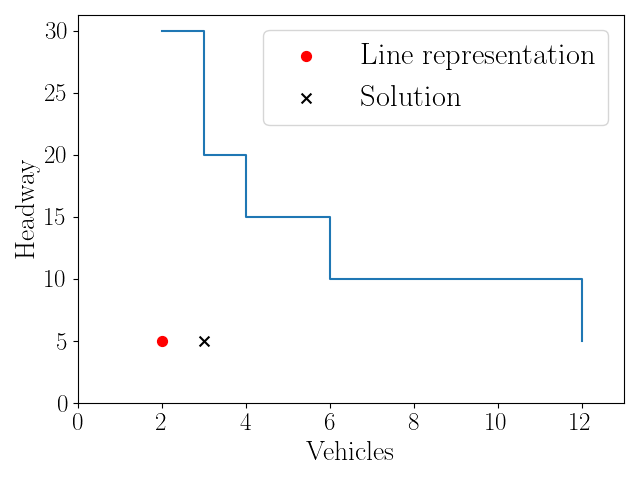}}
 \subfloat[Iteration 1]
 {\includegraphics[width=0.25\textwidth]{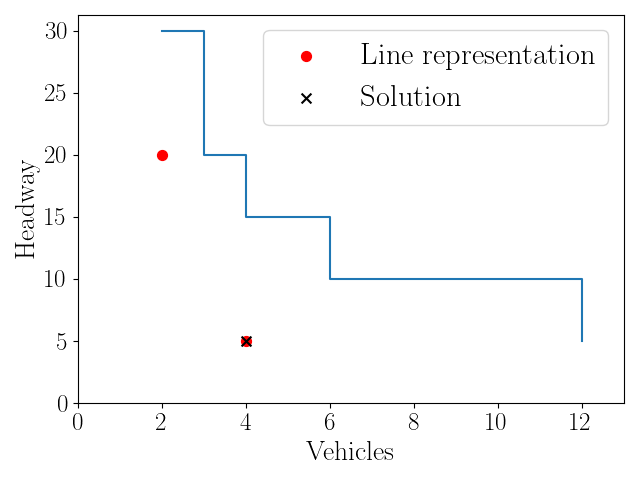}} 
 \subfloat[Iteration 2]{\includegraphics[width=0.25\textwidth]{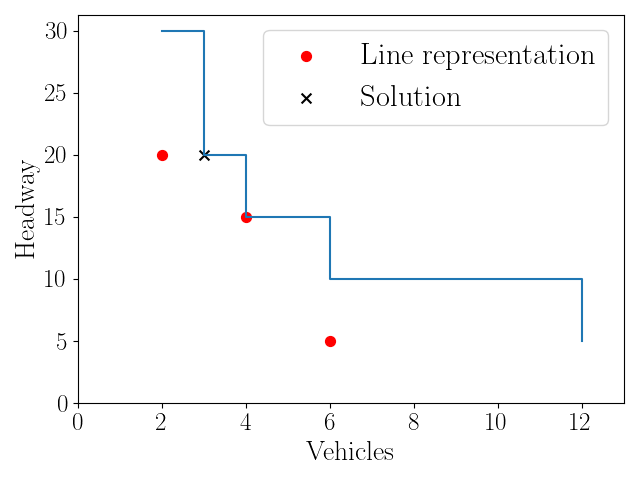}} \caption{Algorithm iterations for solving example in \autoref{tab:smallExample}}
 \label{fig:60WPH}
 \end{figure}

For this small example, it can be determined that the optimal solution is to open the line at $h_4$ of 20 minutes, which allows for just enough vehicles to cover demand. \autoref{tab:algIterations} shows how the line is represented in each iteration of the algorithm and the solution obtained from solving the \problem with the reduced headway set. At initialization, the problem contains only the minimum headway ($h_1$ = 5 minutes), which is set to require the minimum number of vehicles ($\underline{\phi}_{lh_1}=2$). 
The model can, therefore, only select this headway, and it assigns three vehicles, as this is the minimum required to transport the demand. 
However, this solution is not valid as 12 vehicles are required to operate the headway, i.e., $z^*_l<f_v(h_1)$.  In the next iteration, the line specification is updated to cut off the previous solution. Now, at least four vehicles need to be assigned to operate the headway $h_1$. The optimal choice is to open the line at the minimum headway and assign the required vehicles ($z^*_l=4$). This is still not enough to operate the promised service level, so the line representation is updated again. In iteration two, operating $h_1$ requires at least six vehicles, while  $h_3$ and $h_4$ require at least four and two vehicles, respectively. In this iteration, it is no longer cost-effective to open the minimum headway. Instead, the model selects $h_4$ and assigns three vehicles. Since there are sufficient vehicles to operate this frequency, $3\geq f_v(h_4)$, the solution is feasible, and the algorithm terminates. The algorithm has found the optimal solution while only including three out of five candidate headways. \autoref{fig:60WPH} illustrates the iterative process of the \frameworkNameShortLP.

\begin{table}[htbp]
\centering
\caption{Overview of algorithm iterations showing the line representation, computed solution, and feasibility with respect to \problem{($\headwaySet$)} at each step.}
\label{tab:algIterations}
\begin{tabular}{@{}llll@{}}
\toprule
 & Representation of line $l$ & Solution & Feasible solution \\ \midrule
Initialization and first solve & $\{\{h_1, \underline{\phi}_{lh_1} = 2 \}\}$ & $y_{lh_1} = 1, z_l = 3$ & No \\ \midrule
Iteration 1 & \begin{tabular}[c]{@{}l@{}}\{\{$h_1, \underline{\phi}_{lh_1} = 4 $\},\\ \ \{$h_4, \underline{\phi}_{lh_4} = 2$\}\}\end{tabular} & $y_{lh_1} = 1, z_l = 4$ & No \\ \midrule
Iteration 2 & \begin{tabular}[c]{@{}l@{}}\{\{$h_1, \underline{\phi}_{lh_1} = 6 $\},\\ \ \{$h_3, \underline{\phi}_{lh_3} = 4$\},\\ \ \{$h_4, \underline{\phi}_{lh_4}=2 $\}\}\end{tabular} & $y_{lh_4} = 1, z_l = 3$ & Yes \\ \bottomrule
\end{tabular}
\end{table}

\begin{figure}[H]
 \centering
 \subfloat[Initialization and first solve] {\includegraphics[width=0.25\textwidth]{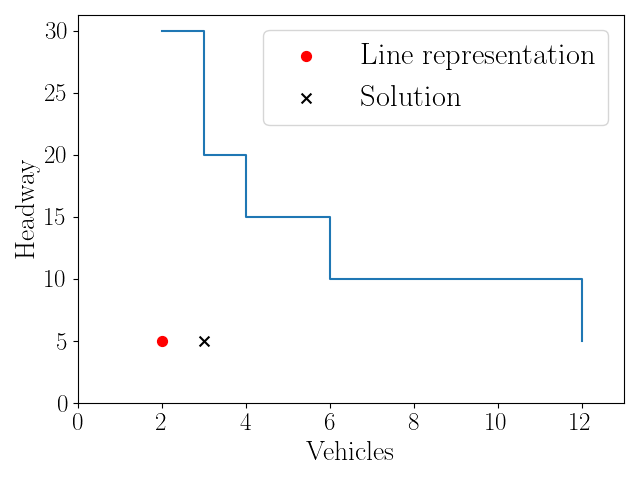}}
 \subfloat[Iteration 1]
 {\includegraphics[width=0.25\textwidth]{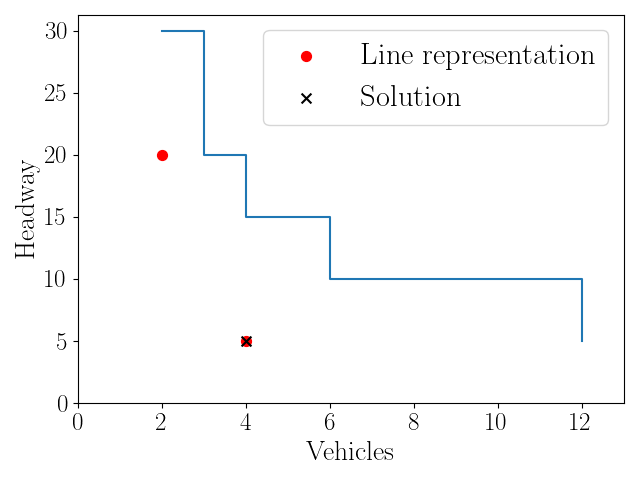}} 
 \subfloat[Iteration 2]{\includegraphics[width=0.25\textwidth]{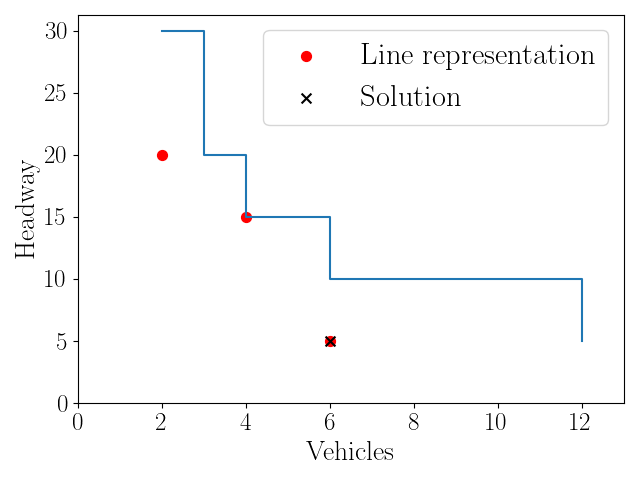}} 
 \subfloat[Iteration 3]
 {\includegraphics[width=0.25\textwidth]{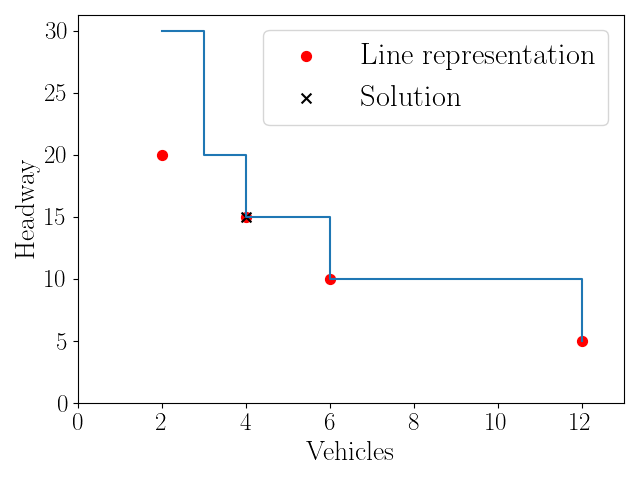}} 
 \caption{Algorithm iterations for solving example in \autoref{tab:smallExample} with increased value of time}
 \label{fig:180WPH}
 \end{figure}

The weighing of operator costs relative to the passenger costs affects the solution procedure. 
Assume that passengers perceive the travel time as more inconvenient; say, the path cost is three times the weighted travel time. 
Now, the optimal solution is to upgrade the service level to $h_3$ (15 minutes) at the cost of assigning one additional vehicle. \autoref{fig:180WPH} shows the iterative process given this setting. The algorithm converges to the optimal solution after four iterations, selecting a headway of 15 minutes. 
This is caused by the objective weighing passenger costs higher, and the maximum frequency thus remains the optimal choice until the representation is updated to require the true number of vehicles required. The objective's cost structure hence impacts not only the final solution but also how fast the algorithm converges. The algorithm uses the capacity assignment to guide the refinement, so this step is most efficient when the capacity is important for the optimal solution. 

%% file: Appendix/Proofs.tex
\section{Proofs}\label{sec:Proofs}

\subsection{Theorem~\ref{thm:lb}}
\label{sec:proof1}

\begin{proof}
Let $(\mathbf{y}^*,\mathbf{x}^*,\mathbf{z}^*)$ be an optimal solution to \problem{$(\mathcal{F})$} with objective value $\theta^*$. We define the associated optimal set of lines that are selected in this solution as \added{$\mathcal{L}^*=\left\{l\in \mathcal{L} \mid \exists f \in \mathcal{F}_l : y_{lf}^{*} =1 \right\}$} and assume that the optimal frequency for line $l\in\mathcal{L}^*$ is given by $f_l^*$. We construct the feasible solution $(\mathbf{\bar{y}},\mathbf{\bar{x}}, \mathbf{z}^*)^{(i)}$ with objective function value $\bar{\theta}_i$ at iteration $i$ of the algorithm using the lines specified in $\mathcal{L}^*$ and the same vehicle assignment $\mathbf{z}^*$. Specifically, we let $\bar{y}_{lf^\prime}=1$ for $l\in \mathcal{L}^*$ for some $f^\prime \geq f_l^*$ with $\underline{\phi}_{lf^\prime}^{i} \leq z^*_l$. The existence of $f^\prime$ is guaranteed by design of the \frameworkNameShortLP; it iteratively refines the problem starting with the highest frequencies. The refinement strategy always ensures that $z^*_l$ remains feasible for the headway $f_h(z^*_l)$ and possibly shorter headways due to the underestimation of the number of vehicles necessary. Solution $(\mathbf{\bar{y}},\mathbf{\bar{x}}, \mathbf{z}^*)^{(i)}$ uses the same vehicle assignment and lines as $(\mathbf{y}^*,\mathbf{x}^*,\mathbf{z}^*)$. The respective objective function values can therefore only differ in passenger routing costs. Passenger costs are, however, non-increasing in the frequency. As the same lines are used and the provided vehicle capacity is the same, alternative, potentially cheaper, paths exist for all passengers that are routed in $(\mathbf{y}^*,\mathbf{x}^*,\mathbf{z}^*)$; they are just assigned higher frequency arcs in the CGN. We conclude that $\bar{\theta_i}\leq \theta^*$. Since the optimal solution at iteration $i$, $\vectorSol{i}$ has an objective value $\theta_i^* \leq \bar{\theta}_i$, we conclude that $\theta_i^*$ provides a lower bound on $\theta^*$.
\end{proof}

\subsection{Theorem~\ref{thm:finite}}
\label{sec:proof2}
\begin{proof}
Iterations of the algorithm continue as long as the solution \vectorSol{i} is infeasible for \problem{$(\frequencySet)$}. In such a solution, there must exist a line $l$ for which $y_{lf^\prime}=1$ for some $f^\prime \in \mathcal{F}^i_l$ such that $z_l^*<f_v(f^\prime)$, i.e., not enough vehicles have been assigned to operate frequency $f^\prime$. The algorithm corrects for this by introducing the frequency that can be operated with $z_l^*$ vehicles and increases the minimum number of vehicles needed to operate frequency $f^\prime$. At each iteration, at least one frequency and the corresponding paths are therefore added to the current model. For a given line, we can therefore iterate at most $|\mathcal{F}_l|$ times.
\end{proof}

%% file: Appendix/Tables.tex
\section{Computational results for different values of \texorpdfstring{$\lambda$}{lambda}}\label{sec:TablesSup}

The following tables report the computational effort and solution quality for the line planning instances solved using both the \frameworkNameShortLP (\emph{A}) and CPLEX (\emph{M}) for different values of $\lambda$. We report the total time spent solving the \problem, the reported optimality gap at termination, the number of iterations for \emph{A}, the relative improvement in lower and upper bounds (RLB and RUB) reported by \emph{A} compared to \emph{M}, and finally the number of variables and constraints generated by \emph{A} in the final iteration relative to the model size in \emph{M}. The tables highlight the best-performing method in \textcolor{blue}{blue}, based on solution time, then upper bound, and finally lower bound. An asterisk indicates the method with the smallest optimality gap. If the lower bound results in an optimality gap exceeding 100\%, the RLB is reported as “-”. The relative improvement in bound is calculated as follows: 
\setlength{\jot}{2pt}  
\begin{align*}
\text{Relative Lower Bound Improvement (RLB)} &= \frac{z_{LB}^A - z_{LB}^M}{z_{LB}^M} \\
\text{Relative Upper Bound Improvement (RUB)} &= \frac{z_{UB}^M - z_{UB}^A}{z_{UB}^M}
\end{align*}
RLB and RUB values larger than zero mean that $A$ provides stronger lower and upper bounds, respectively.

\begin{table}[H]
	\centering
	\caption{
		Computational results for solving the \problem using either \frameworkNameShortLP (\emph{A}) or CPLEX (\emph{M}) with $\lambda = 0.1$.}
	\resizebox{\textwidth}{!}{
		\begin{tabular}{lrrrrrrrrrrr}\toprule
			\multicolumn{1}{l}{\textbf{Network}} &
			\multicolumn{1}{c}{\textbf{Pool}} &
			\multicolumn{1}{c}{\textbf{ODs}} &
			\multicolumn{2}{c}{\textbf{Time (min)}} &
			\multicolumn{2}{c}{\textbf{Gap (\%)}} &
			\multicolumn{1}{c}{\textbf{Iter.}} &
			\multicolumn{1}{c}{$\textbf{RLB (\%)}$} &
			\multicolumn{1}{c}{$\textbf{RUB (\%)}$} &
			\multicolumn{1}{c}{\textbf{Var. (\%)}} &
			\multicolumn{1}{c}{\textbf{Const. (\%)}} \\
			\cmidrule(lr){4-5}
			\cmidrule(lr){6-7} 
			\cmidrule(lr){8-8}
			\cmidrule(lr){9-9}
			\cmidrule(lr){10-10}
			\cmidrule(lr){11-11}
			\cmidrule(lr){12-12}
			& $|\linePool|$
			& $|\ODset|$
			& $A$ & M
			& $A$ & M
			& $A$ & $A$ & $A$ & $A$ & $A$ \\ \midrule 
			\midrule 
\sioux & 23 & 100 &  \textcolor{blue}{0.9} & 61.5 & 0.00 & 0.00 & 12 & 0.0 & 0.0 & 4.4 & 24.1 \\
\sioux & 23 & 200 & \textcolor{blue}{134.6} & 600.0 & 0.00$^*$ & $>100$ & 30 & $-$ & 8.5 & 11.9 & 38.3 \\
\sioux & 23 & 300 & \textcolor{blue}{317.4} & 500.3 & 0.00 & 0.00 & 30 & 0.0 & 0.0 & 13.7 & 39.6 \\
\sioux & 46 & 100 & \textcolor{blue}{5.3} & 419.4 & 0.00 & 0.00 & 15 & 0.0 & 0.0 & 4.5 & 25.8 \\
\sioux & 46 & 200 & \textcolor{blue}{494.8} & 600.0 & 0.00$^*$ & $>100$ & 36 & $-$ & 59.9 & 11.3 & 38.4 \\
\sioux & 46 & 300 & \textcolor{blue}{600.0} & 600.0 & 18.59$^*$ & $>100$ & 27 & $-$ & 145.8 & 10.0 & 34.5 \\
\sioux & 92 & 100 & \textcolor{blue}{600.0} & 600.0 & 2.17$^*$ & $>100$ & 40 & $-$ & 41.9 & 4.4 & 27.1 \\
\sioux & 92 & 200 & \textcolor{blue}{600.0} & 600.0 & $>100$ & $>100$ & 12 & $-$ & 73.3 & 2.8 & 20.5 \\
\sioux & 92 & 300 & \textcolor{blue}{600.0} & 600.0 & $>100$ & $>100$ & 8 & $-$ & 99.6 & 2.1 & 17.4 \\
\midrule
\grid & 30 & 100 & \textcolor{blue}{0.1} & 37.9 & 0.00 & 0.00 & 9 & 0.0 & 0.0 & 2.4 & 17.4 \\
\grid & 30 & 200 & \textcolor{blue}{2.5} & 233.1 & 0.00 & 0.00 & 16 & 0.0 & 0.0 & 4.2 & 21.8 \\
\grid & 30 & 300 & \textcolor{blue}{9.8} & 419.1 & 0.00 & 0.00 & 14 & 0.0 & 0.0 & 4.0 & 20.9 \\
\grid & 60 & 100 & \textcolor{blue}{15.4} & 600.0 & 0.00$^*$ & $>100$ & 25 & $-$ & 24.1 & 2.6 & 19.3 \\
\grid & 60 & 200 & \textcolor{blue}{600.0} & 600.0 & 19.24$^*$ & $>100$ & 38 & $-$ & 67.0 & 4.0 & 22.5 \\
\grid & 60 & 300 & \textcolor{blue}{600.0} & 600.0 & 46.85$^*$ & $>100$ & 19 & $-$ & 63.5 & 3.0 & 18.4 \\
\grid & 123 & 100 & \textcolor{blue}{600.0} & 600.0 & 12.76$^*$ & $>100$ & 40 & $-$ & 54.5 & 2.2 & 17.9 \\
\grid & 123 & 200 & \textcolor{blue}{600.0} & 600.0 & 80.42$^*$ & $>100$ & 15 & $-$ & 64.6 & 1.4 & 13.3 \\
\grid & 123 & 300 & \textcolor{blue}{600.0} & 600.0 & $>100$ & $>100$ & 9 & $-$ & 74.6 & 1.1 & 11.3 \\
\midrule
\lowersaxony & 41 & 100 & \textcolor{blue}{0.0} & 13.9 & 0.00 & 0.00 & 2 & 0.0 & 0.0 & 1.9 & 28.8 \\
\lowersaxony & 41 & 200 & \textcolor{blue}{0.9} & 309.1 & 0.00 & 0.00 & 6 & 0.0 & 0.0 & 1.8 & 25.4 \\
\lowersaxony & 41 & 300 & \textcolor{blue}{3.1} & 600.0 & 0.00$^*$ & $>100$ & 8 & $-$ & 0.0 & 1.7 & 22.9 \\
\lowersaxony & 82 & 100 & \textcolor{blue}{0.3} & 308.1 & 0.00 & 0.00 & 3 & 0.0 & 0.0 & 1.8 & 39.4 \\
\lowersaxony & 82 & 200 & \textcolor{blue}{8.0} & 600.0 & 0.00$^*$ & $>100$ & 8 & $-$ & 0.0 & 1.7 & 34.0 \\
\lowersaxony & 82 & 300 & \textcolor{blue}{23.4} & 600.0 & 0.00$^*$ & $>100$ & 9 & $-$ & 72.8 & 1.6 & 29.6 \\
\lowersaxony & 120 & 100 & \textcolor{blue}{1.3} & 600.0 & 0.00$^*$ & $>100$ & 3 & $-$ & 0.0 & 1.6 & 52.1 \\
\lowersaxony & 120 & 200 & \textcolor{blue}{264.7} & 600.0 & 0.00$^*$ & $>100$ & 21 & $-$ & 0.0 & 1.7 & 48.0 \\
\lowersaxony & 120 & 300 & \textcolor{blue}{600.0} & 600.0 & $>100$ & $>100$ & 16 & $-$ & 85.5 & 1.5 & 39.0 \\
\midrule
\athens & 32 & 100 & \textcolor{blue}{0.0} & 0.1 & 0.00 & 0.00 & 3 & 0.0 & 0.0 & 4.6 & 51.3 \\
\athens & 32 & 200 & \textcolor{blue}{0.1} & 4.5 & 0.00 & 0.00 & 3 & 0.0 & 0.0 & 3.0 & 40.8 \\
\athens & 32 & 300 & \textcolor{blue}{1.2} & 27.6 & 0.00 & 0.00 & 8 & 0.0 & 0.0 & 3.5 & 42.7 \\
\athens & 64 & 100 & \textcolor{blue}{0.1} & 0.7 & 0.00 & 0.00 & 3 & 0.0 & 0.0 & 4.0 & 72.7 \\
\athens & 64 & 200 & \textcolor{blue}{0.6} & 89.1 & 0.00 & 0.00 & 3 & 0.0 & 0.0 & 2.9 & 58.2 \\
\athens & 64 & 300 & \textcolor{blue}{47.0} & 398.2 & 0.00 & 0.00 & 14 & 0.0 & 0.0 & 3.3 & 61.8 \\
\athens & 129 & 100 & \textcolor{blue}{0.3} & 1.9 & 0.00 & 0.00 & 3 & 0.0 & 0.0 & 4.9 & 114.1 \\
\athens & 129 & 200 & \textcolor{blue}{2.4} & 238.4 & 0.00 & 0.00 & 3 & 0.0 & 0.0 & 3.4 & 92.1 \\
\athens & 129 & 300 & \textcolor{blue}{374.8} & 600.0 & 0.00$^*$ & 70.57 & 14 & 194.7 & 13.3 & 3.8 & 98.9 
\\	\bottomrule
		\end{tabular}}
\label{tab:TOwlambda0.1} 
\end{table}

\begin{table}[H]
	\centering
	\caption{
		Computational results for solving the \problem using either \frameworkNameShortLP (\emph{A}) or CPLEX (\emph{M}) with $\lambda = 0.25$.}
	\resizebox{\textwidth}{!}{
		\begin{tabular}{lrrrrrrrrrrr}\toprule
			\multicolumn{1}{l}{\textbf{Network}} &
			\multicolumn{1}{c}{\textbf{Pool}} &
			\multicolumn{1}{c}{\textbf{ODs}} &
			\multicolumn{2}{c}{\textbf{Time (min)}} &
			\multicolumn{2}{c}{\textbf{Gap (\%)}} &
			\multicolumn{1}{c}{\textbf{Iter.}} &
			\multicolumn{1}{c}{$\textbf{RLB (\%)}$} &
			\multicolumn{1}{c}{$\textbf{RUB (\%)}$} &
			\multicolumn{1}{c}{\textbf{Var. (\%)}} &
			\multicolumn{1}{c}{\textbf{Const. (\%)}} \\
			\cmidrule(lr){4-5}
			\cmidrule(lr){6-7} 
			\cmidrule(lr){8-8}
			\cmidrule(lr){9-9}
			\cmidrule(lr){10-10}
			\cmidrule(lr){11-11}
			\cmidrule(lr){12-12}
			& $|\linePool|$
			& $|\ODset|$
			& $A$ & M
			& $A$ & M
			& $A$ & $A$ & $A$ & $A$ & $A$ \\ \midrule 
			\midrule 
\sioux & 23 & 100 & \textcolor{blue}{10.8} & 88.8 & 0.00 & 0.00 &  25 & 0.0 & 0.0 & 10.6 & 37.5 \\
\sioux & 23 & 200 & \textcolor{blue}{173.0} & 234.7 & 0.00 & 0.00 &  29 & 0.0 & 0.0 & 17.7 & 46.1 \\
\sioux & 23 & 300 & 467.1 & \textcolor{blue}{420.9} & 0.00 & 0.00 &  30 & 0.0 & 0.0 & 19.7 & 47.5 \\
\sioux & 46 & 100 & \textcolor{blue}{30.3} & 600.0 & 0.00$^*$ & 17.25 & 27 & 20.9 & 0.0 & 9.8 & 37.3 \\
\sioux & 46 & 200 & \textcolor{blue}{479.7} & 600.0 & 0.00$^*$ & 76.74 & 32 & 317.0 & 3.0 & 17.5 & 47.1 \\
\sioux & 46 & 300 & \textcolor{blue}{600.0} & 600.0 & 11.19$^*$ & 80.19 & 25 & 331.2 & 3.8 & 15.5 & 42.5 \\
\sioux & 92 & 100 & \textcolor{blue}{600.0} & 600.0 & 13.73$^*$ & 83.99 & 32 & 417.5 & 3.9 & 5.6 & 30.2 \\
\sioux & 92 & 200 & \textcolor{blue}{600.0} & 600.0 & 51.80$^*$ & 94.93 & 14 & 780.6 & 7.3 & 3.4 & 22.7 \\
\sioux & 92 & 300 & \textcolor{blue}{600.0} & 600.0 & 65.84$^*$ & $>100$ & 8 & $-$ & 47.5 & 2.2 & 17.4 \\
\midrule
\grid & 30 & 100 & \textcolor{blue}{0.6} & 44.2 & 0.00 & 0.00 &  15 & 0.0 & 0.0 & 4.9 & 25.1 \\
\grid & 30 & 200 & \textcolor{blue}{21.1} & 140.8 & 0.00 & 0.00 &  21 & 0.0 & 0.0 & 8.3 & 30.2 \\
\grid & 30 & 300 & \textcolor{blue}{79.8} & 483.3 & 0.00 & 0.00 &  22 & 0.0 & 0.0 & 9.8 & 32.2 \\
\grid & 60 & 100 & \textcolor{blue}{231.3} & 600.0 & 0.00$^*$ & 63.50 & 47 & 146.3 & 10.1 & 6.5 & 29.9 \\
\grid & 60 & 200 & \textcolor{blue}{600.0} & 600.0 & 15.48$^*$ & 78.41 & 33 & 186.8 & 26.7 & 6.1 & 27.4 \\
\grid & 60 & 300 & \textcolor{blue}{600.0} & 600.0 & 24.57$^*$ & 87.54 & 20 & 290.5 & 35.5 & 4.1 & 21.6 \\
\grid & 123 & 100 & \textcolor{blue}{600.0} & 600.0 & 15.84$^*$ & 76.04 & 38 & 148.0 & 29.4 & 3.0 & 20.5 \\
\grid & 123 & 200 & \textcolor{blue}{600.0} & 600.0 & 38.99$^*$ & $>100$ & 17 & $-$ & 37.2 & 1.7 & 14.5 \\
\grid & 123 & 300 & \textcolor{blue}{600.0} & 600.0 & 47.12$^*$ & $>100$ & 10 & $-$ & 40.7 & 1.3 & 12.0 \\
\midrule
\lowersaxony & 41 & 100 & \textcolor{blue}{0.2} & 16.8 & 0.00 & 0.00 &  4 & 0.0 & 0.0 & 2.2 & 30.7 \\
\lowersaxony & 41 & 200 & \textcolor{blue}{11.8} & 310.7 & 0.00 & 0.00 &  18 & 0.0 & 0.0 & 2.9 & 29.9 \\
\lowersaxony & 41 & 300 & \textcolor{blue}{84.0} & 600.0 & 0.00$^*$ & $>100$ & 26 & $-$ & 10.0 & 3.5 & 29.5 \\
\lowersaxony & 82 & 100 & \textcolor{blue}{1.8} & 313.9 & 0.00 & 0.00 &  7 & 0.0 & 0.0 & 2.1 & 41.9 \\
\lowersaxony & 82 & 200 & \textcolor{blue}{90.1} & 600.0 & 0.00$^*$ & $>100$ & 21 & $-$ & 0.0 & 2.6 & 38.3 \\
\lowersaxony & 82 & 300 & \textcolor{blue}{520.7} & 600.0 & 0.00$^*$ & $>100$ & 30 & $-$ & 0.0 & 2.9 & 35.3 \\
\lowersaxony & 120 & 100 & \textcolor{blue}{13.6} & 600.0 & 0.00$^*$ & $>100$ & 10 & $-$ & 0.0 & 1.9 & 54.8 \\
\lowersaxony & 120 & 200 & 600.0 & \textcolor{blue}{600.0} & 24.29$^*$ & $>100$ & 20 & $-$ & -11.9 & 1.9 & 49.1 \\
\lowersaxony & 120 & 300 & \textcolor{blue}{600.0} & 600.0 & 35.05$^*$ & $>100$ & 10 & $-$ & 40.5 & 1.5 & 40.0 \\
\midrule
\athens & 32 & 100 & \textcolor{blue}{$<0.1$} & 0.1 & 0.00 & 0.00 &  6 & 0.0 & 0.0 & 5.1 & 53.9 \\
\athens & 32 & 200 & \textcolor{blue}{0.7} & 6.3 & 0.00 & 0.00 &  15 & 0.0 & 0.0 & 4.5 & 48.9 \\
\athens & 32 & 300 & \textcolor{blue}{10.2} & 69.9 & 0.00 & 0.00 &  25 & 0.0 & 0.0 & 5.7 & 54.5 \\
\athens & 64 & 100 & \textcolor{blue}{0.3} & 0.8 & 0.00 & 0.00 &  6 & 0.0 & 0.0 & 4.3 & 74.3 \\
\athens & 64 & 200 & \textcolor{blue}{34.4} & 121.8 & 0.00 & 0.00 &  29 & 0.0 & 0.0 & 4.1 & 65.6 \\
\athens & 64 & 300 & 600.0 & \textcolor{blue}{600.0} & 7.44 & 4.31$^*$ & 31 & -3.3 & 0.0 & 4.1 & 66.8 \\
\athens & 129 & 100 & \textcolor{blue}{1.3} & 2.9 & 0.00 & 0.00 &  6 & 0.0 & 0.0 & 5.1 & 116.6 \\
\athens & 129 & 200 & \textcolor{blue}{177.7} & 388.0 & 0.00 & 0.00 &  29 & 0.0 & 0.0 & 4.4 & 99.8 \\
\athens & 129 & 300 & \textcolor{blue}{600.0} & 600.0 & 12.43$^*$ & 44.36 & 8 & 56.2 & 0.8 & 3.6 & 98.2 \\	\bottomrule
		\end{tabular}
	}
	\label{tab:TOwlambda0.25} 
\end{table}

\begin{table}[H]
	\centering
	\caption{
		Computational results for solving the \problem using either \frameworkNameShortLP (\emph{A}) or CPLEX (\emph{M}) with $\lambda = 0.5$.}
	\resizebox{\textwidth}{!}{
		\begin{tabular}{lrrrrrrrrrrr}\toprule
			\multicolumn{1}{l}{\textbf{Network}} &
			\multicolumn{1}{c}{\textbf{Pool}} &
			\multicolumn{1}{c}{\textbf{ODs}} &
			\multicolumn{2}{c}{\textbf{Time (min)}} &
			\multicolumn{2}{c}{\textbf{Gap (\%)}} &
			\multicolumn{1}{c}{\textbf{Iter.}} &
			\multicolumn{1}{c}{$\textbf{RLB (\%)}$} &
			\multicolumn{1}{c}{$\textbf{RUB (\%)}$} &
			\multicolumn{1}{c}{\textbf{Var. (\%)}} &
			\multicolumn{1}{c}{\textbf{Const. (\%)}} \\
			\cmidrule(lr){4-5}
			\cmidrule(lr){6-7} 
			\cmidrule(lr){8-8}
			\cmidrule(lr){9-9}
			\cmidrule(lr){10-10}
			\cmidrule(lr){11-11}
			\cmidrule(lr){12-12}
			& $|\linePool|$
			& $|\ODset|$
			& $A$ & M
			& $A$ & M
			& $A$ & $A$ & $A$ & $A$ & $A$ \\ \midrule 
			\midrule 
\sioux & 23 & 100 & \textcolor{blue}{36.7} & 264.2 & 0.00 & 0.00 & 37 & 0.0 & 0.0 & 21.5 & 53.1 \\
\sioux & 23 & 200 & 178.7 & \textcolor{blue}{137.7} & 0.00 & 0.00 & 31 & 0.0 & 0.0 & 27.1 & 57.2 \\
\sioux & 23 & 300 & 392.2 & \textcolor{blue}{185.8} & 0.00 & 0.00 & 27 & 0.0 & 0.0 & 27.7 & 56.3 \\
\sioux & 46 & 100 & \textcolor{blue}{105.2} & 600.0 & 0.00$^*$ & 26.17 & 42 & 35.5 & 0.0 & 19.9 & 52.4 \\
\sioux & 46 & 200 & 567.9 & \textcolor{blue}{232.7} & 0.00 & 0.00 & 34 & 0.0 & 0.0 & 25.3 & 56.9 \\
\sioux & 46 & 300 & 600.0 & \textcolor{blue}{599.0} & 6.92 & 0.00$^*$ & 27 & -3.7 & -3.5 & 23.5 & 52.7 \\
\sioux & 92 & 100 & \textcolor{blue}{600.0} & 600.0 & 18.53$^*$ & 51.17 & 31 & 59.3 & 4.5 & 6.5 & 32.0 \\
\sioux & 92 & 200 & \textcolor{blue}{600.0} & 600.0 & 33.46$^*$ & 55.64 & 18 & 40.4 & 6.4 & 4.5 & 25.2 \\
\sioux & 92 & 300 & \textcolor{blue}{600.0} & 600.0 & 43.10$^*$ & 72.55 & 12 & 36.5 & 34.2 & 3.2 & 20.5 \\
\midrule
\grid & 30 & 100 & \textcolor{blue}{9.7} & 44.5 & 0.00 & 0.00 & 28 & 0.0 & 0.0 & 12.2 & 37.9 \\
\grid & 30 & 200 & \textcolor{blue}{76.3} & 79.4 & 0.00 & 0.00 & 31 & 0.0 & 0.0 & 19.1 & 45.1 \\
\grid & 30 & 300 & \textcolor{blue}{339.2} & 345.1 & 0.00 & 0.00 & 34 & 0.0 & 0.0 & 22.8 & 48.4 \\
\grid & 60 & 100 & 600.0 & \textcolor{blue}{600.0} & 4.56$^*$ & 41.14 & 63 & 62.4 & -0.1 & 13.6 & 41.4 \\
\grid & 60 & 200 & \textcolor{blue}{600.0} & 600.0 & 19.06$^*$ & 52.70 & 34 & 55.0 & 9.5 & 8.2 & 31.4 \\
\grid & 60 & 300 & \textcolor{blue}{600.0} & 600.0 & 22.80$^*$ & 59.62 & 22 & 53.0 & 20.0 & 5.5 & 25.0 \\
\grid & 123 & 100 & \textcolor{blue}{600.0} & 600.0 & 25.45$^*$ & 54.02 & 40 & 43.8 & 11.3 & 3.5 & 22.3 \\
\grid & 123 & 200 & \textcolor{blue}{600.0} & 600.0 & 33.13$^*$ & $>100$ & 19 & $-$ & 30.1 & 2.1 & 16.3 \\
\grid & 123 & 300 & \textcolor{blue}{600.0} & 600.0 & 39.02$^*$ & $>100$ & 11 & $-$ & 30.1 & 1.5 & 13.0 \\
\midrule
\lowersaxony & 41 & 100 & \textcolor{blue}{2.5} & 13.9 & 0.00 & 0.00 & 19 & 0.0 & 0.0 & 3.6 & 38.3 \\
\lowersaxony & 41 & 200 & \textcolor{blue}{333.8} & 600.0 & 0.00$^*$ & 56.16 & 52 & 128.1 & 0.0 & 8.0 & 44.6 \\
\lowersaxony & 41 & 300 & \textcolor{blue}{592.0} & 600.0 & 0.00$^*$ & 51.02 & 47 & 97.3 & 3.3 & 7.3 & 41.0 \\
\lowersaxony & 82 & 100 & \textcolor{blue}{28.5} & 600.0 & 0.00$^*$ & 47.50 & 25 & 90.5 & 0.0 & 3.4 & 51.1 \\
\lowersaxony & 82 & 200 & 600.0 & \textcolor{blue}{600.0} & 7.75$^*$ & 62.05 & 30 & 147.9 & -1.9 & 3.9 & 45.4 \\
\lowersaxony & 82 & 300 & \textcolor{blue}{600.0} & 600.0 & 10.46$^*$ & 64.60 & 23 & 149.4 & 1.4 & 3.1 & 37.5 \\
\lowersaxony & 120 & 100 & \textcolor{blue}{191.9} & 600.0 & 0.00$^*$ & 54.27 & 49 & 118.7 & 0.0 & 3.1 & 65.3 \\
\lowersaxony & 120 & 200 & 600.0 & \textcolor{blue}{600.0} & 23.82$^*$ & 71.63 & 17 & 192.2 & -8.8 & 2.1 & 52.7 \\
\lowersaxony & 120 & 300 & \textcolor{blue}{600.0} & 600.0 & 27.86$^*$ & $>100$ & 9 & $-$ & 29.1 & 1.5 & 40.5 \\
\midrule
\athens & 32 & 100 & \textcolor{blue}{0.0} & 0.1 & 0.00 & 0.00 & 10 & 0.0 & 0.0 & 6.1 & 62.0 \\
\athens & 32 & 200 & \textcolor{blue}{6.5} & 9.7 & 0.00 & 0.00 & 29 & 0.0 & 0.0 & 7.2 & 63.4 \\
\athens & 32 & 300 & \textcolor{blue}{36.4} & 80.5 & 0.00 & 0.00 & 33 & 0.0 & 0.0 & 9.8 & 71.3 \\
\athens & 64 & 100 & \textcolor{blue}{1.1} & 2.2 & 0.00 & 0.00 & 11 & 0.0 & 0.0 & 4.7 & 80.0 \\
\athens & 64 & 200 & 424.8 & \textcolor{blue}{156.4} & 0.00 & 0.00 & 70 & 0.0 & 0.0 & 7.1 & 86.3 \\
\athens & 64 & 300 & \textcolor{blue}{600.0} & 600.0 & 8.96$^*$ & 17.68 & 19 & 9.9 & 0.6 & 4.2 & 68.2 \\
\athens & 129 & 100 & \textcolor{blue}{4.7} & 5.8 & 0.00 & 0.00 & 11 & 0.0 & 0.0 & 5.6 & 125.0 \\
\athens & 129 & 200 & 600.0 & \textcolor{blue}{519.8} & 6.81 & 0.00$^*$ & 30 & -6.8 & 0.0 & 4.6 & 100.3 \\
\athens & 129 & 300 & \textcolor{blue}{600.0} & 600.0 & 12.45$^*$ & 22.33 & 6 & 12.7 & 0.0 & 3.7 & 98.6 \\
\bottomrule
		\end{tabular}
	}
	\label{tab:TOwlambda0.5} 
\end{table}

\begin{table}[H]
	\centering
	\caption{
		Computational results for solving the \problem using either \frameworkNameShortLP (\emph{A}) or CPLEX (\emph{M}) with $\lambda = 1$.}
	\resizebox{\textwidth}{!}{
		\begin{tabular}{lrrrrrrrrrrr}\toprule
			\multicolumn{1}{l}{\textbf{Network}} &
			\multicolumn{1}{c}{\textbf{Pool}} &
			\multicolumn{1}{c}{\textbf{ODs}} &
			\multicolumn{2}{c}{\textbf{Time (min)}} &
			\multicolumn{2}{c}{\textbf{Gap (\%)}} &
			\multicolumn{1}{c}{\textbf{Iter.}} &
			\multicolumn{1}{c}{$\textbf{RLB (\%)}$} &
			\multicolumn{1}{c}{$\textbf{RUB (\%)}$} &
			\multicolumn{1}{c}{\textbf{Var. (\%)}} &
			\multicolumn{1}{c}{\textbf{Const. (\%)}} \\
			\cmidrule(lr){4-5}
			\cmidrule(lr){6-7} 
			\cmidrule(lr){8-8}
			\cmidrule(lr){9-9}
			\cmidrule(lr){10-10}
			\cmidrule(lr){11-11}
			\cmidrule(lr){12-12}
			& $|\linePool|$
			& $|\ODset|$
			& $A$ & M
			& $A$ & M
			& $A$ & $A$ & $A$ & $A$ & $A$ \\ \midrule 
			\midrule 
\sioux & 23 & 100 & 34.0 & \textcolor{blue}{13.1} & 0.00 & 0.00 &  36 & 0.0 & 0.0 & 32.6 & 65.4 \\
\sioux & 23 & 200 & 167.4 & \textcolor{blue}{58.4} & 0.00 & 0.00 &  32 & 0.0 & 0.0 & 42.6 & 72.4 \\
\sioux & 23 & 300 & 433.0 & \textcolor{blue}{116.9} & 0.00 & 0.00 &  29 & 0.0 & 0.0 & 43.9 & 71.5 \\
\sioux & 46 & 100 & \textcolor{blue}{93.1} & 503.9 & 0.00 & 0.00 &  41 & 0.0 & 0.0 & 29.3 & 64.7 \\
\sioux & 46 & 200 & 418.1 & \textcolor{blue}{131.9} & 0.00 & 0.00 &  35 & 0.0 & 0.0 & 37.9 & 71.5 \\
\sioux & 46 & 300 & 600.0 & \textcolor{blue}{288.1} & 7.92 & 0.00$^*$ & 31 & -2.8 & -5.6 & 36.4 & 67.1 \\
\sioux & 92 & 100 & 600.0 & \textcolor{blue}{600.0} & 15.30$^*$ & 34.08 & 45 & 32.9 & -3.4 & 13.1 & 45.7 \\
\sioux & 92 & 200 & \textcolor{blue}{600.0} & 600.0 & 25.36$^*$ & 36.22 & 22 & 13.6 & 2.9 & 7.3 & 32.7 \\
\sioux & 92 & 300 & 600.0 & \textcolor{blue}{600.0} & 32.48$^*$ & 34.45 & 18 & 10.4 & -7.2 & 5.9 & 27.9 \\
\midrule
\grid & 30 & 100 & 30.3 & \textcolor{blue}{23.4} & 0.00 & 0.00 &  37 & 0.0 & 0.0 & 30.2 & 57.6 \\
\grid & 30 & 200 & 161.3 & \textcolor{blue}{40.4} & 0.00 & 0.00 &  39 & 0.0 & 0.0 & 36.9 & 61.3 \\
\grid & 30 & 300 & 509.7 & \textcolor{blue}{220.8} & 0.00 & 0.00 &  37 & 0.0 & 0.0 & 40.9 & 64.2 \\
\grid & 60 & 100 & 600.0 & \textcolor{blue}{600.0} & 9.84$^*$ & 34.84 & 69 & 41.2 & -2.0 & 22.6 & 51.8 \\
\grid & 60 & 200 & \textcolor{blue}{600.0} & 600.0 & 19.36$^*$ & 40.08 & 38 & 28.2 & 4.8 & 12.4 & 37.8 \\
\grid & 60 & 300 & \textcolor{blue}{600.0} & 600.0 & 25.93$^*$ & 47.32 & 25 & 21.0 & 13.9 & 8.2 & 30.3 \\
\grid & 123 & 100 & \textcolor{blue}{600.0} & 600.0 & 25.17$^*$ & 41.99 & 42 & 19.3 & 7.5 & 4.7 & 25.5 \\
\grid & 123 & 200 & \textcolor{blue}{600.0} & 600.0 & 32.90$^*$ & 58.75 & 22 & 19.3 & 26.7 & 2.7 & 18.3 \\
\grid & 123 & 300 & \textcolor{blue}{600.0} & 600.0 & 35.78$^*$ & $>100$ & 15 & $-$ & 27.4 & 2.0 & 15.1 \\
\midrule
\lowersaxony & 41 & 100 & 105.4 & \textcolor{blue}{95.5} & 0.00 & 0.00 &  59 & 0.0 & 0.0 & 12.2 & 65.9 \\
\lowersaxony & 41 & 200 & 600.0 & \textcolor{blue}{600.0} & 10.49$^*$ & 31.96 & 45 & 36.8 & -4.0 & 11.0 & 56.1 \\
\lowersaxony & 41 & 300 & 600.0 & \textcolor{blue}{600.0} & 13.45$^*$ & 34.35 & 30 & 34.9 & -2.3 & 6.9 & 41.4 \\
\lowersaxony & 82 & 100 & 600.0 & \textcolor{blue}{600.0} & 10.74$^*$ & 32.67 & 37 & 40.3 & -5.8 & 5.5 & 63.0 \\
\lowersaxony & 82 & 200 & 600.0 & \textcolor{blue}{600.0} & 16.83$^*$ & 34.81 & 23 & 32.0 & -3.5 & 4.0 & 48.0 \\
\lowersaxony & 82 & 300 & 600.0 & \textcolor{blue}{600.0} & 20.92$^*$ & 42.09 & 17 & 38.6 & -1.5 & 3.1 & 39.0 \\
\lowersaxony & 120 & 100 & 600.0 & \textcolor{blue}{600.0} & 19.89$^*$ & 40.82 & 27 & 46.4 & -8.1 & 3.2 & 70.4 \\
\lowersaxony & 120 & 200 & 600.0 & \textcolor{blue}{600.0} & 32.25$^*$ & 50.55 & 10 & 40.5 & -2.5 & 1.8 & 49.8 \\
\lowersaxony & 120 & 300 & \textcolor{blue}{600.0} & 600.0 & 35.39$^*$ & 64.44 & 6 & 46.3 & 19.5 & 1.4 & 39.3 \\
\midrule
\athens & 32 & 100 & \textcolor{blue}{0.3} & \textcolor{blue}{0.3} & 0.00 & 0.00 & 16 & 0.0 & 0.0 & 7.1 & 68.1 \\
\athens & 32 & 200 & \textcolor{blue}{12.5} & 17.6 & 0.00 & 0.00 &  32 & 0.0 & 0.0 & 11.7 & 79.8 \\
\athens & 32 & 300 & 114.7 & \textcolor{blue}{90.6} & 0.00 & 0.00 &  33 & 0.0 & 0.0 & 16.0 & 97.8 \\
\athens & 64 & 100 & 6.7 & \textcolor{blue}{3.1} & 0.00 & 0.00 &  25 & 0.0 & 0.0 & 5.5 & 85.6 \\
\athens & 64 & 200 & \textcolor{blue}{600.0} & 600.0 & 8.15$^*$ & 13.54 & 42 & 6.2 & 0.0 & 6.2 & 81.8 \\
\athens & 64 & 300 & 600.0 & \textcolor{blue}{600.0} & 14.53$^*$ & 20.79 & 14 & 10.1 & -2.1 & 4.5 & 69.5 \\
\athens & 129 & 100 & 42.0 & \textcolor{blue}{17.0} & 0.00 & 0.00 &  25 & 0.0 & 0.0 & 6.4 & 132.5 \\
\athens & 129 & 200 & 600.0 & \textcolor{blue}{600.0} & 12.42$^*$ & 19.70 & 12 & 9.7 & -0.6 & 4.2 & 98.4 \\
\athens & 129 & 300 & 600.0 & \textcolor{blue}{600.0} & 20.15$^*$ & 22.23 & 4 & 5.1 & -2.3 & 3.6 & 98.6 \\		\bottomrule
		\end{tabular}
	}
	\label{tab:TOwlambda1.0} 
\end{table}

%% file: Appendix/SocialMaximum.tex
\section{Objective value components over \texorpdfstring{$\lambda$}{lambda}}\label{sec:appSocialOptimum}

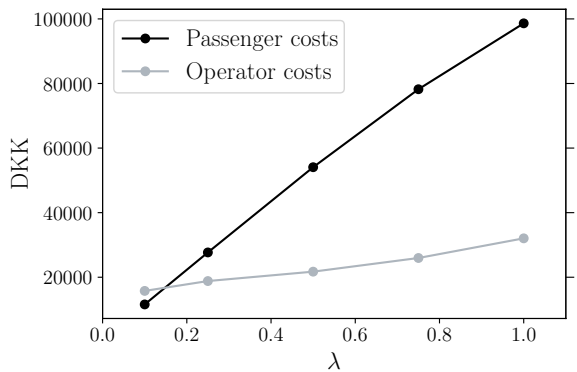
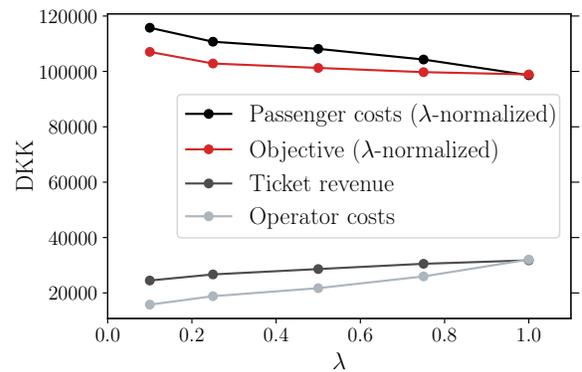
\begin{figure}[H]
	\centering
		\subfloat[t][Passenger and operator costs]{\resizebox{!}{0.32\textwidth}
{\input{figures/Comp/PIandOCVSLambda.pgf}}
			\label{fig:cost1}} \hfill
            		\subfloat[t][Objective value components normalized to a fixed value of $\lambda$]{\resizebox{!}{0.32\textwidth}
{\input{figures/Comp/CostsVSLambda.pgf}}
			\label{fig:cost2}} 
		\caption{Objective value and cost components averaged over all instances solved to optimality over varying values of $\lambda$.}
	\label{fig:Objective_VS_lambda}
\end{figure}

\autoref{fig:Objective_VS_lambda} illustrates how the objective function balances operator and passenger costs as the trade-off parameter $\lambda$ varies. In \autoref{fig:cost1}, the costs are shown scaled by the corresponding $\lambda$ values, reflecting their contribution to the actual objective function. In contrast, \autoref{fig:cost2} presents the same components normalized to a fixed $\lambda$ value to allow for consistent comparison across different $\lambda$ settings.
The results are averaged over the subset of instances that were solved to optimality for all $\lambda$ values. As shown in \autoref{fig:cost2}, the passenger and operator cost components operate on different magnitudes, and the total objective value decreases as $\lambda$ increases. 

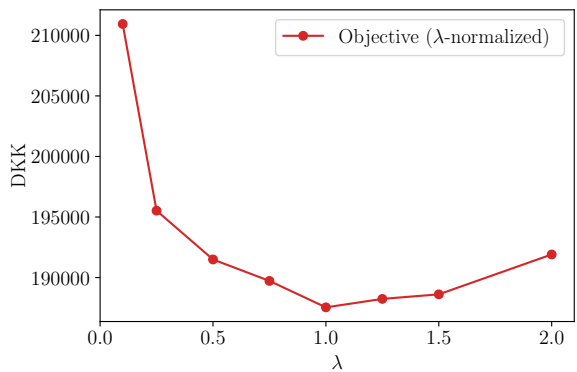
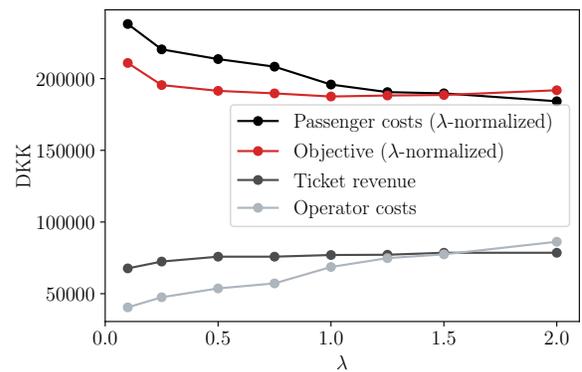
\begin{figure}[H]
	\centering
		\subfloat[t][Objective value]{\resizebox{!}{0.32\textwidth}
{\input{figures/Comp/Inst_ObjVSLambda.pgf}}
			\label{fig:InstObjVsL}} \hfill
            		\subfloat[t][Objective value components]{\resizebox{!}{0.32\textwidth}
{\input{figures/Comp/Inst_CostsVSLambda.pgf}}
			\label{fig:InstCostsVsL}} 
		\caption{Objective value and cost components for a representative instance using the \sioux network considering a line pool with 23 lines and 300 OD pairs. The objective value and passenger costs are normalized to a fixed value of $\lambda$.}
	\label{fig:siouxL23OD300}
\end{figure}

We consider the range of $\lambda \leq 1$ to be the most relevant for this study. To illustrate how the objective evolves beyond this range, \autoref{fig:InstCostsVsL} shows the $\lambda$-normalized objective values for a selected instance for the \sioux network as $\lambda$ varies from 0.1 to 2. We observe that $\lambda = 1$ corresponds to a socially optimal solution where passenger and operator costs are weighted proportionally to their actual magnitudes. 

%% file: figures/Comp/CostsVSLambda.pgf
\begingroup%
\makeatletter%
\begin{pgfpicture}%
\pgfpathrectangle{\pgfpointorigin}{\pgfqpoint{6.000000in}{4.000000in}}%
\pgfusepath{use as bounding box, clip}%
\begin{pgfscope}%
\pgfsetbuttcap%
\pgfsetmiterjoin%
\definecolor{currentfill}{rgb}{1.000000,1.000000,1.000000}%
\pgfsetfillcolor{currentfill}%
\pgfsetlinewidth{0.000000pt}%
\definecolor{currentstroke}{rgb}{1.000000,1.000000,1.000000}%
\pgfsetstrokecolor{currentstroke}%
\pgfsetdash{}{0pt}%
\pgfpathmoveto{\pgfqpoint{0.000000in}{0.000000in}}%
\pgfpathlineto{\pgfqpoint{6.000000in}{0.000000in}}%
\pgfpathlineto{\pgfqpoint{6.000000in}{4.000000in}}%
\pgfpathlineto{\pgfqpoint{0.000000in}{4.000000in}}%
\pgfpathlineto{\pgfqpoint{0.000000in}{0.000000in}}%
\pgfpathclose%
\pgfusepath{fill}%
\end{pgfscope}%
\begin{pgfscope}%
\pgfsetbuttcap%
\pgfsetmiterjoin%
\definecolor{currentfill}{rgb}{1.000000,1.000000,1.000000}%
\pgfsetfillcolor{currentfill}%
\pgfsetlinewidth{0.000000pt}%
\definecolor{currentstroke}{rgb}{0.000000,0.000000,0.000000}%
\pgfsetstrokecolor{currentstroke}%
\pgfsetstrokeopacity{0.000000}%
\pgfsetdash{}{0pt}%
\pgfpathmoveto{\pgfqpoint{1.128611in}{0.731763in}}%
\pgfpathlineto{\pgfqpoint{5.706667in}{0.731763in}}%
\pgfpathlineto{\pgfqpoint{5.706667in}{3.744191in}}%
\pgfpathlineto{\pgfqpoint{1.128611in}{3.744191in}}%
\pgfpathlineto{\pgfqpoint{1.128611in}{0.731763in}}%
\pgfpathclose%
\pgfusepath{fill}%
\end{pgfscope}%
\begin{pgfscope}%
\pgfsetbuttcap%
\pgfsetroundjoin%
\definecolor{currentfill}{rgb}{0.000000,0.000000,0.000000}%
\pgfsetfillcolor{currentfill}%
\pgfsetlinewidth{0.803000pt}%
\definecolor{currentstroke}{rgb}{0.000000,0.000000,0.000000}%
\pgfsetstrokecolor{currentstroke}%
\pgfsetdash{}{0pt}%
\pgfsys@defobject{currentmarker}{\pgfqpoint{0.000000in}{-0.048611in}}{\pgfqpoint{0.000000in}{0.000000in}}{%
\pgfpathmoveto{\pgfqpoint{0.000000in}{0.000000in}}%
\pgfpathlineto{\pgfqpoint{0.000000in}{-0.048611in}}%
\pgfusepath{stroke,fill}%
}%
\begin{pgfscope}%
\pgfsys@transformshift{1.128611in}{0.731763in}%
\pgfsys@useobject{currentmarker}{}%
\end{pgfscope}%
\end{pgfscope}%
\begin{pgfscope}%
\definecolor{textcolor}{rgb}{0.000000,0.000000,0.000000}%
\pgfsetstrokecolor{textcolor}%
\pgfsetfillcolor{textcolor}%
\pgftext[x=1.128611in,y=0.634541in,,top]{\color{textcolor}\rmfamily\fontsize{14.000000}{16.800000}\selectfont \(\displaystyle {0.0}\)}%
\end{pgfscope}%
\begin{pgfscope}%
\pgfsetbuttcap%
\pgfsetroundjoin%
\definecolor{currentfill}{rgb}{0.000000,0.000000,0.000000}%
\pgfsetfillcolor{currentfill}%
\pgfsetlinewidth{0.803000pt}%
\definecolor{currentstroke}{rgb}{0.000000,0.000000,0.000000}%
\pgfsetstrokecolor{currentstroke}%
\pgfsetdash{}{0pt}%
\pgfsys@defobject{currentmarker}{\pgfqpoint{0.000000in}{-0.048611in}}{\pgfqpoint{0.000000in}{0.000000in}}{%
\pgfpathmoveto{\pgfqpoint{0.000000in}{0.000000in}}%
\pgfpathlineto{\pgfqpoint{0.000000in}{-0.048611in}}%
\pgfusepath{stroke,fill}%
}%
\begin{pgfscope}%
\pgfsys@transformshift{1.960985in}{0.731763in}%
\pgfsys@useobject{currentmarker}{}%
\end{pgfscope}%
\end{pgfscope}%
\begin{pgfscope}%
\definecolor{textcolor}{rgb}{0.000000,0.000000,0.000000}%
\pgfsetstrokecolor{textcolor}%
\pgfsetfillcolor{textcolor}%
\pgftext[x=1.960985in,y=0.634541in,,top]{\color{textcolor}\rmfamily\fontsize{14.000000}{16.800000}\selectfont \(\displaystyle {0.2}\)}%
\end{pgfscope}%
\begin{pgfscope}%
\pgfsetbuttcap%
\pgfsetroundjoin%
\definecolor{currentfill}{rgb}{0.000000,0.000000,0.000000}%
\pgfsetfillcolor{currentfill}%
\pgfsetlinewidth{0.803000pt}%
\definecolor{currentstroke}{rgb}{0.000000,0.000000,0.000000}%
\pgfsetstrokecolor{currentstroke}%
\pgfsetdash{}{0pt}%
\pgfsys@defobject{currentmarker}{\pgfqpoint{0.000000in}{-0.048611in}}{\pgfqpoint{0.000000in}{0.000000in}}{%
\pgfpathmoveto{\pgfqpoint{0.000000in}{0.000000in}}%
\pgfpathlineto{\pgfqpoint{0.000000in}{-0.048611in}}%
\pgfusepath{stroke,fill}%
}%
\begin{pgfscope}%
\pgfsys@transformshift{2.793359in}{0.731763in}%
\pgfsys@useobject{currentmarker}{}%
\end{pgfscope}%
\end{pgfscope}%
\begin{pgfscope}%
\definecolor{textcolor}{rgb}{0.000000,0.000000,0.000000}%
\pgfsetstrokecolor{textcolor}%
\pgfsetfillcolor{textcolor}%
\pgftext[x=2.793359in,y=0.634541in,,top]{\color{textcolor}\rmfamily\fontsize{14.000000}{16.800000}\selectfont \(\displaystyle {0.4}\)}%
\end{pgfscope}%
\begin{pgfscope}%
\pgfsetbuttcap%
\pgfsetroundjoin%
\definecolor{currentfill}{rgb}{0.000000,0.000000,0.000000}%
\pgfsetfillcolor{currentfill}%
\pgfsetlinewidth{0.803000pt}%
\definecolor{currentstroke}{rgb}{0.000000,0.000000,0.000000}%
\pgfsetstrokecolor{currentstroke}%
\pgfsetdash{}{0pt}%
\pgfsys@defobject{currentmarker}{\pgfqpoint{0.000000in}{-0.048611in}}{\pgfqpoint{0.000000in}{0.000000in}}{%
\pgfpathmoveto{\pgfqpoint{0.000000in}{0.000000in}}%
\pgfpathlineto{\pgfqpoint{0.000000in}{-0.048611in}}%
\pgfusepath{stroke,fill}%
}%
\begin{pgfscope}%
\pgfsys@transformshift{3.625732in}{0.731763in}%
\pgfsys@useobject{currentmarker}{}%
\end{pgfscope}%
\end{pgfscope}%
\begin{pgfscope}%
\definecolor{textcolor}{rgb}{0.000000,0.000000,0.000000}%
\pgfsetstrokecolor{textcolor}%
\pgfsetfillcolor{textcolor}%
\pgftext[x=3.625732in,y=0.634541in,,top]{\color{textcolor}\rmfamily\fontsize{14.000000}{16.800000}\selectfont \(\displaystyle {0.6}\)}%
\end{pgfscope}%
\begin{pgfscope}%
\pgfsetbuttcap%
\pgfsetroundjoin%
\definecolor{currentfill}{rgb}{0.000000,0.000000,0.000000}%
\pgfsetfillcolor{currentfill}%
\pgfsetlinewidth{0.803000pt}%
\definecolor{currentstroke}{rgb}{0.000000,0.000000,0.000000}%
\pgfsetstrokecolor{currentstroke}%
\pgfsetdash{}{0pt}%
\pgfsys@defobject{currentmarker}{\pgfqpoint{0.000000in}{-0.048611in}}{\pgfqpoint{0.000000in}{0.000000in}}{%
\pgfpathmoveto{\pgfqpoint{0.000000in}{0.000000in}}%
\pgfpathlineto{\pgfqpoint{0.000000in}{-0.048611in}}%
\pgfusepath{stroke,fill}%
}%
\begin{pgfscope}%
\pgfsys@transformshift{4.458106in}{0.731763in}%
\pgfsys@useobject{currentmarker}{}%
\end{pgfscope}%
\end{pgfscope}%
\begin{pgfscope}%
\definecolor{textcolor}{rgb}{0.000000,0.000000,0.000000}%
\pgfsetstrokecolor{textcolor}%
\pgfsetfillcolor{textcolor}%
\pgftext[x=4.458106in,y=0.634541in,,top]{\color{textcolor}\rmfamily\fontsize{14.000000}{16.800000}\selectfont \(\displaystyle {0.8}\)}%
\end{pgfscope}%
\begin{pgfscope}%
\pgfsetbuttcap%
\pgfsetroundjoin%
\definecolor{currentfill}{rgb}{0.000000,0.000000,0.000000}%
\pgfsetfillcolor{currentfill}%
\pgfsetlinewidth{0.803000pt}%
\definecolor{currentstroke}{rgb}{0.000000,0.000000,0.000000}%
\pgfsetstrokecolor{currentstroke}%
\pgfsetdash{}{0pt}%
\pgfsys@defobject{currentmarker}{\pgfqpoint{0.000000in}{-0.048611in}}{\pgfqpoint{0.000000in}{0.000000in}}{%
\pgfpathmoveto{\pgfqpoint{0.000000in}{0.000000in}}%
\pgfpathlineto{\pgfqpoint{0.000000in}{-0.048611in}}%
\pgfusepath{stroke,fill}%
}%
\begin{pgfscope}%
\pgfsys@transformshift{5.290480in}{0.731763in}%
\pgfsys@useobject{currentmarker}{}%
\end{pgfscope}%
\end{pgfscope}%
\begin{pgfscope}%
\definecolor{textcolor}{rgb}{0.000000,0.000000,0.000000}%
\pgfsetstrokecolor{textcolor}%
\pgfsetfillcolor{textcolor}%
\pgftext[x=5.290480in,y=0.634541in,,top]{\color{textcolor}\rmfamily\fontsize{14.000000}{16.800000}\selectfont \(\displaystyle {1.0}\)}%
\end{pgfscope}%
\begin{pgfscope}%
\definecolor{textcolor}{rgb}{0.000000,0.000000,0.000000}%
\pgfsetstrokecolor{textcolor}%
\pgfsetfillcolor{textcolor}%
\pgftext[x=3.417639in,y=0.401208in,,top]{\color{textcolor}\rmfamily\fontsize{16.000000}{19.200000}\selectfont \(\displaystyle \lambda\)}%
\end{pgfscope}%
\begin{pgfscope}%
\pgfsetbuttcap%
\pgfsetroundjoin%
\definecolor{currentfill}{rgb}{0.000000,0.000000,0.000000}%
\pgfsetfillcolor{currentfill}%
\pgfsetlinewidth{0.803000pt}%
\definecolor{currentstroke}{rgb}{0.000000,0.000000,0.000000}%
\pgfsetstrokecolor{currentstroke}%
\pgfsetdash{}{0pt}%
\pgfsys@defobject{currentmarker}{\pgfqpoint{-0.048611in}{0.000000in}}{\pgfqpoint{-0.000000in}{0.000000in}}{%
\pgfpathmoveto{\pgfqpoint{-0.000000in}{0.000000in}}%
\pgfpathlineto{\pgfqpoint{-0.048611in}{0.000000in}}%
\pgfusepath{stroke,fill}%
}%
\begin{pgfscope}%
\pgfsys@transformshift{1.128611in}{0.984625in}%
\pgfsys@useobject{currentmarker}{}%
\end{pgfscope}%
\end{pgfscope}%
\begin{pgfscope}%
\definecolor{textcolor}{rgb}{0.000000,0.000000,0.000000}%
\pgfsetstrokecolor{textcolor}%
\pgfsetfillcolor{textcolor}%
\pgftext[x=0.541812in, y=0.915181in, left, base]{\color{textcolor}\rmfamily\fontsize{14.000000}{16.800000}\selectfont \(\displaystyle {20000}\)}%
\end{pgfscope}%
\begin{pgfscope}%
\pgfsetbuttcap%
\pgfsetroundjoin%
\definecolor{currentfill}{rgb}{0.000000,0.000000,0.000000}%
\pgfsetfillcolor{currentfill}%
\pgfsetlinewidth{0.803000pt}%
\definecolor{currentstroke}{rgb}{0.000000,0.000000,0.000000}%
\pgfsetstrokecolor{currentstroke}%
\pgfsetdash{}{0pt}%
\pgfsys@defobject{currentmarker}{\pgfqpoint{-0.048611in}{0.000000in}}{\pgfqpoint{-0.000000in}{0.000000in}}{%
\pgfpathmoveto{\pgfqpoint{-0.000000in}{0.000000in}}%
\pgfpathlineto{\pgfqpoint{-0.048611in}{0.000000in}}%
\pgfusepath{stroke,fill}%
}%
\begin{pgfscope}%
\pgfsys@transformshift{1.128611in}{1.532342in}%
\pgfsys@useobject{currentmarker}{}%
\end{pgfscope}%
\end{pgfscope}%
\begin{pgfscope}%
\definecolor{textcolor}{rgb}{0.000000,0.000000,0.000000}%
\pgfsetstrokecolor{textcolor}%
\pgfsetfillcolor{textcolor}%
\pgftext[x=0.541812in, y=1.462897in, left, base]{\color{textcolor}\rmfamily\fontsize{14.000000}{16.800000}\selectfont \(\displaystyle {40000}\)}%
\end{pgfscope}%
\begin{pgfscope}%
\pgfsetbuttcap%
\pgfsetroundjoin%
\definecolor{currentfill}{rgb}{0.000000,0.000000,0.000000}%
\pgfsetfillcolor{currentfill}%
\pgfsetlinewidth{0.803000pt}%
\definecolor{currentstroke}{rgb}{0.000000,0.000000,0.000000}%
\pgfsetstrokecolor{currentstroke}%
\pgfsetdash{}{0pt}%
\pgfsys@defobject{currentmarker}{\pgfqpoint{-0.048611in}{0.000000in}}{\pgfqpoint{-0.000000in}{0.000000in}}{%
\pgfpathmoveto{\pgfqpoint{-0.000000in}{0.000000in}}%
\pgfpathlineto{\pgfqpoint{-0.048611in}{0.000000in}}%
\pgfusepath{stroke,fill}%
}%
\begin{pgfscope}%
\pgfsys@transformshift{1.128611in}{2.080058in}%
\pgfsys@useobject{currentmarker}{}%
\end{pgfscope}%
\end{pgfscope}%
\begin{pgfscope}%
\definecolor{textcolor}{rgb}{0.000000,0.000000,0.000000}%
\pgfsetstrokecolor{textcolor}%
\pgfsetfillcolor{textcolor}%
\pgftext[x=0.541812in, y=2.010614in, left, base]{\color{textcolor}\rmfamily\fontsize{14.000000}{16.800000}\selectfont \(\displaystyle {60000}\)}%
\end{pgfscope}%
\begin{pgfscope}%
\pgfsetbuttcap%
\pgfsetroundjoin%
\definecolor{currentfill}{rgb}{0.000000,0.000000,0.000000}%
\pgfsetfillcolor{currentfill}%
\pgfsetlinewidth{0.803000pt}%
\definecolor{currentstroke}{rgb}{0.000000,0.000000,0.000000}%
\pgfsetstrokecolor{currentstroke}%
\pgfsetdash{}{0pt}%
\pgfsys@defobject{currentmarker}{\pgfqpoint{-0.048611in}{0.000000in}}{\pgfqpoint{-0.000000in}{0.000000in}}{%
\pgfpathmoveto{\pgfqpoint{-0.000000in}{0.000000in}}%
\pgfpathlineto{\pgfqpoint{-0.048611in}{0.000000in}}%
\pgfusepath{stroke,fill}%
}%
\begin{pgfscope}%
\pgfsys@transformshift{1.128611in}{2.627775in}%
\pgfsys@useobject{currentmarker}{}%
\end{pgfscope}%
\end{pgfscope}%
\begin{pgfscope}%
\definecolor{textcolor}{rgb}{0.000000,0.000000,0.000000}%
\pgfsetstrokecolor{textcolor}%
\pgfsetfillcolor{textcolor}%
\pgftext[x=0.541812in, y=2.558330in, left, base]{\color{textcolor}\rmfamily\fontsize{14.000000}{16.800000}\selectfont \(\displaystyle {80000}\)}%
\end{pgfscope}%
\begin{pgfscope}%
\pgfsetbuttcap%
\pgfsetroundjoin%
\definecolor{currentfill}{rgb}{0.000000,0.000000,0.000000}%
\pgfsetfillcolor{currentfill}%
\pgfsetlinewidth{0.803000pt}%
\definecolor{currentstroke}{rgb}{0.000000,0.000000,0.000000}%
\pgfsetstrokecolor{currentstroke}%
\pgfsetdash{}{0pt}%
\pgfsys@defobject{currentmarker}{\pgfqpoint{-0.048611in}{0.000000in}}{\pgfqpoint{-0.000000in}{0.000000in}}{%
\pgfpathmoveto{\pgfqpoint{-0.000000in}{0.000000in}}%
\pgfpathlineto{\pgfqpoint{-0.048611in}{0.000000in}}%
\pgfusepath{stroke,fill}%
}%
\begin{pgfscope}%
\pgfsys@transformshift{1.128611in}{3.175491in}%
\pgfsys@useobject{currentmarker}{}%
\end{pgfscope}%
\end{pgfscope}%
\begin{pgfscope}%
\definecolor{textcolor}{rgb}{0.000000,0.000000,0.000000}%
\pgfsetstrokecolor{textcolor}%
\pgfsetfillcolor{textcolor}%
\pgftext[x=0.443896in, y=3.106047in, left, base]{\color{textcolor}\rmfamily\fontsize{14.000000}{16.800000}\selectfont \(\displaystyle {100000}\)}%
\end{pgfscope}%
\begin{pgfscope}%
\pgfsetbuttcap%
\pgfsetroundjoin%
\definecolor{currentfill}{rgb}{0.000000,0.000000,0.000000}%
\pgfsetfillcolor{currentfill}%
\pgfsetlinewidth{0.803000pt}%
\definecolor{currentstroke}{rgb}{0.000000,0.000000,0.000000}%
\pgfsetstrokecolor{currentstroke}%
\pgfsetdash{}{0pt}%
\pgfsys@defobject{currentmarker}{\pgfqpoint{-0.048611in}{0.000000in}}{\pgfqpoint{-0.000000in}{0.000000in}}{%
\pgfpathmoveto{\pgfqpoint{-0.000000in}{0.000000in}}%
\pgfpathlineto{\pgfqpoint{-0.048611in}{0.000000in}}%
\pgfusepath{stroke,fill}%
}%
\begin{pgfscope}%
\pgfsys@transformshift{1.128611in}{3.723208in}%
\pgfsys@useobject{currentmarker}{}%
\end{pgfscope}%
\end{pgfscope}%
\begin{pgfscope}%
\definecolor{textcolor}{rgb}{0.000000,0.000000,0.000000}%
\pgfsetstrokecolor{textcolor}%
\pgfsetfillcolor{textcolor}%
\pgftext[x=0.443896in, y=3.653763in, left, base]{\color{textcolor}\rmfamily\fontsize{14.000000}{16.800000}\selectfont \(\displaystyle {120000}\)}%
\end{pgfscope}%
\begin{pgfscope}%
\definecolor{textcolor}{rgb}{0.000000,0.000000,0.000000}%
\pgfsetstrokecolor{textcolor}%
\pgfsetfillcolor{textcolor}%
\pgftext[x=0.388341in,y=2.237977in,,bottom,rotate=90.000000]{\color{textcolor}\rmfamily\fontsize{16.000000}{19.200000}\selectfont DKK}%
\end{pgfscope}%
\begin{pgfscope}%
\pgfpathrectangle{\pgfqpoint{1.128611in}{0.731763in}}{\pgfqpoint{4.578056in}{3.012428in}}%
\pgfusepath{clip}%
\pgfsetrectcap%
\pgfsetroundjoin%
\pgfsetlinewidth{1.505625pt}%
\definecolor{currentstroke}{rgb}{0.000000,0.000000,0.000000}%
\pgfsetstrokecolor{currentstroke}%
\pgfsetdash{}{0pt}%
\pgfpathmoveto{\pgfqpoint{1.544798in}{3.607263in}}%
\pgfpathlineto{\pgfqpoint{2.169078in}{3.469090in}}%
\pgfpathlineto{\pgfqpoint{3.209545in}{3.398695in}}%
\pgfpathlineto{\pgfqpoint{4.250013in}{3.293174in}}%
\pgfpathlineto{\pgfqpoint{5.290480in}{3.137993in}}%
\pgfusepath{stroke}%
\end{pgfscope}%
\begin{pgfscope}%
\pgfpathrectangle{\pgfqpoint{1.128611in}{0.731763in}}{\pgfqpoint{4.578056in}{3.012428in}}%
\pgfusepath{clip}%
\pgfsetbuttcap%
\pgfsetroundjoin%
\definecolor{currentfill}{rgb}{0.000000,0.000000,0.000000}%
\pgfsetfillcolor{currentfill}%
\pgfsetlinewidth{1.003750pt}%
\definecolor{currentstroke}{rgb}{0.000000,0.000000,0.000000}%
\pgfsetstrokecolor{currentstroke}%
\pgfsetdash{}{0pt}%
\pgfsys@defobject{currentmarker}{\pgfqpoint{-0.041667in}{-0.041667in}}{\pgfqpoint{0.041667in}{0.041667in}}{%
\pgfpathmoveto{\pgfqpoint{0.000000in}{-0.041667in}}%
\pgfpathcurveto{\pgfqpoint{0.011050in}{-0.041667in}}{\pgfqpoint{0.021649in}{-0.037276in}}{\pgfqpoint{0.029463in}{-0.029463in}}%
\pgfpathcurveto{\pgfqpoint{0.037276in}{-0.021649in}}{\pgfqpoint{0.041667in}{-0.011050in}}{\pgfqpoint{0.041667in}{0.000000in}}%
\pgfpathcurveto{\pgfqpoint{0.041667in}{0.011050in}}{\pgfqpoint{0.037276in}{0.021649in}}{\pgfqpoint{0.029463in}{0.029463in}}%
\pgfpathcurveto{\pgfqpoint{0.021649in}{0.037276in}}{\pgfqpoint{0.011050in}{0.041667in}}{\pgfqpoint{0.000000in}{0.041667in}}%
\pgfpathcurveto{\pgfqpoint{-0.011050in}{0.041667in}}{\pgfqpoint{-0.021649in}{0.037276in}}{\pgfqpoint{-0.029463in}{0.029463in}}%
\pgfpathcurveto{\pgfqpoint{-0.037276in}{0.021649in}}{\pgfqpoint{-0.041667in}{0.011050in}}{\pgfqpoint{-0.041667in}{0.000000in}}%
\pgfpathcurveto{\pgfqpoint{-0.041667in}{-0.011050in}}{\pgfqpoint{-0.037276in}{-0.021649in}}{\pgfqpoint{-0.029463in}{-0.029463in}}%
\pgfpathcurveto{\pgfqpoint{-0.021649in}{-0.037276in}}{\pgfqpoint{-0.011050in}{-0.041667in}}{\pgfqpoint{0.000000in}{-0.041667in}}%
\pgfpathlineto{\pgfqpoint{0.000000in}{-0.041667in}}%
\pgfpathclose%
\pgfusepath{stroke,fill}%
}%
\begin{pgfscope}%
\pgfsys@transformshift{1.544798in}{3.607263in}%
\pgfsys@useobject{currentmarker}{}%
\end{pgfscope}%
\begin{pgfscope}%
\pgfsys@transformshift{2.169078in}{3.469090in}%
\pgfsys@useobject{currentmarker}{}%
\end{pgfscope}%
\begin{pgfscope}%
\pgfsys@transformshift{3.209545in}{3.398695in}%
\pgfsys@useobject{currentmarker}{}%
\end{pgfscope}%
\begin{pgfscope}%
\pgfsys@transformshift{4.250013in}{3.293174in}%
\pgfsys@useobject{currentmarker}{}%
\end{pgfscope}%
\begin{pgfscope}%
\pgfsys@transformshift{5.290480in}{3.137993in}%
\pgfsys@useobject{currentmarker}{}%
\end{pgfscope}%
\end{pgfscope}%
\begin{pgfscope}%
\pgfpathrectangle{\pgfqpoint{1.128611in}{0.731763in}}{\pgfqpoint{4.578056in}{3.012428in}}%
\pgfusepath{clip}%
\pgfsetrectcap%
\pgfsetroundjoin%
\pgfsetlinewidth{1.505625pt}%
\definecolor{currentstroke}{rgb}{0.839216,0.152941,0.156863}%
\pgfsetstrokecolor{currentstroke}%
\pgfsetdash{}{0pt}%
\pgfpathmoveto{\pgfqpoint{1.544798in}{3.368518in}}%
\pgfpathlineto{\pgfqpoint{2.169078in}{3.253697in}}%
\pgfpathlineto{\pgfqpoint{3.209545in}{3.209969in}}%
\pgfpathlineto{\pgfqpoint{4.250013in}{3.168277in}}%
\pgfpathlineto{\pgfqpoint{5.290480in}{3.145682in}}%
\pgfusepath{stroke}%
\end{pgfscope}%
\begin{pgfscope}%
\pgfpathrectangle{\pgfqpoint{1.128611in}{0.731763in}}{\pgfqpoint{4.578056in}{3.012428in}}%
\pgfusepath{clip}%
\pgfsetbuttcap%
\pgfsetroundjoin%
\definecolor{currentfill}{rgb}{0.839216,0.152941,0.156863}%
\pgfsetfillcolor{currentfill}%
\pgfsetlinewidth{1.003750pt}%
\definecolor{currentstroke}{rgb}{0.839216,0.152941,0.156863}%
\pgfsetstrokecolor{currentstroke}%
\pgfsetdash{}{0pt}%
\pgfsys@defobject{currentmarker}{\pgfqpoint{-0.041667in}{-0.041667in}}{\pgfqpoint{0.041667in}{0.041667in}}{%
\pgfpathmoveto{\pgfqpoint{0.000000in}{-0.041667in}}%
\pgfpathcurveto{\pgfqpoint{0.011050in}{-0.041667in}}{\pgfqpoint{0.021649in}{-0.037276in}}{\pgfqpoint{0.029463in}{-0.029463in}}%
\pgfpathcurveto{\pgfqpoint{0.037276in}{-0.021649in}}{\pgfqpoint{0.041667in}{-0.011050in}}{\pgfqpoint{0.041667in}{0.000000in}}%
\pgfpathcurveto{\pgfqpoint{0.041667in}{0.011050in}}{\pgfqpoint{0.037276in}{0.021649in}}{\pgfqpoint{0.029463in}{0.029463in}}%
\pgfpathcurveto{\pgfqpoint{0.021649in}{0.037276in}}{\pgfqpoint{0.011050in}{0.041667in}}{\pgfqpoint{0.000000in}{0.041667in}}%
\pgfpathcurveto{\pgfqpoint{-0.011050in}{0.041667in}}{\pgfqpoint{-0.021649in}{0.037276in}}{\pgfqpoint{-0.029463in}{0.029463in}}%
\pgfpathcurveto{\pgfqpoint{-0.037276in}{0.021649in}}{\pgfqpoint{-0.041667in}{0.011050in}}{\pgfqpoint{-0.041667in}{0.000000in}}%
\pgfpathcurveto{\pgfqpoint{-0.041667in}{-0.011050in}}{\pgfqpoint{-0.037276in}{-0.021649in}}{\pgfqpoint{-0.029463in}{-0.029463in}}%
\pgfpathcurveto{\pgfqpoint{-0.021649in}{-0.037276in}}{\pgfqpoint{-0.011050in}{-0.041667in}}{\pgfqpoint{0.000000in}{-0.041667in}}%
\pgfpathlineto{\pgfqpoint{0.000000in}{-0.041667in}}%
\pgfpathclose%
\pgfusepath{stroke,fill}%
}%
\begin{pgfscope}%
\pgfsys@transformshift{1.544798in}{3.368518in}%
\pgfsys@useobject{currentmarker}{}%
\end{pgfscope}%
\begin{pgfscope}%
\pgfsys@transformshift{2.169078in}{3.253697in}%
\pgfsys@useobject{currentmarker}{}%
\end{pgfscope}%
\begin{pgfscope}%
\pgfsys@transformshift{3.209545in}{3.209969in}%
\pgfsys@useobject{currentmarker}{}%
\end{pgfscope}%
\begin{pgfscope}%
\pgfsys@transformshift{4.250013in}{3.168277in}%
\pgfsys@useobject{currentmarker}{}%
\end{pgfscope}%
\begin{pgfscope}%
\pgfsys@transformshift{5.290480in}{3.145682in}%
\pgfsys@useobject{currentmarker}{}%
\end{pgfscope}%
\end{pgfscope}%
\begin{pgfscope}%
\pgfpathrectangle{\pgfqpoint{1.128611in}{0.731763in}}{\pgfqpoint{4.578056in}{3.012428in}}%
\pgfusepath{clip}%
\pgfsetrectcap%
\pgfsetroundjoin%
\pgfsetlinewidth{1.505625pt}%
\definecolor{currentstroke}{rgb}{0.301961,0.301961,0.301961}%
\pgfsetstrokecolor{currentstroke}%
\pgfsetdash{}{0pt}%
\pgfpathmoveto{\pgfqpoint{1.544798in}{1.107437in}}%
\pgfpathlineto{\pgfqpoint{2.169078in}{1.167678in}}%
\pgfpathlineto{\pgfqpoint{3.209545in}{1.220652in}}%
\pgfpathlineto{\pgfqpoint{4.250013in}{1.272742in}}%
\pgfpathlineto{\pgfqpoint{5.290480in}{1.306845in}}%
\pgfusepath{stroke}%
\end{pgfscope}%
\begin{pgfscope}%
\pgfpathrectangle{\pgfqpoint{1.128611in}{0.731763in}}{\pgfqpoint{4.578056in}{3.012428in}}%
\pgfusepath{clip}%
\pgfsetbuttcap%
\pgfsetroundjoin%
\definecolor{currentfill}{rgb}{0.301961,0.301961,0.301961}%
\pgfsetfillcolor{currentfill}%
\pgfsetlinewidth{1.003750pt}%
\definecolor{currentstroke}{rgb}{0.301961,0.301961,0.301961}%
\pgfsetstrokecolor{currentstroke}%
\pgfsetdash{}{0pt}%
\pgfsys@defobject{currentmarker}{\pgfqpoint{-0.041667in}{-0.041667in}}{\pgfqpoint{0.041667in}{0.041667in}}{%
\pgfpathmoveto{\pgfqpoint{0.000000in}{-0.041667in}}%
\pgfpathcurveto{\pgfqpoint{0.011050in}{-0.041667in}}{\pgfqpoint{0.021649in}{-0.037276in}}{\pgfqpoint{0.029463in}{-0.029463in}}%
\pgfpathcurveto{\pgfqpoint{0.037276in}{-0.021649in}}{\pgfqpoint{0.041667in}{-0.011050in}}{\pgfqpoint{0.041667in}{0.000000in}}%
\pgfpathcurveto{\pgfqpoint{0.041667in}{0.011050in}}{\pgfqpoint{0.037276in}{0.021649in}}{\pgfqpoint{0.029463in}{0.029463in}}%
\pgfpathcurveto{\pgfqpoint{0.021649in}{0.037276in}}{\pgfqpoint{0.011050in}{0.041667in}}{\pgfqpoint{0.000000in}{0.041667in}}%
\pgfpathcurveto{\pgfqpoint{-0.011050in}{0.041667in}}{\pgfqpoint{-0.021649in}{0.037276in}}{\pgfqpoint{-0.029463in}{0.029463in}}%
\pgfpathcurveto{\pgfqpoint{-0.037276in}{0.021649in}}{\pgfqpoint{-0.041667in}{0.011050in}}{\pgfqpoint{-0.041667in}{0.000000in}}%
\pgfpathcurveto{\pgfqpoint{-0.041667in}{-0.011050in}}{\pgfqpoint{-0.037276in}{-0.021649in}}{\pgfqpoint{-0.029463in}{-0.029463in}}%
\pgfpathcurveto{\pgfqpoint{-0.021649in}{-0.037276in}}{\pgfqpoint{-0.011050in}{-0.041667in}}{\pgfqpoint{0.000000in}{-0.041667in}}%
\pgfpathlineto{\pgfqpoint{0.000000in}{-0.041667in}}%
\pgfpathclose%
\pgfusepath{stroke,fill}%
}%
\begin{pgfscope}%
\pgfsys@transformshift{1.544798in}{1.107437in}%
\pgfsys@useobject{currentmarker}{}%
\end{pgfscope}%
\begin{pgfscope}%
\pgfsys@transformshift{2.169078in}{1.167678in}%
\pgfsys@useobject{currentmarker}{}%
\end{pgfscope}%
\begin{pgfscope}%
\pgfsys@transformshift{3.209545in}{1.220652in}%
\pgfsys@useobject{currentmarker}{}%
\end{pgfscope}%
\begin{pgfscope}%
\pgfsys@transformshift{4.250013in}{1.272742in}%
\pgfsys@useobject{currentmarker}{}%
\end{pgfscope}%
\begin{pgfscope}%
\pgfsys@transformshift{5.290480in}{1.306845in}%
\pgfsys@useobject{currentmarker}{}%
\end{pgfscope}%
\end{pgfscope}%
\begin{pgfscope}%
\pgfpathrectangle{\pgfqpoint{1.128611in}{0.731763in}}{\pgfqpoint{4.578056in}{3.012428in}}%
\pgfusepath{clip}%
\pgfsetrectcap%
\pgfsetroundjoin%
\pgfsetlinewidth{1.505625pt}%
\definecolor{currentstroke}{rgb}{0.678431,0.709804,0.741176}%
\pgfsetstrokecolor{currentstroke}%
\pgfsetdash{}{0pt}%
\pgfpathmoveto{\pgfqpoint{1.544798in}{0.868692in}}%
\pgfpathlineto{\pgfqpoint{2.169078in}{0.952285in}}%
\pgfpathlineto{\pgfqpoint{3.209545in}{1.031926in}}%
\pgfpathlineto{\pgfqpoint{4.250013in}{1.147845in}}%
\pgfpathlineto{\pgfqpoint{5.290480in}{1.314533in}}%
\pgfusepath{stroke}%
\end{pgfscope}%
\begin{pgfscope}%
\pgfpathrectangle{\pgfqpoint{1.128611in}{0.731763in}}{\pgfqpoint{4.578056in}{3.012428in}}%
\pgfusepath{clip}%
\pgfsetbuttcap%
\pgfsetroundjoin%
\definecolor{currentfill}{rgb}{0.678431,0.709804,0.741176}%
\pgfsetfillcolor{currentfill}%
\pgfsetlinewidth{1.003750pt}%
\definecolor{currentstroke}{rgb}{0.678431,0.709804,0.741176}%
\pgfsetstrokecolor{currentstroke}%
\pgfsetdash{}{0pt}%
\pgfsys@defobject{currentmarker}{\pgfqpoint{-0.041667in}{-0.041667in}}{\pgfqpoint{0.041667in}{0.041667in}}{%
\pgfpathmoveto{\pgfqpoint{0.000000in}{-0.041667in}}%
\pgfpathcurveto{\pgfqpoint{0.011050in}{-0.041667in}}{\pgfqpoint{0.021649in}{-0.037276in}}{\pgfqpoint{0.029463in}{-0.029463in}}%
\pgfpathcurveto{\pgfqpoint{0.037276in}{-0.021649in}}{\pgfqpoint{0.041667in}{-0.011050in}}{\pgfqpoint{0.041667in}{0.000000in}}%
\pgfpathcurveto{\pgfqpoint{0.041667in}{0.011050in}}{\pgfqpoint{0.037276in}{0.021649in}}{\pgfqpoint{0.029463in}{0.029463in}}%
\pgfpathcurveto{\pgfqpoint{0.021649in}{0.037276in}}{\pgfqpoint{0.011050in}{0.041667in}}{\pgfqpoint{0.000000in}{0.041667in}}%
\pgfpathcurveto{\pgfqpoint{-0.011050in}{0.041667in}}{\pgfqpoint{-0.021649in}{0.037276in}}{\pgfqpoint{-0.029463in}{0.029463in}}%
\pgfpathcurveto{\pgfqpoint{-0.037276in}{0.021649in}}{\pgfqpoint{-0.041667in}{0.011050in}}{\pgfqpoint{-0.041667in}{0.000000in}}%
\pgfpathcurveto{\pgfqpoint{-0.041667in}{-0.011050in}}{\pgfqpoint{-0.037276in}{-0.021649in}}{\pgfqpoint{-0.029463in}{-0.029463in}}%
\pgfpathcurveto{\pgfqpoint{-0.021649in}{-0.037276in}}{\pgfqpoint{-0.011050in}{-0.041667in}}{\pgfqpoint{0.000000in}{-0.041667in}}%
\pgfpathlineto{\pgfqpoint{0.000000in}{-0.041667in}}%
\pgfpathclose%
\pgfusepath{stroke,fill}%
}%
\begin{pgfscope}%
\pgfsys@transformshift{1.544798in}{0.868692in}%
\pgfsys@useobject{currentmarker}{}%
\end{pgfscope}%
\begin{pgfscope}%
\pgfsys@transformshift{2.169078in}{0.952285in}%
\pgfsys@useobject{currentmarker}{}%
\end{pgfscope}%
\begin{pgfscope}%
\pgfsys@transformshift{3.209545in}{1.031926in}%
\pgfsys@useobject{currentmarker}{}%
\end{pgfscope}%
\begin{pgfscope}%
\pgfsys@transformshift{4.250013in}{1.147845in}%
\pgfsys@useobject{currentmarker}{}%
\end{pgfscope}%
\begin{pgfscope}%
\pgfsys@transformshift{5.290480in}{1.314533in}%
\pgfsys@useobject{currentmarker}{}%
\end{pgfscope}%
\end{pgfscope}%
\begin{pgfscope}%
\pgfsetrectcap%
\pgfsetmiterjoin%
\pgfsetlinewidth{0.803000pt}%
\definecolor{currentstroke}{rgb}{0.000000,0.000000,0.000000}%
\pgfsetstrokecolor{currentstroke}%
\pgfsetdash{}{0pt}%
\pgfpathmoveto{\pgfqpoint{1.128611in}{0.731763in}}%
\pgfpathlineto{\pgfqpoint{1.128611in}{3.744191in}}%
\pgfusepath{stroke}%
\end{pgfscope}%
\begin{pgfscope}%
\pgfsetrectcap%
\pgfsetmiterjoin%
\pgfsetlinewidth{0.803000pt}%
\definecolor{currentstroke}{rgb}{0.000000,0.000000,0.000000}%
\pgfsetstrokecolor{currentstroke}%
\pgfsetdash{}{0pt}%
\pgfpathmoveto{\pgfqpoint{5.706667in}{0.731763in}}%
\pgfpathlineto{\pgfqpoint{5.706667in}{3.744191in}}%
\pgfusepath{stroke}%
\end{pgfscope}%
\begin{pgfscope}%
\pgfsetrectcap%
\pgfsetmiterjoin%
\pgfsetlinewidth{0.803000pt}%
\definecolor{currentstroke}{rgb}{0.000000,0.000000,0.000000}%
\pgfsetstrokecolor{currentstroke}%
\pgfsetdash{}{0pt}%
\pgfpathmoveto{\pgfqpoint{1.128611in}{0.731763in}}%
\pgfpathlineto{\pgfqpoint{5.706667in}{0.731763in}}%
\pgfusepath{stroke}%
\end{pgfscope}%
\begin{pgfscope}%
\pgfsetrectcap%
\pgfsetmiterjoin%
\pgfsetlinewidth{0.803000pt}%
\definecolor{currentstroke}{rgb}{0.000000,0.000000,0.000000}%
\pgfsetstrokecolor{currentstroke}%
\pgfsetdash{}{0pt}%
\pgfpathmoveto{\pgfqpoint{1.128611in}{3.744191in}}%
\pgfpathlineto{\pgfqpoint{5.706667in}{3.744191in}}%
\pgfusepath{stroke}%
\end{pgfscope}%
\begin{pgfscope}%
\pgfsetbuttcap%
\pgfsetmiterjoin%
\definecolor{currentfill}{rgb}{1.000000,1.000000,1.000000}%
\pgfsetfillcolor{currentfill}%
\pgfsetfillopacity{0.800000}%
\pgfsetlinewidth{1.003750pt}%
\definecolor{currentstroke}{rgb}{0.800000,0.800000,0.800000}%
\pgfsetstrokecolor{currentstroke}%
\pgfsetstrokeopacity{0.800000}%
\pgfsetdash{}{0pt}%
\pgfpathmoveto{\pgfqpoint{1.860182in}{1.529490in}}%
\pgfpathlineto{\pgfqpoint{5.551111in}{1.529490in}}%
\pgfpathquadraticcurveto{\pgfqpoint{5.595556in}{1.529490in}}{\pgfqpoint{5.595556in}{1.573935in}}%
\pgfpathlineto{\pgfqpoint{5.595556in}{2.902020in}}%
\pgfpathquadraticcurveto{\pgfqpoint{5.595556in}{2.946465in}}{\pgfqpoint{5.551111in}{2.946465in}}%
\pgfpathlineto{\pgfqpoint{1.860182in}{2.946465in}}%
\pgfpathquadraticcurveto{\pgfqpoint{1.815738in}{2.946465in}}{\pgfqpoint{1.815738in}{2.902020in}}%
\pgfpathlineto{\pgfqpoint{1.815738in}{1.573935in}}%
\pgfpathquadraticcurveto{\pgfqpoint{1.815738in}{1.529490in}}{\pgfqpoint{1.860182in}{1.529490in}}%
\pgfpathlineto{\pgfqpoint{1.860182in}{1.529490in}}%
\pgfpathclose%
\pgfusepath{stroke,fill}%
\end{pgfscope}%
\begin{pgfscope}%
\pgfsetrectcap%
\pgfsetroundjoin%
\pgfsetlinewidth{1.505625pt}%
\definecolor{currentstroke}{rgb}{0.000000,0.000000,0.000000}%
\pgfsetstrokecolor{currentstroke}%
\pgfsetdash{}{0pt}%
\pgfpathmoveto{\pgfqpoint{1.904627in}{2.755570in}}%
\pgfpathlineto{\pgfqpoint{2.126849in}{2.755570in}}%
\pgfpathlineto{\pgfqpoint{2.349071in}{2.755570in}}%
\pgfusepath{stroke}%
\end{pgfscope}%
\begin{pgfscope}%
\pgfsetbuttcap%
\pgfsetroundjoin%
\definecolor{currentfill}{rgb}{0.000000,0.000000,0.000000}%
\pgfsetfillcolor{currentfill}%
\pgfsetlinewidth{1.003750pt}%
\definecolor{currentstroke}{rgb}{0.000000,0.000000,0.000000}%
\pgfsetstrokecolor{currentstroke}%
\pgfsetdash{}{0pt}%
\pgfsys@defobject{currentmarker}{\pgfqpoint{-0.041667in}{-0.041667in}}{\pgfqpoint{0.041667in}{0.041667in}}{%
\pgfpathmoveto{\pgfqpoint{0.000000in}{-0.041667in}}%
\pgfpathcurveto{\pgfqpoint{0.011050in}{-0.041667in}}{\pgfqpoint{0.021649in}{-0.037276in}}{\pgfqpoint{0.029463in}{-0.029463in}}%
\pgfpathcurveto{\pgfqpoint{0.037276in}{-0.021649in}}{\pgfqpoint{0.041667in}{-0.011050in}}{\pgfqpoint{0.041667in}{0.000000in}}%
\pgfpathcurveto{\pgfqpoint{0.041667in}{0.011050in}}{\pgfqpoint{0.037276in}{0.021649in}}{\pgfqpoint{0.029463in}{0.029463in}}%
\pgfpathcurveto{\pgfqpoint{0.021649in}{0.037276in}}{\pgfqpoint{0.011050in}{0.041667in}}{\pgfqpoint{0.000000in}{0.041667in}}%
\pgfpathcurveto{\pgfqpoint{-0.011050in}{0.041667in}}{\pgfqpoint{-0.021649in}{0.037276in}}{\pgfqpoint{-0.029463in}{0.029463in}}%
\pgfpathcurveto{\pgfqpoint{-0.037276in}{0.021649in}}{\pgfqpoint{-0.041667in}{0.011050in}}{\pgfqpoint{-0.041667in}{0.000000in}}%
\pgfpathcurveto{\pgfqpoint{-0.041667in}{-0.011050in}}{\pgfqpoint{-0.037276in}{-0.021649in}}{\pgfqpoint{-0.029463in}{-0.029463in}}%
\pgfpathcurveto{\pgfqpoint{-0.021649in}{-0.037276in}}{\pgfqpoint{-0.011050in}{-0.041667in}}{\pgfqpoint{0.000000in}{-0.041667in}}%
\pgfpathlineto{\pgfqpoint{0.000000in}{-0.041667in}}%
\pgfpathclose%
\pgfusepath{stroke,fill}%
}%
\begin{pgfscope}%
\pgfsys@transformshift{2.126849in}{2.755570in}%
\pgfsys@useobject{currentmarker}{}%
\end{pgfscope}%
\end{pgfscope}%
\begin{pgfscope}%
\definecolor{textcolor}{rgb}{0.000000,0.000000,0.000000}%
\pgfsetstrokecolor{textcolor}%
\pgfsetfillcolor{textcolor}%
\pgftext[x=2.526849in,y=2.677792in,left,base]{\color{textcolor}\rmfamily\fontsize{16.000000}{19.200000}\selectfont Passenger costs (\(\displaystyle \lambda\)-normalized)}%
\end{pgfscope}%
\begin{pgfscope}%
\pgfsetrectcap%
\pgfsetroundjoin%
\pgfsetlinewidth{1.505625pt}%
\definecolor{currentstroke}{rgb}{0.839216,0.152941,0.156863}%
\pgfsetstrokecolor{currentstroke}%
\pgfsetdash{}{0pt}%
\pgfpathmoveto{\pgfqpoint{1.904627in}{2.404876in}}%
\pgfpathlineto{\pgfqpoint{2.126849in}{2.404876in}}%
\pgfpathlineto{\pgfqpoint{2.349071in}{2.404876in}}%
\pgfusepath{stroke}%
\end{pgfscope}%
\begin{pgfscope}%
\pgfsetbuttcap%
\pgfsetroundjoin%
\definecolor{currentfill}{rgb}{0.839216,0.152941,0.156863}%
\pgfsetfillcolor{currentfill}%
\pgfsetlinewidth{1.003750pt}%
\definecolor{currentstroke}{rgb}{0.839216,0.152941,0.156863}%
\pgfsetstrokecolor{currentstroke}%
\pgfsetdash{}{0pt}%
\pgfsys@defobject{currentmarker}{\pgfqpoint{-0.041667in}{-0.041667in}}{\pgfqpoint{0.041667in}{0.041667in}}{%
\pgfpathmoveto{\pgfqpoint{0.000000in}{-0.041667in}}%
\pgfpathcurveto{\pgfqpoint{0.011050in}{-0.041667in}}{\pgfqpoint{0.021649in}{-0.037276in}}{\pgfqpoint{0.029463in}{-0.029463in}}%
\pgfpathcurveto{\pgfqpoint{0.037276in}{-0.021649in}}{\pgfqpoint{0.041667in}{-0.011050in}}{\pgfqpoint{0.041667in}{0.000000in}}%
\pgfpathcurveto{\pgfqpoint{0.041667in}{0.011050in}}{\pgfqpoint{0.037276in}{0.021649in}}{\pgfqpoint{0.029463in}{0.029463in}}%
\pgfpathcurveto{\pgfqpoint{0.021649in}{0.037276in}}{\pgfqpoint{0.011050in}{0.041667in}}{\pgfqpoint{0.000000in}{0.041667in}}%
\pgfpathcurveto{\pgfqpoint{-0.011050in}{0.041667in}}{\pgfqpoint{-0.021649in}{0.037276in}}{\pgfqpoint{-0.029463in}{0.029463in}}%
\pgfpathcurveto{\pgfqpoint{-0.037276in}{0.021649in}}{\pgfqpoint{-0.041667in}{0.011050in}}{\pgfqpoint{-0.041667in}{0.000000in}}%
\pgfpathcurveto{\pgfqpoint{-0.041667in}{-0.011050in}}{\pgfqpoint{-0.037276in}{-0.021649in}}{\pgfqpoint{-0.029463in}{-0.029463in}}%
\pgfpathcurveto{\pgfqpoint{-0.021649in}{-0.037276in}}{\pgfqpoint{-0.011050in}{-0.041667in}}{\pgfqpoint{0.000000in}{-0.041667in}}%
\pgfpathlineto{\pgfqpoint{0.000000in}{-0.041667in}}%
\pgfpathclose%
\pgfusepath{stroke,fill}%
}%
\begin{pgfscope}%
\pgfsys@transformshift{2.126849in}{2.404876in}%
\pgfsys@useobject{currentmarker}{}%
\end{pgfscope}%
\end{pgfscope}%
\begin{pgfscope}%
\definecolor{textcolor}{rgb}{0.000000,0.000000,0.000000}%
\pgfsetstrokecolor{textcolor}%
\pgfsetfillcolor{textcolor}%
\pgftext[x=2.526849in,y=2.327098in,left,base]{\color{textcolor}\rmfamily\fontsize{16.000000}{19.200000}\selectfont Objective (\(\displaystyle \lambda\)-normalized)}%
\end{pgfscope}%
\begin{pgfscope}%
\pgfsetrectcap%
\pgfsetroundjoin%
\pgfsetlinewidth{1.505625pt}%
\definecolor{currentstroke}{rgb}{0.301961,0.301961,0.301961}%
\pgfsetstrokecolor{currentstroke}%
\pgfsetdash{}{0pt}%
\pgfpathmoveto{\pgfqpoint{1.904627in}{2.067299in}}%
\pgfpathlineto{\pgfqpoint{2.126849in}{2.067299in}}%
\pgfpathlineto{\pgfqpoint{2.349071in}{2.067299in}}%
\pgfusepath{stroke}%
\end{pgfscope}%
\begin{pgfscope}%
\pgfsetbuttcap%
\pgfsetroundjoin%
\definecolor{currentfill}{rgb}{0.301961,0.301961,0.301961}%
\pgfsetfillcolor{currentfill}%
\pgfsetlinewidth{1.003750pt}%
\definecolor{currentstroke}{rgb}{0.301961,0.301961,0.301961}%
\pgfsetstrokecolor{currentstroke}%
\pgfsetdash{}{0pt}%
\pgfsys@defobject{currentmarker}{\pgfqpoint{-0.041667in}{-0.041667in}}{\pgfqpoint{0.041667in}{0.041667in}}{%
\pgfpathmoveto{\pgfqpoint{0.000000in}{-0.041667in}}%
\pgfpathcurveto{\pgfqpoint{0.011050in}{-0.041667in}}{\pgfqpoint{0.021649in}{-0.037276in}}{\pgfqpoint{0.029463in}{-0.029463in}}%
\pgfpathcurveto{\pgfqpoint{0.037276in}{-0.021649in}}{\pgfqpoint{0.041667in}{-0.011050in}}{\pgfqpoint{0.041667in}{0.000000in}}%
\pgfpathcurveto{\pgfqpoint{0.041667in}{0.011050in}}{\pgfqpoint{0.037276in}{0.021649in}}{\pgfqpoint{0.029463in}{0.029463in}}%
\pgfpathcurveto{\pgfqpoint{0.021649in}{0.037276in}}{\pgfqpoint{0.011050in}{0.041667in}}{\pgfqpoint{0.000000in}{0.041667in}}%
\pgfpathcurveto{\pgfqpoint{-0.011050in}{0.041667in}}{\pgfqpoint{-0.021649in}{0.037276in}}{\pgfqpoint{-0.029463in}{0.029463in}}%
\pgfpathcurveto{\pgfqpoint{-0.037276in}{0.021649in}}{\pgfqpoint{-0.041667in}{0.011050in}}{\pgfqpoint{-0.041667in}{0.000000in}}%
\pgfpathcurveto{\pgfqpoint{-0.041667in}{-0.011050in}}{\pgfqpoint{-0.037276in}{-0.021649in}}{\pgfqpoint{-0.029463in}{-0.029463in}}%
\pgfpathcurveto{\pgfqpoint{-0.021649in}{-0.037276in}}{\pgfqpoint{-0.011050in}{-0.041667in}}{\pgfqpoint{0.000000in}{-0.041667in}}%
\pgfpathlineto{\pgfqpoint{0.000000in}{-0.041667in}}%
\pgfpathclose%
\pgfusepath{stroke,fill}%
}%
\begin{pgfscope}%
\pgfsys@transformshift{2.126849in}{2.067299in}%
\pgfsys@useobject{currentmarker}{}%
\end{pgfscope}%
\end{pgfscope}%
\begin{pgfscope}%
\definecolor{textcolor}{rgb}{0.000000,0.000000,0.000000}%
\pgfsetstrokecolor{textcolor}%
\pgfsetfillcolor{textcolor}%
\pgftext[x=2.526849in,y=1.989521in,left,base]{\color{textcolor}\rmfamily\fontsize{16.000000}{19.200000}\selectfont Ticket revenue}%
\end{pgfscope}%
\begin{pgfscope}%
\pgfsetrectcap%
\pgfsetroundjoin%
\pgfsetlinewidth{1.505625pt}%
\definecolor{currentstroke}{rgb}{0.678431,0.709804,0.741176}%
\pgfsetstrokecolor{currentstroke}%
\pgfsetdash{}{0pt}%
\pgfpathmoveto{\pgfqpoint{1.904627in}{1.742839in}}%
\pgfpathlineto{\pgfqpoint{2.126849in}{1.742839in}}%
\pgfpathlineto{\pgfqpoint{2.349071in}{1.742839in}}%
\pgfusepath{stroke}%
\end{pgfscope}%
\begin{pgfscope}%
\pgfsetbuttcap%
\pgfsetroundjoin%
\definecolor{currentfill}{rgb}{0.678431,0.709804,0.741176}%
\pgfsetfillcolor{currentfill}%
\pgfsetlinewidth{1.003750pt}%
\definecolor{currentstroke}{rgb}{0.678431,0.709804,0.741176}%
\pgfsetstrokecolor{currentstroke}%
\pgfsetdash{}{0pt}%
\pgfsys@defobject{currentmarker}{\pgfqpoint{-0.041667in}{-0.041667in}}{\pgfqpoint{0.041667in}{0.041667in}}{%
\pgfpathmoveto{\pgfqpoint{0.000000in}{-0.041667in}}%
\pgfpathcurveto{\pgfqpoint{0.011050in}{-0.041667in}}{\pgfqpoint{0.021649in}{-0.037276in}}{\pgfqpoint{0.029463in}{-0.029463in}}%
\pgfpathcurveto{\pgfqpoint{0.037276in}{-0.021649in}}{\pgfqpoint{0.041667in}{-0.011050in}}{\pgfqpoint{0.041667in}{0.000000in}}%
\pgfpathcurveto{\pgfqpoint{0.041667in}{0.011050in}}{\pgfqpoint{0.037276in}{0.021649in}}{\pgfqpoint{0.029463in}{0.029463in}}%
\pgfpathcurveto{\pgfqpoint{0.021649in}{0.037276in}}{\pgfqpoint{0.011050in}{0.041667in}}{\pgfqpoint{0.000000in}{0.041667in}}%
\pgfpathcurveto{\pgfqpoint{-0.011050in}{0.041667in}}{\pgfqpoint{-0.021649in}{0.037276in}}{\pgfqpoint{-0.029463in}{0.029463in}}%
\pgfpathcurveto{\pgfqpoint{-0.037276in}{0.021649in}}{\pgfqpoint{-0.041667in}{0.011050in}}{\pgfqpoint{-0.041667in}{0.000000in}}%
\pgfpathcurveto{\pgfqpoint{-0.041667in}{-0.011050in}}{\pgfqpoint{-0.037276in}{-0.021649in}}{\pgfqpoint{-0.029463in}{-0.029463in}}%
\pgfpathcurveto{\pgfqpoint{-0.021649in}{-0.037276in}}{\pgfqpoint{-0.011050in}{-0.041667in}}{\pgfqpoint{0.000000in}{-0.041667in}}%
\pgfpathlineto{\pgfqpoint{0.000000in}{-0.041667in}}%
\pgfpathclose%
\pgfusepath{stroke,fill}%
}%
\begin{pgfscope}%
\pgfsys@transformshift{2.126849in}{1.742839in}%
\pgfsys@useobject{currentmarker}{}%
\end{pgfscope}%
\end{pgfscope}%
\begin{pgfscope}%
\definecolor{textcolor}{rgb}{0.000000,0.000000,0.000000}%
\pgfsetstrokecolor{textcolor}%
\pgfsetfillcolor{textcolor}%
\pgftext[x=2.526849in,y=1.665061in,left,base]{\color{textcolor}\rmfamily\fontsize{16.000000}{19.200000}\selectfont Operator costs}%
\end{pgfscope}%
\begin{pgfscope}%
\pgfsetbuttcap%
\pgfsetroundjoin%
\definecolor{currentfill}{rgb}{0.000000,0.000000,0.000000}%
\pgfsetfillcolor{currentfill}%
\pgfsetlinewidth{0.803000pt}%
\definecolor{currentstroke}{rgb}{0.000000,0.000000,0.000000}%
\pgfsetstrokecolor{currentstroke}%
\pgfsetdash{}{0pt}%
\pgfsys@defobject{currentmarker}{\pgfqpoint{0.000000in}{0.000000in}}{\pgfqpoint{0.083333in}{0.000000in}}{%
\pgfpathmoveto{\pgfqpoint{0.000000in}{0.000000in}}%
\pgfpathlineto{\pgfqpoint{0.083333in}{0.000000in}}%
\pgfusepath{stroke,fill}%
}%
\begin{pgfscope}%
\pgfsys@transformshift{5.706667in}{0.984625in}%
\pgfsys@useobject{currentmarker}{}%
\end{pgfscope}%
\end{pgfscope}%
\begin{pgfscope}%
\pgfsetbuttcap%
\pgfsetroundjoin%
\definecolor{currentfill}{rgb}{0.000000,0.000000,0.000000}%
\pgfsetfillcolor{currentfill}%
\pgfsetlinewidth{0.803000pt}%
\definecolor{currentstroke}{rgb}{0.000000,0.000000,0.000000}%
\pgfsetstrokecolor{currentstroke}%
\pgfsetdash{}{0pt}%
\pgfsys@defobject{currentmarker}{\pgfqpoint{0.000000in}{0.000000in}}{\pgfqpoint{0.083333in}{0.000000in}}{%
\pgfpathmoveto{\pgfqpoint{0.000000in}{0.000000in}}%
\pgfpathlineto{\pgfqpoint{0.083333in}{0.000000in}}%
\pgfusepath{stroke,fill}%
}%
\begin{pgfscope}%
\pgfsys@transformshift{5.706667in}{1.532342in}%
\pgfsys@useobject{currentmarker}{}%
\end{pgfscope}%
\end{pgfscope}%
\begin{pgfscope}%
\pgfsetbuttcap%
\pgfsetroundjoin%
\definecolor{currentfill}{rgb}{0.000000,0.000000,0.000000}%
\pgfsetfillcolor{currentfill}%
\pgfsetlinewidth{0.803000pt}%
\definecolor{currentstroke}{rgb}{0.000000,0.000000,0.000000}%
\pgfsetstrokecolor{currentstroke}%
\pgfsetdash{}{0pt}%
\pgfsys@defobject{currentmarker}{\pgfqpoint{0.000000in}{0.000000in}}{\pgfqpoint{0.083333in}{0.000000in}}{%
\pgfpathmoveto{\pgfqpoint{0.000000in}{0.000000in}}%
\pgfpathlineto{\pgfqpoint{0.083333in}{0.000000in}}%
\pgfusepath{stroke,fill}%
}%
\begin{pgfscope}%
\pgfsys@transformshift{5.706667in}{2.080058in}%
\pgfsys@useobject{currentmarker}{}%
\end{pgfscope}%
\end{pgfscope}%
\begin{pgfscope}%
\pgfsetbuttcap%
\pgfsetroundjoin%
\definecolor{currentfill}{rgb}{0.000000,0.000000,0.000000}%
\pgfsetfillcolor{currentfill}%
\pgfsetlinewidth{0.803000pt}%
\definecolor{currentstroke}{rgb}{0.000000,0.000000,0.000000}%
\pgfsetstrokecolor{currentstroke}%
\pgfsetdash{}{0pt}%
\pgfsys@defobject{currentmarker}{\pgfqpoint{0.000000in}{0.000000in}}{\pgfqpoint{0.083333in}{0.000000in}}{%
\pgfpathmoveto{\pgfqpoint{0.000000in}{0.000000in}}%
\pgfpathlineto{\pgfqpoint{0.083333in}{0.000000in}}%
\pgfusepath{stroke,fill}%
}%
\begin{pgfscope}%
\pgfsys@transformshift{5.706667in}{2.627775in}%
\pgfsys@useobject{currentmarker}{}%
\end{pgfscope}%
\end{pgfscope}%
\begin{pgfscope}%
\pgfsetbuttcap%
\pgfsetroundjoin%
\definecolor{currentfill}{rgb}{0.000000,0.000000,0.000000}%
\pgfsetfillcolor{currentfill}%
\pgfsetlinewidth{0.803000pt}%
\definecolor{currentstroke}{rgb}{0.000000,0.000000,0.000000}%
\pgfsetstrokecolor{currentstroke}%
\pgfsetdash{}{0pt}%
\pgfsys@defobject{currentmarker}{\pgfqpoint{0.000000in}{0.000000in}}{\pgfqpoint{0.083333in}{0.000000in}}{%
\pgfpathmoveto{\pgfqpoint{0.000000in}{0.000000in}}%
\pgfpathlineto{\pgfqpoint{0.083333in}{0.000000in}}%
\pgfusepath{stroke,fill}%
}%
\begin{pgfscope}%
\pgfsys@transformshift{5.706667in}{3.175491in}%
\pgfsys@useobject{currentmarker}{}%
\end{pgfscope}%
\end{pgfscope}%
\begin{pgfscope}%
\pgfsetbuttcap%
\pgfsetroundjoin%
\definecolor{currentfill}{rgb}{0.000000,0.000000,0.000000}%
\pgfsetfillcolor{currentfill}%
\pgfsetlinewidth{0.803000pt}%
\definecolor{currentstroke}{rgb}{0.000000,0.000000,0.000000}%
\pgfsetstrokecolor{currentstroke}%
\pgfsetdash{}{0pt}%
\pgfsys@defobject{currentmarker}{\pgfqpoint{0.000000in}{0.000000in}}{\pgfqpoint{0.083333in}{0.000000in}}{%
\pgfpathmoveto{\pgfqpoint{0.000000in}{0.000000in}}%
\pgfpathlineto{\pgfqpoint{0.083333in}{0.000000in}}%
\pgfusepath{stroke,fill}%
}%
\begin{pgfscope}%
\pgfsys@transformshift{5.706667in}{3.723208in}%
\pgfsys@useobject{currentmarker}{}%
\end{pgfscope}%
\end{pgfscope}%
\begin{pgfscope}%
\pgfsetrectcap%
\pgfsetmiterjoin%
\pgfsetlinewidth{0.803000pt}%
\definecolor{currentstroke}{rgb}{0.000000,0.000000,0.000000}%
\pgfsetstrokecolor{currentstroke}%
\pgfsetdash{}{0pt}%
\pgfpathmoveto{\pgfqpoint{1.128611in}{0.731763in}}%
\pgfpathlineto{\pgfqpoint{1.128611in}{3.744191in}}%
\pgfusepath{stroke}%
\end{pgfscope}%
\begin{pgfscope}%
\pgfsetrectcap%
\pgfsetmiterjoin%
\pgfsetlinewidth{0.803000pt}%
\definecolor{currentstroke}{rgb}{0.000000,0.000000,0.000000}%
\pgfsetstrokecolor{currentstroke}%
\pgfsetdash{}{0pt}%
\pgfpathmoveto{\pgfqpoint{5.706667in}{0.731763in}}%
\pgfpathlineto{\pgfqpoint{5.706667in}{3.744191in}}%
\pgfusepath{stroke}%
\end{pgfscope}%
\begin{pgfscope}%
\pgfsetrectcap%
\pgfsetmiterjoin%
\pgfsetlinewidth{0.803000pt}%
\definecolor{currentstroke}{rgb}{0.000000,0.000000,0.000000}%
\pgfsetstrokecolor{currentstroke}%
\pgfsetdash{}{0pt}%
\pgfpathmoveto{\pgfqpoint{1.128611in}{0.731763in}}%
\pgfpathlineto{\pgfqpoint{5.706667in}{0.731763in}}%
\pgfusepath{stroke}%
\end{pgfscope}%
\begin{pgfscope}%
\pgfsetrectcap%
\pgfsetmiterjoin%
\pgfsetlinewidth{0.803000pt}%
\definecolor{currentstroke}{rgb}{0.000000,0.000000,0.000000}%
\pgfsetstrokecolor{currentstroke}%
\pgfsetdash{}{0pt}%
\pgfpathmoveto{\pgfqpoint{1.128611in}{3.744191in}}%
\pgfpathlineto{\pgfqpoint{5.706667in}{3.744191in}}%
\pgfusepath{stroke}%
\end{pgfscope}%
\end{pgfpicture}%
\makeatother%
\endgroup%

%% file: figures/Comp/Inst_ObjVSLambda.pgf
\begingroup%
\makeatletter%
\begin{pgfpicture}%
\pgfpathrectangle{\pgfpointorigin}{\pgfqpoint{6.000000in}{4.000000in}}%
\pgfusepath{use as bounding box, clip}%
\begin{pgfscope}%
\pgfsetbuttcap%
\pgfsetmiterjoin%
\definecolor{currentfill}{rgb}{1.000000,1.000000,1.000000}%
\pgfsetfillcolor{currentfill}%
\pgfsetlinewidth{0.000000pt}%
\definecolor{currentstroke}{rgb}{1.000000,1.000000,1.000000}%
\pgfsetstrokecolor{currentstroke}%
\pgfsetdash{}{0pt}%
\pgfpathmoveto{\pgfqpoint{0.000000in}{0.000000in}}%
\pgfpathlineto{\pgfqpoint{6.000000in}{0.000000in}}%
\pgfpathlineto{\pgfqpoint{6.000000in}{4.000000in}}%
\pgfpathlineto{\pgfqpoint{0.000000in}{4.000000in}}%
\pgfpathlineto{\pgfqpoint{0.000000in}{0.000000in}}%
\pgfpathclose%
\pgfusepath{fill}%
\end{pgfscope}%
\begin{pgfscope}%
\pgfsetbuttcap%
\pgfsetmiterjoin%
\definecolor{currentfill}{rgb}{1.000000,1.000000,1.000000}%
\pgfsetfillcolor{currentfill}%
\pgfsetlinewidth{0.000000pt}%
\definecolor{currentstroke}{rgb}{0.000000,0.000000,0.000000}%
\pgfsetstrokecolor{currentstroke}%
\pgfsetstrokeopacity{0.000000}%
\pgfsetdash{}{0pt}%
\pgfpathmoveto{\pgfqpoint{1.104012in}{0.707164in}}%
\pgfpathlineto{\pgfqpoint{5.790000in}{0.707164in}}%
\pgfpathlineto{\pgfqpoint{5.790000in}{3.790000in}}%
\pgfpathlineto{\pgfqpoint{1.104012in}{3.790000in}}%
\pgfpathlineto{\pgfqpoint{1.104012in}{0.707164in}}%
\pgfpathclose%
\pgfusepath{fill}%
\end{pgfscope}%
\begin{pgfscope}%
\pgfsetbuttcap%
\pgfsetroundjoin%
\definecolor{currentfill}{rgb}{0.000000,0.000000,0.000000}%
\pgfsetfillcolor{currentfill}%
\pgfsetlinewidth{0.803000pt}%
\definecolor{currentstroke}{rgb}{0.000000,0.000000,0.000000}%
\pgfsetstrokecolor{currentstroke}%
\pgfsetdash{}{0pt}%
\pgfsys@defobject{currentmarker}{\pgfqpoint{0.000000in}{-0.048611in}}{\pgfqpoint{0.000000in}{0.000000in}}{%
\pgfpathmoveto{\pgfqpoint{0.000000in}{0.000000in}}%
\pgfpathlineto{\pgfqpoint{0.000000in}{-0.048611in}}%
\pgfusepath{stroke,fill}%
}%
\begin{pgfscope}%
\pgfsys@transformshift{1.104012in}{0.707164in}%
\pgfsys@useobject{currentmarker}{}%
\end{pgfscope}%
\end{pgfscope}%
\begin{pgfscope}%
\definecolor{textcolor}{rgb}{0.000000,0.000000,0.000000}%
\pgfsetstrokecolor{textcolor}%
\pgfsetfillcolor{textcolor}%
\pgftext[x=1.104012in,y=0.609942in,,top]{\color{textcolor}\rmfamily\fontsize{14.000000}{16.800000}\selectfont \(\displaystyle {0.0}\)}%
\end{pgfscope}%
\begin{pgfscope}%
\pgfsetbuttcap%
\pgfsetroundjoin%
\definecolor{currentfill}{rgb}{0.000000,0.000000,0.000000}%
\pgfsetfillcolor{currentfill}%
\pgfsetlinewidth{0.803000pt}%
\definecolor{currentstroke}{rgb}{0.000000,0.000000,0.000000}%
\pgfsetstrokecolor{currentstroke}%
\pgfsetdash{}{0pt}%
\pgfsys@defobject{currentmarker}{\pgfqpoint{0.000000in}{-0.048611in}}{\pgfqpoint{0.000000in}{0.000000in}}{%
\pgfpathmoveto{\pgfqpoint{0.000000in}{0.000000in}}%
\pgfpathlineto{\pgfqpoint{0.000000in}{-0.048611in}}%
\pgfusepath{stroke,fill}%
}%
\begin{pgfscope}%
\pgfsys@transformshift{2.219723in}{0.707164in}%
\pgfsys@useobject{currentmarker}{}%
\end{pgfscope}%
\end{pgfscope}%
\begin{pgfscope}%
\definecolor{textcolor}{rgb}{0.000000,0.000000,0.000000}%
\pgfsetstrokecolor{textcolor}%
\pgfsetfillcolor{textcolor}%
\pgftext[x=2.219723in,y=0.609942in,,top]{\color{textcolor}\rmfamily\fontsize{14.000000}{16.800000}\selectfont \(\displaystyle {0.5}\)}%
\end{pgfscope}%
\begin{pgfscope}%
\pgfsetbuttcap%
\pgfsetroundjoin%
\definecolor{currentfill}{rgb}{0.000000,0.000000,0.000000}%
\pgfsetfillcolor{currentfill}%
\pgfsetlinewidth{0.803000pt}%
\definecolor{currentstroke}{rgb}{0.000000,0.000000,0.000000}%
\pgfsetstrokecolor{currentstroke}%
\pgfsetdash{}{0pt}%
\pgfsys@defobject{currentmarker}{\pgfqpoint{0.000000in}{-0.048611in}}{\pgfqpoint{0.000000in}{0.000000in}}{%
\pgfpathmoveto{\pgfqpoint{0.000000in}{0.000000in}}%
\pgfpathlineto{\pgfqpoint{0.000000in}{-0.048611in}}%
\pgfusepath{stroke,fill}%
}%
\begin{pgfscope}%
\pgfsys@transformshift{3.335435in}{0.707164in}%
\pgfsys@useobject{currentmarker}{}%
\end{pgfscope}%
\end{pgfscope}%
\begin{pgfscope}%
\definecolor{textcolor}{rgb}{0.000000,0.000000,0.000000}%
\pgfsetstrokecolor{textcolor}%
\pgfsetfillcolor{textcolor}%
\pgftext[x=3.335435in,y=0.609942in,,top]{\color{textcolor}\rmfamily\fontsize{14.000000}{16.800000}\selectfont \(\displaystyle {1.0}\)}%
\end{pgfscope}%
\begin{pgfscope}%
\pgfsetbuttcap%
\pgfsetroundjoin%
\definecolor{currentfill}{rgb}{0.000000,0.000000,0.000000}%
\pgfsetfillcolor{currentfill}%
\pgfsetlinewidth{0.803000pt}%
\definecolor{currentstroke}{rgb}{0.000000,0.000000,0.000000}%
\pgfsetstrokecolor{currentstroke}%
\pgfsetdash{}{0pt}%
\pgfsys@defobject{currentmarker}{\pgfqpoint{0.000000in}{-0.048611in}}{\pgfqpoint{0.000000in}{0.000000in}}{%
\pgfpathmoveto{\pgfqpoint{0.000000in}{0.000000in}}%
\pgfpathlineto{\pgfqpoint{0.000000in}{-0.048611in}}%
\pgfusepath{stroke,fill}%
}%
\begin{pgfscope}%
\pgfsys@transformshift{4.451146in}{0.707164in}%
\pgfsys@useobject{currentmarker}{}%
\end{pgfscope}%
\end{pgfscope}%
\begin{pgfscope}%
\definecolor{textcolor}{rgb}{0.000000,0.000000,0.000000}%
\pgfsetstrokecolor{textcolor}%
\pgfsetfillcolor{textcolor}%
\pgftext[x=4.451146in,y=0.609942in,,top]{\color{textcolor}\rmfamily\fontsize{14.000000}{16.800000}\selectfont \(\displaystyle {1.5}\)}%
\end{pgfscope}%
\begin{pgfscope}%
\pgfsetbuttcap%
\pgfsetroundjoin%
\definecolor{currentfill}{rgb}{0.000000,0.000000,0.000000}%
\pgfsetfillcolor{currentfill}%
\pgfsetlinewidth{0.803000pt}%
\definecolor{currentstroke}{rgb}{0.000000,0.000000,0.000000}%
\pgfsetstrokecolor{currentstroke}%
\pgfsetdash{}{0pt}%
\pgfsys@defobject{currentmarker}{\pgfqpoint{0.000000in}{-0.048611in}}{\pgfqpoint{0.000000in}{0.000000in}}{%
\pgfpathmoveto{\pgfqpoint{0.000000in}{0.000000in}}%
\pgfpathlineto{\pgfqpoint{0.000000in}{-0.048611in}}%
\pgfusepath{stroke,fill}%
}%
\begin{pgfscope}%
\pgfsys@transformshift{5.566858in}{0.707164in}%
\pgfsys@useobject{currentmarker}{}%
\end{pgfscope}%
\end{pgfscope}%
\begin{pgfscope}%
\definecolor{textcolor}{rgb}{0.000000,0.000000,0.000000}%
\pgfsetstrokecolor{textcolor}%
\pgfsetfillcolor{textcolor}%
\pgftext[x=5.566858in,y=0.609942in,,top]{\color{textcolor}\rmfamily\fontsize{14.000000}{16.800000}\selectfont \(\displaystyle {2.0}\)}%
\end{pgfscope}%
\begin{pgfscope}%
\definecolor{textcolor}{rgb}{0.000000,0.000000,0.000000}%
\pgfsetstrokecolor{textcolor}%
\pgfsetfillcolor{textcolor}%
\pgftext[x=3.447006in,y=0.376609in,,top]{\color{textcolor}\rmfamily\fontsize{14.000000}{16.800000}\selectfont \(\displaystyle \lambda\)}%
\end{pgfscope}%
\begin{pgfscope}%
\pgfsetbuttcap%
\pgfsetroundjoin%
\definecolor{currentfill}{rgb}{0.000000,0.000000,0.000000}%
\pgfsetfillcolor{currentfill}%
\pgfsetlinewidth{0.803000pt}%
\definecolor{currentstroke}{rgb}{0.000000,0.000000,0.000000}%
\pgfsetstrokecolor{currentstroke}%
\pgfsetdash{}{0pt}%
\pgfsys@defobject{currentmarker}{\pgfqpoint{-0.048611in}{0.000000in}}{\pgfqpoint{-0.000000in}{0.000000in}}{%
\pgfpathmoveto{\pgfqpoint{-0.000000in}{0.000000in}}%
\pgfpathlineto{\pgfqpoint{-0.048611in}{0.000000in}}%
\pgfusepath{stroke,fill}%
}%
\begin{pgfscope}%
\pgfsys@transformshift{1.104012in}{1.142055in}%
\pgfsys@useobject{currentmarker}{}%
\end{pgfscope}%
\end{pgfscope}%
\begin{pgfscope}%
\definecolor{textcolor}{rgb}{0.000000,0.000000,0.000000}%
\pgfsetstrokecolor{textcolor}%
\pgfsetfillcolor{textcolor}%
\pgftext[x=0.419297in, y=1.072611in, left, base]{\color{textcolor}\rmfamily\fontsize{14.000000}{16.800000}\selectfont \(\displaystyle {190000}\)}%
\end{pgfscope}%
\begin{pgfscope}%
\pgfsetbuttcap%
\pgfsetroundjoin%
\definecolor{currentfill}{rgb}{0.000000,0.000000,0.000000}%
\pgfsetfillcolor{currentfill}%
\pgfsetlinewidth{0.803000pt}%
\definecolor{currentstroke}{rgb}{0.000000,0.000000,0.000000}%
\pgfsetstrokecolor{currentstroke}%
\pgfsetdash{}{0pt}%
\pgfsys@defobject{currentmarker}{\pgfqpoint{-0.048611in}{0.000000in}}{\pgfqpoint{-0.000000in}{0.000000in}}{%
\pgfpathmoveto{\pgfqpoint{-0.000000in}{0.000000in}}%
\pgfpathlineto{\pgfqpoint{-0.048611in}{0.000000in}}%
\pgfusepath{stroke,fill}%
}%
\begin{pgfscope}%
\pgfsys@transformshift{1.104012in}{1.740922in}%
\pgfsys@useobject{currentmarker}{}%
\end{pgfscope}%
\end{pgfscope}%
\begin{pgfscope}%
\definecolor{textcolor}{rgb}{0.000000,0.000000,0.000000}%
\pgfsetstrokecolor{textcolor}%
\pgfsetfillcolor{textcolor}%
\pgftext[x=0.419297in, y=1.671478in, left, base]{\color{textcolor}\rmfamily\fontsize{14.000000}{16.800000}\selectfont \(\displaystyle {195000}\)}%
\end{pgfscope}%
\begin{pgfscope}%
\pgfsetbuttcap%
\pgfsetroundjoin%
\definecolor{currentfill}{rgb}{0.000000,0.000000,0.000000}%
\pgfsetfillcolor{currentfill}%
\pgfsetlinewidth{0.803000pt}%
\definecolor{currentstroke}{rgb}{0.000000,0.000000,0.000000}%
\pgfsetstrokecolor{currentstroke}%
\pgfsetdash{}{0pt}%
\pgfsys@defobject{currentmarker}{\pgfqpoint{-0.048611in}{0.000000in}}{\pgfqpoint{-0.000000in}{0.000000in}}{%
\pgfpathmoveto{\pgfqpoint{-0.000000in}{0.000000in}}%
\pgfpathlineto{\pgfqpoint{-0.048611in}{0.000000in}}%
\pgfusepath{stroke,fill}%
}%
\begin{pgfscope}%
\pgfsys@transformshift{1.104012in}{2.339789in}%
\pgfsys@useobject{currentmarker}{}%
\end{pgfscope}%
\end{pgfscope}%
\begin{pgfscope}%
\definecolor{textcolor}{rgb}{0.000000,0.000000,0.000000}%
\pgfsetstrokecolor{textcolor}%
\pgfsetfillcolor{textcolor}%
\pgftext[x=0.419297in, y=2.270345in, left, base]{\color{textcolor}\rmfamily\fontsize{14.000000}{16.800000}\selectfont \(\displaystyle {200000}\)}%
\end{pgfscope}%
\begin{pgfscope}%
\pgfsetbuttcap%
\pgfsetroundjoin%
\definecolor{currentfill}{rgb}{0.000000,0.000000,0.000000}%
\pgfsetfillcolor{currentfill}%
\pgfsetlinewidth{0.803000pt}%
\definecolor{currentstroke}{rgb}{0.000000,0.000000,0.000000}%
\pgfsetstrokecolor{currentstroke}%
\pgfsetdash{}{0pt}%
\pgfsys@defobject{currentmarker}{\pgfqpoint{-0.048611in}{0.000000in}}{\pgfqpoint{-0.000000in}{0.000000in}}{%
\pgfpathmoveto{\pgfqpoint{-0.000000in}{0.000000in}}%
\pgfpathlineto{\pgfqpoint{-0.048611in}{0.000000in}}%
\pgfusepath{stroke,fill}%
}%
\begin{pgfscope}%
\pgfsys@transformshift{1.104012in}{2.938657in}%
\pgfsys@useobject{currentmarker}{}%
\end{pgfscope}%
\end{pgfscope}%
\begin{pgfscope}%
\definecolor{textcolor}{rgb}{0.000000,0.000000,0.000000}%
\pgfsetstrokecolor{textcolor}%
\pgfsetfillcolor{textcolor}%
\pgftext[x=0.419297in, y=2.869212in, left, base]{\color{textcolor}\rmfamily\fontsize{14.000000}{16.800000}\selectfont \(\displaystyle {205000}\)}%
\end{pgfscope}%
\begin{pgfscope}%
\pgfsetbuttcap%
\pgfsetroundjoin%
\definecolor{currentfill}{rgb}{0.000000,0.000000,0.000000}%
\pgfsetfillcolor{currentfill}%
\pgfsetlinewidth{0.803000pt}%
\definecolor{currentstroke}{rgb}{0.000000,0.000000,0.000000}%
\pgfsetstrokecolor{currentstroke}%
\pgfsetdash{}{0pt}%
\pgfsys@defobject{currentmarker}{\pgfqpoint{-0.048611in}{0.000000in}}{\pgfqpoint{-0.000000in}{0.000000in}}{%
\pgfpathmoveto{\pgfqpoint{-0.000000in}{0.000000in}}%
\pgfpathlineto{\pgfqpoint{-0.048611in}{0.000000in}}%
\pgfusepath{stroke,fill}%
}%
\begin{pgfscope}%
\pgfsys@transformshift{1.104012in}{3.537524in}%
\pgfsys@useobject{currentmarker}{}%
\end{pgfscope}%
\end{pgfscope}%
\begin{pgfscope}%
\definecolor{textcolor}{rgb}{0.000000,0.000000,0.000000}%
\pgfsetstrokecolor{textcolor}%
\pgfsetfillcolor{textcolor}%
\pgftext[x=0.419297in, y=3.468079in, left, base]{\color{textcolor}\rmfamily\fontsize{14.000000}{16.800000}\selectfont \(\displaystyle {210000}\)}%
\end{pgfscope}%
\begin{pgfscope}%
\definecolor{textcolor}{rgb}{0.000000,0.000000,0.000000}%
\pgfsetstrokecolor{textcolor}%
\pgfsetfillcolor{textcolor}%
\pgftext[x=0.363741in,y=2.248582in,,bottom,rotate=90.000000]{\color{textcolor}\rmfamily\fontsize{14.000000}{16.800000}\selectfont DKK}%
\end{pgfscope}%
\begin{pgfscope}%
\pgfpathrectangle{\pgfqpoint{1.104012in}{0.707164in}}{\pgfqpoint{4.685988in}{3.082836in}}%
\pgfusepath{clip}%
\pgfsetrectcap%
\pgfsetroundjoin%
\pgfsetlinewidth{1.505625pt}%
\definecolor{currentstroke}{rgb}{0.839216,0.152941,0.156863}%
\pgfsetstrokecolor{currentstroke}%
\pgfsetdash{}{0pt}%
\pgfpathmoveto{\pgfqpoint{1.327154in}{3.649871in}}%
\pgfpathlineto{\pgfqpoint{1.661867in}{1.803564in}}%
\pgfpathlineto{\pgfqpoint{2.219723in}{1.320877in}}%
\pgfpathlineto{\pgfqpoint{2.777579in}{1.108718in}}%
\pgfpathlineto{\pgfqpoint{3.335435in}{0.847293in}}%
\pgfpathlineto{\pgfqpoint{3.893290in}{0.930895in}}%
\pgfpathlineto{\pgfqpoint{4.451146in}{0.976009in}}%
\pgfpathlineto{\pgfqpoint{5.566858in}{1.370403in}}%
\pgfusepath{stroke}%
\end{pgfscope}%
\begin{pgfscope}%
\pgfpathrectangle{\pgfqpoint{1.104012in}{0.707164in}}{\pgfqpoint{4.685988in}{3.082836in}}%
\pgfusepath{clip}%
\pgfsetbuttcap%
\pgfsetroundjoin%
\definecolor{currentfill}{rgb}{0.839216,0.152941,0.156863}%
\pgfsetfillcolor{currentfill}%
\pgfsetlinewidth{1.003750pt}%
\definecolor{currentstroke}{rgb}{0.839216,0.152941,0.156863}%
\pgfsetstrokecolor{currentstroke}%
\pgfsetdash{}{0pt}%
\pgfsys@defobject{currentmarker}{\pgfqpoint{-0.041667in}{-0.041667in}}{\pgfqpoint{0.041667in}{0.041667in}}{%
\pgfpathmoveto{\pgfqpoint{0.000000in}{-0.041667in}}%
\pgfpathcurveto{\pgfqpoint{0.011050in}{-0.041667in}}{\pgfqpoint{0.021649in}{-0.037276in}}{\pgfqpoint{0.029463in}{-0.029463in}}%
\pgfpathcurveto{\pgfqpoint{0.037276in}{-0.021649in}}{\pgfqpoint{0.041667in}{-0.011050in}}{\pgfqpoint{0.041667in}{0.000000in}}%
\pgfpathcurveto{\pgfqpoint{0.041667in}{0.011050in}}{\pgfqpoint{0.037276in}{0.021649in}}{\pgfqpoint{0.029463in}{0.029463in}}%
\pgfpathcurveto{\pgfqpoint{0.021649in}{0.037276in}}{\pgfqpoint{0.011050in}{0.041667in}}{\pgfqpoint{0.000000in}{0.041667in}}%
\pgfpathcurveto{\pgfqpoint{-0.011050in}{0.041667in}}{\pgfqpoint{-0.021649in}{0.037276in}}{\pgfqpoint{-0.029463in}{0.029463in}}%
\pgfpathcurveto{\pgfqpoint{-0.037276in}{0.021649in}}{\pgfqpoint{-0.041667in}{0.011050in}}{\pgfqpoint{-0.041667in}{0.000000in}}%
\pgfpathcurveto{\pgfqpoint{-0.041667in}{-0.011050in}}{\pgfqpoint{-0.037276in}{-0.021649in}}{\pgfqpoint{-0.029463in}{-0.029463in}}%
\pgfpathcurveto{\pgfqpoint{-0.021649in}{-0.037276in}}{\pgfqpoint{-0.011050in}{-0.041667in}}{\pgfqpoint{0.000000in}{-0.041667in}}%
\pgfpathlineto{\pgfqpoint{0.000000in}{-0.041667in}}%
\pgfpathclose%
\pgfusepath{stroke,fill}%
}%
\begin{pgfscope}%
\pgfsys@transformshift{1.327154in}{3.649871in}%
\pgfsys@useobject{currentmarker}{}%
\end{pgfscope}%
\begin{pgfscope}%
\pgfsys@transformshift{1.661867in}{1.803564in}%
\pgfsys@useobject{currentmarker}{}%
\end{pgfscope}%
\begin{pgfscope}%
\pgfsys@transformshift{2.219723in}{1.320877in}%
\pgfsys@useobject{currentmarker}{}%
\end{pgfscope}%
\begin{pgfscope}%
\pgfsys@transformshift{2.777579in}{1.108718in}%
\pgfsys@useobject{currentmarker}{}%
\end{pgfscope}%
\begin{pgfscope}%
\pgfsys@transformshift{3.335435in}{0.847293in}%
\pgfsys@useobject{currentmarker}{}%
\end{pgfscope}%
\begin{pgfscope}%
\pgfsys@transformshift{3.893290in}{0.930895in}%
\pgfsys@useobject{currentmarker}{}%
\end{pgfscope}%
\begin{pgfscope}%
\pgfsys@transformshift{4.451146in}{0.976009in}%
\pgfsys@useobject{currentmarker}{}%
\end{pgfscope}%
\begin{pgfscope}%
\pgfsys@transformshift{5.566858in}{1.370403in}%
\pgfsys@useobject{currentmarker}{}%
\end{pgfscope}%
\end{pgfscope}%
\begin{pgfscope}%
\pgfsetrectcap%
\pgfsetmiterjoin%
\pgfsetlinewidth{0.803000pt}%
\definecolor{currentstroke}{rgb}{0.000000,0.000000,0.000000}%
\pgfsetstrokecolor{currentstroke}%
\pgfsetdash{}{0pt}%
\pgfpathmoveto{\pgfqpoint{1.104012in}{0.707164in}}%
\pgfpathlineto{\pgfqpoint{1.104012in}{3.790000in}}%
\pgfusepath{stroke}%
\end{pgfscope}%
\begin{pgfscope}%
\pgfsetrectcap%
\pgfsetmiterjoin%
\pgfsetlinewidth{0.803000pt}%
\definecolor{currentstroke}{rgb}{0.000000,0.000000,0.000000}%
\pgfsetstrokecolor{currentstroke}%
\pgfsetdash{}{0pt}%
\pgfpathmoveto{\pgfqpoint{5.790000in}{0.707164in}}%
\pgfpathlineto{\pgfqpoint{5.790000in}{3.790000in}}%
\pgfusepath{stroke}%
\end{pgfscope}%
\begin{pgfscope}%
\pgfsetrectcap%
\pgfsetmiterjoin%
\pgfsetlinewidth{0.803000pt}%
\definecolor{currentstroke}{rgb}{0.000000,0.000000,0.000000}%
\pgfsetstrokecolor{currentstroke}%
\pgfsetdash{}{0pt}%
\pgfpathmoveto{\pgfqpoint{1.104012in}{0.707164in}}%
\pgfpathlineto{\pgfqpoint{5.790000in}{0.707164in}}%
\pgfusepath{stroke}%
\end{pgfscope}%
\begin{pgfscope}%
\pgfsetrectcap%
\pgfsetmiterjoin%
\pgfsetlinewidth{0.803000pt}%
\definecolor{currentstroke}{rgb}{0.000000,0.000000,0.000000}%
\pgfsetstrokecolor{currentstroke}%
\pgfsetdash{}{0pt}%
\pgfpathmoveto{\pgfqpoint{1.104012in}{3.790000in}}%
\pgfpathlineto{\pgfqpoint{5.790000in}{3.790000in}}%
\pgfusepath{stroke}%
\end{pgfscope}%
\begin{pgfscope}%
\pgfsetbuttcap%
\pgfsetmiterjoin%
\definecolor{currentfill}{rgb}{1.000000,1.000000,1.000000}%
\pgfsetfillcolor{currentfill}%
\pgfsetfillopacity{0.800000}%
\pgfsetlinewidth{1.003750pt}%
\definecolor{currentstroke}{rgb}{0.800000,0.800000,0.800000}%
\pgfsetstrokecolor{currentstroke}%
\pgfsetstrokeopacity{0.800000}%
\pgfsetdash{}{0pt}%
\pgfpathmoveto{\pgfqpoint{2.877528in}{3.337223in}}%
\pgfpathlineto{\pgfqpoint{5.653889in}{3.337223in}}%
\pgfpathquadraticcurveto{\pgfqpoint{5.692778in}{3.337223in}}{\pgfqpoint{5.692778in}{3.376111in}}%
\pgfpathlineto{\pgfqpoint{5.692778in}{3.653889in}}%
\pgfpathquadraticcurveto{\pgfqpoint{5.692778in}{3.692778in}}{\pgfqpoint{5.653889in}{3.692778in}}%
\pgfpathlineto{\pgfqpoint{2.877528in}{3.692778in}}%
\pgfpathquadraticcurveto{\pgfqpoint{2.838639in}{3.692778in}}{\pgfqpoint{2.838639in}{3.653889in}}%
\pgfpathlineto{\pgfqpoint{2.838639in}{3.376111in}}%
\pgfpathquadraticcurveto{\pgfqpoint{2.838639in}{3.337223in}}{\pgfqpoint{2.877528in}{3.337223in}}%
\pgfpathlineto{\pgfqpoint{2.877528in}{3.337223in}}%
\pgfpathclose%
\pgfusepath{stroke,fill}%
\end{pgfscope}%
\begin{pgfscope}%
\pgfsetrectcap%
\pgfsetroundjoin%
\pgfsetlinewidth{1.505625pt}%
\definecolor{currentstroke}{rgb}{0.839216,0.152941,0.156863}%
\pgfsetstrokecolor{currentstroke}%
\pgfsetdash{}{0pt}%
\pgfpathmoveto{\pgfqpoint{2.916417in}{3.533056in}}%
\pgfpathlineto{\pgfqpoint{3.110861in}{3.533056in}}%
\pgfpathlineto{\pgfqpoint{3.305306in}{3.533056in}}%
\pgfusepath{stroke}%
\end{pgfscope}%
\begin{pgfscope}%
\pgfsetbuttcap%
\pgfsetroundjoin%
\definecolor{currentfill}{rgb}{0.839216,0.152941,0.156863}%
\pgfsetfillcolor{currentfill}%
\pgfsetlinewidth{1.003750pt}%
\definecolor{currentstroke}{rgb}{0.839216,0.152941,0.156863}%
\pgfsetstrokecolor{currentstroke}%
\pgfsetdash{}{0pt}%
\pgfsys@defobject{currentmarker}{\pgfqpoint{-0.041667in}{-0.041667in}}{\pgfqpoint{0.041667in}{0.041667in}}{%
\pgfpathmoveto{\pgfqpoint{0.000000in}{-0.041667in}}%
\pgfpathcurveto{\pgfqpoint{0.011050in}{-0.041667in}}{\pgfqpoint{0.021649in}{-0.037276in}}{\pgfqpoint{0.029463in}{-0.029463in}}%
\pgfpathcurveto{\pgfqpoint{0.037276in}{-0.021649in}}{\pgfqpoint{0.041667in}{-0.011050in}}{\pgfqpoint{0.041667in}{0.000000in}}%
\pgfpathcurveto{\pgfqpoint{0.041667in}{0.011050in}}{\pgfqpoint{0.037276in}{0.021649in}}{\pgfqpoint{0.029463in}{0.029463in}}%
\pgfpathcurveto{\pgfqpoint{0.021649in}{0.037276in}}{\pgfqpoint{0.011050in}{0.041667in}}{\pgfqpoint{0.000000in}{0.041667in}}%
\pgfpathcurveto{\pgfqpoint{-0.011050in}{0.041667in}}{\pgfqpoint{-0.021649in}{0.037276in}}{\pgfqpoint{-0.029463in}{0.029463in}}%
\pgfpathcurveto{\pgfqpoint{-0.037276in}{0.021649in}}{\pgfqpoint{-0.041667in}{0.011050in}}{\pgfqpoint{-0.041667in}{0.000000in}}%
\pgfpathcurveto{\pgfqpoint{-0.041667in}{-0.011050in}}{\pgfqpoint{-0.037276in}{-0.021649in}}{\pgfqpoint{-0.029463in}{-0.029463in}}%
\pgfpathcurveto{\pgfqpoint{-0.021649in}{-0.037276in}}{\pgfqpoint{-0.011050in}{-0.041667in}}{\pgfqpoint{0.000000in}{-0.041667in}}%
\pgfpathlineto{\pgfqpoint{0.000000in}{-0.041667in}}%
\pgfpathclose%
\pgfusepath{stroke,fill}%
}%
\begin{pgfscope}%
\pgfsys@transformshift{3.110861in}{3.533056in}%
\pgfsys@useobject{currentmarker}{}%
\end{pgfscope}%
\end{pgfscope}%
\begin{pgfscope}%
\definecolor{textcolor}{rgb}{0.000000,0.000000,0.000000}%
\pgfsetstrokecolor{textcolor}%
\pgfsetfillcolor{textcolor}%
\pgftext[x=3.460861in,y=3.465000in,left,base]{\color{textcolor}\rmfamily\fontsize{14.000000}{16.800000}\selectfont Objective (\(\displaystyle \lambda\)-normalized)}%
\end{pgfscope}%
\end{pgfpicture}%
\makeatother%
\endgroup%

%% file: figures/Comp/Inst_CostsVSLambda.pgf
\begingroup%
\makeatletter%
\begin{pgfpicture}%
\pgfpathrectangle{\pgfpointorigin}{\pgfqpoint{6.000000in}{4.000000in}}%
\pgfusepath{use as bounding box, clip}%
\begin{pgfscope}%
\pgfsetbuttcap%
\pgfsetmiterjoin%
\definecolor{currentfill}{rgb}{1.000000,1.000000,1.000000}%
\pgfsetfillcolor{currentfill}%
\pgfsetlinewidth{0.000000pt}%
\definecolor{currentstroke}{rgb}{1.000000,1.000000,1.000000}%
\pgfsetstrokecolor{currentstroke}%
\pgfsetdash{}{0pt}%
\pgfpathmoveto{\pgfqpoint{0.000000in}{0.000000in}}%
\pgfpathlineto{\pgfqpoint{6.000000in}{0.000000in}}%
\pgfpathlineto{\pgfqpoint{6.000000in}{4.000000in}}%
\pgfpathlineto{\pgfqpoint{0.000000in}{4.000000in}}%
\pgfpathlineto{\pgfqpoint{0.000000in}{0.000000in}}%
\pgfpathclose%
\pgfusepath{fill}%
\end{pgfscope}%
\begin{pgfscope}%
\pgfsetbuttcap%
\pgfsetmiterjoin%
\definecolor{currentfill}{rgb}{1.000000,1.000000,1.000000}%
\pgfsetfillcolor{currentfill}%
\pgfsetlinewidth{0.000000pt}%
\definecolor{currentstroke}{rgb}{0.000000,0.000000,0.000000}%
\pgfsetstrokecolor{currentstroke}%
\pgfsetstrokeopacity{0.000000}%
\pgfsetdash{}{0pt}%
\pgfpathmoveto{\pgfqpoint{1.104012in}{0.707164in}}%
\pgfpathlineto{\pgfqpoint{5.790000in}{0.707164in}}%
\pgfpathlineto{\pgfqpoint{5.790000in}{3.790000in}}%
\pgfpathlineto{\pgfqpoint{1.104012in}{3.790000in}}%
\pgfpathlineto{\pgfqpoint{1.104012in}{0.707164in}}%
\pgfpathclose%
\pgfusepath{fill}%
\end{pgfscope}%
\begin{pgfscope}%
\pgfsetbuttcap%
\pgfsetroundjoin%
\definecolor{currentfill}{rgb}{0.000000,0.000000,0.000000}%
\pgfsetfillcolor{currentfill}%
\pgfsetlinewidth{0.803000pt}%
\definecolor{currentstroke}{rgb}{0.000000,0.000000,0.000000}%
\pgfsetstrokecolor{currentstroke}%
\pgfsetdash{}{0pt}%
\pgfsys@defobject{currentmarker}{\pgfqpoint{0.000000in}{-0.048611in}}{\pgfqpoint{0.000000in}{0.000000in}}{%
\pgfpathmoveto{\pgfqpoint{0.000000in}{0.000000in}}%
\pgfpathlineto{\pgfqpoint{0.000000in}{-0.048611in}}%
\pgfusepath{stroke,fill}%
}%
\begin{pgfscope}%
\pgfsys@transformshift{1.104012in}{0.707164in}%
\pgfsys@useobject{currentmarker}{}%
\end{pgfscope}%
\end{pgfscope}%
\begin{pgfscope}%
\definecolor{textcolor}{rgb}{0.000000,0.000000,0.000000}%
\pgfsetstrokecolor{textcolor}%
\pgfsetfillcolor{textcolor}%
\pgftext[x=1.104012in,y=0.609942in,,top]{\color{textcolor}\rmfamily\fontsize{14.000000}{16.800000}\selectfont \(\displaystyle {0.0}\)}%
\end{pgfscope}%
\begin{pgfscope}%
\pgfsetbuttcap%
\pgfsetroundjoin%
\definecolor{currentfill}{rgb}{0.000000,0.000000,0.000000}%
\pgfsetfillcolor{currentfill}%
\pgfsetlinewidth{0.803000pt}%
\definecolor{currentstroke}{rgb}{0.000000,0.000000,0.000000}%
\pgfsetstrokecolor{currentstroke}%
\pgfsetdash{}{0pt}%
\pgfsys@defobject{currentmarker}{\pgfqpoint{0.000000in}{-0.048611in}}{\pgfqpoint{0.000000in}{0.000000in}}{%
\pgfpathmoveto{\pgfqpoint{0.000000in}{0.000000in}}%
\pgfpathlineto{\pgfqpoint{0.000000in}{-0.048611in}}%
\pgfusepath{stroke,fill}%
}%
\begin{pgfscope}%
\pgfsys@transformshift{2.219723in}{0.707164in}%
\pgfsys@useobject{currentmarker}{}%
\end{pgfscope}%
\end{pgfscope}%
\begin{pgfscope}%
\definecolor{textcolor}{rgb}{0.000000,0.000000,0.000000}%
\pgfsetstrokecolor{textcolor}%
\pgfsetfillcolor{textcolor}%
\pgftext[x=2.219723in,y=0.609942in,,top]{\color{textcolor}\rmfamily\fontsize{14.000000}{16.800000}\selectfont \(\displaystyle {0.5}\)}%
\end{pgfscope}%
\begin{pgfscope}%
\pgfsetbuttcap%
\pgfsetroundjoin%
\definecolor{currentfill}{rgb}{0.000000,0.000000,0.000000}%
\pgfsetfillcolor{currentfill}%
\pgfsetlinewidth{0.803000pt}%
\definecolor{currentstroke}{rgb}{0.000000,0.000000,0.000000}%
\pgfsetstrokecolor{currentstroke}%
\pgfsetdash{}{0pt}%
\pgfsys@defobject{currentmarker}{\pgfqpoint{0.000000in}{-0.048611in}}{\pgfqpoint{0.000000in}{0.000000in}}{%
\pgfpathmoveto{\pgfqpoint{0.000000in}{0.000000in}}%
\pgfpathlineto{\pgfqpoint{0.000000in}{-0.048611in}}%
\pgfusepath{stroke,fill}%
}%
\begin{pgfscope}%
\pgfsys@transformshift{3.335435in}{0.707164in}%
\pgfsys@useobject{currentmarker}{}%
\end{pgfscope}%
\end{pgfscope}%
\begin{pgfscope}%
\definecolor{textcolor}{rgb}{0.000000,0.000000,0.000000}%
\pgfsetstrokecolor{textcolor}%
\pgfsetfillcolor{textcolor}%
\pgftext[x=3.335435in,y=0.609942in,,top]{\color{textcolor}\rmfamily\fontsize{14.000000}{16.800000}\selectfont \(\displaystyle {1.0}\)}%
\end{pgfscope}%
\begin{pgfscope}%
\pgfsetbuttcap%
\pgfsetroundjoin%
\definecolor{currentfill}{rgb}{0.000000,0.000000,0.000000}%
\pgfsetfillcolor{currentfill}%
\pgfsetlinewidth{0.803000pt}%
\definecolor{currentstroke}{rgb}{0.000000,0.000000,0.000000}%
\pgfsetstrokecolor{currentstroke}%
\pgfsetdash{}{0pt}%
\pgfsys@defobject{currentmarker}{\pgfqpoint{0.000000in}{-0.048611in}}{\pgfqpoint{0.000000in}{0.000000in}}{%
\pgfpathmoveto{\pgfqpoint{0.000000in}{0.000000in}}%
\pgfpathlineto{\pgfqpoint{0.000000in}{-0.048611in}}%
\pgfusepath{stroke,fill}%
}%
\begin{pgfscope}%
\pgfsys@transformshift{4.451146in}{0.707164in}%
\pgfsys@useobject{currentmarker}{}%
\end{pgfscope}%
\end{pgfscope}%
\begin{pgfscope}%
\definecolor{textcolor}{rgb}{0.000000,0.000000,0.000000}%
\pgfsetstrokecolor{textcolor}%
\pgfsetfillcolor{textcolor}%
\pgftext[x=4.451146in,y=0.609942in,,top]{\color{textcolor}\rmfamily\fontsize{14.000000}{16.800000}\selectfont \(\displaystyle {1.5}\)}%
\end{pgfscope}%
\begin{pgfscope}%
\pgfsetbuttcap%
\pgfsetroundjoin%
\definecolor{currentfill}{rgb}{0.000000,0.000000,0.000000}%
\pgfsetfillcolor{currentfill}%
\pgfsetlinewidth{0.803000pt}%
\definecolor{currentstroke}{rgb}{0.000000,0.000000,0.000000}%
\pgfsetstrokecolor{currentstroke}%
\pgfsetdash{}{0pt}%
\pgfsys@defobject{currentmarker}{\pgfqpoint{0.000000in}{-0.048611in}}{\pgfqpoint{0.000000in}{0.000000in}}{%
\pgfpathmoveto{\pgfqpoint{0.000000in}{0.000000in}}%
\pgfpathlineto{\pgfqpoint{0.000000in}{-0.048611in}}%
\pgfusepath{stroke,fill}%
}%
\begin{pgfscope}%
\pgfsys@transformshift{5.566858in}{0.707164in}%
\pgfsys@useobject{currentmarker}{}%
\end{pgfscope}%
\end{pgfscope}%
\begin{pgfscope}%
\definecolor{textcolor}{rgb}{0.000000,0.000000,0.000000}%
\pgfsetstrokecolor{textcolor}%
\pgfsetfillcolor{textcolor}%
\pgftext[x=5.566858in,y=0.609942in,,top]{\color{textcolor}\rmfamily\fontsize{14.000000}{16.800000}\selectfont \(\displaystyle {2.0}\)}%
\end{pgfscope}%
\begin{pgfscope}%
\definecolor{textcolor}{rgb}{0.000000,0.000000,0.000000}%
\pgfsetstrokecolor{textcolor}%
\pgfsetfillcolor{textcolor}%
\pgftext[x=3.447006in,y=0.376609in,,top]{\color{textcolor}\rmfamily\fontsize{14.000000}{16.800000}\selectfont \(\displaystyle \lambda\)}%
\end{pgfscope}%
\begin{pgfscope}%
\pgfsetbuttcap%
\pgfsetroundjoin%
\definecolor{currentfill}{rgb}{0.000000,0.000000,0.000000}%
\pgfsetfillcolor{currentfill}%
\pgfsetlinewidth{0.803000pt}%
\definecolor{currentstroke}{rgb}{0.000000,0.000000,0.000000}%
\pgfsetstrokecolor{currentstroke}%
\pgfsetdash{}{0pt}%
\pgfsys@defobject{currentmarker}{\pgfqpoint{-0.048611in}{0.000000in}}{\pgfqpoint{-0.000000in}{0.000000in}}{%
\pgfpathmoveto{\pgfqpoint{-0.000000in}{0.000000in}}%
\pgfpathlineto{\pgfqpoint{-0.048611in}{0.000000in}}%
\pgfusepath{stroke,fill}%
}%
\begin{pgfscope}%
\pgfsys@transformshift{1.104012in}{0.982289in}%
\pgfsys@useobject{currentmarker}{}%
\end{pgfscope}%
\end{pgfscope}%
\begin{pgfscope}%
\definecolor{textcolor}{rgb}{0.000000,0.000000,0.000000}%
\pgfsetstrokecolor{textcolor}%
\pgfsetfillcolor{textcolor}%
\pgftext[x=0.517212in, y=0.912844in, left, base]{\color{textcolor}\rmfamily\fontsize{14.000000}{16.800000}\selectfont \(\displaystyle {50000}\)}%
\end{pgfscope}%
\begin{pgfscope}%
\pgfsetbuttcap%
\pgfsetroundjoin%
\definecolor{currentfill}{rgb}{0.000000,0.000000,0.000000}%
\pgfsetfillcolor{currentfill}%
\pgfsetlinewidth{0.803000pt}%
\definecolor{currentstroke}{rgb}{0.000000,0.000000,0.000000}%
\pgfsetstrokecolor{currentstroke}%
\pgfsetdash{}{0pt}%
\pgfsys@defobject{currentmarker}{\pgfqpoint{-0.048611in}{0.000000in}}{\pgfqpoint{-0.000000in}{0.000000in}}{%
\pgfpathmoveto{\pgfqpoint{-0.000000in}{0.000000in}}%
\pgfpathlineto{\pgfqpoint{-0.048611in}{0.000000in}}%
\pgfusepath{stroke,fill}%
}%
\begin{pgfscope}%
\pgfsys@transformshift{1.104012in}{1.691299in}%
\pgfsys@useobject{currentmarker}{}%
\end{pgfscope}%
\end{pgfscope}%
\begin{pgfscope}%
\definecolor{textcolor}{rgb}{0.000000,0.000000,0.000000}%
\pgfsetstrokecolor{textcolor}%
\pgfsetfillcolor{textcolor}%
\pgftext[x=0.419297in, y=1.621855in, left, base]{\color{textcolor}\rmfamily\fontsize{14.000000}{16.800000}\selectfont \(\displaystyle {100000}\)}%
\end{pgfscope}%
\begin{pgfscope}%
\pgfsetbuttcap%
\pgfsetroundjoin%
\definecolor{currentfill}{rgb}{0.000000,0.000000,0.000000}%
\pgfsetfillcolor{currentfill}%
\pgfsetlinewidth{0.803000pt}%
\definecolor{currentstroke}{rgb}{0.000000,0.000000,0.000000}%
\pgfsetstrokecolor{currentstroke}%
\pgfsetdash{}{0pt}%
\pgfsys@defobject{currentmarker}{\pgfqpoint{-0.048611in}{0.000000in}}{\pgfqpoint{-0.000000in}{0.000000in}}{%
\pgfpathmoveto{\pgfqpoint{-0.000000in}{0.000000in}}%
\pgfpathlineto{\pgfqpoint{-0.048611in}{0.000000in}}%
\pgfusepath{stroke,fill}%
}%
\begin{pgfscope}%
\pgfsys@transformshift{1.104012in}{2.400310in}%
\pgfsys@useobject{currentmarker}{}%
\end{pgfscope}%
\end{pgfscope}%
\begin{pgfscope}%
\definecolor{textcolor}{rgb}{0.000000,0.000000,0.000000}%
\pgfsetstrokecolor{textcolor}%
\pgfsetfillcolor{textcolor}%
\pgftext[x=0.419297in, y=2.330866in, left, base]{\color{textcolor}\rmfamily\fontsize{14.000000}{16.800000}\selectfont \(\displaystyle {150000}\)}%
\end{pgfscope}%
\begin{pgfscope}%
\pgfsetbuttcap%
\pgfsetroundjoin%
\definecolor{currentfill}{rgb}{0.000000,0.000000,0.000000}%
\pgfsetfillcolor{currentfill}%
\pgfsetlinewidth{0.803000pt}%
\definecolor{currentstroke}{rgb}{0.000000,0.000000,0.000000}%
\pgfsetstrokecolor{currentstroke}%
\pgfsetdash{}{0pt}%
\pgfsys@defobject{currentmarker}{\pgfqpoint{-0.048611in}{0.000000in}}{\pgfqpoint{-0.000000in}{0.000000in}}{%
\pgfpathmoveto{\pgfqpoint{-0.000000in}{0.000000in}}%
\pgfpathlineto{\pgfqpoint{-0.048611in}{0.000000in}}%
\pgfusepath{stroke,fill}%
}%
\begin{pgfscope}%
\pgfsys@transformshift{1.104012in}{3.109321in}%
\pgfsys@useobject{currentmarker}{}%
\end{pgfscope}%
\end{pgfscope}%
\begin{pgfscope}%
\definecolor{textcolor}{rgb}{0.000000,0.000000,0.000000}%
\pgfsetstrokecolor{textcolor}%
\pgfsetfillcolor{textcolor}%
\pgftext[x=0.419297in, y=3.039877in, left, base]{\color{textcolor}\rmfamily\fontsize{14.000000}{16.800000}\selectfont \(\displaystyle {200000}\)}%
\end{pgfscope}%
\begin{pgfscope}%
\definecolor{textcolor}{rgb}{0.000000,0.000000,0.000000}%
\pgfsetstrokecolor{textcolor}%
\pgfsetfillcolor{textcolor}%
\pgftext[x=0.363741in,y=2.248582in,,bottom,rotate=90.000000]{\color{textcolor}\rmfamily\fontsize{14.000000}{16.800000}\selectfont DKK}%
\end{pgfscope}%
\begin{pgfscope}%
\pgfpathrectangle{\pgfqpoint{1.104012in}{0.707164in}}{\pgfqpoint{4.685988in}{3.082836in}}%
\pgfusepath{clip}%
\pgfsetrectcap%
\pgfsetroundjoin%
\pgfsetlinewidth{1.505625pt}%
\definecolor{currentstroke}{rgb}{0.000000,0.000000,0.000000}%
\pgfsetstrokecolor{currentstroke}%
\pgfsetdash{}{0pt}%
\pgfpathmoveto{\pgfqpoint{1.327154in}{3.649871in}}%
\pgfpathlineto{\pgfqpoint{1.661867in}{3.398881in}}%
\pgfpathlineto{\pgfqpoint{2.219723in}{3.302768in}}%
\pgfpathlineto{\pgfqpoint{2.777579in}{3.227735in}}%
\pgfpathlineto{\pgfqpoint{3.335435in}{3.051012in}}%
\pgfpathlineto{\pgfqpoint{3.893290in}{2.976141in}}%
\pgfpathlineto{\pgfqpoint{4.451146in}{2.963757in}}%
\pgfpathlineto{\pgfqpoint{5.566858in}{2.885678in}}%
\pgfusepath{stroke}%
\end{pgfscope}%
\begin{pgfscope}%
\pgfpathrectangle{\pgfqpoint{1.104012in}{0.707164in}}{\pgfqpoint{4.685988in}{3.082836in}}%
\pgfusepath{clip}%
\pgfsetbuttcap%
\pgfsetroundjoin%
\definecolor{currentfill}{rgb}{0.000000,0.000000,0.000000}%
\pgfsetfillcolor{currentfill}%
\pgfsetlinewidth{1.003750pt}%
\definecolor{currentstroke}{rgb}{0.000000,0.000000,0.000000}%
\pgfsetstrokecolor{currentstroke}%
\pgfsetdash{}{0pt}%
\pgfsys@defobject{currentmarker}{\pgfqpoint{-0.041667in}{-0.041667in}}{\pgfqpoint{0.041667in}{0.041667in}}{%
\pgfpathmoveto{\pgfqpoint{0.000000in}{-0.041667in}}%
\pgfpathcurveto{\pgfqpoint{0.011050in}{-0.041667in}}{\pgfqpoint{0.021649in}{-0.037276in}}{\pgfqpoint{0.029463in}{-0.029463in}}%
\pgfpathcurveto{\pgfqpoint{0.037276in}{-0.021649in}}{\pgfqpoint{0.041667in}{-0.011050in}}{\pgfqpoint{0.041667in}{0.000000in}}%
\pgfpathcurveto{\pgfqpoint{0.041667in}{0.011050in}}{\pgfqpoint{0.037276in}{0.021649in}}{\pgfqpoint{0.029463in}{0.029463in}}%
\pgfpathcurveto{\pgfqpoint{0.021649in}{0.037276in}}{\pgfqpoint{0.011050in}{0.041667in}}{\pgfqpoint{0.000000in}{0.041667in}}%
\pgfpathcurveto{\pgfqpoint{-0.011050in}{0.041667in}}{\pgfqpoint{-0.021649in}{0.037276in}}{\pgfqpoint{-0.029463in}{0.029463in}}%
\pgfpathcurveto{\pgfqpoint{-0.037276in}{0.021649in}}{\pgfqpoint{-0.041667in}{0.011050in}}{\pgfqpoint{-0.041667in}{0.000000in}}%
\pgfpathcurveto{\pgfqpoint{-0.041667in}{-0.011050in}}{\pgfqpoint{-0.037276in}{-0.021649in}}{\pgfqpoint{-0.029463in}{-0.029463in}}%
\pgfpathcurveto{\pgfqpoint{-0.021649in}{-0.037276in}}{\pgfqpoint{-0.011050in}{-0.041667in}}{\pgfqpoint{0.000000in}{-0.041667in}}%
\pgfpathlineto{\pgfqpoint{0.000000in}{-0.041667in}}%
\pgfpathclose%
\pgfusepath{stroke,fill}%
}%
\begin{pgfscope}%
\pgfsys@transformshift{1.327154in}{3.649871in}%
\pgfsys@useobject{currentmarker}{}%
\end{pgfscope}%
\begin{pgfscope}%
\pgfsys@transformshift{1.661867in}{3.398881in}%
\pgfsys@useobject{currentmarker}{}%
\end{pgfscope}%
\begin{pgfscope}%
\pgfsys@transformshift{2.219723in}{3.302768in}%
\pgfsys@useobject{currentmarker}{}%
\end{pgfscope}%
\begin{pgfscope}%
\pgfsys@transformshift{2.777579in}{3.227735in}%
\pgfsys@useobject{currentmarker}{}%
\end{pgfscope}%
\begin{pgfscope}%
\pgfsys@transformshift{3.335435in}{3.051012in}%
\pgfsys@useobject{currentmarker}{}%
\end{pgfscope}%
\begin{pgfscope}%
\pgfsys@transformshift{3.893290in}{2.976141in}%
\pgfsys@useobject{currentmarker}{}%
\end{pgfscope}%
\begin{pgfscope}%
\pgfsys@transformshift{4.451146in}{2.963757in}%
\pgfsys@useobject{currentmarker}{}%
\end{pgfscope}%
\begin{pgfscope}%
\pgfsys@transformshift{5.566858in}{2.885678in}%
\pgfsys@useobject{currentmarker}{}%
\end{pgfscope}%
\end{pgfscope}%
\begin{pgfscope}%
\pgfpathrectangle{\pgfqpoint{1.104012in}{0.707164in}}{\pgfqpoint{4.685988in}{3.082836in}}%
\pgfusepath{clip}%
\pgfsetrectcap%
\pgfsetroundjoin%
\pgfsetlinewidth{1.505625pt}%
\definecolor{currentstroke}{rgb}{0.839216,0.152941,0.156863}%
\pgfsetstrokecolor{currentstroke}%
\pgfsetdash{}{0pt}%
\pgfpathmoveto{\pgfqpoint{1.327154in}{3.264424in}}%
\pgfpathlineto{\pgfqpoint{1.661867in}{3.045836in}}%
\pgfpathlineto{\pgfqpoint{2.219723in}{2.988690in}}%
\pgfpathlineto{\pgfqpoint{2.777579in}{2.963572in}}%
\pgfpathlineto{\pgfqpoint{3.335435in}{2.932622in}}%
\pgfpathlineto{\pgfqpoint{3.893290in}{2.942519in}}%
\pgfpathlineto{\pgfqpoint{4.451146in}{2.947861in}}%
\pgfpathlineto{\pgfqpoint{5.566858in}{2.994554in}}%
\pgfusepath{stroke}%
\end{pgfscope}%
\begin{pgfscope}%
\pgfpathrectangle{\pgfqpoint{1.104012in}{0.707164in}}{\pgfqpoint{4.685988in}{3.082836in}}%
\pgfusepath{clip}%
\pgfsetbuttcap%
\pgfsetroundjoin%
\definecolor{currentfill}{rgb}{0.839216,0.152941,0.156863}%
\pgfsetfillcolor{currentfill}%
\pgfsetlinewidth{1.003750pt}%
\definecolor{currentstroke}{rgb}{0.839216,0.152941,0.156863}%
\pgfsetstrokecolor{currentstroke}%
\pgfsetdash{}{0pt}%
\pgfsys@defobject{currentmarker}{\pgfqpoint{-0.041667in}{-0.041667in}}{\pgfqpoint{0.041667in}{0.041667in}}{%
\pgfpathmoveto{\pgfqpoint{0.000000in}{-0.041667in}}%
\pgfpathcurveto{\pgfqpoint{0.011050in}{-0.041667in}}{\pgfqpoint{0.021649in}{-0.037276in}}{\pgfqpoint{0.029463in}{-0.029463in}}%
\pgfpathcurveto{\pgfqpoint{0.037276in}{-0.021649in}}{\pgfqpoint{0.041667in}{-0.011050in}}{\pgfqpoint{0.041667in}{0.000000in}}%
\pgfpathcurveto{\pgfqpoint{0.041667in}{0.011050in}}{\pgfqpoint{0.037276in}{0.021649in}}{\pgfqpoint{0.029463in}{0.029463in}}%
\pgfpathcurveto{\pgfqpoint{0.021649in}{0.037276in}}{\pgfqpoint{0.011050in}{0.041667in}}{\pgfqpoint{0.000000in}{0.041667in}}%
\pgfpathcurveto{\pgfqpoint{-0.011050in}{0.041667in}}{\pgfqpoint{-0.021649in}{0.037276in}}{\pgfqpoint{-0.029463in}{0.029463in}}%
\pgfpathcurveto{\pgfqpoint{-0.037276in}{0.021649in}}{\pgfqpoint{-0.041667in}{0.011050in}}{\pgfqpoint{-0.041667in}{0.000000in}}%
\pgfpathcurveto{\pgfqpoint{-0.041667in}{-0.011050in}}{\pgfqpoint{-0.037276in}{-0.021649in}}{\pgfqpoint{-0.029463in}{-0.029463in}}%
\pgfpathcurveto{\pgfqpoint{-0.021649in}{-0.037276in}}{\pgfqpoint{-0.011050in}{-0.041667in}}{\pgfqpoint{0.000000in}{-0.041667in}}%
\pgfpathlineto{\pgfqpoint{0.000000in}{-0.041667in}}%
\pgfpathclose%
\pgfusepath{stroke,fill}%
}%
\begin{pgfscope}%
\pgfsys@transformshift{1.327154in}{3.264424in}%
\pgfsys@useobject{currentmarker}{}%
\end{pgfscope}%
\begin{pgfscope}%
\pgfsys@transformshift{1.661867in}{3.045836in}%
\pgfsys@useobject{currentmarker}{}%
\end{pgfscope}%
\begin{pgfscope}%
\pgfsys@transformshift{2.219723in}{2.988690in}%
\pgfsys@useobject{currentmarker}{}%
\end{pgfscope}%
\begin{pgfscope}%
\pgfsys@transformshift{2.777579in}{2.963572in}%
\pgfsys@useobject{currentmarker}{}%
\end{pgfscope}%
\begin{pgfscope}%
\pgfsys@transformshift{3.335435in}{2.932622in}%
\pgfsys@useobject{currentmarker}{}%
\end{pgfscope}%
\begin{pgfscope}%
\pgfsys@transformshift{3.893290in}{2.942519in}%
\pgfsys@useobject{currentmarker}{}%
\end{pgfscope}%
\begin{pgfscope}%
\pgfsys@transformshift{4.451146in}{2.947861in}%
\pgfsys@useobject{currentmarker}{}%
\end{pgfscope}%
\begin{pgfscope}%
\pgfsys@transformshift{5.566858in}{2.994554in}%
\pgfsys@useobject{currentmarker}{}%
\end{pgfscope}%
\end{pgfscope}%
\begin{pgfscope}%
\pgfpathrectangle{\pgfqpoint{1.104012in}{0.707164in}}{\pgfqpoint{4.685988in}{3.082836in}}%
\pgfusepath{clip}%
\pgfsetrectcap%
\pgfsetroundjoin%
\pgfsetlinewidth{1.505625pt}%
\definecolor{currentstroke}{rgb}{0.301961,0.301961,0.301961}%
\pgfsetstrokecolor{currentstroke}%
\pgfsetdash{}{0pt}%
\pgfpathmoveto{\pgfqpoint{1.327154in}{1.232740in}}%
\pgfpathlineto{\pgfqpoint{1.661867in}{1.300167in}}%
\pgfpathlineto{\pgfqpoint{2.219723in}{1.348549in}}%
\pgfpathlineto{\pgfqpoint{2.777579in}{1.348549in}}%
\pgfpathlineto{\pgfqpoint{3.335435in}{1.365027in}}%
\pgfpathlineto{\pgfqpoint{3.893290in}{1.367579in}}%
\pgfpathlineto{\pgfqpoint{4.451146in}{1.387304in}}%
\pgfpathlineto{\pgfqpoint{5.566858in}{1.387304in}}%
\pgfusepath{stroke}%
\end{pgfscope}%
\begin{pgfscope}%
\pgfpathrectangle{\pgfqpoint{1.104012in}{0.707164in}}{\pgfqpoint{4.685988in}{3.082836in}}%
\pgfusepath{clip}%
\pgfsetbuttcap%
\pgfsetroundjoin%
\definecolor{currentfill}{rgb}{0.301961,0.301961,0.301961}%
\pgfsetfillcolor{currentfill}%
\pgfsetlinewidth{1.003750pt}%
\definecolor{currentstroke}{rgb}{0.301961,0.301961,0.301961}%
\pgfsetstrokecolor{currentstroke}%
\pgfsetdash{}{0pt}%
\pgfsys@defobject{currentmarker}{\pgfqpoint{-0.041667in}{-0.041667in}}{\pgfqpoint{0.041667in}{0.041667in}}{%
\pgfpathmoveto{\pgfqpoint{0.000000in}{-0.041667in}}%
\pgfpathcurveto{\pgfqpoint{0.011050in}{-0.041667in}}{\pgfqpoint{0.021649in}{-0.037276in}}{\pgfqpoint{0.029463in}{-0.029463in}}%
\pgfpathcurveto{\pgfqpoint{0.037276in}{-0.021649in}}{\pgfqpoint{0.041667in}{-0.011050in}}{\pgfqpoint{0.041667in}{0.000000in}}%
\pgfpathcurveto{\pgfqpoint{0.041667in}{0.011050in}}{\pgfqpoint{0.037276in}{0.021649in}}{\pgfqpoint{0.029463in}{0.029463in}}%
\pgfpathcurveto{\pgfqpoint{0.021649in}{0.037276in}}{\pgfqpoint{0.011050in}{0.041667in}}{\pgfqpoint{0.000000in}{0.041667in}}%
\pgfpathcurveto{\pgfqpoint{-0.011050in}{0.041667in}}{\pgfqpoint{-0.021649in}{0.037276in}}{\pgfqpoint{-0.029463in}{0.029463in}}%
\pgfpathcurveto{\pgfqpoint{-0.037276in}{0.021649in}}{\pgfqpoint{-0.041667in}{0.011050in}}{\pgfqpoint{-0.041667in}{0.000000in}}%
\pgfpathcurveto{\pgfqpoint{-0.041667in}{-0.011050in}}{\pgfqpoint{-0.037276in}{-0.021649in}}{\pgfqpoint{-0.029463in}{-0.029463in}}%
\pgfpathcurveto{\pgfqpoint{-0.021649in}{-0.037276in}}{\pgfqpoint{-0.011050in}{-0.041667in}}{\pgfqpoint{0.000000in}{-0.041667in}}%
\pgfpathlineto{\pgfqpoint{0.000000in}{-0.041667in}}%
\pgfpathclose%
\pgfusepath{stroke,fill}%
}%
\begin{pgfscope}%
\pgfsys@transformshift{1.327154in}{1.232740in}%
\pgfsys@useobject{currentmarker}{}%
\end{pgfscope}%
\begin{pgfscope}%
\pgfsys@transformshift{1.661867in}{1.300167in}%
\pgfsys@useobject{currentmarker}{}%
\end{pgfscope}%
\begin{pgfscope}%
\pgfsys@transformshift{2.219723in}{1.348549in}%
\pgfsys@useobject{currentmarker}{}%
\end{pgfscope}%
\begin{pgfscope}%
\pgfsys@transformshift{2.777579in}{1.348549in}%
\pgfsys@useobject{currentmarker}{}%
\end{pgfscope}%
\begin{pgfscope}%
\pgfsys@transformshift{3.335435in}{1.365027in}%
\pgfsys@useobject{currentmarker}{}%
\end{pgfscope}%
\begin{pgfscope}%
\pgfsys@transformshift{3.893290in}{1.367579in}%
\pgfsys@useobject{currentmarker}{}%
\end{pgfscope}%
\begin{pgfscope}%
\pgfsys@transformshift{4.451146in}{1.387304in}%
\pgfsys@useobject{currentmarker}{}%
\end{pgfscope}%
\begin{pgfscope}%
\pgfsys@transformshift{5.566858in}{1.387304in}%
\pgfsys@useobject{currentmarker}{}%
\end{pgfscope}%
\end{pgfscope}%
\begin{pgfscope}%
\pgfpathrectangle{\pgfqpoint{1.104012in}{0.707164in}}{\pgfqpoint{4.685988in}{3.082836in}}%
\pgfusepath{clip}%
\pgfsetrectcap%
\pgfsetroundjoin%
\pgfsetlinewidth{1.505625pt}%
\definecolor{currentstroke}{rgb}{0.678431,0.709804,0.741176}%
\pgfsetstrokecolor{currentstroke}%
\pgfsetdash{}{0pt}%
\pgfpathmoveto{\pgfqpoint{1.327154in}{0.847293in}}%
\pgfpathlineto{\pgfqpoint{1.661867in}{0.947122in}}%
\pgfpathlineto{\pgfqpoint{2.219723in}{1.034472in}}%
\pgfpathlineto{\pgfqpoint{2.777579in}{1.084386in}}%
\pgfpathlineto{\pgfqpoint{3.335435in}{1.246608in}}%
\pgfpathlineto{\pgfqpoint{3.893290in}{1.333958in}}%
\pgfpathlineto{\pgfqpoint{4.451146in}{1.371394in}}%
\pgfpathlineto{\pgfqpoint{5.566858in}{1.496180in}}%
\pgfusepath{stroke}%
\end{pgfscope}%
\begin{pgfscope}%
\pgfpathrectangle{\pgfqpoint{1.104012in}{0.707164in}}{\pgfqpoint{4.685988in}{3.082836in}}%
\pgfusepath{clip}%
\pgfsetbuttcap%
\pgfsetroundjoin%
\definecolor{currentfill}{rgb}{0.678431,0.709804,0.741176}%
\pgfsetfillcolor{currentfill}%
\pgfsetlinewidth{1.003750pt}%
\definecolor{currentstroke}{rgb}{0.678431,0.709804,0.741176}%
\pgfsetstrokecolor{currentstroke}%
\pgfsetdash{}{0pt}%
\pgfsys@defobject{currentmarker}{\pgfqpoint{-0.041667in}{-0.041667in}}{\pgfqpoint{0.041667in}{0.041667in}}{%
\pgfpathmoveto{\pgfqpoint{0.000000in}{-0.041667in}}%
\pgfpathcurveto{\pgfqpoint{0.011050in}{-0.041667in}}{\pgfqpoint{0.021649in}{-0.037276in}}{\pgfqpoint{0.029463in}{-0.029463in}}%
\pgfpathcurveto{\pgfqpoint{0.037276in}{-0.021649in}}{\pgfqpoint{0.041667in}{-0.011050in}}{\pgfqpoint{0.041667in}{0.000000in}}%
\pgfpathcurveto{\pgfqpoint{0.041667in}{0.011050in}}{\pgfqpoint{0.037276in}{0.021649in}}{\pgfqpoint{0.029463in}{0.029463in}}%
\pgfpathcurveto{\pgfqpoint{0.021649in}{0.037276in}}{\pgfqpoint{0.011050in}{0.041667in}}{\pgfqpoint{0.000000in}{0.041667in}}%
\pgfpathcurveto{\pgfqpoint{-0.011050in}{0.041667in}}{\pgfqpoint{-0.021649in}{0.037276in}}{\pgfqpoint{-0.029463in}{0.029463in}}%
\pgfpathcurveto{\pgfqpoint{-0.037276in}{0.021649in}}{\pgfqpoint{-0.041667in}{0.011050in}}{\pgfqpoint{-0.041667in}{0.000000in}}%
\pgfpathcurveto{\pgfqpoint{-0.041667in}{-0.011050in}}{\pgfqpoint{-0.037276in}{-0.021649in}}{\pgfqpoint{-0.029463in}{-0.029463in}}%
\pgfpathcurveto{\pgfqpoint{-0.021649in}{-0.037276in}}{\pgfqpoint{-0.011050in}{-0.041667in}}{\pgfqpoint{0.000000in}{-0.041667in}}%
\pgfpathlineto{\pgfqpoint{0.000000in}{-0.041667in}}%
\pgfpathclose%
\pgfusepath{stroke,fill}%
}%
\begin{pgfscope}%
\pgfsys@transformshift{1.327154in}{0.847293in}%
\pgfsys@useobject{currentmarker}{}%
\end{pgfscope}%
\begin{pgfscope}%
\pgfsys@transformshift{1.661867in}{0.947122in}%
\pgfsys@useobject{currentmarker}{}%
\end{pgfscope}%
\begin{pgfscope}%
\pgfsys@transformshift{2.219723in}{1.034472in}%
\pgfsys@useobject{currentmarker}{}%
\end{pgfscope}%
\begin{pgfscope}%
\pgfsys@transformshift{2.777579in}{1.084386in}%
\pgfsys@useobject{currentmarker}{}%
\end{pgfscope}%
\begin{pgfscope}%
\pgfsys@transformshift{3.335435in}{1.246608in}%
\pgfsys@useobject{currentmarker}{}%
\end{pgfscope}%
\begin{pgfscope}%
\pgfsys@transformshift{3.893290in}{1.333958in}%
\pgfsys@useobject{currentmarker}{}%
\end{pgfscope}%
\begin{pgfscope}%
\pgfsys@transformshift{4.451146in}{1.371394in}%
\pgfsys@useobject{currentmarker}{}%
\end{pgfscope}%
\begin{pgfscope}%
\pgfsys@transformshift{5.566858in}{1.496180in}%
\pgfsys@useobject{currentmarker}{}%
\end{pgfscope}%
\end{pgfscope}%
\begin{pgfscope}%
\pgfsetrectcap%
\pgfsetmiterjoin%
\pgfsetlinewidth{0.803000pt}%
\definecolor{currentstroke}{rgb}{0.000000,0.000000,0.000000}%
\pgfsetstrokecolor{currentstroke}%
\pgfsetdash{}{0pt}%
\pgfpathmoveto{\pgfqpoint{1.104012in}{0.707164in}}%
\pgfpathlineto{\pgfqpoint{1.104012in}{3.790000in}}%
\pgfusepath{stroke}%
\end{pgfscope}%
\begin{pgfscope}%
\pgfsetrectcap%
\pgfsetmiterjoin%
\pgfsetlinewidth{0.803000pt}%
\definecolor{currentstroke}{rgb}{0.000000,0.000000,0.000000}%
\pgfsetstrokecolor{currentstroke}%
\pgfsetdash{}{0pt}%
\pgfpathmoveto{\pgfqpoint{5.790000in}{0.707164in}}%
\pgfpathlineto{\pgfqpoint{5.790000in}{3.790000in}}%
\pgfusepath{stroke}%
\end{pgfscope}%
\begin{pgfscope}%
\pgfsetrectcap%
\pgfsetmiterjoin%
\pgfsetlinewidth{0.803000pt}%
\definecolor{currentstroke}{rgb}{0.000000,0.000000,0.000000}%
\pgfsetstrokecolor{currentstroke}%
\pgfsetdash{}{0pt}%
\pgfpathmoveto{\pgfqpoint{1.104012in}{0.707164in}}%
\pgfpathlineto{\pgfqpoint{5.790000in}{0.707164in}}%
\pgfusepath{stroke}%
\end{pgfscope}%
\begin{pgfscope}%
\pgfsetrectcap%
\pgfsetmiterjoin%
\pgfsetlinewidth{0.803000pt}%
\definecolor{currentstroke}{rgb}{0.000000,0.000000,0.000000}%
\pgfsetstrokecolor{currentstroke}%
\pgfsetdash{}{0pt}%
\pgfpathmoveto{\pgfqpoint{1.104012in}{3.790000in}}%
\pgfpathlineto{\pgfqpoint{5.790000in}{3.790000in}}%
\pgfusepath{stroke}%
\end{pgfscope}%
\begin{pgfscope}%
\pgfsetbuttcap%
\pgfsetmiterjoin%
\definecolor{currentfill}{rgb}{1.000000,1.000000,1.000000}%
\pgfsetfillcolor{currentfill}%
\pgfsetfillopacity{0.800000}%
\pgfsetlinewidth{1.003750pt}%
\definecolor{currentstroke}{rgb}{0.800000,0.800000,0.800000}%
\pgfsetstrokecolor{currentstroke}%
\pgfsetstrokeopacity{0.800000}%
\pgfsetdash{}{0pt}%
\pgfpathmoveto{\pgfqpoint{2.380913in}{1.647194in}}%
\pgfpathlineto{\pgfqpoint{5.653889in}{1.647194in}}%
\pgfpathquadraticcurveto{\pgfqpoint{5.692778in}{1.647194in}}{\pgfqpoint{5.692778in}{1.686083in}}%
\pgfpathlineto{\pgfqpoint{5.692778in}{2.811081in}}%
\pgfpathquadraticcurveto{\pgfqpoint{5.692778in}{2.849970in}}{\pgfqpoint{5.653889in}{2.849970in}}%
\pgfpathlineto{\pgfqpoint{2.380913in}{2.849970in}}%
\pgfpathquadraticcurveto{\pgfqpoint{2.342024in}{2.849970in}}{\pgfqpoint{2.342024in}{2.811081in}}%
\pgfpathlineto{\pgfqpoint{2.342024in}{1.686083in}}%
\pgfpathquadraticcurveto{\pgfqpoint{2.342024in}{1.647194in}}{\pgfqpoint{2.380913in}{1.647194in}}%
\pgfpathlineto{\pgfqpoint{2.380913in}{1.647194in}}%
\pgfpathclose%
\pgfusepath{stroke,fill}%
\end{pgfscope}%
\begin{pgfscope}%
\pgfsetrectcap%
\pgfsetroundjoin%
\pgfsetlinewidth{1.505625pt}%
\definecolor{currentstroke}{rgb}{0.000000,0.000000,0.000000}%
\pgfsetstrokecolor{currentstroke}%
\pgfsetdash{}{0pt}%
\pgfpathmoveto{\pgfqpoint{2.419802in}{2.690248in}}%
\pgfpathlineto{\pgfqpoint{2.614246in}{2.690248in}}%
\pgfpathlineto{\pgfqpoint{2.808691in}{2.690248in}}%
\pgfusepath{stroke}%
\end{pgfscope}%
\begin{pgfscope}%
\pgfsetbuttcap%
\pgfsetroundjoin%
\definecolor{currentfill}{rgb}{0.000000,0.000000,0.000000}%
\pgfsetfillcolor{currentfill}%
\pgfsetlinewidth{1.003750pt}%
\definecolor{currentstroke}{rgb}{0.000000,0.000000,0.000000}%
\pgfsetstrokecolor{currentstroke}%
\pgfsetdash{}{0pt}%
\pgfsys@defobject{currentmarker}{\pgfqpoint{-0.041667in}{-0.041667in}}{\pgfqpoint{0.041667in}{0.041667in}}{%
\pgfpathmoveto{\pgfqpoint{0.000000in}{-0.041667in}}%
\pgfpathcurveto{\pgfqpoint{0.011050in}{-0.041667in}}{\pgfqpoint{0.021649in}{-0.037276in}}{\pgfqpoint{0.029463in}{-0.029463in}}%
\pgfpathcurveto{\pgfqpoint{0.037276in}{-0.021649in}}{\pgfqpoint{0.041667in}{-0.011050in}}{\pgfqpoint{0.041667in}{0.000000in}}%
\pgfpathcurveto{\pgfqpoint{0.041667in}{0.011050in}}{\pgfqpoint{0.037276in}{0.021649in}}{\pgfqpoint{0.029463in}{0.029463in}}%
\pgfpathcurveto{\pgfqpoint{0.021649in}{0.037276in}}{\pgfqpoint{0.011050in}{0.041667in}}{\pgfqpoint{0.000000in}{0.041667in}}%
\pgfpathcurveto{\pgfqpoint{-0.011050in}{0.041667in}}{\pgfqpoint{-0.021649in}{0.037276in}}{\pgfqpoint{-0.029463in}{0.029463in}}%
\pgfpathcurveto{\pgfqpoint{-0.037276in}{0.021649in}}{\pgfqpoint{-0.041667in}{0.011050in}}{\pgfqpoint{-0.041667in}{0.000000in}}%
\pgfpathcurveto{\pgfqpoint{-0.041667in}{-0.011050in}}{\pgfqpoint{-0.037276in}{-0.021649in}}{\pgfqpoint{-0.029463in}{-0.029463in}}%
\pgfpathcurveto{\pgfqpoint{-0.021649in}{-0.037276in}}{\pgfqpoint{-0.011050in}{-0.041667in}}{\pgfqpoint{0.000000in}{-0.041667in}}%
\pgfpathlineto{\pgfqpoint{0.000000in}{-0.041667in}}%
\pgfpathclose%
\pgfusepath{stroke,fill}%
}%
\begin{pgfscope}%
\pgfsys@transformshift{2.614246in}{2.690248in}%
\pgfsys@useobject{currentmarker}{}%
\end{pgfscope}%
\end{pgfscope}%
\begin{pgfscope}%
\definecolor{textcolor}{rgb}{0.000000,0.000000,0.000000}%
\pgfsetstrokecolor{textcolor}%
\pgfsetfillcolor{textcolor}%
\pgftext[x=2.964246in,y=2.622193in,left,base]{\color{textcolor}\rmfamily\fontsize{14.000000}{16.800000}\selectfont Passenger costs (\(\displaystyle \lambda\)-normalized)}%
\end{pgfscope}%
\begin{pgfscope}%
\pgfsetrectcap%
\pgfsetroundjoin%
\pgfsetlinewidth{1.505625pt}%
\definecolor{currentstroke}{rgb}{0.839216,0.152941,0.156863}%
\pgfsetstrokecolor{currentstroke}%
\pgfsetdash{}{0pt}%
\pgfpathmoveto{\pgfqpoint{2.419802in}{2.393026in}}%
\pgfpathlineto{\pgfqpoint{2.614246in}{2.393026in}}%
\pgfpathlineto{\pgfqpoint{2.808691in}{2.393026in}}%
\pgfusepath{stroke}%
\end{pgfscope}%
\begin{pgfscope}%
\pgfsetbuttcap%
\pgfsetroundjoin%
\definecolor{currentfill}{rgb}{0.839216,0.152941,0.156863}%
\pgfsetfillcolor{currentfill}%
\pgfsetlinewidth{1.003750pt}%
\definecolor{currentstroke}{rgb}{0.839216,0.152941,0.156863}%
\pgfsetstrokecolor{currentstroke}%
\pgfsetdash{}{0pt}%
\pgfsys@defobject{currentmarker}{\pgfqpoint{-0.041667in}{-0.041667in}}{\pgfqpoint{0.041667in}{0.041667in}}{%
\pgfpathmoveto{\pgfqpoint{0.000000in}{-0.041667in}}%
\pgfpathcurveto{\pgfqpoint{0.011050in}{-0.041667in}}{\pgfqpoint{0.021649in}{-0.037276in}}{\pgfqpoint{0.029463in}{-0.029463in}}%
\pgfpathcurveto{\pgfqpoint{0.037276in}{-0.021649in}}{\pgfqpoint{0.041667in}{-0.011050in}}{\pgfqpoint{0.041667in}{0.000000in}}%
\pgfpathcurveto{\pgfqpoint{0.041667in}{0.011050in}}{\pgfqpoint{0.037276in}{0.021649in}}{\pgfqpoint{0.029463in}{0.029463in}}%
\pgfpathcurveto{\pgfqpoint{0.021649in}{0.037276in}}{\pgfqpoint{0.011050in}{0.041667in}}{\pgfqpoint{0.000000in}{0.041667in}}%
\pgfpathcurveto{\pgfqpoint{-0.011050in}{0.041667in}}{\pgfqpoint{-0.021649in}{0.037276in}}{\pgfqpoint{-0.029463in}{0.029463in}}%
\pgfpathcurveto{\pgfqpoint{-0.037276in}{0.021649in}}{\pgfqpoint{-0.041667in}{0.011050in}}{\pgfqpoint{-0.041667in}{0.000000in}}%
\pgfpathcurveto{\pgfqpoint{-0.041667in}{-0.011050in}}{\pgfqpoint{-0.037276in}{-0.021649in}}{\pgfqpoint{-0.029463in}{-0.029463in}}%
\pgfpathcurveto{\pgfqpoint{-0.021649in}{-0.037276in}}{\pgfqpoint{-0.011050in}{-0.041667in}}{\pgfqpoint{0.000000in}{-0.041667in}}%
\pgfpathlineto{\pgfqpoint{0.000000in}{-0.041667in}}%
\pgfpathclose%
\pgfusepath{stroke,fill}%
}%
\begin{pgfscope}%
\pgfsys@transformshift{2.614246in}{2.393026in}%
\pgfsys@useobject{currentmarker}{}%
\end{pgfscope}%
\end{pgfscope}%
\begin{pgfscope}%
\definecolor{textcolor}{rgb}{0.000000,0.000000,0.000000}%
\pgfsetstrokecolor{textcolor}%
\pgfsetfillcolor{textcolor}%
\pgftext[x=2.964246in,y=2.324971in,left,base]{\color{textcolor}\rmfamily\fontsize{14.000000}{16.800000}\selectfont Objective (\(\displaystyle \lambda\)-normalized)}%
\end{pgfscope}%
\begin{pgfscope}%
\pgfsetrectcap%
\pgfsetroundjoin%
\pgfsetlinewidth{1.505625pt}%
\definecolor{currentstroke}{rgb}{0.301961,0.301961,0.301961}%
\pgfsetstrokecolor{currentstroke}%
\pgfsetdash{}{0pt}%
\pgfpathmoveto{\pgfqpoint{2.419802in}{2.106915in}}%
\pgfpathlineto{\pgfqpoint{2.614246in}{2.106915in}}%
\pgfpathlineto{\pgfqpoint{2.808691in}{2.106915in}}%
\pgfusepath{stroke}%
\end{pgfscope}%
\begin{pgfscope}%
\pgfsetbuttcap%
\pgfsetroundjoin%
\definecolor{currentfill}{rgb}{0.301961,0.301961,0.301961}%
\pgfsetfillcolor{currentfill}%
\pgfsetlinewidth{1.003750pt}%
\definecolor{currentstroke}{rgb}{0.301961,0.301961,0.301961}%
\pgfsetstrokecolor{currentstroke}%
\pgfsetdash{}{0pt}%
\pgfsys@defobject{currentmarker}{\pgfqpoint{-0.041667in}{-0.041667in}}{\pgfqpoint{0.041667in}{0.041667in}}{%
\pgfpathmoveto{\pgfqpoint{0.000000in}{-0.041667in}}%
\pgfpathcurveto{\pgfqpoint{0.011050in}{-0.041667in}}{\pgfqpoint{0.021649in}{-0.037276in}}{\pgfqpoint{0.029463in}{-0.029463in}}%
\pgfpathcurveto{\pgfqpoint{0.037276in}{-0.021649in}}{\pgfqpoint{0.041667in}{-0.011050in}}{\pgfqpoint{0.041667in}{0.000000in}}%
\pgfpathcurveto{\pgfqpoint{0.041667in}{0.011050in}}{\pgfqpoint{0.037276in}{0.021649in}}{\pgfqpoint{0.029463in}{0.029463in}}%
\pgfpathcurveto{\pgfqpoint{0.021649in}{0.037276in}}{\pgfqpoint{0.011050in}{0.041667in}}{\pgfqpoint{0.000000in}{0.041667in}}%
\pgfpathcurveto{\pgfqpoint{-0.011050in}{0.041667in}}{\pgfqpoint{-0.021649in}{0.037276in}}{\pgfqpoint{-0.029463in}{0.029463in}}%
\pgfpathcurveto{\pgfqpoint{-0.037276in}{0.021649in}}{\pgfqpoint{-0.041667in}{0.011050in}}{\pgfqpoint{-0.041667in}{0.000000in}}%
\pgfpathcurveto{\pgfqpoint{-0.041667in}{-0.011050in}}{\pgfqpoint{-0.037276in}{-0.021649in}}{\pgfqpoint{-0.029463in}{-0.029463in}}%
\pgfpathcurveto{\pgfqpoint{-0.021649in}{-0.037276in}}{\pgfqpoint{-0.011050in}{-0.041667in}}{\pgfqpoint{0.000000in}{-0.041667in}}%
\pgfpathlineto{\pgfqpoint{0.000000in}{-0.041667in}}%
\pgfpathclose%
\pgfusepath{stroke,fill}%
}%
\begin{pgfscope}%
\pgfsys@transformshift{2.614246in}{2.106915in}%
\pgfsys@useobject{currentmarker}{}%
\end{pgfscope}%
\end{pgfscope}%
\begin{pgfscope}%
\definecolor{textcolor}{rgb}{0.000000,0.000000,0.000000}%
\pgfsetstrokecolor{textcolor}%
\pgfsetfillcolor{textcolor}%
\pgftext[x=2.964246in,y=2.038860in,left,base]{\color{textcolor}\rmfamily\fontsize{14.000000}{16.800000}\selectfont Ticket revenue}%
\end{pgfscope}%
\begin{pgfscope}%
\pgfsetrectcap%
\pgfsetroundjoin%
\pgfsetlinewidth{1.505625pt}%
\definecolor{currentstroke}{rgb}{0.678431,0.709804,0.741176}%
\pgfsetstrokecolor{currentstroke}%
\pgfsetdash{}{0pt}%
\pgfpathmoveto{\pgfqpoint{2.419802in}{1.831916in}}%
\pgfpathlineto{\pgfqpoint{2.614246in}{1.831916in}}%
\pgfpathlineto{\pgfqpoint{2.808691in}{1.831916in}}%
\pgfusepath{stroke}%
\end{pgfscope}%
\begin{pgfscope}%
\pgfsetbuttcap%
\pgfsetroundjoin%
\definecolor{currentfill}{rgb}{0.678431,0.709804,0.741176}%
\pgfsetfillcolor{currentfill}%
\pgfsetlinewidth{1.003750pt}%
\definecolor{currentstroke}{rgb}{0.678431,0.709804,0.741176}%
\pgfsetstrokecolor{currentstroke}%
\pgfsetdash{}{0pt}%
\pgfsys@defobject{currentmarker}{\pgfqpoint{-0.041667in}{-0.041667in}}{\pgfqpoint{0.041667in}{0.041667in}}{%
\pgfpathmoveto{\pgfqpoint{0.000000in}{-0.041667in}}%
\pgfpathcurveto{\pgfqpoint{0.011050in}{-0.041667in}}{\pgfqpoint{0.021649in}{-0.037276in}}{\pgfqpoint{0.029463in}{-0.029463in}}%
\pgfpathcurveto{\pgfqpoint{0.037276in}{-0.021649in}}{\pgfqpoint{0.041667in}{-0.011050in}}{\pgfqpoint{0.041667in}{0.000000in}}%
\pgfpathcurveto{\pgfqpoint{0.041667in}{0.011050in}}{\pgfqpoint{0.037276in}{0.021649in}}{\pgfqpoint{0.029463in}{0.029463in}}%
\pgfpathcurveto{\pgfqpoint{0.021649in}{0.037276in}}{\pgfqpoint{0.011050in}{0.041667in}}{\pgfqpoint{0.000000in}{0.041667in}}%
\pgfpathcurveto{\pgfqpoint{-0.011050in}{0.041667in}}{\pgfqpoint{-0.021649in}{0.037276in}}{\pgfqpoint{-0.029463in}{0.029463in}}%
\pgfpathcurveto{\pgfqpoint{-0.037276in}{0.021649in}}{\pgfqpoint{-0.041667in}{0.011050in}}{\pgfqpoint{-0.041667in}{0.000000in}}%
\pgfpathcurveto{\pgfqpoint{-0.041667in}{-0.011050in}}{\pgfqpoint{-0.037276in}{-0.021649in}}{\pgfqpoint{-0.029463in}{-0.029463in}}%
\pgfpathcurveto{\pgfqpoint{-0.021649in}{-0.037276in}}{\pgfqpoint{-0.011050in}{-0.041667in}}{\pgfqpoint{0.000000in}{-0.041667in}}%
\pgfpathlineto{\pgfqpoint{0.000000in}{-0.041667in}}%
\pgfpathclose%
\pgfusepath{stroke,fill}%
}%
\begin{pgfscope}%
\pgfsys@transformshift{2.614246in}{1.831916in}%
\pgfsys@useobject{currentmarker}{}%
\end{pgfscope}%
\end{pgfscope}%
\begin{pgfscope}%
\definecolor{textcolor}{rgb}{0.000000,0.000000,0.000000}%
\pgfsetstrokecolor{textcolor}%
\pgfsetfillcolor{textcolor}%
\pgftext[x=2.964246in,y=1.763860in,left,base]{\color{textcolor}\rmfamily\fontsize{14.000000}{16.800000}\selectfont Operator costs}%
\end{pgfscope}%
\end{pgfpicture}%
\makeatother%
\endgroup%

%% file: Appendix/Plots.tex
\section{Computational performance for different values of \texorpdfstring{$\lambda$ }{lambda}}\label{sec:PlotSup}

\begin{figure}[H]
	\centering
		\subfloat[{\added{$\lambda = 0.1$}}]{\includegraphics[width=0.35\textwidth]{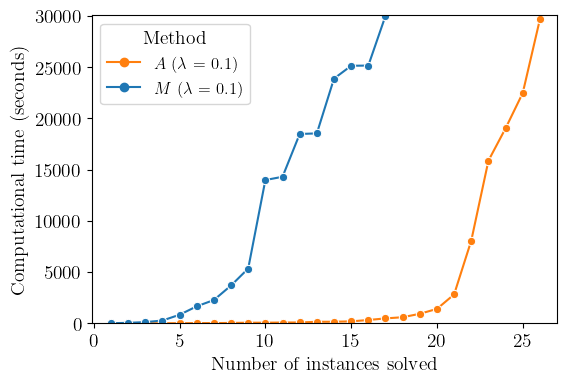}\label{fig:A}}
		\subfloat[{\added{$\lambda = 0.1$}}]{\includegraphics[width=0.35\textwidth]{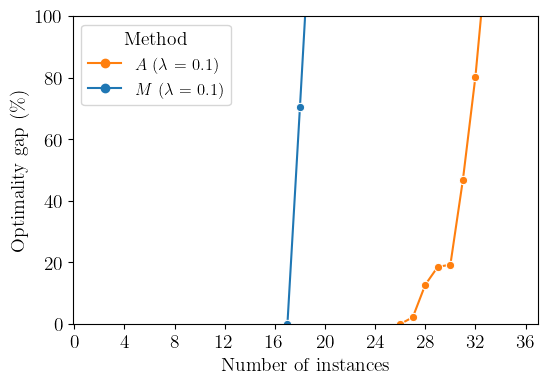}\label{fig:B}} \\
        \subfloat[{\added{$\lambda = 0.5$}}]{\includegraphics[width=0.35\textwidth]{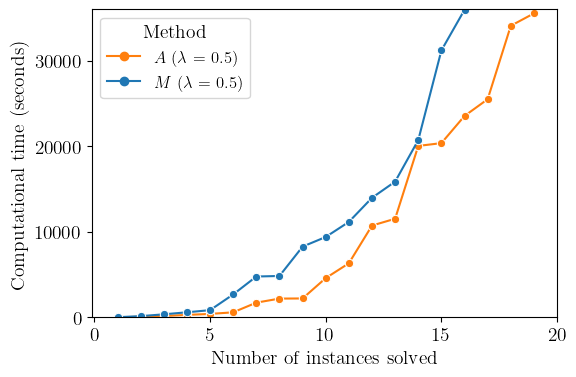}\label{fig:C}}
		\subfloat[{\added{$\lambda = 0.5$}}]{\includegraphics[width=0.35\textwidth]{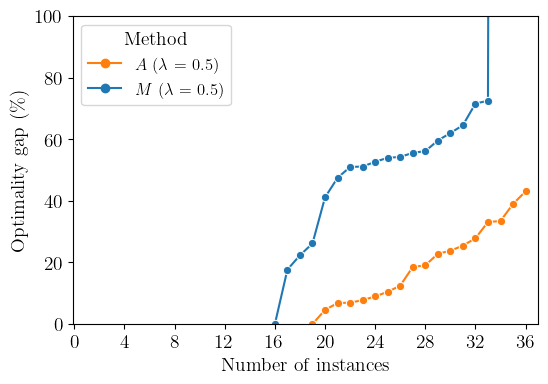}\label{fig:D}} \\
            \subfloat[{\added{$\lambda = 0.75$}}]{\includegraphics[width=0.35\textwidth]{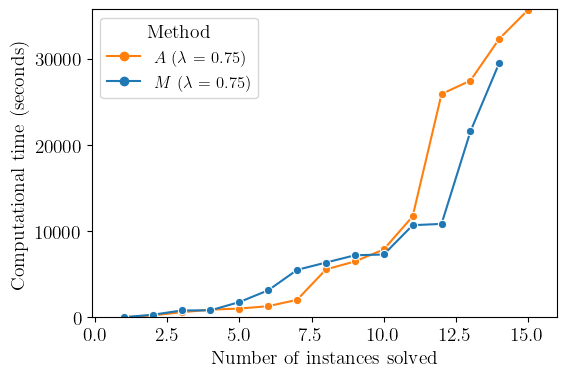}\label{fig:E}}
		\subfloat[{\added{$\lambda = 0.75$}}]{\includegraphics[width=0.35\textwidth]{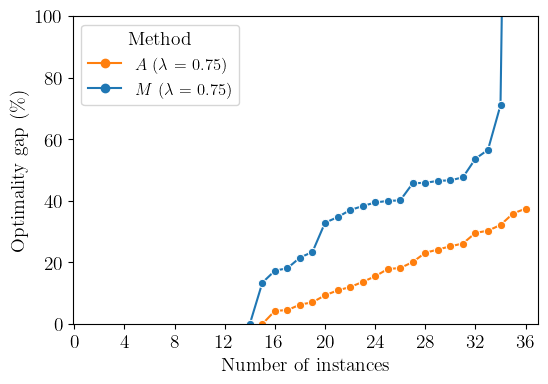}\label{fig:F}} \\
        \subfloat[{\added{$\lambda = 1$}}]{\includegraphics[width=0.35\textwidth]{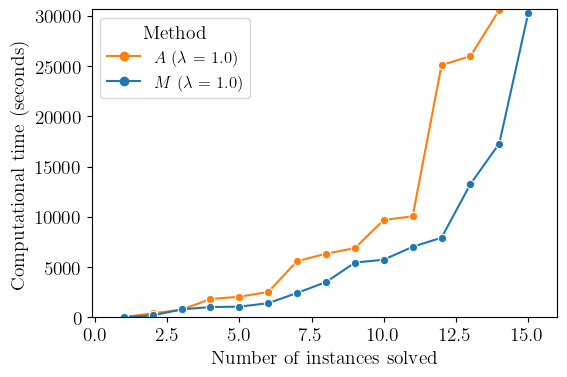}\label{fig:G}}
		\subfloat[{\added{$\lambda = 1$}}]{\includegraphics[width=0.35\textwidth]{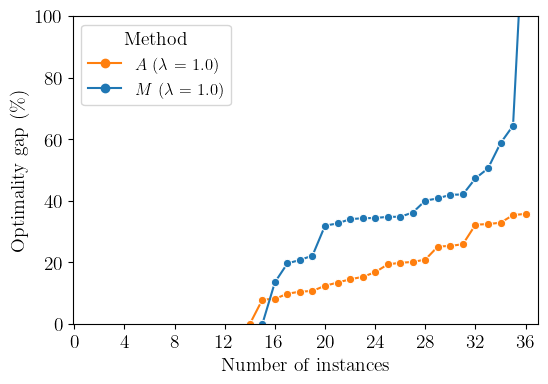}\label{fig:H}}
		\caption{Results using \frameworkNameShortLP ($A$) and CPLEX ($M$). }
\label{fig:AvsM_InstancesSolvedAndGap}
\end{figure}

\autoref{fig:AvsM_InstancesSolvedAndGap} illustrates how many instances are solved to optimality over time, as well as how many achieve solutions within a given optimality gap, for $\lambda \in \{0.1, 0.5, 0.75, 1\}$. We note that while the improvement in computation time is limited for $\lambda \in \{0.5, 0.75, 1\}$, the gaps reported by the \frameworkNameShortLP are stronger even for $\lambda = 1$. 

\section{Computational performance for different problem dimension and values of \texorpdfstring{$\lambda$}{lambda}}\label{sec:problemDims}

The following figures illustrate the number of instances solved to optimality by the \frameworkNameShortLP ($A$) and CPLEX ($M$) grouped by network, $|\ODset|$, and $|\linePool|$. The total number of instances in each group is marked by a dashed line. 

\begin{figure}[H]
    \centering
    \includegraphics[width=1.0\linewidth]{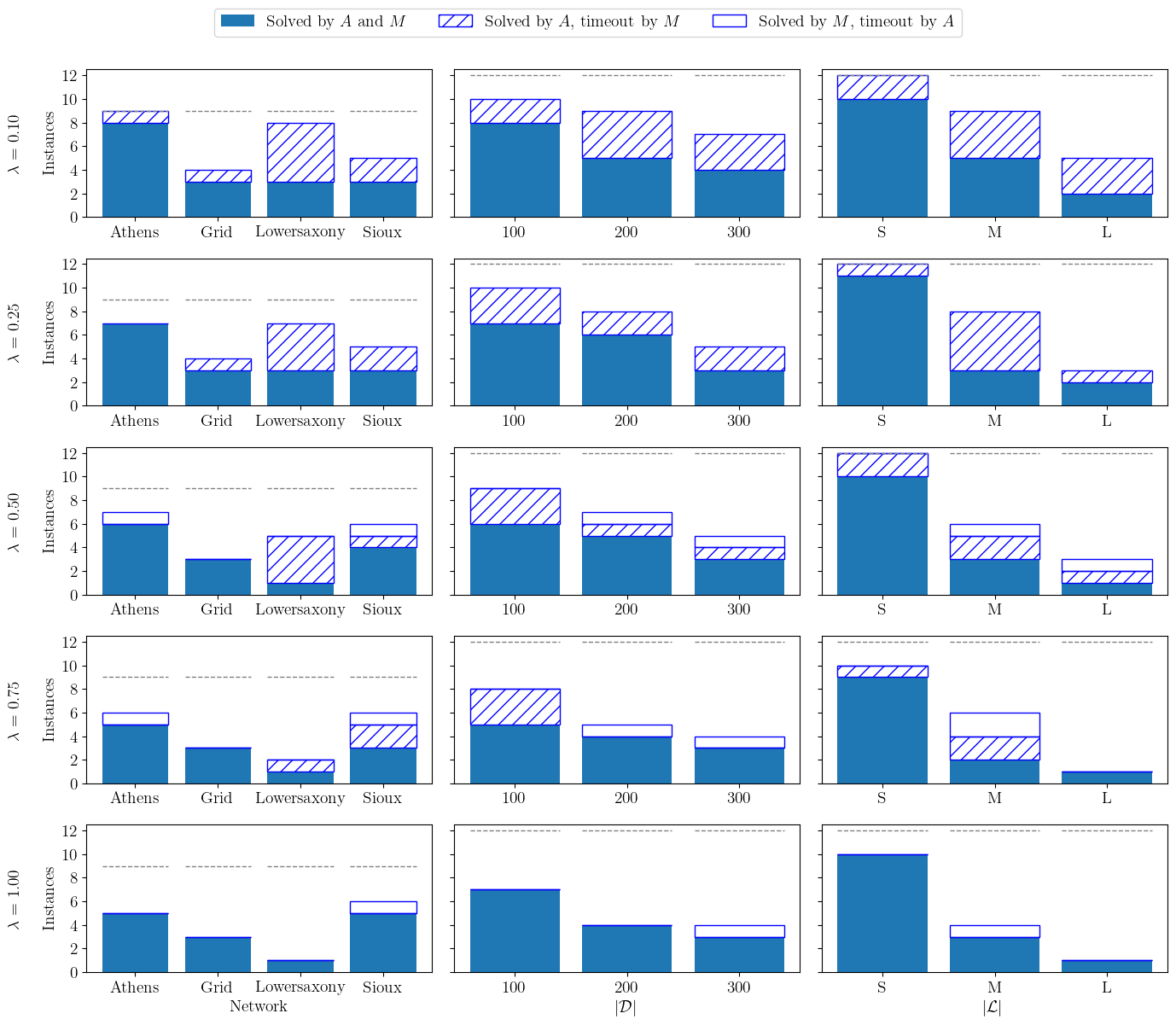}
    \caption{Number of instances solved to optimality}
    \label{fig:DimsALL}
\end{figure}

%% file: Appendix/TableSummary.tex
\section{Aggregated computational results for different values of \texorpdfstring{$\lambda$}{lambda}}\label{sec:AggTables}

\begin{table}[H]
\centering
\caption{Average computational performance for instances solved to optimality by both \frameworkNameShortLP\ (\emph{A}) and CPLEX (\emph{M}), by either, or by both methods. The table reports the number of instances in each group, average computation time, average per-instance speed-up factor, average optimality gap at termination, and the relative improvement in upper bound quality.}
\label{tab:CompAverage01}
\begin{tabular}{@{}lrrrrrrrr@{}}
\toprule
 & & \multicolumn{1}{c}{} & \multicolumn{3}{c}{Time (minutes)} & \multicolumn{2}{c}{Gap (\%)}& \multicolumn{1}{c}{RUB ($\%$)} \\   \cmidrule(l){3-3} \cmidrule(l){4-6} \cmidrule(l){7-8} \cmidrule(l){9-9} 
 &  & \multicolumn{1}{c}{\#Inst} &  \multicolumn{1}{c}{$A$} & \multicolumn{1}{c}{$M$} & \multicolumn{1}{c}{Speed-up} & \multicolumn{1}{c}{$A$} & \multicolumn{1}{c}{$M$} & \multicolumn{1}{c}{$A$} \\  \midrule
& Both optimal & 17  & 22.9 & 180.2 & 192.0 & 0 & 0 & 0 \\
$\lambda $= 0.1 & $A$ optimal  & 9 & 146.7 & $>600$  & 87.8 & 0 &$ >100$ & 19.8 \\
& Neither optimal & 10 & $>600$ &$>600$  & 1.00 & $>100$ & $>100$ & 77.0 \\
\midrule
& Both optimal  & 16  & 61.9 & 165.3 & 33.1 & 0 & 0 & 0 \\
$\lambda $= 0.25& $A$ optimal & 7  & 207.1 & $>600$& 11.8 & 0 & 91.51 & 3.3 \\
& Neither optimal &13  & $>600$ &$>600$ & 1.0 & 28.0 & $>100$ & 20.1 \\ 
\midrule
& Both optimal  & 14 & 148.3 & 111.3 & 11.6 &0 & 0 & 0  \\
$\lambda $= 0.5& $A$ optimal  &5  & 250.3 & $>600$ & 6.5 & 0 & 47.0 & 0.7  \\
&Neither optimal & 15  &$>600$&$>600$ & 1.0 & 22.0 & 661.2 & 11.1  \\
&$M$ optimal  & 2 & $>600$& 559.4 & 0.9 & 6.9 & 0 & -1.8 \\ 
\midrule
& Both optimal  & 12 & 174.5 & 76.0 & 1.2 & 0 & 0 & 0  \\
$\lambda $= 0.75& $A$ optimal  &3  & 186.2 & $>600$ & 8.5 & 0 & 21.5 & $<0.1$  \\
&Neither optimal & 19  &$>600$&$>600$ & 1.0 & 20.5 & 76.0 & 4.3  \\
&$M$ optimal  & 2 & $>600$& 426.0 & 0.7 & 6.9 & 0 & -3.3 \\ 
\midrule
& Both optimal  & 14 & 152.0 & 95.2 & 0.9 & 0 & 0 & 0 \\
$\lambda $= 1.0 & Neither optimal & 21 & $>600$ & $>600$& 1.0 & 20.8 & 41.8 & 2.7 \\
& $M$ optimal & 1 & $>600$ & 288.1 & 0.48 & 7.9 & 0 & -5.6 \\
\bottomrule
\end{tabular}
\end{table}

%% file: Appendix/Demand.tex
\section{Impact of considering \elasticDemand}\label{sec:DE_appendix}

\subsection{Results averaged over a consistent subset of instances solved to optimality}
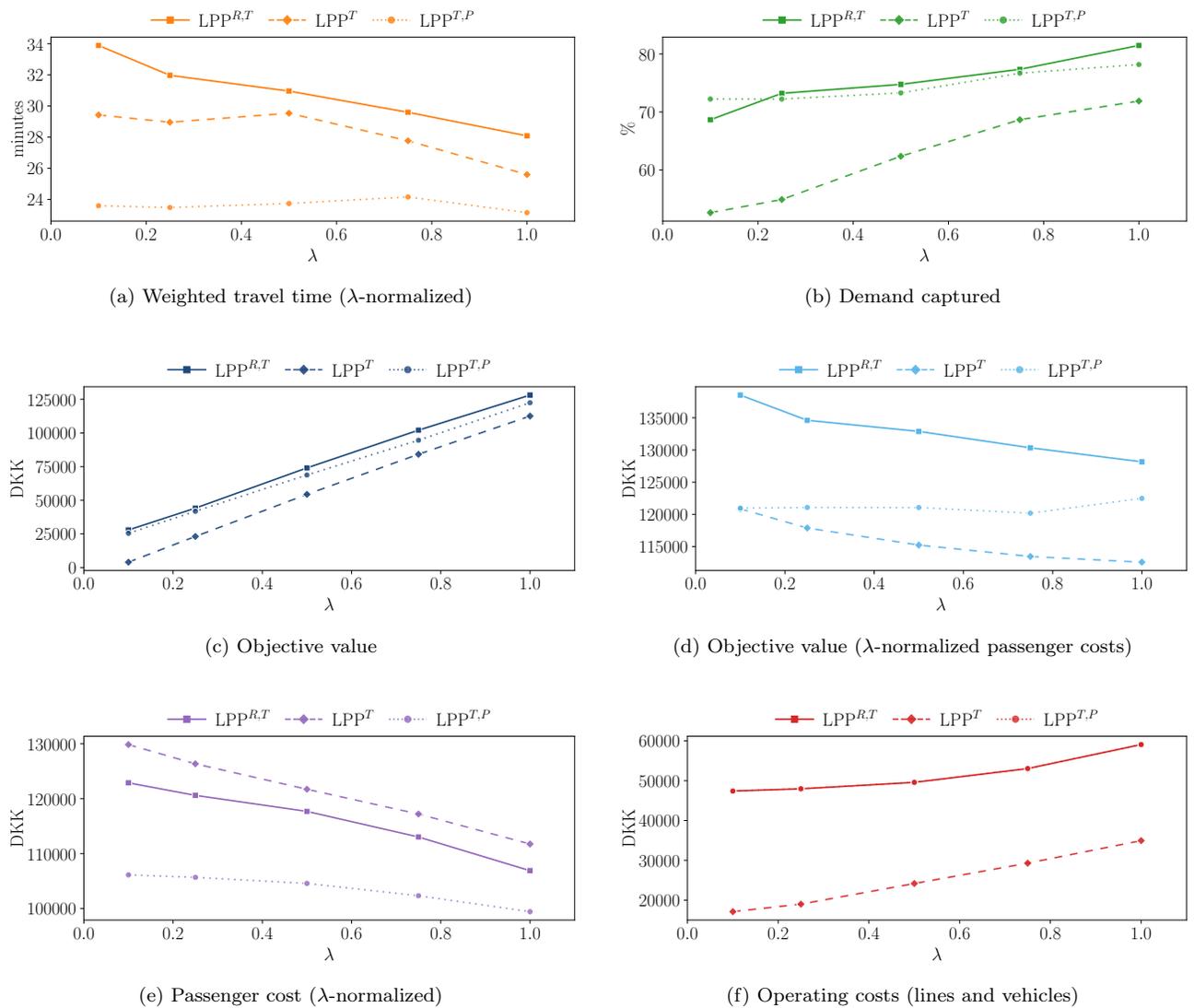
\begin{figure}[H]
  \centering
  \begin{tabular}{cc}
    \begin{subfigure}[b]{0.5\textwidth}
      \centering
      \resizebox{\linewidth}{!}{\input{figures/DemandPlots/TightLayout/DE_avgWTT_rescaled.pgf}}
      \caption{Weighted travel time ($\lambda$-normalized)}
    \end{subfigure} &
    \begin{subfigure}[b]{0.5\textwidth}
      \centering
      \resizebox{\linewidth}{!}{\input{figures/DemandPlots/TightLayout/DE_DemandCapturedPct.pgf}}
      \caption{Demand captured}
    \end{subfigure}  \\[1em]
    \begin{subfigure}[b]{0.5\textwidth}
      \centering
      \resizebox{\linewidth}{!}{\input{figures/DemandPlots/TightLayout/DE_UB.pgf}}
      \caption{Objective value}
    \end{subfigure} & 
        \begin{subfigure}[b]{0.5\textwidth}
      \centering
      \resizebox{\linewidth}{!}{\input{figures/DemandPlots/TightLayout/DE_UB_rescaled.pgf}}
      \caption{Objective value ($\lambda$-normalized passenger costs)}
    \end{subfigure}\\[1em]
    \begin{subfigure}[b]{0.5\textwidth}
      \centering
      \resizebox{\linewidth}{!}{\input{figures/DemandPlots/TightLayout/DE_PI_rescaled.pgf}}
      \caption{Passenger cost ($\lambda$-normalized)}
    \end{subfigure} &
    \begin{subfigure}[b]{0.5\textwidth}
      \centering
      \resizebox{\linewidth}{!}{\input{figures/DemandPlots/TightLayout/DE_OC.pgf}}
      \caption{Operating costs (lines and vehicles)}
    \end{subfigure} 
  \end{tabular}
  \caption{Solution quality for different values of $\lambda$ considering the 7 instances that can be solved to optimality}
  \label{fig:SolQuality}
\end{figure}

\subsection{Results averaged over instances solved within 10\% of optimality}

\begin{figure}[H]
  \centering
  \begin{tabular}{cc}
    \begin{subfigure}[b]{0.5\textwidth}
      \centering
      \resizebox{\linewidth}{!}{\input{figures/DemandPlots/TightLayout/DE_avgWTT_rescaled_Gap10pct.pgf}}
      \caption{Weighted travel time ($\lambda$-normalized)}
    \end{subfigure} &
    \begin{subfigure}[b]{0.5\textwidth}
      \centering
      \resizebox{\linewidth}{!}{\input{figures/DemandPlots/TightLayout/DE_DemandCapturedPct_Gap10pct.pgf}}
      \caption{Demand captured}
    \end{subfigure}  \\[1em]
     \begin{subfigure}[b]{0.5\textwidth}
      \centering
      \resizebox{\linewidth}{!}{\input{figures/DemandPlots/TightLayout/DE_UB_Gap10pct.pgf}}
      \caption{Objective value}
    \end{subfigure} & 
        \begin{subfigure}[b]{0.5\textwidth}
      \centering
      \resizebox{\linewidth}{!}{\input{figures/DemandPlots/TightLayout/DE_UB_rescaled_Gap10pct.pgf}}
      \caption{Objective value ($\lambda$-normalized passenger costs)}
    \end{subfigure}\\[1em]
    \begin{subfigure}[b]{0.5\textwidth}
      \centering
      \resizebox{\linewidth}{!}{\input{figures/DemandPlots/TightLayout/DE_PI_rescaled_Gap10pct.pgf}}
      \caption{Passenger cost ($\lambda$-normalized)}
    \end{subfigure} &
    \begin{subfigure}[b]{0.5\textwidth}
      \centering
      \resizebox{\linewidth}{!}{\input{figures/DemandPlots/TightLayout/DE_OC_Gap10pct.pgf}}
      \caption{Operating costs (lines and vehicles)}
    \end{subfigure} 
  \end{tabular}
  \caption{Solution quality for different values of $\lambda$ considering the 12 instances that can be solved within 10\% of optimality}
  \label{fig:SolQuality10PctGap}
\end{figure}
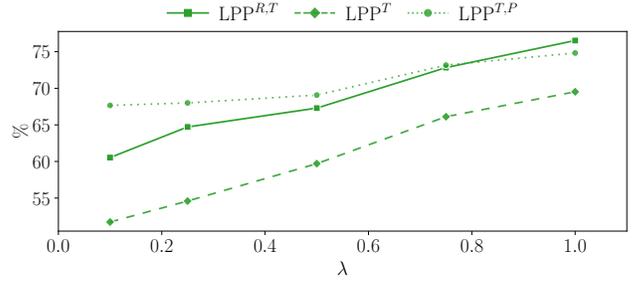
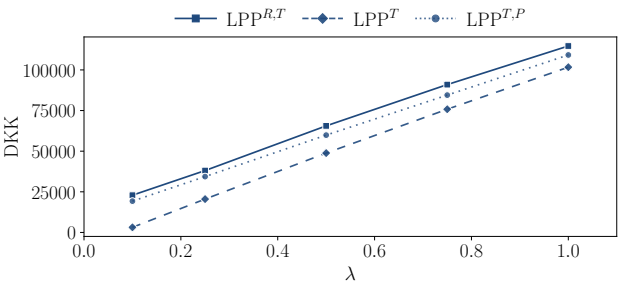
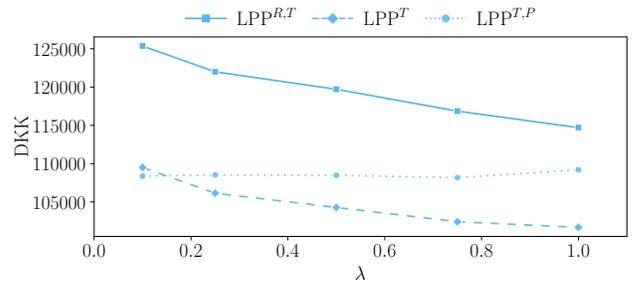
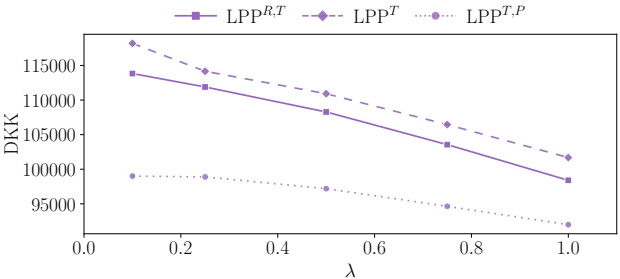
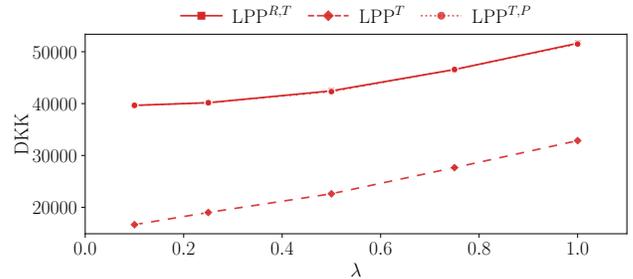

%% file: figures/DemandPlots/TightLayout/DE_UB.pgf
\begingroup%
\makeatletter%
\begin{pgfpicture}%
\pgfpathrectangle{\pgfpointorigin}{\pgfqpoint{7.782018in}{3.651310in}}%
\pgfusepath{use as bounding box, clip}%
\begin{pgfscope}%
\pgfsetbuttcap%
\pgfsetmiterjoin%
\definecolor{currentfill}{rgb}{1.000000,1.000000,1.000000}%
\pgfsetfillcolor{currentfill}%
\pgfsetlinewidth{0.000000pt}%
\definecolor{currentstroke}{rgb}{1.000000,1.000000,1.000000}%
\pgfsetstrokecolor{currentstroke}%
\pgfsetdash{}{0pt}%
\pgfpathmoveto{\pgfqpoint{0.000000in}{0.000000in}}%
\pgfpathlineto{\pgfqpoint{7.782018in}{0.000000in}}%
\pgfpathlineto{\pgfqpoint{7.782018in}{3.651310in}}%
\pgfpathlineto{\pgfqpoint{0.000000in}{3.651310in}}%
\pgfpathlineto{\pgfqpoint{0.000000in}{0.000000in}}%
\pgfpathclose%
\pgfusepath{fill}%
\end{pgfscope}%
\begin{pgfscope}%
\pgfsetbuttcap%
\pgfsetmiterjoin%
\definecolor{currentfill}{rgb}{1.000000,1.000000,1.000000}%
\pgfsetfillcolor{currentfill}%
\pgfsetlinewidth{0.000000pt}%
\definecolor{currentstroke}{rgb}{0.000000,0.000000,0.000000}%
\pgfsetstrokecolor{currentstroke}%
\pgfsetstrokeopacity{0.000000}%
\pgfsetdash{}{0pt}%
\pgfpathmoveto{\pgfqpoint{1.126536in}{0.679475in}}%
\pgfpathlineto{\pgfqpoint{7.682018in}{0.679475in}}%
\pgfpathlineto{\pgfqpoint{7.682018in}{3.135897in}}%
\pgfpathlineto{\pgfqpoint{1.126536in}{3.135897in}}%
\pgfpathlineto{\pgfqpoint{1.126536in}{0.679475in}}%
\pgfpathclose%
\pgfusepath{fill}%
\end{pgfscope}%
\begin{pgfscope}%
\pgfsetbuttcap%
\pgfsetroundjoin%
\definecolor{currentfill}{rgb}{0.000000,0.000000,0.000000}%
\pgfsetfillcolor{currentfill}%
\pgfsetlinewidth{0.803000pt}%
\definecolor{currentstroke}{rgb}{0.000000,0.000000,0.000000}%
\pgfsetstrokecolor{currentstroke}%
\pgfsetdash{}{0pt}%
\pgfsys@defobject{currentmarker}{\pgfqpoint{0.000000in}{-0.048611in}}{\pgfqpoint{0.000000in}{0.000000in}}{%
\pgfpathmoveto{\pgfqpoint{0.000000in}{0.000000in}}%
\pgfpathlineto{\pgfqpoint{0.000000in}{-0.048611in}}%
\pgfusepath{stroke,fill}%
}%
\begin{pgfscope}%
\pgfsys@transformshift{1.126536in}{0.679475in}%
\pgfsys@useobject{currentmarker}{}%
\end{pgfscope}%
\end{pgfscope}%
\begin{pgfscope}%
\definecolor{textcolor}{rgb}{0.000000,0.000000,0.000000}%
\pgfsetstrokecolor{textcolor}%
\pgfsetfillcolor{textcolor}%
\pgftext[x=1.126536in,y=0.582253in,,top]{\color{textcolor}\rmfamily\fontsize{16.000000}{19.200000}\selectfont \(\displaystyle {0.0}\)}%
\end{pgfscope}%
\begin{pgfscope}%
\pgfsetbuttcap%
\pgfsetroundjoin%
\definecolor{currentfill}{rgb}{0.000000,0.000000,0.000000}%
\pgfsetfillcolor{currentfill}%
\pgfsetlinewidth{0.803000pt}%
\definecolor{currentstroke}{rgb}{0.000000,0.000000,0.000000}%
\pgfsetstrokecolor{currentstroke}%
\pgfsetdash{}{0pt}%
\pgfsys@defobject{currentmarker}{\pgfqpoint{0.000000in}{-0.048611in}}{\pgfqpoint{0.000000in}{0.000000in}}{%
\pgfpathmoveto{\pgfqpoint{0.000000in}{0.000000in}}%
\pgfpathlineto{\pgfqpoint{0.000000in}{-0.048611in}}%
\pgfusepath{stroke,fill}%
}%
\begin{pgfscope}%
\pgfsys@transformshift{2.318442in}{0.679475in}%
\pgfsys@useobject{currentmarker}{}%
\end{pgfscope}%
\end{pgfscope}%
\begin{pgfscope}%
\definecolor{textcolor}{rgb}{0.000000,0.000000,0.000000}%
\pgfsetstrokecolor{textcolor}%
\pgfsetfillcolor{textcolor}%
\pgftext[x=2.318442in,y=0.582253in,,top]{\color{textcolor}\rmfamily\fontsize{16.000000}{19.200000}\selectfont \(\displaystyle {0.2}\)}%
\end{pgfscope}%
\begin{pgfscope}%
\pgfsetbuttcap%
\pgfsetroundjoin%
\definecolor{currentfill}{rgb}{0.000000,0.000000,0.000000}%
\pgfsetfillcolor{currentfill}%
\pgfsetlinewidth{0.803000pt}%
\definecolor{currentstroke}{rgb}{0.000000,0.000000,0.000000}%
\pgfsetstrokecolor{currentstroke}%
\pgfsetdash{}{0pt}%
\pgfsys@defobject{currentmarker}{\pgfqpoint{0.000000in}{-0.048611in}}{\pgfqpoint{0.000000in}{0.000000in}}{%
\pgfpathmoveto{\pgfqpoint{0.000000in}{0.000000in}}%
\pgfpathlineto{\pgfqpoint{0.000000in}{-0.048611in}}%
\pgfusepath{stroke,fill}%
}%
\begin{pgfscope}%
\pgfsys@transformshift{3.510347in}{0.679475in}%
\pgfsys@useobject{currentmarker}{}%
\end{pgfscope}%
\end{pgfscope}%
\begin{pgfscope}%
\definecolor{textcolor}{rgb}{0.000000,0.000000,0.000000}%
\pgfsetstrokecolor{textcolor}%
\pgfsetfillcolor{textcolor}%
\pgftext[x=3.510347in,y=0.582253in,,top]{\color{textcolor}\rmfamily\fontsize{16.000000}{19.200000}\selectfont \(\displaystyle {0.4}\)}%
\end{pgfscope}%
\begin{pgfscope}%
\pgfsetbuttcap%
\pgfsetroundjoin%
\definecolor{currentfill}{rgb}{0.000000,0.000000,0.000000}%
\pgfsetfillcolor{currentfill}%
\pgfsetlinewidth{0.803000pt}%
\definecolor{currentstroke}{rgb}{0.000000,0.000000,0.000000}%
\pgfsetstrokecolor{currentstroke}%
\pgfsetdash{}{0pt}%
\pgfsys@defobject{currentmarker}{\pgfqpoint{0.000000in}{-0.048611in}}{\pgfqpoint{0.000000in}{0.000000in}}{%
\pgfpathmoveto{\pgfqpoint{0.000000in}{0.000000in}}%
\pgfpathlineto{\pgfqpoint{0.000000in}{-0.048611in}}%
\pgfusepath{stroke,fill}%
}%
\begin{pgfscope}%
\pgfsys@transformshift{4.702253in}{0.679475in}%
\pgfsys@useobject{currentmarker}{}%
\end{pgfscope}%
\end{pgfscope}%
\begin{pgfscope}%
\definecolor{textcolor}{rgb}{0.000000,0.000000,0.000000}%
\pgfsetstrokecolor{textcolor}%
\pgfsetfillcolor{textcolor}%
\pgftext[x=4.702253in,y=0.582253in,,top]{\color{textcolor}\rmfamily\fontsize{16.000000}{19.200000}\selectfont \(\displaystyle {0.6}\)}%
\end{pgfscope}%
\begin{pgfscope}%
\pgfsetbuttcap%
\pgfsetroundjoin%
\definecolor{currentfill}{rgb}{0.000000,0.000000,0.000000}%
\pgfsetfillcolor{currentfill}%
\pgfsetlinewidth{0.803000pt}%
\definecolor{currentstroke}{rgb}{0.000000,0.000000,0.000000}%
\pgfsetstrokecolor{currentstroke}%
\pgfsetdash{}{0pt}%
\pgfsys@defobject{currentmarker}{\pgfqpoint{0.000000in}{-0.048611in}}{\pgfqpoint{0.000000in}{0.000000in}}{%
\pgfpathmoveto{\pgfqpoint{0.000000in}{0.000000in}}%
\pgfpathlineto{\pgfqpoint{0.000000in}{-0.048611in}}%
\pgfusepath{stroke,fill}%
}%
\begin{pgfscope}%
\pgfsys@transformshift{5.894159in}{0.679475in}%
\pgfsys@useobject{currentmarker}{}%
\end{pgfscope}%
\end{pgfscope}%
\begin{pgfscope}%
\definecolor{textcolor}{rgb}{0.000000,0.000000,0.000000}%
\pgfsetstrokecolor{textcolor}%
\pgfsetfillcolor{textcolor}%
\pgftext[x=5.894159in,y=0.582253in,,top]{\color{textcolor}\rmfamily\fontsize{16.000000}{19.200000}\selectfont \(\displaystyle {0.8}\)}%
\end{pgfscope}%
\begin{pgfscope}%
\pgfsetbuttcap%
\pgfsetroundjoin%
\definecolor{currentfill}{rgb}{0.000000,0.000000,0.000000}%
\pgfsetfillcolor{currentfill}%
\pgfsetlinewidth{0.803000pt}%
\definecolor{currentstroke}{rgb}{0.000000,0.000000,0.000000}%
\pgfsetstrokecolor{currentstroke}%
\pgfsetdash{}{0pt}%
\pgfsys@defobject{currentmarker}{\pgfqpoint{0.000000in}{-0.048611in}}{\pgfqpoint{0.000000in}{0.000000in}}{%
\pgfpathmoveto{\pgfqpoint{0.000000in}{0.000000in}}%
\pgfpathlineto{\pgfqpoint{0.000000in}{-0.048611in}}%
\pgfusepath{stroke,fill}%
}%
\begin{pgfscope}%
\pgfsys@transformshift{7.086065in}{0.679475in}%
\pgfsys@useobject{currentmarker}{}%
\end{pgfscope}%
\end{pgfscope}%
\begin{pgfscope}%
\definecolor{textcolor}{rgb}{0.000000,0.000000,0.000000}%
\pgfsetstrokecolor{textcolor}%
\pgfsetfillcolor{textcolor}%
\pgftext[x=7.086065in,y=0.582253in,,top]{\color{textcolor}\rmfamily\fontsize{16.000000}{19.200000}\selectfont \(\displaystyle {1.0}\)}%
\end{pgfscope}%
\begin{pgfscope}%
\definecolor{textcolor}{rgb}{0.000000,0.000000,0.000000}%
\pgfsetstrokecolor{textcolor}%
\pgfsetfillcolor{textcolor}%
\pgftext[x=4.404277in,y=0.313349in,,top]{\color{textcolor}\rmfamily\fontsize{18.000000}{21.600000}\selectfont \(\displaystyle \lambda\)}%
\end{pgfscope}%
\begin{pgfscope}%
\pgfsetbuttcap%
\pgfsetroundjoin%
\definecolor{currentfill}{rgb}{0.000000,0.000000,0.000000}%
\pgfsetfillcolor{currentfill}%
\pgfsetlinewidth{0.803000pt}%
\definecolor{currentstroke}{rgb}{0.000000,0.000000,0.000000}%
\pgfsetstrokecolor{currentstroke}%
\pgfsetdash{}{0pt}%
\pgfsys@defobject{currentmarker}{\pgfqpoint{-0.048611in}{0.000000in}}{\pgfqpoint{-0.000000in}{0.000000in}}{%
\pgfpathmoveto{\pgfqpoint{-0.000000in}{0.000000in}}%
\pgfpathlineto{\pgfqpoint{-0.048611in}{0.000000in}}%
\pgfusepath{stroke,fill}%
}%
\begin{pgfscope}%
\pgfsys@transformshift{1.126536in}{0.719599in}%
\pgfsys@useobject{currentmarker}{}%
\end{pgfscope}%
\end{pgfscope}%
\begin{pgfscope}%
\definecolor{textcolor}{rgb}{0.000000,0.000000,0.000000}%
\pgfsetstrokecolor{textcolor}%
\pgfsetfillcolor{textcolor}%
\pgftext[x=0.919245in, y=0.636266in, left, base]{\color{textcolor}\rmfamily\fontsize{16.000000}{19.200000}\selectfont \(\displaystyle {0}\)}%
\end{pgfscope}%
\begin{pgfscope}%
\pgfsetbuttcap%
\pgfsetroundjoin%
\definecolor{currentfill}{rgb}{0.000000,0.000000,0.000000}%
\pgfsetfillcolor{currentfill}%
\pgfsetlinewidth{0.803000pt}%
\definecolor{currentstroke}{rgb}{0.000000,0.000000,0.000000}%
\pgfsetstrokecolor{currentstroke}%
\pgfsetdash{}{0pt}%
\pgfsys@defobject{currentmarker}{\pgfqpoint{-0.048611in}{0.000000in}}{\pgfqpoint{-0.000000in}{0.000000in}}{%
\pgfpathmoveto{\pgfqpoint{-0.000000in}{0.000000in}}%
\pgfpathlineto{\pgfqpoint{-0.048611in}{0.000000in}}%
\pgfusepath{stroke,fill}%
}%
\begin{pgfscope}%
\pgfsys@transformshift{1.126536in}{1.169175in}%
\pgfsys@useobject{currentmarker}{}%
\end{pgfscope}%
\end{pgfscope}%
\begin{pgfscope}%
\definecolor{textcolor}{rgb}{0.000000,0.000000,0.000000}%
\pgfsetstrokecolor{textcolor}%
\pgfsetfillcolor{textcolor}%
\pgftext[x=0.478972in, y=1.085842in, left, base]{\color{textcolor}\rmfamily\fontsize{16.000000}{19.200000}\selectfont \(\displaystyle {25000}\)}%
\end{pgfscope}%
\begin{pgfscope}%
\pgfsetbuttcap%
\pgfsetroundjoin%
\definecolor{currentfill}{rgb}{0.000000,0.000000,0.000000}%
\pgfsetfillcolor{currentfill}%
\pgfsetlinewidth{0.803000pt}%
\definecolor{currentstroke}{rgb}{0.000000,0.000000,0.000000}%
\pgfsetstrokecolor{currentstroke}%
\pgfsetdash{}{0pt}%
\pgfsys@defobject{currentmarker}{\pgfqpoint{-0.048611in}{0.000000in}}{\pgfqpoint{-0.000000in}{0.000000in}}{%
\pgfpathmoveto{\pgfqpoint{-0.000000in}{0.000000in}}%
\pgfpathlineto{\pgfqpoint{-0.048611in}{0.000000in}}%
\pgfusepath{stroke,fill}%
}%
\begin{pgfscope}%
\pgfsys@transformshift{1.126536in}{1.618751in}%
\pgfsys@useobject{currentmarker}{}%
\end{pgfscope}%
\end{pgfscope}%
\begin{pgfscope}%
\definecolor{textcolor}{rgb}{0.000000,0.000000,0.000000}%
\pgfsetstrokecolor{textcolor}%
\pgfsetfillcolor{textcolor}%
\pgftext[x=0.478972in, y=1.535418in, left, base]{\color{textcolor}\rmfamily\fontsize{16.000000}{19.200000}\selectfont \(\displaystyle {50000}\)}%
\end{pgfscope}%
\begin{pgfscope}%
\pgfsetbuttcap%
\pgfsetroundjoin%
\definecolor{currentfill}{rgb}{0.000000,0.000000,0.000000}%
\pgfsetfillcolor{currentfill}%
\pgfsetlinewidth{0.803000pt}%
\definecolor{currentstroke}{rgb}{0.000000,0.000000,0.000000}%
\pgfsetstrokecolor{currentstroke}%
\pgfsetdash{}{0pt}%
\pgfsys@defobject{currentmarker}{\pgfqpoint{-0.048611in}{0.000000in}}{\pgfqpoint{-0.000000in}{0.000000in}}{%
\pgfpathmoveto{\pgfqpoint{-0.000000in}{0.000000in}}%
\pgfpathlineto{\pgfqpoint{-0.048611in}{0.000000in}}%
\pgfusepath{stroke,fill}%
}%
\begin{pgfscope}%
\pgfsys@transformshift{1.126536in}{2.068327in}%
\pgfsys@useobject{currentmarker}{}%
\end{pgfscope}%
\end{pgfscope}%
\begin{pgfscope}%
\definecolor{textcolor}{rgb}{0.000000,0.000000,0.000000}%
\pgfsetstrokecolor{textcolor}%
\pgfsetfillcolor{textcolor}%
\pgftext[x=0.478972in, y=1.984994in, left, base]{\color{textcolor}\rmfamily\fontsize{16.000000}{19.200000}\selectfont \(\displaystyle {75000}\)}%
\end{pgfscope}%
\begin{pgfscope}%
\pgfsetbuttcap%
\pgfsetroundjoin%
\definecolor{currentfill}{rgb}{0.000000,0.000000,0.000000}%
\pgfsetfillcolor{currentfill}%
\pgfsetlinewidth{0.803000pt}%
\definecolor{currentstroke}{rgb}{0.000000,0.000000,0.000000}%
\pgfsetstrokecolor{currentstroke}%
\pgfsetdash{}{0pt}%
\pgfsys@defobject{currentmarker}{\pgfqpoint{-0.048611in}{0.000000in}}{\pgfqpoint{-0.000000in}{0.000000in}}{%
\pgfpathmoveto{\pgfqpoint{-0.000000in}{0.000000in}}%
\pgfpathlineto{\pgfqpoint{-0.048611in}{0.000000in}}%
\pgfusepath{stroke,fill}%
}%
\begin{pgfscope}%
\pgfsys@transformshift{1.126536in}{2.517903in}%
\pgfsys@useobject{currentmarker}{}%
\end{pgfscope}%
\end{pgfscope}%
\begin{pgfscope}%
\definecolor{textcolor}{rgb}{0.000000,0.000000,0.000000}%
\pgfsetstrokecolor{textcolor}%
\pgfsetfillcolor{textcolor}%
\pgftext[x=0.368904in, y=2.434570in, left, base]{\color{textcolor}\rmfamily\fontsize{16.000000}{19.200000}\selectfont \(\displaystyle {100000}\)}%
\end{pgfscope}%
\begin{pgfscope}%
\pgfsetbuttcap%
\pgfsetroundjoin%
\definecolor{currentfill}{rgb}{0.000000,0.000000,0.000000}%
\pgfsetfillcolor{currentfill}%
\pgfsetlinewidth{0.803000pt}%
\definecolor{currentstroke}{rgb}{0.000000,0.000000,0.000000}%
\pgfsetstrokecolor{currentstroke}%
\pgfsetdash{}{0pt}%
\pgfsys@defobject{currentmarker}{\pgfqpoint{-0.048611in}{0.000000in}}{\pgfqpoint{-0.000000in}{0.000000in}}{%
\pgfpathmoveto{\pgfqpoint{-0.000000in}{0.000000in}}%
\pgfpathlineto{\pgfqpoint{-0.048611in}{0.000000in}}%
\pgfusepath{stroke,fill}%
}%
\begin{pgfscope}%
\pgfsys@transformshift{1.126536in}{2.967479in}%
\pgfsys@useobject{currentmarker}{}%
\end{pgfscope}%
\end{pgfscope}%
\begin{pgfscope}%
\definecolor{textcolor}{rgb}{0.000000,0.000000,0.000000}%
\pgfsetstrokecolor{textcolor}%
\pgfsetfillcolor{textcolor}%
\pgftext[x=0.368904in, y=2.884146in, left, base]{\color{textcolor}\rmfamily\fontsize{16.000000}{19.200000}\selectfont \(\displaystyle {125000}\)}%
\end{pgfscope}%
\begin{pgfscope}%
\definecolor{textcolor}{rgb}{0.000000,0.000000,0.000000}%
\pgfsetstrokecolor{textcolor}%
\pgfsetfillcolor{textcolor}%
\pgftext[x=0.313349in,y=1.907686in,,bottom,rotate=90.000000]{\color{textcolor}\rmfamily\fontsize{18.000000}{21.600000}\selectfont DKK}%
\end{pgfscope}%
\begin{pgfscope}%
\pgfpathrectangle{\pgfqpoint{1.126536in}{0.679475in}}{\pgfqpoint{6.555482in}{2.456422in}}%
\pgfusepath{clip}%
\pgfsetrectcap%
\pgfsetroundjoin%
\pgfsetlinewidth{1.505625pt}%
\definecolor{currentstroke}{rgb}{0.121569,0.286275,0.490196}%
\pgfsetstrokecolor{currentstroke}%
\pgfsetdash{}{0pt}%
\pgfpathmoveto{\pgfqpoint{1.722489in}{1.221637in}}%
\pgfpathlineto{\pgfqpoint{2.616418in}{1.513499in}}%
\pgfpathlineto{\pgfqpoint{4.106300in}{2.050981in}}%
\pgfpathlineto{\pgfqpoint{5.596183in}{2.555241in}}%
\pgfpathlineto{\pgfqpoint{7.086065in}{3.024241in}}%
\pgfusepath{stroke}%
\end{pgfscope}%
\begin{pgfscope}%
\pgfpathrectangle{\pgfqpoint{1.126536in}{0.679475in}}{\pgfqpoint{6.555482in}{2.456422in}}%
\pgfusepath{clip}%
\pgfsetbuttcap%
\pgfsetmiterjoin%
\definecolor{currentfill}{rgb}{0.121569,0.286275,0.490196}%
\pgfsetfillcolor{currentfill}%
\pgfsetlinewidth{0.752812pt}%
\definecolor{currentstroke}{rgb}{1.000000,1.000000,1.000000}%
\pgfsetstrokecolor{currentstroke}%
\pgfsetdash{}{0pt}%
\pgfsys@defobject{currentmarker}{\pgfqpoint{-0.041667in}{-0.041667in}}{\pgfqpoint{0.041667in}{0.041667in}}{%
\pgfpathmoveto{\pgfqpoint{-0.041667in}{-0.041667in}}%
\pgfpathlineto{\pgfqpoint{0.041667in}{-0.041667in}}%
\pgfpathlineto{\pgfqpoint{0.041667in}{0.041667in}}%
\pgfpathlineto{\pgfqpoint{-0.041667in}{0.041667in}}%
\pgfpathlineto{\pgfqpoint{-0.041667in}{-0.041667in}}%
\pgfpathclose%
\pgfusepath{stroke,fill}%
}%
\begin{pgfscope}%
\pgfsys@transformshift{1.722489in}{1.221637in}%
\pgfsys@useobject{currentmarker}{}%
\end{pgfscope}%
\begin{pgfscope}%
\pgfsys@transformshift{2.616418in}{1.513499in}%
\pgfsys@useobject{currentmarker}{}%
\end{pgfscope}%
\begin{pgfscope}%
\pgfsys@transformshift{4.106300in}{2.050981in}%
\pgfsys@useobject{currentmarker}{}%
\end{pgfscope}%
\begin{pgfscope}%
\pgfsys@transformshift{5.596183in}{2.555241in}%
\pgfsys@useobject{currentmarker}{}%
\end{pgfscope}%
\begin{pgfscope}%
\pgfsys@transformshift{7.086065in}{3.024241in}%
\pgfsys@useobject{currentmarker}{}%
\end{pgfscope}%
\end{pgfscope}%
\begin{pgfscope}%
\pgfpathrectangle{\pgfqpoint{1.126536in}{0.679475in}}{\pgfqpoint{6.555482in}{2.456422in}}%
\pgfusepath{clip}%
\pgfsetbuttcap%
\pgfsetroundjoin%
\pgfsetlinewidth{1.505625pt}%
\definecolor{currentstroke}{rgb}{0.121569,0.286275,0.490196}%
\pgfsetstrokecolor{currentstroke}%
\pgfsetstrokeopacity{0.900000}%
\pgfsetdash{{7.500000pt}{7.500000pt}}{0.000000pt}%
\pgfpathmoveto{\pgfqpoint{1.722489in}{0.791131in}}%
\pgfpathlineto{\pgfqpoint{2.616418in}{1.135154in}}%
\pgfpathlineto{\pgfqpoint{4.106300in}{1.697211in}}%
\pgfpathlineto{\pgfqpoint{5.596183in}{2.232875in}}%
\pgfpathlineto{\pgfqpoint{7.086065in}{2.743922in}}%
\pgfusepath{stroke}%
\end{pgfscope}%
\begin{pgfscope}%
\pgfpathrectangle{\pgfqpoint{1.126536in}{0.679475in}}{\pgfqpoint{6.555482in}{2.456422in}}%
\pgfusepath{clip}%
\pgfsetbuttcap%
\pgfsetmiterjoin%
\definecolor{currentfill}{rgb}{0.121569,0.286275,0.490196}%
\pgfsetfillcolor{currentfill}%
\pgfsetfillopacity{0.900000}%
\pgfsetlinewidth{0.752812pt}%
\definecolor{currentstroke}{rgb}{1.000000,1.000000,1.000000}%
\pgfsetstrokecolor{currentstroke}%
\pgfsetdash{}{0pt}%
\pgfsys@defobject{currentmarker}{\pgfqpoint{-0.058926in}{-0.058926in}}{\pgfqpoint{0.058926in}{0.058926in}}{%
\pgfpathmoveto{\pgfqpoint{0.000000in}{-0.058926in}}%
\pgfpathlineto{\pgfqpoint{0.058926in}{0.000000in}}%
\pgfpathlineto{\pgfqpoint{0.000000in}{0.058926in}}%
\pgfpathlineto{\pgfqpoint{-0.058926in}{0.000000in}}%
\pgfpathlineto{\pgfqpoint{0.000000in}{-0.058926in}}%
\pgfpathclose%
\pgfusepath{stroke,fill}%
}%
\begin{pgfscope}%
\pgfsys@transformshift{1.722489in}{0.791131in}%
\pgfsys@useobject{currentmarker}{}%
\end{pgfscope}%
\begin{pgfscope}%
\pgfsys@transformshift{2.616418in}{1.135154in}%
\pgfsys@useobject{currentmarker}{}%
\end{pgfscope}%
\begin{pgfscope}%
\pgfsys@transformshift{4.106300in}{1.697211in}%
\pgfsys@useobject{currentmarker}{}%
\end{pgfscope}%
\begin{pgfscope}%
\pgfsys@transformshift{5.596183in}{2.232875in}%
\pgfsys@useobject{currentmarker}{}%
\end{pgfscope}%
\begin{pgfscope}%
\pgfsys@transformshift{7.086065in}{2.743922in}%
\pgfsys@useobject{currentmarker}{}%
\end{pgfscope}%
\end{pgfscope}%
\begin{pgfscope}%
\pgfpathrectangle{\pgfqpoint{1.126536in}{0.679475in}}{\pgfqpoint{6.555482in}{2.456422in}}%
\pgfusepath{clip}%
\pgfsetbuttcap%
\pgfsetroundjoin%
\pgfsetlinewidth{1.505625pt}%
\definecolor{currentstroke}{rgb}{0.121569,0.286275,0.490196}%
\pgfsetstrokecolor{currentstroke}%
\pgfsetstrokeopacity{0.800000}%
\pgfsetdash{{1.500000pt}{4.500000pt}}{0.000000pt}%
\pgfpathmoveto{\pgfqpoint{1.722489in}{1.176944in}}%
\pgfpathlineto{\pgfqpoint{2.616418in}{1.471259in}}%
\pgfpathlineto{\pgfqpoint{4.106300in}{1.956118in}}%
\pgfpathlineto{\pgfqpoint{5.596183in}{2.420764in}}%
\pgfpathlineto{\pgfqpoint{7.086065in}{2.922123in}}%
\pgfusepath{stroke}%
\end{pgfscope}%
\begin{pgfscope}%
\pgfpathrectangle{\pgfqpoint{1.126536in}{0.679475in}}{\pgfqpoint{6.555482in}{2.456422in}}%
\pgfusepath{clip}%
\pgfsetbuttcap%
\pgfsetroundjoin%
\definecolor{currentfill}{rgb}{0.121569,0.286275,0.490196}%
\pgfsetfillcolor{currentfill}%
\pgfsetfillopacity{0.800000}%
\pgfsetlinewidth{0.752812pt}%
\definecolor{currentstroke}{rgb}{1.000000,1.000000,1.000000}%
\pgfsetstrokecolor{currentstroke}%
\pgfsetdash{}{0pt}%
\pgfsys@defobject{currentmarker}{\pgfqpoint{-0.041667in}{-0.041667in}}{\pgfqpoint{0.041667in}{0.041667in}}{%
\pgfpathmoveto{\pgfqpoint{0.000000in}{-0.041667in}}%
\pgfpathcurveto{\pgfqpoint{0.011050in}{-0.041667in}}{\pgfqpoint{0.021649in}{-0.037276in}}{\pgfqpoint{0.029463in}{-0.029463in}}%
\pgfpathcurveto{\pgfqpoint{0.037276in}{-0.021649in}}{\pgfqpoint{0.041667in}{-0.011050in}}{\pgfqpoint{0.041667in}{0.000000in}}%
\pgfpathcurveto{\pgfqpoint{0.041667in}{0.011050in}}{\pgfqpoint{0.037276in}{0.021649in}}{\pgfqpoint{0.029463in}{0.029463in}}%
\pgfpathcurveto{\pgfqpoint{0.021649in}{0.037276in}}{\pgfqpoint{0.011050in}{0.041667in}}{\pgfqpoint{0.000000in}{0.041667in}}%
\pgfpathcurveto{\pgfqpoint{-0.011050in}{0.041667in}}{\pgfqpoint{-0.021649in}{0.037276in}}{\pgfqpoint{-0.029463in}{0.029463in}}%
\pgfpathcurveto{\pgfqpoint{-0.037276in}{0.021649in}}{\pgfqpoint{-0.041667in}{0.011050in}}{\pgfqpoint{-0.041667in}{0.000000in}}%
\pgfpathcurveto{\pgfqpoint{-0.041667in}{-0.011050in}}{\pgfqpoint{-0.037276in}{-0.021649in}}{\pgfqpoint{-0.029463in}{-0.029463in}}%
\pgfpathcurveto{\pgfqpoint{-0.021649in}{-0.037276in}}{\pgfqpoint{-0.011050in}{-0.041667in}}{\pgfqpoint{0.000000in}{-0.041667in}}%
\pgfpathlineto{\pgfqpoint{0.000000in}{-0.041667in}}%
\pgfpathclose%
\pgfusepath{stroke,fill}%
}%
\begin{pgfscope}%
\pgfsys@transformshift{1.722489in}{1.176944in}%
\pgfsys@useobject{currentmarker}{}%
\end{pgfscope}%
\begin{pgfscope}%
\pgfsys@transformshift{2.616418in}{1.471259in}%
\pgfsys@useobject{currentmarker}{}%
\end{pgfscope}%
\begin{pgfscope}%
\pgfsys@transformshift{4.106300in}{1.956118in}%
\pgfsys@useobject{currentmarker}{}%
\end{pgfscope}%
\begin{pgfscope}%
\pgfsys@transformshift{5.596183in}{2.420764in}%
\pgfsys@useobject{currentmarker}{}%
\end{pgfscope}%
\begin{pgfscope}%
\pgfsys@transformshift{7.086065in}{2.922123in}%
\pgfsys@useobject{currentmarker}{}%
\end{pgfscope}%
\end{pgfscope}%
\begin{pgfscope}%
\pgfsetrectcap%
\pgfsetmiterjoin%
\pgfsetlinewidth{0.803000pt}%
\definecolor{currentstroke}{rgb}{0.000000,0.000000,0.000000}%
\pgfsetstrokecolor{currentstroke}%
\pgfsetdash{}{0pt}%
\pgfpathmoveto{\pgfqpoint{1.126536in}{0.679475in}}%
\pgfpathlineto{\pgfqpoint{1.126536in}{3.135897in}}%
\pgfusepath{stroke}%
\end{pgfscope}%
\begin{pgfscope}%
\pgfsetrectcap%
\pgfsetmiterjoin%
\pgfsetlinewidth{0.803000pt}%
\definecolor{currentstroke}{rgb}{0.000000,0.000000,0.000000}%
\pgfsetstrokecolor{currentstroke}%
\pgfsetdash{}{0pt}%
\pgfpathmoveto{\pgfqpoint{7.682018in}{0.679475in}}%
\pgfpathlineto{\pgfqpoint{7.682018in}{3.135897in}}%
\pgfusepath{stroke}%
\end{pgfscope}%
\begin{pgfscope}%
\pgfsetrectcap%
\pgfsetmiterjoin%
\pgfsetlinewidth{0.803000pt}%
\definecolor{currentstroke}{rgb}{0.000000,0.000000,0.000000}%
\pgfsetstrokecolor{currentstroke}%
\pgfsetdash{}{0pt}%
\pgfpathmoveto{\pgfqpoint{1.126536in}{0.679475in}}%
\pgfpathlineto{\pgfqpoint{7.682018in}{0.679475in}}%
\pgfusepath{stroke}%
\end{pgfscope}%
\begin{pgfscope}%
\pgfsetrectcap%
\pgfsetmiterjoin%
\pgfsetlinewidth{0.803000pt}%
\definecolor{currentstroke}{rgb}{0.000000,0.000000,0.000000}%
\pgfsetstrokecolor{currentstroke}%
\pgfsetdash{}{0pt}%
\pgfpathmoveto{\pgfqpoint{1.126536in}{3.135897in}}%
\pgfpathlineto{\pgfqpoint{7.682018in}{3.135897in}}%
\pgfusepath{stroke}%
\end{pgfscope}%
\begin{pgfscope}%
\pgfsetrectcap%
\pgfsetroundjoin%
\pgfsetlinewidth{1.505625pt}%
\definecolor{currentstroke}{rgb}{0.121569,0.286275,0.490196}%
\pgfsetstrokecolor{currentstroke}%
\pgfsetdash{}{0pt}%
\pgfpathmoveto{\pgfqpoint{2.251710in}{3.363810in}}%
\pgfpathlineto{\pgfqpoint{2.501710in}{3.363810in}}%
\pgfpathlineto{\pgfqpoint{2.751710in}{3.363810in}}%
\pgfusepath{stroke}%
\end{pgfscope}%
\begin{pgfscope}%
\pgfsetbuttcap%
\pgfsetmiterjoin%
\definecolor{currentfill}{rgb}{0.121569,0.286275,0.490196}%
\pgfsetfillcolor{currentfill}%
\pgfsetlinewidth{1.003750pt}%
\definecolor{currentstroke}{rgb}{0.121569,0.286275,0.490196}%
\pgfsetstrokecolor{currentstroke}%
\pgfsetdash{}{0pt}%
\pgfsys@defobject{currentmarker}{\pgfqpoint{-0.041667in}{-0.041667in}}{\pgfqpoint{0.041667in}{0.041667in}}{%
\pgfpathmoveto{\pgfqpoint{-0.041667in}{-0.041667in}}%
\pgfpathlineto{\pgfqpoint{0.041667in}{-0.041667in}}%
\pgfpathlineto{\pgfqpoint{0.041667in}{0.041667in}}%
\pgfpathlineto{\pgfqpoint{-0.041667in}{0.041667in}}%
\pgfpathlineto{\pgfqpoint{-0.041667in}{-0.041667in}}%
\pgfpathclose%
\pgfusepath{stroke,fill}%
}%
\begin{pgfscope}%
\pgfsys@transformshift{2.501710in}{3.363810in}%
\pgfsys@useobject{currentmarker}{}%
\end{pgfscope}%
\end{pgfscope}%
\begin{pgfscope}%
\definecolor{textcolor}{rgb}{0.000000,0.000000,0.000000}%
\pgfsetstrokecolor{textcolor}%
\pgfsetfillcolor{textcolor}%
\pgftext[x=2.876710in,y=3.276310in,left,base]{\color{textcolor}\rmfamily\fontsize{18.000000}{21.600000}\selectfont LPP\(\displaystyle ^{R,T}\)}%
\end{pgfscope}%
\begin{pgfscope}%
\pgfsetbuttcap%
\pgfsetroundjoin%
\pgfsetlinewidth{1.505625pt}%
\definecolor{currentstroke}{rgb}{0.121569,0.286275,0.490196}%
\pgfsetstrokecolor{currentstroke}%
\pgfsetstrokeopacity{0.900000}%
\pgfsetdash{{5.550000pt}{2.400000pt}}{0.000000pt}%
\pgfpathmoveto{\pgfqpoint{3.812661in}{3.363810in}}%
\pgfpathlineto{\pgfqpoint{4.062661in}{3.363810in}}%
\pgfpathlineto{\pgfqpoint{4.312661in}{3.363810in}}%
\pgfusepath{stroke}%
\end{pgfscope}%
\begin{pgfscope}%
\pgfsetbuttcap%
\pgfsetmiterjoin%
\definecolor{currentfill}{rgb}{0.121569,0.286275,0.490196}%
\pgfsetfillcolor{currentfill}%
\pgfsetfillopacity{0.900000}%
\pgfsetlinewidth{1.003750pt}%
\definecolor{currentstroke}{rgb}{0.121569,0.286275,0.490196}%
\pgfsetstrokecolor{currentstroke}%
\pgfsetstrokeopacity{0.900000}%
\pgfsetdash{}{0pt}%
\pgfsys@defobject{currentmarker}{\pgfqpoint{-0.058926in}{-0.058926in}}{\pgfqpoint{0.058926in}{0.058926in}}{%
\pgfpathmoveto{\pgfqpoint{0.000000in}{-0.058926in}}%
\pgfpathlineto{\pgfqpoint{0.058926in}{0.000000in}}%
\pgfpathlineto{\pgfqpoint{0.000000in}{0.058926in}}%
\pgfpathlineto{\pgfqpoint{-0.058926in}{0.000000in}}%
\pgfpathlineto{\pgfqpoint{0.000000in}{-0.058926in}}%
\pgfpathclose%
\pgfusepath{stroke,fill}%
}%
\begin{pgfscope}%
\pgfsys@transformshift{4.062661in}{3.363810in}%
\pgfsys@useobject{currentmarker}{}%
\end{pgfscope}%
\end{pgfscope}%
\begin{pgfscope}%
\definecolor{textcolor}{rgb}{0.000000,0.000000,0.000000}%
\pgfsetstrokecolor{textcolor}%
\pgfsetfillcolor{textcolor}%
\pgftext[x=4.437661in,y=3.276310in,left,base]{\color{textcolor}\rmfamily\fontsize{18.000000}{21.600000}\selectfont LPP\(\displaystyle ^{T}\)}%
\end{pgfscope}%
\begin{pgfscope}%
\pgfsetbuttcap%
\pgfsetroundjoin%
\pgfsetlinewidth{1.505625pt}%
\definecolor{currentstroke}{rgb}{0.121569,0.286275,0.490196}%
\pgfsetstrokecolor{currentstroke}%
\pgfsetstrokeopacity{0.800000}%
\pgfsetdash{{1.500000pt}{2.475000pt}}{0.000000pt}%
\pgfpathmoveto{\pgfqpoint{5.202692in}{3.363810in}}%
\pgfpathlineto{\pgfqpoint{5.452692in}{3.363810in}}%
\pgfpathlineto{\pgfqpoint{5.702692in}{3.363810in}}%
\pgfusepath{stroke}%
\end{pgfscope}%
\begin{pgfscope}%
\pgfsetbuttcap%
\pgfsetroundjoin%
\definecolor{currentfill}{rgb}{0.121569,0.286275,0.490196}%
\pgfsetfillcolor{currentfill}%
\pgfsetfillopacity{0.800000}%
\pgfsetlinewidth{1.003750pt}%
\definecolor{currentstroke}{rgb}{0.121569,0.286275,0.490196}%
\pgfsetstrokecolor{currentstroke}%
\pgfsetstrokeopacity{0.800000}%
\pgfsetdash{}{0pt}%
\pgfsys@defobject{currentmarker}{\pgfqpoint{-0.041667in}{-0.041667in}}{\pgfqpoint{0.041667in}{0.041667in}}{%
\pgfpathmoveto{\pgfqpoint{0.000000in}{-0.041667in}}%
\pgfpathcurveto{\pgfqpoint{0.011050in}{-0.041667in}}{\pgfqpoint{0.021649in}{-0.037276in}}{\pgfqpoint{0.029463in}{-0.029463in}}%
\pgfpathcurveto{\pgfqpoint{0.037276in}{-0.021649in}}{\pgfqpoint{0.041667in}{-0.011050in}}{\pgfqpoint{0.041667in}{0.000000in}}%
\pgfpathcurveto{\pgfqpoint{0.041667in}{0.011050in}}{\pgfqpoint{0.037276in}{0.021649in}}{\pgfqpoint{0.029463in}{0.029463in}}%
\pgfpathcurveto{\pgfqpoint{0.021649in}{0.037276in}}{\pgfqpoint{0.011050in}{0.041667in}}{\pgfqpoint{0.000000in}{0.041667in}}%
\pgfpathcurveto{\pgfqpoint{-0.011050in}{0.041667in}}{\pgfqpoint{-0.021649in}{0.037276in}}{\pgfqpoint{-0.029463in}{0.029463in}}%
\pgfpathcurveto{\pgfqpoint{-0.037276in}{0.021649in}}{\pgfqpoint{-0.041667in}{0.011050in}}{\pgfqpoint{-0.041667in}{0.000000in}}%
\pgfpathcurveto{\pgfqpoint{-0.041667in}{-0.011050in}}{\pgfqpoint{-0.037276in}{-0.021649in}}{\pgfqpoint{-0.029463in}{-0.029463in}}%
\pgfpathcurveto{\pgfqpoint{-0.021649in}{-0.037276in}}{\pgfqpoint{-0.011050in}{-0.041667in}}{\pgfqpoint{0.000000in}{-0.041667in}}%
\pgfpathlineto{\pgfqpoint{0.000000in}{-0.041667in}}%
\pgfpathclose%
\pgfusepath{stroke,fill}%
}%
\begin{pgfscope}%
\pgfsys@transformshift{5.452692in}{3.363810in}%
\pgfsys@useobject{currentmarker}{}%
\end{pgfscope}%
\end{pgfscope}%
\begin{pgfscope}%
\definecolor{textcolor}{rgb}{0.000000,0.000000,0.000000}%
\pgfsetstrokecolor{textcolor}%
\pgfsetfillcolor{textcolor}%
\pgftext[x=5.827692in,y=3.276310in,left,base]{\color{textcolor}\rmfamily\fontsize{18.000000}{21.600000}\selectfont LPP\(\displaystyle ^{T, P}\)}%
\end{pgfscope}%
\end{pgfpicture}%
\makeatother%
\endgroup%

%% file: figures/DemandPlots/TightLayout/DE_OC.pgf
\begingroup%
\makeatletter%
\begin{pgfpicture}%
\pgfpathrectangle{\pgfpointorigin}{\pgfqpoint{7.780338in}{3.651310in}}%
\pgfusepath{use as bounding box, clip}%
\begin{pgfscope}%
\pgfsetbuttcap%
\pgfsetmiterjoin%
\definecolor{currentfill}{rgb}{1.000000,1.000000,1.000000}%
\pgfsetfillcolor{currentfill}%
\pgfsetlinewidth{0.000000pt}%
\definecolor{currentstroke}{rgb}{1.000000,1.000000,1.000000}%
\pgfsetstrokecolor{currentstroke}%
\pgfsetdash{}{0pt}%
\pgfpathmoveto{\pgfqpoint{0.000000in}{0.000000in}}%
\pgfpathlineto{\pgfqpoint{7.780338in}{0.000000in}}%
\pgfpathlineto{\pgfqpoint{7.780338in}{3.651310in}}%
\pgfpathlineto{\pgfqpoint{0.000000in}{3.651310in}}%
\pgfpathlineto{\pgfqpoint{0.000000in}{0.000000in}}%
\pgfpathclose%
\pgfusepath{fill}%
\end{pgfscope}%
\begin{pgfscope}%
\pgfsetbuttcap%
\pgfsetmiterjoin%
\definecolor{currentfill}{rgb}{1.000000,1.000000,1.000000}%
\pgfsetfillcolor{currentfill}%
\pgfsetlinewidth{0.000000pt}%
\definecolor{currentstroke}{rgb}{0.000000,0.000000,0.000000}%
\pgfsetstrokecolor{currentstroke}%
\pgfsetstrokeopacity{0.000000}%
\pgfsetdash{}{0pt}%
\pgfpathmoveto{\pgfqpoint{1.016468in}{0.679475in}}%
\pgfpathlineto{\pgfqpoint{7.680338in}{0.679475in}}%
\pgfpathlineto{\pgfqpoint{7.680338in}{3.135897in}}%
\pgfpathlineto{\pgfqpoint{1.016468in}{3.135897in}}%
\pgfpathlineto{\pgfqpoint{1.016468in}{0.679475in}}%
\pgfpathclose%
\pgfusepath{fill}%
\end{pgfscope}%
\begin{pgfscope}%
\pgfsetbuttcap%
\pgfsetroundjoin%
\definecolor{currentfill}{rgb}{0.000000,0.000000,0.000000}%
\pgfsetfillcolor{currentfill}%
\pgfsetlinewidth{0.803000pt}%
\definecolor{currentstroke}{rgb}{0.000000,0.000000,0.000000}%
\pgfsetstrokecolor{currentstroke}%
\pgfsetdash{}{0pt}%
\pgfsys@defobject{currentmarker}{\pgfqpoint{0.000000in}{-0.048611in}}{\pgfqpoint{0.000000in}{0.000000in}}{%
\pgfpathmoveto{\pgfqpoint{0.000000in}{0.000000in}}%
\pgfpathlineto{\pgfqpoint{0.000000in}{-0.048611in}}%
\pgfusepath{stroke,fill}%
}%
\begin{pgfscope}%
\pgfsys@transformshift{1.016468in}{0.679475in}%
\pgfsys@useobject{currentmarker}{}%
\end{pgfscope}%
\end{pgfscope}%
\begin{pgfscope}%
\definecolor{textcolor}{rgb}{0.000000,0.000000,0.000000}%
\pgfsetstrokecolor{textcolor}%
\pgfsetfillcolor{textcolor}%
\pgftext[x=1.016468in,y=0.582253in,,top]{\color{textcolor}\rmfamily\fontsize{16.000000}{19.200000}\selectfont \(\displaystyle {0.0}\)}%
\end{pgfscope}%
\begin{pgfscope}%
\pgfsetbuttcap%
\pgfsetroundjoin%
\definecolor{currentfill}{rgb}{0.000000,0.000000,0.000000}%
\pgfsetfillcolor{currentfill}%
\pgfsetlinewidth{0.803000pt}%
\definecolor{currentstroke}{rgb}{0.000000,0.000000,0.000000}%
\pgfsetstrokecolor{currentstroke}%
\pgfsetdash{}{0pt}%
\pgfsys@defobject{currentmarker}{\pgfqpoint{0.000000in}{-0.048611in}}{\pgfqpoint{0.000000in}{0.000000in}}{%
\pgfpathmoveto{\pgfqpoint{0.000000in}{0.000000in}}%
\pgfpathlineto{\pgfqpoint{0.000000in}{-0.048611in}}%
\pgfusepath{stroke,fill}%
}%
\begin{pgfscope}%
\pgfsys@transformshift{2.228080in}{0.679475in}%
\pgfsys@useobject{currentmarker}{}%
\end{pgfscope}%
\end{pgfscope}%
\begin{pgfscope}%
\definecolor{textcolor}{rgb}{0.000000,0.000000,0.000000}%
\pgfsetstrokecolor{textcolor}%
\pgfsetfillcolor{textcolor}%
\pgftext[x=2.228080in,y=0.582253in,,top]{\color{textcolor}\rmfamily\fontsize{16.000000}{19.200000}\selectfont \(\displaystyle {0.2}\)}%
\end{pgfscope}%
\begin{pgfscope}%
\pgfsetbuttcap%
\pgfsetroundjoin%
\definecolor{currentfill}{rgb}{0.000000,0.000000,0.000000}%
\pgfsetfillcolor{currentfill}%
\pgfsetlinewidth{0.803000pt}%
\definecolor{currentstroke}{rgb}{0.000000,0.000000,0.000000}%
\pgfsetstrokecolor{currentstroke}%
\pgfsetdash{}{0pt}%
\pgfsys@defobject{currentmarker}{\pgfqpoint{0.000000in}{-0.048611in}}{\pgfqpoint{0.000000in}{0.000000in}}{%
\pgfpathmoveto{\pgfqpoint{0.000000in}{0.000000in}}%
\pgfpathlineto{\pgfqpoint{0.000000in}{-0.048611in}}%
\pgfusepath{stroke,fill}%
}%
\begin{pgfscope}%
\pgfsys@transformshift{3.439693in}{0.679475in}%
\pgfsys@useobject{currentmarker}{}%
\end{pgfscope}%
\end{pgfscope}%
\begin{pgfscope}%
\definecolor{textcolor}{rgb}{0.000000,0.000000,0.000000}%
\pgfsetstrokecolor{textcolor}%
\pgfsetfillcolor{textcolor}%
\pgftext[x=3.439693in,y=0.582253in,,top]{\color{textcolor}\rmfamily\fontsize{16.000000}{19.200000}\selectfont \(\displaystyle {0.4}\)}%
\end{pgfscope}%
\begin{pgfscope}%
\pgfsetbuttcap%
\pgfsetroundjoin%
\definecolor{currentfill}{rgb}{0.000000,0.000000,0.000000}%
\pgfsetfillcolor{currentfill}%
\pgfsetlinewidth{0.803000pt}%
\definecolor{currentstroke}{rgb}{0.000000,0.000000,0.000000}%
\pgfsetstrokecolor{currentstroke}%
\pgfsetdash{}{0pt}%
\pgfsys@defobject{currentmarker}{\pgfqpoint{0.000000in}{-0.048611in}}{\pgfqpoint{0.000000in}{0.000000in}}{%
\pgfpathmoveto{\pgfqpoint{0.000000in}{0.000000in}}%
\pgfpathlineto{\pgfqpoint{0.000000in}{-0.048611in}}%
\pgfusepath{stroke,fill}%
}%
\begin{pgfscope}%
\pgfsys@transformshift{4.651306in}{0.679475in}%
\pgfsys@useobject{currentmarker}{}%
\end{pgfscope}%
\end{pgfscope}%
\begin{pgfscope}%
\definecolor{textcolor}{rgb}{0.000000,0.000000,0.000000}%
\pgfsetstrokecolor{textcolor}%
\pgfsetfillcolor{textcolor}%
\pgftext[x=4.651306in,y=0.582253in,,top]{\color{textcolor}\rmfamily\fontsize{16.000000}{19.200000}\selectfont \(\displaystyle {0.6}\)}%
\end{pgfscope}%
\begin{pgfscope}%
\pgfsetbuttcap%
\pgfsetroundjoin%
\definecolor{currentfill}{rgb}{0.000000,0.000000,0.000000}%
\pgfsetfillcolor{currentfill}%
\pgfsetlinewidth{0.803000pt}%
\definecolor{currentstroke}{rgb}{0.000000,0.000000,0.000000}%
\pgfsetstrokecolor{currentstroke}%
\pgfsetdash{}{0pt}%
\pgfsys@defobject{currentmarker}{\pgfqpoint{0.000000in}{-0.048611in}}{\pgfqpoint{0.000000in}{0.000000in}}{%
\pgfpathmoveto{\pgfqpoint{0.000000in}{0.000000in}}%
\pgfpathlineto{\pgfqpoint{0.000000in}{-0.048611in}}%
\pgfusepath{stroke,fill}%
}%
\begin{pgfscope}%
\pgfsys@transformshift{5.862919in}{0.679475in}%
\pgfsys@useobject{currentmarker}{}%
\end{pgfscope}%
\end{pgfscope}%
\begin{pgfscope}%
\definecolor{textcolor}{rgb}{0.000000,0.000000,0.000000}%
\pgfsetstrokecolor{textcolor}%
\pgfsetfillcolor{textcolor}%
\pgftext[x=5.862919in,y=0.582253in,,top]{\color{textcolor}\rmfamily\fontsize{16.000000}{19.200000}\selectfont \(\displaystyle {0.8}\)}%
\end{pgfscope}%
\begin{pgfscope}%
\pgfsetbuttcap%
\pgfsetroundjoin%
\definecolor{currentfill}{rgb}{0.000000,0.000000,0.000000}%
\pgfsetfillcolor{currentfill}%
\pgfsetlinewidth{0.803000pt}%
\definecolor{currentstroke}{rgb}{0.000000,0.000000,0.000000}%
\pgfsetstrokecolor{currentstroke}%
\pgfsetdash{}{0pt}%
\pgfsys@defobject{currentmarker}{\pgfqpoint{0.000000in}{-0.048611in}}{\pgfqpoint{0.000000in}{0.000000in}}{%
\pgfpathmoveto{\pgfqpoint{0.000000in}{0.000000in}}%
\pgfpathlineto{\pgfqpoint{0.000000in}{-0.048611in}}%
\pgfusepath{stroke,fill}%
}%
\begin{pgfscope}%
\pgfsys@transformshift{7.074532in}{0.679475in}%
\pgfsys@useobject{currentmarker}{}%
\end{pgfscope}%
\end{pgfscope}%
\begin{pgfscope}%
\definecolor{textcolor}{rgb}{0.000000,0.000000,0.000000}%
\pgfsetstrokecolor{textcolor}%
\pgfsetfillcolor{textcolor}%
\pgftext[x=7.074532in,y=0.582253in,,top]{\color{textcolor}\rmfamily\fontsize{16.000000}{19.200000}\selectfont \(\displaystyle {1.0}\)}%
\end{pgfscope}%
\begin{pgfscope}%
\definecolor{textcolor}{rgb}{0.000000,0.000000,0.000000}%
\pgfsetstrokecolor{textcolor}%
\pgfsetfillcolor{textcolor}%
\pgftext[x=4.348403in,y=0.313349in,,top]{\color{textcolor}\rmfamily\fontsize{18.000000}{21.600000}\selectfont \(\displaystyle \lambda\)}%
\end{pgfscope}%
\begin{pgfscope}%
\pgfsetbuttcap%
\pgfsetroundjoin%
\definecolor{currentfill}{rgb}{0.000000,0.000000,0.000000}%
\pgfsetfillcolor{currentfill}%
\pgfsetlinewidth{0.803000pt}%
\definecolor{currentstroke}{rgb}{0.000000,0.000000,0.000000}%
\pgfsetstrokecolor{currentstroke}%
\pgfsetdash{}{0pt}%
\pgfsys@defobject{currentmarker}{\pgfqpoint{-0.048611in}{0.000000in}}{\pgfqpoint{-0.000000in}{0.000000in}}{%
\pgfpathmoveto{\pgfqpoint{-0.000000in}{0.000000in}}%
\pgfpathlineto{\pgfqpoint{-0.048611in}{0.000000in}}%
\pgfusepath{stroke,fill}%
}%
\begin{pgfscope}%
\pgfsys@transformshift{1.016468in}{0.945515in}%
\pgfsys@useobject{currentmarker}{}%
\end{pgfscope}%
\end{pgfscope}%
\begin{pgfscope}%
\definecolor{textcolor}{rgb}{0.000000,0.000000,0.000000}%
\pgfsetstrokecolor{textcolor}%
\pgfsetfillcolor{textcolor}%
\pgftext[x=0.368904in, y=0.862182in, left, base]{\color{textcolor}\rmfamily\fontsize{16.000000}{19.200000}\selectfont \(\displaystyle {20000}\)}%
\end{pgfscope}%
\begin{pgfscope}%
\pgfsetbuttcap%
\pgfsetroundjoin%
\definecolor{currentfill}{rgb}{0.000000,0.000000,0.000000}%
\pgfsetfillcolor{currentfill}%
\pgfsetlinewidth{0.803000pt}%
\definecolor{currentstroke}{rgb}{0.000000,0.000000,0.000000}%
\pgfsetstrokecolor{currentstroke}%
\pgfsetdash{}{0pt}%
\pgfsys@defobject{currentmarker}{\pgfqpoint{-0.048611in}{0.000000in}}{\pgfqpoint{-0.000000in}{0.000000in}}{%
\pgfpathmoveto{\pgfqpoint{-0.000000in}{0.000000in}}%
\pgfpathlineto{\pgfqpoint{-0.048611in}{0.000000in}}%
\pgfusepath{stroke,fill}%
}%
\begin{pgfscope}%
\pgfsys@transformshift{1.016468in}{1.477353in}%
\pgfsys@useobject{currentmarker}{}%
\end{pgfscope}%
\end{pgfscope}%
\begin{pgfscope}%
\definecolor{textcolor}{rgb}{0.000000,0.000000,0.000000}%
\pgfsetstrokecolor{textcolor}%
\pgfsetfillcolor{textcolor}%
\pgftext[x=0.368904in, y=1.394020in, left, base]{\color{textcolor}\rmfamily\fontsize{16.000000}{19.200000}\selectfont \(\displaystyle {30000}\)}%
\end{pgfscope}%
\begin{pgfscope}%
\pgfsetbuttcap%
\pgfsetroundjoin%
\definecolor{currentfill}{rgb}{0.000000,0.000000,0.000000}%
\pgfsetfillcolor{currentfill}%
\pgfsetlinewidth{0.803000pt}%
\definecolor{currentstroke}{rgb}{0.000000,0.000000,0.000000}%
\pgfsetstrokecolor{currentstroke}%
\pgfsetdash{}{0pt}%
\pgfsys@defobject{currentmarker}{\pgfqpoint{-0.048611in}{0.000000in}}{\pgfqpoint{-0.000000in}{0.000000in}}{%
\pgfpathmoveto{\pgfqpoint{-0.000000in}{0.000000in}}%
\pgfpathlineto{\pgfqpoint{-0.048611in}{0.000000in}}%
\pgfusepath{stroke,fill}%
}%
\begin{pgfscope}%
\pgfsys@transformshift{1.016468in}{2.009191in}%
\pgfsys@useobject{currentmarker}{}%
\end{pgfscope}%
\end{pgfscope}%
\begin{pgfscope}%
\definecolor{textcolor}{rgb}{0.000000,0.000000,0.000000}%
\pgfsetstrokecolor{textcolor}%
\pgfsetfillcolor{textcolor}%
\pgftext[x=0.368904in, y=1.925858in, left, base]{\color{textcolor}\rmfamily\fontsize{16.000000}{19.200000}\selectfont \(\displaystyle {40000}\)}%
\end{pgfscope}%
\begin{pgfscope}%
\pgfsetbuttcap%
\pgfsetroundjoin%
\definecolor{currentfill}{rgb}{0.000000,0.000000,0.000000}%
\pgfsetfillcolor{currentfill}%
\pgfsetlinewidth{0.803000pt}%
\definecolor{currentstroke}{rgb}{0.000000,0.000000,0.000000}%
\pgfsetstrokecolor{currentstroke}%
\pgfsetdash{}{0pt}%
\pgfsys@defobject{currentmarker}{\pgfqpoint{-0.048611in}{0.000000in}}{\pgfqpoint{-0.000000in}{0.000000in}}{%
\pgfpathmoveto{\pgfqpoint{-0.000000in}{0.000000in}}%
\pgfpathlineto{\pgfqpoint{-0.048611in}{0.000000in}}%
\pgfusepath{stroke,fill}%
}%
\begin{pgfscope}%
\pgfsys@transformshift{1.016468in}{2.541029in}%
\pgfsys@useobject{currentmarker}{}%
\end{pgfscope}%
\end{pgfscope}%
\begin{pgfscope}%
\definecolor{textcolor}{rgb}{0.000000,0.000000,0.000000}%
\pgfsetstrokecolor{textcolor}%
\pgfsetfillcolor{textcolor}%
\pgftext[x=0.368904in, y=2.457695in, left, base]{\color{textcolor}\rmfamily\fontsize{16.000000}{19.200000}\selectfont \(\displaystyle {50000}\)}%
\end{pgfscope}%
\begin{pgfscope}%
\pgfsetbuttcap%
\pgfsetroundjoin%
\definecolor{currentfill}{rgb}{0.000000,0.000000,0.000000}%
\pgfsetfillcolor{currentfill}%
\pgfsetlinewidth{0.803000pt}%
\definecolor{currentstroke}{rgb}{0.000000,0.000000,0.000000}%
\pgfsetstrokecolor{currentstroke}%
\pgfsetdash{}{0pt}%
\pgfsys@defobject{currentmarker}{\pgfqpoint{-0.048611in}{0.000000in}}{\pgfqpoint{-0.000000in}{0.000000in}}{%
\pgfpathmoveto{\pgfqpoint{-0.000000in}{0.000000in}}%
\pgfpathlineto{\pgfqpoint{-0.048611in}{0.000000in}}%
\pgfusepath{stroke,fill}%
}%
\begin{pgfscope}%
\pgfsys@transformshift{1.016468in}{3.072866in}%
\pgfsys@useobject{currentmarker}{}%
\end{pgfscope}%
\end{pgfscope}%
\begin{pgfscope}%
\definecolor{textcolor}{rgb}{0.000000,0.000000,0.000000}%
\pgfsetstrokecolor{textcolor}%
\pgfsetfillcolor{textcolor}%
\pgftext[x=0.368904in, y=2.989533in, left, base]{\color{textcolor}\rmfamily\fontsize{16.000000}{19.200000}\selectfont \(\displaystyle {60000}\)}%
\end{pgfscope}%
\begin{pgfscope}%
\definecolor{textcolor}{rgb}{0.000000,0.000000,0.000000}%
\pgfsetstrokecolor{textcolor}%
\pgfsetfillcolor{textcolor}%
\pgftext[x=0.313349in,y=1.907686in,,bottom,rotate=90.000000]{\color{textcolor}\rmfamily\fontsize{18.000000}{21.600000}\selectfont DKK}%
\end{pgfscope}%
\begin{pgfscope}%
\pgfpathrectangle{\pgfqpoint{1.016468in}{0.679475in}}{\pgfqpoint{6.663871in}{2.456422in}}%
\pgfusepath{clip}%
\pgfsetrectcap%
\pgfsetroundjoin%
\pgfsetlinewidth{1.505625pt}%
\definecolor{currentstroke}{rgb}{0.839216,0.152941,0.156863}%
\pgfsetstrokecolor{currentstroke}%
\pgfsetdash{}{0pt}%
\pgfpathmoveto{\pgfqpoint{1.622274in}{2.403174in}}%
\pgfpathlineto{\pgfqpoint{2.530984in}{2.432883in}}%
\pgfpathlineto{\pgfqpoint{4.045500in}{2.518983in}}%
\pgfpathlineto{\pgfqpoint{5.560016in}{2.703315in}}%
\pgfpathlineto{\pgfqpoint{7.074532in}{3.024241in}}%
\pgfusepath{stroke}%
\end{pgfscope}%
\begin{pgfscope}%
\pgfpathrectangle{\pgfqpoint{1.016468in}{0.679475in}}{\pgfqpoint{6.663871in}{2.456422in}}%
\pgfusepath{clip}%
\pgfsetbuttcap%
\pgfsetmiterjoin%
\definecolor{currentfill}{rgb}{0.839216,0.152941,0.156863}%
\pgfsetfillcolor{currentfill}%
\pgfsetlinewidth{0.752812pt}%
\definecolor{currentstroke}{rgb}{1.000000,1.000000,1.000000}%
\pgfsetstrokecolor{currentstroke}%
\pgfsetdash{}{0pt}%
\pgfsys@defobject{currentmarker}{\pgfqpoint{-0.041667in}{-0.041667in}}{\pgfqpoint{0.041667in}{0.041667in}}{%
\pgfpathmoveto{\pgfqpoint{-0.041667in}{-0.041667in}}%
\pgfpathlineto{\pgfqpoint{0.041667in}{-0.041667in}}%
\pgfpathlineto{\pgfqpoint{0.041667in}{0.041667in}}%
\pgfpathlineto{\pgfqpoint{-0.041667in}{0.041667in}}%
\pgfpathlineto{\pgfqpoint{-0.041667in}{-0.041667in}}%
\pgfpathclose%
\pgfusepath{stroke,fill}%
}%
\begin{pgfscope}%
\pgfsys@transformshift{1.622274in}{2.403174in}%
\pgfsys@useobject{currentmarker}{}%
\end{pgfscope}%
\begin{pgfscope}%
\pgfsys@transformshift{2.530984in}{2.432883in}%
\pgfsys@useobject{currentmarker}{}%
\end{pgfscope}%
\begin{pgfscope}%
\pgfsys@transformshift{4.045500in}{2.518983in}%
\pgfsys@useobject{currentmarker}{}%
\end{pgfscope}%
\begin{pgfscope}%
\pgfsys@transformshift{5.560016in}{2.703315in}%
\pgfsys@useobject{currentmarker}{}%
\end{pgfscope}%
\begin{pgfscope}%
\pgfsys@transformshift{7.074532in}{3.024241in}%
\pgfsys@useobject{currentmarker}{}%
\end{pgfscope}%
\end{pgfscope}%
\begin{pgfscope}%
\pgfpathrectangle{\pgfqpoint{1.016468in}{0.679475in}}{\pgfqpoint{6.663871in}{2.456422in}}%
\pgfusepath{clip}%
\pgfsetbuttcap%
\pgfsetroundjoin%
\pgfsetlinewidth{1.505625pt}%
\definecolor{currentstroke}{rgb}{0.839216,0.152941,0.156863}%
\pgfsetstrokecolor{currentstroke}%
\pgfsetstrokeopacity{0.900000}%
\pgfsetdash{{7.500000pt}{7.500000pt}}{0.000000pt}%
\pgfpathmoveto{\pgfqpoint{1.622274in}{0.791131in}}%
\pgfpathlineto{\pgfqpoint{2.530984in}{0.892249in}}%
\pgfpathlineto{\pgfqpoint{4.045500in}{1.167416in}}%
\pgfpathlineto{\pgfqpoint{5.560016in}{1.439669in}}%
\pgfpathlineto{\pgfqpoint{7.074532in}{1.740537in}}%
\pgfusepath{stroke}%
\end{pgfscope}%
\begin{pgfscope}%
\pgfpathrectangle{\pgfqpoint{1.016468in}{0.679475in}}{\pgfqpoint{6.663871in}{2.456422in}}%
\pgfusepath{clip}%
\pgfsetbuttcap%
\pgfsetmiterjoin%
\definecolor{currentfill}{rgb}{0.839216,0.152941,0.156863}%
\pgfsetfillcolor{currentfill}%
\pgfsetfillopacity{0.900000}%
\pgfsetlinewidth{0.752812pt}%
\definecolor{currentstroke}{rgb}{1.000000,1.000000,1.000000}%
\pgfsetstrokecolor{currentstroke}%
\pgfsetdash{}{0pt}%
\pgfsys@defobject{currentmarker}{\pgfqpoint{-0.058926in}{-0.058926in}}{\pgfqpoint{0.058926in}{0.058926in}}{%
\pgfpathmoveto{\pgfqpoint{0.000000in}{-0.058926in}}%
\pgfpathlineto{\pgfqpoint{0.058926in}{0.000000in}}%
\pgfpathlineto{\pgfqpoint{0.000000in}{0.058926in}}%
\pgfpathlineto{\pgfqpoint{-0.058926in}{0.000000in}}%
\pgfpathlineto{\pgfqpoint{0.000000in}{-0.058926in}}%
\pgfpathclose%
\pgfusepath{stroke,fill}%
}%
\begin{pgfscope}%
\pgfsys@transformshift{1.622274in}{0.791131in}%
\pgfsys@useobject{currentmarker}{}%
\end{pgfscope}%
\begin{pgfscope}%
\pgfsys@transformshift{2.530984in}{0.892249in}%
\pgfsys@useobject{currentmarker}{}%
\end{pgfscope}%
\begin{pgfscope}%
\pgfsys@transformshift{4.045500in}{1.167416in}%
\pgfsys@useobject{currentmarker}{}%
\end{pgfscope}%
\begin{pgfscope}%
\pgfsys@transformshift{5.560016in}{1.439669in}%
\pgfsys@useobject{currentmarker}{}%
\end{pgfscope}%
\begin{pgfscope}%
\pgfsys@transformshift{7.074532in}{1.740537in}%
\pgfsys@useobject{currentmarker}{}%
\end{pgfscope}%
\end{pgfscope}%
\begin{pgfscope}%
\pgfpathrectangle{\pgfqpoint{1.016468in}{0.679475in}}{\pgfqpoint{6.663871in}{2.456422in}}%
\pgfusepath{clip}%
\pgfsetbuttcap%
\pgfsetroundjoin%
\pgfsetlinewidth{1.505625pt}%
\definecolor{currentstroke}{rgb}{0.839216,0.152941,0.156863}%
\pgfsetstrokecolor{currentstroke}%
\pgfsetstrokeopacity{0.800000}%
\pgfsetdash{{1.500000pt}{4.500000pt}}{0.000000pt}%
\pgfpathmoveto{\pgfqpoint{1.622274in}{2.403169in}}%
\pgfpathlineto{\pgfqpoint{2.530984in}{2.432883in}}%
\pgfpathlineto{\pgfqpoint{4.045500in}{2.518980in}}%
\pgfpathlineto{\pgfqpoint{5.560016in}{2.703315in}}%
\pgfpathlineto{\pgfqpoint{7.074532in}{3.024241in}}%
\pgfusepath{stroke}%
\end{pgfscope}%
\begin{pgfscope}%
\pgfpathrectangle{\pgfqpoint{1.016468in}{0.679475in}}{\pgfqpoint{6.663871in}{2.456422in}}%
\pgfusepath{clip}%
\pgfsetbuttcap%
\pgfsetroundjoin%
\definecolor{currentfill}{rgb}{0.839216,0.152941,0.156863}%
\pgfsetfillcolor{currentfill}%
\pgfsetfillopacity{0.800000}%
\pgfsetlinewidth{0.752812pt}%
\definecolor{currentstroke}{rgb}{1.000000,1.000000,1.000000}%
\pgfsetstrokecolor{currentstroke}%
\pgfsetdash{}{0pt}%
\pgfsys@defobject{currentmarker}{\pgfqpoint{-0.041667in}{-0.041667in}}{\pgfqpoint{0.041667in}{0.041667in}}{%
\pgfpathmoveto{\pgfqpoint{0.000000in}{-0.041667in}}%
\pgfpathcurveto{\pgfqpoint{0.011050in}{-0.041667in}}{\pgfqpoint{0.021649in}{-0.037276in}}{\pgfqpoint{0.029463in}{-0.029463in}}%
\pgfpathcurveto{\pgfqpoint{0.037276in}{-0.021649in}}{\pgfqpoint{0.041667in}{-0.011050in}}{\pgfqpoint{0.041667in}{0.000000in}}%
\pgfpathcurveto{\pgfqpoint{0.041667in}{0.011050in}}{\pgfqpoint{0.037276in}{0.021649in}}{\pgfqpoint{0.029463in}{0.029463in}}%
\pgfpathcurveto{\pgfqpoint{0.021649in}{0.037276in}}{\pgfqpoint{0.011050in}{0.041667in}}{\pgfqpoint{0.000000in}{0.041667in}}%
\pgfpathcurveto{\pgfqpoint{-0.011050in}{0.041667in}}{\pgfqpoint{-0.021649in}{0.037276in}}{\pgfqpoint{-0.029463in}{0.029463in}}%
\pgfpathcurveto{\pgfqpoint{-0.037276in}{0.021649in}}{\pgfqpoint{-0.041667in}{0.011050in}}{\pgfqpoint{-0.041667in}{0.000000in}}%
\pgfpathcurveto{\pgfqpoint{-0.041667in}{-0.011050in}}{\pgfqpoint{-0.037276in}{-0.021649in}}{\pgfqpoint{-0.029463in}{-0.029463in}}%
\pgfpathcurveto{\pgfqpoint{-0.021649in}{-0.037276in}}{\pgfqpoint{-0.011050in}{-0.041667in}}{\pgfqpoint{0.000000in}{-0.041667in}}%
\pgfpathlineto{\pgfqpoint{0.000000in}{-0.041667in}}%
\pgfpathclose%
\pgfusepath{stroke,fill}%
}%
\begin{pgfscope}%
\pgfsys@transformshift{1.622274in}{2.403169in}%
\pgfsys@useobject{currentmarker}{}%
\end{pgfscope}%
\begin{pgfscope}%
\pgfsys@transformshift{2.530984in}{2.432883in}%
\pgfsys@useobject{currentmarker}{}%
\end{pgfscope}%
\begin{pgfscope}%
\pgfsys@transformshift{4.045500in}{2.518980in}%
\pgfsys@useobject{currentmarker}{}%
\end{pgfscope}%
\begin{pgfscope}%
\pgfsys@transformshift{5.560016in}{2.703315in}%
\pgfsys@useobject{currentmarker}{}%
\end{pgfscope}%
\begin{pgfscope}%
\pgfsys@transformshift{7.074532in}{3.024241in}%
\pgfsys@useobject{currentmarker}{}%
\end{pgfscope}%
\end{pgfscope}%
\begin{pgfscope}%
\pgfsetrectcap%
\pgfsetmiterjoin%
\pgfsetlinewidth{0.803000pt}%
\definecolor{currentstroke}{rgb}{0.000000,0.000000,0.000000}%
\pgfsetstrokecolor{currentstroke}%
\pgfsetdash{}{0pt}%
\pgfpathmoveto{\pgfqpoint{1.016468in}{0.679475in}}%
\pgfpathlineto{\pgfqpoint{1.016468in}{3.135897in}}%
\pgfusepath{stroke}%
\end{pgfscope}%
\begin{pgfscope}%
\pgfsetrectcap%
\pgfsetmiterjoin%
\pgfsetlinewidth{0.803000pt}%
\definecolor{currentstroke}{rgb}{0.000000,0.000000,0.000000}%
\pgfsetstrokecolor{currentstroke}%
\pgfsetdash{}{0pt}%
\pgfpathmoveto{\pgfqpoint{7.680338in}{0.679475in}}%
\pgfpathlineto{\pgfqpoint{7.680338in}{3.135897in}}%
\pgfusepath{stroke}%
\end{pgfscope}%
\begin{pgfscope}%
\pgfsetrectcap%
\pgfsetmiterjoin%
\pgfsetlinewidth{0.803000pt}%
\definecolor{currentstroke}{rgb}{0.000000,0.000000,0.000000}%
\pgfsetstrokecolor{currentstroke}%
\pgfsetdash{}{0pt}%
\pgfpathmoveto{\pgfqpoint{1.016468in}{0.679475in}}%
\pgfpathlineto{\pgfqpoint{7.680338in}{0.679475in}}%
\pgfusepath{stroke}%
\end{pgfscope}%
\begin{pgfscope}%
\pgfsetrectcap%
\pgfsetmiterjoin%
\pgfsetlinewidth{0.803000pt}%
\definecolor{currentstroke}{rgb}{0.000000,0.000000,0.000000}%
\pgfsetstrokecolor{currentstroke}%
\pgfsetdash{}{0pt}%
\pgfpathmoveto{\pgfqpoint{1.016468in}{3.135897in}}%
\pgfpathlineto{\pgfqpoint{7.680338in}{3.135897in}}%
\pgfusepath{stroke}%
\end{pgfscope}%
\begin{pgfscope}%
\pgfsetrectcap%
\pgfsetroundjoin%
\pgfsetlinewidth{1.505625pt}%
\definecolor{currentstroke}{rgb}{0.839216,0.152941,0.156863}%
\pgfsetstrokecolor{currentstroke}%
\pgfsetdash{}{0pt}%
\pgfpathmoveto{\pgfqpoint{2.195836in}{3.363810in}}%
\pgfpathlineto{\pgfqpoint{2.445836in}{3.363810in}}%
\pgfpathlineto{\pgfqpoint{2.695836in}{3.363810in}}%
\pgfusepath{stroke}%
\end{pgfscope}%
\begin{pgfscope}%
\pgfsetbuttcap%
\pgfsetmiterjoin%
\definecolor{currentfill}{rgb}{0.839216,0.152941,0.156863}%
\pgfsetfillcolor{currentfill}%
\pgfsetlinewidth{1.003750pt}%
\definecolor{currentstroke}{rgb}{0.839216,0.152941,0.156863}%
\pgfsetstrokecolor{currentstroke}%
\pgfsetdash{}{0pt}%
\pgfsys@defobject{currentmarker}{\pgfqpoint{-0.041667in}{-0.041667in}}{\pgfqpoint{0.041667in}{0.041667in}}{%
\pgfpathmoveto{\pgfqpoint{-0.041667in}{-0.041667in}}%
\pgfpathlineto{\pgfqpoint{0.041667in}{-0.041667in}}%
\pgfpathlineto{\pgfqpoint{0.041667in}{0.041667in}}%
\pgfpathlineto{\pgfqpoint{-0.041667in}{0.041667in}}%
\pgfpathlineto{\pgfqpoint{-0.041667in}{-0.041667in}}%
\pgfpathclose%
\pgfusepath{stroke,fill}%
}%
\begin{pgfscope}%
\pgfsys@transformshift{2.445836in}{3.363810in}%
\pgfsys@useobject{currentmarker}{}%
\end{pgfscope}%
\end{pgfscope}%
\begin{pgfscope}%
\definecolor{textcolor}{rgb}{0.000000,0.000000,0.000000}%
\pgfsetstrokecolor{textcolor}%
\pgfsetfillcolor{textcolor}%
\pgftext[x=2.820836in,y=3.276310in,left,base]{\color{textcolor}\rmfamily\fontsize{18.000000}{21.600000}\selectfont LPP\(\displaystyle ^{R,T}\)}%
\end{pgfscope}%
\begin{pgfscope}%
\pgfsetbuttcap%
\pgfsetroundjoin%
\pgfsetlinewidth{1.505625pt}%
\definecolor{currentstroke}{rgb}{0.839216,0.152941,0.156863}%
\pgfsetstrokecolor{currentstroke}%
\pgfsetstrokeopacity{0.900000}%
\pgfsetdash{{5.550000pt}{2.400000pt}}{0.000000pt}%
\pgfpathmoveto{\pgfqpoint{3.756787in}{3.363810in}}%
\pgfpathlineto{\pgfqpoint{4.006787in}{3.363810in}}%
\pgfpathlineto{\pgfqpoint{4.256787in}{3.363810in}}%
\pgfusepath{stroke}%
\end{pgfscope}%
\begin{pgfscope}%
\pgfsetbuttcap%
\pgfsetmiterjoin%
\definecolor{currentfill}{rgb}{0.839216,0.152941,0.156863}%
\pgfsetfillcolor{currentfill}%
\pgfsetfillopacity{0.900000}%
\pgfsetlinewidth{1.003750pt}%
\definecolor{currentstroke}{rgb}{0.839216,0.152941,0.156863}%
\pgfsetstrokecolor{currentstroke}%
\pgfsetstrokeopacity{0.900000}%
\pgfsetdash{}{0pt}%
\pgfsys@defobject{currentmarker}{\pgfqpoint{-0.058926in}{-0.058926in}}{\pgfqpoint{0.058926in}{0.058926in}}{%
\pgfpathmoveto{\pgfqpoint{0.000000in}{-0.058926in}}%
\pgfpathlineto{\pgfqpoint{0.058926in}{0.000000in}}%
\pgfpathlineto{\pgfqpoint{0.000000in}{0.058926in}}%
\pgfpathlineto{\pgfqpoint{-0.058926in}{0.000000in}}%
\pgfpathlineto{\pgfqpoint{0.000000in}{-0.058926in}}%
\pgfpathclose%
\pgfusepath{stroke,fill}%
}%
\begin{pgfscope}%
\pgfsys@transformshift{4.006787in}{3.363810in}%
\pgfsys@useobject{currentmarker}{}%
\end{pgfscope}%
\end{pgfscope}%
\begin{pgfscope}%
\definecolor{textcolor}{rgb}{0.000000,0.000000,0.000000}%
\pgfsetstrokecolor{textcolor}%
\pgfsetfillcolor{textcolor}%
\pgftext[x=4.381787in,y=3.276310in,left,base]{\color{textcolor}\rmfamily\fontsize{18.000000}{21.600000}\selectfont LPP\(\displaystyle ^{T}\)}%
\end{pgfscope}%
\begin{pgfscope}%
\pgfsetbuttcap%
\pgfsetroundjoin%
\pgfsetlinewidth{1.505625pt}%
\definecolor{currentstroke}{rgb}{0.839216,0.152941,0.156863}%
\pgfsetstrokecolor{currentstroke}%
\pgfsetstrokeopacity{0.800000}%
\pgfsetdash{{1.500000pt}{2.475000pt}}{0.000000pt}%
\pgfpathmoveto{\pgfqpoint{5.146819in}{3.363810in}}%
\pgfpathlineto{\pgfqpoint{5.396819in}{3.363810in}}%
\pgfpathlineto{\pgfqpoint{5.646819in}{3.363810in}}%
\pgfusepath{stroke}%
\end{pgfscope}%
\begin{pgfscope}%
\pgfsetbuttcap%
\pgfsetroundjoin%
\definecolor{currentfill}{rgb}{0.839216,0.152941,0.156863}%
\pgfsetfillcolor{currentfill}%
\pgfsetfillopacity{0.800000}%
\pgfsetlinewidth{1.003750pt}%
\definecolor{currentstroke}{rgb}{0.839216,0.152941,0.156863}%
\pgfsetstrokecolor{currentstroke}%
\pgfsetstrokeopacity{0.800000}%
\pgfsetdash{}{0pt}%
\pgfsys@defobject{currentmarker}{\pgfqpoint{-0.041667in}{-0.041667in}}{\pgfqpoint{0.041667in}{0.041667in}}{%
\pgfpathmoveto{\pgfqpoint{0.000000in}{-0.041667in}}%
\pgfpathcurveto{\pgfqpoint{0.011050in}{-0.041667in}}{\pgfqpoint{0.021649in}{-0.037276in}}{\pgfqpoint{0.029463in}{-0.029463in}}%
\pgfpathcurveto{\pgfqpoint{0.037276in}{-0.021649in}}{\pgfqpoint{0.041667in}{-0.011050in}}{\pgfqpoint{0.041667in}{0.000000in}}%
\pgfpathcurveto{\pgfqpoint{0.041667in}{0.011050in}}{\pgfqpoint{0.037276in}{0.021649in}}{\pgfqpoint{0.029463in}{0.029463in}}%
\pgfpathcurveto{\pgfqpoint{0.021649in}{0.037276in}}{\pgfqpoint{0.011050in}{0.041667in}}{\pgfqpoint{0.000000in}{0.041667in}}%
\pgfpathcurveto{\pgfqpoint{-0.011050in}{0.041667in}}{\pgfqpoint{-0.021649in}{0.037276in}}{\pgfqpoint{-0.029463in}{0.029463in}}%
\pgfpathcurveto{\pgfqpoint{-0.037276in}{0.021649in}}{\pgfqpoint{-0.041667in}{0.011050in}}{\pgfqpoint{-0.041667in}{0.000000in}}%
\pgfpathcurveto{\pgfqpoint{-0.041667in}{-0.011050in}}{\pgfqpoint{-0.037276in}{-0.021649in}}{\pgfqpoint{-0.029463in}{-0.029463in}}%
\pgfpathcurveto{\pgfqpoint{-0.021649in}{-0.037276in}}{\pgfqpoint{-0.011050in}{-0.041667in}}{\pgfqpoint{0.000000in}{-0.041667in}}%
\pgfpathlineto{\pgfqpoint{0.000000in}{-0.041667in}}%
\pgfpathclose%
\pgfusepath{stroke,fill}%
}%
\begin{pgfscope}%
\pgfsys@transformshift{5.396819in}{3.363810in}%
\pgfsys@useobject{currentmarker}{}%
\end{pgfscope}%
\end{pgfscope}%
\begin{pgfscope}%
\definecolor{textcolor}{rgb}{0.000000,0.000000,0.000000}%
\pgfsetstrokecolor{textcolor}%
\pgfsetfillcolor{textcolor}%
\pgftext[x=5.771819in,y=3.276310in,left,base]{\color{textcolor}\rmfamily\fontsize{18.000000}{21.600000}\selectfont LPP\(\displaystyle ^{T, P}\)}%
\end{pgfscope}%
\end{pgfpicture}%
\makeatother%
\endgroup%

%% file: figures/DemandPlots/TightLayout/DE_avgWTT_rescaled_Gap10pct.pgf
\begingroup%
\makeatletter%
\begin{pgfpicture}%
\pgfpathrectangle{\pgfpointorigin}{\pgfqpoint{7.950645in}{3.651310in}}%
\pgfusepath{use as bounding box, clip}%
\begin{pgfscope}%
\pgfsetbuttcap%
\pgfsetmiterjoin%
\definecolor{currentfill}{rgb}{1.000000,1.000000,1.000000}%
\pgfsetfillcolor{currentfill}%
\pgfsetlinewidth{0.000000pt}%
\definecolor{currentstroke}{rgb}{1.000000,1.000000,1.000000}%
\pgfsetstrokecolor{currentstroke}%
\pgfsetdash{}{0pt}%
\pgfpathmoveto{\pgfqpoint{0.000000in}{0.000000in}}%
\pgfpathlineto{\pgfqpoint{7.950645in}{0.000000in}}%
\pgfpathlineto{\pgfqpoint{7.950645in}{3.651310in}}%
\pgfpathlineto{\pgfqpoint{0.000000in}{3.651310in}}%
\pgfpathlineto{\pgfqpoint{0.000000in}{0.000000in}}%
\pgfpathclose%
\pgfusepath{fill}%
\end{pgfscope}%
\begin{pgfscope}%
\pgfsetbuttcap%
\pgfsetmiterjoin%
\definecolor{currentfill}{rgb}{1.000000,1.000000,1.000000}%
\pgfsetfillcolor{currentfill}%
\pgfsetlinewidth{0.000000pt}%
\definecolor{currentstroke}{rgb}{0.000000,0.000000,0.000000}%
\pgfsetstrokecolor{currentstroke}%
\pgfsetstrokeopacity{0.000000}%
\pgfsetdash{}{0pt}%
\pgfpathmoveto{\pgfqpoint{0.861608in}{0.679475in}}%
\pgfpathlineto{\pgfqpoint{7.850645in}{0.679475in}}%
\pgfpathlineto{\pgfqpoint{7.850645in}{3.135897in}}%
\pgfpathlineto{\pgfqpoint{0.861608in}{3.135897in}}%
\pgfpathlineto{\pgfqpoint{0.861608in}{0.679475in}}%
\pgfpathclose%
\pgfusepath{fill}%
\end{pgfscope}%
\begin{pgfscope}%
\pgfsetbuttcap%
\pgfsetroundjoin%
\definecolor{currentfill}{rgb}{0.000000,0.000000,0.000000}%
\pgfsetfillcolor{currentfill}%
\pgfsetlinewidth{0.803000pt}%
\definecolor{currentstroke}{rgb}{0.000000,0.000000,0.000000}%
\pgfsetstrokecolor{currentstroke}%
\pgfsetdash{}{0pt}%
\pgfsys@defobject{currentmarker}{\pgfqpoint{0.000000in}{-0.048611in}}{\pgfqpoint{0.000000in}{0.000000in}}{%
\pgfpathmoveto{\pgfqpoint{0.000000in}{0.000000in}}%
\pgfpathlineto{\pgfqpoint{0.000000in}{-0.048611in}}%
\pgfusepath{stroke,fill}%
}%
\begin{pgfscope}%
\pgfsys@transformshift{0.861608in}{0.679475in}%
\pgfsys@useobject{currentmarker}{}%
\end{pgfscope}%
\end{pgfscope}%
\begin{pgfscope}%
\definecolor{textcolor}{rgb}{0.000000,0.000000,0.000000}%
\pgfsetstrokecolor{textcolor}%
\pgfsetfillcolor{textcolor}%
\pgftext[x=0.861608in,y=0.582253in,,top]{\color{textcolor}\rmfamily\fontsize{16.000000}{19.200000}\selectfont \(\displaystyle {0.0}\)}%
\end{pgfscope}%
\begin{pgfscope}%
\pgfsetbuttcap%
\pgfsetroundjoin%
\definecolor{currentfill}{rgb}{0.000000,0.000000,0.000000}%
\pgfsetfillcolor{currentfill}%
\pgfsetlinewidth{0.803000pt}%
\definecolor{currentstroke}{rgb}{0.000000,0.000000,0.000000}%
\pgfsetstrokecolor{currentstroke}%
\pgfsetdash{}{0pt}%
\pgfsys@defobject{currentmarker}{\pgfqpoint{0.000000in}{-0.048611in}}{\pgfqpoint{0.000000in}{0.000000in}}{%
\pgfpathmoveto{\pgfqpoint{0.000000in}{0.000000in}}%
\pgfpathlineto{\pgfqpoint{0.000000in}{-0.048611in}}%
\pgfusepath{stroke,fill}%
}%
\begin{pgfscope}%
\pgfsys@transformshift{2.132342in}{0.679475in}%
\pgfsys@useobject{currentmarker}{}%
\end{pgfscope}%
\end{pgfscope}%
\begin{pgfscope}%
\definecolor{textcolor}{rgb}{0.000000,0.000000,0.000000}%
\pgfsetstrokecolor{textcolor}%
\pgfsetfillcolor{textcolor}%
\pgftext[x=2.132342in,y=0.582253in,,top]{\color{textcolor}\rmfamily\fontsize{16.000000}{19.200000}\selectfont \(\displaystyle {0.2}\)}%
\end{pgfscope}%
\begin{pgfscope}%
\pgfsetbuttcap%
\pgfsetroundjoin%
\definecolor{currentfill}{rgb}{0.000000,0.000000,0.000000}%
\pgfsetfillcolor{currentfill}%
\pgfsetlinewidth{0.803000pt}%
\definecolor{currentstroke}{rgb}{0.000000,0.000000,0.000000}%
\pgfsetstrokecolor{currentstroke}%
\pgfsetdash{}{0pt}%
\pgfsys@defobject{currentmarker}{\pgfqpoint{0.000000in}{-0.048611in}}{\pgfqpoint{0.000000in}{0.000000in}}{%
\pgfpathmoveto{\pgfqpoint{0.000000in}{0.000000in}}%
\pgfpathlineto{\pgfqpoint{0.000000in}{-0.048611in}}%
\pgfusepath{stroke,fill}%
}%
\begin{pgfscope}%
\pgfsys@transformshift{3.403076in}{0.679475in}%
\pgfsys@useobject{currentmarker}{}%
\end{pgfscope}%
\end{pgfscope}%
\begin{pgfscope}%
\definecolor{textcolor}{rgb}{0.000000,0.000000,0.000000}%
\pgfsetstrokecolor{textcolor}%
\pgfsetfillcolor{textcolor}%
\pgftext[x=3.403076in,y=0.582253in,,top]{\color{textcolor}\rmfamily\fontsize{16.000000}{19.200000}\selectfont \(\displaystyle {0.4}\)}%
\end{pgfscope}%
\begin{pgfscope}%
\pgfsetbuttcap%
\pgfsetroundjoin%
\definecolor{currentfill}{rgb}{0.000000,0.000000,0.000000}%
\pgfsetfillcolor{currentfill}%
\pgfsetlinewidth{0.803000pt}%
\definecolor{currentstroke}{rgb}{0.000000,0.000000,0.000000}%
\pgfsetstrokecolor{currentstroke}%
\pgfsetdash{}{0pt}%
\pgfsys@defobject{currentmarker}{\pgfqpoint{0.000000in}{-0.048611in}}{\pgfqpoint{0.000000in}{0.000000in}}{%
\pgfpathmoveto{\pgfqpoint{0.000000in}{0.000000in}}%
\pgfpathlineto{\pgfqpoint{0.000000in}{-0.048611in}}%
\pgfusepath{stroke,fill}%
}%
\begin{pgfscope}%
\pgfsys@transformshift{4.673810in}{0.679475in}%
\pgfsys@useobject{currentmarker}{}%
\end{pgfscope}%
\end{pgfscope}%
\begin{pgfscope}%
\definecolor{textcolor}{rgb}{0.000000,0.000000,0.000000}%
\pgfsetstrokecolor{textcolor}%
\pgfsetfillcolor{textcolor}%
\pgftext[x=4.673810in,y=0.582253in,,top]{\color{textcolor}\rmfamily\fontsize{16.000000}{19.200000}\selectfont \(\displaystyle {0.6}\)}%
\end{pgfscope}%
\begin{pgfscope}%
\pgfsetbuttcap%
\pgfsetroundjoin%
\definecolor{currentfill}{rgb}{0.000000,0.000000,0.000000}%
\pgfsetfillcolor{currentfill}%
\pgfsetlinewidth{0.803000pt}%
\definecolor{currentstroke}{rgb}{0.000000,0.000000,0.000000}%
\pgfsetstrokecolor{currentstroke}%
\pgfsetdash{}{0pt}%
\pgfsys@defobject{currentmarker}{\pgfqpoint{0.000000in}{-0.048611in}}{\pgfqpoint{0.000000in}{0.000000in}}{%
\pgfpathmoveto{\pgfqpoint{0.000000in}{0.000000in}}%
\pgfpathlineto{\pgfqpoint{0.000000in}{-0.048611in}}%
\pgfusepath{stroke,fill}%
}%
\begin{pgfscope}%
\pgfsys@transformshift{5.944544in}{0.679475in}%
\pgfsys@useobject{currentmarker}{}%
\end{pgfscope}%
\end{pgfscope}%
\begin{pgfscope}%
\definecolor{textcolor}{rgb}{0.000000,0.000000,0.000000}%
\pgfsetstrokecolor{textcolor}%
\pgfsetfillcolor{textcolor}%
\pgftext[x=5.944544in,y=0.582253in,,top]{\color{textcolor}\rmfamily\fontsize{16.000000}{19.200000}\selectfont \(\displaystyle {0.8}\)}%
\end{pgfscope}%
\begin{pgfscope}%
\pgfsetbuttcap%
\pgfsetroundjoin%
\definecolor{currentfill}{rgb}{0.000000,0.000000,0.000000}%
\pgfsetfillcolor{currentfill}%
\pgfsetlinewidth{0.803000pt}%
\definecolor{currentstroke}{rgb}{0.000000,0.000000,0.000000}%
\pgfsetstrokecolor{currentstroke}%
\pgfsetdash{}{0pt}%
\pgfsys@defobject{currentmarker}{\pgfqpoint{0.000000in}{-0.048611in}}{\pgfqpoint{0.000000in}{0.000000in}}{%
\pgfpathmoveto{\pgfqpoint{0.000000in}{0.000000in}}%
\pgfpathlineto{\pgfqpoint{0.000000in}{-0.048611in}}%
\pgfusepath{stroke,fill}%
}%
\begin{pgfscope}%
\pgfsys@transformshift{7.215278in}{0.679475in}%
\pgfsys@useobject{currentmarker}{}%
\end{pgfscope}%
\end{pgfscope}%
\begin{pgfscope}%
\definecolor{textcolor}{rgb}{0.000000,0.000000,0.000000}%
\pgfsetstrokecolor{textcolor}%
\pgfsetfillcolor{textcolor}%
\pgftext[x=7.215278in,y=0.582253in,,top]{\color{textcolor}\rmfamily\fontsize{16.000000}{19.200000}\selectfont \(\displaystyle {1.0}\)}%
\end{pgfscope}%
\begin{pgfscope}%
\definecolor{textcolor}{rgb}{0.000000,0.000000,0.000000}%
\pgfsetstrokecolor{textcolor}%
\pgfsetfillcolor{textcolor}%
\pgftext[x=4.356126in,y=0.313349in,,top]{\color{textcolor}\rmfamily\fontsize{18.000000}{21.600000}\selectfont \(\displaystyle \lambda\)}%
\end{pgfscope}%
\begin{pgfscope}%
\pgfsetbuttcap%
\pgfsetroundjoin%
\definecolor{currentfill}{rgb}{0.000000,0.000000,0.000000}%
\pgfsetfillcolor{currentfill}%
\pgfsetlinewidth{0.803000pt}%
\definecolor{currentstroke}{rgb}{0.000000,0.000000,0.000000}%
\pgfsetstrokecolor{currentstroke}%
\pgfsetdash{}{0pt}%
\pgfsys@defobject{currentmarker}{\pgfqpoint{-0.048611in}{0.000000in}}{\pgfqpoint{-0.000000in}{0.000000in}}{%
\pgfpathmoveto{\pgfqpoint{-0.000000in}{0.000000in}}%
\pgfpathlineto{\pgfqpoint{-0.048611in}{0.000000in}}%
\pgfusepath{stroke,fill}%
}%
\begin{pgfscope}%
\pgfsys@transformshift{0.861608in}{0.881292in}%
\pgfsys@useobject{currentmarker}{}%
\end{pgfscope}%
\end{pgfscope}%
\begin{pgfscope}%
\definecolor{textcolor}{rgb}{0.000000,0.000000,0.000000}%
\pgfsetstrokecolor{textcolor}%
\pgfsetfillcolor{textcolor}%
\pgftext[x=0.368904in, y=0.797959in, left, base]{\color{textcolor}\rmfamily\fontsize{16.000000}{19.200000}\selectfont \(\displaystyle {22.5}\)}%
\end{pgfscope}%
\begin{pgfscope}%
\pgfsetbuttcap%
\pgfsetroundjoin%
\definecolor{currentfill}{rgb}{0.000000,0.000000,0.000000}%
\pgfsetfillcolor{currentfill}%
\pgfsetlinewidth{0.803000pt}%
\definecolor{currentstroke}{rgb}{0.000000,0.000000,0.000000}%
\pgfsetstrokecolor{currentstroke}%
\pgfsetdash{}{0pt}%
\pgfsys@defobject{currentmarker}{\pgfqpoint{-0.048611in}{0.000000in}}{\pgfqpoint{-0.000000in}{0.000000in}}{%
\pgfpathmoveto{\pgfqpoint{-0.000000in}{0.000000in}}%
\pgfpathlineto{\pgfqpoint{-0.048611in}{0.000000in}}%
\pgfusepath{stroke,fill}%
}%
\begin{pgfscope}%
\pgfsys@transformshift{0.861608in}{1.370945in}%
\pgfsys@useobject{currentmarker}{}%
\end{pgfscope}%
\end{pgfscope}%
\begin{pgfscope}%
\definecolor{textcolor}{rgb}{0.000000,0.000000,0.000000}%
\pgfsetstrokecolor{textcolor}%
\pgfsetfillcolor{textcolor}%
\pgftext[x=0.368904in, y=1.287611in, left, base]{\color{textcolor}\rmfamily\fontsize{16.000000}{19.200000}\selectfont \(\displaystyle {25.0}\)}%
\end{pgfscope}%
\begin{pgfscope}%
\pgfsetbuttcap%
\pgfsetroundjoin%
\definecolor{currentfill}{rgb}{0.000000,0.000000,0.000000}%
\pgfsetfillcolor{currentfill}%
\pgfsetlinewidth{0.803000pt}%
\definecolor{currentstroke}{rgb}{0.000000,0.000000,0.000000}%
\pgfsetstrokecolor{currentstroke}%
\pgfsetdash{}{0pt}%
\pgfsys@defobject{currentmarker}{\pgfqpoint{-0.048611in}{0.000000in}}{\pgfqpoint{-0.000000in}{0.000000in}}{%
\pgfpathmoveto{\pgfqpoint{-0.000000in}{0.000000in}}%
\pgfpathlineto{\pgfqpoint{-0.048611in}{0.000000in}}%
\pgfusepath{stroke,fill}%
}%
\begin{pgfscope}%
\pgfsys@transformshift{0.861608in}{1.860598in}%
\pgfsys@useobject{currentmarker}{}%
\end{pgfscope}%
\end{pgfscope}%
\begin{pgfscope}%
\definecolor{textcolor}{rgb}{0.000000,0.000000,0.000000}%
\pgfsetstrokecolor{textcolor}%
\pgfsetfillcolor{textcolor}%
\pgftext[x=0.368904in, y=1.777264in, left, base]{\color{textcolor}\rmfamily\fontsize{16.000000}{19.200000}\selectfont \(\displaystyle {27.5}\)}%
\end{pgfscope}%
\begin{pgfscope}%
\pgfsetbuttcap%
\pgfsetroundjoin%
\definecolor{currentfill}{rgb}{0.000000,0.000000,0.000000}%
\pgfsetfillcolor{currentfill}%
\pgfsetlinewidth{0.803000pt}%
\definecolor{currentstroke}{rgb}{0.000000,0.000000,0.000000}%
\pgfsetstrokecolor{currentstroke}%
\pgfsetdash{}{0pt}%
\pgfsys@defobject{currentmarker}{\pgfqpoint{-0.048611in}{0.000000in}}{\pgfqpoint{-0.000000in}{0.000000in}}{%
\pgfpathmoveto{\pgfqpoint{-0.000000in}{0.000000in}}%
\pgfpathlineto{\pgfqpoint{-0.048611in}{0.000000in}}%
\pgfusepath{stroke,fill}%
}%
\begin{pgfscope}%
\pgfsys@transformshift{0.861608in}{2.350250in}%
\pgfsys@useobject{currentmarker}{}%
\end{pgfscope}%
\end{pgfscope}%
\begin{pgfscope}%
\definecolor{textcolor}{rgb}{0.000000,0.000000,0.000000}%
\pgfsetstrokecolor{textcolor}%
\pgfsetfillcolor{textcolor}%
\pgftext[x=0.368904in, y=2.266917in, left, base]{\color{textcolor}\rmfamily\fontsize{16.000000}{19.200000}\selectfont \(\displaystyle {30.0}\)}%
\end{pgfscope}%
\begin{pgfscope}%
\pgfsetbuttcap%
\pgfsetroundjoin%
\definecolor{currentfill}{rgb}{0.000000,0.000000,0.000000}%
\pgfsetfillcolor{currentfill}%
\pgfsetlinewidth{0.803000pt}%
\definecolor{currentstroke}{rgb}{0.000000,0.000000,0.000000}%
\pgfsetstrokecolor{currentstroke}%
\pgfsetdash{}{0pt}%
\pgfsys@defobject{currentmarker}{\pgfqpoint{-0.048611in}{0.000000in}}{\pgfqpoint{-0.000000in}{0.000000in}}{%
\pgfpathmoveto{\pgfqpoint{-0.000000in}{0.000000in}}%
\pgfpathlineto{\pgfqpoint{-0.048611in}{0.000000in}}%
\pgfusepath{stroke,fill}%
}%
\begin{pgfscope}%
\pgfsys@transformshift{0.861608in}{2.839903in}%
\pgfsys@useobject{currentmarker}{}%
\end{pgfscope}%
\end{pgfscope}%
\begin{pgfscope}%
\definecolor{textcolor}{rgb}{0.000000,0.000000,0.000000}%
\pgfsetstrokecolor{textcolor}%
\pgfsetfillcolor{textcolor}%
\pgftext[x=0.368904in, y=2.756570in, left, base]{\color{textcolor}\rmfamily\fontsize{16.000000}{19.200000}\selectfont \(\displaystyle {32.5}\)}%
\end{pgfscope}%
\begin{pgfscope}%
\definecolor{textcolor}{rgb}{0.000000,0.000000,0.000000}%
\pgfsetstrokecolor{textcolor}%
\pgfsetfillcolor{textcolor}%
\pgftext[x=0.313349in,y=1.907686in,,bottom,rotate=90.000000]{\color{textcolor}\rmfamily\fontsize{18.000000}{21.600000}\selectfont minutes}%
\end{pgfscope}%
\begin{pgfscope}%
\pgfpathrectangle{\pgfqpoint{0.861608in}{0.679475in}}{\pgfqpoint{6.989037in}{2.456422in}}%
\pgfusepath{clip}%
\pgfsetrectcap%
\pgfsetroundjoin%
\pgfsetlinewidth{1.505625pt}%
\definecolor{currentstroke}{rgb}{1.000000,0.498039,0.054902}%
\pgfsetstrokecolor{currentstroke}%
\pgfsetdash{}{0pt}%
\pgfpathmoveto{\pgfqpoint{1.496975in}{3.024241in}}%
\pgfpathlineto{\pgfqpoint{2.450025in}{2.685336in}}%
\pgfpathlineto{\pgfqpoint{4.038443in}{2.353074in}}%
\pgfpathlineto{\pgfqpoint{5.626860in}{1.943079in}}%
\pgfpathlineto{\pgfqpoint{7.215278in}{1.689985in}}%
\pgfusepath{stroke}%
\end{pgfscope}%
\begin{pgfscope}%
\pgfpathrectangle{\pgfqpoint{0.861608in}{0.679475in}}{\pgfqpoint{6.989037in}{2.456422in}}%
\pgfusepath{clip}%
\pgfsetbuttcap%
\pgfsetmiterjoin%
\definecolor{currentfill}{rgb}{1.000000,0.498039,0.054902}%
\pgfsetfillcolor{currentfill}%
\pgfsetlinewidth{0.752812pt}%
\definecolor{currentstroke}{rgb}{1.000000,1.000000,1.000000}%
\pgfsetstrokecolor{currentstroke}%
\pgfsetdash{}{0pt}%
\pgfsys@defobject{currentmarker}{\pgfqpoint{-0.041667in}{-0.041667in}}{\pgfqpoint{0.041667in}{0.041667in}}{%
\pgfpathmoveto{\pgfqpoint{-0.041667in}{-0.041667in}}%
\pgfpathlineto{\pgfqpoint{0.041667in}{-0.041667in}}%
\pgfpathlineto{\pgfqpoint{0.041667in}{0.041667in}}%
\pgfpathlineto{\pgfqpoint{-0.041667in}{0.041667in}}%
\pgfpathlineto{\pgfqpoint{-0.041667in}{-0.041667in}}%
\pgfpathclose%
\pgfusepath{stroke,fill}%
}%
\begin{pgfscope}%
\pgfsys@transformshift{1.496975in}{3.024241in}%
\pgfsys@useobject{currentmarker}{}%
\end{pgfscope}%
\begin{pgfscope}%
\pgfsys@transformshift{2.450025in}{2.685336in}%
\pgfsys@useobject{currentmarker}{}%
\end{pgfscope}%
\begin{pgfscope}%
\pgfsys@transformshift{4.038443in}{2.353074in}%
\pgfsys@useobject{currentmarker}{}%
\end{pgfscope}%
\begin{pgfscope}%
\pgfsys@transformshift{5.626860in}{1.943079in}%
\pgfsys@useobject{currentmarker}{}%
\end{pgfscope}%
\begin{pgfscope}%
\pgfsys@transformshift{7.215278in}{1.689985in}%
\pgfsys@useobject{currentmarker}{}%
\end{pgfscope}%
\end{pgfscope}%
\begin{pgfscope}%
\pgfpathrectangle{\pgfqpoint{0.861608in}{0.679475in}}{\pgfqpoint{6.989037in}{2.456422in}}%
\pgfusepath{clip}%
\pgfsetbuttcap%
\pgfsetroundjoin%
\pgfsetlinewidth{1.505625pt}%
\definecolor{currentstroke}{rgb}{1.000000,0.498039,0.054902}%
\pgfsetstrokecolor{currentstroke}%
\pgfsetstrokeopacity{0.900000}%
\pgfsetdash{{7.500000pt}{7.500000pt}}{0.000000pt}%
\pgfpathmoveto{\pgfqpoint{1.496975in}{1.860369in}}%
\pgfpathlineto{\pgfqpoint{2.450025in}{1.832759in}}%
\pgfpathlineto{\pgfqpoint{4.038443in}{1.893662in}}%
\pgfpathlineto{\pgfqpoint{5.626860in}{1.494731in}}%
\pgfpathlineto{\pgfqpoint{7.215278in}{1.132437in}}%
\pgfusepath{stroke}%
\end{pgfscope}%
\begin{pgfscope}%
\pgfpathrectangle{\pgfqpoint{0.861608in}{0.679475in}}{\pgfqpoint{6.989037in}{2.456422in}}%
\pgfusepath{clip}%
\pgfsetbuttcap%
\pgfsetmiterjoin%
\definecolor{currentfill}{rgb}{1.000000,0.498039,0.054902}%
\pgfsetfillcolor{currentfill}%
\pgfsetfillopacity{0.900000}%
\pgfsetlinewidth{0.752812pt}%
\definecolor{currentstroke}{rgb}{1.000000,1.000000,1.000000}%
\pgfsetstrokecolor{currentstroke}%
\pgfsetdash{}{0pt}%
\pgfsys@defobject{currentmarker}{\pgfqpoint{-0.058926in}{-0.058926in}}{\pgfqpoint{0.058926in}{0.058926in}}{%
\pgfpathmoveto{\pgfqpoint{0.000000in}{-0.058926in}}%
\pgfpathlineto{\pgfqpoint{0.058926in}{0.000000in}}%
\pgfpathlineto{\pgfqpoint{0.000000in}{0.058926in}}%
\pgfpathlineto{\pgfqpoint{-0.058926in}{0.000000in}}%
\pgfpathlineto{\pgfqpoint{0.000000in}{-0.058926in}}%
\pgfpathclose%
\pgfusepath{stroke,fill}%
}%
\begin{pgfscope}%
\pgfsys@transformshift{1.496975in}{1.860369in}%
\pgfsys@useobject{currentmarker}{}%
\end{pgfscope}%
\begin{pgfscope}%
\pgfsys@transformshift{2.450025in}{1.832759in}%
\pgfsys@useobject{currentmarker}{}%
\end{pgfscope}%
\begin{pgfscope}%
\pgfsys@transformshift{4.038443in}{1.893662in}%
\pgfsys@useobject{currentmarker}{}%
\end{pgfscope}%
\begin{pgfscope}%
\pgfsys@transformshift{5.626860in}{1.494731in}%
\pgfsys@useobject{currentmarker}{}%
\end{pgfscope}%
\begin{pgfscope}%
\pgfsys@transformshift{7.215278in}{1.132437in}%
\pgfsys@useobject{currentmarker}{}%
\end{pgfscope}%
\end{pgfscope}%
\begin{pgfscope}%
\pgfpathrectangle{\pgfqpoint{0.861608in}{0.679475in}}{\pgfqpoint{6.989037in}{2.456422in}}%
\pgfusepath{clip}%
\pgfsetbuttcap%
\pgfsetroundjoin%
\pgfsetlinewidth{1.505625pt}%
\definecolor{currentstroke}{rgb}{1.000000,0.498039,0.054902}%
\pgfsetstrokecolor{currentstroke}%
\pgfsetstrokeopacity{0.800000}%
\pgfsetdash{{1.500000pt}{4.500000pt}}{0.000000pt}%
\pgfpathmoveto{\pgfqpoint{1.496975in}{0.866700in}}%
\pgfpathlineto{\pgfqpoint{2.450025in}{0.791131in}}%
\pgfpathlineto{\pgfqpoint{4.038443in}{0.838016in}}%
\pgfpathlineto{\pgfqpoint{5.626860in}{0.867926in}}%
\pgfpathlineto{\pgfqpoint{7.215278in}{0.795193in}}%
\pgfusepath{stroke}%
\end{pgfscope}%
\begin{pgfscope}%
\pgfpathrectangle{\pgfqpoint{0.861608in}{0.679475in}}{\pgfqpoint{6.989037in}{2.456422in}}%
\pgfusepath{clip}%
\pgfsetbuttcap%
\pgfsetroundjoin%
\definecolor{currentfill}{rgb}{1.000000,0.498039,0.054902}%
\pgfsetfillcolor{currentfill}%
\pgfsetfillopacity{0.800000}%
\pgfsetlinewidth{0.752812pt}%
\definecolor{currentstroke}{rgb}{1.000000,1.000000,1.000000}%
\pgfsetstrokecolor{currentstroke}%
\pgfsetdash{}{0pt}%
\pgfsys@defobject{currentmarker}{\pgfqpoint{-0.041667in}{-0.041667in}}{\pgfqpoint{0.041667in}{0.041667in}}{%
\pgfpathmoveto{\pgfqpoint{0.000000in}{-0.041667in}}%
\pgfpathcurveto{\pgfqpoint{0.011050in}{-0.041667in}}{\pgfqpoint{0.021649in}{-0.037276in}}{\pgfqpoint{0.029463in}{-0.029463in}}%
\pgfpathcurveto{\pgfqpoint{0.037276in}{-0.021649in}}{\pgfqpoint{0.041667in}{-0.011050in}}{\pgfqpoint{0.041667in}{0.000000in}}%
\pgfpathcurveto{\pgfqpoint{0.041667in}{0.011050in}}{\pgfqpoint{0.037276in}{0.021649in}}{\pgfqpoint{0.029463in}{0.029463in}}%
\pgfpathcurveto{\pgfqpoint{0.021649in}{0.037276in}}{\pgfqpoint{0.011050in}{0.041667in}}{\pgfqpoint{0.000000in}{0.041667in}}%
\pgfpathcurveto{\pgfqpoint{-0.011050in}{0.041667in}}{\pgfqpoint{-0.021649in}{0.037276in}}{\pgfqpoint{-0.029463in}{0.029463in}}%
\pgfpathcurveto{\pgfqpoint{-0.037276in}{0.021649in}}{\pgfqpoint{-0.041667in}{0.011050in}}{\pgfqpoint{-0.041667in}{0.000000in}}%
\pgfpathcurveto{\pgfqpoint{-0.041667in}{-0.011050in}}{\pgfqpoint{-0.037276in}{-0.021649in}}{\pgfqpoint{-0.029463in}{-0.029463in}}%
\pgfpathcurveto{\pgfqpoint{-0.021649in}{-0.037276in}}{\pgfqpoint{-0.011050in}{-0.041667in}}{\pgfqpoint{0.000000in}{-0.041667in}}%
\pgfpathlineto{\pgfqpoint{0.000000in}{-0.041667in}}%
\pgfpathclose%
\pgfusepath{stroke,fill}%
}%
\begin{pgfscope}%
\pgfsys@transformshift{1.496975in}{0.866700in}%
\pgfsys@useobject{currentmarker}{}%
\end{pgfscope}%
\begin{pgfscope}%
\pgfsys@transformshift{2.450025in}{0.791131in}%
\pgfsys@useobject{currentmarker}{}%
\end{pgfscope}%
\begin{pgfscope}%
\pgfsys@transformshift{4.038443in}{0.838016in}%
\pgfsys@useobject{currentmarker}{}%
\end{pgfscope}%
\begin{pgfscope}%
\pgfsys@transformshift{5.626860in}{0.867926in}%
\pgfsys@useobject{currentmarker}{}%
\end{pgfscope}%
\begin{pgfscope}%
\pgfsys@transformshift{7.215278in}{0.795193in}%
\pgfsys@useobject{currentmarker}{}%
\end{pgfscope}%
\end{pgfscope}%
\begin{pgfscope}%
\pgfsetrectcap%
\pgfsetmiterjoin%
\pgfsetlinewidth{0.803000pt}%
\definecolor{currentstroke}{rgb}{0.000000,0.000000,0.000000}%
\pgfsetstrokecolor{currentstroke}%
\pgfsetdash{}{0pt}%
\pgfpathmoveto{\pgfqpoint{0.861608in}{0.679475in}}%
\pgfpathlineto{\pgfqpoint{0.861608in}{3.135897in}}%
\pgfusepath{stroke}%
\end{pgfscope}%
\begin{pgfscope}%
\pgfsetrectcap%
\pgfsetmiterjoin%
\pgfsetlinewidth{0.803000pt}%
\definecolor{currentstroke}{rgb}{0.000000,0.000000,0.000000}%
\pgfsetstrokecolor{currentstroke}%
\pgfsetdash{}{0pt}%
\pgfpathmoveto{\pgfqpoint{7.850645in}{0.679475in}}%
\pgfpathlineto{\pgfqpoint{7.850645in}{3.135897in}}%
\pgfusepath{stroke}%
\end{pgfscope}%
\begin{pgfscope}%
\pgfsetrectcap%
\pgfsetmiterjoin%
\pgfsetlinewidth{0.803000pt}%
\definecolor{currentstroke}{rgb}{0.000000,0.000000,0.000000}%
\pgfsetstrokecolor{currentstroke}%
\pgfsetdash{}{0pt}%
\pgfpathmoveto{\pgfqpoint{0.861608in}{0.679475in}}%
\pgfpathlineto{\pgfqpoint{7.850645in}{0.679475in}}%
\pgfusepath{stroke}%
\end{pgfscope}%
\begin{pgfscope}%
\pgfsetrectcap%
\pgfsetmiterjoin%
\pgfsetlinewidth{0.803000pt}%
\definecolor{currentstroke}{rgb}{0.000000,0.000000,0.000000}%
\pgfsetstrokecolor{currentstroke}%
\pgfsetdash{}{0pt}%
\pgfpathmoveto{\pgfqpoint{0.861608in}{3.135897in}}%
\pgfpathlineto{\pgfqpoint{7.850645in}{3.135897in}}%
\pgfusepath{stroke}%
\end{pgfscope}%
\begin{pgfscope}%
\pgfsetrectcap%
\pgfsetroundjoin%
\pgfsetlinewidth{1.505625pt}%
\definecolor{currentstroke}{rgb}{1.000000,0.498039,0.054902}%
\pgfsetstrokecolor{currentstroke}%
\pgfsetdash{}{0pt}%
\pgfpathmoveto{\pgfqpoint{2.203560in}{3.363810in}}%
\pgfpathlineto{\pgfqpoint{2.453560in}{3.363810in}}%
\pgfpathlineto{\pgfqpoint{2.703560in}{3.363810in}}%
\pgfusepath{stroke}%
\end{pgfscope}%
\begin{pgfscope}%
\pgfsetbuttcap%
\pgfsetmiterjoin%
\definecolor{currentfill}{rgb}{1.000000,0.498039,0.054902}%
\pgfsetfillcolor{currentfill}%
\pgfsetlinewidth{1.003750pt}%
\definecolor{currentstroke}{rgb}{1.000000,0.498039,0.054902}%
\pgfsetstrokecolor{currentstroke}%
\pgfsetdash{}{0pt}%
\pgfsys@defobject{currentmarker}{\pgfqpoint{-0.041667in}{-0.041667in}}{\pgfqpoint{0.041667in}{0.041667in}}{%
\pgfpathmoveto{\pgfqpoint{-0.041667in}{-0.041667in}}%
\pgfpathlineto{\pgfqpoint{0.041667in}{-0.041667in}}%
\pgfpathlineto{\pgfqpoint{0.041667in}{0.041667in}}%
\pgfpathlineto{\pgfqpoint{-0.041667in}{0.041667in}}%
\pgfpathlineto{\pgfqpoint{-0.041667in}{-0.041667in}}%
\pgfpathclose%
\pgfusepath{stroke,fill}%
}%
\begin{pgfscope}%
\pgfsys@transformshift{2.453560in}{3.363810in}%
\pgfsys@useobject{currentmarker}{}%
\end{pgfscope}%
\end{pgfscope}%
\begin{pgfscope}%
\definecolor{textcolor}{rgb}{0.000000,0.000000,0.000000}%
\pgfsetstrokecolor{textcolor}%
\pgfsetfillcolor{textcolor}%
\pgftext[x=2.828560in,y=3.276310in,left,base]{\color{textcolor}\rmfamily\fontsize{18.000000}{21.600000}\selectfont LPP\(\displaystyle ^{R,T}\)}%
\end{pgfscope}%
\begin{pgfscope}%
\pgfsetbuttcap%
\pgfsetroundjoin%
\pgfsetlinewidth{1.505625pt}%
\definecolor{currentstroke}{rgb}{1.000000,0.498039,0.054902}%
\pgfsetstrokecolor{currentstroke}%
\pgfsetstrokeopacity{0.900000}%
\pgfsetdash{{5.550000pt}{2.400000pt}}{0.000000pt}%
\pgfpathmoveto{\pgfqpoint{3.764510in}{3.363810in}}%
\pgfpathlineto{\pgfqpoint{4.014510in}{3.363810in}}%
\pgfpathlineto{\pgfqpoint{4.264510in}{3.363810in}}%
\pgfusepath{stroke}%
\end{pgfscope}%
\begin{pgfscope}%
\pgfsetbuttcap%
\pgfsetmiterjoin%
\definecolor{currentfill}{rgb}{1.000000,0.498039,0.054902}%
\pgfsetfillcolor{currentfill}%
\pgfsetfillopacity{0.900000}%
\pgfsetlinewidth{1.003750pt}%
\definecolor{currentstroke}{rgb}{1.000000,0.498039,0.054902}%
\pgfsetstrokecolor{currentstroke}%
\pgfsetstrokeopacity{0.900000}%
\pgfsetdash{}{0pt}%
\pgfsys@defobject{currentmarker}{\pgfqpoint{-0.058926in}{-0.058926in}}{\pgfqpoint{0.058926in}{0.058926in}}{%
\pgfpathmoveto{\pgfqpoint{0.000000in}{-0.058926in}}%
\pgfpathlineto{\pgfqpoint{0.058926in}{0.000000in}}%
\pgfpathlineto{\pgfqpoint{0.000000in}{0.058926in}}%
\pgfpathlineto{\pgfqpoint{-0.058926in}{0.000000in}}%
\pgfpathlineto{\pgfqpoint{0.000000in}{-0.058926in}}%
\pgfpathclose%
\pgfusepath{stroke,fill}%
}%
\begin{pgfscope}%
\pgfsys@transformshift{4.014510in}{3.363810in}%
\pgfsys@useobject{currentmarker}{}%
\end{pgfscope}%
\end{pgfscope}%
\begin{pgfscope}%
\definecolor{textcolor}{rgb}{0.000000,0.000000,0.000000}%
\pgfsetstrokecolor{textcolor}%
\pgfsetfillcolor{textcolor}%
\pgftext[x=4.389510in,y=3.276310in,left,base]{\color{textcolor}\rmfamily\fontsize{18.000000}{21.600000}\selectfont LPP\(\displaystyle ^{T}\)}%
\end{pgfscope}%
\begin{pgfscope}%
\pgfsetbuttcap%
\pgfsetroundjoin%
\pgfsetlinewidth{1.505625pt}%
\definecolor{currentstroke}{rgb}{1.000000,0.498039,0.054902}%
\pgfsetstrokecolor{currentstroke}%
\pgfsetstrokeopacity{0.800000}%
\pgfsetdash{{1.500000pt}{2.475000pt}}{0.000000pt}%
\pgfpathmoveto{\pgfqpoint{5.154542in}{3.363810in}}%
\pgfpathlineto{\pgfqpoint{5.404542in}{3.363810in}}%
\pgfpathlineto{\pgfqpoint{5.654542in}{3.363810in}}%
\pgfusepath{stroke}%
\end{pgfscope}%
\begin{pgfscope}%
\pgfsetbuttcap%
\pgfsetroundjoin%
\definecolor{currentfill}{rgb}{1.000000,0.498039,0.054902}%
\pgfsetfillcolor{currentfill}%
\pgfsetfillopacity{0.800000}%
\pgfsetlinewidth{1.003750pt}%
\definecolor{currentstroke}{rgb}{1.000000,0.498039,0.054902}%
\pgfsetstrokecolor{currentstroke}%
\pgfsetstrokeopacity{0.800000}%
\pgfsetdash{}{0pt}%
\pgfsys@defobject{currentmarker}{\pgfqpoint{-0.041667in}{-0.041667in}}{\pgfqpoint{0.041667in}{0.041667in}}{%
\pgfpathmoveto{\pgfqpoint{0.000000in}{-0.041667in}}%
\pgfpathcurveto{\pgfqpoint{0.011050in}{-0.041667in}}{\pgfqpoint{0.021649in}{-0.037276in}}{\pgfqpoint{0.029463in}{-0.029463in}}%
\pgfpathcurveto{\pgfqpoint{0.037276in}{-0.021649in}}{\pgfqpoint{0.041667in}{-0.011050in}}{\pgfqpoint{0.041667in}{0.000000in}}%
\pgfpathcurveto{\pgfqpoint{0.041667in}{0.011050in}}{\pgfqpoint{0.037276in}{0.021649in}}{\pgfqpoint{0.029463in}{0.029463in}}%
\pgfpathcurveto{\pgfqpoint{0.021649in}{0.037276in}}{\pgfqpoint{0.011050in}{0.041667in}}{\pgfqpoint{0.000000in}{0.041667in}}%
\pgfpathcurveto{\pgfqpoint{-0.011050in}{0.041667in}}{\pgfqpoint{-0.021649in}{0.037276in}}{\pgfqpoint{-0.029463in}{0.029463in}}%
\pgfpathcurveto{\pgfqpoint{-0.037276in}{0.021649in}}{\pgfqpoint{-0.041667in}{0.011050in}}{\pgfqpoint{-0.041667in}{0.000000in}}%
\pgfpathcurveto{\pgfqpoint{-0.041667in}{-0.011050in}}{\pgfqpoint{-0.037276in}{-0.021649in}}{\pgfqpoint{-0.029463in}{-0.029463in}}%
\pgfpathcurveto{\pgfqpoint{-0.021649in}{-0.037276in}}{\pgfqpoint{-0.011050in}{-0.041667in}}{\pgfqpoint{0.000000in}{-0.041667in}}%
\pgfpathlineto{\pgfqpoint{0.000000in}{-0.041667in}}%
\pgfpathclose%
\pgfusepath{stroke,fill}%
}%
\begin{pgfscope}%
\pgfsys@transformshift{5.404542in}{3.363810in}%
\pgfsys@useobject{currentmarker}{}%
\end{pgfscope}%
\end{pgfscope}%
\begin{pgfscope}%
\definecolor{textcolor}{rgb}{0.000000,0.000000,0.000000}%
\pgfsetstrokecolor{textcolor}%
\pgfsetfillcolor{textcolor}%
\pgftext[x=5.779542in,y=3.276310in,left,base]{\color{textcolor}\rmfamily\fontsize{18.000000}{21.600000}\selectfont LPP\(\displaystyle ^{T, P}\)}%
\end{pgfscope}%
\end{pgfpicture}%
\makeatother%
\endgroup%

%% file: figures/DemandPlots/TightLayout/DE_DemandCapturedPct_Gap10pct.pgf
\begingroup%
\makeatletter%
\begin{pgfpicture}%
\pgfpathrectangle{\pgfpointorigin}{\pgfqpoint{7.775300in}{3.651310in}}%
\pgfusepath{use as bounding box, clip}%
\begin{pgfscope}%
\pgfsetbuttcap%
\pgfsetmiterjoin%
\definecolor{currentfill}{rgb}{1.000000,1.000000,1.000000}%
\pgfsetfillcolor{currentfill}%
\pgfsetlinewidth{0.000000pt}%
\definecolor{currentstroke}{rgb}{1.000000,1.000000,1.000000}%
\pgfsetstrokecolor{currentstroke}%
\pgfsetdash{}{0pt}%
\pgfpathmoveto{\pgfqpoint{0.000000in}{0.000000in}}%
\pgfpathlineto{\pgfqpoint{7.775300in}{0.000000in}}%
\pgfpathlineto{\pgfqpoint{7.775300in}{3.651310in}}%
\pgfpathlineto{\pgfqpoint{0.000000in}{3.651310in}}%
\pgfpathlineto{\pgfqpoint{0.000000in}{0.000000in}}%
\pgfpathclose%
\pgfusepath{fill}%
\end{pgfscope}%
\begin{pgfscope}%
\pgfsetbuttcap%
\pgfsetmiterjoin%
\definecolor{currentfill}{rgb}{1.000000,1.000000,1.000000}%
\pgfsetfillcolor{currentfill}%
\pgfsetlinewidth{0.000000pt}%
\definecolor{currentstroke}{rgb}{0.000000,0.000000,0.000000}%
\pgfsetstrokecolor{currentstroke}%
\pgfsetstrokeopacity{0.000000}%
\pgfsetdash{}{0pt}%
\pgfpathmoveto{\pgfqpoint{0.686263in}{0.679475in}}%
\pgfpathlineto{\pgfqpoint{7.675300in}{0.679475in}}%
\pgfpathlineto{\pgfqpoint{7.675300in}{3.135897in}}%
\pgfpathlineto{\pgfqpoint{0.686263in}{3.135897in}}%
\pgfpathlineto{\pgfqpoint{0.686263in}{0.679475in}}%
\pgfpathclose%
\pgfusepath{fill}%
\end{pgfscope}%
\begin{pgfscope}%
\pgfsetbuttcap%
\pgfsetroundjoin%
\definecolor{currentfill}{rgb}{0.000000,0.000000,0.000000}%
\pgfsetfillcolor{currentfill}%
\pgfsetlinewidth{0.803000pt}%
\definecolor{currentstroke}{rgb}{0.000000,0.000000,0.000000}%
\pgfsetstrokecolor{currentstroke}%
\pgfsetdash{}{0pt}%
\pgfsys@defobject{currentmarker}{\pgfqpoint{0.000000in}{-0.048611in}}{\pgfqpoint{0.000000in}{0.000000in}}{%
\pgfpathmoveto{\pgfqpoint{0.000000in}{0.000000in}}%
\pgfpathlineto{\pgfqpoint{0.000000in}{-0.048611in}}%
\pgfusepath{stroke,fill}%
}%
\begin{pgfscope}%
\pgfsys@transformshift{0.686263in}{0.679475in}%
\pgfsys@useobject{currentmarker}{}%
\end{pgfscope}%
\end{pgfscope}%
\begin{pgfscope}%
\definecolor{textcolor}{rgb}{0.000000,0.000000,0.000000}%
\pgfsetstrokecolor{textcolor}%
\pgfsetfillcolor{textcolor}%
\pgftext[x=0.686263in,y=0.582253in,,top]{\color{textcolor}\rmfamily\fontsize{16.000000}{19.200000}\selectfont \(\displaystyle {0.0}\)}%
\end{pgfscope}%
\begin{pgfscope}%
\pgfsetbuttcap%
\pgfsetroundjoin%
\definecolor{currentfill}{rgb}{0.000000,0.000000,0.000000}%
\pgfsetfillcolor{currentfill}%
\pgfsetlinewidth{0.803000pt}%
\definecolor{currentstroke}{rgb}{0.000000,0.000000,0.000000}%
\pgfsetstrokecolor{currentstroke}%
\pgfsetdash{}{0pt}%
\pgfsys@defobject{currentmarker}{\pgfqpoint{0.000000in}{-0.048611in}}{\pgfqpoint{0.000000in}{0.000000in}}{%
\pgfpathmoveto{\pgfqpoint{0.000000in}{0.000000in}}%
\pgfpathlineto{\pgfqpoint{0.000000in}{-0.048611in}}%
\pgfusepath{stroke,fill}%
}%
\begin{pgfscope}%
\pgfsys@transformshift{1.956997in}{0.679475in}%
\pgfsys@useobject{currentmarker}{}%
\end{pgfscope}%
\end{pgfscope}%
\begin{pgfscope}%
\definecolor{textcolor}{rgb}{0.000000,0.000000,0.000000}%
\pgfsetstrokecolor{textcolor}%
\pgfsetfillcolor{textcolor}%
\pgftext[x=1.956997in,y=0.582253in,,top]{\color{textcolor}\rmfamily\fontsize{16.000000}{19.200000}\selectfont \(\displaystyle {0.2}\)}%
\end{pgfscope}%
\begin{pgfscope}%
\pgfsetbuttcap%
\pgfsetroundjoin%
\definecolor{currentfill}{rgb}{0.000000,0.000000,0.000000}%
\pgfsetfillcolor{currentfill}%
\pgfsetlinewidth{0.803000pt}%
\definecolor{currentstroke}{rgb}{0.000000,0.000000,0.000000}%
\pgfsetstrokecolor{currentstroke}%
\pgfsetdash{}{0pt}%
\pgfsys@defobject{currentmarker}{\pgfqpoint{0.000000in}{-0.048611in}}{\pgfqpoint{0.000000in}{0.000000in}}{%
\pgfpathmoveto{\pgfqpoint{0.000000in}{0.000000in}}%
\pgfpathlineto{\pgfqpoint{0.000000in}{-0.048611in}}%
\pgfusepath{stroke,fill}%
}%
\begin{pgfscope}%
\pgfsys@transformshift{3.227731in}{0.679475in}%
\pgfsys@useobject{currentmarker}{}%
\end{pgfscope}%
\end{pgfscope}%
\begin{pgfscope}%
\definecolor{textcolor}{rgb}{0.000000,0.000000,0.000000}%
\pgfsetstrokecolor{textcolor}%
\pgfsetfillcolor{textcolor}%
\pgftext[x=3.227731in,y=0.582253in,,top]{\color{textcolor}\rmfamily\fontsize{16.000000}{19.200000}\selectfont \(\displaystyle {0.4}\)}%
\end{pgfscope}%
\begin{pgfscope}%
\pgfsetbuttcap%
\pgfsetroundjoin%
\definecolor{currentfill}{rgb}{0.000000,0.000000,0.000000}%
\pgfsetfillcolor{currentfill}%
\pgfsetlinewidth{0.803000pt}%
\definecolor{currentstroke}{rgb}{0.000000,0.000000,0.000000}%
\pgfsetstrokecolor{currentstroke}%
\pgfsetdash{}{0pt}%
\pgfsys@defobject{currentmarker}{\pgfqpoint{0.000000in}{-0.048611in}}{\pgfqpoint{0.000000in}{0.000000in}}{%
\pgfpathmoveto{\pgfqpoint{0.000000in}{0.000000in}}%
\pgfpathlineto{\pgfqpoint{0.000000in}{-0.048611in}}%
\pgfusepath{stroke,fill}%
}%
\begin{pgfscope}%
\pgfsys@transformshift{4.498465in}{0.679475in}%
\pgfsys@useobject{currentmarker}{}%
\end{pgfscope}%
\end{pgfscope}%
\begin{pgfscope}%
\definecolor{textcolor}{rgb}{0.000000,0.000000,0.000000}%
\pgfsetstrokecolor{textcolor}%
\pgfsetfillcolor{textcolor}%
\pgftext[x=4.498465in,y=0.582253in,,top]{\color{textcolor}\rmfamily\fontsize{16.000000}{19.200000}\selectfont \(\displaystyle {0.6}\)}%
\end{pgfscope}%
\begin{pgfscope}%
\pgfsetbuttcap%
\pgfsetroundjoin%
\definecolor{currentfill}{rgb}{0.000000,0.000000,0.000000}%
\pgfsetfillcolor{currentfill}%
\pgfsetlinewidth{0.803000pt}%
\definecolor{currentstroke}{rgb}{0.000000,0.000000,0.000000}%
\pgfsetstrokecolor{currentstroke}%
\pgfsetdash{}{0pt}%
\pgfsys@defobject{currentmarker}{\pgfqpoint{0.000000in}{-0.048611in}}{\pgfqpoint{0.000000in}{0.000000in}}{%
\pgfpathmoveto{\pgfqpoint{0.000000in}{0.000000in}}%
\pgfpathlineto{\pgfqpoint{0.000000in}{-0.048611in}}%
\pgfusepath{stroke,fill}%
}%
\begin{pgfscope}%
\pgfsys@transformshift{5.769199in}{0.679475in}%
\pgfsys@useobject{currentmarker}{}%
\end{pgfscope}%
\end{pgfscope}%
\begin{pgfscope}%
\definecolor{textcolor}{rgb}{0.000000,0.000000,0.000000}%
\pgfsetstrokecolor{textcolor}%
\pgfsetfillcolor{textcolor}%
\pgftext[x=5.769199in,y=0.582253in,,top]{\color{textcolor}\rmfamily\fontsize{16.000000}{19.200000}\selectfont \(\displaystyle {0.8}\)}%
\end{pgfscope}%
\begin{pgfscope}%
\pgfsetbuttcap%
\pgfsetroundjoin%
\definecolor{currentfill}{rgb}{0.000000,0.000000,0.000000}%
\pgfsetfillcolor{currentfill}%
\pgfsetlinewidth{0.803000pt}%
\definecolor{currentstroke}{rgb}{0.000000,0.000000,0.000000}%
\pgfsetstrokecolor{currentstroke}%
\pgfsetdash{}{0pt}%
\pgfsys@defobject{currentmarker}{\pgfqpoint{0.000000in}{-0.048611in}}{\pgfqpoint{0.000000in}{0.000000in}}{%
\pgfpathmoveto{\pgfqpoint{0.000000in}{0.000000in}}%
\pgfpathlineto{\pgfqpoint{0.000000in}{-0.048611in}}%
\pgfusepath{stroke,fill}%
}%
\begin{pgfscope}%
\pgfsys@transformshift{7.039933in}{0.679475in}%
\pgfsys@useobject{currentmarker}{}%
\end{pgfscope}%
\end{pgfscope}%
\begin{pgfscope}%
\definecolor{textcolor}{rgb}{0.000000,0.000000,0.000000}%
\pgfsetstrokecolor{textcolor}%
\pgfsetfillcolor{textcolor}%
\pgftext[x=7.039933in,y=0.582253in,,top]{\color{textcolor}\rmfamily\fontsize{16.000000}{19.200000}\selectfont \(\displaystyle {1.0}\)}%
\end{pgfscope}%
\begin{pgfscope}%
\definecolor{textcolor}{rgb}{0.000000,0.000000,0.000000}%
\pgfsetstrokecolor{textcolor}%
\pgfsetfillcolor{textcolor}%
\pgftext[x=4.180781in,y=0.313349in,,top]{\color{textcolor}\rmfamily\fontsize{18.000000}{21.600000}\selectfont \(\displaystyle \lambda\)}%
\end{pgfscope}%
\begin{pgfscope}%
\pgfsetbuttcap%
\pgfsetroundjoin%
\definecolor{currentfill}{rgb}{0.000000,0.000000,0.000000}%
\pgfsetfillcolor{currentfill}%
\pgfsetlinewidth{0.803000pt}%
\definecolor{currentstroke}{rgb}{0.000000,0.000000,0.000000}%
\pgfsetstrokecolor{currentstroke}%
\pgfsetdash{}{0pt}%
\pgfsys@defobject{currentmarker}{\pgfqpoint{-0.048611in}{0.000000in}}{\pgfqpoint{-0.000000in}{0.000000in}}{%
\pgfpathmoveto{\pgfqpoint{-0.000000in}{0.000000in}}%
\pgfpathlineto{\pgfqpoint{-0.048611in}{0.000000in}}%
\pgfusepath{stroke,fill}%
}%
\begin{pgfscope}%
\pgfsys@transformshift{0.686263in}{1.086027in}%
\pgfsys@useobject{currentmarker}{}%
\end{pgfscope}%
\end{pgfscope}%
\begin{pgfscope}%
\definecolor{textcolor}{rgb}{0.000000,0.000000,0.000000}%
\pgfsetstrokecolor{textcolor}%
\pgfsetfillcolor{textcolor}%
\pgftext[x=0.368904in, y=1.002694in, left, base]{\color{textcolor}\rmfamily\fontsize{16.000000}{19.200000}\selectfont \(\displaystyle {55}\)}%
\end{pgfscope}%
\begin{pgfscope}%
\pgfsetbuttcap%
\pgfsetroundjoin%
\definecolor{currentfill}{rgb}{0.000000,0.000000,0.000000}%
\pgfsetfillcolor{currentfill}%
\pgfsetlinewidth{0.803000pt}%
\definecolor{currentstroke}{rgb}{0.000000,0.000000,0.000000}%
\pgfsetstrokecolor{currentstroke}%
\pgfsetdash{}{0pt}%
\pgfsys@defobject{currentmarker}{\pgfqpoint{-0.048611in}{0.000000in}}{\pgfqpoint{-0.000000in}{0.000000in}}{%
\pgfpathmoveto{\pgfqpoint{-0.000000in}{0.000000in}}%
\pgfpathlineto{\pgfqpoint{-0.048611in}{0.000000in}}%
\pgfusepath{stroke,fill}%
}%
\begin{pgfscope}%
\pgfsys@transformshift{0.686263in}{1.536251in}%
\pgfsys@useobject{currentmarker}{}%
\end{pgfscope}%
\end{pgfscope}%
\begin{pgfscope}%
\definecolor{textcolor}{rgb}{0.000000,0.000000,0.000000}%
\pgfsetstrokecolor{textcolor}%
\pgfsetfillcolor{textcolor}%
\pgftext[x=0.368904in, y=1.452918in, left, base]{\color{textcolor}\rmfamily\fontsize{16.000000}{19.200000}\selectfont \(\displaystyle {60}\)}%
\end{pgfscope}%
\begin{pgfscope}%
\pgfsetbuttcap%
\pgfsetroundjoin%
\definecolor{currentfill}{rgb}{0.000000,0.000000,0.000000}%
\pgfsetfillcolor{currentfill}%
\pgfsetlinewidth{0.803000pt}%
\definecolor{currentstroke}{rgb}{0.000000,0.000000,0.000000}%
\pgfsetstrokecolor{currentstroke}%
\pgfsetdash{}{0pt}%
\pgfsys@defobject{currentmarker}{\pgfqpoint{-0.048611in}{0.000000in}}{\pgfqpoint{-0.000000in}{0.000000in}}{%
\pgfpathmoveto{\pgfqpoint{-0.000000in}{0.000000in}}%
\pgfpathlineto{\pgfqpoint{-0.048611in}{0.000000in}}%
\pgfusepath{stroke,fill}%
}%
\begin{pgfscope}%
\pgfsys@transformshift{0.686263in}{1.986475in}%
\pgfsys@useobject{currentmarker}{}%
\end{pgfscope}%
\end{pgfscope}%
\begin{pgfscope}%
\definecolor{textcolor}{rgb}{0.000000,0.000000,0.000000}%
\pgfsetstrokecolor{textcolor}%
\pgfsetfillcolor{textcolor}%
\pgftext[x=0.368904in, y=1.903142in, left, base]{\color{textcolor}\rmfamily\fontsize{16.000000}{19.200000}\selectfont \(\displaystyle {65}\)}%
\end{pgfscope}%
\begin{pgfscope}%
\pgfsetbuttcap%
\pgfsetroundjoin%
\definecolor{currentfill}{rgb}{0.000000,0.000000,0.000000}%
\pgfsetfillcolor{currentfill}%
\pgfsetlinewidth{0.803000pt}%
\definecolor{currentstroke}{rgb}{0.000000,0.000000,0.000000}%
\pgfsetstrokecolor{currentstroke}%
\pgfsetdash{}{0pt}%
\pgfsys@defobject{currentmarker}{\pgfqpoint{-0.048611in}{0.000000in}}{\pgfqpoint{-0.000000in}{0.000000in}}{%
\pgfpathmoveto{\pgfqpoint{-0.000000in}{0.000000in}}%
\pgfpathlineto{\pgfqpoint{-0.048611in}{0.000000in}}%
\pgfusepath{stroke,fill}%
}%
\begin{pgfscope}%
\pgfsys@transformshift{0.686263in}{2.436699in}%
\pgfsys@useobject{currentmarker}{}%
\end{pgfscope}%
\end{pgfscope}%
\begin{pgfscope}%
\definecolor{textcolor}{rgb}{0.000000,0.000000,0.000000}%
\pgfsetstrokecolor{textcolor}%
\pgfsetfillcolor{textcolor}%
\pgftext[x=0.368904in, y=2.353366in, left, base]{\color{textcolor}\rmfamily\fontsize{16.000000}{19.200000}\selectfont \(\displaystyle {70}\)}%
\end{pgfscope}%
\begin{pgfscope}%
\pgfsetbuttcap%
\pgfsetroundjoin%
\definecolor{currentfill}{rgb}{0.000000,0.000000,0.000000}%
\pgfsetfillcolor{currentfill}%
\pgfsetlinewidth{0.803000pt}%
\definecolor{currentstroke}{rgb}{0.000000,0.000000,0.000000}%
\pgfsetstrokecolor{currentstroke}%
\pgfsetdash{}{0pt}%
\pgfsys@defobject{currentmarker}{\pgfqpoint{-0.048611in}{0.000000in}}{\pgfqpoint{-0.000000in}{0.000000in}}{%
\pgfpathmoveto{\pgfqpoint{-0.000000in}{0.000000in}}%
\pgfpathlineto{\pgfqpoint{-0.048611in}{0.000000in}}%
\pgfusepath{stroke,fill}%
}%
\begin{pgfscope}%
\pgfsys@transformshift{0.686263in}{2.886923in}%
\pgfsys@useobject{currentmarker}{}%
\end{pgfscope}%
\end{pgfscope}%
\begin{pgfscope}%
\definecolor{textcolor}{rgb}{0.000000,0.000000,0.000000}%
\pgfsetstrokecolor{textcolor}%
\pgfsetfillcolor{textcolor}%
\pgftext[x=0.368904in, y=2.803590in, left, base]{\color{textcolor}\rmfamily\fontsize{16.000000}{19.200000}\selectfont \(\displaystyle {75}\)}%
\end{pgfscope}%
\begin{pgfscope}%
\definecolor{textcolor}{rgb}{0.000000,0.000000,0.000000}%
\pgfsetstrokecolor{textcolor}%
\pgfsetfillcolor{textcolor}%
\pgftext[x=0.313349in,y=1.907686in,,bottom,rotate=90.000000]{\color{textcolor}\rmfamily\fontsize{18.000000}{21.600000}\selectfont \%}%
\end{pgfscope}%
\begin{pgfscope}%
\pgfpathrectangle{\pgfqpoint{0.686263in}{0.679475in}}{\pgfqpoint{6.989037in}{2.456422in}}%
\pgfusepath{clip}%
\pgfsetrectcap%
\pgfsetroundjoin%
\pgfsetlinewidth{1.505625pt}%
\definecolor{currentstroke}{rgb}{0.172549,0.627451,0.172549}%
\pgfsetstrokecolor{currentstroke}%
\pgfsetdash{}{0pt}%
\pgfpathmoveto{\pgfqpoint{1.321630in}{1.584275in}}%
\pgfpathlineto{\pgfqpoint{2.274680in}{1.960962in}}%
\pgfpathlineto{\pgfqpoint{3.863098in}{2.192828in}}%
\pgfpathlineto{\pgfqpoint{5.451515in}{2.691076in}}%
\pgfpathlineto{\pgfqpoint{7.039933in}{3.024241in}}%
\pgfusepath{stroke}%
\end{pgfscope}%
\begin{pgfscope}%
\pgfpathrectangle{\pgfqpoint{0.686263in}{0.679475in}}{\pgfqpoint{6.989037in}{2.456422in}}%
\pgfusepath{clip}%
\pgfsetbuttcap%
\pgfsetmiterjoin%
\definecolor{currentfill}{rgb}{0.172549,0.627451,0.172549}%
\pgfsetfillcolor{currentfill}%
\pgfsetlinewidth{0.752812pt}%
\definecolor{currentstroke}{rgb}{1.000000,1.000000,1.000000}%
\pgfsetstrokecolor{currentstroke}%
\pgfsetdash{}{0pt}%
\pgfsys@defobject{currentmarker}{\pgfqpoint{-0.041667in}{-0.041667in}}{\pgfqpoint{0.041667in}{0.041667in}}{%
\pgfpathmoveto{\pgfqpoint{-0.041667in}{-0.041667in}}%
\pgfpathlineto{\pgfqpoint{0.041667in}{-0.041667in}}%
\pgfpathlineto{\pgfqpoint{0.041667in}{0.041667in}}%
\pgfpathlineto{\pgfqpoint{-0.041667in}{0.041667in}}%
\pgfpathlineto{\pgfqpoint{-0.041667in}{-0.041667in}}%
\pgfpathclose%
\pgfusepath{stroke,fill}%
}%
\begin{pgfscope}%
\pgfsys@transformshift{1.321630in}{1.584275in}%
\pgfsys@useobject{currentmarker}{}%
\end{pgfscope}%
\begin{pgfscope}%
\pgfsys@transformshift{2.274680in}{1.960962in}%
\pgfsys@useobject{currentmarker}{}%
\end{pgfscope}%
\begin{pgfscope}%
\pgfsys@transformshift{3.863098in}{2.192828in}%
\pgfsys@useobject{currentmarker}{}%
\end{pgfscope}%
\begin{pgfscope}%
\pgfsys@transformshift{5.451515in}{2.691076in}%
\pgfsys@useobject{currentmarker}{}%
\end{pgfscope}%
\begin{pgfscope}%
\pgfsys@transformshift{7.039933in}{3.024241in}%
\pgfsys@useobject{currentmarker}{}%
\end{pgfscope}%
\end{pgfscope}%
\begin{pgfscope}%
\pgfpathrectangle{\pgfqpoint{0.686263in}{0.679475in}}{\pgfqpoint{6.989037in}{2.456422in}}%
\pgfusepath{clip}%
\pgfsetbuttcap%
\pgfsetroundjoin%
\pgfsetlinewidth{1.505625pt}%
\definecolor{currentstroke}{rgb}{0.172549,0.627451,0.172549}%
\pgfsetstrokecolor{currentstroke}%
\pgfsetstrokeopacity{0.900000}%
\pgfsetdash{{7.500000pt}{7.500000pt}}{0.000000pt}%
\pgfpathmoveto{\pgfqpoint{1.321630in}{0.791131in}}%
\pgfpathlineto{\pgfqpoint{2.274680in}{1.049259in}}%
\pgfpathlineto{\pgfqpoint{3.863098in}{1.509238in}}%
\pgfpathlineto{\pgfqpoint{5.451515in}{2.087025in}}%
\pgfpathlineto{\pgfqpoint{7.039933in}{2.393928in}}%
\pgfusepath{stroke}%
\end{pgfscope}%
\begin{pgfscope}%
\pgfpathrectangle{\pgfqpoint{0.686263in}{0.679475in}}{\pgfqpoint{6.989037in}{2.456422in}}%
\pgfusepath{clip}%
\pgfsetbuttcap%
\pgfsetmiterjoin%
\definecolor{currentfill}{rgb}{0.172549,0.627451,0.172549}%
\pgfsetfillcolor{currentfill}%
\pgfsetfillopacity{0.900000}%
\pgfsetlinewidth{0.752812pt}%
\definecolor{currentstroke}{rgb}{1.000000,1.000000,1.000000}%
\pgfsetstrokecolor{currentstroke}%
\pgfsetdash{}{0pt}%
\pgfsys@defobject{currentmarker}{\pgfqpoint{-0.058926in}{-0.058926in}}{\pgfqpoint{0.058926in}{0.058926in}}{%
\pgfpathmoveto{\pgfqpoint{0.000000in}{-0.058926in}}%
\pgfpathlineto{\pgfqpoint{0.058926in}{0.000000in}}%
\pgfpathlineto{\pgfqpoint{0.000000in}{0.058926in}}%
\pgfpathlineto{\pgfqpoint{-0.058926in}{0.000000in}}%
\pgfpathlineto{\pgfqpoint{0.000000in}{-0.058926in}}%
\pgfpathclose%
\pgfusepath{stroke,fill}%
}%
\begin{pgfscope}%
\pgfsys@transformshift{1.321630in}{0.791131in}%
\pgfsys@useobject{currentmarker}{}%
\end{pgfscope}%
\begin{pgfscope}%
\pgfsys@transformshift{2.274680in}{1.049259in}%
\pgfsys@useobject{currentmarker}{}%
\end{pgfscope}%
\begin{pgfscope}%
\pgfsys@transformshift{3.863098in}{1.509238in}%
\pgfsys@useobject{currentmarker}{}%
\end{pgfscope}%
\begin{pgfscope}%
\pgfsys@transformshift{5.451515in}{2.087025in}%
\pgfsys@useobject{currentmarker}{}%
\end{pgfscope}%
\begin{pgfscope}%
\pgfsys@transformshift{7.039933in}{2.393928in}%
\pgfsys@useobject{currentmarker}{}%
\end{pgfscope}%
\end{pgfscope}%
\begin{pgfscope}%
\pgfpathrectangle{\pgfqpoint{0.686263in}{0.679475in}}{\pgfqpoint{6.989037in}{2.456422in}}%
\pgfusepath{clip}%
\pgfsetbuttcap%
\pgfsetroundjoin%
\pgfsetlinewidth{1.505625pt}%
\definecolor{currentstroke}{rgb}{0.172549,0.627451,0.172549}%
\pgfsetstrokecolor{currentstroke}%
\pgfsetstrokeopacity{0.800000}%
\pgfsetdash{{1.500000pt}{4.500000pt}}{0.000000pt}%
\pgfpathmoveto{\pgfqpoint{1.321630in}{2.225844in}}%
\pgfpathlineto{\pgfqpoint{2.274680in}{2.256609in}}%
\pgfpathlineto{\pgfqpoint{3.863098in}{2.353408in}}%
\pgfpathlineto{\pgfqpoint{5.451515in}{2.721841in}}%
\pgfpathlineto{\pgfqpoint{7.039933in}{2.871165in}}%
\pgfusepath{stroke}%
\end{pgfscope}%
\begin{pgfscope}%
\pgfpathrectangle{\pgfqpoint{0.686263in}{0.679475in}}{\pgfqpoint{6.989037in}{2.456422in}}%
\pgfusepath{clip}%
\pgfsetbuttcap%
\pgfsetroundjoin%
\definecolor{currentfill}{rgb}{0.172549,0.627451,0.172549}%
\pgfsetfillcolor{currentfill}%
\pgfsetfillopacity{0.800000}%
\pgfsetlinewidth{0.752812pt}%
\definecolor{currentstroke}{rgb}{1.000000,1.000000,1.000000}%
\pgfsetstrokecolor{currentstroke}%
\pgfsetdash{}{0pt}%
\pgfsys@defobject{currentmarker}{\pgfqpoint{-0.041667in}{-0.041667in}}{\pgfqpoint{0.041667in}{0.041667in}}{%
\pgfpathmoveto{\pgfqpoint{0.000000in}{-0.041667in}}%
\pgfpathcurveto{\pgfqpoint{0.011050in}{-0.041667in}}{\pgfqpoint{0.021649in}{-0.037276in}}{\pgfqpoint{0.029463in}{-0.029463in}}%
\pgfpathcurveto{\pgfqpoint{0.037276in}{-0.021649in}}{\pgfqpoint{0.041667in}{-0.011050in}}{\pgfqpoint{0.041667in}{0.000000in}}%
\pgfpathcurveto{\pgfqpoint{0.041667in}{0.011050in}}{\pgfqpoint{0.037276in}{0.021649in}}{\pgfqpoint{0.029463in}{0.029463in}}%
\pgfpathcurveto{\pgfqpoint{0.021649in}{0.037276in}}{\pgfqpoint{0.011050in}{0.041667in}}{\pgfqpoint{0.000000in}{0.041667in}}%
\pgfpathcurveto{\pgfqpoint{-0.011050in}{0.041667in}}{\pgfqpoint{-0.021649in}{0.037276in}}{\pgfqpoint{-0.029463in}{0.029463in}}%
\pgfpathcurveto{\pgfqpoint{-0.037276in}{0.021649in}}{\pgfqpoint{-0.041667in}{0.011050in}}{\pgfqpoint{-0.041667in}{0.000000in}}%
\pgfpathcurveto{\pgfqpoint{-0.041667in}{-0.011050in}}{\pgfqpoint{-0.037276in}{-0.021649in}}{\pgfqpoint{-0.029463in}{-0.029463in}}%
\pgfpathcurveto{\pgfqpoint{-0.021649in}{-0.037276in}}{\pgfqpoint{-0.011050in}{-0.041667in}}{\pgfqpoint{0.000000in}{-0.041667in}}%
\pgfpathlineto{\pgfqpoint{0.000000in}{-0.041667in}}%
\pgfpathclose%
\pgfusepath{stroke,fill}%
}%
\begin{pgfscope}%
\pgfsys@transformshift{1.321630in}{2.225844in}%
\pgfsys@useobject{currentmarker}{}%
\end{pgfscope}%
\begin{pgfscope}%
\pgfsys@transformshift{2.274680in}{2.256609in}%
\pgfsys@useobject{currentmarker}{}%
\end{pgfscope}%
\begin{pgfscope}%
\pgfsys@transformshift{3.863098in}{2.353408in}%
\pgfsys@useobject{currentmarker}{}%
\end{pgfscope}%
\begin{pgfscope}%
\pgfsys@transformshift{5.451515in}{2.721841in}%
\pgfsys@useobject{currentmarker}{}%
\end{pgfscope}%
\begin{pgfscope}%
\pgfsys@transformshift{7.039933in}{2.871165in}%
\pgfsys@useobject{currentmarker}{}%
\end{pgfscope}%
\end{pgfscope}%
\begin{pgfscope}%
\pgfsetrectcap%
\pgfsetmiterjoin%
\pgfsetlinewidth{0.803000pt}%
\definecolor{currentstroke}{rgb}{0.000000,0.000000,0.000000}%
\pgfsetstrokecolor{currentstroke}%
\pgfsetdash{}{0pt}%
\pgfpathmoveto{\pgfqpoint{0.686263in}{0.679475in}}%
\pgfpathlineto{\pgfqpoint{0.686263in}{3.135897in}}%
\pgfusepath{stroke}%
\end{pgfscope}%
\begin{pgfscope}%
\pgfsetrectcap%
\pgfsetmiterjoin%
\pgfsetlinewidth{0.803000pt}%
\definecolor{currentstroke}{rgb}{0.000000,0.000000,0.000000}%
\pgfsetstrokecolor{currentstroke}%
\pgfsetdash{}{0pt}%
\pgfpathmoveto{\pgfqpoint{7.675300in}{0.679475in}}%
\pgfpathlineto{\pgfqpoint{7.675300in}{3.135897in}}%
\pgfusepath{stroke}%
\end{pgfscope}%
\begin{pgfscope}%
\pgfsetrectcap%
\pgfsetmiterjoin%
\pgfsetlinewidth{0.803000pt}%
\definecolor{currentstroke}{rgb}{0.000000,0.000000,0.000000}%
\pgfsetstrokecolor{currentstroke}%
\pgfsetdash{}{0pt}%
\pgfpathmoveto{\pgfqpoint{0.686263in}{0.679475in}}%
\pgfpathlineto{\pgfqpoint{7.675300in}{0.679475in}}%
\pgfusepath{stroke}%
\end{pgfscope}%
\begin{pgfscope}%
\pgfsetrectcap%
\pgfsetmiterjoin%
\pgfsetlinewidth{0.803000pt}%
\definecolor{currentstroke}{rgb}{0.000000,0.000000,0.000000}%
\pgfsetstrokecolor{currentstroke}%
\pgfsetdash{}{0pt}%
\pgfpathmoveto{\pgfqpoint{0.686263in}{3.135897in}}%
\pgfpathlineto{\pgfqpoint{7.675300in}{3.135897in}}%
\pgfusepath{stroke}%
\end{pgfscope}%
\begin{pgfscope}%
\pgfsetrectcap%
\pgfsetroundjoin%
\pgfsetlinewidth{1.505625pt}%
\definecolor{currentstroke}{rgb}{0.172549,0.627451,0.172549}%
\pgfsetstrokecolor{currentstroke}%
\pgfsetdash{}{0pt}%
\pgfpathmoveto{\pgfqpoint{2.028215in}{3.363810in}}%
\pgfpathlineto{\pgfqpoint{2.278215in}{3.363810in}}%
\pgfpathlineto{\pgfqpoint{2.528215in}{3.363810in}}%
\pgfusepath{stroke}%
\end{pgfscope}%
\begin{pgfscope}%
\pgfsetbuttcap%
\pgfsetmiterjoin%
\definecolor{currentfill}{rgb}{0.172549,0.627451,0.172549}%
\pgfsetfillcolor{currentfill}%
\pgfsetlinewidth{1.003750pt}%
\definecolor{currentstroke}{rgb}{0.172549,0.627451,0.172549}%
\pgfsetstrokecolor{currentstroke}%
\pgfsetdash{}{0pt}%
\pgfsys@defobject{currentmarker}{\pgfqpoint{-0.041667in}{-0.041667in}}{\pgfqpoint{0.041667in}{0.041667in}}{%
\pgfpathmoveto{\pgfqpoint{-0.041667in}{-0.041667in}}%
\pgfpathlineto{\pgfqpoint{0.041667in}{-0.041667in}}%
\pgfpathlineto{\pgfqpoint{0.041667in}{0.041667in}}%
\pgfpathlineto{\pgfqpoint{-0.041667in}{0.041667in}}%
\pgfpathlineto{\pgfqpoint{-0.041667in}{-0.041667in}}%
\pgfpathclose%
\pgfusepath{stroke,fill}%
}%
\begin{pgfscope}%
\pgfsys@transformshift{2.278215in}{3.363810in}%
\pgfsys@useobject{currentmarker}{}%
\end{pgfscope}%
\end{pgfscope}%
\begin{pgfscope}%
\definecolor{textcolor}{rgb}{0.000000,0.000000,0.000000}%
\pgfsetstrokecolor{textcolor}%
\pgfsetfillcolor{textcolor}%
\pgftext[x=2.653215in,y=3.276310in,left,base]{\color{textcolor}\rmfamily\fontsize{18.000000}{21.600000}\selectfont LPP\(\displaystyle ^{R,T}\)}%
\end{pgfscope}%
\begin{pgfscope}%
\pgfsetbuttcap%
\pgfsetroundjoin%
\pgfsetlinewidth{1.505625pt}%
\definecolor{currentstroke}{rgb}{0.172549,0.627451,0.172549}%
\pgfsetstrokecolor{currentstroke}%
\pgfsetstrokeopacity{0.900000}%
\pgfsetdash{{5.550000pt}{2.400000pt}}{0.000000pt}%
\pgfpathmoveto{\pgfqpoint{3.589165in}{3.363810in}}%
\pgfpathlineto{\pgfqpoint{3.839165in}{3.363810in}}%
\pgfpathlineto{\pgfqpoint{4.089165in}{3.363810in}}%
\pgfusepath{stroke}%
\end{pgfscope}%
\begin{pgfscope}%
\pgfsetbuttcap%
\pgfsetmiterjoin%
\definecolor{currentfill}{rgb}{0.172549,0.627451,0.172549}%
\pgfsetfillcolor{currentfill}%
\pgfsetfillopacity{0.900000}%
\pgfsetlinewidth{1.003750pt}%
\definecolor{currentstroke}{rgb}{0.172549,0.627451,0.172549}%
\pgfsetstrokecolor{currentstroke}%
\pgfsetstrokeopacity{0.900000}%
\pgfsetdash{}{0pt}%
\pgfsys@defobject{currentmarker}{\pgfqpoint{-0.058926in}{-0.058926in}}{\pgfqpoint{0.058926in}{0.058926in}}{%
\pgfpathmoveto{\pgfqpoint{0.000000in}{-0.058926in}}%
\pgfpathlineto{\pgfqpoint{0.058926in}{0.000000in}}%
\pgfpathlineto{\pgfqpoint{0.000000in}{0.058926in}}%
\pgfpathlineto{\pgfqpoint{-0.058926in}{0.000000in}}%
\pgfpathlineto{\pgfqpoint{0.000000in}{-0.058926in}}%
\pgfpathclose%
\pgfusepath{stroke,fill}%
}%
\begin{pgfscope}%
\pgfsys@transformshift{3.839165in}{3.363810in}%
\pgfsys@useobject{currentmarker}{}%
\end{pgfscope}%
\end{pgfscope}%
\begin{pgfscope}%
\definecolor{textcolor}{rgb}{0.000000,0.000000,0.000000}%
\pgfsetstrokecolor{textcolor}%
\pgfsetfillcolor{textcolor}%
\pgftext[x=4.214165in,y=3.276310in,left,base]{\color{textcolor}\rmfamily\fontsize{18.000000}{21.600000}\selectfont LPP\(\displaystyle ^{T}\)}%
\end{pgfscope}%
\begin{pgfscope}%
\pgfsetbuttcap%
\pgfsetroundjoin%
\pgfsetlinewidth{1.505625pt}%
\definecolor{currentstroke}{rgb}{0.172549,0.627451,0.172549}%
\pgfsetstrokecolor{currentstroke}%
\pgfsetstrokeopacity{0.800000}%
\pgfsetdash{{1.500000pt}{2.475000pt}}{0.000000pt}%
\pgfpathmoveto{\pgfqpoint{4.979197in}{3.363810in}}%
\pgfpathlineto{\pgfqpoint{5.229197in}{3.363810in}}%
\pgfpathlineto{\pgfqpoint{5.479197in}{3.363810in}}%
\pgfusepath{stroke}%
\end{pgfscope}%
\begin{pgfscope}%
\pgfsetbuttcap%
\pgfsetroundjoin%
\definecolor{currentfill}{rgb}{0.172549,0.627451,0.172549}%
\pgfsetfillcolor{currentfill}%
\pgfsetfillopacity{0.800000}%
\pgfsetlinewidth{1.003750pt}%
\definecolor{currentstroke}{rgb}{0.172549,0.627451,0.172549}%
\pgfsetstrokecolor{currentstroke}%
\pgfsetstrokeopacity{0.800000}%
\pgfsetdash{}{0pt}%
\pgfsys@defobject{currentmarker}{\pgfqpoint{-0.041667in}{-0.041667in}}{\pgfqpoint{0.041667in}{0.041667in}}{%
\pgfpathmoveto{\pgfqpoint{0.000000in}{-0.041667in}}%
\pgfpathcurveto{\pgfqpoint{0.011050in}{-0.041667in}}{\pgfqpoint{0.021649in}{-0.037276in}}{\pgfqpoint{0.029463in}{-0.029463in}}%
\pgfpathcurveto{\pgfqpoint{0.037276in}{-0.021649in}}{\pgfqpoint{0.041667in}{-0.011050in}}{\pgfqpoint{0.041667in}{0.000000in}}%
\pgfpathcurveto{\pgfqpoint{0.041667in}{0.011050in}}{\pgfqpoint{0.037276in}{0.021649in}}{\pgfqpoint{0.029463in}{0.029463in}}%
\pgfpathcurveto{\pgfqpoint{0.021649in}{0.037276in}}{\pgfqpoint{0.011050in}{0.041667in}}{\pgfqpoint{0.000000in}{0.041667in}}%
\pgfpathcurveto{\pgfqpoint{-0.011050in}{0.041667in}}{\pgfqpoint{-0.021649in}{0.037276in}}{\pgfqpoint{-0.029463in}{0.029463in}}%
\pgfpathcurveto{\pgfqpoint{-0.037276in}{0.021649in}}{\pgfqpoint{-0.041667in}{0.011050in}}{\pgfqpoint{-0.041667in}{0.000000in}}%
\pgfpathcurveto{\pgfqpoint{-0.041667in}{-0.011050in}}{\pgfqpoint{-0.037276in}{-0.021649in}}{\pgfqpoint{-0.029463in}{-0.029463in}}%
\pgfpathcurveto{\pgfqpoint{-0.021649in}{-0.037276in}}{\pgfqpoint{-0.011050in}{-0.041667in}}{\pgfqpoint{0.000000in}{-0.041667in}}%
\pgfpathlineto{\pgfqpoint{0.000000in}{-0.041667in}}%
\pgfpathclose%
\pgfusepath{stroke,fill}%
}%
\begin{pgfscope}%
\pgfsys@transformshift{5.229197in}{3.363810in}%
\pgfsys@useobject{currentmarker}{}%
\end{pgfscope}%
\end{pgfscope}%
\begin{pgfscope}%
\definecolor{textcolor}{rgb}{0.000000,0.000000,0.000000}%
\pgfsetstrokecolor{textcolor}%
\pgfsetfillcolor{textcolor}%
\pgftext[x=5.604197in,y=3.276310in,left,base]{\color{textcolor}\rmfamily\fontsize{18.000000}{21.600000}\selectfont LPP\(\displaystyle ^{T, P}\)}%
\end{pgfscope}%
\end{pgfpicture}%
\makeatother%
\endgroup%

%% file: figures/DemandPlots/TightLayout/DE_UB_Gap10pct.pgf
\begingroup%
\makeatletter%
\begin{pgfpicture}%
\pgfpathrectangle{\pgfpointorigin}{\pgfqpoint{7.782018in}{3.651310in}}%
\pgfusepath{use as bounding box, clip}%
\begin{pgfscope}%
\pgfsetbuttcap%
\pgfsetmiterjoin%
\definecolor{currentfill}{rgb}{1.000000,1.000000,1.000000}%
\pgfsetfillcolor{currentfill}%
\pgfsetlinewidth{0.000000pt}%
\definecolor{currentstroke}{rgb}{1.000000,1.000000,1.000000}%
\pgfsetstrokecolor{currentstroke}%
\pgfsetdash{}{0pt}%
\pgfpathmoveto{\pgfqpoint{0.000000in}{0.000000in}}%
\pgfpathlineto{\pgfqpoint{7.782018in}{0.000000in}}%
\pgfpathlineto{\pgfqpoint{7.782018in}{3.651310in}}%
\pgfpathlineto{\pgfqpoint{0.000000in}{3.651310in}}%
\pgfpathlineto{\pgfqpoint{0.000000in}{0.000000in}}%
\pgfpathclose%
\pgfusepath{fill}%
\end{pgfscope}%
\begin{pgfscope}%
\pgfsetbuttcap%
\pgfsetmiterjoin%
\definecolor{currentfill}{rgb}{1.000000,1.000000,1.000000}%
\pgfsetfillcolor{currentfill}%
\pgfsetlinewidth{0.000000pt}%
\definecolor{currentstroke}{rgb}{0.000000,0.000000,0.000000}%
\pgfsetstrokecolor{currentstroke}%
\pgfsetstrokeopacity{0.000000}%
\pgfsetdash{}{0pt}%
\pgfpathmoveto{\pgfqpoint{1.126536in}{0.679475in}}%
\pgfpathlineto{\pgfqpoint{7.682018in}{0.679475in}}%
\pgfpathlineto{\pgfqpoint{7.682018in}{3.135897in}}%
\pgfpathlineto{\pgfqpoint{1.126536in}{3.135897in}}%
\pgfpathlineto{\pgfqpoint{1.126536in}{0.679475in}}%
\pgfpathclose%
\pgfusepath{fill}%
\end{pgfscope}%
\begin{pgfscope}%
\pgfsetbuttcap%
\pgfsetroundjoin%
\definecolor{currentfill}{rgb}{0.000000,0.000000,0.000000}%
\pgfsetfillcolor{currentfill}%
\pgfsetlinewidth{0.803000pt}%
\definecolor{currentstroke}{rgb}{0.000000,0.000000,0.000000}%
\pgfsetstrokecolor{currentstroke}%
\pgfsetdash{}{0pt}%
\pgfsys@defobject{currentmarker}{\pgfqpoint{0.000000in}{-0.048611in}}{\pgfqpoint{0.000000in}{0.000000in}}{%
\pgfpathmoveto{\pgfqpoint{0.000000in}{0.000000in}}%
\pgfpathlineto{\pgfqpoint{0.000000in}{-0.048611in}}%
\pgfusepath{stroke,fill}%
}%
\begin{pgfscope}%
\pgfsys@transformshift{1.126536in}{0.679475in}%
\pgfsys@useobject{currentmarker}{}%
\end{pgfscope}%
\end{pgfscope}%
\begin{pgfscope}%
\definecolor{textcolor}{rgb}{0.000000,0.000000,0.000000}%
\pgfsetstrokecolor{textcolor}%
\pgfsetfillcolor{textcolor}%
\pgftext[x=1.126536in,y=0.582253in,,top]{\color{textcolor}\rmfamily\fontsize{16.000000}{19.200000}\selectfont \(\displaystyle {0.0}\)}%
\end{pgfscope}%
\begin{pgfscope}%
\pgfsetbuttcap%
\pgfsetroundjoin%
\definecolor{currentfill}{rgb}{0.000000,0.000000,0.000000}%
\pgfsetfillcolor{currentfill}%
\pgfsetlinewidth{0.803000pt}%
\definecolor{currentstroke}{rgb}{0.000000,0.000000,0.000000}%
\pgfsetstrokecolor{currentstroke}%
\pgfsetdash{}{0pt}%
\pgfsys@defobject{currentmarker}{\pgfqpoint{0.000000in}{-0.048611in}}{\pgfqpoint{0.000000in}{0.000000in}}{%
\pgfpathmoveto{\pgfqpoint{0.000000in}{0.000000in}}%
\pgfpathlineto{\pgfqpoint{0.000000in}{-0.048611in}}%
\pgfusepath{stroke,fill}%
}%
\begin{pgfscope}%
\pgfsys@transformshift{2.318442in}{0.679475in}%
\pgfsys@useobject{currentmarker}{}%
\end{pgfscope}%
\end{pgfscope}%
\begin{pgfscope}%
\definecolor{textcolor}{rgb}{0.000000,0.000000,0.000000}%
\pgfsetstrokecolor{textcolor}%
\pgfsetfillcolor{textcolor}%
\pgftext[x=2.318442in,y=0.582253in,,top]{\color{textcolor}\rmfamily\fontsize{16.000000}{19.200000}\selectfont \(\displaystyle {0.2}\)}%
\end{pgfscope}%
\begin{pgfscope}%
\pgfsetbuttcap%
\pgfsetroundjoin%
\definecolor{currentfill}{rgb}{0.000000,0.000000,0.000000}%
\pgfsetfillcolor{currentfill}%
\pgfsetlinewidth{0.803000pt}%
\definecolor{currentstroke}{rgb}{0.000000,0.000000,0.000000}%
\pgfsetstrokecolor{currentstroke}%
\pgfsetdash{}{0pt}%
\pgfsys@defobject{currentmarker}{\pgfqpoint{0.000000in}{-0.048611in}}{\pgfqpoint{0.000000in}{0.000000in}}{%
\pgfpathmoveto{\pgfqpoint{0.000000in}{0.000000in}}%
\pgfpathlineto{\pgfqpoint{0.000000in}{-0.048611in}}%
\pgfusepath{stroke,fill}%
}%
\begin{pgfscope}%
\pgfsys@transformshift{3.510347in}{0.679475in}%
\pgfsys@useobject{currentmarker}{}%
\end{pgfscope}%
\end{pgfscope}%
\begin{pgfscope}%
\definecolor{textcolor}{rgb}{0.000000,0.000000,0.000000}%
\pgfsetstrokecolor{textcolor}%
\pgfsetfillcolor{textcolor}%
\pgftext[x=3.510347in,y=0.582253in,,top]{\color{textcolor}\rmfamily\fontsize{16.000000}{19.200000}\selectfont \(\displaystyle {0.4}\)}%
\end{pgfscope}%
\begin{pgfscope}%
\pgfsetbuttcap%
\pgfsetroundjoin%
\definecolor{currentfill}{rgb}{0.000000,0.000000,0.000000}%
\pgfsetfillcolor{currentfill}%
\pgfsetlinewidth{0.803000pt}%
\definecolor{currentstroke}{rgb}{0.000000,0.000000,0.000000}%
\pgfsetstrokecolor{currentstroke}%
\pgfsetdash{}{0pt}%
\pgfsys@defobject{currentmarker}{\pgfqpoint{0.000000in}{-0.048611in}}{\pgfqpoint{0.000000in}{0.000000in}}{%
\pgfpathmoveto{\pgfqpoint{0.000000in}{0.000000in}}%
\pgfpathlineto{\pgfqpoint{0.000000in}{-0.048611in}}%
\pgfusepath{stroke,fill}%
}%
\begin{pgfscope}%
\pgfsys@transformshift{4.702253in}{0.679475in}%
\pgfsys@useobject{currentmarker}{}%
\end{pgfscope}%
\end{pgfscope}%
\begin{pgfscope}%
\definecolor{textcolor}{rgb}{0.000000,0.000000,0.000000}%
\pgfsetstrokecolor{textcolor}%
\pgfsetfillcolor{textcolor}%
\pgftext[x=4.702253in,y=0.582253in,,top]{\color{textcolor}\rmfamily\fontsize{16.000000}{19.200000}\selectfont \(\displaystyle {0.6}\)}%
\end{pgfscope}%
\begin{pgfscope}%
\pgfsetbuttcap%
\pgfsetroundjoin%
\definecolor{currentfill}{rgb}{0.000000,0.000000,0.000000}%
\pgfsetfillcolor{currentfill}%
\pgfsetlinewidth{0.803000pt}%
\definecolor{currentstroke}{rgb}{0.000000,0.000000,0.000000}%
\pgfsetstrokecolor{currentstroke}%
\pgfsetdash{}{0pt}%
\pgfsys@defobject{currentmarker}{\pgfqpoint{0.000000in}{-0.048611in}}{\pgfqpoint{0.000000in}{0.000000in}}{%
\pgfpathmoveto{\pgfqpoint{0.000000in}{0.000000in}}%
\pgfpathlineto{\pgfqpoint{0.000000in}{-0.048611in}}%
\pgfusepath{stroke,fill}%
}%
\begin{pgfscope}%
\pgfsys@transformshift{5.894159in}{0.679475in}%
\pgfsys@useobject{currentmarker}{}%
\end{pgfscope}%
\end{pgfscope}%
\begin{pgfscope}%
\definecolor{textcolor}{rgb}{0.000000,0.000000,0.000000}%
\pgfsetstrokecolor{textcolor}%
\pgfsetfillcolor{textcolor}%
\pgftext[x=5.894159in,y=0.582253in,,top]{\color{textcolor}\rmfamily\fontsize{16.000000}{19.200000}\selectfont \(\displaystyle {0.8}\)}%
\end{pgfscope}%
\begin{pgfscope}%
\pgfsetbuttcap%
\pgfsetroundjoin%
\definecolor{currentfill}{rgb}{0.000000,0.000000,0.000000}%
\pgfsetfillcolor{currentfill}%
\pgfsetlinewidth{0.803000pt}%
\definecolor{currentstroke}{rgb}{0.000000,0.000000,0.000000}%
\pgfsetstrokecolor{currentstroke}%
\pgfsetdash{}{0pt}%
\pgfsys@defobject{currentmarker}{\pgfqpoint{0.000000in}{-0.048611in}}{\pgfqpoint{0.000000in}{0.000000in}}{%
\pgfpathmoveto{\pgfqpoint{0.000000in}{0.000000in}}%
\pgfpathlineto{\pgfqpoint{0.000000in}{-0.048611in}}%
\pgfusepath{stroke,fill}%
}%
\begin{pgfscope}%
\pgfsys@transformshift{7.086065in}{0.679475in}%
\pgfsys@useobject{currentmarker}{}%
\end{pgfscope}%
\end{pgfscope}%
\begin{pgfscope}%
\definecolor{textcolor}{rgb}{0.000000,0.000000,0.000000}%
\pgfsetstrokecolor{textcolor}%
\pgfsetfillcolor{textcolor}%
\pgftext[x=7.086065in,y=0.582253in,,top]{\color{textcolor}\rmfamily\fontsize{16.000000}{19.200000}\selectfont \(\displaystyle {1.0}\)}%
\end{pgfscope}%
\begin{pgfscope}%
\definecolor{textcolor}{rgb}{0.000000,0.000000,0.000000}%
\pgfsetstrokecolor{textcolor}%
\pgfsetfillcolor{textcolor}%
\pgftext[x=4.404277in,y=0.313349in,,top]{\color{textcolor}\rmfamily\fontsize{18.000000}{21.600000}\selectfont \(\displaystyle \lambda\)}%
\end{pgfscope}%
\begin{pgfscope}%
\pgfsetbuttcap%
\pgfsetroundjoin%
\definecolor{currentfill}{rgb}{0.000000,0.000000,0.000000}%
\pgfsetfillcolor{currentfill}%
\pgfsetlinewidth{0.803000pt}%
\definecolor{currentstroke}{rgb}{0.000000,0.000000,0.000000}%
\pgfsetstrokecolor{currentstroke}%
\pgfsetdash{}{0pt}%
\pgfsys@defobject{currentmarker}{\pgfqpoint{-0.048611in}{0.000000in}}{\pgfqpoint{-0.000000in}{0.000000in}}{%
\pgfpathmoveto{\pgfqpoint{-0.000000in}{0.000000in}}%
\pgfpathlineto{\pgfqpoint{-0.048611in}{0.000000in}}%
\pgfusepath{stroke,fill}%
}%
\begin{pgfscope}%
\pgfsys@transformshift{1.126536in}{0.728437in}%
\pgfsys@useobject{currentmarker}{}%
\end{pgfscope}%
\end{pgfscope}%
\begin{pgfscope}%
\definecolor{textcolor}{rgb}{0.000000,0.000000,0.000000}%
\pgfsetstrokecolor{textcolor}%
\pgfsetfillcolor{textcolor}%
\pgftext[x=0.919245in, y=0.645104in, left, base]{\color{textcolor}\rmfamily\fontsize{16.000000}{19.200000}\selectfont \(\displaystyle {0}\)}%
\end{pgfscope}%
\begin{pgfscope}%
\pgfsetbuttcap%
\pgfsetroundjoin%
\definecolor{currentfill}{rgb}{0.000000,0.000000,0.000000}%
\pgfsetfillcolor{currentfill}%
\pgfsetlinewidth{0.803000pt}%
\definecolor{currentstroke}{rgb}{0.000000,0.000000,0.000000}%
\pgfsetstrokecolor{currentstroke}%
\pgfsetdash{}{0pt}%
\pgfsys@defobject{currentmarker}{\pgfqpoint{-0.048611in}{0.000000in}}{\pgfqpoint{-0.000000in}{0.000000in}}{%
\pgfpathmoveto{\pgfqpoint{-0.000000in}{0.000000in}}%
\pgfpathlineto{\pgfqpoint{-0.048611in}{0.000000in}}%
\pgfusepath{stroke,fill}%
}%
\begin{pgfscope}%
\pgfsys@transformshift{1.126536in}{1.228810in}%
\pgfsys@useobject{currentmarker}{}%
\end{pgfscope}%
\end{pgfscope}%
\begin{pgfscope}%
\definecolor{textcolor}{rgb}{0.000000,0.000000,0.000000}%
\pgfsetstrokecolor{textcolor}%
\pgfsetfillcolor{textcolor}%
\pgftext[x=0.478972in, y=1.145477in, left, base]{\color{textcolor}\rmfamily\fontsize{16.000000}{19.200000}\selectfont \(\displaystyle {25000}\)}%
\end{pgfscope}%
\begin{pgfscope}%
\pgfsetbuttcap%
\pgfsetroundjoin%
\definecolor{currentfill}{rgb}{0.000000,0.000000,0.000000}%
\pgfsetfillcolor{currentfill}%
\pgfsetlinewidth{0.803000pt}%
\definecolor{currentstroke}{rgb}{0.000000,0.000000,0.000000}%
\pgfsetstrokecolor{currentstroke}%
\pgfsetdash{}{0pt}%
\pgfsys@defobject{currentmarker}{\pgfqpoint{-0.048611in}{0.000000in}}{\pgfqpoint{-0.000000in}{0.000000in}}{%
\pgfpathmoveto{\pgfqpoint{-0.000000in}{0.000000in}}%
\pgfpathlineto{\pgfqpoint{-0.048611in}{0.000000in}}%
\pgfusepath{stroke,fill}%
}%
\begin{pgfscope}%
\pgfsys@transformshift{1.126536in}{1.729184in}%
\pgfsys@useobject{currentmarker}{}%
\end{pgfscope}%
\end{pgfscope}%
\begin{pgfscope}%
\definecolor{textcolor}{rgb}{0.000000,0.000000,0.000000}%
\pgfsetstrokecolor{textcolor}%
\pgfsetfillcolor{textcolor}%
\pgftext[x=0.478972in, y=1.645850in, left, base]{\color{textcolor}\rmfamily\fontsize{16.000000}{19.200000}\selectfont \(\displaystyle {50000}\)}%
\end{pgfscope}%
\begin{pgfscope}%
\pgfsetbuttcap%
\pgfsetroundjoin%
\definecolor{currentfill}{rgb}{0.000000,0.000000,0.000000}%
\pgfsetfillcolor{currentfill}%
\pgfsetlinewidth{0.803000pt}%
\definecolor{currentstroke}{rgb}{0.000000,0.000000,0.000000}%
\pgfsetstrokecolor{currentstroke}%
\pgfsetdash{}{0pt}%
\pgfsys@defobject{currentmarker}{\pgfqpoint{-0.048611in}{0.000000in}}{\pgfqpoint{-0.000000in}{0.000000in}}{%
\pgfpathmoveto{\pgfqpoint{-0.000000in}{0.000000in}}%
\pgfpathlineto{\pgfqpoint{-0.048611in}{0.000000in}}%
\pgfusepath{stroke,fill}%
}%
\begin{pgfscope}%
\pgfsys@transformshift{1.126536in}{2.229557in}%
\pgfsys@useobject{currentmarker}{}%
\end{pgfscope}%
\end{pgfscope}%
\begin{pgfscope}%
\definecolor{textcolor}{rgb}{0.000000,0.000000,0.000000}%
\pgfsetstrokecolor{textcolor}%
\pgfsetfillcolor{textcolor}%
\pgftext[x=0.478972in, y=2.146224in, left, base]{\color{textcolor}\rmfamily\fontsize{16.000000}{19.200000}\selectfont \(\displaystyle {75000}\)}%
\end{pgfscope}%
\begin{pgfscope}%
\pgfsetbuttcap%
\pgfsetroundjoin%
\definecolor{currentfill}{rgb}{0.000000,0.000000,0.000000}%
\pgfsetfillcolor{currentfill}%
\pgfsetlinewidth{0.803000pt}%
\definecolor{currentstroke}{rgb}{0.000000,0.000000,0.000000}%
\pgfsetstrokecolor{currentstroke}%
\pgfsetdash{}{0pt}%
\pgfsys@defobject{currentmarker}{\pgfqpoint{-0.048611in}{0.000000in}}{\pgfqpoint{-0.000000in}{0.000000in}}{%
\pgfpathmoveto{\pgfqpoint{-0.000000in}{0.000000in}}%
\pgfpathlineto{\pgfqpoint{-0.048611in}{0.000000in}}%
\pgfusepath{stroke,fill}%
}%
\begin{pgfscope}%
\pgfsys@transformshift{1.126536in}{2.729930in}%
\pgfsys@useobject{currentmarker}{}%
\end{pgfscope}%
\end{pgfscope}%
\begin{pgfscope}%
\definecolor{textcolor}{rgb}{0.000000,0.000000,0.000000}%
\pgfsetstrokecolor{textcolor}%
\pgfsetfillcolor{textcolor}%
\pgftext[x=0.368904in, y=2.646597in, left, base]{\color{textcolor}\rmfamily\fontsize{16.000000}{19.200000}\selectfont \(\displaystyle {100000}\)}%
\end{pgfscope}%
\begin{pgfscope}%
\definecolor{textcolor}{rgb}{0.000000,0.000000,0.000000}%
\pgfsetstrokecolor{textcolor}%
\pgfsetfillcolor{textcolor}%
\pgftext[x=0.313349in,y=1.907686in,,bottom,rotate=90.000000]{\color{textcolor}\rmfamily\fontsize{18.000000}{21.600000}\selectfont DKK}%
\end{pgfscope}%
\begin{pgfscope}%
\pgfpathrectangle{\pgfqpoint{1.126536in}{0.679475in}}{\pgfqpoint{6.555482in}{2.456422in}}%
\pgfusepath{clip}%
\pgfsetrectcap%
\pgfsetroundjoin%
\pgfsetlinewidth{1.505625pt}%
\definecolor{currentstroke}{rgb}{0.121569,0.286275,0.490196}%
\pgfsetstrokecolor{currentstroke}%
\pgfsetdash{}{0pt}%
\pgfpathmoveto{\pgfqpoint{1.722489in}{1.187039in}}%
\pgfpathlineto{\pgfqpoint{2.616418in}{1.490456in}}%
\pgfpathlineto{\pgfqpoint{4.106300in}{2.040789in}}%
\pgfpathlineto{\pgfqpoint{5.596183in}{2.549220in}}%
\pgfpathlineto{\pgfqpoint{7.086065in}{3.024241in}}%
\pgfusepath{stroke}%
\end{pgfscope}%
\begin{pgfscope}%
\pgfpathrectangle{\pgfqpoint{1.126536in}{0.679475in}}{\pgfqpoint{6.555482in}{2.456422in}}%
\pgfusepath{clip}%
\pgfsetbuttcap%
\pgfsetmiterjoin%
\definecolor{currentfill}{rgb}{0.121569,0.286275,0.490196}%
\pgfsetfillcolor{currentfill}%
\pgfsetlinewidth{0.752812pt}%
\definecolor{currentstroke}{rgb}{1.000000,1.000000,1.000000}%
\pgfsetstrokecolor{currentstroke}%
\pgfsetdash{}{0pt}%
\pgfsys@defobject{currentmarker}{\pgfqpoint{-0.041667in}{-0.041667in}}{\pgfqpoint{0.041667in}{0.041667in}}{%
\pgfpathmoveto{\pgfqpoint{-0.041667in}{-0.041667in}}%
\pgfpathlineto{\pgfqpoint{0.041667in}{-0.041667in}}%
\pgfpathlineto{\pgfqpoint{0.041667in}{0.041667in}}%
\pgfpathlineto{\pgfqpoint{-0.041667in}{0.041667in}}%
\pgfpathlineto{\pgfqpoint{-0.041667in}{-0.041667in}}%
\pgfpathclose%
\pgfusepath{stroke,fill}%
}%
\begin{pgfscope}%
\pgfsys@transformshift{1.722489in}{1.187039in}%
\pgfsys@useobject{currentmarker}{}%
\end{pgfscope}%
\begin{pgfscope}%
\pgfsys@transformshift{2.616418in}{1.490456in}%
\pgfsys@useobject{currentmarker}{}%
\end{pgfscope}%
\begin{pgfscope}%
\pgfsys@transformshift{4.106300in}{2.040789in}%
\pgfsys@useobject{currentmarker}{}%
\end{pgfscope}%
\begin{pgfscope}%
\pgfsys@transformshift{5.596183in}{2.549220in}%
\pgfsys@useobject{currentmarker}{}%
\end{pgfscope}%
\begin{pgfscope}%
\pgfsys@transformshift{7.086065in}{3.024241in}%
\pgfsys@useobject{currentmarker}{}%
\end{pgfscope}%
\end{pgfscope}%
\begin{pgfscope}%
\pgfpathrectangle{\pgfqpoint{1.126536in}{0.679475in}}{\pgfqpoint{6.555482in}{2.456422in}}%
\pgfusepath{clip}%
\pgfsetbuttcap%
\pgfsetroundjoin%
\pgfsetlinewidth{1.505625pt}%
\definecolor{currentstroke}{rgb}{0.121569,0.286275,0.490196}%
\pgfsetstrokecolor{currentstroke}%
\pgfsetstrokeopacity{0.900000}%
\pgfsetdash{{7.500000pt}{7.500000pt}}{0.000000pt}%
\pgfpathmoveto{\pgfqpoint{1.722489in}{0.791131in}}%
\pgfpathlineto{\pgfqpoint{2.616418in}{1.139117in}}%
\pgfpathlineto{\pgfqpoint{4.106300in}{1.705699in}}%
\pgfpathlineto{\pgfqpoint{5.596183in}{2.245600in}}%
\pgfpathlineto{\pgfqpoint{7.086065in}{2.763533in}}%
\pgfusepath{stroke}%
\end{pgfscope}%
\begin{pgfscope}%
\pgfpathrectangle{\pgfqpoint{1.126536in}{0.679475in}}{\pgfqpoint{6.555482in}{2.456422in}}%
\pgfusepath{clip}%
\pgfsetbuttcap%
\pgfsetmiterjoin%
\definecolor{currentfill}{rgb}{0.121569,0.286275,0.490196}%
\pgfsetfillcolor{currentfill}%
\pgfsetfillopacity{0.900000}%
\pgfsetlinewidth{0.752812pt}%
\definecolor{currentstroke}{rgb}{1.000000,1.000000,1.000000}%
\pgfsetstrokecolor{currentstroke}%
\pgfsetdash{}{0pt}%
\pgfsys@defobject{currentmarker}{\pgfqpoint{-0.058926in}{-0.058926in}}{\pgfqpoint{0.058926in}{0.058926in}}{%
\pgfpathmoveto{\pgfqpoint{0.000000in}{-0.058926in}}%
\pgfpathlineto{\pgfqpoint{0.058926in}{0.000000in}}%
\pgfpathlineto{\pgfqpoint{0.000000in}{0.058926in}}%
\pgfpathlineto{\pgfqpoint{-0.058926in}{0.000000in}}%
\pgfpathlineto{\pgfqpoint{0.000000in}{-0.058926in}}%
\pgfpathclose%
\pgfusepath{stroke,fill}%
}%
\begin{pgfscope}%
\pgfsys@transformshift{1.722489in}{0.791131in}%
\pgfsys@useobject{currentmarker}{}%
\end{pgfscope}%
\begin{pgfscope}%
\pgfsys@transformshift{2.616418in}{1.139117in}%
\pgfsys@useobject{currentmarker}{}%
\end{pgfscope}%
\begin{pgfscope}%
\pgfsys@transformshift{4.106300in}{1.705699in}%
\pgfsys@useobject{currentmarker}{}%
\end{pgfscope}%
\begin{pgfscope}%
\pgfsys@transformshift{5.596183in}{2.245600in}%
\pgfsys@useobject{currentmarker}{}%
\end{pgfscope}%
\begin{pgfscope}%
\pgfsys@transformshift{7.086065in}{2.763533in}%
\pgfsys@useobject{currentmarker}{}%
\end{pgfscope}%
\end{pgfscope}%
\begin{pgfscope}%
\pgfpathrectangle{\pgfqpoint{1.126536in}{0.679475in}}{\pgfqpoint{6.555482in}{2.456422in}}%
\pgfusepath{clip}%
\pgfsetbuttcap%
\pgfsetroundjoin%
\pgfsetlinewidth{1.505625pt}%
\definecolor{currentstroke}{rgb}{0.121569,0.286275,0.490196}%
\pgfsetstrokecolor{currentstroke}%
\pgfsetstrokeopacity{0.800000}%
\pgfsetdash{{1.500000pt}{4.500000pt}}{0.000000pt}%
\pgfpathmoveto{\pgfqpoint{1.722489in}{1.114058in}}%
\pgfpathlineto{\pgfqpoint{2.616418in}{1.416123in}}%
\pgfpathlineto{\pgfqpoint{4.106300in}{1.926948in}}%
\pgfpathlineto{\pgfqpoint{5.596183in}{2.419785in}}%
\pgfpathlineto{\pgfqpoint{7.086065in}{2.913976in}}%
\pgfusepath{stroke}%
\end{pgfscope}%
\begin{pgfscope}%
\pgfpathrectangle{\pgfqpoint{1.126536in}{0.679475in}}{\pgfqpoint{6.555482in}{2.456422in}}%
\pgfusepath{clip}%
\pgfsetbuttcap%
\pgfsetroundjoin%
\definecolor{currentfill}{rgb}{0.121569,0.286275,0.490196}%
\pgfsetfillcolor{currentfill}%
\pgfsetfillopacity{0.800000}%
\pgfsetlinewidth{0.752812pt}%
\definecolor{currentstroke}{rgb}{1.000000,1.000000,1.000000}%
\pgfsetstrokecolor{currentstroke}%
\pgfsetdash{}{0pt}%
\pgfsys@defobject{currentmarker}{\pgfqpoint{-0.041667in}{-0.041667in}}{\pgfqpoint{0.041667in}{0.041667in}}{%
\pgfpathmoveto{\pgfqpoint{0.000000in}{-0.041667in}}%
\pgfpathcurveto{\pgfqpoint{0.011050in}{-0.041667in}}{\pgfqpoint{0.021649in}{-0.037276in}}{\pgfqpoint{0.029463in}{-0.029463in}}%
\pgfpathcurveto{\pgfqpoint{0.037276in}{-0.021649in}}{\pgfqpoint{0.041667in}{-0.011050in}}{\pgfqpoint{0.041667in}{0.000000in}}%
\pgfpathcurveto{\pgfqpoint{0.041667in}{0.011050in}}{\pgfqpoint{0.037276in}{0.021649in}}{\pgfqpoint{0.029463in}{0.029463in}}%
\pgfpathcurveto{\pgfqpoint{0.021649in}{0.037276in}}{\pgfqpoint{0.011050in}{0.041667in}}{\pgfqpoint{0.000000in}{0.041667in}}%
\pgfpathcurveto{\pgfqpoint{-0.011050in}{0.041667in}}{\pgfqpoint{-0.021649in}{0.037276in}}{\pgfqpoint{-0.029463in}{0.029463in}}%
\pgfpathcurveto{\pgfqpoint{-0.037276in}{0.021649in}}{\pgfqpoint{-0.041667in}{0.011050in}}{\pgfqpoint{-0.041667in}{0.000000in}}%
\pgfpathcurveto{\pgfqpoint{-0.041667in}{-0.011050in}}{\pgfqpoint{-0.037276in}{-0.021649in}}{\pgfqpoint{-0.029463in}{-0.029463in}}%
\pgfpathcurveto{\pgfqpoint{-0.021649in}{-0.037276in}}{\pgfqpoint{-0.011050in}{-0.041667in}}{\pgfqpoint{0.000000in}{-0.041667in}}%
\pgfpathlineto{\pgfqpoint{0.000000in}{-0.041667in}}%
\pgfpathclose%
\pgfusepath{stroke,fill}%
}%
\begin{pgfscope}%
\pgfsys@transformshift{1.722489in}{1.114058in}%
\pgfsys@useobject{currentmarker}{}%
\end{pgfscope}%
\begin{pgfscope}%
\pgfsys@transformshift{2.616418in}{1.416123in}%
\pgfsys@useobject{currentmarker}{}%
\end{pgfscope}%
\begin{pgfscope}%
\pgfsys@transformshift{4.106300in}{1.926948in}%
\pgfsys@useobject{currentmarker}{}%
\end{pgfscope}%
\begin{pgfscope}%
\pgfsys@transformshift{5.596183in}{2.419785in}%
\pgfsys@useobject{currentmarker}{}%
\end{pgfscope}%
\begin{pgfscope}%
\pgfsys@transformshift{7.086065in}{2.913976in}%
\pgfsys@useobject{currentmarker}{}%
\end{pgfscope}%
\end{pgfscope}%
\begin{pgfscope}%
\pgfsetrectcap%
\pgfsetmiterjoin%
\pgfsetlinewidth{0.803000pt}%
\definecolor{currentstroke}{rgb}{0.000000,0.000000,0.000000}%
\pgfsetstrokecolor{currentstroke}%
\pgfsetdash{}{0pt}%
\pgfpathmoveto{\pgfqpoint{1.126536in}{0.679475in}}%
\pgfpathlineto{\pgfqpoint{1.126536in}{3.135897in}}%
\pgfusepath{stroke}%
\end{pgfscope}%
\begin{pgfscope}%
\pgfsetrectcap%
\pgfsetmiterjoin%
\pgfsetlinewidth{0.803000pt}%
\definecolor{currentstroke}{rgb}{0.000000,0.000000,0.000000}%
\pgfsetstrokecolor{currentstroke}%
\pgfsetdash{}{0pt}%
\pgfpathmoveto{\pgfqpoint{7.682018in}{0.679475in}}%
\pgfpathlineto{\pgfqpoint{7.682018in}{3.135897in}}%
\pgfusepath{stroke}%
\end{pgfscope}%
\begin{pgfscope}%
\pgfsetrectcap%
\pgfsetmiterjoin%
\pgfsetlinewidth{0.803000pt}%
\definecolor{currentstroke}{rgb}{0.000000,0.000000,0.000000}%
\pgfsetstrokecolor{currentstroke}%
\pgfsetdash{}{0pt}%
\pgfpathmoveto{\pgfqpoint{1.126536in}{0.679475in}}%
\pgfpathlineto{\pgfqpoint{7.682018in}{0.679475in}}%
\pgfusepath{stroke}%
\end{pgfscope}%
\begin{pgfscope}%
\pgfsetrectcap%
\pgfsetmiterjoin%
\pgfsetlinewidth{0.803000pt}%
\definecolor{currentstroke}{rgb}{0.000000,0.000000,0.000000}%
\pgfsetstrokecolor{currentstroke}%
\pgfsetdash{}{0pt}%
\pgfpathmoveto{\pgfqpoint{1.126536in}{3.135897in}}%
\pgfpathlineto{\pgfqpoint{7.682018in}{3.135897in}}%
\pgfusepath{stroke}%
\end{pgfscope}%
\begin{pgfscope}%
\pgfsetrectcap%
\pgfsetroundjoin%
\pgfsetlinewidth{1.505625pt}%
\definecolor{currentstroke}{rgb}{0.121569,0.286275,0.490196}%
\pgfsetstrokecolor{currentstroke}%
\pgfsetdash{}{0pt}%
\pgfpathmoveto{\pgfqpoint{2.251710in}{3.363810in}}%
\pgfpathlineto{\pgfqpoint{2.501710in}{3.363810in}}%
\pgfpathlineto{\pgfqpoint{2.751710in}{3.363810in}}%
\pgfusepath{stroke}%
\end{pgfscope}%
\begin{pgfscope}%
\pgfsetbuttcap%
\pgfsetmiterjoin%
\definecolor{currentfill}{rgb}{0.121569,0.286275,0.490196}%
\pgfsetfillcolor{currentfill}%
\pgfsetlinewidth{1.003750pt}%
\definecolor{currentstroke}{rgb}{0.121569,0.286275,0.490196}%
\pgfsetstrokecolor{currentstroke}%
\pgfsetdash{}{0pt}%
\pgfsys@defobject{currentmarker}{\pgfqpoint{-0.041667in}{-0.041667in}}{\pgfqpoint{0.041667in}{0.041667in}}{%
\pgfpathmoveto{\pgfqpoint{-0.041667in}{-0.041667in}}%
\pgfpathlineto{\pgfqpoint{0.041667in}{-0.041667in}}%
\pgfpathlineto{\pgfqpoint{0.041667in}{0.041667in}}%
\pgfpathlineto{\pgfqpoint{-0.041667in}{0.041667in}}%
\pgfpathlineto{\pgfqpoint{-0.041667in}{-0.041667in}}%
\pgfpathclose%
\pgfusepath{stroke,fill}%
}%
\begin{pgfscope}%
\pgfsys@transformshift{2.501710in}{3.363810in}%
\pgfsys@useobject{currentmarker}{}%
\end{pgfscope}%
\end{pgfscope}%
\begin{pgfscope}%
\definecolor{textcolor}{rgb}{0.000000,0.000000,0.000000}%
\pgfsetstrokecolor{textcolor}%
\pgfsetfillcolor{textcolor}%
\pgftext[x=2.876710in,y=3.276310in,left,base]{\color{textcolor}\rmfamily\fontsize{18.000000}{21.600000}\selectfont LPP\(\displaystyle ^{R,T}\)}%
\end{pgfscope}%
\begin{pgfscope}%
\pgfsetbuttcap%
\pgfsetroundjoin%
\pgfsetlinewidth{1.505625pt}%
\definecolor{currentstroke}{rgb}{0.121569,0.286275,0.490196}%
\pgfsetstrokecolor{currentstroke}%
\pgfsetstrokeopacity{0.900000}%
\pgfsetdash{{5.550000pt}{2.400000pt}}{0.000000pt}%
\pgfpathmoveto{\pgfqpoint{3.812661in}{3.363810in}}%
\pgfpathlineto{\pgfqpoint{4.062661in}{3.363810in}}%
\pgfpathlineto{\pgfqpoint{4.312661in}{3.363810in}}%
\pgfusepath{stroke}%
\end{pgfscope}%
\begin{pgfscope}%
\pgfsetbuttcap%
\pgfsetmiterjoin%
\definecolor{currentfill}{rgb}{0.121569,0.286275,0.490196}%
\pgfsetfillcolor{currentfill}%
\pgfsetfillopacity{0.900000}%
\pgfsetlinewidth{1.003750pt}%
\definecolor{currentstroke}{rgb}{0.121569,0.286275,0.490196}%
\pgfsetstrokecolor{currentstroke}%
\pgfsetstrokeopacity{0.900000}%
\pgfsetdash{}{0pt}%
\pgfsys@defobject{currentmarker}{\pgfqpoint{-0.058926in}{-0.058926in}}{\pgfqpoint{0.058926in}{0.058926in}}{%
\pgfpathmoveto{\pgfqpoint{0.000000in}{-0.058926in}}%
\pgfpathlineto{\pgfqpoint{0.058926in}{0.000000in}}%
\pgfpathlineto{\pgfqpoint{0.000000in}{0.058926in}}%
\pgfpathlineto{\pgfqpoint{-0.058926in}{0.000000in}}%
\pgfpathlineto{\pgfqpoint{0.000000in}{-0.058926in}}%
\pgfpathclose%
\pgfusepath{stroke,fill}%
}%
\begin{pgfscope}%
\pgfsys@transformshift{4.062661in}{3.363810in}%
\pgfsys@useobject{currentmarker}{}%
\end{pgfscope}%
\end{pgfscope}%
\begin{pgfscope}%
\definecolor{textcolor}{rgb}{0.000000,0.000000,0.000000}%
\pgfsetstrokecolor{textcolor}%
\pgfsetfillcolor{textcolor}%
\pgftext[x=4.437661in,y=3.276310in,left,base]{\color{textcolor}\rmfamily\fontsize{18.000000}{21.600000}\selectfont LPP\(\displaystyle ^{T}\)}%
\end{pgfscope}%
\begin{pgfscope}%
\pgfsetbuttcap%
\pgfsetroundjoin%
\pgfsetlinewidth{1.505625pt}%
\definecolor{currentstroke}{rgb}{0.121569,0.286275,0.490196}%
\pgfsetstrokecolor{currentstroke}%
\pgfsetstrokeopacity{0.800000}%
\pgfsetdash{{1.500000pt}{2.475000pt}}{0.000000pt}%
\pgfpathmoveto{\pgfqpoint{5.202692in}{3.363810in}}%
\pgfpathlineto{\pgfqpoint{5.452692in}{3.363810in}}%
\pgfpathlineto{\pgfqpoint{5.702692in}{3.363810in}}%
\pgfusepath{stroke}%
\end{pgfscope}%
\begin{pgfscope}%
\pgfsetbuttcap%
\pgfsetroundjoin%
\definecolor{currentfill}{rgb}{0.121569,0.286275,0.490196}%
\pgfsetfillcolor{currentfill}%
\pgfsetfillopacity{0.800000}%
\pgfsetlinewidth{1.003750pt}%
\definecolor{currentstroke}{rgb}{0.121569,0.286275,0.490196}%
\pgfsetstrokecolor{currentstroke}%
\pgfsetstrokeopacity{0.800000}%
\pgfsetdash{}{0pt}%
\pgfsys@defobject{currentmarker}{\pgfqpoint{-0.041667in}{-0.041667in}}{\pgfqpoint{0.041667in}{0.041667in}}{%
\pgfpathmoveto{\pgfqpoint{0.000000in}{-0.041667in}}%
\pgfpathcurveto{\pgfqpoint{0.011050in}{-0.041667in}}{\pgfqpoint{0.021649in}{-0.037276in}}{\pgfqpoint{0.029463in}{-0.029463in}}%
\pgfpathcurveto{\pgfqpoint{0.037276in}{-0.021649in}}{\pgfqpoint{0.041667in}{-0.011050in}}{\pgfqpoint{0.041667in}{0.000000in}}%
\pgfpathcurveto{\pgfqpoint{0.041667in}{0.011050in}}{\pgfqpoint{0.037276in}{0.021649in}}{\pgfqpoint{0.029463in}{0.029463in}}%
\pgfpathcurveto{\pgfqpoint{0.021649in}{0.037276in}}{\pgfqpoint{0.011050in}{0.041667in}}{\pgfqpoint{0.000000in}{0.041667in}}%
\pgfpathcurveto{\pgfqpoint{-0.011050in}{0.041667in}}{\pgfqpoint{-0.021649in}{0.037276in}}{\pgfqpoint{-0.029463in}{0.029463in}}%
\pgfpathcurveto{\pgfqpoint{-0.037276in}{0.021649in}}{\pgfqpoint{-0.041667in}{0.011050in}}{\pgfqpoint{-0.041667in}{0.000000in}}%
\pgfpathcurveto{\pgfqpoint{-0.041667in}{-0.011050in}}{\pgfqpoint{-0.037276in}{-0.021649in}}{\pgfqpoint{-0.029463in}{-0.029463in}}%
\pgfpathcurveto{\pgfqpoint{-0.021649in}{-0.037276in}}{\pgfqpoint{-0.011050in}{-0.041667in}}{\pgfqpoint{0.000000in}{-0.041667in}}%
\pgfpathlineto{\pgfqpoint{0.000000in}{-0.041667in}}%
\pgfpathclose%
\pgfusepath{stroke,fill}%
}%
\begin{pgfscope}%
\pgfsys@transformshift{5.452692in}{3.363810in}%
\pgfsys@useobject{currentmarker}{}%
\end{pgfscope}%
\end{pgfscope}%
\begin{pgfscope}%
\definecolor{textcolor}{rgb}{0.000000,0.000000,0.000000}%
\pgfsetstrokecolor{textcolor}%
\pgfsetfillcolor{textcolor}%
\pgftext[x=5.827692in,y=3.276310in,left,base]{\color{textcolor}\rmfamily\fontsize{18.000000}{21.600000}\selectfont LPP\(\displaystyle ^{T, P}\)}%
\end{pgfscope}%
\end{pgfpicture}%
\makeatother%
\endgroup%

%% file: figures/DemandPlots/TightLayout/DE_UB_rescaled_Gap10pct.pgf
\begingroup%
\makeatletter%
\begin{pgfpicture}%
\pgfpathrectangle{\pgfpointorigin}{\pgfqpoint{7.782018in}{3.651310in}}%
\pgfusepath{use as bounding box, clip}%
\begin{pgfscope}%
\pgfsetbuttcap%
\pgfsetmiterjoin%
\definecolor{currentfill}{rgb}{1.000000,1.000000,1.000000}%
\pgfsetfillcolor{currentfill}%
\pgfsetlinewidth{0.000000pt}%
\definecolor{currentstroke}{rgb}{1.000000,1.000000,1.000000}%
\pgfsetstrokecolor{currentstroke}%
\pgfsetdash{}{0pt}%
\pgfpathmoveto{\pgfqpoint{0.000000in}{0.000000in}}%
\pgfpathlineto{\pgfqpoint{7.782018in}{0.000000in}}%
\pgfpathlineto{\pgfqpoint{7.782018in}{3.651310in}}%
\pgfpathlineto{\pgfqpoint{0.000000in}{3.651310in}}%
\pgfpathlineto{\pgfqpoint{0.000000in}{0.000000in}}%
\pgfpathclose%
\pgfusepath{fill}%
\end{pgfscope}%
\begin{pgfscope}%
\pgfsetbuttcap%
\pgfsetmiterjoin%
\definecolor{currentfill}{rgb}{1.000000,1.000000,1.000000}%
\pgfsetfillcolor{currentfill}%
\pgfsetlinewidth{0.000000pt}%
\definecolor{currentstroke}{rgb}{0.000000,0.000000,0.000000}%
\pgfsetstrokecolor{currentstroke}%
\pgfsetstrokeopacity{0.000000}%
\pgfsetdash{}{0pt}%
\pgfpathmoveto{\pgfqpoint{1.126536in}{0.679475in}}%
\pgfpathlineto{\pgfqpoint{7.682018in}{0.679475in}}%
\pgfpathlineto{\pgfqpoint{7.682018in}{3.135897in}}%
\pgfpathlineto{\pgfqpoint{1.126536in}{3.135897in}}%
\pgfpathlineto{\pgfqpoint{1.126536in}{0.679475in}}%
\pgfpathclose%
\pgfusepath{fill}%
\end{pgfscope}%
\begin{pgfscope}%
\pgfsetbuttcap%
\pgfsetroundjoin%
\definecolor{currentfill}{rgb}{0.000000,0.000000,0.000000}%
\pgfsetfillcolor{currentfill}%
\pgfsetlinewidth{0.803000pt}%
\definecolor{currentstroke}{rgb}{0.000000,0.000000,0.000000}%
\pgfsetstrokecolor{currentstroke}%
\pgfsetdash{}{0pt}%
\pgfsys@defobject{currentmarker}{\pgfqpoint{0.000000in}{-0.048611in}}{\pgfqpoint{0.000000in}{0.000000in}}{%
\pgfpathmoveto{\pgfqpoint{0.000000in}{0.000000in}}%
\pgfpathlineto{\pgfqpoint{0.000000in}{-0.048611in}}%
\pgfusepath{stroke,fill}%
}%
\begin{pgfscope}%
\pgfsys@transformshift{1.126536in}{0.679475in}%
\pgfsys@useobject{currentmarker}{}%
\end{pgfscope}%
\end{pgfscope}%
\begin{pgfscope}%
\definecolor{textcolor}{rgb}{0.000000,0.000000,0.000000}%
\pgfsetstrokecolor{textcolor}%
\pgfsetfillcolor{textcolor}%
\pgftext[x=1.126536in,y=0.582253in,,top]{\color{textcolor}\rmfamily\fontsize{16.000000}{19.200000}\selectfont \(\displaystyle {0.0}\)}%
\end{pgfscope}%
\begin{pgfscope}%
\pgfsetbuttcap%
\pgfsetroundjoin%
\definecolor{currentfill}{rgb}{0.000000,0.000000,0.000000}%
\pgfsetfillcolor{currentfill}%
\pgfsetlinewidth{0.803000pt}%
\definecolor{currentstroke}{rgb}{0.000000,0.000000,0.000000}%
\pgfsetstrokecolor{currentstroke}%
\pgfsetdash{}{0pt}%
\pgfsys@defobject{currentmarker}{\pgfqpoint{0.000000in}{-0.048611in}}{\pgfqpoint{0.000000in}{0.000000in}}{%
\pgfpathmoveto{\pgfqpoint{0.000000in}{0.000000in}}%
\pgfpathlineto{\pgfqpoint{0.000000in}{-0.048611in}}%
\pgfusepath{stroke,fill}%
}%
\begin{pgfscope}%
\pgfsys@transformshift{2.318442in}{0.679475in}%
\pgfsys@useobject{currentmarker}{}%
\end{pgfscope}%
\end{pgfscope}%
\begin{pgfscope}%
\definecolor{textcolor}{rgb}{0.000000,0.000000,0.000000}%
\pgfsetstrokecolor{textcolor}%
\pgfsetfillcolor{textcolor}%
\pgftext[x=2.318442in,y=0.582253in,,top]{\color{textcolor}\rmfamily\fontsize{16.000000}{19.200000}\selectfont \(\displaystyle {0.2}\)}%
\end{pgfscope}%
\begin{pgfscope}%
\pgfsetbuttcap%
\pgfsetroundjoin%
\definecolor{currentfill}{rgb}{0.000000,0.000000,0.000000}%
\pgfsetfillcolor{currentfill}%
\pgfsetlinewidth{0.803000pt}%
\definecolor{currentstroke}{rgb}{0.000000,0.000000,0.000000}%
\pgfsetstrokecolor{currentstroke}%
\pgfsetdash{}{0pt}%
\pgfsys@defobject{currentmarker}{\pgfqpoint{0.000000in}{-0.048611in}}{\pgfqpoint{0.000000in}{0.000000in}}{%
\pgfpathmoveto{\pgfqpoint{0.000000in}{0.000000in}}%
\pgfpathlineto{\pgfqpoint{0.000000in}{-0.048611in}}%
\pgfusepath{stroke,fill}%
}%
\begin{pgfscope}%
\pgfsys@transformshift{3.510347in}{0.679475in}%
\pgfsys@useobject{currentmarker}{}%
\end{pgfscope}%
\end{pgfscope}%
\begin{pgfscope}%
\definecolor{textcolor}{rgb}{0.000000,0.000000,0.000000}%
\pgfsetstrokecolor{textcolor}%
\pgfsetfillcolor{textcolor}%
\pgftext[x=3.510347in,y=0.582253in,,top]{\color{textcolor}\rmfamily\fontsize{16.000000}{19.200000}\selectfont \(\displaystyle {0.4}\)}%
\end{pgfscope}%
\begin{pgfscope}%
\pgfsetbuttcap%
\pgfsetroundjoin%
\definecolor{currentfill}{rgb}{0.000000,0.000000,0.000000}%
\pgfsetfillcolor{currentfill}%
\pgfsetlinewidth{0.803000pt}%
\definecolor{currentstroke}{rgb}{0.000000,0.000000,0.000000}%
\pgfsetstrokecolor{currentstroke}%
\pgfsetdash{}{0pt}%
\pgfsys@defobject{currentmarker}{\pgfqpoint{0.000000in}{-0.048611in}}{\pgfqpoint{0.000000in}{0.000000in}}{%
\pgfpathmoveto{\pgfqpoint{0.000000in}{0.000000in}}%
\pgfpathlineto{\pgfqpoint{0.000000in}{-0.048611in}}%
\pgfusepath{stroke,fill}%
}%
\begin{pgfscope}%
\pgfsys@transformshift{4.702253in}{0.679475in}%
\pgfsys@useobject{currentmarker}{}%
\end{pgfscope}%
\end{pgfscope}%
\begin{pgfscope}%
\definecolor{textcolor}{rgb}{0.000000,0.000000,0.000000}%
\pgfsetstrokecolor{textcolor}%
\pgfsetfillcolor{textcolor}%
\pgftext[x=4.702253in,y=0.582253in,,top]{\color{textcolor}\rmfamily\fontsize{16.000000}{19.200000}\selectfont \(\displaystyle {0.6}\)}%
\end{pgfscope}%
\begin{pgfscope}%
\pgfsetbuttcap%
\pgfsetroundjoin%
\definecolor{currentfill}{rgb}{0.000000,0.000000,0.000000}%
\pgfsetfillcolor{currentfill}%
\pgfsetlinewidth{0.803000pt}%
\definecolor{currentstroke}{rgb}{0.000000,0.000000,0.000000}%
\pgfsetstrokecolor{currentstroke}%
\pgfsetdash{}{0pt}%
\pgfsys@defobject{currentmarker}{\pgfqpoint{0.000000in}{-0.048611in}}{\pgfqpoint{0.000000in}{0.000000in}}{%
\pgfpathmoveto{\pgfqpoint{0.000000in}{0.000000in}}%
\pgfpathlineto{\pgfqpoint{0.000000in}{-0.048611in}}%
\pgfusepath{stroke,fill}%
}%
\begin{pgfscope}%
\pgfsys@transformshift{5.894159in}{0.679475in}%
\pgfsys@useobject{currentmarker}{}%
\end{pgfscope}%
\end{pgfscope}%
\begin{pgfscope}%
\definecolor{textcolor}{rgb}{0.000000,0.000000,0.000000}%
\pgfsetstrokecolor{textcolor}%
\pgfsetfillcolor{textcolor}%
\pgftext[x=5.894159in,y=0.582253in,,top]{\color{textcolor}\rmfamily\fontsize{16.000000}{19.200000}\selectfont \(\displaystyle {0.8}\)}%
\end{pgfscope}%
\begin{pgfscope}%
\pgfsetbuttcap%
\pgfsetroundjoin%
\definecolor{currentfill}{rgb}{0.000000,0.000000,0.000000}%
\pgfsetfillcolor{currentfill}%
\pgfsetlinewidth{0.803000pt}%
\definecolor{currentstroke}{rgb}{0.000000,0.000000,0.000000}%
\pgfsetstrokecolor{currentstroke}%
\pgfsetdash{}{0pt}%
\pgfsys@defobject{currentmarker}{\pgfqpoint{0.000000in}{-0.048611in}}{\pgfqpoint{0.000000in}{0.000000in}}{%
\pgfpathmoveto{\pgfqpoint{0.000000in}{0.000000in}}%
\pgfpathlineto{\pgfqpoint{0.000000in}{-0.048611in}}%
\pgfusepath{stroke,fill}%
}%
\begin{pgfscope}%
\pgfsys@transformshift{7.086065in}{0.679475in}%
\pgfsys@useobject{currentmarker}{}%
\end{pgfscope}%
\end{pgfscope}%
\begin{pgfscope}%
\definecolor{textcolor}{rgb}{0.000000,0.000000,0.000000}%
\pgfsetstrokecolor{textcolor}%
\pgfsetfillcolor{textcolor}%
\pgftext[x=7.086065in,y=0.582253in,,top]{\color{textcolor}\rmfamily\fontsize{16.000000}{19.200000}\selectfont \(\displaystyle {1.0}\)}%
\end{pgfscope}%
\begin{pgfscope}%
\definecolor{textcolor}{rgb}{0.000000,0.000000,0.000000}%
\pgfsetstrokecolor{textcolor}%
\pgfsetfillcolor{textcolor}%
\pgftext[x=4.404277in,y=0.313349in,,top]{\color{textcolor}\rmfamily\fontsize{18.000000}{21.600000}\selectfont \(\displaystyle \lambda\)}%
\end{pgfscope}%
\begin{pgfscope}%
\pgfsetbuttcap%
\pgfsetroundjoin%
\definecolor{currentfill}{rgb}{0.000000,0.000000,0.000000}%
\pgfsetfillcolor{currentfill}%
\pgfsetlinewidth{0.803000pt}%
\definecolor{currentstroke}{rgb}{0.000000,0.000000,0.000000}%
\pgfsetstrokecolor{currentstroke}%
\pgfsetdash{}{0pt}%
\pgfsys@defobject{currentmarker}{\pgfqpoint{-0.048611in}{0.000000in}}{\pgfqpoint{-0.000000in}{0.000000in}}{%
\pgfpathmoveto{\pgfqpoint{-0.000000in}{0.000000in}}%
\pgfpathlineto{\pgfqpoint{-0.048611in}{0.000000in}}%
\pgfusepath{stroke,fill}%
}%
\begin{pgfscope}%
\pgfsys@transformshift{1.126536in}{1.104297in}%
\pgfsys@useobject{currentmarker}{}%
\end{pgfscope}%
\end{pgfscope}%
\begin{pgfscope}%
\definecolor{textcolor}{rgb}{0.000000,0.000000,0.000000}%
\pgfsetstrokecolor{textcolor}%
\pgfsetfillcolor{textcolor}%
\pgftext[x=0.368904in, y=1.020963in, left, base]{\color{textcolor}\rmfamily\fontsize{16.000000}{19.200000}\selectfont \(\displaystyle {105000}\)}%
\end{pgfscope}%
\begin{pgfscope}%
\pgfsetbuttcap%
\pgfsetroundjoin%
\definecolor{currentfill}{rgb}{0.000000,0.000000,0.000000}%
\pgfsetfillcolor{currentfill}%
\pgfsetlinewidth{0.803000pt}%
\definecolor{currentstroke}{rgb}{0.000000,0.000000,0.000000}%
\pgfsetstrokecolor{currentstroke}%
\pgfsetdash{}{0pt}%
\pgfsys@defobject{currentmarker}{\pgfqpoint{-0.048611in}{0.000000in}}{\pgfqpoint{-0.000000in}{0.000000in}}{%
\pgfpathmoveto{\pgfqpoint{-0.000000in}{0.000000in}}%
\pgfpathlineto{\pgfqpoint{-0.048611in}{0.000000in}}%
\pgfusepath{stroke,fill}%
}%
\begin{pgfscope}%
\pgfsys@transformshift{1.126536in}{1.575778in}%
\pgfsys@useobject{currentmarker}{}%
\end{pgfscope}%
\end{pgfscope}%
\begin{pgfscope}%
\definecolor{textcolor}{rgb}{0.000000,0.000000,0.000000}%
\pgfsetstrokecolor{textcolor}%
\pgfsetfillcolor{textcolor}%
\pgftext[x=0.368904in, y=1.492445in, left, base]{\color{textcolor}\rmfamily\fontsize{16.000000}{19.200000}\selectfont \(\displaystyle {110000}\)}%
\end{pgfscope}%
\begin{pgfscope}%
\pgfsetbuttcap%
\pgfsetroundjoin%
\definecolor{currentfill}{rgb}{0.000000,0.000000,0.000000}%
\pgfsetfillcolor{currentfill}%
\pgfsetlinewidth{0.803000pt}%
\definecolor{currentstroke}{rgb}{0.000000,0.000000,0.000000}%
\pgfsetstrokecolor{currentstroke}%
\pgfsetdash{}{0pt}%
\pgfsys@defobject{currentmarker}{\pgfqpoint{-0.048611in}{0.000000in}}{\pgfqpoint{-0.000000in}{0.000000in}}{%
\pgfpathmoveto{\pgfqpoint{-0.000000in}{0.000000in}}%
\pgfpathlineto{\pgfqpoint{-0.048611in}{0.000000in}}%
\pgfusepath{stroke,fill}%
}%
\begin{pgfscope}%
\pgfsys@transformshift{1.126536in}{2.047260in}%
\pgfsys@useobject{currentmarker}{}%
\end{pgfscope}%
\end{pgfscope}%
\begin{pgfscope}%
\definecolor{textcolor}{rgb}{0.000000,0.000000,0.000000}%
\pgfsetstrokecolor{textcolor}%
\pgfsetfillcolor{textcolor}%
\pgftext[x=0.368904in, y=1.963927in, left, base]{\color{textcolor}\rmfamily\fontsize{16.000000}{19.200000}\selectfont \(\displaystyle {115000}\)}%
\end{pgfscope}%
\begin{pgfscope}%
\pgfsetbuttcap%
\pgfsetroundjoin%
\definecolor{currentfill}{rgb}{0.000000,0.000000,0.000000}%
\pgfsetfillcolor{currentfill}%
\pgfsetlinewidth{0.803000pt}%
\definecolor{currentstroke}{rgb}{0.000000,0.000000,0.000000}%
\pgfsetstrokecolor{currentstroke}%
\pgfsetdash{}{0pt}%
\pgfsys@defobject{currentmarker}{\pgfqpoint{-0.048611in}{0.000000in}}{\pgfqpoint{-0.000000in}{0.000000in}}{%
\pgfpathmoveto{\pgfqpoint{-0.000000in}{0.000000in}}%
\pgfpathlineto{\pgfqpoint{-0.048611in}{0.000000in}}%
\pgfusepath{stroke,fill}%
}%
\begin{pgfscope}%
\pgfsys@transformshift{1.126536in}{2.518742in}%
\pgfsys@useobject{currentmarker}{}%
\end{pgfscope}%
\end{pgfscope}%
\begin{pgfscope}%
\definecolor{textcolor}{rgb}{0.000000,0.000000,0.000000}%
\pgfsetstrokecolor{textcolor}%
\pgfsetfillcolor{textcolor}%
\pgftext[x=0.368904in, y=2.435409in, left, base]{\color{textcolor}\rmfamily\fontsize{16.000000}{19.200000}\selectfont \(\displaystyle {120000}\)}%
\end{pgfscope}%
\begin{pgfscope}%
\pgfsetbuttcap%
\pgfsetroundjoin%
\definecolor{currentfill}{rgb}{0.000000,0.000000,0.000000}%
\pgfsetfillcolor{currentfill}%
\pgfsetlinewidth{0.803000pt}%
\definecolor{currentstroke}{rgb}{0.000000,0.000000,0.000000}%
\pgfsetstrokecolor{currentstroke}%
\pgfsetdash{}{0pt}%
\pgfsys@defobject{currentmarker}{\pgfqpoint{-0.048611in}{0.000000in}}{\pgfqpoint{-0.000000in}{0.000000in}}{%
\pgfpathmoveto{\pgfqpoint{-0.000000in}{0.000000in}}%
\pgfpathlineto{\pgfqpoint{-0.048611in}{0.000000in}}%
\pgfusepath{stroke,fill}%
}%
\begin{pgfscope}%
\pgfsys@transformshift{1.126536in}{2.990224in}%
\pgfsys@useobject{currentmarker}{}%
\end{pgfscope}%
\end{pgfscope}%
\begin{pgfscope}%
\definecolor{textcolor}{rgb}{0.000000,0.000000,0.000000}%
\pgfsetstrokecolor{textcolor}%
\pgfsetfillcolor{textcolor}%
\pgftext[x=0.368904in, y=2.906891in, left, base]{\color{textcolor}\rmfamily\fontsize{16.000000}{19.200000}\selectfont \(\displaystyle {125000}\)}%
\end{pgfscope}%
\begin{pgfscope}%
\definecolor{textcolor}{rgb}{0.000000,0.000000,0.000000}%
\pgfsetstrokecolor{textcolor}%
\pgfsetfillcolor{textcolor}%
\pgftext[x=0.313349in,y=1.907686in,,bottom,rotate=90.000000]{\color{textcolor}\rmfamily\fontsize{18.000000}{21.600000}\selectfont DKK}%
\end{pgfscope}%
\begin{pgfscope}%
\pgfpathrectangle{\pgfqpoint{1.126536in}{0.679475in}}{\pgfqpoint{6.555482in}{2.456422in}}%
\pgfusepath{clip}%
\pgfsetrectcap%
\pgfsetroundjoin%
\pgfsetlinewidth{1.505625pt}%
\definecolor{currentstroke}{rgb}{0.337255,0.705882,0.913725}%
\pgfsetstrokecolor{currentstroke}%
\pgfsetdash{}{0pt}%
\pgfpathmoveto{\pgfqpoint{1.722489in}{3.024241in}}%
\pgfpathlineto{\pgfqpoint{2.616418in}{2.706345in}}%
\pgfpathlineto{\pgfqpoint{4.106300in}{2.491200in}}%
\pgfpathlineto{\pgfqpoint{5.596183in}{2.222696in}}%
\pgfpathlineto{\pgfqpoint{7.086065in}{2.019404in}}%
\pgfusepath{stroke}%
\end{pgfscope}%
\begin{pgfscope}%
\pgfpathrectangle{\pgfqpoint{1.126536in}{0.679475in}}{\pgfqpoint{6.555482in}{2.456422in}}%
\pgfusepath{clip}%
\pgfsetbuttcap%
\pgfsetmiterjoin%
\definecolor{currentfill}{rgb}{0.337255,0.705882,0.913725}%
\pgfsetfillcolor{currentfill}%
\pgfsetlinewidth{0.752812pt}%
\definecolor{currentstroke}{rgb}{1.000000,1.000000,1.000000}%
\pgfsetstrokecolor{currentstroke}%
\pgfsetdash{}{0pt}%
\pgfsys@defobject{currentmarker}{\pgfqpoint{-0.041667in}{-0.041667in}}{\pgfqpoint{0.041667in}{0.041667in}}{%
\pgfpathmoveto{\pgfqpoint{-0.041667in}{-0.041667in}}%
\pgfpathlineto{\pgfqpoint{0.041667in}{-0.041667in}}%
\pgfpathlineto{\pgfqpoint{0.041667in}{0.041667in}}%
\pgfpathlineto{\pgfqpoint{-0.041667in}{0.041667in}}%
\pgfpathlineto{\pgfqpoint{-0.041667in}{-0.041667in}}%
\pgfpathclose%
\pgfusepath{stroke,fill}%
}%
\begin{pgfscope}%
\pgfsys@transformshift{1.722489in}{3.024241in}%
\pgfsys@useobject{currentmarker}{}%
\end{pgfscope}%
\begin{pgfscope}%
\pgfsys@transformshift{2.616418in}{2.706345in}%
\pgfsys@useobject{currentmarker}{}%
\end{pgfscope}%
\begin{pgfscope}%
\pgfsys@transformshift{4.106300in}{2.491200in}%
\pgfsys@useobject{currentmarker}{}%
\end{pgfscope}%
\begin{pgfscope}%
\pgfsys@transformshift{5.596183in}{2.222696in}%
\pgfsys@useobject{currentmarker}{}%
\end{pgfscope}%
\begin{pgfscope}%
\pgfsys@transformshift{7.086065in}{2.019404in}%
\pgfsys@useobject{currentmarker}{}%
\end{pgfscope}%
\end{pgfscope}%
\begin{pgfscope}%
\pgfpathrectangle{\pgfqpoint{1.126536in}{0.679475in}}{\pgfqpoint{6.555482in}{2.456422in}}%
\pgfusepath{clip}%
\pgfsetbuttcap%
\pgfsetroundjoin%
\pgfsetlinewidth{1.505625pt}%
\definecolor{currentstroke}{rgb}{0.337255,0.705882,0.913725}%
\pgfsetstrokecolor{currentstroke}%
\pgfsetstrokeopacity{0.900000}%
\pgfsetdash{{7.500000pt}{7.500000pt}}{0.000000pt}%
\pgfpathmoveto{\pgfqpoint{1.722489in}{1.529227in}}%
\pgfpathlineto{\pgfqpoint{2.616418in}{1.211457in}}%
\pgfpathlineto{\pgfqpoint{4.106300in}{1.036647in}}%
\pgfpathlineto{\pgfqpoint{5.596183in}{0.860554in}}%
\pgfpathlineto{\pgfqpoint{7.086065in}{0.791131in}}%
\pgfusepath{stroke}%
\end{pgfscope}%
\begin{pgfscope}%
\pgfpathrectangle{\pgfqpoint{1.126536in}{0.679475in}}{\pgfqpoint{6.555482in}{2.456422in}}%
\pgfusepath{clip}%
\pgfsetbuttcap%
\pgfsetmiterjoin%
\definecolor{currentfill}{rgb}{0.337255,0.705882,0.913725}%
\pgfsetfillcolor{currentfill}%
\pgfsetfillopacity{0.900000}%
\pgfsetlinewidth{0.752812pt}%
\definecolor{currentstroke}{rgb}{1.000000,1.000000,1.000000}%
\pgfsetstrokecolor{currentstroke}%
\pgfsetdash{}{0pt}%
\pgfsys@defobject{currentmarker}{\pgfqpoint{-0.058926in}{-0.058926in}}{\pgfqpoint{0.058926in}{0.058926in}}{%
\pgfpathmoveto{\pgfqpoint{0.000000in}{-0.058926in}}%
\pgfpathlineto{\pgfqpoint{0.058926in}{0.000000in}}%
\pgfpathlineto{\pgfqpoint{0.000000in}{0.058926in}}%
\pgfpathlineto{\pgfqpoint{-0.058926in}{0.000000in}}%
\pgfpathlineto{\pgfqpoint{0.000000in}{-0.058926in}}%
\pgfpathclose%
\pgfusepath{stroke,fill}%
}%
\begin{pgfscope}%
\pgfsys@transformshift{1.722489in}{1.529227in}%
\pgfsys@useobject{currentmarker}{}%
\end{pgfscope}%
\begin{pgfscope}%
\pgfsys@transformshift{2.616418in}{1.211457in}%
\pgfsys@useobject{currentmarker}{}%
\end{pgfscope}%
\begin{pgfscope}%
\pgfsys@transformshift{4.106300in}{1.036647in}%
\pgfsys@useobject{currentmarker}{}%
\end{pgfscope}%
\begin{pgfscope}%
\pgfsys@transformshift{5.596183in}{0.860554in}%
\pgfsys@useobject{currentmarker}{}%
\end{pgfscope}%
\begin{pgfscope}%
\pgfsys@transformshift{7.086065in}{0.791131in}%
\pgfsys@useobject{currentmarker}{}%
\end{pgfscope}%
\end{pgfscope}%
\begin{pgfscope}%
\pgfpathrectangle{\pgfqpoint{1.126536in}{0.679475in}}{\pgfqpoint{6.555482in}{2.456422in}}%
\pgfusepath{clip}%
\pgfsetbuttcap%
\pgfsetroundjoin%
\pgfsetlinewidth{1.505625pt}%
\definecolor{currentstroke}{rgb}{0.337255,0.705882,0.913725}%
\pgfsetstrokecolor{currentstroke}%
\pgfsetstrokeopacity{0.800000}%
\pgfsetdash{{1.500000pt}{4.500000pt}}{0.000000pt}%
\pgfpathmoveto{\pgfqpoint{1.722489in}{1.422680in}}%
\pgfpathlineto{\pgfqpoint{2.616418in}{1.436896in}}%
\pgfpathlineto{\pgfqpoint{4.106300in}{1.432275in}}%
\pgfpathlineto{\pgfqpoint{5.596183in}{1.402802in}}%
\pgfpathlineto{\pgfqpoint{7.086065in}{1.499909in}}%
\pgfusepath{stroke}%
\end{pgfscope}%
\begin{pgfscope}%
\pgfpathrectangle{\pgfqpoint{1.126536in}{0.679475in}}{\pgfqpoint{6.555482in}{2.456422in}}%
\pgfusepath{clip}%
\pgfsetbuttcap%
\pgfsetroundjoin%
\definecolor{currentfill}{rgb}{0.337255,0.705882,0.913725}%
\pgfsetfillcolor{currentfill}%
\pgfsetfillopacity{0.800000}%
\pgfsetlinewidth{0.752812pt}%
\definecolor{currentstroke}{rgb}{1.000000,1.000000,1.000000}%
\pgfsetstrokecolor{currentstroke}%
\pgfsetdash{}{0pt}%
\pgfsys@defobject{currentmarker}{\pgfqpoint{-0.041667in}{-0.041667in}}{\pgfqpoint{0.041667in}{0.041667in}}{%
\pgfpathmoveto{\pgfqpoint{0.000000in}{-0.041667in}}%
\pgfpathcurveto{\pgfqpoint{0.011050in}{-0.041667in}}{\pgfqpoint{0.021649in}{-0.037276in}}{\pgfqpoint{0.029463in}{-0.029463in}}%
\pgfpathcurveto{\pgfqpoint{0.037276in}{-0.021649in}}{\pgfqpoint{0.041667in}{-0.011050in}}{\pgfqpoint{0.041667in}{0.000000in}}%
\pgfpathcurveto{\pgfqpoint{0.041667in}{0.011050in}}{\pgfqpoint{0.037276in}{0.021649in}}{\pgfqpoint{0.029463in}{0.029463in}}%
\pgfpathcurveto{\pgfqpoint{0.021649in}{0.037276in}}{\pgfqpoint{0.011050in}{0.041667in}}{\pgfqpoint{0.000000in}{0.041667in}}%
\pgfpathcurveto{\pgfqpoint{-0.011050in}{0.041667in}}{\pgfqpoint{-0.021649in}{0.037276in}}{\pgfqpoint{-0.029463in}{0.029463in}}%
\pgfpathcurveto{\pgfqpoint{-0.037276in}{0.021649in}}{\pgfqpoint{-0.041667in}{0.011050in}}{\pgfqpoint{-0.041667in}{0.000000in}}%
\pgfpathcurveto{\pgfqpoint{-0.041667in}{-0.011050in}}{\pgfqpoint{-0.037276in}{-0.021649in}}{\pgfqpoint{-0.029463in}{-0.029463in}}%
\pgfpathcurveto{\pgfqpoint{-0.021649in}{-0.037276in}}{\pgfqpoint{-0.011050in}{-0.041667in}}{\pgfqpoint{0.000000in}{-0.041667in}}%
\pgfpathlineto{\pgfqpoint{0.000000in}{-0.041667in}}%
\pgfpathclose%
\pgfusepath{stroke,fill}%
}%
\begin{pgfscope}%
\pgfsys@transformshift{1.722489in}{1.422680in}%
\pgfsys@useobject{currentmarker}{}%
\end{pgfscope}%
\begin{pgfscope}%
\pgfsys@transformshift{2.616418in}{1.436896in}%
\pgfsys@useobject{currentmarker}{}%
\end{pgfscope}%
\begin{pgfscope}%
\pgfsys@transformshift{4.106300in}{1.432275in}%
\pgfsys@useobject{currentmarker}{}%
\end{pgfscope}%
\begin{pgfscope}%
\pgfsys@transformshift{5.596183in}{1.402802in}%
\pgfsys@useobject{currentmarker}{}%
\end{pgfscope}%
\begin{pgfscope}%
\pgfsys@transformshift{7.086065in}{1.499909in}%
\pgfsys@useobject{currentmarker}{}%
\end{pgfscope}%
\end{pgfscope}%
\begin{pgfscope}%
\pgfsetrectcap%
\pgfsetmiterjoin%
\pgfsetlinewidth{0.803000pt}%
\definecolor{currentstroke}{rgb}{0.000000,0.000000,0.000000}%
\pgfsetstrokecolor{currentstroke}%
\pgfsetdash{}{0pt}%
\pgfpathmoveto{\pgfqpoint{1.126536in}{0.679475in}}%
\pgfpathlineto{\pgfqpoint{1.126536in}{3.135897in}}%
\pgfusepath{stroke}%
\end{pgfscope}%
\begin{pgfscope}%
\pgfsetrectcap%
\pgfsetmiterjoin%
\pgfsetlinewidth{0.803000pt}%
\definecolor{currentstroke}{rgb}{0.000000,0.000000,0.000000}%
\pgfsetstrokecolor{currentstroke}%
\pgfsetdash{}{0pt}%
\pgfpathmoveto{\pgfqpoint{7.682018in}{0.679475in}}%
\pgfpathlineto{\pgfqpoint{7.682018in}{3.135897in}}%
\pgfusepath{stroke}%
\end{pgfscope}%
\begin{pgfscope}%
\pgfsetrectcap%
\pgfsetmiterjoin%
\pgfsetlinewidth{0.803000pt}%
\definecolor{currentstroke}{rgb}{0.000000,0.000000,0.000000}%
\pgfsetstrokecolor{currentstroke}%
\pgfsetdash{}{0pt}%
\pgfpathmoveto{\pgfqpoint{1.126536in}{0.679475in}}%
\pgfpathlineto{\pgfqpoint{7.682018in}{0.679475in}}%
\pgfusepath{stroke}%
\end{pgfscope}%
\begin{pgfscope}%
\pgfsetrectcap%
\pgfsetmiterjoin%
\pgfsetlinewidth{0.803000pt}%
\definecolor{currentstroke}{rgb}{0.000000,0.000000,0.000000}%
\pgfsetstrokecolor{currentstroke}%
\pgfsetdash{}{0pt}%
\pgfpathmoveto{\pgfqpoint{1.126536in}{3.135897in}}%
\pgfpathlineto{\pgfqpoint{7.682018in}{3.135897in}}%
\pgfusepath{stroke}%
\end{pgfscope}%
\begin{pgfscope}%
\pgfsetrectcap%
\pgfsetroundjoin%
\pgfsetlinewidth{1.505625pt}%
\definecolor{currentstroke}{rgb}{0.337255,0.705882,0.913725}%
\pgfsetstrokecolor{currentstroke}%
\pgfsetdash{}{0pt}%
\pgfpathmoveto{\pgfqpoint{2.251710in}{3.363810in}}%
\pgfpathlineto{\pgfqpoint{2.501710in}{3.363810in}}%
\pgfpathlineto{\pgfqpoint{2.751710in}{3.363810in}}%
\pgfusepath{stroke}%
\end{pgfscope}%
\begin{pgfscope}%
\pgfsetbuttcap%
\pgfsetmiterjoin%
\definecolor{currentfill}{rgb}{0.337255,0.705882,0.913725}%
\pgfsetfillcolor{currentfill}%
\pgfsetlinewidth{1.003750pt}%
\definecolor{currentstroke}{rgb}{0.337255,0.705882,0.913725}%
\pgfsetstrokecolor{currentstroke}%
\pgfsetdash{}{0pt}%
\pgfsys@defobject{currentmarker}{\pgfqpoint{-0.041667in}{-0.041667in}}{\pgfqpoint{0.041667in}{0.041667in}}{%
\pgfpathmoveto{\pgfqpoint{-0.041667in}{-0.041667in}}%
\pgfpathlineto{\pgfqpoint{0.041667in}{-0.041667in}}%
\pgfpathlineto{\pgfqpoint{0.041667in}{0.041667in}}%
\pgfpathlineto{\pgfqpoint{-0.041667in}{0.041667in}}%
\pgfpathlineto{\pgfqpoint{-0.041667in}{-0.041667in}}%
\pgfpathclose%
\pgfusepath{stroke,fill}%
}%
\begin{pgfscope}%
\pgfsys@transformshift{2.501710in}{3.363810in}%
\pgfsys@useobject{currentmarker}{}%
\end{pgfscope}%
\end{pgfscope}%
\begin{pgfscope}%
\definecolor{textcolor}{rgb}{0.000000,0.000000,0.000000}%
\pgfsetstrokecolor{textcolor}%
\pgfsetfillcolor{textcolor}%
\pgftext[x=2.876710in,y=3.276310in,left,base]{\color{textcolor}\rmfamily\fontsize{18.000000}{21.600000}\selectfont LPP\(\displaystyle ^{R,T}\)}%
\end{pgfscope}%
\begin{pgfscope}%
\pgfsetbuttcap%
\pgfsetroundjoin%
\pgfsetlinewidth{1.505625pt}%
\definecolor{currentstroke}{rgb}{0.337255,0.705882,0.913725}%
\pgfsetstrokecolor{currentstroke}%
\pgfsetstrokeopacity{0.900000}%
\pgfsetdash{{5.550000pt}{2.400000pt}}{0.000000pt}%
\pgfpathmoveto{\pgfqpoint{3.812661in}{3.363810in}}%
\pgfpathlineto{\pgfqpoint{4.062661in}{3.363810in}}%
\pgfpathlineto{\pgfqpoint{4.312661in}{3.363810in}}%
\pgfusepath{stroke}%
\end{pgfscope}%
\begin{pgfscope}%
\pgfsetbuttcap%
\pgfsetmiterjoin%
\definecolor{currentfill}{rgb}{0.337255,0.705882,0.913725}%
\pgfsetfillcolor{currentfill}%
\pgfsetfillopacity{0.900000}%
\pgfsetlinewidth{1.003750pt}%
\definecolor{currentstroke}{rgb}{0.337255,0.705882,0.913725}%
\pgfsetstrokecolor{currentstroke}%
\pgfsetstrokeopacity{0.900000}%
\pgfsetdash{}{0pt}%
\pgfsys@defobject{currentmarker}{\pgfqpoint{-0.058926in}{-0.058926in}}{\pgfqpoint{0.058926in}{0.058926in}}{%
\pgfpathmoveto{\pgfqpoint{0.000000in}{-0.058926in}}%
\pgfpathlineto{\pgfqpoint{0.058926in}{0.000000in}}%
\pgfpathlineto{\pgfqpoint{0.000000in}{0.058926in}}%
\pgfpathlineto{\pgfqpoint{-0.058926in}{0.000000in}}%
\pgfpathlineto{\pgfqpoint{0.000000in}{-0.058926in}}%
\pgfpathclose%
\pgfusepath{stroke,fill}%
}%
\begin{pgfscope}%
\pgfsys@transformshift{4.062661in}{3.363810in}%
\pgfsys@useobject{currentmarker}{}%
\end{pgfscope}%
\end{pgfscope}%
\begin{pgfscope}%
\definecolor{textcolor}{rgb}{0.000000,0.000000,0.000000}%
\pgfsetstrokecolor{textcolor}%
\pgfsetfillcolor{textcolor}%
\pgftext[x=4.437661in,y=3.276310in,left,base]{\color{textcolor}\rmfamily\fontsize{18.000000}{21.600000}\selectfont LPP\(\displaystyle ^{T}\)}%
\end{pgfscope}%
\begin{pgfscope}%
\pgfsetbuttcap%
\pgfsetroundjoin%
\pgfsetlinewidth{1.505625pt}%
\definecolor{currentstroke}{rgb}{0.337255,0.705882,0.913725}%
\pgfsetstrokecolor{currentstroke}%
\pgfsetstrokeopacity{0.800000}%
\pgfsetdash{{1.500000pt}{2.475000pt}}{0.000000pt}%
\pgfpathmoveto{\pgfqpoint{5.202692in}{3.363810in}}%
\pgfpathlineto{\pgfqpoint{5.452692in}{3.363810in}}%
\pgfpathlineto{\pgfqpoint{5.702692in}{3.363810in}}%
\pgfusepath{stroke}%
\end{pgfscope}%
\begin{pgfscope}%
\pgfsetbuttcap%
\pgfsetroundjoin%
\definecolor{currentfill}{rgb}{0.337255,0.705882,0.913725}%
\pgfsetfillcolor{currentfill}%
\pgfsetfillopacity{0.800000}%
\pgfsetlinewidth{1.003750pt}%
\definecolor{currentstroke}{rgb}{0.337255,0.705882,0.913725}%
\pgfsetstrokecolor{currentstroke}%
\pgfsetstrokeopacity{0.800000}%
\pgfsetdash{}{0pt}%
\pgfsys@defobject{currentmarker}{\pgfqpoint{-0.041667in}{-0.041667in}}{\pgfqpoint{0.041667in}{0.041667in}}{%
\pgfpathmoveto{\pgfqpoint{0.000000in}{-0.041667in}}%
\pgfpathcurveto{\pgfqpoint{0.011050in}{-0.041667in}}{\pgfqpoint{0.021649in}{-0.037276in}}{\pgfqpoint{0.029463in}{-0.029463in}}%
\pgfpathcurveto{\pgfqpoint{0.037276in}{-0.021649in}}{\pgfqpoint{0.041667in}{-0.011050in}}{\pgfqpoint{0.041667in}{0.000000in}}%
\pgfpathcurveto{\pgfqpoint{0.041667in}{0.011050in}}{\pgfqpoint{0.037276in}{0.021649in}}{\pgfqpoint{0.029463in}{0.029463in}}%
\pgfpathcurveto{\pgfqpoint{0.021649in}{0.037276in}}{\pgfqpoint{0.011050in}{0.041667in}}{\pgfqpoint{0.000000in}{0.041667in}}%
\pgfpathcurveto{\pgfqpoint{-0.011050in}{0.041667in}}{\pgfqpoint{-0.021649in}{0.037276in}}{\pgfqpoint{-0.029463in}{0.029463in}}%
\pgfpathcurveto{\pgfqpoint{-0.037276in}{0.021649in}}{\pgfqpoint{-0.041667in}{0.011050in}}{\pgfqpoint{-0.041667in}{0.000000in}}%
\pgfpathcurveto{\pgfqpoint{-0.041667in}{-0.011050in}}{\pgfqpoint{-0.037276in}{-0.021649in}}{\pgfqpoint{-0.029463in}{-0.029463in}}%
\pgfpathcurveto{\pgfqpoint{-0.021649in}{-0.037276in}}{\pgfqpoint{-0.011050in}{-0.041667in}}{\pgfqpoint{0.000000in}{-0.041667in}}%
\pgfpathlineto{\pgfqpoint{0.000000in}{-0.041667in}}%
\pgfpathclose%
\pgfusepath{stroke,fill}%
}%
\begin{pgfscope}%
\pgfsys@transformshift{5.452692in}{3.363810in}%
\pgfsys@useobject{currentmarker}{}%
\end{pgfscope}%
\end{pgfscope}%
\begin{pgfscope}%
\definecolor{textcolor}{rgb}{0.000000,0.000000,0.000000}%
\pgfsetstrokecolor{textcolor}%
\pgfsetfillcolor{textcolor}%
\pgftext[x=5.827692in,y=3.276310in,left,base]{\color{textcolor}\rmfamily\fontsize{18.000000}{21.600000}\selectfont LPP\(\displaystyle ^{T, P}\)}%
\end{pgfscope}%
\end{pgfpicture}%
\makeatother%
\endgroup%

%% file: figures/DemandPlots/TightLayout/DE_PI_rescaled_Gap10pct.pgf
\begingroup%
\makeatletter%
\begin{pgfpicture}%
\pgfpathrectangle{\pgfpointorigin}{\pgfqpoint{7.782018in}{3.651310in}}%
\pgfusepath{use as bounding box, clip}%
\begin{pgfscope}%
\pgfsetbuttcap%
\pgfsetmiterjoin%
\definecolor{currentfill}{rgb}{1.000000,1.000000,1.000000}%
\pgfsetfillcolor{currentfill}%
\pgfsetlinewidth{0.000000pt}%
\definecolor{currentstroke}{rgb}{1.000000,1.000000,1.000000}%
\pgfsetstrokecolor{currentstroke}%
\pgfsetdash{}{0pt}%
\pgfpathmoveto{\pgfqpoint{0.000000in}{0.000000in}}%
\pgfpathlineto{\pgfqpoint{7.782018in}{0.000000in}}%
\pgfpathlineto{\pgfqpoint{7.782018in}{3.651310in}}%
\pgfpathlineto{\pgfqpoint{0.000000in}{3.651310in}}%
\pgfpathlineto{\pgfqpoint{0.000000in}{0.000000in}}%
\pgfpathclose%
\pgfusepath{fill}%
\end{pgfscope}%
\begin{pgfscope}%
\pgfsetbuttcap%
\pgfsetmiterjoin%
\definecolor{currentfill}{rgb}{1.000000,1.000000,1.000000}%
\pgfsetfillcolor{currentfill}%
\pgfsetlinewidth{0.000000pt}%
\definecolor{currentstroke}{rgb}{0.000000,0.000000,0.000000}%
\pgfsetstrokecolor{currentstroke}%
\pgfsetstrokeopacity{0.000000}%
\pgfsetdash{}{0pt}%
\pgfpathmoveto{\pgfqpoint{1.126536in}{0.679475in}}%
\pgfpathlineto{\pgfqpoint{7.682018in}{0.679475in}}%
\pgfpathlineto{\pgfqpoint{7.682018in}{3.135897in}}%
\pgfpathlineto{\pgfqpoint{1.126536in}{3.135897in}}%
\pgfpathlineto{\pgfqpoint{1.126536in}{0.679475in}}%
\pgfpathclose%
\pgfusepath{fill}%
\end{pgfscope}%
\begin{pgfscope}%
\pgfsetbuttcap%
\pgfsetroundjoin%
\definecolor{currentfill}{rgb}{0.000000,0.000000,0.000000}%
\pgfsetfillcolor{currentfill}%
\pgfsetlinewidth{0.803000pt}%
\definecolor{currentstroke}{rgb}{0.000000,0.000000,0.000000}%
\pgfsetstrokecolor{currentstroke}%
\pgfsetdash{}{0pt}%
\pgfsys@defobject{currentmarker}{\pgfqpoint{0.000000in}{-0.048611in}}{\pgfqpoint{0.000000in}{0.000000in}}{%
\pgfpathmoveto{\pgfqpoint{0.000000in}{0.000000in}}%
\pgfpathlineto{\pgfqpoint{0.000000in}{-0.048611in}}%
\pgfusepath{stroke,fill}%
}%
\begin{pgfscope}%
\pgfsys@transformshift{1.126536in}{0.679475in}%
\pgfsys@useobject{currentmarker}{}%
\end{pgfscope}%
\end{pgfscope}%
\begin{pgfscope}%
\definecolor{textcolor}{rgb}{0.000000,0.000000,0.000000}%
\pgfsetstrokecolor{textcolor}%
\pgfsetfillcolor{textcolor}%
\pgftext[x=1.126536in,y=0.582253in,,top]{\color{textcolor}\rmfamily\fontsize{16.000000}{19.200000}\selectfont \(\displaystyle {0.0}\)}%
\end{pgfscope}%
\begin{pgfscope}%
\pgfsetbuttcap%
\pgfsetroundjoin%
\definecolor{currentfill}{rgb}{0.000000,0.000000,0.000000}%
\pgfsetfillcolor{currentfill}%
\pgfsetlinewidth{0.803000pt}%
\definecolor{currentstroke}{rgb}{0.000000,0.000000,0.000000}%
\pgfsetstrokecolor{currentstroke}%
\pgfsetdash{}{0pt}%
\pgfsys@defobject{currentmarker}{\pgfqpoint{0.000000in}{-0.048611in}}{\pgfqpoint{0.000000in}{0.000000in}}{%
\pgfpathmoveto{\pgfqpoint{0.000000in}{0.000000in}}%
\pgfpathlineto{\pgfqpoint{0.000000in}{-0.048611in}}%
\pgfusepath{stroke,fill}%
}%
\begin{pgfscope}%
\pgfsys@transformshift{2.318442in}{0.679475in}%
\pgfsys@useobject{currentmarker}{}%
\end{pgfscope}%
\end{pgfscope}%
\begin{pgfscope}%
\definecolor{textcolor}{rgb}{0.000000,0.000000,0.000000}%
\pgfsetstrokecolor{textcolor}%
\pgfsetfillcolor{textcolor}%
\pgftext[x=2.318442in,y=0.582253in,,top]{\color{textcolor}\rmfamily\fontsize{16.000000}{19.200000}\selectfont \(\displaystyle {0.2}\)}%
\end{pgfscope}%
\begin{pgfscope}%
\pgfsetbuttcap%
\pgfsetroundjoin%
\definecolor{currentfill}{rgb}{0.000000,0.000000,0.000000}%
\pgfsetfillcolor{currentfill}%
\pgfsetlinewidth{0.803000pt}%
\definecolor{currentstroke}{rgb}{0.000000,0.000000,0.000000}%
\pgfsetstrokecolor{currentstroke}%
\pgfsetdash{}{0pt}%
\pgfsys@defobject{currentmarker}{\pgfqpoint{0.000000in}{-0.048611in}}{\pgfqpoint{0.000000in}{0.000000in}}{%
\pgfpathmoveto{\pgfqpoint{0.000000in}{0.000000in}}%
\pgfpathlineto{\pgfqpoint{0.000000in}{-0.048611in}}%
\pgfusepath{stroke,fill}%
}%
\begin{pgfscope}%
\pgfsys@transformshift{3.510347in}{0.679475in}%
\pgfsys@useobject{currentmarker}{}%
\end{pgfscope}%
\end{pgfscope}%
\begin{pgfscope}%
\definecolor{textcolor}{rgb}{0.000000,0.000000,0.000000}%
\pgfsetstrokecolor{textcolor}%
\pgfsetfillcolor{textcolor}%
\pgftext[x=3.510347in,y=0.582253in,,top]{\color{textcolor}\rmfamily\fontsize{16.000000}{19.200000}\selectfont \(\displaystyle {0.4}\)}%
\end{pgfscope}%
\begin{pgfscope}%
\pgfsetbuttcap%
\pgfsetroundjoin%
\definecolor{currentfill}{rgb}{0.000000,0.000000,0.000000}%
\pgfsetfillcolor{currentfill}%
\pgfsetlinewidth{0.803000pt}%
\definecolor{currentstroke}{rgb}{0.000000,0.000000,0.000000}%
\pgfsetstrokecolor{currentstroke}%
\pgfsetdash{}{0pt}%
\pgfsys@defobject{currentmarker}{\pgfqpoint{0.000000in}{-0.048611in}}{\pgfqpoint{0.000000in}{0.000000in}}{%
\pgfpathmoveto{\pgfqpoint{0.000000in}{0.000000in}}%
\pgfpathlineto{\pgfqpoint{0.000000in}{-0.048611in}}%
\pgfusepath{stroke,fill}%
}%
\begin{pgfscope}%
\pgfsys@transformshift{4.702253in}{0.679475in}%
\pgfsys@useobject{currentmarker}{}%
\end{pgfscope}%
\end{pgfscope}%
\begin{pgfscope}%
\definecolor{textcolor}{rgb}{0.000000,0.000000,0.000000}%
\pgfsetstrokecolor{textcolor}%
\pgfsetfillcolor{textcolor}%
\pgftext[x=4.702253in,y=0.582253in,,top]{\color{textcolor}\rmfamily\fontsize{16.000000}{19.200000}\selectfont \(\displaystyle {0.6}\)}%
\end{pgfscope}%
\begin{pgfscope}%
\pgfsetbuttcap%
\pgfsetroundjoin%
\definecolor{currentfill}{rgb}{0.000000,0.000000,0.000000}%
\pgfsetfillcolor{currentfill}%
\pgfsetlinewidth{0.803000pt}%
\definecolor{currentstroke}{rgb}{0.000000,0.000000,0.000000}%
\pgfsetstrokecolor{currentstroke}%
\pgfsetdash{}{0pt}%
\pgfsys@defobject{currentmarker}{\pgfqpoint{0.000000in}{-0.048611in}}{\pgfqpoint{0.000000in}{0.000000in}}{%
\pgfpathmoveto{\pgfqpoint{0.000000in}{0.000000in}}%
\pgfpathlineto{\pgfqpoint{0.000000in}{-0.048611in}}%
\pgfusepath{stroke,fill}%
}%
\begin{pgfscope}%
\pgfsys@transformshift{5.894159in}{0.679475in}%
\pgfsys@useobject{currentmarker}{}%
\end{pgfscope}%
\end{pgfscope}%
\begin{pgfscope}%
\definecolor{textcolor}{rgb}{0.000000,0.000000,0.000000}%
\pgfsetstrokecolor{textcolor}%
\pgfsetfillcolor{textcolor}%
\pgftext[x=5.894159in,y=0.582253in,,top]{\color{textcolor}\rmfamily\fontsize{16.000000}{19.200000}\selectfont \(\displaystyle {0.8}\)}%
\end{pgfscope}%
\begin{pgfscope}%
\pgfsetbuttcap%
\pgfsetroundjoin%
\definecolor{currentfill}{rgb}{0.000000,0.000000,0.000000}%
\pgfsetfillcolor{currentfill}%
\pgfsetlinewidth{0.803000pt}%
\definecolor{currentstroke}{rgb}{0.000000,0.000000,0.000000}%
\pgfsetstrokecolor{currentstroke}%
\pgfsetdash{}{0pt}%
\pgfsys@defobject{currentmarker}{\pgfqpoint{0.000000in}{-0.048611in}}{\pgfqpoint{0.000000in}{0.000000in}}{%
\pgfpathmoveto{\pgfqpoint{0.000000in}{0.000000in}}%
\pgfpathlineto{\pgfqpoint{0.000000in}{-0.048611in}}%
\pgfusepath{stroke,fill}%
}%
\begin{pgfscope}%
\pgfsys@transformshift{7.086065in}{0.679475in}%
\pgfsys@useobject{currentmarker}{}%
\end{pgfscope}%
\end{pgfscope}%
\begin{pgfscope}%
\definecolor{textcolor}{rgb}{0.000000,0.000000,0.000000}%
\pgfsetstrokecolor{textcolor}%
\pgfsetfillcolor{textcolor}%
\pgftext[x=7.086065in,y=0.582253in,,top]{\color{textcolor}\rmfamily\fontsize{16.000000}{19.200000}\selectfont \(\displaystyle {1.0}\)}%
\end{pgfscope}%
\begin{pgfscope}%
\definecolor{textcolor}{rgb}{0.000000,0.000000,0.000000}%
\pgfsetstrokecolor{textcolor}%
\pgfsetfillcolor{textcolor}%
\pgftext[x=4.404277in,y=0.313349in,,top]{\color{textcolor}\rmfamily\fontsize{18.000000}{21.600000}\selectfont \(\displaystyle \lambda\)}%
\end{pgfscope}%
\begin{pgfscope}%
\pgfsetbuttcap%
\pgfsetroundjoin%
\definecolor{currentfill}{rgb}{0.000000,0.000000,0.000000}%
\pgfsetfillcolor{currentfill}%
\pgfsetlinewidth{0.803000pt}%
\definecolor{currentstroke}{rgb}{0.000000,0.000000,0.000000}%
\pgfsetstrokecolor{currentstroke}%
\pgfsetdash{}{0pt}%
\pgfsys@defobject{currentmarker}{\pgfqpoint{-0.048611in}{0.000000in}}{\pgfqpoint{-0.000000in}{0.000000in}}{%
\pgfpathmoveto{\pgfqpoint{-0.000000in}{0.000000in}}%
\pgfpathlineto{\pgfqpoint{-0.048611in}{0.000000in}}%
\pgfusepath{stroke,fill}%
}%
\begin{pgfscope}%
\pgfsys@transformshift{1.126536in}{1.047889in}%
\pgfsys@useobject{currentmarker}{}%
\end{pgfscope}%
\end{pgfscope}%
\begin{pgfscope}%
\definecolor{textcolor}{rgb}{0.000000,0.000000,0.000000}%
\pgfsetstrokecolor{textcolor}%
\pgfsetfillcolor{textcolor}%
\pgftext[x=0.478972in, y=0.964555in, left, base]{\color{textcolor}\rmfamily\fontsize{16.000000}{19.200000}\selectfont \(\displaystyle {95000}\)}%
\end{pgfscope}%
\begin{pgfscope}%
\pgfsetbuttcap%
\pgfsetroundjoin%
\definecolor{currentfill}{rgb}{0.000000,0.000000,0.000000}%
\pgfsetfillcolor{currentfill}%
\pgfsetlinewidth{0.803000pt}%
\definecolor{currentstroke}{rgb}{0.000000,0.000000,0.000000}%
\pgfsetstrokecolor{currentstroke}%
\pgfsetdash{}{0pt}%
\pgfsys@defobject{currentmarker}{\pgfqpoint{-0.048611in}{0.000000in}}{\pgfqpoint{-0.000000in}{0.000000in}}{%
\pgfpathmoveto{\pgfqpoint{-0.000000in}{0.000000in}}%
\pgfpathlineto{\pgfqpoint{-0.048611in}{0.000000in}}%
\pgfusepath{stroke,fill}%
}%
\begin{pgfscope}%
\pgfsys@transformshift{1.126536in}{1.473949in}%
\pgfsys@useobject{currentmarker}{}%
\end{pgfscope}%
\end{pgfscope}%
\begin{pgfscope}%
\definecolor{textcolor}{rgb}{0.000000,0.000000,0.000000}%
\pgfsetstrokecolor{textcolor}%
\pgfsetfillcolor{textcolor}%
\pgftext[x=0.368904in, y=1.390616in, left, base]{\color{textcolor}\rmfamily\fontsize{16.000000}{19.200000}\selectfont \(\displaystyle {100000}\)}%
\end{pgfscope}%
\begin{pgfscope}%
\pgfsetbuttcap%
\pgfsetroundjoin%
\definecolor{currentfill}{rgb}{0.000000,0.000000,0.000000}%
\pgfsetfillcolor{currentfill}%
\pgfsetlinewidth{0.803000pt}%
\definecolor{currentstroke}{rgb}{0.000000,0.000000,0.000000}%
\pgfsetstrokecolor{currentstroke}%
\pgfsetdash{}{0pt}%
\pgfsys@defobject{currentmarker}{\pgfqpoint{-0.048611in}{0.000000in}}{\pgfqpoint{-0.000000in}{0.000000in}}{%
\pgfpathmoveto{\pgfqpoint{-0.000000in}{0.000000in}}%
\pgfpathlineto{\pgfqpoint{-0.048611in}{0.000000in}}%
\pgfusepath{stroke,fill}%
}%
\begin{pgfscope}%
\pgfsys@transformshift{1.126536in}{1.900010in}%
\pgfsys@useobject{currentmarker}{}%
\end{pgfscope}%
\end{pgfscope}%
\begin{pgfscope}%
\definecolor{textcolor}{rgb}{0.000000,0.000000,0.000000}%
\pgfsetstrokecolor{textcolor}%
\pgfsetfillcolor{textcolor}%
\pgftext[x=0.368904in, y=1.816676in, left, base]{\color{textcolor}\rmfamily\fontsize{16.000000}{19.200000}\selectfont \(\displaystyle {105000}\)}%
\end{pgfscope}%
\begin{pgfscope}%
\pgfsetbuttcap%
\pgfsetroundjoin%
\definecolor{currentfill}{rgb}{0.000000,0.000000,0.000000}%
\pgfsetfillcolor{currentfill}%
\pgfsetlinewidth{0.803000pt}%
\definecolor{currentstroke}{rgb}{0.000000,0.000000,0.000000}%
\pgfsetstrokecolor{currentstroke}%
\pgfsetdash{}{0pt}%
\pgfsys@defobject{currentmarker}{\pgfqpoint{-0.048611in}{0.000000in}}{\pgfqpoint{-0.000000in}{0.000000in}}{%
\pgfpathmoveto{\pgfqpoint{-0.000000in}{0.000000in}}%
\pgfpathlineto{\pgfqpoint{-0.048611in}{0.000000in}}%
\pgfusepath{stroke,fill}%
}%
\begin{pgfscope}%
\pgfsys@transformshift{1.126536in}{2.326070in}%
\pgfsys@useobject{currentmarker}{}%
\end{pgfscope}%
\end{pgfscope}%
\begin{pgfscope}%
\definecolor{textcolor}{rgb}{0.000000,0.000000,0.000000}%
\pgfsetstrokecolor{textcolor}%
\pgfsetfillcolor{textcolor}%
\pgftext[x=0.368904in, y=2.242737in, left, base]{\color{textcolor}\rmfamily\fontsize{16.000000}{19.200000}\selectfont \(\displaystyle {110000}\)}%
\end{pgfscope}%
\begin{pgfscope}%
\pgfsetbuttcap%
\pgfsetroundjoin%
\definecolor{currentfill}{rgb}{0.000000,0.000000,0.000000}%
\pgfsetfillcolor{currentfill}%
\pgfsetlinewidth{0.803000pt}%
\definecolor{currentstroke}{rgb}{0.000000,0.000000,0.000000}%
\pgfsetstrokecolor{currentstroke}%
\pgfsetdash{}{0pt}%
\pgfsys@defobject{currentmarker}{\pgfqpoint{-0.048611in}{0.000000in}}{\pgfqpoint{-0.000000in}{0.000000in}}{%
\pgfpathmoveto{\pgfqpoint{-0.000000in}{0.000000in}}%
\pgfpathlineto{\pgfqpoint{-0.048611in}{0.000000in}}%
\pgfusepath{stroke,fill}%
}%
\begin{pgfscope}%
\pgfsys@transformshift{1.126536in}{2.752131in}%
\pgfsys@useobject{currentmarker}{}%
\end{pgfscope}%
\end{pgfscope}%
\begin{pgfscope}%
\definecolor{textcolor}{rgb}{0.000000,0.000000,0.000000}%
\pgfsetstrokecolor{textcolor}%
\pgfsetfillcolor{textcolor}%
\pgftext[x=0.368904in, y=2.668797in, left, base]{\color{textcolor}\rmfamily\fontsize{16.000000}{19.200000}\selectfont \(\displaystyle {115000}\)}%
\end{pgfscope}%
\begin{pgfscope}%
\definecolor{textcolor}{rgb}{0.000000,0.000000,0.000000}%
\pgfsetstrokecolor{textcolor}%
\pgfsetfillcolor{textcolor}%
\pgftext[x=0.313349in,y=1.907686in,,bottom,rotate=90.000000]{\color{textcolor}\rmfamily\fontsize{18.000000}{21.600000}\selectfont DKK}%
\end{pgfscope}%
\begin{pgfscope}%
\pgfpathrectangle{\pgfqpoint{1.126536in}{0.679475in}}{\pgfqpoint{6.555482in}{2.456422in}}%
\pgfusepath{clip}%
\pgfsetrectcap%
\pgfsetroundjoin%
\pgfsetlinewidth{1.505625pt}%
\definecolor{currentstroke}{rgb}{0.580392,0.403922,0.741176}%
\pgfsetstrokecolor{currentstroke}%
\pgfsetdash{}{0pt}%
\pgfpathmoveto{\pgfqpoint{1.722489in}{2.652504in}}%
\pgfpathlineto{\pgfqpoint{2.616418in}{2.487064in}}%
\pgfpathlineto{\pgfqpoint{4.106300in}{2.179378in}}%
\pgfpathlineto{\pgfqpoint{5.596183in}{1.777039in}}%
\pgfpathlineto{\pgfqpoint{7.086065in}{1.337312in}}%
\pgfusepath{stroke}%
\end{pgfscope}%
\begin{pgfscope}%
\pgfpathrectangle{\pgfqpoint{1.126536in}{0.679475in}}{\pgfqpoint{6.555482in}{2.456422in}}%
\pgfusepath{clip}%
\pgfsetbuttcap%
\pgfsetmiterjoin%
\definecolor{currentfill}{rgb}{0.580392,0.403922,0.741176}%
\pgfsetfillcolor{currentfill}%
\pgfsetlinewidth{0.752812pt}%
\definecolor{currentstroke}{rgb}{1.000000,1.000000,1.000000}%
\pgfsetstrokecolor{currentstroke}%
\pgfsetdash{}{0pt}%
\pgfsys@defobject{currentmarker}{\pgfqpoint{-0.041667in}{-0.041667in}}{\pgfqpoint{0.041667in}{0.041667in}}{%
\pgfpathmoveto{\pgfqpoint{-0.041667in}{-0.041667in}}%
\pgfpathlineto{\pgfqpoint{0.041667in}{-0.041667in}}%
\pgfpathlineto{\pgfqpoint{0.041667in}{0.041667in}}%
\pgfpathlineto{\pgfqpoint{-0.041667in}{0.041667in}}%
\pgfpathlineto{\pgfqpoint{-0.041667in}{-0.041667in}}%
\pgfpathclose%
\pgfusepath{stroke,fill}%
}%
\begin{pgfscope}%
\pgfsys@transformshift{1.722489in}{2.652504in}%
\pgfsys@useobject{currentmarker}{}%
\end{pgfscope}%
\begin{pgfscope}%
\pgfsys@transformshift{2.616418in}{2.487064in}%
\pgfsys@useobject{currentmarker}{}%
\end{pgfscope}%
\begin{pgfscope}%
\pgfsys@transformshift{4.106300in}{2.179378in}%
\pgfsys@useobject{currentmarker}{}%
\end{pgfscope}%
\begin{pgfscope}%
\pgfsys@transformshift{5.596183in}{1.777039in}%
\pgfsys@useobject{currentmarker}{}%
\end{pgfscope}%
\begin{pgfscope}%
\pgfsys@transformshift{7.086065in}{1.337312in}%
\pgfsys@useobject{currentmarker}{}%
\end{pgfscope}%
\end{pgfscope}%
\begin{pgfscope}%
\pgfpathrectangle{\pgfqpoint{1.126536in}{0.679475in}}{\pgfqpoint{6.555482in}{2.456422in}}%
\pgfusepath{clip}%
\pgfsetbuttcap%
\pgfsetroundjoin%
\pgfsetlinewidth{1.505625pt}%
\definecolor{currentstroke}{rgb}{0.580392,0.403922,0.741176}%
\pgfsetstrokecolor{currentstroke}%
\pgfsetstrokeopacity{0.900000}%
\pgfsetdash{{7.500000pt}{7.500000pt}}{0.000000pt}%
\pgfpathmoveto{\pgfqpoint{1.722489in}{3.024241in}}%
\pgfpathlineto{\pgfqpoint{2.616418in}{2.680297in}}%
\pgfpathlineto{\pgfqpoint{4.106300in}{2.403769in}}%
\pgfpathlineto{\pgfqpoint{5.596183in}{2.023927in}}%
\pgfpathlineto{\pgfqpoint{7.086065in}{1.617596in}}%
\pgfusepath{stroke}%
\end{pgfscope}%
\begin{pgfscope}%
\pgfpathrectangle{\pgfqpoint{1.126536in}{0.679475in}}{\pgfqpoint{6.555482in}{2.456422in}}%
\pgfusepath{clip}%
\pgfsetbuttcap%
\pgfsetmiterjoin%
\definecolor{currentfill}{rgb}{0.580392,0.403922,0.741176}%
\pgfsetfillcolor{currentfill}%
\pgfsetfillopacity{0.900000}%
\pgfsetlinewidth{0.752812pt}%
\definecolor{currentstroke}{rgb}{1.000000,1.000000,1.000000}%
\pgfsetstrokecolor{currentstroke}%
\pgfsetdash{}{0pt}%
\pgfsys@defobject{currentmarker}{\pgfqpoint{-0.058926in}{-0.058926in}}{\pgfqpoint{0.058926in}{0.058926in}}{%
\pgfpathmoveto{\pgfqpoint{0.000000in}{-0.058926in}}%
\pgfpathlineto{\pgfqpoint{0.058926in}{0.000000in}}%
\pgfpathlineto{\pgfqpoint{0.000000in}{0.058926in}}%
\pgfpathlineto{\pgfqpoint{-0.058926in}{0.000000in}}%
\pgfpathlineto{\pgfqpoint{0.000000in}{-0.058926in}}%
\pgfpathclose%
\pgfusepath{stroke,fill}%
}%
\begin{pgfscope}%
\pgfsys@transformshift{1.722489in}{3.024241in}%
\pgfsys@useobject{currentmarker}{}%
\end{pgfscope}%
\begin{pgfscope}%
\pgfsys@transformshift{2.616418in}{2.680297in}%
\pgfsys@useobject{currentmarker}{}%
\end{pgfscope}%
\begin{pgfscope}%
\pgfsys@transformshift{4.106300in}{2.403769in}%
\pgfsys@useobject{currentmarker}{}%
\end{pgfscope}%
\begin{pgfscope}%
\pgfsys@transformshift{5.596183in}{2.023927in}%
\pgfsys@useobject{currentmarker}{}%
\end{pgfscope}%
\begin{pgfscope}%
\pgfsys@transformshift{7.086065in}{1.617596in}%
\pgfsys@useobject{currentmarker}{}%
\end{pgfscope}%
\end{pgfscope}%
\begin{pgfscope}%
\pgfpathrectangle{\pgfqpoint{1.126536in}{0.679475in}}{\pgfqpoint{6.555482in}{2.456422in}}%
\pgfusepath{clip}%
\pgfsetbuttcap%
\pgfsetroundjoin%
\pgfsetlinewidth{1.505625pt}%
\definecolor{currentstroke}{rgb}{0.580392,0.403922,0.741176}%
\pgfsetstrokecolor{currentstroke}%
\pgfsetstrokeopacity{0.800000}%
\pgfsetdash{{1.500000pt}{4.500000pt}}{0.000000pt}%
\pgfpathmoveto{\pgfqpoint{1.722489in}{1.389660in}}%
\pgfpathlineto{\pgfqpoint{2.616418in}{1.379477in}}%
\pgfpathlineto{\pgfqpoint{4.106300in}{1.234901in}}%
\pgfpathlineto{\pgfqpoint{5.596183in}{1.017648in}}%
\pgfpathlineto{\pgfqpoint{7.086065in}{0.791131in}}%
\pgfusepath{stroke}%
\end{pgfscope}%
\begin{pgfscope}%
\pgfpathrectangle{\pgfqpoint{1.126536in}{0.679475in}}{\pgfqpoint{6.555482in}{2.456422in}}%
\pgfusepath{clip}%
\pgfsetbuttcap%
\pgfsetroundjoin%
\definecolor{currentfill}{rgb}{0.580392,0.403922,0.741176}%
\pgfsetfillcolor{currentfill}%
\pgfsetfillopacity{0.800000}%
\pgfsetlinewidth{0.752812pt}%
\definecolor{currentstroke}{rgb}{1.000000,1.000000,1.000000}%
\pgfsetstrokecolor{currentstroke}%
\pgfsetdash{}{0pt}%
\pgfsys@defobject{currentmarker}{\pgfqpoint{-0.041667in}{-0.041667in}}{\pgfqpoint{0.041667in}{0.041667in}}{%
\pgfpathmoveto{\pgfqpoint{0.000000in}{-0.041667in}}%
\pgfpathcurveto{\pgfqpoint{0.011050in}{-0.041667in}}{\pgfqpoint{0.021649in}{-0.037276in}}{\pgfqpoint{0.029463in}{-0.029463in}}%
\pgfpathcurveto{\pgfqpoint{0.037276in}{-0.021649in}}{\pgfqpoint{0.041667in}{-0.011050in}}{\pgfqpoint{0.041667in}{0.000000in}}%
\pgfpathcurveto{\pgfqpoint{0.041667in}{0.011050in}}{\pgfqpoint{0.037276in}{0.021649in}}{\pgfqpoint{0.029463in}{0.029463in}}%
\pgfpathcurveto{\pgfqpoint{0.021649in}{0.037276in}}{\pgfqpoint{0.011050in}{0.041667in}}{\pgfqpoint{0.000000in}{0.041667in}}%
\pgfpathcurveto{\pgfqpoint{-0.011050in}{0.041667in}}{\pgfqpoint{-0.021649in}{0.037276in}}{\pgfqpoint{-0.029463in}{0.029463in}}%
\pgfpathcurveto{\pgfqpoint{-0.037276in}{0.021649in}}{\pgfqpoint{-0.041667in}{0.011050in}}{\pgfqpoint{-0.041667in}{0.000000in}}%
\pgfpathcurveto{\pgfqpoint{-0.041667in}{-0.011050in}}{\pgfqpoint{-0.037276in}{-0.021649in}}{\pgfqpoint{-0.029463in}{-0.029463in}}%
\pgfpathcurveto{\pgfqpoint{-0.021649in}{-0.037276in}}{\pgfqpoint{-0.011050in}{-0.041667in}}{\pgfqpoint{0.000000in}{-0.041667in}}%
\pgfpathlineto{\pgfqpoint{0.000000in}{-0.041667in}}%
\pgfpathclose%
\pgfusepath{stroke,fill}%
}%
\begin{pgfscope}%
\pgfsys@transformshift{1.722489in}{1.389660in}%
\pgfsys@useobject{currentmarker}{}%
\end{pgfscope}%
\begin{pgfscope}%
\pgfsys@transformshift{2.616418in}{1.379477in}%
\pgfsys@useobject{currentmarker}{}%
\end{pgfscope}%
\begin{pgfscope}%
\pgfsys@transformshift{4.106300in}{1.234901in}%
\pgfsys@useobject{currentmarker}{}%
\end{pgfscope}%
\begin{pgfscope}%
\pgfsys@transformshift{5.596183in}{1.017648in}%
\pgfsys@useobject{currentmarker}{}%
\end{pgfscope}%
\begin{pgfscope}%
\pgfsys@transformshift{7.086065in}{0.791131in}%
\pgfsys@useobject{currentmarker}{}%
\end{pgfscope}%
\end{pgfscope}%
\begin{pgfscope}%
\pgfsetrectcap%
\pgfsetmiterjoin%
\pgfsetlinewidth{0.803000pt}%
\definecolor{currentstroke}{rgb}{0.000000,0.000000,0.000000}%
\pgfsetstrokecolor{currentstroke}%
\pgfsetdash{}{0pt}%
\pgfpathmoveto{\pgfqpoint{1.126536in}{0.679475in}}%
\pgfpathlineto{\pgfqpoint{1.126536in}{3.135897in}}%
\pgfusepath{stroke}%
\end{pgfscope}%
\begin{pgfscope}%
\pgfsetrectcap%
\pgfsetmiterjoin%
\pgfsetlinewidth{0.803000pt}%
\definecolor{currentstroke}{rgb}{0.000000,0.000000,0.000000}%
\pgfsetstrokecolor{currentstroke}%
\pgfsetdash{}{0pt}%
\pgfpathmoveto{\pgfqpoint{7.682018in}{0.679475in}}%
\pgfpathlineto{\pgfqpoint{7.682018in}{3.135897in}}%
\pgfusepath{stroke}%
\end{pgfscope}%
\begin{pgfscope}%
\pgfsetrectcap%
\pgfsetmiterjoin%
\pgfsetlinewidth{0.803000pt}%
\definecolor{currentstroke}{rgb}{0.000000,0.000000,0.000000}%
\pgfsetstrokecolor{currentstroke}%
\pgfsetdash{}{0pt}%
\pgfpathmoveto{\pgfqpoint{1.126536in}{0.679475in}}%
\pgfpathlineto{\pgfqpoint{7.682018in}{0.679475in}}%
\pgfusepath{stroke}%
\end{pgfscope}%
\begin{pgfscope}%
\pgfsetrectcap%
\pgfsetmiterjoin%
\pgfsetlinewidth{0.803000pt}%
\definecolor{currentstroke}{rgb}{0.000000,0.000000,0.000000}%
\pgfsetstrokecolor{currentstroke}%
\pgfsetdash{}{0pt}%
\pgfpathmoveto{\pgfqpoint{1.126536in}{3.135897in}}%
\pgfpathlineto{\pgfqpoint{7.682018in}{3.135897in}}%
\pgfusepath{stroke}%
\end{pgfscope}%
\begin{pgfscope}%
\pgfsetrectcap%
\pgfsetroundjoin%
\pgfsetlinewidth{1.505625pt}%
\definecolor{currentstroke}{rgb}{0.580392,0.403922,0.741176}%
\pgfsetstrokecolor{currentstroke}%
\pgfsetdash{}{0pt}%
\pgfpathmoveto{\pgfqpoint{2.251710in}{3.363810in}}%
\pgfpathlineto{\pgfqpoint{2.501710in}{3.363810in}}%
\pgfpathlineto{\pgfqpoint{2.751710in}{3.363810in}}%
\pgfusepath{stroke}%
\end{pgfscope}%
\begin{pgfscope}%
\pgfsetbuttcap%
\pgfsetmiterjoin%
\definecolor{currentfill}{rgb}{0.580392,0.403922,0.741176}%
\pgfsetfillcolor{currentfill}%
\pgfsetlinewidth{1.003750pt}%
\definecolor{currentstroke}{rgb}{0.580392,0.403922,0.741176}%
\pgfsetstrokecolor{currentstroke}%
\pgfsetdash{}{0pt}%
\pgfsys@defobject{currentmarker}{\pgfqpoint{-0.041667in}{-0.041667in}}{\pgfqpoint{0.041667in}{0.041667in}}{%
\pgfpathmoveto{\pgfqpoint{-0.041667in}{-0.041667in}}%
\pgfpathlineto{\pgfqpoint{0.041667in}{-0.041667in}}%
\pgfpathlineto{\pgfqpoint{0.041667in}{0.041667in}}%
\pgfpathlineto{\pgfqpoint{-0.041667in}{0.041667in}}%
\pgfpathlineto{\pgfqpoint{-0.041667in}{-0.041667in}}%
\pgfpathclose%
\pgfusepath{stroke,fill}%
}%
\begin{pgfscope}%
\pgfsys@transformshift{2.501710in}{3.363810in}%
\pgfsys@useobject{currentmarker}{}%
\end{pgfscope}%
\end{pgfscope}%
\begin{pgfscope}%
\definecolor{textcolor}{rgb}{0.000000,0.000000,0.000000}%
\pgfsetstrokecolor{textcolor}%
\pgfsetfillcolor{textcolor}%
\pgftext[x=2.876710in,y=3.276310in,left,base]{\color{textcolor}\rmfamily\fontsize{18.000000}{21.600000}\selectfont LPP\(\displaystyle ^{R,T}\)}%
\end{pgfscope}%
\begin{pgfscope}%
\pgfsetbuttcap%
\pgfsetroundjoin%
\pgfsetlinewidth{1.505625pt}%
\definecolor{currentstroke}{rgb}{0.580392,0.403922,0.741176}%
\pgfsetstrokecolor{currentstroke}%
\pgfsetstrokeopacity{0.900000}%
\pgfsetdash{{5.550000pt}{2.400000pt}}{0.000000pt}%
\pgfpathmoveto{\pgfqpoint{3.812661in}{3.363810in}}%
\pgfpathlineto{\pgfqpoint{4.062661in}{3.363810in}}%
\pgfpathlineto{\pgfqpoint{4.312661in}{3.363810in}}%
\pgfusepath{stroke}%
\end{pgfscope}%
\begin{pgfscope}%
\pgfsetbuttcap%
\pgfsetmiterjoin%
\definecolor{currentfill}{rgb}{0.580392,0.403922,0.741176}%
\pgfsetfillcolor{currentfill}%
\pgfsetfillopacity{0.900000}%
\pgfsetlinewidth{1.003750pt}%
\definecolor{currentstroke}{rgb}{0.580392,0.403922,0.741176}%
\pgfsetstrokecolor{currentstroke}%
\pgfsetstrokeopacity{0.900000}%
\pgfsetdash{}{0pt}%
\pgfsys@defobject{currentmarker}{\pgfqpoint{-0.058926in}{-0.058926in}}{\pgfqpoint{0.058926in}{0.058926in}}{%
\pgfpathmoveto{\pgfqpoint{0.000000in}{-0.058926in}}%
\pgfpathlineto{\pgfqpoint{0.058926in}{0.000000in}}%
\pgfpathlineto{\pgfqpoint{0.000000in}{0.058926in}}%
\pgfpathlineto{\pgfqpoint{-0.058926in}{0.000000in}}%
\pgfpathlineto{\pgfqpoint{0.000000in}{-0.058926in}}%
\pgfpathclose%
\pgfusepath{stroke,fill}%
}%
\begin{pgfscope}%
\pgfsys@transformshift{4.062661in}{3.363810in}%
\pgfsys@useobject{currentmarker}{}%
\end{pgfscope}%
\end{pgfscope}%
\begin{pgfscope}%
\definecolor{textcolor}{rgb}{0.000000,0.000000,0.000000}%
\pgfsetstrokecolor{textcolor}%
\pgfsetfillcolor{textcolor}%
\pgftext[x=4.437661in,y=3.276310in,left,base]{\color{textcolor}\rmfamily\fontsize{18.000000}{21.600000}\selectfont LPP\(\displaystyle ^{T}\)}%
\end{pgfscope}%
\begin{pgfscope}%
\pgfsetbuttcap%
\pgfsetroundjoin%
\pgfsetlinewidth{1.505625pt}%
\definecolor{currentstroke}{rgb}{0.580392,0.403922,0.741176}%
\pgfsetstrokecolor{currentstroke}%
\pgfsetstrokeopacity{0.800000}%
\pgfsetdash{{1.500000pt}{2.475000pt}}{0.000000pt}%
\pgfpathmoveto{\pgfqpoint{5.202692in}{3.363810in}}%
\pgfpathlineto{\pgfqpoint{5.452692in}{3.363810in}}%
\pgfpathlineto{\pgfqpoint{5.702692in}{3.363810in}}%
\pgfusepath{stroke}%
\end{pgfscope}%
\begin{pgfscope}%
\pgfsetbuttcap%
\pgfsetroundjoin%
\definecolor{currentfill}{rgb}{0.580392,0.403922,0.741176}%
\pgfsetfillcolor{currentfill}%
\pgfsetfillopacity{0.800000}%
\pgfsetlinewidth{1.003750pt}%
\definecolor{currentstroke}{rgb}{0.580392,0.403922,0.741176}%
\pgfsetstrokecolor{currentstroke}%
\pgfsetstrokeopacity{0.800000}%
\pgfsetdash{}{0pt}%
\pgfsys@defobject{currentmarker}{\pgfqpoint{-0.041667in}{-0.041667in}}{\pgfqpoint{0.041667in}{0.041667in}}{%
\pgfpathmoveto{\pgfqpoint{0.000000in}{-0.041667in}}%
\pgfpathcurveto{\pgfqpoint{0.011050in}{-0.041667in}}{\pgfqpoint{0.021649in}{-0.037276in}}{\pgfqpoint{0.029463in}{-0.029463in}}%
\pgfpathcurveto{\pgfqpoint{0.037276in}{-0.021649in}}{\pgfqpoint{0.041667in}{-0.011050in}}{\pgfqpoint{0.041667in}{0.000000in}}%
\pgfpathcurveto{\pgfqpoint{0.041667in}{0.011050in}}{\pgfqpoint{0.037276in}{0.021649in}}{\pgfqpoint{0.029463in}{0.029463in}}%
\pgfpathcurveto{\pgfqpoint{0.021649in}{0.037276in}}{\pgfqpoint{0.011050in}{0.041667in}}{\pgfqpoint{0.000000in}{0.041667in}}%
\pgfpathcurveto{\pgfqpoint{-0.011050in}{0.041667in}}{\pgfqpoint{-0.021649in}{0.037276in}}{\pgfqpoint{-0.029463in}{0.029463in}}%
\pgfpathcurveto{\pgfqpoint{-0.037276in}{0.021649in}}{\pgfqpoint{-0.041667in}{0.011050in}}{\pgfqpoint{-0.041667in}{0.000000in}}%
\pgfpathcurveto{\pgfqpoint{-0.041667in}{-0.011050in}}{\pgfqpoint{-0.037276in}{-0.021649in}}{\pgfqpoint{-0.029463in}{-0.029463in}}%
\pgfpathcurveto{\pgfqpoint{-0.021649in}{-0.037276in}}{\pgfqpoint{-0.011050in}{-0.041667in}}{\pgfqpoint{0.000000in}{-0.041667in}}%
\pgfpathlineto{\pgfqpoint{0.000000in}{-0.041667in}}%
\pgfpathclose%
\pgfusepath{stroke,fill}%
}%
\begin{pgfscope}%
\pgfsys@transformshift{5.452692in}{3.363810in}%
\pgfsys@useobject{currentmarker}{}%
\end{pgfscope}%
\end{pgfscope}%
\begin{pgfscope}%
\definecolor{textcolor}{rgb}{0.000000,0.000000,0.000000}%
\pgfsetstrokecolor{textcolor}%
\pgfsetfillcolor{textcolor}%
\pgftext[x=5.827692in,y=3.276310in,left,base]{\color{textcolor}\rmfamily\fontsize{18.000000}{21.600000}\selectfont LPP\(\displaystyle ^{T, P}\)}%
\end{pgfscope}%
\end{pgfpicture}%
\makeatother%
\endgroup%

%% file: figures/DemandPlots/TightLayout/DE_OC_Gap10pct.pgf
\begingroup%
\makeatletter%
\begin{pgfpicture}%
\pgfpathrectangle{\pgfpointorigin}{\pgfqpoint{7.780338in}{3.651310in}}%
\pgfusepath{use as bounding box, clip}%
\begin{pgfscope}%
\pgfsetbuttcap%
\pgfsetmiterjoin%
\definecolor{currentfill}{rgb}{1.000000,1.000000,1.000000}%
\pgfsetfillcolor{currentfill}%
\pgfsetlinewidth{0.000000pt}%
\definecolor{currentstroke}{rgb}{1.000000,1.000000,1.000000}%
\pgfsetstrokecolor{currentstroke}%
\pgfsetdash{}{0pt}%
\pgfpathmoveto{\pgfqpoint{0.000000in}{0.000000in}}%
\pgfpathlineto{\pgfqpoint{7.780338in}{0.000000in}}%
\pgfpathlineto{\pgfqpoint{7.780338in}{3.651310in}}%
\pgfpathlineto{\pgfqpoint{0.000000in}{3.651310in}}%
\pgfpathlineto{\pgfqpoint{0.000000in}{0.000000in}}%
\pgfpathclose%
\pgfusepath{fill}%
\end{pgfscope}%
\begin{pgfscope}%
\pgfsetbuttcap%
\pgfsetmiterjoin%
\definecolor{currentfill}{rgb}{1.000000,1.000000,1.000000}%
\pgfsetfillcolor{currentfill}%
\pgfsetlinewidth{0.000000pt}%
\definecolor{currentstroke}{rgb}{0.000000,0.000000,0.000000}%
\pgfsetstrokecolor{currentstroke}%
\pgfsetstrokeopacity{0.000000}%
\pgfsetdash{}{0pt}%
\pgfpathmoveto{\pgfqpoint{1.016468in}{0.679475in}}%
\pgfpathlineto{\pgfqpoint{7.680338in}{0.679475in}}%
\pgfpathlineto{\pgfqpoint{7.680338in}{3.135897in}}%
\pgfpathlineto{\pgfqpoint{1.016468in}{3.135897in}}%
\pgfpathlineto{\pgfqpoint{1.016468in}{0.679475in}}%
\pgfpathclose%
\pgfusepath{fill}%
\end{pgfscope}%
\begin{pgfscope}%
\pgfsetbuttcap%
\pgfsetroundjoin%
\definecolor{currentfill}{rgb}{0.000000,0.000000,0.000000}%
\pgfsetfillcolor{currentfill}%
\pgfsetlinewidth{0.803000pt}%
\definecolor{currentstroke}{rgb}{0.000000,0.000000,0.000000}%
\pgfsetstrokecolor{currentstroke}%
\pgfsetdash{}{0pt}%
\pgfsys@defobject{currentmarker}{\pgfqpoint{0.000000in}{-0.048611in}}{\pgfqpoint{0.000000in}{0.000000in}}{%
\pgfpathmoveto{\pgfqpoint{0.000000in}{0.000000in}}%
\pgfpathlineto{\pgfqpoint{0.000000in}{-0.048611in}}%
\pgfusepath{stroke,fill}%
}%
\begin{pgfscope}%
\pgfsys@transformshift{1.016468in}{0.679475in}%
\pgfsys@useobject{currentmarker}{}%
\end{pgfscope}%
\end{pgfscope}%
\begin{pgfscope}%
\definecolor{textcolor}{rgb}{0.000000,0.000000,0.000000}%
\pgfsetstrokecolor{textcolor}%
\pgfsetfillcolor{textcolor}%
\pgftext[x=1.016468in,y=0.582253in,,top]{\color{textcolor}\rmfamily\fontsize{16.000000}{19.200000}\selectfont \(\displaystyle {0.0}\)}%
\end{pgfscope}%
\begin{pgfscope}%
\pgfsetbuttcap%
\pgfsetroundjoin%
\definecolor{currentfill}{rgb}{0.000000,0.000000,0.000000}%
\pgfsetfillcolor{currentfill}%
\pgfsetlinewidth{0.803000pt}%
\definecolor{currentstroke}{rgb}{0.000000,0.000000,0.000000}%
\pgfsetstrokecolor{currentstroke}%
\pgfsetdash{}{0pt}%
\pgfsys@defobject{currentmarker}{\pgfqpoint{0.000000in}{-0.048611in}}{\pgfqpoint{0.000000in}{0.000000in}}{%
\pgfpathmoveto{\pgfqpoint{0.000000in}{0.000000in}}%
\pgfpathlineto{\pgfqpoint{0.000000in}{-0.048611in}}%
\pgfusepath{stroke,fill}%
}%
\begin{pgfscope}%
\pgfsys@transformshift{2.228080in}{0.679475in}%
\pgfsys@useobject{currentmarker}{}%
\end{pgfscope}%
\end{pgfscope}%
\begin{pgfscope}%
\definecolor{textcolor}{rgb}{0.000000,0.000000,0.000000}%
\pgfsetstrokecolor{textcolor}%
\pgfsetfillcolor{textcolor}%
\pgftext[x=2.228080in,y=0.582253in,,top]{\color{textcolor}\rmfamily\fontsize{16.000000}{19.200000}\selectfont \(\displaystyle {0.2}\)}%
\end{pgfscope}%
\begin{pgfscope}%
\pgfsetbuttcap%
\pgfsetroundjoin%
\definecolor{currentfill}{rgb}{0.000000,0.000000,0.000000}%
\pgfsetfillcolor{currentfill}%
\pgfsetlinewidth{0.803000pt}%
\definecolor{currentstroke}{rgb}{0.000000,0.000000,0.000000}%
\pgfsetstrokecolor{currentstroke}%
\pgfsetdash{}{0pt}%
\pgfsys@defobject{currentmarker}{\pgfqpoint{0.000000in}{-0.048611in}}{\pgfqpoint{0.000000in}{0.000000in}}{%
\pgfpathmoveto{\pgfqpoint{0.000000in}{0.000000in}}%
\pgfpathlineto{\pgfqpoint{0.000000in}{-0.048611in}}%
\pgfusepath{stroke,fill}%
}%
\begin{pgfscope}%
\pgfsys@transformshift{3.439693in}{0.679475in}%
\pgfsys@useobject{currentmarker}{}%
\end{pgfscope}%
\end{pgfscope}%
\begin{pgfscope}%
\definecolor{textcolor}{rgb}{0.000000,0.000000,0.000000}%
\pgfsetstrokecolor{textcolor}%
\pgfsetfillcolor{textcolor}%
\pgftext[x=3.439693in,y=0.582253in,,top]{\color{textcolor}\rmfamily\fontsize{16.000000}{19.200000}\selectfont \(\displaystyle {0.4}\)}%
\end{pgfscope}%
\begin{pgfscope}%
\pgfsetbuttcap%
\pgfsetroundjoin%
\definecolor{currentfill}{rgb}{0.000000,0.000000,0.000000}%
\pgfsetfillcolor{currentfill}%
\pgfsetlinewidth{0.803000pt}%
\definecolor{currentstroke}{rgb}{0.000000,0.000000,0.000000}%
\pgfsetstrokecolor{currentstroke}%
\pgfsetdash{}{0pt}%
\pgfsys@defobject{currentmarker}{\pgfqpoint{0.000000in}{-0.048611in}}{\pgfqpoint{0.000000in}{0.000000in}}{%
\pgfpathmoveto{\pgfqpoint{0.000000in}{0.000000in}}%
\pgfpathlineto{\pgfqpoint{0.000000in}{-0.048611in}}%
\pgfusepath{stroke,fill}%
}%
\begin{pgfscope}%
\pgfsys@transformshift{4.651306in}{0.679475in}%
\pgfsys@useobject{currentmarker}{}%
\end{pgfscope}%
\end{pgfscope}%
\begin{pgfscope}%
\definecolor{textcolor}{rgb}{0.000000,0.000000,0.000000}%
\pgfsetstrokecolor{textcolor}%
\pgfsetfillcolor{textcolor}%
\pgftext[x=4.651306in,y=0.582253in,,top]{\color{textcolor}\rmfamily\fontsize{16.000000}{19.200000}\selectfont \(\displaystyle {0.6}\)}%
\end{pgfscope}%
\begin{pgfscope}%
\pgfsetbuttcap%
\pgfsetroundjoin%
\definecolor{currentfill}{rgb}{0.000000,0.000000,0.000000}%
\pgfsetfillcolor{currentfill}%
\pgfsetlinewidth{0.803000pt}%
\definecolor{currentstroke}{rgb}{0.000000,0.000000,0.000000}%
\pgfsetstrokecolor{currentstroke}%
\pgfsetdash{}{0pt}%
\pgfsys@defobject{currentmarker}{\pgfqpoint{0.000000in}{-0.048611in}}{\pgfqpoint{0.000000in}{0.000000in}}{%
\pgfpathmoveto{\pgfqpoint{0.000000in}{0.000000in}}%
\pgfpathlineto{\pgfqpoint{0.000000in}{-0.048611in}}%
\pgfusepath{stroke,fill}%
}%
\begin{pgfscope}%
\pgfsys@transformshift{5.862919in}{0.679475in}%
\pgfsys@useobject{currentmarker}{}%
\end{pgfscope}%
\end{pgfscope}%
\begin{pgfscope}%
\definecolor{textcolor}{rgb}{0.000000,0.000000,0.000000}%
\pgfsetstrokecolor{textcolor}%
\pgfsetfillcolor{textcolor}%
\pgftext[x=5.862919in,y=0.582253in,,top]{\color{textcolor}\rmfamily\fontsize{16.000000}{19.200000}\selectfont \(\displaystyle {0.8}\)}%
\end{pgfscope}%
\begin{pgfscope}%
\pgfsetbuttcap%
\pgfsetroundjoin%
\definecolor{currentfill}{rgb}{0.000000,0.000000,0.000000}%
\pgfsetfillcolor{currentfill}%
\pgfsetlinewidth{0.803000pt}%
\definecolor{currentstroke}{rgb}{0.000000,0.000000,0.000000}%
\pgfsetstrokecolor{currentstroke}%
\pgfsetdash{}{0pt}%
\pgfsys@defobject{currentmarker}{\pgfqpoint{0.000000in}{-0.048611in}}{\pgfqpoint{0.000000in}{0.000000in}}{%
\pgfpathmoveto{\pgfqpoint{0.000000in}{0.000000in}}%
\pgfpathlineto{\pgfqpoint{0.000000in}{-0.048611in}}%
\pgfusepath{stroke,fill}%
}%
\begin{pgfscope}%
\pgfsys@transformshift{7.074532in}{0.679475in}%
\pgfsys@useobject{currentmarker}{}%
\end{pgfscope}%
\end{pgfscope}%
\begin{pgfscope}%
\definecolor{textcolor}{rgb}{0.000000,0.000000,0.000000}%
\pgfsetstrokecolor{textcolor}%
\pgfsetfillcolor{textcolor}%
\pgftext[x=7.074532in,y=0.582253in,,top]{\color{textcolor}\rmfamily\fontsize{16.000000}{19.200000}\selectfont \(\displaystyle {1.0}\)}%
\end{pgfscope}%
\begin{pgfscope}%
\definecolor{textcolor}{rgb}{0.000000,0.000000,0.000000}%
\pgfsetstrokecolor{textcolor}%
\pgfsetfillcolor{textcolor}%
\pgftext[x=4.348403in,y=0.313349in,,top]{\color{textcolor}\rmfamily\fontsize{18.000000}{21.600000}\selectfont \(\displaystyle \lambda\)}%
\end{pgfscope}%
\begin{pgfscope}%
\pgfsetbuttcap%
\pgfsetroundjoin%
\definecolor{currentfill}{rgb}{0.000000,0.000000,0.000000}%
\pgfsetfillcolor{currentfill}%
\pgfsetlinewidth{0.803000pt}%
\definecolor{currentstroke}{rgb}{0.000000,0.000000,0.000000}%
\pgfsetstrokecolor{currentstroke}%
\pgfsetdash{}{0pt}%
\pgfsys@defobject{currentmarker}{\pgfqpoint{-0.048611in}{0.000000in}}{\pgfqpoint{-0.000000in}{0.000000in}}{%
\pgfpathmoveto{\pgfqpoint{-0.000000in}{0.000000in}}%
\pgfpathlineto{\pgfqpoint{-0.048611in}{0.000000in}}%
\pgfusepath{stroke,fill}%
}%
\begin{pgfscope}%
\pgfsys@transformshift{1.016468in}{1.005206in}%
\pgfsys@useobject{currentmarker}{}%
\end{pgfscope}%
\end{pgfscope}%
\begin{pgfscope}%
\definecolor{textcolor}{rgb}{0.000000,0.000000,0.000000}%
\pgfsetstrokecolor{textcolor}%
\pgfsetfillcolor{textcolor}%
\pgftext[x=0.368904in, y=0.921873in, left, base]{\color{textcolor}\rmfamily\fontsize{16.000000}{19.200000}\selectfont \(\displaystyle {20000}\)}%
\end{pgfscope}%
\begin{pgfscope}%
\pgfsetbuttcap%
\pgfsetroundjoin%
\definecolor{currentfill}{rgb}{0.000000,0.000000,0.000000}%
\pgfsetfillcolor{currentfill}%
\pgfsetlinewidth{0.803000pt}%
\definecolor{currentstroke}{rgb}{0.000000,0.000000,0.000000}%
\pgfsetstrokecolor{currentstroke}%
\pgfsetdash{}{0pt}%
\pgfsys@defobject{currentmarker}{\pgfqpoint{-0.048611in}{0.000000in}}{\pgfqpoint{-0.000000in}{0.000000in}}{%
\pgfpathmoveto{\pgfqpoint{-0.000000in}{0.000000in}}%
\pgfpathlineto{\pgfqpoint{-0.048611in}{0.000000in}}%
\pgfusepath{stroke,fill}%
}%
\begin{pgfscope}%
\pgfsys@transformshift{1.016468in}{1.643603in}%
\pgfsys@useobject{currentmarker}{}%
\end{pgfscope}%
\end{pgfscope}%
\begin{pgfscope}%
\definecolor{textcolor}{rgb}{0.000000,0.000000,0.000000}%
\pgfsetstrokecolor{textcolor}%
\pgfsetfillcolor{textcolor}%
\pgftext[x=0.368904in, y=1.560269in, left, base]{\color{textcolor}\rmfamily\fontsize{16.000000}{19.200000}\selectfont \(\displaystyle {30000}\)}%
\end{pgfscope}%
\begin{pgfscope}%
\pgfsetbuttcap%
\pgfsetroundjoin%
\definecolor{currentfill}{rgb}{0.000000,0.000000,0.000000}%
\pgfsetfillcolor{currentfill}%
\pgfsetlinewidth{0.803000pt}%
\definecolor{currentstroke}{rgb}{0.000000,0.000000,0.000000}%
\pgfsetstrokecolor{currentstroke}%
\pgfsetdash{}{0pt}%
\pgfsys@defobject{currentmarker}{\pgfqpoint{-0.048611in}{0.000000in}}{\pgfqpoint{-0.000000in}{0.000000in}}{%
\pgfpathmoveto{\pgfqpoint{-0.000000in}{0.000000in}}%
\pgfpathlineto{\pgfqpoint{-0.048611in}{0.000000in}}%
\pgfusepath{stroke,fill}%
}%
\begin{pgfscope}%
\pgfsys@transformshift{1.016468in}{2.281999in}%
\pgfsys@useobject{currentmarker}{}%
\end{pgfscope}%
\end{pgfscope}%
\begin{pgfscope}%
\definecolor{textcolor}{rgb}{0.000000,0.000000,0.000000}%
\pgfsetstrokecolor{textcolor}%
\pgfsetfillcolor{textcolor}%
\pgftext[x=0.368904in, y=2.198666in, left, base]{\color{textcolor}\rmfamily\fontsize{16.000000}{19.200000}\selectfont \(\displaystyle {40000}\)}%
\end{pgfscope}%
\begin{pgfscope}%
\pgfsetbuttcap%
\pgfsetroundjoin%
\definecolor{currentfill}{rgb}{0.000000,0.000000,0.000000}%
\pgfsetfillcolor{currentfill}%
\pgfsetlinewidth{0.803000pt}%
\definecolor{currentstroke}{rgb}{0.000000,0.000000,0.000000}%
\pgfsetstrokecolor{currentstroke}%
\pgfsetdash{}{0pt}%
\pgfsys@defobject{currentmarker}{\pgfqpoint{-0.048611in}{0.000000in}}{\pgfqpoint{-0.000000in}{0.000000in}}{%
\pgfpathmoveto{\pgfqpoint{-0.000000in}{0.000000in}}%
\pgfpathlineto{\pgfqpoint{-0.048611in}{0.000000in}}%
\pgfusepath{stroke,fill}%
}%
\begin{pgfscope}%
\pgfsys@transformshift{1.016468in}{2.920395in}%
\pgfsys@useobject{currentmarker}{}%
\end{pgfscope}%
\end{pgfscope}%
\begin{pgfscope}%
\definecolor{textcolor}{rgb}{0.000000,0.000000,0.000000}%
\pgfsetstrokecolor{textcolor}%
\pgfsetfillcolor{textcolor}%
\pgftext[x=0.368904in, y=2.837062in, left, base]{\color{textcolor}\rmfamily\fontsize{16.000000}{19.200000}\selectfont \(\displaystyle {50000}\)}%
\end{pgfscope}%
\begin{pgfscope}%
\definecolor{textcolor}{rgb}{0.000000,0.000000,0.000000}%
\pgfsetstrokecolor{textcolor}%
\pgfsetfillcolor{textcolor}%
\pgftext[x=0.313349in,y=1.907686in,,bottom,rotate=90.000000]{\color{textcolor}\rmfamily\fontsize{18.000000}{21.600000}\selectfont DKK}%
\end{pgfscope}%
\begin{pgfscope}%
\pgfpathrectangle{\pgfqpoint{1.016468in}{0.679475in}}{\pgfqpoint{6.663871in}{2.456422in}}%
\pgfusepath{clip}%
\pgfsetrectcap%
\pgfsetroundjoin%
\pgfsetlinewidth{1.505625pt}%
\definecolor{currentstroke}{rgb}{0.839216,0.152941,0.156863}%
\pgfsetstrokecolor{currentstroke}%
\pgfsetdash{}{0pt}%
\pgfpathmoveto{\pgfqpoint{1.622274in}{2.259885in}}%
\pgfpathlineto{\pgfqpoint{2.530984in}{2.292597in}}%
\pgfpathlineto{\pgfqpoint{4.045500in}{2.438389in}}%
\pgfpathlineto{\pgfqpoint{5.560016in}{2.701213in}}%
\pgfpathlineto{\pgfqpoint{7.074532in}{3.024241in}}%
\pgfusepath{stroke}%
\end{pgfscope}%
\begin{pgfscope}%
\pgfpathrectangle{\pgfqpoint{1.016468in}{0.679475in}}{\pgfqpoint{6.663871in}{2.456422in}}%
\pgfusepath{clip}%
\pgfsetbuttcap%
\pgfsetmiterjoin%
\definecolor{currentfill}{rgb}{0.839216,0.152941,0.156863}%
\pgfsetfillcolor{currentfill}%
\pgfsetlinewidth{0.752812pt}%
\definecolor{currentstroke}{rgb}{1.000000,1.000000,1.000000}%
\pgfsetstrokecolor{currentstroke}%
\pgfsetdash{}{0pt}%
\pgfsys@defobject{currentmarker}{\pgfqpoint{-0.041667in}{-0.041667in}}{\pgfqpoint{0.041667in}{0.041667in}}{%
\pgfpathmoveto{\pgfqpoint{-0.041667in}{-0.041667in}}%
\pgfpathlineto{\pgfqpoint{0.041667in}{-0.041667in}}%
\pgfpathlineto{\pgfqpoint{0.041667in}{0.041667in}}%
\pgfpathlineto{\pgfqpoint{-0.041667in}{0.041667in}}%
\pgfpathlineto{\pgfqpoint{-0.041667in}{-0.041667in}}%
\pgfpathclose%
\pgfusepath{stroke,fill}%
}%
\begin{pgfscope}%
\pgfsys@transformshift{1.622274in}{2.259885in}%
\pgfsys@useobject{currentmarker}{}%
\end{pgfscope}%
\begin{pgfscope}%
\pgfsys@transformshift{2.530984in}{2.292597in}%
\pgfsys@useobject{currentmarker}{}%
\end{pgfscope}%
\begin{pgfscope}%
\pgfsys@transformshift{4.045500in}{2.438389in}%
\pgfsys@useobject{currentmarker}{}%
\end{pgfscope}%
\begin{pgfscope}%
\pgfsys@transformshift{5.560016in}{2.701213in}%
\pgfsys@useobject{currentmarker}{}%
\end{pgfscope}%
\begin{pgfscope}%
\pgfsys@transformshift{7.074532in}{3.024241in}%
\pgfsys@useobject{currentmarker}{}%
\end{pgfscope}%
\end{pgfscope}%
\begin{pgfscope}%
\pgfpathrectangle{\pgfqpoint{1.016468in}{0.679475in}}{\pgfqpoint{6.663871in}{2.456422in}}%
\pgfusepath{clip}%
\pgfsetbuttcap%
\pgfsetroundjoin%
\pgfsetlinewidth{1.505625pt}%
\definecolor{currentstroke}{rgb}{0.839216,0.152941,0.156863}%
\pgfsetstrokecolor{currentstroke}%
\pgfsetstrokeopacity{0.900000}%
\pgfsetdash{{7.500000pt}{7.500000pt}}{0.000000pt}%
\pgfpathmoveto{\pgfqpoint{1.622274in}{0.791131in}}%
\pgfpathlineto{\pgfqpoint{2.530984in}{0.941521in}}%
\pgfpathlineto{\pgfqpoint{4.045500in}{1.171649in}}%
\pgfpathlineto{\pgfqpoint{5.560016in}{1.493367in}}%
\pgfpathlineto{\pgfqpoint{7.074532in}{1.825758in}}%
\pgfusepath{stroke}%
\end{pgfscope}%
\begin{pgfscope}%
\pgfpathrectangle{\pgfqpoint{1.016468in}{0.679475in}}{\pgfqpoint{6.663871in}{2.456422in}}%
\pgfusepath{clip}%
\pgfsetbuttcap%
\pgfsetmiterjoin%
\definecolor{currentfill}{rgb}{0.839216,0.152941,0.156863}%
\pgfsetfillcolor{currentfill}%
\pgfsetfillopacity{0.900000}%
\pgfsetlinewidth{0.752812pt}%
\definecolor{currentstroke}{rgb}{1.000000,1.000000,1.000000}%
\pgfsetstrokecolor{currentstroke}%
\pgfsetdash{}{0pt}%
\pgfsys@defobject{currentmarker}{\pgfqpoint{-0.058926in}{-0.058926in}}{\pgfqpoint{0.058926in}{0.058926in}}{%
\pgfpathmoveto{\pgfqpoint{0.000000in}{-0.058926in}}%
\pgfpathlineto{\pgfqpoint{0.058926in}{0.000000in}}%
\pgfpathlineto{\pgfqpoint{0.000000in}{0.058926in}}%
\pgfpathlineto{\pgfqpoint{-0.058926in}{0.000000in}}%
\pgfpathlineto{\pgfqpoint{0.000000in}{-0.058926in}}%
\pgfpathclose%
\pgfusepath{stroke,fill}%
}%
\begin{pgfscope}%
\pgfsys@transformshift{1.622274in}{0.791131in}%
\pgfsys@useobject{currentmarker}{}%
\end{pgfscope}%
\begin{pgfscope}%
\pgfsys@transformshift{2.530984in}{0.941521in}%
\pgfsys@useobject{currentmarker}{}%
\end{pgfscope}%
\begin{pgfscope}%
\pgfsys@transformshift{4.045500in}{1.171649in}%
\pgfsys@useobject{currentmarker}{}%
\end{pgfscope}%
\begin{pgfscope}%
\pgfsys@transformshift{5.560016in}{1.493367in}%
\pgfsys@useobject{currentmarker}{}%
\end{pgfscope}%
\begin{pgfscope}%
\pgfsys@transformshift{7.074532in}{1.825758in}%
\pgfsys@useobject{currentmarker}{}%
\end{pgfscope}%
\end{pgfscope}%
\begin{pgfscope}%
\pgfpathrectangle{\pgfqpoint{1.016468in}{0.679475in}}{\pgfqpoint{6.663871in}{2.456422in}}%
\pgfusepath{clip}%
\pgfsetbuttcap%
\pgfsetroundjoin%
\pgfsetlinewidth{1.505625pt}%
\definecolor{currentstroke}{rgb}{0.839216,0.152941,0.156863}%
\pgfsetstrokecolor{currentstroke}%
\pgfsetstrokeopacity{0.800000}%
\pgfsetdash{{1.500000pt}{4.500000pt}}{0.000000pt}%
\pgfpathmoveto{\pgfqpoint{1.622274in}{2.259884in}}%
\pgfpathlineto{\pgfqpoint{2.530984in}{2.292596in}}%
\pgfpathlineto{\pgfqpoint{4.045500in}{2.429022in}}%
\pgfpathlineto{\pgfqpoint{5.560016in}{2.701213in}}%
\pgfpathlineto{\pgfqpoint{7.074532in}{3.014878in}}%
\pgfusepath{stroke}%
\end{pgfscope}%
\begin{pgfscope}%
\pgfpathrectangle{\pgfqpoint{1.016468in}{0.679475in}}{\pgfqpoint{6.663871in}{2.456422in}}%
\pgfusepath{clip}%
\pgfsetbuttcap%
\pgfsetroundjoin%
\definecolor{currentfill}{rgb}{0.839216,0.152941,0.156863}%
\pgfsetfillcolor{currentfill}%
\pgfsetfillopacity{0.800000}%
\pgfsetlinewidth{0.752812pt}%
\definecolor{currentstroke}{rgb}{1.000000,1.000000,1.000000}%
\pgfsetstrokecolor{currentstroke}%
\pgfsetdash{}{0pt}%
\pgfsys@defobject{currentmarker}{\pgfqpoint{-0.041667in}{-0.041667in}}{\pgfqpoint{0.041667in}{0.041667in}}{%
\pgfpathmoveto{\pgfqpoint{0.000000in}{-0.041667in}}%
\pgfpathcurveto{\pgfqpoint{0.011050in}{-0.041667in}}{\pgfqpoint{0.021649in}{-0.037276in}}{\pgfqpoint{0.029463in}{-0.029463in}}%
\pgfpathcurveto{\pgfqpoint{0.037276in}{-0.021649in}}{\pgfqpoint{0.041667in}{-0.011050in}}{\pgfqpoint{0.041667in}{0.000000in}}%
\pgfpathcurveto{\pgfqpoint{0.041667in}{0.011050in}}{\pgfqpoint{0.037276in}{0.021649in}}{\pgfqpoint{0.029463in}{0.029463in}}%
\pgfpathcurveto{\pgfqpoint{0.021649in}{0.037276in}}{\pgfqpoint{0.011050in}{0.041667in}}{\pgfqpoint{0.000000in}{0.041667in}}%
\pgfpathcurveto{\pgfqpoint{-0.011050in}{0.041667in}}{\pgfqpoint{-0.021649in}{0.037276in}}{\pgfqpoint{-0.029463in}{0.029463in}}%
\pgfpathcurveto{\pgfqpoint{-0.037276in}{0.021649in}}{\pgfqpoint{-0.041667in}{0.011050in}}{\pgfqpoint{-0.041667in}{0.000000in}}%
\pgfpathcurveto{\pgfqpoint{-0.041667in}{-0.011050in}}{\pgfqpoint{-0.037276in}{-0.021649in}}{\pgfqpoint{-0.029463in}{-0.029463in}}%
\pgfpathcurveto{\pgfqpoint{-0.021649in}{-0.037276in}}{\pgfqpoint{-0.011050in}{-0.041667in}}{\pgfqpoint{0.000000in}{-0.041667in}}%
\pgfpathlineto{\pgfqpoint{0.000000in}{-0.041667in}}%
\pgfpathclose%
\pgfusepath{stroke,fill}%
}%
\begin{pgfscope}%
\pgfsys@transformshift{1.622274in}{2.259884in}%
\pgfsys@useobject{currentmarker}{}%
\end{pgfscope}%
\begin{pgfscope}%
\pgfsys@transformshift{2.530984in}{2.292596in}%
\pgfsys@useobject{currentmarker}{}%
\end{pgfscope}%
\begin{pgfscope}%
\pgfsys@transformshift{4.045500in}{2.429022in}%
\pgfsys@useobject{currentmarker}{}%
\end{pgfscope}%
\begin{pgfscope}%
\pgfsys@transformshift{5.560016in}{2.701213in}%
\pgfsys@useobject{currentmarker}{}%
\end{pgfscope}%
\begin{pgfscope}%
\pgfsys@transformshift{7.074532in}{3.014878in}%
\pgfsys@useobject{currentmarker}{}%
\end{pgfscope}%
\end{pgfscope}%
\begin{pgfscope}%
\pgfsetrectcap%
\pgfsetmiterjoin%
\pgfsetlinewidth{0.803000pt}%
\definecolor{currentstroke}{rgb}{0.000000,0.000000,0.000000}%
\pgfsetstrokecolor{currentstroke}%
\pgfsetdash{}{0pt}%
\pgfpathmoveto{\pgfqpoint{1.016468in}{0.679475in}}%
\pgfpathlineto{\pgfqpoint{1.016468in}{3.135897in}}%
\pgfusepath{stroke}%
\end{pgfscope}%
\begin{pgfscope}%
\pgfsetrectcap%
\pgfsetmiterjoin%
\pgfsetlinewidth{0.803000pt}%
\definecolor{currentstroke}{rgb}{0.000000,0.000000,0.000000}%
\pgfsetstrokecolor{currentstroke}%
\pgfsetdash{}{0pt}%
\pgfpathmoveto{\pgfqpoint{7.680338in}{0.679475in}}%
\pgfpathlineto{\pgfqpoint{7.680338in}{3.135897in}}%
\pgfusepath{stroke}%
\end{pgfscope}%
\begin{pgfscope}%
\pgfsetrectcap%
\pgfsetmiterjoin%
\pgfsetlinewidth{0.803000pt}%
\definecolor{currentstroke}{rgb}{0.000000,0.000000,0.000000}%
\pgfsetstrokecolor{currentstroke}%
\pgfsetdash{}{0pt}%
\pgfpathmoveto{\pgfqpoint{1.016468in}{0.679475in}}%
\pgfpathlineto{\pgfqpoint{7.680338in}{0.679475in}}%
\pgfusepath{stroke}%
\end{pgfscope}%
\begin{pgfscope}%
\pgfsetrectcap%
\pgfsetmiterjoin%
\pgfsetlinewidth{0.803000pt}%
\definecolor{currentstroke}{rgb}{0.000000,0.000000,0.000000}%
\pgfsetstrokecolor{currentstroke}%
\pgfsetdash{}{0pt}%
\pgfpathmoveto{\pgfqpoint{1.016468in}{3.135897in}}%
\pgfpathlineto{\pgfqpoint{7.680338in}{3.135897in}}%
\pgfusepath{stroke}%
\end{pgfscope}%
\begin{pgfscope}%
\pgfsetrectcap%
\pgfsetroundjoin%
\pgfsetlinewidth{1.505625pt}%
\definecolor{currentstroke}{rgb}{0.839216,0.152941,0.156863}%
\pgfsetstrokecolor{currentstroke}%
\pgfsetdash{}{0pt}%
\pgfpathmoveto{\pgfqpoint{2.195836in}{3.363810in}}%
\pgfpathlineto{\pgfqpoint{2.445836in}{3.363810in}}%
\pgfpathlineto{\pgfqpoint{2.695836in}{3.363810in}}%
\pgfusepath{stroke}%
\end{pgfscope}%
\begin{pgfscope}%
\pgfsetbuttcap%
\pgfsetmiterjoin%
\definecolor{currentfill}{rgb}{0.839216,0.152941,0.156863}%
\pgfsetfillcolor{currentfill}%
\pgfsetlinewidth{1.003750pt}%
\definecolor{currentstroke}{rgb}{0.839216,0.152941,0.156863}%
\pgfsetstrokecolor{currentstroke}%
\pgfsetdash{}{0pt}%
\pgfsys@defobject{currentmarker}{\pgfqpoint{-0.041667in}{-0.041667in}}{\pgfqpoint{0.041667in}{0.041667in}}{%
\pgfpathmoveto{\pgfqpoint{-0.041667in}{-0.041667in}}%
\pgfpathlineto{\pgfqpoint{0.041667in}{-0.041667in}}%
\pgfpathlineto{\pgfqpoint{0.041667in}{0.041667in}}%
\pgfpathlineto{\pgfqpoint{-0.041667in}{0.041667in}}%
\pgfpathlineto{\pgfqpoint{-0.041667in}{-0.041667in}}%
\pgfpathclose%
\pgfusepath{stroke,fill}%
}%
\begin{pgfscope}%
\pgfsys@transformshift{2.445836in}{3.363810in}%
\pgfsys@useobject{currentmarker}{}%
\end{pgfscope}%
\end{pgfscope}%
\begin{pgfscope}%
\definecolor{textcolor}{rgb}{0.000000,0.000000,0.000000}%
\pgfsetstrokecolor{textcolor}%
\pgfsetfillcolor{textcolor}%
\pgftext[x=2.820836in,y=3.276310in,left,base]{\color{textcolor}\rmfamily\fontsize{18.000000}{21.600000}\selectfont LPP\(\displaystyle ^{R,T}\)}%
\end{pgfscope}%
\begin{pgfscope}%
\pgfsetbuttcap%
\pgfsetroundjoin%
\pgfsetlinewidth{1.505625pt}%
\definecolor{currentstroke}{rgb}{0.839216,0.152941,0.156863}%
\pgfsetstrokecolor{currentstroke}%
\pgfsetstrokeopacity{0.900000}%
\pgfsetdash{{5.550000pt}{2.400000pt}}{0.000000pt}%
\pgfpathmoveto{\pgfqpoint{3.756787in}{3.363810in}}%
\pgfpathlineto{\pgfqpoint{4.006787in}{3.363810in}}%
\pgfpathlineto{\pgfqpoint{4.256787in}{3.363810in}}%
\pgfusepath{stroke}%
\end{pgfscope}%
\begin{pgfscope}%
\pgfsetbuttcap%
\pgfsetmiterjoin%
\definecolor{currentfill}{rgb}{0.839216,0.152941,0.156863}%
\pgfsetfillcolor{currentfill}%
\pgfsetfillopacity{0.900000}%
\pgfsetlinewidth{1.003750pt}%
\definecolor{currentstroke}{rgb}{0.839216,0.152941,0.156863}%
\pgfsetstrokecolor{currentstroke}%
\pgfsetstrokeopacity{0.900000}%
\pgfsetdash{}{0pt}%
\pgfsys@defobject{currentmarker}{\pgfqpoint{-0.058926in}{-0.058926in}}{\pgfqpoint{0.058926in}{0.058926in}}{%
\pgfpathmoveto{\pgfqpoint{0.000000in}{-0.058926in}}%
\pgfpathlineto{\pgfqpoint{0.058926in}{0.000000in}}%
\pgfpathlineto{\pgfqpoint{0.000000in}{0.058926in}}%
\pgfpathlineto{\pgfqpoint{-0.058926in}{0.000000in}}%
\pgfpathlineto{\pgfqpoint{0.000000in}{-0.058926in}}%
\pgfpathclose%
\pgfusepath{stroke,fill}%
}%
\begin{pgfscope}%
\pgfsys@transformshift{4.006787in}{3.363810in}%
\pgfsys@useobject{currentmarker}{}%
\end{pgfscope}%
\end{pgfscope}%
\begin{pgfscope}%
\definecolor{textcolor}{rgb}{0.000000,0.000000,0.000000}%
\pgfsetstrokecolor{textcolor}%
\pgfsetfillcolor{textcolor}%
\pgftext[x=4.381787in,y=3.276310in,left,base]{\color{textcolor}\rmfamily\fontsize{18.000000}{21.600000}\selectfont LPP\(\displaystyle ^{T}\)}%
\end{pgfscope}%
\begin{pgfscope}%
\pgfsetbuttcap%
\pgfsetroundjoin%
\pgfsetlinewidth{1.505625pt}%
\definecolor{currentstroke}{rgb}{0.839216,0.152941,0.156863}%
\pgfsetstrokecolor{currentstroke}%
\pgfsetstrokeopacity{0.800000}%
\pgfsetdash{{1.500000pt}{2.475000pt}}{0.000000pt}%
\pgfpathmoveto{\pgfqpoint{5.146819in}{3.363810in}}%
\pgfpathlineto{\pgfqpoint{5.396819in}{3.363810in}}%
\pgfpathlineto{\pgfqpoint{5.646819in}{3.363810in}}%
\pgfusepath{stroke}%
\end{pgfscope}%
\begin{pgfscope}%
\pgfsetbuttcap%
\pgfsetroundjoin%
\definecolor{currentfill}{rgb}{0.839216,0.152941,0.156863}%
\pgfsetfillcolor{currentfill}%
\pgfsetfillopacity{0.800000}%
\pgfsetlinewidth{1.003750pt}%
\definecolor{currentstroke}{rgb}{0.839216,0.152941,0.156863}%
\pgfsetstrokecolor{currentstroke}%
\pgfsetstrokeopacity{0.800000}%
\pgfsetdash{}{0pt}%
\pgfsys@defobject{currentmarker}{\pgfqpoint{-0.041667in}{-0.041667in}}{\pgfqpoint{0.041667in}{0.041667in}}{%
\pgfpathmoveto{\pgfqpoint{0.000000in}{-0.041667in}}%
\pgfpathcurveto{\pgfqpoint{0.011050in}{-0.041667in}}{\pgfqpoint{0.021649in}{-0.037276in}}{\pgfqpoint{0.029463in}{-0.029463in}}%
\pgfpathcurveto{\pgfqpoint{0.037276in}{-0.021649in}}{\pgfqpoint{0.041667in}{-0.011050in}}{\pgfqpoint{0.041667in}{0.000000in}}%
\pgfpathcurveto{\pgfqpoint{0.041667in}{0.011050in}}{\pgfqpoint{0.037276in}{0.021649in}}{\pgfqpoint{0.029463in}{0.029463in}}%
\pgfpathcurveto{\pgfqpoint{0.021649in}{0.037276in}}{\pgfqpoint{0.011050in}{0.041667in}}{\pgfqpoint{0.000000in}{0.041667in}}%
\pgfpathcurveto{\pgfqpoint{-0.011050in}{0.041667in}}{\pgfqpoint{-0.021649in}{0.037276in}}{\pgfqpoint{-0.029463in}{0.029463in}}%
\pgfpathcurveto{\pgfqpoint{-0.037276in}{0.021649in}}{\pgfqpoint{-0.041667in}{0.011050in}}{\pgfqpoint{-0.041667in}{0.000000in}}%
\pgfpathcurveto{\pgfqpoint{-0.041667in}{-0.011050in}}{\pgfqpoint{-0.037276in}{-0.021649in}}{\pgfqpoint{-0.029463in}{-0.029463in}}%
\pgfpathcurveto{\pgfqpoint{-0.021649in}{-0.037276in}}{\pgfqpoint{-0.011050in}{-0.041667in}}{\pgfqpoint{0.000000in}{-0.041667in}}%
\pgfpathlineto{\pgfqpoint{0.000000in}{-0.041667in}}%
\pgfpathclose%
\pgfusepath{stroke,fill}%
}%
\begin{pgfscope}%
\pgfsys@transformshift{5.396819in}{3.363810in}%
\pgfsys@useobject{currentmarker}{}%
\end{pgfscope}%
\end{pgfscope}%
\begin{pgfscope}%
\definecolor{textcolor}{rgb}{0.000000,0.000000,0.000000}%
\pgfsetstrokecolor{textcolor}%
\pgfsetfillcolor{textcolor}%
\pgftext[x=5.771819in,y=3.276310in,left,base]{\color{textcolor}\rmfamily\fontsize{18.000000}{21.600000}\selectfont LPP\(\displaystyle ^{T, P}\)}%
\end{pgfscope}%
\end{pgfpicture}%
\makeatother%
\endgroup%

%% file: Appendix/ValidInequalities.tex
\section{Valid Inequalities}\label{sec:AppendixValidIneq}
\begin{figure}[H]
\centering
\subfloat[Instances solved over time]{\includegraphics[width=0.45\textwidth]{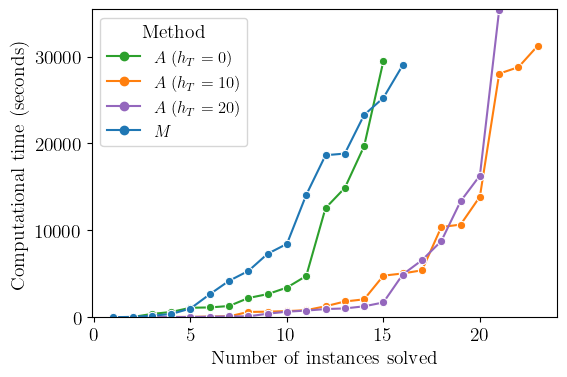}}
\subfloat[{Optimality gap at termination}]{\includegraphics[width=0.45\textwidth]{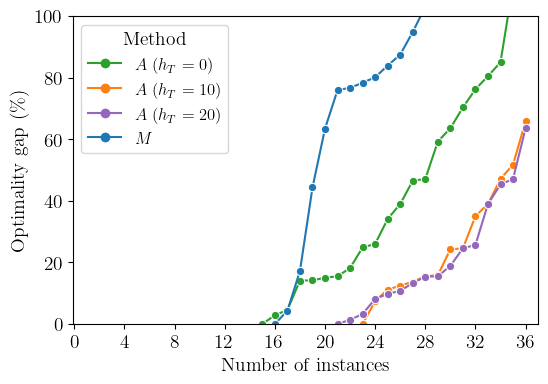}} \\
\caption{Computational performance ($\lambda = 0.25)$ }
\label{fig:SolvedAndGap_VI}
\end{figure}

\autoref{fig:SolvedAndGap_VI} presents the cumulative number of instances solved over time, along with the number of instances that reach a certain optimality gap at termination. The results compare the performance of \frameworkNameShortLP (\emph{A}) with varying threshold values for valid inequalities $h_T$, and CPLEX ($M$), all evaluated with $\lambda = 0.25$. 

%% file: Appendix/FrequencySet.tex
\section{Impact of the candidate frequency set} \label{sec:FreqAnalysis}
In the model, each line may operate at any integer headway between two and 20 minutes. Since line frequency directly influences the weighted travel time (WTT), the frequency options are important in modeling \elasticDemand. To ensure that the model appropriately balances service quality with operational cost, it is necessary to provide sufficient resolution in the available frequency options.

We evaluate four alternative headway sets defined by their step size:
$H_1=\{2,3,4,..,20\}$,\ $H_2=\{2,4,6,..,20\}$,\ $H_3=\{2,5,8,..,20\}$, and $H_4=\{2,8,14,20\}$
These sets span the same range but vary in granularity, from fine (small step size) to coarse (large step size). The analysis is based on 19 instances that are solvable to proven optimality for all headway sets.

\begin{table}[H]
\centering
\caption{Results averaged over 19 instances solved to optimality for {\problem{(\Pset)}} for $\lambda=0.25$ }
\label{tab:FreqSet}
\begin{tabular}{@{}lccc@{}}
\toprule
 & Objective (DKK) & WTT (minutes) & Demand (\%) \\
\midrule
$H_1$ & 91,814.1 & 25.1 & 50.7 \\
$H_2$ & 93,007.2 & 22.0 & 49.7 \\
$H_3$ & 94,009.2 & 22.6 & 48.2 \\
$H_4$ & 96,362.3 & 23.9 & 41.0 \\
\bottomrule
\end{tabular}
\end{table}

The results show that finer headway granularity generally yields better outcomes across all metrics. As the step size increases, the average objective value (i.e., total cost) rises, and the ridership declines. Interestingly, the minimum WTT is achieved with $H_2$, indicating that very fine-grained options ($H_1$) do not necessarily improve travel time but may contribute to more demand captured. 

%% file: Appendix/Capacity.tex
\section{Load profiles} \label{apx:Capacity}
\begin{figure}[htbp]
    \centering
    \begin{subfigure}{0.49\textwidth}
        \centering
        \begin{tikzpicture}
\begin{axis}[
    xlabel={$\lambda$},
    ylabel={Arcs operating at capacity (\%)},
    legend pos=north west,
    grid=both,
    width=\textwidth,
    height=6cm
]

\addplot[
    only marks,
    mark=*,
    mark size=1.5pt,
        opacity=0.5,
    color=blue, forget plot
] 
table [x=lambdaVal, y=capacityPctSelectedArcs, col sep=comma] {Appendix/routing_data/capacity_utilization_10pct_ALL.csv};

Line plot of average per lambda
\addplot[
    thick,
    color=red
]
table [x=lambdaVal, y=capacityPctSelectedArcs, col sep=comma] {Appendix/routing_data/AVG_capacity_utilization_10pct.csv};
\addlegendentry{Average}
\end{axis}
\end{tikzpicture}
\caption{Fraction of capacity-constrained arcs among arcs selected in the solution}
\label{fig:loads_selected}
\end{subfigure}%
\hfill
   \begin{subfigure}{0.49\textwidth}
        \centering
        \begin{tikzpicture}
\begin{axis}[
    xlabel={$\lambda$},
    ylabel={Arcs operating at capacity (\%)},
    legend pos=north west,
    grid=both,
    width=\textwidth,
    height=6cm
]

\addplot[
    only marks,
    mark=*,
    mark size=1.5pt,
    opacity=0.5,
    color=blue, forget plot
] 
table [x=lambdaVal, y=capacityPctAllArcs, col sep=comma] {Appendix/routing_data/capacity_utilization_10pct_ALL.csv};

Line plot of average per lambda
\addplot[
    thick,
    color=red
]
table [x=lambdaVal, y=capacityPctAllArcs, col sep=comma] {Appendix/routing_data/avg_capacity_utilization_10pct.csv};
\addlegendentry{Average}

\end{axis}
\end{tikzpicture}
\caption{Fraction of capacity-constrained arcs among all directed arcs in the PTN}
\label{fig:loads_all}
\end{subfigure}
\caption{Arcs operating at capacity (lines show means; points show individual instances). Only solutions
with an optimality gap below 10\% are included.}
\label{fig:loads}
\end{figure}
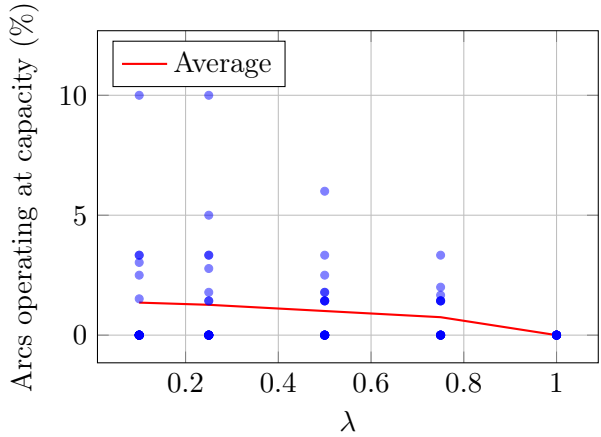
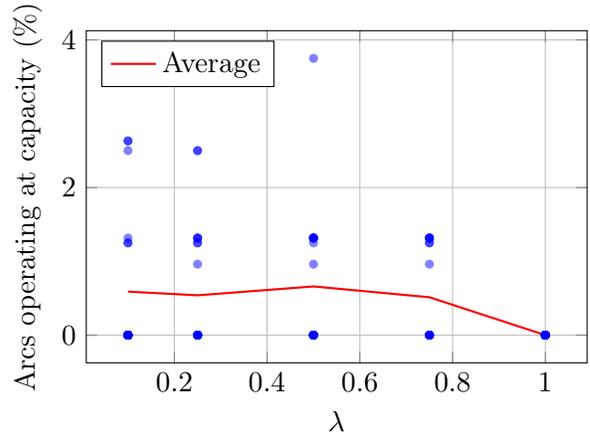

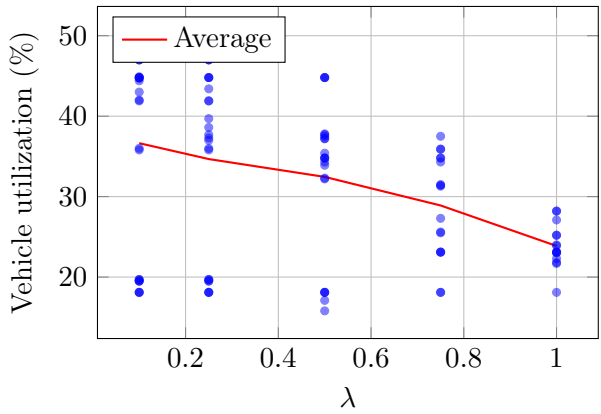
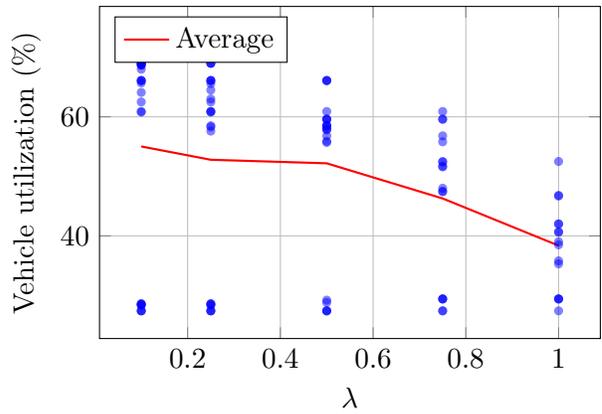
\begin{figure}[htbp]
    \centering
    \begin{subfigure}{0.49\textwidth}
        \centering
        \begin{tikzpicture}
\begin{axis}[
    xlabel={$\lambda$},
    ylabel={Vehicle utilization (\%)},
    legend pos=north west,
    grid=both,
    width=\textwidth,
    height=6cm
]

\addplot[
    only marks,
    mark=*,
    mark size=1.5pt,
    opacity=0.5,
    color=blue, forget plot
] 
table [x=lambdaVal, y=UTIL, col sep=comma] {Appendix/routing_data/capacity_utilization_10pct_ALL.csv};

Line plot of average per lambda
\addplot[
    thick,
    color=red
]
table [x=lambdaVal, y=UTIL, col sep=comma] {Appendix/routing_data/avg_capacity_utilization_10pct.csv};
\addlegendentry{Average}

\end{axis}
\end{tikzpicture}
\caption{Arc-level vehicle utilization}
\label{fig:util_arc}
\end{subfigure}%
\hfill
   \begin{subfigure}{0.49\textwidth}
        \centering
        \begin{tikzpicture}
\begin{axis}[
    xlabel={$\lambda$},
    ylabel={Vehicle utilization (\%)},
    legend pos=north west,
    grid=both,
    width=\textwidth,
    height=6cm
]

\addplot[
    only marks,
    mark=*,
    mark size=1.5pt,
    opacity=0.5,
    color=blue, forget plot
] 
table [x=lambdaVal, y=UTIL_WTTP, col sep=comma] {Appendix/routing_data/capacity_utilization_10pct_ALL.csv};

Line plot of average per lambda
\addplot[
    thick,
    color=red
]
table [x=lambdaVal, y=UTIL_WTTP, col sep=comma] {Appendix/routing_data/avg_capacity_utilization_10pct.csv};
\addlegendentry{Average}

\end{axis}
\end{tikzpicture}
\caption{Passenger- and time-weighted utilization}
\label{fig:util_wttp}
\end{subfigure}
\caption{Vehicle utilization rates across all solutions (lines show means, scatter points show individual instances). Only solutions with an optimality gap below 10\% are included.}
\label{fig:util}
\end{figure}

\added{\autoref{fig:loads} reports the fraction of arcs operating at maximum capacity, i.e., the assigned passenger flow to the segment of line equals the provided seating capacity. 
These figures have been computed based on those instances solved with an optimality gap below 10\%.
The number of at-capacity arcs ranges, on average, from 0.5 for $\lambda = 0.1$ ($1.4\%$ of the selected arcs) to $0$ at $\lambda = 1$. This trend reflects the structure of the model: as $\lambda$ increases, the objective increasingly favors operating more frequent services, thereby reducing the likelihood that arcs will reach the capacity limits. Similarly, \autoref{fig:util} shows that average vehicle utilization remains modest, ranging from 36\% to 24\% when measured at the arc level, and from 55\% to 38\% when measured as the “experienced congestion level”, i.e., utilization weighted by passenger-minutes.}

